\documentclass[reqno]{amsart}

\address{Department of Mathematics, Kyoto
University, Kyoto, Japan } \email{fukaya@math.kyoto-u.ac.jp}
\address{Department of Mathematics, University of
Wisconsin, Madison, WI, USA \& Department of Mathematics,
POSTSECH, Pohang, Korea} \email{oh@math.wisc.edu}
\address{Graduate School of Mathematics,
Nagoya University, Nagoya, Japan \&
Korea Institute for Advanced Study, Seoul, Korea} \email{ohta@math.nagoya-u.ac.jp}
\address{Research Institute for Mathematical Sciences, Kyoto University, Kyoto, Japan \&
Korea Institute for Advanced Study, Seoul, Korea}
\email{ono@kurims.kyoto-u.ac.jp}

\usepackage{graphicx}
\usepackage{amsmath}
\usepackage{amscd}
\usepackage{amssymb}
\usepackage{amstext}
\usepackage{amsmath}
\usepackage[all]{xy}

\def\E{\ifmmode{\mathbb E}\else{$\mathbb E$}\fi} %natural numbers
\def\N{\ifmmode{\mathbb N}\else{$\mathbb N$}\fi} %natural numbers%
\def\R{\ifmmode{\mathbb R}\else{$\mathbb R$}\fi} %real numbers
\def\Q{\ifmmode{\mathbb Q}\else{$\mathbb Q$}\fi} %rational numbers
\def\C{\ifmmode{\mathbb C}\else{$\mathbb C$}\fi} %complex numbers
\def\H{\ifmmode{\mathbb H}\else{$\mathbb H$}\fi} %complex numbers
\def\Z{\ifmmode{\mathbb Z}\else{$\mathbb Z$}\fi} %integers
\def\P{\ifmmode{\mathbb P}\else{$\mathbb P$}\fi} %real numbers
\def\T{\ifmmode{\mathbb T}\else{$\mathbb T$}\fi} %real numbers
\def\SS{\ifmmode{\mathbb S}\else{$\mathbb S$}\fi} %real numbers
\def\DD{\ifmmode{\mathbb D}\else{$\mathbb D$}\fi} %real numbers
\def\K{\ifmmode{\mathbb K}\else{$\mathbb K$}\fi}

\def\EE{{\mathcal E}}
\def\KK{{\mathcal K}}

\theoremstyle{theorem}
\newtheorem{thm}{Theorem}[section]
\newtheorem{cor}[thm]{Corollary}
\newtheorem{lem}[thm]{Lemma}
\newtheorem{sublem}[thm]{Sublemma}
\newtheorem{prop}[thm]{Proposition}

\newtheorem{lemdef}[thm]{Lemma-Definition}

\theoremstyle{definition}
\newtheorem{defn}[thm]{Definition}
\newtheorem{rem}[thm]{Remark}

\newtheorem{exm}[thm]{Example}
\newtheorem{conds}[thm]{Condition}
\newtheorem{proper}[thm]{Property}

\newtheorem{assump}[thm]{Assumption}

\numberwithin{equation}{section}
\begin{document}
\title[details on Kuranishi structure]{Technical details on Kuranishi structure and virtual fundamental chain}
\author{Kenji Fukaya, Yong-Geun Oh, Hiroshi Ohta, Kaoru Ono}
\begin{abstract}
This is an expository article describing the theory of Kuranishi structure in great detail.
The origin of this article lies in a series of pdf
files \cite{Fu1,FOn2,fooo:ans3,fooo:ans34,fooo:ans5} that
the authors of the present article uploaded for the discussion of the google
group named `Kuranishi' (with its administrator H. Hofer). In these pdf files,
the present authors replied to several questions concerning the foundation of
the Kuranishi structure which  were raised by K. Wehrheim.
At this stage we submit this article to the e-print arXiv, all the
questions or objections asked in the google group `Kuranishi'
were answered, supplemented or confuted by us.
\par
We provide (in Part \ref{origin}) our confutations against several criticisms we found in
\cite{MW1} (arXiv:1208.1304v1). We have seen a few instances in which the public display of
negative opinions, such as those written in \cite{MW1}, on the soundness of
virtual fundamental chain/cycle  technique  of Kuranishi  structure
has caused serious trouble for young  mathematicians to publish
his/her papers. Due to this reason, we feel obliged not
to escape from our duty  of  confuting  the article \cite{MW1} and
providing more thorough explanations both on various technical points of
Kuranishi structure and on its implementation in the study of
pseudoholomorphic curves.
\par
We would like to mention that before \cite{MW1} was uploaded,
all the objections in \cite{MW1} had  been already discussed and responded in our files mentioned above
which had been sent to the members of google group
Kuranishi (that include the authors of \cite{MW1}).
\par
In the first part of this article (Part \ref{Part2}) we
discuss the abstract theory of Kuranishi structure
and virtual fundamental chain/cycle.
We review the definition of
Kuranishi structure in \cite{fooo:book1} and
explain in detail the way how we define a
virtual fundamental chain/cycle for a space with
Kuranishi structure.
This part can be read independently from other parts.
\par
In the second part (Parts \ref{secsimple} and \ref{generalcase}) of this article
we describe the construction of Kuranishi structure
on the moduli space of pseudo-holomorphic curves, in great detail.
We include the complete analytic detail of the gluing (stretching the neck) construction
as well as several other issues e.g., the smoothness of the
resulting Kuranishi structure.
The case of $S^1$ equivariant Kuranishi structure which appears
in the study of time independent Hamiltonian and the moduli space of
Floer's equation is included in Part \ref{S1equivariant}.
\end{abstract}

\maketitle
\par\newpage

\tableofcontents
\par\medskip
\par\newpage

\part{Introduction}
\section{Introduction}

This is an expository article describing the theory of Kuranishi structure,
its construction for the moduli space of pseudo-holomorphic curves of various kinds
and the way how to use it to define and study virtual fundamental class.
Our purpose is to provide technical details much enough for mathematicians
who want to apply Kuranishi structure for various purposes can feel confident of
its foundation. We believe that by now (16 years after its discovery), the methodology of Kuranishi
structure can be used as a `black-box' meaning that mathematicians can freely use
it without checking its details in the level we provide in this article, once they understand
the general methodology and basic ideas. To ensure them of the
preciseness and correctness of all the basic results stated in \cite{FOn}, \cite{fooo:book1},
we provide thorough details so that people can use them without doubt.
\par
The origin of this article lies in a series of pdf
files \cite{Fu1,FOn2,fooo:ans3,fooo:ans34,fooo:ans5} that
the authors of the present article uploaded for the discussion of the google
group named `Kuranishi' (with its administrator H. Hofer),
started around March 14, 2012. In these pdf files,
the present authors prepared the replies to several questions concerning the foundation of
the Kuranishi structure which  were raised by K. Wehrheim.
(We mention about the discussion in this google more in Part \ref{origin}.)
The pdf files themselves that we uploaded can be obtained in the home page of
the second named author (http://www.math.wisc.edu/~oh/).
\par
The theory of Kuranishi structure first appeared in 1996 January in a paper by first and fourth named authors,
which was published as \cite{FuOn99I}. More details thereof were published in \cite{FOn}.
These papers contained some technical inaccurate points which were corrected in the book
written by the present authors \cite{fooo:book1}.
In \cite{fooo:book1} the same methodology as the one used in \cite{FuOn99I,FOn} is applied systematically
for the construction of filtered $A_{\infty}$ structure
on the singular (co)homology of a Lagrangian submanifold of a symplectic manifolds.
\par
The construction of the Kuranishi structure of the moduli space of
pseudo-holomorphic curves \cite{FOn} is written in the way suitable
for the purpose of \cite{FOn} (especially \cite[Theorem 1.1, Corollary 1.4]{FOn}).
(See Subsection \ref{subsec342} for more discussion on this point.)
In \cite{fooo:book1} (especially in its A1.4)
we provide more detail so that it can be used for the
purpose of chain level construction we used there.
The present article contains even more details of
the various parts of this construction.
\par
Meanwhile there are several articles which
describe the story of Kuranishi structure,
for example \cite{joyce}, \cite{joyce2}.
\par
Several other versions of the construction of the virtual
fundamental chain or cycle via Kuranishi structure are
included in some of the papers of the present
authors (\cite{fooo:toric1,fooo:toric2,fooo:toricmir,cyclic}) aimed for various applications.
\par
This article is {\it not} a research paper {\it but} an expository article.
All the results together with the basic idea of its proof had been
published in the references we mention above.
\par
Kuranishi structure is one of the various versions of the technique so called
virtual fundamental cycle/chains.
Several other versions of the same technique appeared
in the year 1996 (the same year as \cite{FOn} appeared.)
\cite{LiTi98, LiuTi98,Rua99, Sie96}.
Some more detail of \cite{LiTi98, LiuTi98} was appeared in \cite{LuT}.
\par
Later on, other versions
of virtual technique appeared \cite{HWZ,ChiMo}.
Also there are various expository articles, e.g., \cite{Salamon,McDuff07,operad,suugakuexp,MW1},
which one can obtain in various places.
\par
Because of its origin, this article is written in the style
so that it will serve as a reference that confirms the
solidness of the foundation of the theory.
Therefore the preciseness and rigor are our major concern,
while writing this article.
We are planning to provide a text in the future, which is more easily
accessible
to non-experts, such as graduate students of the area or
researchers from the other related fields.

\section{The technique of virtual fundamental cycle/chain}

We start with a very brief review of the technique of virtual fundamental cycles/chains.
In differential geometry, various moduli spaces appear as `the set of solutions of nonlinear partial
differential equations'. Here we put the word in the quote because there are several important issues to be
taken care of.
\begin{enumerate}
\item[(A)]
The moduli space is in general very much singular.
\item[(B)]
We need to take appropriate equivalence classes of the solution set
to obtain a moduli space.
\item[(C)]
We need an appropriate compactification to obtain a useful
moduli space.
\end{enumerate}
In the case of moduli spaces appearing in differential geometry,
the point (A) was studied by Kuranishi in Kodaira-Spencer theory of deformation of complex structures.
For each compact manifold $(X,J)$, Kuranishi constructed a finite dimensional
complex manifold $V$ on which the group of automorphisms $\text{Aut}(X,J)$ acts
and a holomorphic map $s : V \to E$ to a complex vector space,
such that $E$ is acted by $\text{Aut}(X,J)$ and $s$ is $\text{Aut}(X,J)$-equivariant,
and the moduli space of complex structure of $J$ locally
is described as
\begin{equation}\label{kuranishimodel}
s^{-1}(0)/\text{Aut}(X,J).
\end{equation}
This is called the Kuranishi model.
The map $s$ is called the {\it Kuranishi map}.
\par
In 1980' first by Donaldson in gauge theory and then by Gromov in the theory of
pseudo-holomorphic curves, the idea to use the fundamental homology class
of various moduli spaces to obtain an invariant was discovered.
In the theory of pseudo-holomorphic curves, Gromov-Witten invariant
was obtained in this way. Such a theory was rigorously built
in the case of semi-positive symplectic manifold by
Ruan \cite{ruanD} and Ruan-Tian \cite{RuTi95}. (See also \cite{McSa94}.)
\par
Around the same time, Floer used the moduli space of solutions
of pseudo-holomorphic curve equation with an extra term defined by a
Hamiltonian vector field and succeeded in rigorously defining a homology theory,
that is now called Floer homology of periodic Hamiltonian system.
In \cite{Flo89I} Floer assumed that the symplectic manifold is monotone.
This assumption is weakened to semi-positivity in \cite{HoSa95} and \cite{Ono95}.
\par
In both cases, we need to study the moduli spaces of virtual dimensions $0$ and $1$ (see e.g., page 1020 \cite{FOn}).
We construct multi-valued perturbation (multi-section) inductively from the smallest energy
(thanks to Gromov's compactness theorem) and
can arrange that no zeros appear in the strata of negative virtual dimension.  Combined with Gromov's compactness
theorem, we can also arrange that there are finitely many zeros in the strata of virtual dimension $0$.
When the virtual dimension of the moduli space is $1$, we pick a Kuranishi neighborhood for each zero of the
multi-section such that they are mutually disjoint.  Then we extend the multi-section so that
the weighted number of zeros in the virtual dimension $0$ stratum coincides with the one with
sufficiently large fixed gluing parameter $T$ (see Subsection \ref{subsec:Smoothness of coordinate changes} below).
We also note that in the case of Gromov-Witten invariants, what we need to study is not the detailed geometric data
which moduli space carries but only its homology class, which is significantly weaker information.
For the purpose in \cite{fooo:book1}, we have to work with {\it not} homology classes {\it but} chains.
\par
Roughly speaking the moduli space can be regarded as the zero set of a section of
certain infinite dimensional bundle over an infinite dimensional space.
Thus its fundamental class is nothing but the `Poincar\'e dual' to the `Euler class' of the bundle.
It is well-known that in the finite dimensional case the Euler class is a topological invariant
of the bundles and so in particular we can take any section to study it, when the base space is a closed
manifold.
In the infinite dimensional situation we need to take sections so that
they satisfy appropriate Fredholm properties but still exist so abundantly
that one has much freedom to perturb.
\par
Thus, in a situation when the automorphism group of the objects are trivial,
it is very easy to find a perturbation of the equation in an abstract way and
find a perturbed moduli space that is smooth.
This was actually the way taken by Donaldson  \cite[II.3]{Don83} in gauge theory.
\par
The problem becomes nontrivial when the points (B),(C) enter.
Let us restrict our discussion below to the case of
the moduli space of pseudo-holomorphic curves.
\par
The point (B) causes the most serious trouble in case the group of
automorphisms is noncompact.
In fact in such a case the moduli space is not Hausdorff in general.
This point is studied in the work of Mumford who introduced the
notion of stability and used it systematically to study
moduli space of algebraic varieties.
The case of curves (Riemann surfaces) is an important case of
it.
\par
Gromov-Witten theory or the theory of pseudo-holomorphic curves
is a natural generalization thereof where we consider the pair
$(\Sigma,u)$ where $\Sigma$ is a Riemann surface (which includes the
case of complex curves with nodal singularities when we study compactification),
and $u : \Sigma \to X$ is a pseudo-holomorphic map.
(We may includes a finite number of points $\vec z \in \Sigma$, which are nonsingular
and disjoint.)
The group $\text{Aut}(\Sigma,\vec z,u)$ of automorphisms  consists of
biholomorphic maps $v : \Sigma \to \Sigma$ such that
$u\circ v = u$ and $v$ fixes every marked point in $\vec z$.
\par
The notion of stable maps due coined by Kontsevitch clarifies the issue here.
He called the triple $(\Sigma,\vec z,u)$ stable when $\text{Aut}(\Sigma,\vec z,u)$
is of finite order. This is a natural generalization of the notion of stability due to
Mumford defined for the case of stable curves $(\Sigma,\vec z)$.
Kontsevitch observed that the moduli space of
stable maps $(\Sigma,\vec z,u)$ is Hausdorff.
The first and fourth named author gave a precise
definition of the relevant topology and gave the proof of Hausdorff property in \cite[Definition 10.3, Lemma 10.4]{FOn}.
Hausdorff property is discussed also in \cite{Sie96}.
\par
On the other hand, though the stability implies that the group $\text{Aut}(\Sigma,\vec z,u)$
is finite,  it may still be nontrivial.
It means that in the local description as in (\ref{kuranishimodel}),
the group of automorphisms can still be nontrivial.
In other words our situation is closer to that of an
orbibundle on an orbifold rather than to a vector bundle on a manifold.
An orbibundle is a vector bundle in the category of orbifold.
\par
It is well-known classical fact that a generic equivariant section
of an equivariant vector bundle is not necessarily transversal to zero,
even when the group is finite.
After taking the quotient it means that an orbibundle over an orbifold may not have transversal
sections.
A simple explanation of this fact is that the Euler class of the orbibundle
is not necessarily an element of  cohomology group over the integer
coefficients but is defined only in  cohomology group over the rational coefficients.
\par
At the year 1996 three approaches are proposed and worked out
on this point and applied to the study of the moduli space of
pseudo-holomorphic curves.
\par
One approach due to J. Li and G. Tian is an analytic
version of `locally free resolution of tangent complex'.
The approach by Ruan \cite{Rua99} uses de-Rham
theory and uses the representative of Euler class in
de Rham cohomology.
\par
The approach of \cite{FOn} is based on multi-valued perturbations,
which was called multisections.
\par
Before explaining more about the method of \cite{FOn},
we explain the point (C).
As usual, compactification of the moduli space of
geometric object is obtained by adding certain kinds of
`singular objects'.
In the case of the moduli space of
pseudo-holomorphic curves, such a singular object consists of the triple $(\Sigma,\vec z,u)$
where $\Sigma$ is a curve with only double points, i.e. nodal singularities, as its singularities,
$\vec z$ are marked points and $u : \Sigma \to X$
is a pseudo-holomorphic map. We can define stability condition as Kontsevich did.
\par
The very important point of the story is we can
define a coordinate chart in a neighborhood of such
objects $(\Sigma,\vec z,u)$ in the same way as the case when $\Sigma$ is smooth.
In other words, the moduli space of the triples can be also described locally
as the Kuranishi model (\ref{kuranishimodel}).
This is a consequence of the process called gluing or stretching the neck.
Such a process has its origin in the work of Taubes in gauge theory.
By now it very much became a standard practice also in the case of
pseudo-holomorphic curves.
\par
We thus find that each point of the compactification of our
moduli space has a local description by  Kuranishi model (\ref{kuranishimodel}).
\par
We can then find a multi-valued perturbation (= mutisection)
on each of the Kuranishi model so that $0$ is a regular value for
each branch of our multivalued perturbation.
Then the task is to formulate the way how those local constructions
(perturbation) are globally compatible.
\par
The notion of Kuranishi structures was introduced for this purpose.
Namely a Kuranishi structure by definition provides a way of describing
the moduli space locally as
$
s^{-1}(0)/\Gamma
$
where $s : V \to E$ is a $\Gamma$-equivariant
map from a space $V$ equipped with an action of a
finite group to a vector space $E$ on which $\Gamma$ acts linearly.
We say such a local description as a Kuranishi neighborhood.
\par
A Kuranishi structure also involves coordinate changes between
Kuranishi neighborhoods and requires certain compatibility between  coordinate changes.
\par
Thus the main idea of this story is as follows.
\begin{enumerate}
\item To define some general notion of `spaces' that contain
various moduli spaces of pseudo-holomorphic curves as examples and work out transversality
issue in that abstract setting.
\item  Use multivalued abstract perturbations, that we call
multisection, to achieve relevant transversality.
\end{enumerate}
In this article we describe
technical details of this method.

\section{The outline of this article}

The theory of Kuranishi structures is divided into two parts.
One is the abstract theory in which we first define the notion of Kuranishi
structure and then we describe how we obtain a virtual fundamental chain/cycle
or its homology class in that abstract setting.
The other is the methodology of implementing the abstract theory
of Kuranishi structure in the study of concrete moduli problem, especially
that of the moduli space of pseudo-holomorphic curves.
\par
We discuss the first point in Part \ref{Part2}
and the second point in Parts \ref{secsimple} and \ref{generalcase}.
\par
The definition of the Kuranishi structure is
given in Section \ref{secdefnkura}.
\par
To construct a perturbation (multisection) of a given
Kuranishi map that is transversal to $0$, we work in a local chart (Kuranishi neighborhood)
and apply an appropriate induction process. \par
This is very much similar to Thom's original proof of transversality theorem
in differential topology. Later on, the proof of transversality theorem via Sard's theorem
combined with Baire category theorem gained more popularity. However the latter
approach, which uses functional analysis, meets some trouble in working out for the case of
multisections. In fact, the sum of multisections is rather delicate to define.
There seems to be no way to define the sum so that an additive inverse exists.
Because of this it is unlikely that the totality of the multisections becomes
a vector space of infinite dimension.
\par
To work out the way to inductively define multisection,
we need to find a clever choice of a system of Kuranishi neighborhoods.
We called such a system of Kuranishi neighborhoods `a good coordinate system'.
\par
We remark that in the local description:
$
s^{-1}(0)/\Gamma
$
($s : V \to E$)
of our moduli space, the number $\dim V - \text{rank}\, E$
is the `virtual' dimension of our moduli space and is a
well-defined number. In other words, it is independent of the Kuranish neighborhood.
On the other hand the dimension of the base $V$ may depend on the Kuranishi neighborhood.
As its consequence,
the coordinate change exists sometimes only in one direction, namely from the Kuranishi neighborhood with
smaller $\dim V$ to the one with bigger $\dim V$.
This makes the proof of the existence of a {\it consistent} system of Kuranishi neighborhoods much more
nontrivial compared to the case of ordinary manifolds.
Recall that already in the case of orbifold (that is the case when obstruction space $E$
is always zero), the order of the group $\Gamma$ may vary and so
the natural procedure of constructing a transversal multisection
of an orbibundle over an orbifold is via the induction over the order of $\Gamma$.
The case of Kuranishi structure is slightly more nontrivial
since the dimension of $V$ may also vary.
\par
The definition of a good coordinate system is given in Section \ref{defgoodcoordsec}.
Existence of such a good coordinate system is proved in Section
\ref{sec:existenceofGCS}.
\par
We alert the readers that in this article
more conditions are required for our definition of good coordinate system
compared to that of our earlier papers \cite{FOn}, \cite{fooo:book1}.
The reason is because it is more convenient for the purpose of writing the technical details
of a part of the construction of the virtual fundamental chain/cycle.
This detail was asked recently by a few people\footnote{It includes D. Yang, K. Wehrheim,
D. McDuff. We thank them for asking this question.}.
For example, a question about how we restrict the domains (of the Kuranishi neighborhoods) of the
perturbations so that the zero sets of the Kuranishi maps that are
defined in each Kuranishi neighborhoods can be glued together to define a
Hausdorff and compact space.
\par
We emphasize that this problem of Hausdorff-ness
and compactness is of very much different nature from, for example, that of
Hausdorff-ness or compactness of the moduli spaces itself.
The latter problem is related to some key geometric notion
such as stability and requires to study certain fundamental
points of the story like the Gromov compactness and the
removable singularity results. On the other hand, the former problem,
though it is rather tedious and complicated to write down a precise way to resolve, is of
technical nature. It boils down to finding a right way of restricting various domains
(Kuranishi neighborhoods) with much care and patience. It, however, goes without saying that
writing down this tedious technicality at least once
is certainly a nontrivial and meaningful labor to do for the salience of the field,
which is the main purpose of Part \ref{Part2}. This makes Part \ref{Part2} rather lengthy
and complicated. Section \ref{defgoodcoordsec} contains some of those
technical arguments.
For the purpose above, we use a general lemma in general topology
which we prove in  Section \ref{gentoplem}. This lemma (Proposition \ref{metrizable})
is in principle well-known, we suspect. We include its proof here for the sake of completeness
since we could not locate an appropriate reference.
We gather well-known facts on orbifold, its embedding, and a bundle on it,
in Section \ref{ofd} for reader's convenience.

These technical points, however, should not be confused with the
basic and conceptional points of the theory of Kuranishi structures.
We believe that the readers, especially with geometric applications in their minds,
can safely forget most of those technical arguments
once they go through and convince themselves of
the soundness of the foundation. The bottom line of the Kuranishi methodology
is to make sure the existence of Kuranishi structure on the \emph{compactified}
moduli spaces in the relevant moduli problems.  (This step is {\it not} automatic.)
Then the rest automatically follows
by the general theory of Kuranishi structures.
\par
In Section \ref{sec:existenceofGCS},  we prove the existence of the
good coordinate system in the sense defined in this paper,
(which is more restrictive compared to the one in \cite{FOn}
or \cite{fooo:book1}.)
The proof is based on the idea with its origin in the proof of Lemma 6.3 in \cite[page 957]{FOn}.
We work by a downward induction on the dimension $\dim V$ of the Kuranishi neighborhood and
in each dimension we glue several Kuranishi neighborhoods (of the same
dimension) to obtain a bigger Kuranishi neighborhood.
\par\smallskip
Parts \ref{secsimple} and \ref{generalcase} provide details of
the construction of the Kuranishi structure
on the moduli space of pseudo-holomorphic curves.
There are two main issues in the construction.
\par
One is of analytic nature.
Namely we construct a Kuranishi neighborhood of a given element
of the compactified moduli space.
In the case when the given element is a stable map
from a nonsingular curve (Riemann surface),
the analytic part of the construction is a fairly
standard functional analysis.
\par
In case when the element is a stable map from a curve which has
a nodal singularity, its neighborhood still contains a
stable map from a nonsingular curve.
So we need to study the problem of gluing or of stretching the neck.
Such a problem on gluing solutions of non-linear elliptic partial differential equation
has been studied extensively in gauge theory and in symplectic
geometry during the last decade of the 20th century.
Several methods had been developed to solve the problem which
are also applicable to our case. In this article, following \cite[Section A1.4]{fooo:book1},
we employ the alternating method, which was first exploited by Donaldson \cite{Don86I} in gauge theory.
\par
In this method, one solves the linearization of the given equation on each
piece of the domain (that is the completion of the complement
of the neck region of the source of our pseudo-holomorphic curve.)
Then we construct a sequence of families of maps that converges to a version of
solutions of the Cauchy-Riemann equation, that is,
\begin{equation}\label{mainequation00}
\overline{\partial} u' \equiv 0  \mod E(u')
\end{equation}
and which are parameterized by a manifold (or an orbifold).
Here $E(u')$ is a family of finite dimensional vector spaces
of smooth sections of an appropriate vector bundle
depending on $u'$.  More precisely, $E(u')$ is defined using additional marked points, which makes
the domain curve stable, see Part 3,4.
\par

\subsection{Smoothness of coordinate changes}\label{subsec:Smoothness of coordinate changes}

The authors were sometimes asked a question about the
smoothness of the Kuranishi map $s$ and of the coordinate change of the Kuranishi
neighborhoods\footnote{Among others, Y.B. Ruan, C.C. Liu, J. Solomon, I. Smith
and H. Hofer asked the question. We thank them for asking this question.}.
\par
This problem is described as follows.
Note in our formulation the neck region is
a long cylinder $[-T,T] \times S^1$ or
long rectangle $[-T,T] \times [0,1]$,
and the case when the source curve
is singular corresponds to the case when $T = \infty$.
So a part of the coordinate
of our Kuranishi neighborhood is naturally parametrized by
$(T_0,\infty]$ or its products.
Note $\infty$ is included in $(T_0,\infty]$.
As a topological space $(T_0,\infty]$ has unambiguous meaning.
On the other hand there is  no obvious choice of its smooth structure as a manifold with boundary.
Moreover for several maps such as Kuranishi map, $s$,
it is not obvious whether it is smooth for given
coordinate of $(T_0,\infty]$. (See \cite[Remark 13.16]{FOn}.)
As we will explain in Subsection \ref{subsec342}
there are several ways to resolve this problem.
One approach is rather topological and uses the fact that the chart is smooth
in the $T$-slice where the gluing parameter $T$ above is fixed.
This approach is strong enough to establish all the results of \cite{FOn}.
The method of \cite{McSa94} which is quoted in \cite{FOn}
is strong enough to work out this approach.
However it is not clear to the authors whether it is good enough to establish
smoothness of the Kuranishi map or of the coordinate
changes at $T=\infty$. (This point was mentioned by the first and the fourth named
authors themselves in the year 1996 at \cite[Remark 13.16]{FOn}.)
To prove an existence of the Kuranishi structure that literally satisfies our axioms, we take
a different method in this article.
\par
Using the alternating method described in \cite[Section A1.4]{fooo:book1}
for the same purpose, we can find an appropriate coordinate chart at
$T=\infty$ so that the Kuranishi map and the coordinate changes of our Kuranishi neighborhoods are
of $C^{\infty}$ class. For this purpose, we take the parameter $s = 1/T$.
As we mentioned in \cite[page 771]{fooo:book1}
this parameter $s$ is different from the one taken in algebraic geometry when the
target $X$ is projective. The parameter used in algebraic geometry corresponds to $e^{-T}$.
It seems likely that in our situation either where
almost complex structure is  non-integrable and/or where we include the Lagrangian submanifold
as the boundary condition (the source being the borderded stable curve) the
Kuranishi map or the coordinate changes
are {\it not} smooth with respect to the parameter $e^{-T}$.
But it is smooth in our parameter $s = 1/T$, as is proved in
\cite[Proposition A1.56]{fooo:book1} and Part \ref{generalcase}.
\par
The proof of this smoothness is based on
the exponential decay estimate of the solution of the equation
(\ref{mainequation00}) with respect to $T$, that is, the length of the neck.
The proof of this exponential decay is given
in \cite[Section A1.4, Lemma A1.58]{fooo:book1}.
Because, after the publication of \cite{fooo:book1}, we still heard some demand of providing more
details of this smoothness proof, we provide such detail in Part \ref{secsimple}
and Section \ref{glueing}.
\par
In Part \ref{secsimple}, we restrict ourselves to the case we glue
two (bordered) stable maps such that the source
(without considering the map) is already
stable. By restricting ourselves to this case we can
explain all the analytic issues needed to work out the
general case also without making the notation so much complicated.
\par
We provide the relevant analytic details using the same induction scheme as
\cite[Section A1.4, page 773-776]{fooo:book1}. The only difference is that we use $L^2_{m}$ space (the space of maps
whose derivative up to order $m$ are of $L^2$ class) here,
while we used $L^p_{1}$ space following the tradition of
symplectic geometry community in \cite[Section A1.4]{fooo:book1}.
Actually using $L^2_{m}$ space in place of $L^p_{1}$ space,
it becomes easier to study the derivatives of our solution with respect to the parameter $T$.
See Remark \ref{Abremark}.
\par
In Section \ref{alternatingmethod} we provide the details of the estimate and show that
the induction scheme of \cite[Section A1.4]{fooo:book1} provides
a convergent family of solutions of our equation (\ref{mainequation00}).
(This estimate is actually fairly straightforward to prove although tedious to write down.)
\par
In Section \ref{subsecdecayT}, we provide the detail of the above
mentioned exponential decay estimate of the $T$-derivatives
of our solutions. In Section \ref{surjinj}, we review the well-established classical
proof of the fact that the solutions obtained exhaust all the
nearby solutions and also the map from the parameter
space to the moduli space is injective.
\par\medskip

\subsection{Construction of Kuranishi structure}

In the first half of Part \ref{generalcase}, we discuss the second main issue,
which enters in the construction of a Kuranishi structure.
The problem here is as follows.
To define Kuranishi neighborhoods, we need to take the obstruction spaces
which appear in the right hand side of (\ref{mainequation00}).
We need to choose them so that the Kuranishi neighborhoods
that are the solution spaces of (\ref{mainequation00}) can be glued together.
In other words we need to choose them so that we can define
smooth coordinate changes.
\par
We need to fix a parametrization of the source (the curve) by the following
reason.
Let $\frak p=(\Sigma,\vec z,u)$ be an element of the moduli space of stable maps.
Firstly we consider the case that the domain $(\Sigma, \vec z)$ is stable.
Consider a vector space $E_\frak p$ spanned by finitely many smooth sections of
sections of the bundle $u^*TX \otimes \Lambda^{0,1}$.
For  $\frak q=(\Sigma',\vec z',u')$ close to $\frak p$, we would like to transport $E_\frak p$ to $\frak q$.
Let $K$ be a compact subset of $\Sigma$ where elements of $E_\frak p$ are supported.
(The subset $K$ is chosen so that it is disjoint from the nodal singularities.)
If we fix a diffeomorphism to the image $K \to \Sigma'$ then
we can use the parallel transport along the closed geodesic to transfer $E_\frak p$ to $E_\frak p(\frak q)$.
This family $E_\frak p(\frak q)$ is obviously a smooth family of vector spaces
of smooth sections so we can study the solution space of (\ref{mainequation00})
by using implicit function theorem, for example.
In the case when $(\Sigma,\vec z)$ is stable we can choose
such a diffeomorphism $K \to \Sigma'$
(up to an ambiguity of finite group of automorphisms)
by using the universal family of curves over Deligne-Mumford moduli space.
\par
In case when $(\Sigma,\vec z)$ is not a stable curve  (but $\frak p=(\Sigma,\vec z,u)$
is a stable map) we need some additional argument.
In \cite{FOn} we gave two methods to resolve this trouble.
One is explained in \cite[Section 15]{FOn} and the other
in \cite[Appendix]{FOn}. The first one uses the center of mass technique
from Riemannian geometry. We explain it a bit in
Subsection \ref{comparizon}. In  \cite{fooo:book1}
and several other places we used the second technique (the one in
\cite[Appendix]{FOn}.)
Therefore we use the second method mainly in this article.
In this method we add additional marked points $\vec w$ on $\Sigma$
so that $(\Sigma,\vec z \cup \vec w)$ becomes stable.
We also put the additional marked points
$\vec w'$ on $\Sigma'$
so that $(\Sigma',\vec z' \cup \vec w')$ becomes stable.
Then we can transform $E_\frak p$ to $E_\frak p(\frak q \cup \vec w')$.
The resulting moduli space (in fact, we take a direct sum \eqref{obstruction} below in later argument),
which we call a thickened moduli space,
has too many extra parameters compared to those required by the
virtual dimension obtained by dimension counting. The extra parameters correspond to the
positions of the points $\vec w'$.
We kill these extra freedom as follows.
We take a finite collection of codimension two submanifolds
$D_i \subset X$ for each added marked point in $w_i \in \vec w$ so that
$D_i$  transversally intersects $u(\Sigma)$ at $u(w_i)$, and require
$w'_i$ to satisfy $u'(w'_i) \in D_i$.
\par
This gives a construction of a Kuranishi neighborhood at each point of our moduli space.
To obtain Kuranishi neighborhoods which can be glued to obtain a global structure, we proceed as follows.
\par
We take a sufficiently dense finite subset
$\{\frak p_c = (\Sigma_c,\vec z_c,u_c)\}_c$ in our moduli space.
For each $\frak p_c$ we fix a finite dimensional vector space
of sections $E_c=E_{\frak p_c}$ which will be a part of the obstruction space $E$.
We also fix additional marked points $\vec w_c$ so that $(\Sigma_c,\vec z_c\cup \vec w_c)$
becomes stable.
\par
We next consider  $(\Sigma',\vec z',u')$ for  which we will set up our equation.
We take all $\frak p_c$'s which are `sufficiently close' to $(\Sigma',\vec z',u')$.
For each such $c$, we take additional  $\vec w'_c$ so that $(\Sigma',\vec z'\cup \vec w'_c)$
becomes close to $(\Sigma_c,\vec z_c\cup \vec w_c)$ in an obvious topology.
We then use the diffeomorphism obtained by this closeness, and parallel transport to transfer $E_c$
to a finite dimensional vector space $E_c(u',\vec w'_c)$ of sections on $\Sigma'$.
We take a sum of them over all $c$'s to obtain the fiber of the obstruction bundle

\begin{equation}
E(u';(\vec w'_c)) = \bigoplus_{c} E_c(u';\vec w'_c). \label{obstruction}
\end{equation}
We remark here that this space depends on all the additional marked points $\bigcup_c \vec w'_c$.
We solve the equation \eqref{mainequation00} to obtain the moduli space of larger dimension.
Finally we cut this space by requiring the condition that each of those
additional marked points lies in the corresponding codimension 2 submanifold that
we have chosen at the time as we define $E_c$.
\par
This process of defining $E(u';(\vec w'_c))$ and its solution space is
described in detail in Sections \ref{Graph} -\ref{settin2}.
\par
In the first two sections of those, we discuss the following point.
Note we say that the diffeomorphism $K \to \Sigma'$ is determined
if the source $(\Sigma,\vec z)$ is stable.
More precisely we proceed as follows.
Note $(\Sigma',\vec z')$ is close to $(\Sigma,\vec z)$ in Deligne-Mumford
moduli space.  To identify the set $K \subset \Sigma'$ with a subset of $\Sigma$
we need to fix a trivialization of the universal family.
Actually the universal family is not a smooth fiber bundle even in orbifold sense,
since there is a singular fiber which corresponds to the nodal curve.
So to specify the embedding $K \to \Sigma$ we also need
to fix the way to resolve the node. The notion of
`coordinate at infinity' is introduced for this purpose.
\par
After introducing the notion `coordinate at infinity', we define the notion of
obstruction bundle data in Section \ref{stabilization}.
The obstruction bundle data consist of the finite dimensional
vector space of sections $E_c$ (that will be the obstruction bundle $E_c(u';\vec w'_c)$ nearby)
together with the additional marked points $\vec w_c$ which we use to transform
$E_c$ to the nearby maps as explained above.
Using these data, the thickened moduli space, which is the solution space of the equation
\begin{equation}\label{mainequation01}
\overline{\partial} u' \equiv 0  \mod E(u';(\vec w'_c)),
\end{equation}
is defined in Section \ref{settin2}.
\par
Section \ref{glueing} is a generalization of the analytic
argument of Part \ref{secsimple}. Most of the arguments in
Part \ref{secsimple} can be generalized here without change.
The most important point which is new here is the following:
For a point $\frak p$ in the moduli space, we
construct a thickened moduli space containing it.
To obtain the vector space $E$ at $\frak p$ which is the fiber thereat
of the obstruction bundle, we consider various $\frak p_c$ in a neighborhood
of $\frak p$ and $E_c(u';\vec w'_{c})$ parallel transported from  $E_c$
using the marked points $\vec w'_{c}$ as mentioned before.
On the other hand, we fix a parameterization of the source of the map $u'$ using
a stabilization at $\frak p$ and some other additional marked points
$\vec w'_{\frak p}$ associated to the stabilization.
Therefore the parametrization of $\Sigma'$ used to define
$E_c(u')$ is different from the one that we use
to study our equation (\ref{mainequation01}).
As far as we are working with smooth curves (the curves without nodal singularity)
this is not really an issue since elements of $E_c$ are smooth sections
and they behave nicely under a diffeomorphism (or under the change of variables).
However when we study gluing of solutions (that is, for the case when
$\frak p$ has a node), we need to study the asymptotic behavior
of this coordinate change as the gluing parameter $T$ goes to infinity.
Study of this asymptotic behavior is also needed
when we prove smoothness of the coordinate change at the boundary or at the corner.
The main ingredients that we use for this purpose are Propositions \ref{changeinfcoorprop}, \ref{reparaexpest},
which are generalizations of \cite[Lemma A1.59]{fooo:book1}.
Propositions \ref{changeinfcoorprop}, \ref{reparaexpest}  are proved in Section \ref{changeinfcoorprop}.
\par
In Section \ref{cutting} we discuss the process of putting the condition $u'(w_i) \in D_i$ to cut
the dimension of the thickened moduli space in detail.
In particular we show that after doing this cutting down and taking the quotient by the finite
group of automorphisms, the set of the solutions of the associated Cauchy-Riemann equations
has right dimension.
\par
Now we construct the moduli space of the solutions of the equation (\ref{mainequation01})
this time requiring the left hand side becomes exactly zero, is homeomorphic to the original unperturbed moduli space.
This fact is used in Section \ref{chart} to define a Kuranishi neighborhood at every
point of the moduli space.
\par
In the next three sections, we construct the coordinate changes between Kuranishi neighborhoods
and show they are compatible.

\subsection{$S^1$-equivariant Kuranishi structure and multi-sections}

As we mentioned before Floer studied the pseudo-holomorphic curve equation with
extra term defined by Hamiltonian vector field
and use its moduli space to define Floer homology of
periodic Hamiltonian system. We can define Kuranishi structure on the moduli space of solutions of Floer's
equation in the same way. We can use this to generalize
Floer's definition of  Floer homology of
periodic Hamiltonian system to an arbitrary symplectic manifold.
\par
This part of the generalization is actually fairly straightforward.
The point mainly discussed in Part \ref{S1equivariant} is not
the definition but a calculation of Floer homology of
periodic Hamiltonian system.
Namely it coincides with singular homology.
This fact is used in the proof of the homological version of Arnold's conjecture.
There are two methods to verify this calculation.
One uses the method of Bott-Morse theory, and the other is based on the study
the case of autonomous Hamiltonian that is $C^2$-small and Morse-Smale.
In Part \ref{S1equivariant},  we use the second method in the
present article following \cite{FOn}.
(The approach using Bott-Morse theory is written in \cite[Section 26]{fooospectr}.)
\par
The key point is to use the $S^1$ symmetry of the problem.
Namely when the Hamiltonian is time independent,
the moduli space we study has an extra $S^1$ symmetry arising from domain rotations.
Therefore  contribution of the relevant Floer moduli space to
the matrix elements of the resulting boundary operator is concentrated
to the fixed point set of the $S^1$-action which exactly corresponds to
the moduli space of Morse gradient flows.  This $S^1$ symmetry is used also in the
approach via the Bott-Morse theory.
\par
In Part \ref{S1equivariant} we define the notion of $S^1$-equivariant
Kuranishi structure and prove the $S^1$-equivariant counterparts of the
various results on the Kuranishi structure.
We also construct an $S^1$ equivariant Kuranishi structure
on the moduli space of the solutions of the Floer's equation
and use it to calculate Floer homology of
periodic Hamiltonian system.

\subsection{Epilogue}

The last part is a kind of appendix.
We have already mentioned that the origin of this article is our replies uploaded
for the discussion in the google group `Kuranishi' during which
we replied mainly to the questions raised by K. Wehrheim.
In the first three sections of Part \ref{origin}, we describe
the discussion of that google group and the role of the pdf files we posted there,
from our point of view.
\par
In the arXiv and even in the published literature, we have seen a few articles
which express some negative view on the foundation and raise some doubts
on the solidness of virtual fundamental chain or cycle techniques, although
they have been used for the various purposes successfully in
the published references. In our point of view, many such doubts raised in those
articles are not based on the precise understanding of the virtual
fundamental cycle techniques but based on some prejudice on the mathematical
point of view and on a few minor technical imprecise statements made in the published articles
on the virtual fundamental cycle techniques.
\par
Recently we have seen another instance of such a writing \cite{MW1}
in arXiv that is posted by the very person who have asked us questions
in the
google group and already gotten our answers, which are mostly the same
as Parts
\ref{Part2} - \ref{S1equivariant} in this article except some
polishing of
presentation.  They posed several difficulties, which arise in {\em their}
approach.
For example, Hausdorff property of certain spaces, smoothness of
obstruction bundles,
which we consciously excluded by taking the route via the finite
dimensional reduction.
(They should be taken care of, if one works with infinite dimensional
setting directly.)
We comment on \cite{MW1} more in
Section \ref{HowMWiswrong}.

\subsection{Thanks}
We would like to thank all the participants, especially Wehrheim and
McDuff,
  in the discussions of the `Kuranishi' google group for motivating us to
go through this painstaking labour by their meticulous reading and
questioning of our writings.

KF is supported by
JSPS Grant-in-Aid for Scientific Research \# 19104001, 2322404 and Global
COE Program G08.
YO is supported by
US NSF grant \# 0904197.
HO is supported by
JSPS Grant-in-Aid for Scientific Research \#23340015.
KO is supported by
JSPS Grant-in-Aid for Scientific Research \# 2124402.
\par
KF thanks Simons Center for Geometry and Physics for
hospitality and financial support while this work is done.
\par\newpage
\part{Kuranishi structure and virtual fundamental chain}\label{Part2}

The purpose of this part is to give the definitions of Kuranishi structure
and good coordinate system and to explain the construction of a virtual
fundamental chain of a space with Kuranishi structure
using a good coordinate system. We also provide the details of the proof of the existence of a good coordinate
system on a compact metrizable space with Kuranishi structure and the tangent bundle
in Section \ref{sec:existenceofGCS}.
We take the definition of \cite[Appendix]{fooo:book1}.

\section{Definition of Kuranishi structure}
\label{secdefnkura}

In this section we give the definition of Kuranishi structure.
We mostly follow the exposition and the notations used in \cite{fooo:book1}.
Here is a technical remark. In \cite{fooo:book1} we do not use the notion of germ of
Kuranishi neighborhoods, which was discussed in \cite{FOn}.
The notion of germ is not needed for the proofs of all the results in \cite{FOn}.
See Section \ref{gernkuranishi} about germ of Kuranishi neighborhood etc.

In particular, as in the exposition of \cite{fooo:book1}, the cocycle condition
$$
\underline{\phi}_{pq} \circ \underline{\phi}_{qr}=\underline{\phi}_{pr}
$$
is the exact equality and not the one modulo automorphism of a Kuranishi neighborhood.
(Note that there may be a non trivial automorphism of a Kuranishi neighborhood.)
This is important to avoid usage of 2-category. Here
\begin{equation}\label{eq:barphipq}
\underline{\phi}_{pq}: U_{pq} \to U_p
\end{equation}
is an embedding of the orbifold $U_{pq}=V_{pq}/\Gamma_q \to U_p =V_p/\Gamma_p$, that is
induced by the $h_{pq}$-equivariant map
\begin{equation}\label{eq:phipq}
\phi_{pq}:V_{pq} \to V_p
\end{equation}
where
\begin{equation}\label{eq:hpq}
h_{pq}: \Gamma_q \to \Gamma_p
\end{equation}
is a group homomorphism.
See below for the precise definitions of these notations.
We want to avoid using the language of 2-category unless it is absolutely necessary because
we feel that it makes things unnecessarily complicated and
that it is also harder to use.
(See Section \ref{ofd} where we summarize the notation and definition on orbifold
we use in this article.)
\par
Let $X$ be a compact metrizable space and $p \in X$.
We define a Kuranishi neighborhood at a point $p$ in $X$ as follows.
\begin{defn} (\cite[Definition A1.1]{fooo:book1})\label{Definition A1.1}
A {\it Kuranishi neighborhood} at $p \in X$ is a quintuple
$(V_p, E_p, \Gamma_p, \psi_p, s_p)$ such that:
\begin{enumerate}
\item $V_p$ is a
smooth manifold of finite dimension, which may or may
not have boundary or corner.
\item
$E_p$ is a real vector space of finite dimension. \par
\item $\Gamma_p$ is a finite group acting smoothly and effectively
on $V_p$ and has a linear representation
on $E_p$.
\item
$s_p$ is a $\Gamma_p$ equivariant smooth map $V_p \to E_p$. \par
\item $\psi_p$ is a homeomorphism from
$s_p^{-1}(0)/\Gamma_p$ to a neighborhood of $p$ in $X$.
\end{enumerate}
\begin{rem}
We {\it always} assume orbifolds to be effective.
In our application to the moduli space of pseudo-holomorphic curves
we can take obstruction spaces so that orbifold appearing in its
Kuranishi neighborhood is always effective,
except the case when the target space $X$ is zero dimensional.
\end{rem}
We denote $U_p = V_p/\Gamma_p$ and call $U_p$ a \emph{Kuranishi neighborhood}.
We sometimes also call $V_p$ a {\emph Kuranishi neighborhood} of $p$ by an abuse of
terminology.
\par
We call $E_p \times V_p \to V_p$ an {\it obstruction bundle} and
$s_p$ a {\it Kuranishi map}.
For $x \in V_p$, denote by $(\Gamma_p)_x$ the isotropy subgroup at $x$, i.e.,
$$(\Gamma_p)_x = \{ \gamma \in \Gamma_p \,\vert \,\gamma x = x \}.$$
Let $o_p$ be a point in $V_p$ with $s_p(o_p) = 0$
and $\psi_p([o_p]) = p$. We will assume that $o_p$ is fixed by all elements of
$
\Gamma_p
$.
Therefore $o_p$ is a unique point of $V_p$ which is mapped to $p$ by $\psi_p$.
\end{defn}
\begin{defn}(\cite[Definition A1.3]{fooo:book1})\label{Definition A1.3}
Let $(V_p, E_p, \Gamma_p, \psi_p, s_p)$ and
$(V_q, E_q, \Gamma_q, \psi_q, s_q)$ be a pair of Kuranishi neighborhoods of
$p \in X$ and $q \in \psi_p(s_p^{-1}(0)/\Gamma_p)$, respectively.
We say a triple $\Phi_{pq} = (\hat\phi_{pq},\phi_{pq},h_{pq})$
a {\it coordinate change} if
\begin{enumerate}
\item[(1)]
$h_{pq}$ is an injective homomorphism $\Gamma_q \to \Gamma_p$.
\item[(2)]
$\phi_{pq} : V_{pq} \to V_p$ is an
$h_{pq}$ equivariant smooth embedding
from a $\Gamma_q$ invariant open neighborhood $V_{pq}$  of $o_q$ to $V_p$,
such that the induced map
$\underline{\phi}_{pq}: U_{pq} \to U_p$ is injective. Here and hereafter $\underline{\phi}_{pq} : U_{pq} \to U_q$ is a map induced
by ${\phi}_{pq}$ and $U_{qp} = V_{qp}/\Gamma_p$.
\item[(3)]
$(\hat\phi_{pq},\phi_{pq})$ is an $h_{pq}$-equivariant embedding of vector
bundles $E_q  \times V_{pq} \to E_p \times V_p$.
\end{enumerate}
In other words, the triple $\Phi_{pq}$ induces an embedding of orbibundles
$$
\underline{\hat\phi}_{pq} : \frac{E_q\times V_{pq}}{\Gamma_q} \to \frac{E_p\times V_{p}}{\Gamma_p},
$$
in the sense of Definition \ref{defn:embedding}.
\par
The collections $\Phi_{pq}$ satisfy the following compatibility conditions.
\begin{enumerate}
\item[(4)] $\hat\phi_{pq}\circ s_q =
s_p\circ\phi_{pq}$. Here and hereafter we sometimes regard
$s_p$ as a section $s_p :  V_p \to E_p\times V_p$ of trivial bundle
$E_p \times V_p \to V_p$.
\item[(5)]
$\psi_q =
\psi_p\circ \underline{\phi}_{pq}$ on $(s_q^{-1}(0) \cap V_{pq})/\Gamma_q$.
\item[(6)]
The map $h_{pq}$ restricts to an isomorphism
$(\Gamma_q)_x \to (\Gamma_p)_{\phi_{pq}(x)}$ for any
$x \in V_{pq}$.
Here
$$
(\Gamma_q)_x = \{ \gamma \in \Gamma_{q} \mid \gamma x = x\}.
$$
\end{enumerate}
\end{defn}

\begin{defn}\label{Definition A1.5}(\cite[Definition A1.5]{fooo:book1})
A {\it Kuranishi structure} on $X$ assigns a
Kuranishi neighborhood $(V_p, E_p, \Gamma_p, \psi_p, s_p)$
for each $p \in X$ and a coordinate change
$(\hat\phi_{pq},\phi_{pq},h_{pq})$ for each
$q \in \psi_p(s_p^{-1}(0)/\Gamma_p)$ such that the following holds.
\begin{enumerate}
\item $\dim V_p - \operatorname{rank} E_p$
is independent of $p$.
\par
\item
If $r \in \psi_q((V_{pq}\cap s_q^{-1}(0))/\Gamma_q)$,
$q \in \psi_p(s_p^{-1}(0)/\Gamma_p)$ then there exists $\gamma_{pqr}^{\alpha} \in \Gamma_p$
for each connected component
$(\phi_{qr}^{-1}(V_{pq}) \cap V_{qr} \cap V_{pr})_\alpha$ of $\phi_{qr}^{-1}(V_{pq})
\cap V_{qr} \cap V_{pr}$ such that
$$
h_{pq} \circ h_{qr} = \gamma_{pqr}^{\alpha}
\cdot h_{pr} \cdot (\gamma_{pqr}^{\alpha})^{-1}
, \quad
\phi_{pq} \circ \phi_{qr} = \gamma_{pqr}^{\alpha}
\cdot \phi_{pr}, \quad
\hat\phi_{pq} \circ \hat\phi_{qr} = \gamma_{pqr}^{\alpha}\cdot
\hat\phi_{pr}.
$$
Here the first equality holds on $(\phi_{qr}^{-1}(V_{pq}) \cap V_{qr} \cap V_{pr})_\alpha$
and the second equality holds on $E_r \times (\phi_{qr}^{-1}(V_{pq}) \cap V_{qr} \cap V_{pr})_\alpha$.
\end{enumerate}
\end{defn}

Next we recall that, for a section $s$ of a vector bundle $E$ on a manifold $V$,
the linearization of $s$ induces a canonical map from the restriction of the tangent bundle to the
zero set $s^{-1}(0)$ to $E\vert_{s^{-1}(0)}$.

We note that the differential $d_{\text{fiber}}s_p$ of the Kuranishi map induces a bundle map
\begin{equation}\label{eq:tangent}
d_{\text{fiber}}s_p ~:~ N_{V_{pq}}V_p \to \frac{\hat{\phi}_{pq}^*(E_p \times V_{p})}{E_q\times {V_{pq}}}
\end{equation}
as $\Gamma_q$-equivariant bundles on $V_{pq} \cap s_q^{-1}(0)$,
and a commutative diagram
\begin{equation}
\begin{CD}
0 @ >>> T_{x}V_{pq}  @ > {d_x\phi_{pq}} >>
T_{\phi_{pq}(x)}V_{p}  @ >>>  (N_{V_{pq}}V_p)_x @ > >> 0\\
&& @ V{ds_{q}}VV @ V{ds_p}VV  @VVV\\
0 @ >>> (E_q)_x @ >{\hat{\phi}_{pq}}>> (E_p)_{\phi_{pq}(x)}
@>>> \frac{(E_p)_{\phi_{pq}(x)}}{(E_q)_x}@>>> 0
\end{CD}
\end{equation}

\begin{defn}[Tangent bundle]\label{defn:tangentbundle}
\label{defn:Kura-tangent} We say that a {\it space with Kuranishi structure $(X,\KK)$
has a tangent bundle} if the map \eqref{eq:tangent}
is an $h_{pq}$-equivariant bundle isomorphism on $V_{pq} \cap s_q^{-1}(0)$.
\par
We say it is orientable if the bundles
$$
\det E_q^* \otimes \det TV_q\Big|_{s_q^{-1}(0) \cap \psi_q^{-1}(U_{pq})}
$$
has trivializations compatible with the isomorphisms \eqref{eq:tangent}.
We call a space with Kuranishi structure and tangent bundle a {\it Kuranishi space}.
\end{defn}

\begin{defn}\label{Definition A1.13}
Consider the situation of Definition \ref{Definition A1.5}. Let $Y$ be a topological space.
A family $\{f_p\}$ of $\Gamma_p$-equivariant continuous maps $f_p : V_p \to Y$ is said to be a
{\it strongly continuous map} if
$$
f_p \circ \phi_{pq} = f_q
$$
on $V_{pq}$.  A strongly continuous map induces a continuous map $f: X \to
Y$.
We will ambiguously denote $f = \{f_p\}$ when the meaning is clear.

When $Y$ is a smooth manifold, a strongly continuous map
$f: X \to Y$ is defined to be smooth if all $f_p:V_p \to Y$ are smooth.
We say that it is {\it weakly submersive} if each $f_p$ is a submersion.
\end{defn}

\section{Definition of good coordinate system}
\label{defgoodcoordsec}

The construction of multivalued perturbation (multisection)
is by induction on the coordinate.
For this induction to work we need to take a clever choice of the
coordinates we work with. Such a system of coordinates is
called
a \emph{good coordinate system}, which was introduced in \cite{FOn}.

In this section we define the notion of good coordinate system following \cite{FOn}. But,
as an abstract framework, we require some additional conditions, for example, Condition \ref{Joyce} due to Joyce.
In \cite{FOn} we shrank Kuranishi neighborhoods several times.
To describe this procedure in great detail,
we use these additional conditions.

\begin{defn}\label{ofddef1}
An orbifold is a special case of Kuranishi space where all the obstruction bundles
are trivial, i.e., $E_p=0$ for all $p \in X$.
\end{defn}

\begin{rem}
If we try to define the notion of morphisms between Kuranishi spaces (the space equipped with
Kuranishi structure) it seems necessary to systematically work in 2-category.
When one is interested in Kuranishi spaces on its own, this is
certainly more natural approach to study.
(This is the approach taken by D. Joyce \cite{joyce2} we suppose.)
Our purpose is to use the notion of Kuranishi structure as a method of defining various
invariants by using Kuranishi structure and abstract perturbation,
and implement them to symplectic geometry and/or mirror symmetry etc.
By this reason, we take a way that is as short as possible and also is as general as possible
at the same time, for our particular purpose.
\end{rem}
See Section \ref{ofd} about our terminology on orbifolds.
We will define there the notion of orbifolds, vector bundle on it, and its embedding.

Hereafter we denote
\begin{equation}
\mathcal U_{p} = \psi_{p}(s_{p}^{-1}(0)/\Gamma_{p})
\end{equation}
which defines an open neighborhood of $p$ in $X$ by the assumption on $\psi_{p}$.

We modify the definition of good coordinate system in
\cite[Lemma 6.3]{fooo:book1} as follows.
In our definition of good coordinate system
our `Kuranishi neighborhood' is an orbifold which may
not be a global quotient.
The definition is given
in terms of orbifolds and the obstruction bundles which are
orbibundles in general.
From now on, we call an orbibundle simply a vector bundle.

\begin{defn}\label{goodcoordinatesystem}
Let $X$ be a space with Kuranishi structure.
A \emph{good coordinate system} on it consists
of a partially ordered set  $(\frak P,\le)$ of finite order,
and $(U_{\frak p}, E_{\frak p}, \psi_{\frak p}, s_{\frak p})$
for each $\frak p \in \frak P$, with the following data.
\begin{enumerate}
\item $U_{\frak p}$ is an
orbifold of finite dimension, which may or may
not have boundary or corner.
\item
$E_{\frak p}$ is a real vector bundle over $U_{\frak p}$.
\item
$s_{\frak p}$ is a section of $E_{\frak p} \to U_{\frak p}$.
\item $\psi_{\frak p}$ is an open embedding of $s_{\frak p}^{-1}(0)$ into $X$.
\item If $\frak q \le \frak p$, then there exists an
embedding of vector bundles
$$
\hat{\underline\phi}_{\frak p \frak q}:E_{\frak q}\vert_{U_{\frak p\frak q}} \to E_{\frak p}
$$
over an embedding ${\underline\phi}_{\frak  p\frak q}: U_{\frak p \frak q} \to U_{\frak p}$
of orbifold such that
\begin{enumerate}
\item
$U_{\frak p\frak q}$ is an open subset of $U_{\frak q}$ such that
\begin{equation}\label{eq23}
\psi_{\frak q}(U_{\frak p\frak q}
\cap s_{\frak q}^{-1}(0))
= \psi_{\frak p}(U_{\frak p}\cap s_{\frak p}^{-1}(0)) \cap \psi_{\frak q}(U_{\frak q}\cap s_{\frak q}^{-1}(0)),
\end{equation}
\item $\widehat{\underline\phi}_{\frak p\frak q}\circ s_\frak q = s_\frak p \circ \underline \phi_{\frak p\frak q}$,
\, $\psi_\frak q = \psi_\frak p\circ \underline \phi_{\frak p\frak q}$,
\item
$d_{\text{fiber}}s_{\frak p}$ induces an isomorphism of vector bundles at $s_{\frak q}^{-1}(0)\cap U_{\frak p\frak q}$.
\begin{equation}\label{15tangent}
N_{U_{\frak p\frak q}}U_\frak p \cong \frac{{\underline\phi}_{\frak p\frak q}^*E_\frak p}
{(E_\frak q)\vert_{U_{\frak p\frak q}}}.
\end{equation}
\end{enumerate}
\item
If $\frak r \le \frak q \le \frak p$,
$\psi_{\frak p}(s_{\frak p}^{-1}(0)) \cap \psi_{\frak q}(s_{\frak q}^{-1}(0))
\cap \psi_{\frak r}(s_{\frak r}^{-1}(0)) \ne \emptyset$, we have
$$
\underline\phi_{\frak p\frak q} \circ \underline\phi_{\frak q\frak r} =
\underline\phi_{\frak p\frak r}, \quad
\hat{\underline\phi}_{\frak p\frak q} \circ \hat{\underline\phi}_{\frak q\frak r} =
\hat{\underline\phi}_{\frak p\frak r}.
$$
Here the first equality holds on ${\underline\phi}_{\frak q\frak r}^{-1}(U_{\frak p\frak q}) \cap U_{\frak q\frak r}
\cap U_{\frak p\frak r}$, and the
second equality holds on
$(E_{\frak r})\vert_{({\underline\phi_{\frak q\frak r}^{-1}(U_{\frak p\frak q}) \cap U_{\frak q\frak r}}\cap U_{\frak p\frak r})}$.
\par\medskip
\item
$$
\bigcup_{\frak p\in \frak P} \psi_{\frak p}(s_{\frak p}^{-1}(0)) = X.
$$
\item
If
$\psi_{\frak p}(s_{\frak p}^{-1}(0)) \cap \psi_{\frak q}(s_{\frak q}^{-1}(0)) \ne \emptyset$,
either $\frak p\le \frak q$ or $\frak q \le \frak p$ holds.
\item The Conditions \ref{Joyce}, \ref{plusalpha}, \ref{plusalpha2} and
\ref{proper} below hold.
\end{enumerate}
\end{defn}
\begin{rem}
For the definition of good coordinate system given here,
we assume more conditions than those given in \cite{FOn}.
We use them to more explicitly describe the process
of shrinking the Kuranishi neighborhoods entering in the construction of multisections.
\par
On the other hand, we prove the existence of such a restrictive good coordinate
system assuming the existence of Kuranishi structure with \emph{the same definition} as the one given in
\cite{fooo:book1}.
Therefore the conclusion which is the existence of virtual
fundamental chain or cycle associated to the Kuranishi structure
is the same as \cite{fooo:book1}.
\end{rem}
\begin{conds}[Joyce \cite{joyce}] \label{Joyce}
Suppose $\frak p \ge \frak q \ge \frak r$.
$$
\underline\phi_{\frak p\frak q}(U_{\frak p\frak q})
\cap
\underline\phi_{\frak p\frak r}(U_{\frak p\frak r})
=
\underline\phi_{\frak p\frak r}
(\underline\phi_{\frak q\frak r}^{-1}(U_{\frak p\frak q}))
\cap
U_{\frak p\frak r}).
$$
\end{conds}

\begin{lem}\label{iikaejoyce}
Condition \ref{Joyce} is equivalent to the following statement:
\par
If $\frak p \ge \frak q \ge \frak r$ and $x \in U_{\frak p\frak r}$,
$y \in U_{\frak p\frak q}$ with
$\underline{\phi}_{\frak p\frak r}(x) = \underline{\phi}_{\frak p\frak q}(y)$,
then
\begin{enumerate}
\item $x \in
\underline \phi_{\frak q\frak r}^{-1}(U_{\frak p\frak q}) \cap U_{\frak q\frak r}$,
\item $\underline{\phi}_{\frak q\frak r}(x) = y$.
\end{enumerate}
\end{lem}
\begin{proof} This is obvious.
\end{proof}

\begin{conds}\label{plusalpha}
\begin{enumerate}
\item If
$
\bigcap_{i\in I} U_{\frak p_i\frak q} \ne \emptyset,
$
then
$
\bigcap_{i\in I} \mathcal U_{\frak p_i\frak q} \ne \emptyset.
$
\item If
$
\bigcap_{i\in I} \underline{\phi}_{\frak p\frak q_i}(U_{\frak p\frak q_i}) \ne \emptyset,
$
then
$
\bigcap_{i\in I} \mathcal U_{\frak p\frak q_i} \ne \emptyset.
$
\end{enumerate}
\end{conds}
Here and hereafter we put
\begin{equation}
\mathcal U_{\frak p\frak q} = \psi_{\frak q}(s_{\frak q}^{-1}(0)\cap U_{\frak p\frak q}).
\end{equation}
Condition \ref{plusalpha} and Definition \ref{goodcoordinatesystem} (8)  imply the following:
\begin{lem}\label{intersection}
Suppose $\frak q \le \frak p_j$ for $ j =1,\dots, J$ and
$
\bigcap_{i\in I} U_{\frak p_i\frak q} \ne \emptyset
$.
Then the set $\{\frak q\} \cup \{\frak p_j \mid j =1,\dots, J\}$ are linearly ordered.
(Namely for each  $\frak r,\frak s \in \{\frak q\} \cup \{\frak p_j \mid j =1,\dots, J\}$ at least one of
$\frak r \ge \frak s$ or $\frak s \ge \frak r$ holds.)
\end{lem}

\begin{conds}\label{plusalpha2}
Suppose $U_{\frak p\frak r}\cap U_{\frak q\frak r}
\ne \emptyset$ or
${\underline \phi}_{\frak q\frak r}^{-1}
(U_{\frak p\frak q}) \ne \emptyset$. If
$\frak p \ge \frak q\ge \frak r$ in addition, then we have
$$
{\underline \phi}_{\frak q\frak r}^{-1}
(U_{\frak p\frak q})
=
U_{\frak p\frak r} \cap U_{\frak q\frak r}.
$$
\end{conds}

\begin{conds}\label{proper}
The map
$
U_{\frak p\frak q} \to U_{\frak p} \times U_{\frak q}
$
defined below is proper.
\begin{equation}
x \mapsto ({\underline \phi}_{\frak p\frak q}(x),x).
\end{equation}
\end{conds}
The existence of good coordinate system is proved in Section \ref{sec:existenceofGCS}.
In the rest of this section, we introduce an equivalence relation on
the disjoint union
$$
\widetilde U(X;\frak P)  = \bigcup_{\frak p\in \frak P} U_{\frak p}
$$
and a quotient space $U(X;\frak P)$ thereof. We may use the set $U(X;\frak P)$ as a `global
thickening' of $X$, in which a perturbation of the zero set $s^{-1}(0)$ will reside.\footnote{Actually we do {\it not} need to use such a space to define
virtual fundamental chain $f_*([X])$. We may take simplicial decomposition of the zero set
$s^{-1}_{\frak p}(0)$
(after appropriately shrinking the domain) so that they are compatible
with the coordinate change, by an induction on the partial order
of the set $\frak P$, and can use it to define $f_*([X])$, instead.
Existence of `appropriate shrink' is intuitively clear.
(See \cite[Answer to Question 3]{Fu1}, for example.)
However writing this intuitive picture in detail without introducing
a formal definition is rather cumbersome. (The details we provide
in this article are more  than required in common research papers, according to our
opinion.)
This is the reason why we choose to define the set $U(X)$ explicitly.}
For the simplicity of notations, we omit the dependence on $\frak P$
of $\widetilde U(X;\frak P), \,  U(X;\frak P)$ from their notations and
just denote them by $\widetilde U(X), \, U(X)$ respectively.

\begin{lem}\label{equiv}
The following  $\sim$ is an equivalence
relation on $\widetilde U(X)$.
\par\medskip
Let $x \in U_{\frak p}$ and $y \in U_{\frak q}$.
We say $x\sim y$ if and only if
\begin{enumerate}
\item
$x = y$, or
\item
$\frak p \ge \frak q$ and $\underline\phi_{\frak p\frak q}(y) = x$, or
\item
$\frak q \ge \frak p$ and $\underline\phi_{\frak q\frak p}(x) = y$.
\end{enumerate}
\end{lem}
\begin{proof}
Only transitivity is nontrivial.
Let $x_1 \sim x_2$, $x_2 \sim x_3$, $x_i \in U_{\frak p_i}$.
\par
Suppose $\frak p_1 \le \frak p_2 \le \frak p_3$.
Condition \ref{plusalpha2} implies $U_{\frak p_3\frak p_1} \cap U_{\frak p_2\frak p_1}
= (\underline\phi_{\frak p_2\frak p_1})^{-1}(U_{\frak p_3\frak p_2})$.
Therefore Definition \ref{goodcoordinatesystem} (6) implies
$$
x_3 = \underline\phi_{\frak p_3\frak p_2}(x_2)  =
\underline\phi_{\frak p_3\frak p_2} \underline\phi_{\frak p_2\frak p_1}(x_1)
=  \underline\phi_{\frak p_3\frak p_1}(x_1).
$$
Namely $x_3 \sim x_1$.
The case $\frak p_1 \ge \frak p_2 \ge \frak p_3$ is similar.
\par
Suppose  $\frak p_1 \ge \frak p_2 \le \frak p_3$.
Condition \ref{plusalpha} and Definition \ref{goodcoordinatesystem}
(8) imply either $\frak p_1 \ge \frak p_3$ or $\frak p_1 \le \frak p_3$.
Let us assume $\frak p_1 \le \frak p_3$.
Then Condition \ref{plusalpha2}  implies
$x_2 \in \underline\phi_{\frak p_1\frak p_2}^{-1}(U_{\frak p_3\frak p_1} )$.
Then  Definition \ref{goodcoordinatesystem} (6) implies
$$
\underline\phi_{\frak p_3\frak p_1}(x_1)
=
\underline\phi_{\frak p_3\frak p_1}(\underline\phi_{\frak p_1\frak p_2}(x_2))
= x_3.
$$
Namely $x_1 \sim x_3$.
The case $\frak p_1 \ge \frak p_3$ is similar.
\par
Let us assume  $\frak p_1 \le \frak p_2 \ge \frak p_3$.
By Condition \ref{plusalpha} and Definition \ref{goodcoordinatesystem}
(8), we have either  $\frak p_1 \le \frak p_3$ or $\frak p_1 \ge \frak p_3$.
Then  Condition \ref{Joyce} implies
$
x_3 = \underline\phi_{\frak p_3\frak p_1}(x_1)$ or $x_1 = \underline\phi_{\frak p_1\frak p_3}(x_3),
$
as required.
\end{proof}
\begin{defn}\label{clrelation}
We define $U(X)$ to be the set of $\sim$ equivalence classes.
\par
The map $\Pi_\frak p : U_\frak p \to U(X)$ sends an element of $U_{\frak p}$ to its equivalence class.
The map $\psi_\frak p^{-1}: \psi_\frak p(s_{\frak p}^{-1}(0) \cap U_\frak p) \to s_{\frak p}^{-1}(0) \cap U_\frak p$
followed by the restriction of $\Pi_\frak p$ to $s_{\frak p}^{-1}(0) \cap U_\frak p$
defines an injective map
\begin{equation}\label{eq:iotap}
\iota_\frak p: \psi_\frak p(s_{\frak p}^{-1}(0) \cap U_\frak p) \to U(X).
\end{equation}

\begin{lem}\label{lem:glueiotap} The maps \eqref{eq:iotap} are consistent on
the overlaps. We denote the resulting global map by $I: X \to U(X)$.
\end{lem}
\begin{proof} Let $p \in X$ and suppose
\begin{equation}\label{eq:p=}
p = \psi_\frak p(x) = \psi_{\frak q}(y)
\end{equation}
for $x \in s_{\frak p}^{-1}(0) \cap U_\frak p$ and $y \in s_{\frak q}^{-1}(0) \cap U_{\frak q}$.
It is enough to prove $x \sim y$. By Definition \ref{goodcoordinatesystem} (6), either $\frak p = \frak q$,
$\frak p < \frak q$ or $\frak q < \frak p$.

If $\frak q = \frak p$, we must have $x = y$ since $\psi_\frak p: s_{\frak p}^{-1}(0) \cap U_\frak p
\to X$ is one-one. For the remaining two cases,
we will focus on the case $\frak q < \frak p$ since the other case is the same.
Then we are given an embedding of orbifolds $\underline \phi_{\frak p\frak q}:U_{\frak p\frak q} \to U_\frak p$.
Then Definition \ref{goodcoordinatesystem} (3-b),
$
\psi_\frak q = \psi_\frak p \circ \underline \phi_{\frak p\frak q}
$ on $U_{\frak p\frak q}$.

On the other hand, it follows from \eqref{eq23} and \eqref{eq:p=} that
$p = \psi_\frak q(x')$ with $x' \in s_{\frak q}^{-1}(0) \cap U_{\frak p\frak q} \subset
s_{\frak p}^{-1}(0) \cap U_\frak q$.
Since $\psi_\frak q$ is one-one on $s_{\frak q}^{-1}(0) \cap U_\frak q$, $x' = y$.

Then we derive
$$
\psi_\frak p(\phi_{\frak p\frak q}(y)) = \psi_\frak q(y) = p
$$
from Definition \ref{goodcoordinatesystem} (3-b). Since $\psi_\frak p (x) = p$ as well
and $\psi_\frak p$ is one-one on $s_{\frak p}^{-1}(0) \cap U_\frak q$, we obtain
$x = \phi_{\frak p\frak q}(y)$. This proves $x \sim y$, which finishes the proof of
the lemma.
 \end{proof}

\end{defn}
\begin{prop}\label{goodcoordinateprop2}
Suppose we have a good coordinate system.
Then there exist open subsets $U'_{\frak p} \subset U_{\frak p}$
and $U'_{\frak p\frak q} \subset U_{\frak p\frak q}$
such that the restrictions to $U'_{\frak p}$ and $U'_{\frak p\frak q}$ give a good coordinate system,
and $U'_{\frak p}$ and $U'_{\frak p\frak q}$ are relatively compact in
$U_{\frak p}$ and $U_{\frak p\frak q}$, respectively.
\end{prop}
\begin{proof}
We take an open subset $U'_{\frak p} \subset U_{\frak p}$  for each $\frak p$ that is relatively
compact in $U_{\frak p}$ and
$$
\bigcup_{\frak p\in \frak P}\psi_{\frak p}(s_{\frak p}^{-1}(0) \cap U'_{\frak p}) = X.
$$
We may choose it so that
\begin{equation}\label{for04}
\bigcap_{i\in I} \mathcal U_{\frak p_i} \ne \emptyset
\,\,\,
\Leftrightarrow
\,\,\,
\bigcap_{i\in I} \mathcal U'_{\frak p_i} \ne \emptyset.
\end{equation}
We put
\begin{equation}\label{defUprime}
U'_{\frak p\frak q}
= U_{\frak p\frak q} \cap U'_{\frak q}
\cap \underline\phi_{\frak p\frak q}^{-1}(U'_{\frak p}).
\end{equation}
Condition \ref{proper} implies that $U'_{\frak p\frak q}$ is
relatively compact.
It is straightforward to check that they satisfy the
conditions in Definition \ref{goodcoordinatesystem}.
\end{proof}
\begin{rem}\label{chooseR}
\begin{enumerate}
\item
If a compact subset $\mathcal K_{\frak p}$ of $\mathcal U_\frak p$ is given for each $\frak p$, then
we may choose $U'_\frak p$ etc. in Proposition \ref{goodcoordinateprop2} so that $\mathcal U'_\frak p$
contains $\mathcal K_\frak p$.
\item
On the other hand, we may choose $U'_{\frak p}$ as small as we want
as far as the condition $\bigcup_{\frak p \in \frak P}\mathcal U'_{\frak p}= X$ is satisfied.
In fact at the beginning of the proof we take $U'_{\frak p}$
so that this is satisfied and do not need to change it.
\item
In the case of the good coordinate system we produce in Section \ref{sec:existenceofGCS},
the index set $\frak P$ is a subset of the set of natural number with
obvious order $<$. So it is in fact linearly ordered.
Some of the combinatorial problem we have taken care of above is simpler in that case.
\end{enumerate}
\end{rem}
We define $U'(X)$ from $U'_{\frak p}$ and $U'_{\frak p\frak q}$
in the same way as Definition \ref{clrelation}.
Let $K_{\frak p}$ be the closure of $U'_{\frak p}$ in $U_{\frak p}$
that is compact.
We have:
\begin{equation}
\bigcup_{\frak p \in \frak P}\psi_{\frak p}(s_{\frak p}^{-1}(0) \cap K_{\frak p}) = X.
\end{equation}
For each $\frak p > \frak q$ we put
$$
K_{\frak p\frak q} = \underline\phi_{\frak p\frak q}^{-1}(K_{\frak p}) \cap K_{\frak q}.
$$
Since $\underline\phi_{\frak p\frak q}$ is proper it follows that
$K_{\frak p\frak q}$ is compact.
In the same way as the proof of Lemma \ref{equiv} $K_{\frak p\frak q}$ and
the restriction of $\underline\phi_{\frak p\frak q}$ to it induce an equivalence relation.
So we are in the following situation.
\begin{assump}\label{Kassumption}
\begin{enumerate}
\item
$(\frak P,\le)$ is a finite partial ordered set.
\item
$K_{\frak p}$ is a Hausdorff and compact set for each $\frak p$.
\item
For $\frak p,\frak q \in \frak P$, $\frak p > \frak q$,
$K_{\frak p\frak q} \subset K_{\frak q}$ is a compact subset
and
$\underline\phi_{\frak p\frak q} : K_{\frak p\frak q}  \to K_{\frak p}$
is an embedding.
\item
We define a relation $\sim$ on $\widetilde K(\frak P) = \bigcup_{\frak p \in \frak P} K_{\frak p}$
(disjoint union) as follows.
Let $x \in K_{\frak p}$ and $y \in K_{\frak q}$.
We say $x\sim y$ if and only if
\begin{enumerate}
\item
$x = y$, or
\item
$\frak p \ge \frak q$ and $\underline\phi_{\frak p\frak q}(y) = x$, or
\item
$\frak q \ge \frak p$ and $\underline\phi_{\frak q\frak p}(x) = y$.
\end{enumerate}
Then $\sim$ is an equivalence relation.
\end{enumerate}
\end{assump}
Let $K(\frak P)$ be the set of the $\sim$ equivalence classes of
$\widetilde K(\frak P)$.
We put the quotient topology.
\par
In this situation the following holds.
\begin{prop}\label{metrizable}
In addition to Assumption \ref{Kassumption}, assume that $K_{\frak p}$ satisfies the second axiom of countability
and locally compact. Then $K(\frak P)$ is metrizable.
(In particular it is Hausdorff.)
\end{prop}
The proof of this proposition will be given in Section \ref{gentoplem}. We continue our discussion.
\begin{defn}\label{Imapdef}
We define a map $\frak J_{K(\frak P)U'(X)} : U'(X) \to K(\frak P)$ by sending the ${\sim }$-equivalence class
$[x]$ of $x \in \widetilde U'(X)$ to the equivalence class of $x \in \widetilde U(X)$ in $K(\frak P)$.
\end{defn}
By the definition (\ref{defUprime}) of $U'_{\frak p\frak q}$,
 we find that if $\tilde x \sim \tilde y$ in $\widetilde U(X)$ for
$\tilde x, \tilde y \in \tilde U'(X)$ then $\tilde x \sim \tilde y$ in $\widetilde U'(X)$.
Therefore $\frak J_{K(\frak P)U'(X)}$ is injective.
\begin{defn}
We equip $U'(X)$ with the weakest topology of the map $\frak J_{K(\frak P)U'(X)} : U'(X) \to K(\frak P)$
(with respect to the topology of $K(\frak P)$). We simply call this kind of topology
appearing henceforth \emph{weak topology}.
\end{defn}

The weak topology on $U'(X)$ is nothing but the subspace topology of $K(\frak P)$ if
we identify $U'(X)$ with its image of the injective map $\frak J_{K(\frak P)U'(X)}$
in $K(\frak P)$. Hereafter we use this topology on $U'(X)$ only, not the quotient topology of
$
U'(X) = \widetilde U'(X)/\sim
$,
unless otherwise mentioned explicitly.

\begin{rem}\label{rem520} The weak topology of $U'(X)$ is Hausdorff since $K(\frak P)$ is
Hausdorff and the map $\frak J_{K(\frak P)U'(X)} : U'(X) \to K(\frak P)$ is injective.
This topology is certainly different from the quotient topology thereof as the set of
equivalence classes in $U'(X)$. We sometimes also call this topology the induced topology
as long as it is not ambiguous. In fact the quotient topology does not necessarily satisfy the first axiom of
countability as pointed out by \cite[Example 6.1.14]{MW1}.
\end{rem}
\begin{cor}\label{corollary521}
$U'(X)$ (equipped with the weak topology) is metrizable.
\end{cor}
\begin{proof}
This is a consequence of Proposition \ref{metrizable}.
\end{proof}
We start from $U'_{\frak p}$ and repeat the process.
Namely we take relatively compact subsets $U''_{\frak p}$
and define $U''(X)$. We use $K'_{\frak p} =$ the closure of
$U''_{\frak p}$ and define $K'(\frak P)$.
We define an injective map $\frak J_{U'(X)U''(X)}$ in the same way as
Definition \ref{Imapdef}. We equip
$U''(X)$ with the weak topology of the map $\frak J_{U'(X)U''(X)}$.

\begin{lem}\label{fraIUUhomeo}
$\frak J_{U'(X)U''(X)}$ is a topological embedding, i.e.,
a homeomorphism to its image.
\end{lem}
\begin{proof}
We define $\frak J_{K(\frak P)K'(\frak P)} : K'(\frak P) \to K(\frak P)$
in the same way. It is injective. Since $K'(\frak P)$ is compact
and $K(\frak P)$ is Hausdorff $\frak J_{K(\frak P)K'(\frak P)}$
is a topological embedding. We then have the diagram
$$
\xymatrix{U'(X) \ar[r] & K(\frak P) \\
U''(X) \ar[u] \ar[r] & K'(\frak P) \ar[u]}
$$
in which the top, the bottom and the right column arrows are all
topological embeddings. The lemma immediately follows from this.
\end{proof}

\begin{rem}
Hereafter we further shrink $U'_{\frak p}$ several times
by taking relatively compact subsets. We always equip them with
{\it the weak topology} of the relevant injective map. Lemma \ref{fraIUUhomeo} implies that
the weak topology induced from $K'(\frak P) = (\bigcup \overline U''_{\frak p})/\sim$ on
$U''_{\frak p}$ coincides with the one induced from $K(\frak P)$.
\end{rem}
We define a map $\Pi_{\frak p} : U'_{\frak p} \to U'(X)$ as before
and define the canonical injective map $I': X \to U'(X)$ as in Lemma \ref{lem:glueiotap}
applied to $U'(X)$.
\begin{lem}
\begin{enumerate}
\item $I': X \to U'(X)$ is a topological embedding.
\item $\Pi_{\frak p} : U'_{\frak p} \to U'(X)$ is a topological embedding.
\end{enumerate}
\end{lem}
\begin{proof} The well-definedness of $I$ is proved in Lemma \ref{lem:glueiotap}.
Statement (1) follows from the fact that $I$ is injective, $X$ is compact and $U'(X)$ is Hausdorff.
\par
We consider a map
$\Pi_{\frak p} : K_{\frak p} \to K(\frak P)$, which is defined in the same way.
This map is injective and continuous.
Moreover $K_{\frak p}$ is compact and $K(\frak P)$ is Hausdorff.
Therefore $\Pi_{\frak p} : K_{\frak p} \to K(\frak P)$
is a topological embedding.
Since $\Pi_{\frak p} : U'_{\frak p} \to U'(X)$
is its restriction and the topologies of $U'_{\frak p}$ and of $U'(X)$
are the weak topology, the lemma follows.
\end{proof}
Hereafter we write $I$ for $I' : X \to U'(X)$ also.
\begin{rem}
Hereafter we equip a metric with $U'(X)$ or with similar spaces obtained by restricting
$U'_{\frak p}$ to its relatively compact subsets several times.
We fix a metric on $K(\frak P)$ and the metric we use is always
the restriction of this metric. This metric is compatible with the weak topology on $U'(X)$.
Since $K(\frak P)$ is compact, the metric on it is unique up to equivalence.
(Here two metrics $d$ and $d'$ are said to be equivalent
if there exists a homeomorphism $\Phi_i : \R_{\ge 0} \to \R_{\ge 0}$
for $i=1,2$ such that
$$
\Phi_1(d(x,y)) \le d'(x,y) \le \Phi_2(d(x,y)).
$$
\end{rem}
\begin{lem}\label{openneighborhood}
Let $x = \psi_\frak p(\tilde x) \in X$ and $\tilde x \in s_\frak p^{-1}(0) \subset U'_\frak p$.
Then there exists a neighborhood $\frak O_\frak p(x)$ of $\tilde x$ in $U'_\frak p$ such that
\begin{enumerate}
\item
$\Pi_\frak p : \frak O_\frak p(x) \to U'(X)$ is a topological embedding.
\item
$\Pi_\frak p(\frak O_\frak p(x))$ is an open subset of
$
\bigcup_{\frak q\le \frak p} \Pi_\frak q(U'_\frak q).
$
\end{enumerate}
\end{lem}
\begin{proof}
Choose an open neighborhood $\frak O_\frak p(x) \subset U'_\frak p$ that is relatively compact in $U'_\frak p$.
Clearly the map $\Pi_\frak p : \frak O_\frak p(x) \to U'(X)$ is injective and continuous.
By the choice, it extends so to the closure of $\frak O_\frak p(x)$ that is compact.
Since $U'(X)$ equipped \emph{with weak topology} is Hausdorff,  (1) follows.
\par
We next prove (2). Let $\frak q \le \frak p$. Since $\Pi_\frak p : \frak O_\frak p(x) \to U'(X)$ is an embedding and its image is relatively compact,
we may assume  $\Pi_\frak p (O_\frak p(x)) \subset U''(X)$ by choosing $U''(X)$ appropriately
where $U''(X)$ is as in Lemma \ref{fraIUUhomeo}. Therefore it suffices to show that
$\Pi_\frak p (O_\frak p(x))$ is open in quotient topology of $U'(X)$.
(This is because $K'(\frak P) \to U'(X)$ is continuous with respect to the quotient topologies.)
For this purpose, it suffices to show
$$
(\Pi_\frak q)^{-1}(\Pi_\frak p(\frak O_\frak p(x)))
$$
is open in $U'_\frak q$.
In fact
$$
\underline{\phi}_{\frak p\frak q}^{-1}(\frak O_\frak p(x)) \cap U'_{\frak p\frak q}
=  \underline{\phi}_{\frak p\frak q}^{-1}(\frak O_\frak p(x)) \cap
\left(U_{\frak p\frak q} \cap U'_\frak q \cap  \underline{\phi}_{\frak p\frak q}^{-1}(U'_\frak p)\right)
= \underline{\phi}_{\frak p\frak q}^{-1}(\frak O_\frak p(x)) \cap U'_\frak q.
$$
But we have $\underline{\phi}_{\frak p\frak q}^{-1}(\frak O_\frak p(x)) \cap U'_\frak q
\subset \underline{\phi}_{\frak p\frak q}^{-1}(U'_\frak p) \cap U'_\frak q$
and hence by definition,
$$
\underline{\phi}_{\frak p\frak q}^{-1}(\frak O_\frak p(x)) \cap U'_\frak q
= \Pi_\frak q^{-1}(\Pi_\frak p(\frak O_\frak p(x))).
$$
Combining these, we have finished the proof.
\end{proof}

The next lemma plays a key role in the next section to show basic properties of
the virtual fundamental chain.
\begin{lem}\label{unionofofd}
For any $x \in X$ there exist $\frak q_{1}, \dots, \frak q_{m} \in \frak P$ with $\frak q_{1} \le \dots \le \frak q_{m}$
and open sets $\Omega_{\frak q_i}(x) \subset U''_{\frak q_i}$ with the following properties.
\begin{enumerate}
\item $x \in \mathcal U''_{\frak q_1}$ and $x \in \overline{\mathcal U''_{\frak q_i}}$ for $i=2,\dots,m$.
\item $\psi_{\frak q_1}^{-1}(x) \in \Omega_{\frak q_1}(x)$.
\item
$\psi_{\frak q_i}^{-1}(x) \in  \overline \Omega_{\frak q_i}(x) \setminus \Omega_{\frak q_i}(x)$ for $i>1$.
Here the closure is taken in $U_{\frak q_i}$.
\item
The map $\Pi_{\frak q_i} : \Omega_{\frak q_i}(x) \to U''(X)$ is a topological embedding.
\item
The union of the images of $\Pi_{\frak q_i} : \Omega_{\frak q_i}(x) \to U''(X)$ is a neighborhood of $I(X)$.
\item
$\dim U_{\frak q_1} < \dim U_{\frak q_i}$ for $i\ne 1$.
\end{enumerate}
\end{lem}
\begin{proof}
By Lemma \ref{intersection} there exists a maximal $\frak q_1$ such that
$x \in  \mathcal U''_{\frak q_1}$.
\begin{sublem}\label{linearordercotto}
There exists  $\frak q_m \in \frak P$ such that it is maximal in the set
$
\{ \frak q\in \frak P \mid x \in  \overline{\mathcal U''_\frak q}\}.
$
\end{sublem}
\begin{proof}
Let $\frak q,\frak q' \in
\{ \frak q\in \frak P \mid x \in  \overline{\mathcal U''_\frak q}\}.
$
Since the closure of $\mathcal U''_\frak q$ is contained in $\mathcal U'_\frak q$, it
follows that $\mathcal U_\frak q \cap \mathcal U_{\frak q'}
\ne \emptyset$. Therefore by (\ref{for04})
$\mathcal U''_\frak q \cap \mathcal U''_{\frak q'} \ne \emptyset$.
The sublemma follows from Definition \ref{goodcoordinatesystem} (8).
\end{proof}
By Sublemma \ref{linearordercotto}
we can take $\frak q_1,\dots,\frak q_m$ such that $\frak q_{1} \le \dots \le \frak q_{m}$
and
$$
\{ \frak q\mid x \in  \overline{\mathcal U''_\frak q} \} \cap
\{   \frak q\mid \dim U_{\frak q} > \dim U_{\frak q_1}  \}
=
\{\frak q_i \mid i=2,\dots,m
\}.
$$
Then we can use Lemma \ref{openneighborhood} to find the required $\Omega_{\frak q_i}(x) = \frak O_{\frak q_i} \cap U''_{\frak  q_i}$.
(We do not need $\frak q_i$ ($i\ne 1$) with $\dim U_{\frak q_1} = \dim U_{\frak q_i}$ since
$x$ is in (the interior of) $\Omega_{\frak q_1}(x)$.)
\end{proof}
We put
\begin{equation}\label{frakUdef}
\frak U(x) = \bigcup_{i=1,\dots,m}\Pi_{\frak q_i}(\Omega_{\frak q_i}(x)).
\end{equation}
We take the intersection of $\frak U(x)$ and a sufficiently small open neighborhood of $x$ in $U''(X)$. Then
the property (4) of Definition \ref{goodcoordinatesystem} implies that
it gives an open neighborhood of $x$ in $U''(X)$.
Such $\frak U(x)$ forms a neighborhood basis of $x$.
(See the proof of Sublemmas \ref{9sublem1}, \ref{9sublem2}.)
\begin{exm}\label{omeganeighex}
We consider the following three subsets $U''_i$ ($i=1,2,3$) of $\R^3$.
$
U''_1 =
$ $x$-axis,
$
U''_2 = \{(x,y,0) \mid x > -y^2\}$,
$
U''_3 = \{(x,y,z) \mid x >0\}.
$
We put $U''(X) = U''_1\cup U''_2 \cup U''_3$ and consider
$\vec 0 = (0,0,0) \in U''(X)$.
Then the neighborhood $\frak U(\vec 0)$ we described above
typically is obtained as follows:
We take
$\Omega_1 = \{(x,0,0) \mid \vert x\vert < 3\epsilon\}$,
$\Omega_2 = \{(x,y,0) \mid x > -y^2, \,\, x^2+y^2 < (2\epsilon)^2\}$,
$\Omega_3 = \{(x,y,z) \mid x> 0, x^2+y^2 + z^2< \epsilon^2\}$.
The union of these three sets is $\frak U(x)$.
\end{exm}

\section{Construction of virtual fundamental chain}

We start the construction of perturbation and
virtual fundamental chain.
We use the notion of multisection for this purpose.
(The method to use abstract multivalued perturbation
to define a virtual fundamental chain or cycle
for the moduli space of pseudo-holomorphic curve
was introduced in 1996 January by \cite{FuOn99I}.)
\par
Let us review the definition of multisection here.

We assume that a finite group $\Gamma$ acts on a manifold $V$ and a vector space $E$.
A symmetric group $\frak S_n$ of order $n!$ acts on the
product $E^n$ by
$$
\sigma(x_1,\ldots,x_n) = (x_{\sigma(1)},\ldots,x_{\sigma(n)}).
$$
Let $S^n(E)$ be the quotient space $E^n/\frak S_n$.
Then $\Gamma$ action on $E$ induces an action on $S^n(E)$.
The map $E^n \to E^{nm}$ defined by
$$(x_1,\cdots,x_n)
\mapsto (\underbrace{x_1,\ldots,x_1}_{\text{$m$
times}},\underbrace{x_2,\ldots,x_2}_{\text{$m$
times}},\cdots,\underbrace{x_n,\ldots,x_n}_{\text{$m$
times}})$$
induces a $\Gamma$ equivariant map $S^n(E) \to S^{nm}(E)$.
\par
\begin{defn}
An {\it {$n$-multisection}} $s$ of $\pi : E\times V \to V$ is a
$\Gamma$-equivariant map $V \to S^n(E)$. We say that it is {\it
liftable} if there exists
$\widetilde s =
(\widetilde{s}_1,\ldots,\widetilde{s}_n): V
\to E^n$ such that its composition with $\pi : E^n \to S^n(E)$ is $s$.
(We do not assume $\widetilde s$ to be $\Gamma$ equivariant.) Each
of $\widetilde{s}_1,\ldots,\widetilde{s}_n$ is said to be a {\it branch} of $s$.
\par
If $s : V \to S^n(E)$ is an $n$ multisection, then it induces an $nm$
multisection for each $m$ by composing it with $S^n(E) \to S^{nm}(E)$.
\par
An $n$ multisection $s$ is said to be {\it equivalent} to an $m$
multisection $s'$ if the induced $nm$ multisections coincide to each other.
An equivalence class by this equivalence relation is said to be a
{\it multisection}.
\par
A liftable multisection is said to be {\it transversal} to zero if each of its
branch is transversal to zero.
\par
A family of multisections $s_{\epsilon}$ is said to {\it converge}
to $s$ as $\epsilon \to 0$ if there exists $n$ such that
$s_{\epsilon}$ is represented by an $n$-multisection
$s_{\epsilon}^n$ and $s_{\epsilon}^n$ converges to a representative
\end{defn}
From now on we assume all the multisections are liftable unless
otherwise stated.
\begin{defn}
Let $U$ be an orbifold (which is not necessarily a global quotient)
and $E$ a vector bundle on $U$ in the sense of Definition
\ref{orbibundledef}.
A {\it multisection} of $E$ on $U$ is given
by $\{U_i\}$ and $s_i$ where:
\begin{enumerate}
\item
$U_i = V_i/\Gamma_i$ is a coordinate system of our orbifold $U$
in the sense of Definition \ref{ofddefn}.
\item
$E\vert_{U_i} = (E_i\times V_i)/\Gamma_i$.
\item
$s_i : V_i \to S^{n_i}(E_i)$ is an $n_i$ multisection of the restriction $E_i$ of our bundle $E$
to $U_i$.
\item
The restriction of $s_i$ to $U_i \cap U_j$ is equivalent to
the restriction of $s_j$ to $U_i \cap U_j$ for each $i,j$.
\end{enumerate}
We say  $(\{U_i\},\{s_i\})$ and $(\{U'_i\},\{s'_i\})$ define the same
multisection if the restriction of $s_i$ to $U_i\cap U'_j$
is equivalent to the restriction of $s'_j$ to $U_i\cap U'_j$
for each $i,j$.
\par
Liftability, transversality, and convergence can be defined in the same way
as above.
\end{defn}
We start the construction of perturbation (system of multisections of obstruction bundles.)
Using  Proposition \ref{goodcoordinateprop2}, we shrink the Kuranishi neighborhoods $U_{\frak p}$
as follows.
\par
First we take an extension of the subbundle $\hat{\underline\phi}_{\frak p\frak q}(E_{\frak q}
\vert_{U'_{\frak p\frak q}})$ of
$E_{\frak p}\vert {{\underline\phi}_{\frak p\frak q}(U'_{\frak p\frak q})}$
to its neighborhood in $U_{\frak p}$.\footnote{Here and hereafter we write $E_{\frak p}$ in place of
$(E_{\frak p} \times V_{\frak p})/\Gamma_{\frak p}$ for simplicity.}
We also fix a splitting
\begin{equation}\label{splittingsubbundle}
E_{\frak p}=
E_{\frak q} \oplus E_{\frak q}^{\perp}
\end{equation}
on a neighborhood of ${{\underline\phi}_{\frak p\frak q}(U'_{\frak p\frak q})}$.
We can take such an extension of the subbundle and
splitting since $U'_{\frak p\frak q}$ is a relatively compact open subset of $U_{\frak p\frak q}$ that is
a suborbifold of $U_{\frak p}$.
\par
Using the splitting (\ref{splittingsubbundle}), the normal differential
\begin{equation}\label{dfiber}
d_{\rm fiber}s_{\frak p}
:
N_{U'_{\frak p\frak q}}U'_{\frak p} \to
\frac{{\underline\phi}_{\frak p\frak q}^*E_{\frak p}}{E_{\frak q}}
\end{equation}
is defined. (Note that without fixing the splitting
$d_{\rm fiber}s_{\frak p}$ is well-defined only at $s_{\frak q}^{-1}(0)
\cap U_{\frak p\frak q}$.)
\par
By the definition of tangent bundle, the map (\ref{dfiber}) is a bundle isomorphism on
$s_{\frak q}^{-1}(0)
\cap U'_{\frak p\frak q}$.
We take an open neighborhood $\frak W'_{\frak p\frak q}$ of
$s_{\frak q}^{-1}(0) \cap U'_{\frak p\frak q}$ in $U'_{\frak q}$
so that  (\ref{dfiber}) remains to be a bundle isomorphism on
$\frak W'_{\frak p\frak q}$.
\par
We take $U''_{\frak q}$ for each $\frak q$ so that
$$
U''_{\frak q} \cap U'_{\frak p\frak q}
\subset \frak W'_{\frak p\frak q}
$$
for each $\frak p$.
(We can take such $U''_{\frak q}$ by Remark \ref{chooseR} (2).)
Thus from now on we may assume that (\ref{dfiber}) is an isomorphism
on $U_{\frak p\frak q}$.
\par
We start with this $U_{\frak p}$, $U_{\frak p\frak q}$ and
repeat the construction of the last section.
Namely we take
$U^{(n)}_\frak p, U^{(n)}_{\frak p\frak q}$ such that
\begin{enumerate}
\item
The conclusion of Proposition \ref{goodcoordinateprop2} is satisfied
when we replace
$U_\frak p, U_{\frak p\frak q}$ by $U^{(n-1)}_\frak p, U^{(n-1)}_{\frak p\frak q}$ and
$U^{\prime}_\frak p, U^{\prime}_{\frak p\frak q}$ by $U^{(n)}_\frak p, U^{(n)}_{\frak p\frak q}$.
\item
The conclusions of Lemmas \ref{iikaejoyce}-\ref{unionofofd} hold for
 $U^{(n)}_\frak p, U^{(n)}_{\frak p\frak q}$.
\item
$U^{(1)}_\frak p, U^{(1)}_{\frak p\frak q}$ is
$U'_\frak p, U'_{\frak p\frak q}$, respectively.
\end{enumerate}
Let $U^{(n)}(X)$ be the space obtained from $U^{(n)}_\frak p, U^{(n)}_{\frak p\frak q}$
as in Definition \ref{clrelation}.

Let us consider the good coordinate system $(U_\frak p, E_\frak p, \psi_\frak p, s_\frak p)$
$\frak p\in \frak P$ of our Kuranishi structure. Let $\#\frak P = N$.
We put $\frak P = \{\frak p_1,\ldots,\frak p_N\}$, where $\frak p_i < \frak p_j$
only if $i < j$. We take
$n = 10N^2$ and consider $U^{(n)}_\frak p, U^{(n)}_{\frak p\frak q}$  as above.

\begin{prop}\label{existmulti}
For each $\epsilon>0$, there exists a system of multisections
$s_{\epsilon,\frak p}$ on
$U^{(n)}_\frak p$ for $\frak p \in \frak P$ with the following properties.
\begin{enumerate}
\item
$s_{\epsilon,\frak p}$ is transversal to $0$.
\item
$s_{\epsilon,\frak p}\circ{\underline\phi}_{\frak p\frak q}
= \hat{\underline\phi}_{\frak p\frak q}\circ s_{\epsilon,\frak q}$.
\item
The derivative of (arbitrary branch of) $s_{\epsilon,\frak p}$ induces an isomorphism
\begin{equation}\label{2tangent}
N_{U_{\frak p\frak q}}U_\frak p \cong
\frac{{\underline\phi}_{\frak p\frak q}^*E_\frak p}{E_\frak q\vert_{U_{\frak p\frak q}}}
\end{equation}
that coincides with the isomorphism  (\ref{dfiber}).
\item
The $C^0$ distance of $s_{\epsilon,\frak p}$ from $s_\frak p$ is smaller than
$\epsilon$.
\end{enumerate}
\end{prop}
The proof of Proposition \ref{existmulti} is by double induction.
Since the indices appearing in the proof is rather cumbersome,
we explain the construction in 2 simple cases before giving
the proof of  Proposition \ref{existmulti}.
\par\medskip
We first explain the case where we have two Kuranishi neighborhoods $U_1, U_2$
and the coordinate transformation $\underline{\phi}_{21}
: U_{21} \to U_2$, where $U_{21} \subset U_1$.
\par
The first step is to construct a multisection $s_{1,\epsilon}$ on $U_1$
that is a perturbation of the Kuranishi map $s_1$ and is transversal
to $0$.
Existence of such $s_{1,\epsilon}$ is  a consequence of \cite[Lemma 3.14]{FOn}.
\par
We next extend it to $U_2$ as follows.
We take a relatively compact subsets $U_i^{(1)} \subset U_i$
such that $\psi_i(s_i^{-1}(0) \cap U_i^{(1)})$ ($i=1,2$)
still covers our space $X$.
We put
$U_{21}^{(1)} = U_2^{(1)} \cap \underline{\phi}_{21}^{-1}(U_1^{(1)})$.
This set is relatively compact in $U_{21}$.
\par
We consider the normal bundle $N_{U_{21}^{(1)}}U_2$ and identify
its disk bundle with a tubular neighborhood ${\mathcal N}_{U_{21}^{(1)}}U_2$.
\par
Using the fiber derivative of $s_2$ we have an isomorphism
$$
E_2\vert_{{\mathcal N}_{U_{21}^{(1)}}U_2}
\cong \pi^*(E_1 \oplus E_1^\perp)
\cong \pi^* (E_1
\oplus N_{U_{21}^{(1)}} U_2).
$$
By the implicit function theorem, the $E_1^\perp$-component of the Kuranishi map $s_2$
induces a diffeomorphism from a sufficiently small tubular neighborhood of $U_{21}^{(1)}$ in
$U_2$ onto a neighborhood of the zero section of $E_1^\perp$.

We can extend $s_{1,\epsilon}$ to the tubular neighborhood
such that the first factor is the pull back of $s_{1,\epsilon}$
and the second factor is the $E_1^\perp$-component of the Kuranishi map $s_2$,
which is considered as the inclusion of
${\mathcal N}_{U_{21}^{(1)}}U_2$ to $N_{U_{21}^{(1)}}U_2$.
We denote  the extension we obtain by $s'_{\epsilon,2}$.
\par
Now we use the relative version of \cite[Lemma 3.14]{FOn}
to  extend $s'_{\epsilon,2}$ to $U_2$ further such that the extension coincides
with $s'_{\epsilon,2}$ on a slightly smaller tubular neighborhood.
\par\medskip
We next explain the case where
we have three Kuranishi neighborhoods $U_1, U_2, U_3$,
coordinate changes $\underline{\phi}_{21}, \underline{\phi}_{32}, \underline{\phi}_{31}$ and $U_{21}, U_{32}, U_{31}$.
\par
We take $U_i^{(j)}$ that is relatively compact in $U_i^{(j-1)}$
but $\bigcup_{i=1,2,3}U_i^{(j)}$  still covers $X$ in an obvious sense.
\par
In exactly the same way as the case when we have two Kuranishi neighborhoods,
we obtain $s_{1,\epsilon}^{(2)}$ and $s_{2,\epsilon}^{(2)}$
on $U_{1}^{(2)}$ and $U_{2}^{(2)}$ respectively,
that are compatible on the overlapped part.
\par
We next extend it to $U_3$.
\par
We extend  $s_{1,\epsilon}^{(2)}$ and $s_{2,\epsilon}^{(2)}$
to the tubular neighborhoods of
$U_{31}^{(3)}$ and $U_{32}^{(3)}$ in $U^{(2)}_3$.
We denote them by
$s_{3,\epsilon}^{(3),1}$ and $s_{3,\epsilon}^{(3),2}$.
\par
Let us study them on the part where two tubular neighborhoods intersect.
Let $x \in U_{31}^{(3)} \cap U_{32}^{(3)}$.
We have
$$
(N_{U_{31}^{(3)}}U_3^{(2)})_x \cong (N_{U_{21}^{(3)}}U_2^{(2)})
_{x} \oplus (N_{U_{32}^{(3)}}U_3^{(2)})_{\underline{\phi}_{32}(x)}
$$
This isomorphism is compatible with the isomorphism
$$
\left( \frac{E_3}{E_1}\right)_x
\cong
\left( \frac{E_2}{E_1}\right)_x
\oplus
\left( \frac{E_3}{E_2}\right)_{\underline{\phi}_{32}(x)}.
$$
Therefore we can glue $s_{3,\epsilon}^{(3),1}$ and $s_{3,\epsilon}^{(3),2}$
on the overlapped part by partition of unity.
\par
We then extend them to the whole $U_3^{(3)}$ by
the relative version of \cite[Lemma 3.14]{FOn}.
\par
The general case which is given below, is similar,
though the notation is rather cumbersome.
\begin{proof}
\footnote{The argument below is one written in \cite{Fu1}.
(If a shorter proof is preferable for readers please read page 3-4 of  \cite{Fu1}.)}
We will construct a system
$
(s_{\epsilon,\frak p}^k; \frak p \in \{\frak p_1, \dots, \frak p_k\})
$
where $s_{\epsilon,\frak p}^k$ is a multisection on $U^{(10k^2)}_{\frak p}$,
by upward induction on $k$, so that they satisfy the conditions (1)-(4) above.

When $k=1$, the proof is a standard perturbation argument combined with
the averaging process over the finite isotropy group. (See the proof of Theorem 3.1 \cite{FOn}
for the details.)

So we fix $k$ here and explain the process of constructing $s_{\epsilon,\frak p_i}^k$
for $i = 1, \ldots, k$, under the hypothesis that we have already constructed corresponding
sections $s_{\epsilon,\frak p_i}^{k-1}$ for $i = 1,\ldots, k-1$.
For $s_{\epsilon,\frak p_i}^k$ with $i<k$, the multisection
$s_{\epsilon,\frak p_i}^k$
is obtained by restriction of the domain of $s_{\epsilon,\frak p_i}^{k-1}$
and performing some adjustment around the boundary.
(See (\ref{skpidefform}).)

We identify the image $\underline\phi_{\frak p_k\frak p_i}(U_{\frak p_k\frak p_i}^{(m)})$ in
$U^{(m)}_{\frak p_k}$ with $U_{\frak p_k\frak p_i}^{(m)}$ and regard the latter as a subset
of $U^{(m)}_{\frak p_k}$ for any integer $m$. (Here $i<k$.)
We put
\begin{equation}
\mathcal N_k^i =
\bigcup_{j=i}^{k-1} \mathcal N^i_{U_{\frak p_k\frak p_j}^{(10 (k-1)^2+10(k-i))}}
U_{\frak p_k}^{(10 (k-1)^2+10(k-i))}
\end{equation}
and will define a multisection
$s^{k,i}_{\epsilon,\frak p_k}$ on $\mathcal N_k^i$ by downward induction on $i$.
\par
Here the open subset
$\mathcal N^i_{U_{\frak p_k\frak p_j}^{(10 (k-1)^2+10(k-i))}}U_{\frak p_k}^{(10 (k-1)^2+10(k-i))}$
is a tubular neighborhood of $U_{\frak p_k\frak p_j}^{(10 (k-1)^2+10(k-i))}$, which will be chosen
in the inductive process.
\par
We assume
the closure of
$$\mathcal N^{i}_{U_{\frak p_k\frak p_j}^{(m+1)}}U_{\frak p_k}^{(m+1)}$$
is compact in
$$\mathcal N^{i+1}_{U_{\frak p_k\frak p_j}^{(m)}}U_{\frak p_k}^{(m)}.$$
\par

Let us start the downward induction for $i$ starting from $i=k-1$.
We have an embedding $U_{\frak p_k\frak p_{k-1}}^{(10 (k-1)^2)} \to U_{\frak p_k}^{(10 (k-1)^2)}$.
We take its tubular neighborhood
$\mathcal N^{k-1}_{U_{\frak p_k\frak p_{k-1}}^{(10 (k-1)^2)}}U_{\frak p_k}^{(10 (k-1)^2)}$.
We also have $s_{\epsilon,\frak p_{k-1}}^{k-1}$ on $U_{\frak p_k\frak p_{k-1}}^{(10 (k-1)^2+10(k-i-1))}$
by the induction hypothesis.
We have already taken a splitting
\begin{equation}\label{eq:Epk}
E_{\frak p_k} = E_{\frak p_{k-1}} \oplus E^{\perp}_{\frak p_{k-1}}
\end{equation}
on $U_{\frak p_k\frak p_{k-1}}^{(10 (k-1)^2+10(k-i-1))}$.
(Here we identify $E_{\frak p_{k-1}}$ with its image by $\hat{\phi}_{\frak p_k\frak p_{k-1}}$.)
We also extended the bundles $E_{\frak p_{k-1}}$, $E^{\perp}_{\frak p_{k-1}}$ to
a small tubular neighborhood $\mathcal N^{k-1}_{U_{\frak p_k\frak p_{k-1}}^{(10 (k-1)^2)}}U_{\frak p_k}^{(10 (k-1)^2)}$.\par
We extend $s_{\epsilon,\frak p_{k-1}}^{k-1}$ to
$\mathcal N^{k-1}_{U_{\frak p_k\frak p_{k-1}}^{(10 (k-1)^2)}}U_{\frak p_k}^{(10 (k-1)^2)}$
so that
\begin{enumerate}
\item[(a)] it is $\epsilon/(2^n)$ close to the Kuranishi map $s_{\frak p_k}$,
\item[(b)] it coincides with the perturbation already defined at the zero section of the normal bundle for
the first component $E_{\frak p_{k-1}}$ of \eqref{eq:Epk}, and
\item[(c)] the second component thereof in the splitting \eqref{eq:Epk} is
the same as the given Kuranishi map $s_{\frak p_k}$.
\end{enumerate}
Transversality condition is obviously satisfied and the properties (2)-(4) also hold
true by construction.
We thus obtained $s^{k,k-1}_{\epsilon,\frak p_k}$
\par
Now we go to the inductive step to construct  $s^{k,i}_{\epsilon,\frak p_k}$ assuming we have
$s^{k,i+1}_{\epsilon,\frak p_k}$.
\par
We consider the embedding
$$
U_{\frak p_k\frak p_j}^{(10 (k-1)^2+10(k-i)-10)} \to U_{\frak p_k}^{(10 (k-1)^2+10(k-i)-10)},
$$
for $i+1 \leq j \leq k$.
We identify $U_{\frak p_k\frak p_j}^{(10 (k-1)^2+10(k-i))-10)}$ with the image of embedding
and take its tubular neighborhood
\begin{equation}\label{311}
\mathcal N^{i+1}_{U_{\frak p_k\frak p_j}^{(10 (k-1)^2+10(k-i)-10)}}U_{\frak p_k}^{(10 (k-1)^2+10(k-i)-10)}.
\end{equation}
Note that we already have our section $s^{k,i+1}_{\epsilon,\frak p_k}$
on ${\mathcal N}_k^{i+1}$ which is the union of (\ref{311}), $j=i+1, \dots, k-1$.
\par
We next apply the same argument as the first step
to obtain a $s_{\epsilon,\frak p_k}^{\prime k,i}$ on a tubular neighborhood
$$
\mathcal N^{+,i}_k=\mathcal N^{+,i}_{U_{\frak p_k\frak p_i}^{(10 (k-1)^2+10(k-i)-9)}}U_{\frak p_k}^{(10 (k-1)^2+10(k-i)-9)}.
$$
(Here we use the fact that we have a multisection $s_{\epsilon,\frak p_i}^{k-1}$
by induction hypothesis.)
\par
Note that ${\mathcal N}_k^{i+1}$ and ${\mathcal N}^{+,i}_k$ are open subsets in an orbifold
$U_{\frak p_k}$.
We take a smooth function
$$
\chi : {\mathcal N}_k^{i+1} \cup {\mathcal N}_k^{+,i} \to [0,1]
$$
such that $\{\chi, 1-\chi\}$ is a partition of unity subordinate to $\{ \mathcal N_k^{i+1}, \mathcal N_k^{+,i} \}$
and
$\chi = 1$ on
\begin{equation}\label{315}
\bigcup_{j=i+1}^{k-1} \mathcal N^i_{U_{\frak p_k\frak p_j}^{(10 (k-1)^2+10(k-i)-8)}}
U_{\frak p_k}^{(10 (k-1)^2+10(k-i)-8)} \subset \mathcal N_k^{i+1}.
\end{equation}
We put:
\begin{equation}\label{312}
s^{k,i}_{\epsilon,\frak p_k} = \chi s^{k,i+1}_{\epsilon,\frak p_k} + (1-\chi) s_{\epsilon,\frak p_k}^{\prime k,i}.
\end{equation}
\begin{rem}
The sum of multisections is a bit delicate to define.
In our case $s^{k,i+1}_{\epsilon,\frak p_k}$ and $s_{\epsilon,\frak p_k}^{\prime k,i}$
are defined by extending the multisections to the tubular neighborhood.
This process does not change the number of branches. (Namely they are extended branch-wise.)
So though these two do not coincide as sections, each branch of $s^{k,i+1}_{\epsilon,\frak p_k}$
has a corresponding branch of $s_{\epsilon,\frak p_k}^{\prime k,i}$.
So we can apply the formula (\ref{312}) branch-wise.
\end{rem}
We find that $s^{k,i}_{\epsilon, \frak p_k}$ is transversal to zero on
$$\bigcup_{j=i+1}^{k-1} \mathcal N^i_{U_{\frak p_k\frak p_j}^{(10 (k-1)^2+10(k-i)-8)}}
U_{\frak p_k}^{(10 (k-1)^2+10(k-i)-8)} \subset \mathcal N_k^{i+1},$$
since $\chi =1$ there.
Note that $s^{k,i+1}_{\epsilon, \frak p_k}$ and $s^{\prime k,i}_{\epsilon, \frak p_k}$
coincide on $U_{\frak p_k \frak p_i}^{10(k-1)^2 +10(k-i)-9)} \cap \mathcal N_k^{i+1}$
up to first derivatives (including the normal direction to $U_{\frak p_k \frak p_i}^{10(k-1)^2 +10(k-i)-9)}$
in $U_{\frak p_k}^{10(k-1)^2 +10(k-i)-9)})$.
Hence $s_{\epsilon, \frak p_k}^{k,i}$ is transversal to zero
on
$$\mathcal N^i_{U_{\frak p_k \frak p_i}^{(10(k-1)^2 + 10(k-i) -7}} U_{\frak p_k}^{(10(k-1)^2+10(k-i)-7)},$$
if we take this tubular neighborhood of $U_{\frak p_k \frak p_i}^{(10(k-1)^2 + 10(k-i) -7}$
sufficiently small.

We restrict (\ref{312}) to
$$
\bigcup_{j=i}^{k-1}\mathcal N^i_{U_{\frak p_k\frak p_j}^{(10 (k-1)^2+10(k-i)-7)}}
U_{\frak p_k}^{(10 (k-1)^2+10(k-i)-7)}.
$$
Then
the multi-section  (\ref{312}) is transversal to the zero section.
It also satisfies the properties (2)(3)(4).
\begin{rem}
We use Lemma \ref{intersection} here for the consistency of tubular neighborhoods.
We may use Mather's compatible system of tubular neighborhoods \cite{Math73}.
However the present situation is much simpler because of Lemma \ref{intersection}.
So the compatibility of tubular neighborhoods is obvious.
\end{rem}
The section $s^{k,i}_{\epsilon,\frak p_k}$ coincides with $s^{k,i+1}_{\epsilon,\frak p_k}$  on the overlapped part because
we take $\chi = 1$ on (\ref{315}).
\par\medskip
We continue the induction up to $i=1$. Then we have a required multisection $s^{k,1}_{\epsilon,\frak p_k}$
on
$$
\mathcal N_k^1 =
\bigcup_{j=1}^{k-1} \mathcal N^i_{U_{\frak p_k\frak p_j}^{(10 k^2- 10)}}U_{\frak p_k}^{(10 k^2- 10)}.
$$
($10 (k-1)^2+10(k-1) < 10 k^2- 10$.)
By \cite[Lemma 3.14]{FOn}  we can extend this to $U_{\frak p_k}^{(10 k^2- 9)}$
so that it satisfies (1) and coincides with  $s^{k,1}_{\epsilon,\frak p_k}$ on
$$
\bigcup_{j=1}^{k-1} \mathcal N^i_{U_{\frak p_k\frak p_j}^{(10 k^2- 9)}}U_{\frak p_k}^{(10 k^2- 9)}.
$$
We thus obtain $s^{k}_{\epsilon,\frak p_k}$.
\par
For $i<k$ we define $s^{k}_{\epsilon,\frak p_i}$
as follows
\begin{equation}\label{skpidefform}
s^{k}_{\epsilon,\frak p_i}
=
\begin{cases}
\displaystyle
(\widehat{\underline{\phi}}_{\frak p_k\frak p_i})^{-1} \circ s^{k}_{\epsilon,\frak p_k} \circ \underline{\phi}_{\frak p_k\frak p_i}
&\text{when the right hand side is defined}
\\
s^{k-1}_{\epsilon,\frak p_i}
&\text{otherwise.}
\end{cases}
\end{equation}
\par
The properties (1)-(4) are satisfied.
The proof of Proposition \ref{existmulti} is complete.
\end{proof}

We have thus constructed our perturbation that is a multisection $s_{\epsilon,\frak p}$.
To obtain a virtual fundamental chain and prove its basic
properties we need to restrict the section $s_{\epsilon,\frak p}$ to an appropriate neighborhood
of the union of zero sets of $s_{\frak p}$ and study the properties of $s_{\epsilon,\frak p}$ there.
Namely we will prove the following:
\begin{enumerate}
\item
The zero set $s_{\epsilon,\frak p}^{-1}(0)$ converges to the zero set $s_{\frak p}^{-1}(0)$.
\item
The union of the zero sets $s_{\epsilon,\frak p}^{-1}(0)$ is compact.
\item
The union of the zero sets $s_{\epsilon,\frak p}^{-1}(0)$ carries a fundamental cycle.
\end{enumerate}
We need to shrink the domain to prove this.
\par
Note we have $\lim_{\epsilon\to 0}s_{\epsilon,\frak p} = s_{\frak p}$.
So if the domains of the mutisections are compact,
this implies that the union of zeros of $s_{\epsilon,\frak p}$
converges to the union of zeros of $s_{\frak p}$, which is
nothing but our space $X$.
\par
Since the domains $U_{\frak p}^{(k)}$ of $s_{\epsilon,\frak p}$
is noncompact, the actual proof is slightly more
nontrivial.
Note we took a sequence $U_{\frak p}^{(k)}$ such that
$U_{\frak p}^{(k+1)}$ is relatively compact
in $U_{\frak p}^{(k)}$.
We use this fact in place of the compactness of the domain.
Then, the point to take care of is the fact that
$s_{\frak p}^{-1}(0)$ may intersect with the boundary
$\partial U_{\frak p}^{(k+1)}$.
We use Lemma \ref{unionofofd} for this purpose.
(See Remark \ref{rem68888}.)
\par
We already shrank our good coordinate system several times.
We shrink again below. Let us denote by $U_{\frak p}$, $U(X)$ etc.  the good coordinate system
and a space obtained  when we proved Proposition \ref{existmulti}.
(In other words  $U_{\frak p} =U^{(N)}_{\frak p}$ in the notation we used
in the proof of Proposition \ref{existmulti}.)
During the discussion of the rest of this section, we restart numbering the
shrunken good coordinate system and will write again $U^{(m)}_{\frak p}$.
\par
Lemma \ref{compactmainlemma} below is the key lemma.
Note all the spaces $U^{(m)}(X)$ with the weak topology induced from $K(\frak P)$
by the injective map $\frak J_{K(\frak P)U^{(m)}(X)}$,
can be regarded as subsets of $K(\frak P) = \bigcup \overline U'_{\frak p} / \sim$
(where $K(\frak P)$ is as in Definition \ref{Imapdef}).
Recall that the space $K(\frak P)$ is compact and metrizable. We take and fix a metric
on it. We consider the metrics on $U^{(m)}(X)$ induced from the metric on $K(\frak P)$.
Then by definition the maps $\frak I_{U^{(m')}(X)U^{(m)}(X)}$ are isometries for $m' < m$.
\par
For a given metric space $Z$ and a subset $C \subset Z$ we denote
its $\epsilon$ neighborhood by
\begin{equation}
B_{\epsilon}(C;Z) = \{z \in Z \mid d(z,C) < \epsilon\}.
\end{equation}
\begin{lem}\label{compactmainlemma}
We may choose $U^{(2)}_{\frak p}$, $U^{(3)}_{\frak p}$  and $\delta > 0$, so that
\begin{equation}\label{320formula}
\aligned
 &B_{\delta}(I(X),U^{(2)}(X)) \cap \bigcup_{\frak p_i \in \frak P} \Pi_{\frak p}
 (\overline s_{\epsilon,\frak p}^{-1}(0) \cap U^{(2)}(X))
\\
&=
\frak I_{U^{(2)}(X)U^{(3)}(X)}
\left(B_{\delta}(I(X),U^{(3)}(X)) \cap \bigcup_{\frak p_i \in \frak P} \Pi_{\frak p}
(\overline s_{\epsilon,\frak p}^{-1}(0) \cap U^{(3)}(X))
 \right).
\endaligned
\end{equation}
\end{lem}
\begin{rem}\label{rem67}
Note if we take $U^{(3)}_{\frak p} \subset U^{(2)}_{\frak p}$ for each $\frak p$ then
$U^{(3)}(X)$ is defined by
$$
U^{(3)}(X) = \bigcup_{\frak p} \Pi_{\frak p}(U^{(3)}_{\frak p}) \subset U^{(2)}(X),
$$
where $\Pi_{\frak p} : U^{(2)}_{\frak p} \to U^{(2)}(X)$.
\end{rem}
\begin{rem}\label{rem68888}
We remark that Lemma \ref{compactmainlemma} claims that
the union of zero sets of $s_{\epsilon,\frak p}$
on $U^{(2)}(X)$ is equal
to the union of zero sets of $s_{\epsilon,\frak p}$
on $U^{(3)}(X)$,
if we restrict them to a small neighborhood of $I(X)$.
\end{rem}
\begin{exm}
We consider $U'_i$ in Example \ref{omeganeighex}.
We define an obstruction bundle $E_i$ on it such that
$E_1$ is trivial bundle, $E_2 = U'_2 \times \R$, and $E_3 = U'_3 \times \R^2$.
We define  sections (= Kuranishi maps) $s_i$ by $s_2(x,y,0) = y$, $s_3(x,y,z) = (y,z)$.
It defines Kuranishi structure and a good coordinate system in an obvious way.
The sections $s_i$ are already transversal to $0$ so we do not perturb.
\par
We put $U^{(2)}(X) = U'(X)$, $U_i^{(2)}(X) = U'_i(X)$.
We take
$$
\aligned
&U^{(3)}_1 = \{(x,0,0) \mid x < 2c\}, \\
&U^{(3)}_2 = \{ (x,y,0) \mid x > -y^2 + c\}, \\
& U^{(3)}_3 = \{ (x,y,z) \mid x >c\}.
\endaligned
$$
Here $c>0$.
In this example the zero set of the section $s_i$ lies on the $x$-axis.
So the lemma holds. The key point of the proof of Lemma \ref{compactmainlemma} is in a neighborhood of the border of
$U^{(2)}_{2}$ and $U^{(2)}_3$ the condition $s_i(x) =0$ implies that the
point $x$ lies in $U^{(3)}_1$. This is the consequence of
Proposition \ref{existmulti} (3).
\end{exm}
\begin{proof}[Proof of Lemma \ref{compactmainlemma}]
In Lemma  \ref{unionofofd} we defined
the set $\Omega_{\frak q_i}(x) \subset U''_{\frak q_i}$ for
$x \in X$. We take $U^{(k)}_{\frak p}$ in place of
$U''_{\frak p}$ and obtain
$\Omega_{\frak q_i}^{(k)}(x) \subset U^{(k)}_{\frak q_i}$.
\par
We then define $\frak U^{(k)}(x)$ by (\ref{frakUdef}).
We choose a neighborhood $O_x$ of $x$ in $U^{(2)}(X)$ and $\delta_1$ so that
\begin{equation}\label{Oxsmallerthaq1}
O_x \cap \Pi_{\frak q_1}(\Omega^{(2)}_{\frak q_1}(x))=
O_x \cap  \overline{\Pi_{\frak q_1}(\Omega^{(2)}_{\frak q_1}(x))}
\end{equation}
and
$$
O_x \cap \frak U^{(2)}(x) \supset B_{\delta_1}(I(x);U^{(2)}(X)).
$$
(Note that $x$ is in the interior of $\Omega^{(2)}_{\frak q_1}(x)$. So
we can take $O_x$ small so that (\ref{Oxsmallerthaq1}) is satisfied.)
\par
\begin{sublem}\label{sublem35}
We may take $U^{(3),x}_{\frak p}$ (depending on $x$) and $\delta_3>0$ (independent of $x$)
so that
\begin{equation}\label{sublem69maineq}
\aligned
&O_{x} \cap \frak U^{(2)}(x) \cap B_{\delta_3}(I(x),U^{(2)}(X))
\cap \bigcup_{\frak p \in \frak P} \Pi_{\frak p} (
\overline s_{\epsilon,\frak p}^{-1}(0))  \\
&
\subseteq\frak I_{U^{(2)}(X)U^{(3),x}(X)}
\left(B_{\delta_3}(I(X),U^{(3),x}(X)))\cap \bigcup_{\frak p \in \frak P}
\Pi_{\frak p} (\overline s_{\epsilon,\frak p}^{-1}(0) \cap U^{(3),x}(X)
) \right)
\endaligned
\end{equation}
for each $x \in X$ holds for sufficiently small $\epsilon$.
\end{sublem}
We note that $\{U^{(3),x}_{\frak p} \mid \frak p \in \frak P\}$ determines $U^{(3),x}(X)$ as in Remark \ref{rem67}.
\begin{proof}
We have
$$
U^{(2)}(X) \cap O_x = \bigcup_{i=1,\dots,m}\Pi_{\frak q_i}(\Omega^{(2)}_{\frak q_i}(x))\cap O_x.
$$
By (\ref{Oxsmallerthaq1}) we have
$$
O_x \cap \Pi_{\frak q_1}(\Omega^{(2)}_{\frak q_1}(x))=
O_x \cap  \overline{\Pi_{\frak q_1}(\Omega^{(2)}_{\frak q_1}(x))}.
$$
Therefore we may choose $U^{(3),x}_{\frak p}$ close enough to  $U^{(2)}_{\frak p}$ so that
$$
O_{x} \cap \Pi_{\frak q_1}(\Omega^{(2)}_{\frak q_1}({x}))
\subset
\frak I_{U^{(2)}(X)U^{(3),x}(X)}(\Pi_{\frak q_1}(U^{(3),x}_{\frak q_1})).
$$
\par
On the other hand, in a sufficiently small tubular neighborhood $\mathcal N_{U_{\frak q_1\frak q_i}}U_{\frak q_i}$ ($i>2$)
the zero set $s_{\epsilon,q_i}^{-1}(0)$ is contained in the subset
$U_{\frak q_i\frak q_1} \subset \mathcal N_{U_{\frak q_i\frak q_1}}U_{\frak q_i}$.
(Here $U_{\frak q_i\frak q_1}$ is identified with the zero section of the normal bundle
$\mathcal N_{U_{\frak q_i\frak q_1}}U_{\frak q_i}$.)
This is a consequence of Proposition \ref{existmulti} (3).
\par
We can choose $\delta_3>0$ sufficiently small so that
$O_{x} \cap \frak U^{(2)}(x) \cap B_{\delta_3}(x,U^{(2)}(X))$
is contained in this tubular neighborhood.
(We may choose $\delta_3$ independent of $x$ since $X$ is compact.)
The sublemma follows.
\end{proof}
We find a finite number of $x_i \in I(X)$, $i=1,\dots,\frak I$ and $\delta > 0$, such that
\begin{equation}\label{coverOXI}
\bigcup_{i=1}^{\frak I} O_{x_i} \cap \frak U^{(2)}(x_i) \cap B_{\delta_3}(I(x_i),U^{(2)}(X))
\supset B_{\delta}(I(X);U^{(2)}(X)).
\end{equation}
We choose $U^{(3)}_{\frak p}$ such that
$$
U^{(2)}_{\frak p} \supset \overline U^{(3)}_{\frak p} \supset U^{(3)}_{\frak p} \supset \bigcup_{i=1}^{\frak I} U^{(3),x_i}_{\frak p}.
$$
Then
(\ref{sublem69maineq}) holds for all $x_i$ ($i=1,\dots,\frak I$) when we replace
$U^{(3),x_i}_{\frak p}$ by this $U^{(3)}_{\frak p}$.
By (\ref{coverOXI}) we have:
$$
\aligned
&\bigcup_{i=1}^{\frak I}
\left(
O_{x_i} \cap \frak U^{(2)}(x_i) \cap B_{\delta_3}(x_i,U^{(2)}(X))
\cap \bigcup_{\frak p \in \frak P} \Pi_{\frak p} (
\overline s_{\epsilon,\frak p}^{-1}(0))
\right)\\
&\supset
B_{\delta}(I(X),U^{(2)}(X)) \cap \bigcup_{\frak p \in \frak P} \Pi_{\frak p}
(\overline s_{\epsilon,\frak p}^{-1}(0) \cap U^{(2)}(X)) .
\endaligned$$
Therefore (\ref{sublem69maineq}) implies that the left hand side of (\ref{320formula}) is contained
in the right hand side.
The inclusion of the opposite direction is obvious.
\end{proof}
\begin{lem}\label{cuttedmodulilem}
\begin{equation}\label{cuttedmoduli}
B_{\delta}(I(X),U^{(3)}(X)) \cap \bigcup_{\frak p \in \frak P} \Pi_{\frak p}
(\overline s_{\epsilon,\frak p}^{-1}(0) \cap U^{(3)}(X))
\end{equation}
is compact if $\epsilon > 0$ is sufficiently small.
\end{lem}
\begin{proof}
Lemma \ref{compactmainlemma} implies
\begin{equation}\label{closurehuhen}
\aligned
&B_{\delta}(I(X),U^{(3)}(X)) \cap \bigcup_{\frak p \in \frak P} \Pi_{\frak p}
(\overline s_{\epsilon,\frak p}^{-1}(0) \cap U^{(3)}(X))\\
&=
B_{\delta}(I(X),\overline U^{(3)}(X)) \cap \bigcup_{\frak p \in \frak P}
\Pi_{\frak p} (\overline s_{\epsilon,\frak p}^{-1}(0) \cap \overline U^{(3)}(X)).
\endaligned\end{equation}
(Here  $\overline U^{(3)}(X)$ is the closure of $U^{(3)}(X)$ in $U^{(2)}(X)$.)
We remark $U^{(3)}(X)$ is a relatively compact subspace of $U^{(2)}(X)$.
The right hand side of (\ref{closurehuhen}) is clearly compact.
\end{proof}
Hereafter we fix $\delta$ and write (\ref{cuttedmoduli})
as $s_{\epsilon}^{-1}(0)_{\delta}$. It is a subspace of $U^{(3)}(X)$.
\begin{lem}
$$
\lim_{\epsilon\to 0} \frak I_{U^{(2)}(X)U^{(3)}(X)}(s_{\epsilon}^{-1}(0)_{\delta})
\subseteq
B_{\delta}(I(X),U^{(2)}(X)) \cap \bigcup_{\frak p \in \frak P}
\Pi_{\frak p} (\overline s_{\frak p}^{-1}(0) \cap U^{(2)}(X)).
$$
Here the convergence is by Hausdorff distance.
\end{lem}
\begin{proof}
This is a consequence of the next sublemma applied to the right hand side of
(\ref{closurehuhen}) chartwise and branchwise.
\end{proof}
\begin{sublem}
Let $E \to Z$ be a vector bundle on a compact metric space $Z$ and $s$ its   section.
Suppose that $s_{\epsilon}$ is a family of sections converges to $s$. Then
$$
\lim_{\epsilon \to 0} s_{\epsilon}^{-1}(0) \subseteq s^{-1}(0).
$$
\end{sublem}
\begin{proof}
Let $\rho >0$. We put
$$
\epsilon(\rho) = \inf \{ \vert s(x) \vert \mid x \in Z \setminus B_{\rho}(s^{-1}(0);Z)\}.
$$
Clearly
$$
 s_{\epsilon}^{-1}(0)
 \subseteq
B_{\rho}(s^{-1}(0);Z)
$$
if $\epsilon < \epsilon(\rho)$. The compactness of $Z$ implies $\epsilon(\rho)>0$.
\end{proof}
\begin{lem}
$s^{-1}_{\epsilon}(0)_{\delta}$ has a triangulation.
\end{lem}
This is proved in \cite[Lemma 6.9]{FOn}.
\par
Suppose we have a strongly continuous map $f = \{f_\frak p \mid \frak p \in \frak  P\}$
to a topological space $Z$ and our Kuranishi space $X$ is oriented.
(\cite[A1.17]{fooo:book1}.)
Then we can put a weight to each simplex of top dimension in
$s^{-1}_{\epsilon}(0)_{\delta}$
and $f_*[X]$. (\cite[(6.10)]{FOn}.)
That is a singular chain of $Z$ (with rational coefficients.)
Thus we have constructed a virtual fundamental chain of $X$
\par
We consider the case our Kuranishi structure has no boundary.
\begin{lem}\label{cycleproperties}
We can put the weight on each of the top dimensional simplices
so that $f_*([X])$ is a cycle.
\end{lem}
\begin{proof}
It suffices to prove that for each $x \in s^{-1}_{\epsilon}(0)_{\delta}$ the zero set of each branch
of $s_{\epsilon,\frak p}$ is a smooth manifold in a neighborhood of $x$.
This is again a consequence of the proof of Lemma \ref{compactmainlemma} as follows.
Suppose $x \in O_{x_i} \cap \frak U^{(2)}(x_i)$.
Then we proved that $s^{-1}_{\epsilon}(0)_{\delta}$ in a neighborhood of $x$
coincides with the zero set of $s_{\epsilon,\frak q_1}$.
On the other hand each branch of $s_{\epsilon,\frak q_1}$ is transversal to $0$ by
Proposition \ref{existmulti}.
\end{proof}
\begin{rem}
We remark that we use the metric on $U(X)$ to {\it prove} that $s_{\epsilon}^{-1}(0)_{\delta}$ has
required properties for sufficiently small $\delta$ and $\epsilon$.
In particular we did {\it not} use the metric to {\it define} $s_{\epsilon}$.
Therefore the virtual fundamental cycle (chain) we obtained is obviously independent
of the choice of the metric on $U(X)$.
\end{rem}
In sum we have defined a virtual fundamental cycle $f_*([X])$
using the good coordinate system on the Kuranishi structure of $X$ which has
tangent bundle, oriented and without boundary.
Using the relative version of our construction of this and the next sections,
we can show the cohomology class of $f_*([X])$ is independent of the
choices. (See \cite[Lemmas 17.8,17.9]{FOn2} etc..)

\section{Existence of good coordinate system}\label{sec:existenceofGCS}

The purpose of this section is to prove the next theorem.
\begin{thm}\label{goodcoordinateexists}
Let $X$ be a compact metrizable space with Kuranishi structure in the
sense of Definition \ref{Definition A1.5}.
Then $X$ has a good coordinate system
in the sense of Definition \ref{goodcoordinatesystem}.
They are compatible in the following sense.
\end{thm}
\begin{defn}
Let $X$ be a space with Kuranishi structure.
A good coordinate system is said to be {\it compatible} with this
Kuranishi structure if the following holds.
\par
Let $(U'_\frak p,E'_\frak p,s'_\frak p,\psi'_\frak p)$ be a chart of the given good coordinate
system of $X$. For each $q\in \psi'_\frak p(\tilde{q})$, $\tilde{q} \in U'_\frak p \cap
(s'_\frak p)^{-1}(0)$ there exists $(\hat{\underline\phi}_{\frak p q},{\underline\phi}_{\frak p q})$ such that
\begin{enumerate}
\item[{(1)}]  There exist
an open subset $U_{\frak p q}$ of $U_q$ and an
embedding of vector bundles
$$
\hat{\underline\phi}_{\frak p q}:E_{q}\vert_{U_{\frak p q}} \to E'_{\frak p}
$$
over an embedding ${\underline\phi}_{\frak  p q}: U_{\frak p q} \to U'_{\frak p}$
of orbifold of $U_{\frak p q}$ that satisfy
\begin{enumerate}
\item
${\underline\phi}_{\frak p q}([o_q]) = q$.
\item $\hat{\underline\phi}_{\frak pq}\circ s_q =
s'_\frak p\circ{\underline\phi}_{\frak pq}$.
\item
$\psi_q =
\psi'_\frak p\circ \underline{\phi}_{\frak pq}$ on $s_q^{-1}(0) \cap U_{\frak pq}$.
\end{enumerate}
\item[{(2)}]
$d_{\rm fiber}s'_{\frak p}$ induces an isomorphism of vector bundles
at $s_{q}^{-1}(0)\cap U_{\frak pq}$.
\begin{equation}\label{15tangent22}\nonumber
N_{U_{\frak p q}}U'_\frak p \cong \frac{\underline\phi_{\frak p q}^*E_\frak p}{E_q\vert_{U_{\frak p q}}}.
\end{equation}
\item[{(3)}]
If $r \in \psi_q(U_{\frak pq}\cap s_q^{-1}(0))$ and
$q \in \psi'_\frak p((s'_\frak p)^{-1}(0)\cap U'_{\frak p})$, then
$$
\underline{\phi}_{\frak pq} \circ \underline{\phi}_{qr} = \underline{\phi}_{\frak pr}, \quad
\hat{\underline{\phi}}_{\frak pq} \circ \hat{\underline{\phi}}_{qr} =
\hat{\underline{\phi}}_{\frak pr}.
$$
Here the first equality holds on
$$
\underline\phi_{qr}^{-1}(U_{\frak p q}) \cap U_{qr} \cap U_{\frak p r} = \underline\phi_{qr}^{-1}(U_{\frak p q}) \cap U_{qr} \cap U_{\frak p r}
$$
and the second equality holds on $E_r\vert_{\underline\phi_{qr}^{-1}(U_{\frak p q}) \cap U_{qr} \cap U_{\frak p r}}$.
\item[{(4)}]
Suppose that $\frak o \ge \frak p$ and the coordinate change of good coordinate system is
given by $(U'_{\frak o\frak p},\hat{\underline{\phi}'}_{\frak o\frak p},{\underline{\phi}'}_{\frak o\frak p})$.
Let $q \in \psi'_{\frak p}({s'_\frak p}^{-1}(0)\cap U'_{\frak o\frak p} )$.
Then we have
$$
{\underline{\phi}'}_{\frak o\frak p} \circ {\underline{\phi}}_{\frak pq} =
{\underline{\phi}}_{\frak oq}, \quad
\hat{\underline{\phi}'}_{\frak o\frak p} \circ \hat{\underline{\phi}}_{\frak pq} =
\hat{\underline{\phi}}_{\frak oq}.
$$
Here the first equality holds on ${\underline{\phi}}_{\frak pq}^{-1}(U'_{\frak o\frak p}) \cap U_{\frak pq}
\cap U_{\frak oq}$, and the
second equality holds on $E_q\vert_{{\underline{\phi}}_{\frak pq}^{-1}(U'_{\frak o\frak p}) \cap U_{\frak pq}
\cap U_{\frak oq}}$.
\end{enumerate}
\end{defn}

In the above definition, we use $\{(U_p, E_p, s_p, \psi_p)\}$ for Kuranishi neighborhoods in the definition
of the Kuranishi structure and $\{(U'_{\frak p}, E'_{\frak p},s'_{\frak p}, \psi'_{\frak p})\}$ for the
data in the definition of the good coordinate system in order to distinguish them in the definition of compatibility
between them.
However, we will
write $\{(U_{\frak p}, E_{\frak p}, s_{\frak p}, \psi_{\frak p})\}$ instead of
$\{(U'_{\frak p}, E'_{\frak p},s'_{\frak p}, \psi'_{\frak p})\}$,
in the rest of this section.
\begin{proof}[Proof of Theorem \ref{goodcoordinateexists}]
Any point $p\in X$ carries a well defined dimension, $\dim U_p$, of orbifold.
We put $d_p = \dim U_p$ and
$$
X(\frak d) = \{p \in X \mid d_p = \frak d\}.
$$
\par
The first part of the proof is to construct an orbifold
(plus an obstruction bundle etc.) that is a `neighborhood'
of a compact subset of $X(\frak d)$.
Let us define such a notion precisely.
\begin{defn}\label{pureofdnbd}
Let $K_*$ be a compact subset of $X(\frak d)$.
A {\it pure orbifold neighborhood} of $K_*$ is
$(U_*,E_*,s_*,\psi_*)$ such that the following holds.
\begin{enumerate}
\item
$U_*$ is a $\frak d$-dimensional orbifold.
\item $E_*$ is a vector bundle whose rank is $\frak d - \dim X$.
(Here $\dim X$ is a dimension of $X$ as a Kuranishi space.)
\item
$s_*$ is a section of $E_*$.
\item
$\psi_* : s_*^{-1}(0) \to X$ is a homeomorphism to a neighborhood $\mathcal U_*$ of $K_*$
in $X(\frak d)$.
\end{enumerate}
We also assume the following compatibility condition with Kuranishi structure of $X$.
For any $p \in \psi_*(s_*^{-1}(0)) \subset X$, there exists $(U_{*p},
\hat{\underline{\phi}}_{*p},\underline{\phi}_{*p})$ such that
\begin{enumerate}
\item[(5)] $U_{*p}$ is an open neighborhood of $[o_p]$ in $U_p$.
\item[(6)] $\hat{\underline{\phi}}_{*p} : E_p\vert_{U_{*p}} \to E_*$ is an embedding of vector
bundle over an embedding of orbifold ${\underline{\phi}}_{*p} : U_{*p} \to U_*$ such that
\begin{enumerate}
\item $\hat{\underline{\phi}}_{*p}\circ s_p =
s_*\circ{\underline{\phi}}_{*p}$ on $U_{*p}$
\item
$\psi_p =
\psi_*\circ \underline{\phi}_{*p}$ on $s_p^{-1}(0) \cap U_{*p}$.
\end{enumerate}
\item[(7)]
The restriction of $ds_*$ to the normal direction induces an isomorphism
\begin{equation}\label{tangent**}
N_{U_{*p}}U_* \cong \frac{\hat{\underline{\phi}}_{*p}^*E_*}{E_p\vert_{U_{*p}}}
\end{equation}
as vector bundles on the orbifold $U_{*p}$ at $s_p^{-1}(0)$.
\item[(8)]
If $q  \in \psi_p(s_p^{-1}(0)\cap U_{*p})$, then
$$
\underline\phi_{*p} \circ \underline\phi_{pq} =
\underline\phi_{*q}, \quad
\underline{\hat\phi}_{*p} \circ \underline{\hat\phi}_{pq} =
\underline{\hat\phi}_{*q}.
$$
Here the first equality holds on $\underline{\phi}_{pq}^{-1}(U_{*p}) \cap U_{pq} \cap U_{*q}$,
and the second equality holds on $E_q\vert_{\underline{\phi}_{pq}^{-1}(U_{*p}) \cap U_{pq} \cap U_{*q}}$.
\end{enumerate}
\end{defn}
Hereafter we denote
\begin{equation}\label{calUnasi}
\mathcal U_* = \psi(s_*^{-1}(0)).
\end{equation}
The goal of the first part of the proof of Theorem \ref{goodcoordinateexists} is to prove the following.
\begin{prop}\label{purecover}
For any compact subset $K$ of $X(\frak d)$ there exists a pure orbifold neighborhood of $K$.
\end{prop}
\begin{proof}
We cover $K$ by a finite number of $\mathcal U_{p_j}$'s where $p_j \in K$ and
$$
\psi_{p_j}(s_{p_j}^{-1}(0)) = \mathcal U_{p_j}.
$$
There exist compact subsets $K_j$ of $\mathcal U_{p_j}$ such that
the union of $K_j$ contains $K$.
Thus to prove Proposition \ref{purecover}, by the induction argument
we have only to prove the following lemma.
\end{proof}
\begin{lem}\label{2gluelemma}
Let $K_1, K_2$ be compact subsets of $X(\frak d)$. Suppose $K_1$ and $K_2$ have pure orbifold neighborhoods.
Then $K_1 \cup K_2$ has a pure orbifold neighborhood.
\end{lem}
\begin{proof}
Let $(U_i,E_i,s_i,\psi_i)$ be a pure orbifold neighborhood of $K_i$.
Note that $\dim U_1 = \dim U_2 = \frak d$.
We denote the map $\underline\phi_{*p}$ the open set $U_{*p}$ etc. for
$(U_i,E_i,s_i,\psi_i)$ by $\underline\phi_{ip}$, $U_{ip}$ etc.
(Namely we replace $*$ by $i \in \{1,2\}$.)
The open subset $\mathcal U_i$ is as in (\ref{calUnasi}).
\par
Let $q \in K_1 \cap K_2$. We take an open subset $U_{12;q}$ such that
\begin{equation}
o_q \in U_{12;q} \subset  U_{1q} \cap U_{2q} \subset U_q
\end{equation}
and
\begin{equation}
\mathcal U_{12;q} = \mathcal U_{1q} \cap \mathcal U_{2q} \subset X.
\end{equation}
Here
$
\mathcal U_{iq} = \psi_q(s_q^{-1}(0) \cap U_{iq}).
$
We take $q_1,\dots, q_I$ such that
$$
K_1 \cap K_2 \subseteq \bigcup_{i=1}^I \mathcal U_{12;q_i}.
$$
We take a relatively compact open subset $U^-_{12;q}$ in $U_{12;q}$ such that
$$
K_1 \cap K_2 \subseteq \bigcup_{i=1}^I \mathcal U^-_{12;q_i}; \quad
\mathcal U^-_{12;q} := \psi_q(s_q^{-1}(0) \cap U^-_{12;q}).
$$
By a standard argument in general topology we can choose them
so that the following holds.
\begin{conds}\label{cond46}
If $\mathcal U_{12;q_i} \cap \mathcal U_{12;q_{i'}} \ne \emptyset$,
then $\mathcal U_{12;q_i} \cap \mathcal U_{12;q_{i'}}
\cap X(\frak d) \ne \emptyset$.
\end{conds}
We assume the same condition for $\{\mathcal U^-_{12;q_i}\}$.
\par
For each $r \in K_1 \cap K_2$ we take an open subset $U^0_r$ of $U_r$
containing $o_r$ so that
\begin{conds}\label{cond47}
\begin{enumerate}
\item $U^0_r \subset U_{1r} \cap U_{2r}$.
\item
If $\underline\phi_{1r}(U^0_r) \cap \overline{\underline\phi_{1q_i}(U_{12;q_i}^-)} \ne \emptyset$,
then
$$
U^0_r \subset U_{q_ir} \cap \underline{\phi}_{q_ir}^{-1}(U_{12;q_i}).
$$
\item
If $\underline\phi_{2r}(U^0_r) \cap \overline{\underline\phi_{2q_i}(U_{12; q_i}^-)} \ne \emptyset$,
then
$$
U^0_r \subset U_{q_ir}\cap \underline{\phi}_{q_ir}^{-1}(U_{12;q_i}).
$$
\end{enumerate}
\end{conds}
\begin{lem}\label{existsUr0}
There exists such a choice of $U^0_r$.
\end{lem}
\begin{proof}
We first observe
\begin{equation}\label{75ineq}
s_1^{-1}(0) \cap  \overline{\underline\phi_{1q_i}(U_{12; q_i}^-)}
=
s_1^{-1}(0) \cap  {\underline\phi_{1q_i}(\overline U_{12; q_i}^-)}
\subset
s_1^{-1}(0) \cap  {\underline\phi_{1q_i}(U_{12; q_i})}.
\end{equation}
We also have a similar inclusion when we replace $\phi_{1q_i}$ etc. by $\phi_{2q_i}$ etc..

We next decompose $I$ into disjoint unions
$$
I = I_1^{\rm in} \cup I_1^{\rm out} = I_2^{\rm in} \cup I_2^{\rm out}
$$
respectively where we define
$$
I_1^{\rm in} = \left\{ i \in I \mid
r \in \psi_{q_i} (s_{q_i}^{-1}(0)), \  \underline\phi_{1q_i}(o_r) \in
\overline{\underline{\phi}_{1q_i}(U_{12;q_i}^-)}\right\}, \quad
I_1^{\rm out} = I \setminus I_1^{\rm in}
$$
and similarly for $I_2^{\rm in}$, $I_2^{\rm out}$. Then we put
\begin{equation}\label{eq:Ur0}
\aligned
&U_r^0 =\\
&\left(\bigcap_{i \in I_1^{\rm in}} \left(U_{q_ir}\cap \underline\phi_{q_ir}^{-1}(U_{12;q_i})\right)\cap
\bigcap_{i \in I_1^{\rm out}}\underline\phi_{1r}^{-1}
\left(U_1 \setminus \overline{\underline\phi_{1q_i}(U_{12; q_i}^-)}\right)
\right) \\
&\cap \left(\bigcap_{i \in I_2^{\rm in}} \left(U_{q_ir}
\cap \underline\phi_{q_ir}^{-1}(U_{12;q_i})\right)\cap \bigcap_{i \in I_2^{\rm out}}
\underline\phi_{2r}^{-1}\left(U_2 \setminus \overline{\underline\phi_{2q_i}(U_{12; q_i}^-)}\right)\right).
\endaligned
\end{equation}
Using (\ref{75ineq}), it is easy to see that $U_r^0$ is an open neighborhood of $o_r$ satisfying all
the required properties.
\end{proof}

We choose $r_1,\dots,r_J \in K_1\cap K_2$ such that
\begin{equation}\label{eq:cup1toJ}
\bigcup_{j=1}^J \mathcal U_{r_j}^0 \supset K_1 \cap K_2.
\end{equation}
We put
\begin{equation}
\aligned
U_{21}^{(1)}
&=
\bigcup_{i=1}^I\bigcup_{j=1}^J
\left(
\underline{\phi}_{1r_j}(U^0_{r_j}) \cap
\underline{\phi}_{1q_i}(U^-_{12;q_i})
\right)
\subset U_1, \\
U_{12}^{(1)}
&=
\bigcup_{i=1}^I\bigcup_{j=1}^J
\left(
\underline{\phi}_{2r_j}(U^0_{r_j}) \cap
\underline{\phi}_{2q_i}(U^-_{12;q_i})
\right)
\subset U_2.
\endaligned
\end{equation}
They are open subsets in orbifolds and so are orbifolds.
We note that
$$
\bigcup_{i=1}^I \left( \underline{\phi}_{1r_j}(U^0_{r_j}) \cap
\underline{\phi}_{1q_i}(U_{12;q_i})
\right) \supset \underline{\phi}_{1r_j}(U_{r_j}^0)
$$
and
$\mathcal U_{12;q_i}$, $i=1, \dots, I$ cover $K_1 \cap K_2$
(similarly for $\underline{\phi}_{2r_j}$, $\underline{\phi}_{2q_i}$).
Hence
\begin{equation}
U_{21}^{(1)}  \supset \psi_1^{-1}(K_1\cap K_2),
\qquad
U_{12}^{(1)}  \supset \psi_2^{-1}(K_1\cap K_2)
\end{equation}
by \eqref{eq:cup1toJ}.
\begin{lem}\label{ofdgluemainsublam}
There exists an open embedding of orbifolds $\underline{\phi}_{21} : U_{21}^{(1)} \to U_2$
such that $\underline{\phi}_{21}(U_{21}^{(1)}) = U_{12}^{(1)}$ and satisfies the following:
\begin{enumerate}
\item
If $x = \underline{\phi}_{1r_j}(\tilde x_j)$
then
\begin{equation}\label{phi21defformula}
\underline{\phi}_{21}(x) = \underline{\phi}_{2r_j}(\tilde x_j).
\end{equation}
\item There exists a bundle isomorphism
\begin{equation}
\hat{\underline{\phi}}_{21} : E_1\vert_{U_{21}^{(1)}} \to E_2\vert_{U_{12}^{(1)}}
\nonumber\end{equation}
over $\underline{\phi}_{21}$. On the fiber of $x = \underline{\phi}_{1r_j}(\tilde x_j)$
we have
\begin{equation}\label{hatphi21defformula}
\hat{\underline{\phi}}_{21} = \hat{\underline{\phi}}_{2r_j} \circ \hat{\underline{\phi}}_{1r_j}^{-1}.
\end{equation}
\item
On $U_{21}^{(1)}$ we have:
\begin{equation}\label{compatis428}
s_2 \circ \underline{\phi}_{21} = \hat{\underline{\phi}}_{21} \circ s_1.
\end{equation}
\item On $s_1^{-1}(0) \cap U_{21}^{(1)}$, we have:
\begin{equation}\label{compatis429}
\psi_2 \circ \underline{\phi}_{21} = \psi_1.
\end{equation}
\end{enumerate}
\end{lem}
\begin{proof} For the statement (1),
we note that the right hand side of (\ref{phi21defformula}) is well-defined because of
Condition \ref{cond47} (1).
So to define $\underline{\phi}_{21}$ it suffices to show that the right hand side of (\ref{phi21defformula})
is independent of $j$.
Suppose
$$
x = \underline{\phi}_{1r_j}(\tilde x_j)  = \underline{\phi}_{1r_{j'}}(\tilde x_{j'})
\in
\underline{\phi}_{1q_i}(U^-_{q_i}).
$$
By Condition \ref{cond47}, (2) we have $\tilde x_j \in U_{q_ir_j}$,
$\tilde x_{j'} \in U_{q_ir_{j'}}$ and
$\underline{\phi}_{q_ir_j}(\tilde x_j) \in U_{12;q_i}$,
$\underline{\phi}_{q_ir_{j'}}(\tilde x_{j'}) \in U_{12;q_i}$.
Since
$$
\underline{\phi}_{1q_i}(\underline{\phi}_{q_ir_j}(\tilde x_j))
= x
= \underline{\phi}_{1q_i}(\underline{\phi}_{q_ir_{j'}}(\tilde x_{j'})),
$$
it follows that
$$
\underline{\phi}_{q_ir_j}(\tilde x_j)
=
\underline{\phi}_{q_ir_{j'}}(\tilde x_{j'}).
$$
Therefore
$$
\underline{\phi}_{2r_j}(\tilde x_j)
=
\underline{\phi}_{2q_i}(\underline{\phi}_{q_ir_j}(\tilde x_j))
=
\underline{\phi}_{2q_i}(\underline{\phi}_{q_ir_{j'}}(\tilde x_{j'}))
=
\underline{\phi}_{2r_{j'}}(\tilde x_{j'}),
$$
as required.
\par
We have thus defined $\underline{\phi}_{21}$. We can define $\underline{\phi}_{12}$ in a similar way.
It is easy to see $\underline{\phi}_{21}\circ \underline{\phi}_{12}$
and  $\underline{\phi}_{12}\circ \underline{\phi}_{21}$
are identity maps. Therefore $\underline{\phi}_{21}$ is an isomorphism.
\par
For (2), We define $\hat{\underline{\phi}}_{21}$ by (\ref{hatphi21defformula}).
We can prove that it is well-defined and is an isomorphism in the same way
as the proof of (1).
\par
Finally for (3) and (4), we note that
(\ref{compatis428}) follows from (\ref{phi21defformula}) and (\ref{hatphi21defformula}).
(\ref{compatis429}) follows from (\ref{phi21defformula}).
\end{proof}

We now use the bundle isomorphisms $(\widehat{\underline{\phi}}_{21},\underline{\phi}_{21})$
to glue the pure orbifold neighborhoods $(U_1,E_1,s_1,\psi_1)$ and $(U_2,E_2,s_2,\psi_2)$.
\par
We first shrink $U_1, \, U_2$ so that Definition \ref{pureofdnbd} (4)  holds as follows.
We choose open subsets $V_1, V_2 \subset X$ such that
\begin{equation}
K_1 \subset V_1 \subset \psi_1(s_1^{-1}(0)), \qquad K_2 \subset V_2 \subset \psi_2(s_2^{-1}(0)),  \qquad
V_1 \cap V_2 \subset \psi_1(U_{21}^{(1)} \cap s_1^{-1}(0)).
\end{equation}
(Such $V_1$, $V_2$ exist since $K_1 \cap K_2 \subset \psi_1(U_{21}^{(1)} \cap s_1^{-1}(0))$.)
\par
We take open sets $U'_i \subset U_i$ such that $U'_i \cap s_i^{-1}(0) = \psi_i^{-1}(V_i)$ and
set $U^{\prime (1)}_{21} = U^{(1)}_{21} \cap U'_1 \cap \underline{\phi}_{21}^{-1}(U_2')$,  $U^{\prime (1)}_{12} = U^{(1)}_{12} \cap U'_2 \cap \underline{\phi}_{12}^{-1}(U_1')$.
By restricting the base of the bundle $E_i$,
the maps $\underline{\phi}_{21}$ etc are defined in an obvious way and satisfy the conclusion
of Lemma \ref{ofdgluemainsublam}.
Moreover we have
\begin{equation}\label{intproperty}
\psi_1(U'_{1} \cap s_1^{-1}(0)) \cap \psi_2(U'_{2} \cap s_2^{-1}(0))
\subset
\psi_1(U_{21}^{\prime (1)} \cap s_1^{-1}(0)).
\end{equation}
(\ref{intproperty}) implies that $\psi_1$, $\psi_2$ induce an
{\it injective} map $(U'_{1} \cap s_1^{-1}(0)) \#_{\underline{\phi}_{21}}
(U'_{2} \cap s_2^{-1}(0)) \to X$.
\par
This would have finished the proof of Lemma \ref{2gluelemma}, if we already established that
the glued space $U_1' \#_{\underline\phi_{21}} U_2'$ is Hausdorff, which is
not guaranteed with the current construction.
In order to obtain a Hausdorff space after gluing we need to further shrink the domains as follows.
\par
The argument of this shrinking process will be somewhat similar to that of Section \ref{defgoodcoordsec}.
Let $U_1^0 \subset U'_1$ and $U_2^0 \subset U'_2$ be relatively compact open subsets such that
$$
\psi_1(s_1^{-1}(0) \cap U_1^0) \supset K_1,
\qquad
\psi_2(s_2^{-1}(0) \cap U_2^0) \supset K_2.
$$
Let $W_{21} \subset U_{21}^{(1)}$ be a relatively compact open subset such that
\begin{equation}\label{Wcond}
W_{21} \supset s_1^{-1}(0) \cap \overline {(U_1^0 \cap \underline{\phi}_{21}^{-1}(U_2^0))}.
\end{equation}
The existence of $W_{21}$ with this property follows from the next lemma.
\begin{lem}\label{lemma49}
$s_1^{-1}(0) \cap \overline{(U_1^0 \cap \underline{\phi}_{21}^{-1}(U_2^0))}$ is compact.
\end{lem}
\begin{proof}
$\psi_1(\overline U_{1}^{0} \cap s_1^{-1}(0))$ is a compact subset of $V_1$.
Therefore $\overline U_{1}^{0} \cap s_1^{-1}(0)$ is compact. Since
$s_1^{-1}(0) \cap \overline{(U_1^0 \cap \underline{\phi}_{21}^{-1}(U_2^0))}$ is its closed subset,
the lemma follows.
\end{proof}
We next prove:
\begin{lem}\label{lem410}
There exists an open subset $U_1^{00}$ such that
$$
\psi_1^{-1}(K_1) \subset U_1^{00} \subset  \overline U_1^{00} \subset U_1^0
$$
such that
\begin{equation}\label{415form}
\overline U_1^{00} \cap \overline{\underline{\phi}_{21}^{-1}(U^0_2)} \subset W_{21}.
\end{equation}
\end{lem}
\begin{proof}
For each $x\in K_1$ we define its neighborhood $U_x$'s as follows.
\begin{enumerate}
\item
If $x \notin \overline{\underline{\phi}_{21}^{-1}(U^0_2)}$ then
$\overline U_x \cap \overline{\underline{\phi}_{21}^{-1}(U^0_2)} = \emptyset$.
\item
If $x \in \overline{\underline{\phi}_{21}^{-1}(U^0_2)}$ then $\overline U_x \subset W_{21}$.
\end{enumerate}
Note $K_1 \cap \overline{\underline{\phi}_{12}(U^0_2)} \subset W_{21}$ by
(\ref{Wcond}). Therefore we can find such $U_x$ in Case (2).
We cover $K_1$ by a finite number of such $U_x$'s.
Let $U_1^{00}$ be its union. Then it has the
required properties.
\end{proof}
Let $C_2$ be the closure of $U_2^0$ in $U'_2$
and
$C_1$ be the closure of $U_1^{00}$ in $U'_1$.
We put
$$
C_{21} = C_1 \cap \underline{\phi}_{21}^{-1}(C_2),
\quad
C_{12} = \underline{\phi}_{21}(C_{21}).
$$
We put
\begin{equation}\label{eq:simonC}
C = (C_1 \cup C_2)/\sim
\end{equation}
where $\sim$ is defined as follows $x\sim y$ if and only if one of the
following holds.
\begin{enumerate}
\item $x = y$, or
\item $x \in C_{21}$, $y =  \underline{\phi}_{21}(x) \in C_2$, or
\item $x \in C_{12}$, $y =  \underline{\phi}_{12}(x) \in C_1$.
\end{enumerate}
We put the quotient topology on $C$.
\begin{lem}\label{C21compact}
$C_{21}$ is a compact subset of $C_1$.
\end{lem}
\begin{proof}
If suffices to show that $C_{21}$ is a closed subset of $C_1$.
But this is obvious since $\underline\phi_{21}$ is a continuous map to
$U_2$ and $C_2$ is a closed subset of $U_2$.
\end{proof}
\begin{lem}\label{sublemhaus}
$C$ is Hausdorff.
\end{lem}
\begin{proof}
This is a consequence of  Proposition \ref{metrizable}.
\end{proof}

Let $\Pi_i : C_i \to C$ be the map which sends an element to its
equivalence class. Since $C_i$ is compact and $C$ is Hausdorff, the map
$\Pi_i$ is a topological embedding.
We set
$$
U = \Pi_1(U^{00}_1) \cup \Pi_2(U^0_2) \subset C
$$
and put the subspace topology on it. This set
carries an orbifold structure since each of $U_1^{00}$
and $U_2^0$ does and the gluing is done by (orbifold) diffeomorphism.
Then by restricting $E_i$, $s_i$, $\psi_i$, etc and then gluing them,
we obtain the required gluings $E$, $s$, $\psi$ and etc.

\begin{lem}\label{XcapU} The quadruple $(U,E,s,\psi)$ constructed above
satisfies the conditions of a pure orbifold neighborhood of $K_1 \cup K_2$.
\end{lem}
\begin{proof}
Definition \ref{pureofdnbd} (1),(2) and (3) are obvious.
We can use (\ref{intproperty}), (\ref{415form}) to show that $\psi : s^{-1}(0) \to X$ is injective
where $s^{-1}(0) \subset U$.
Therefore it is an embedding, which is also open.
The proof of Lemma \ref{XcapU} is complete.
\end{proof}
Therefore the proofs of Proposition \ref{purecover} and Lemma  \ref{2gluelemma} are complete.
\end{proof}
\par\medskip
We have thus completed the first part of the proof of Theorem \ref{goodcoordinateexists}
and enter the second part.
\par
Put
$$
\frak D = \{\frak d \in \Z_{>0} \mid X(\frak d) \ne \emptyset \}
$$
and let $(\frak D,\le)$ be an ordered set.
A subset $D \subset \frak D$ is said to be an {\it ideal}
if
$$
\frak d \in D, \,\frak d' \ge \frak d
\quad
\Rightarrow \quad \frak d' \in D.
$$
We call a subset $D' \subset D$ of an ideal $D$ a {\it sub-ideal} thereof if
$D'$ is an ideal with respect to the induced order.

For $D \subset \frak D$ we put
$$
X(D) = \bigcup_{\frak d \in D} X(\frak d).
$$
Then $X(D)$ is a closed subset of $X$ if $D$ is an ideal.
\begin{defn}\label{mixednbd}
Let $D \subset \frak D$ be an ideal.
A \emph{mixed orbifold neighborhood} of $X(D)$ is given by the quadruples
$$
(\mathcal K_{\frak d} \subset U_{\frak d}, E_{\frak d},s_{\frak d}, \psi_{\frak d})
$$
together with embeddings $(\underline{\hat\phi}_{\frak d\frak d'}, \underline\phi_{\frak d\frak d'})$
of vector bundles for each pair $(\frak d,\frak d')$ with $\frak d < \frak d'$.

We assume they have the following properties:
\begin{enumerate}
\item
$\mathcal K_{\frak d}$ is a compact subset of $X(\frak d)$.
\item
$(U_{\frak d}, E_{\frak d},s_{\frak d}, \psi_{\frak d})$ is a pure
orbifold neighborhood of $\mathcal K_{\frak d}$.
\item Let $\psi_{\frak d} : U_{\frak d} \cap s_{\frak d}^{-1}(0) \to X$ be as in
Definition \ref{pureofdnbd} (3). Then we have
\begin{equation}\label{Kislargerthan}
\mathcal K_{\frak d}
\supset
X(\frak d) \setminus \bigcup_{\frak d' > \frak d}
\psi_{\frak d'} \left(U_{\frak d'} \cap s_{\frak d'}^{-1}(0)\right).
\end{equation}
\item
$U_{\frak d'\frak d}$ is an open neighborhood of
$$
\psi_{\frak d}^{-1}\left(
\mathcal K_{\frak d}\cap \psi_{\frak d'} \left(U_{\frak d'} \cap s_{\frak d'}^{-1}(0)\right)
\right)
$$
in $U(\frak d)$.
\item The map
$$
\underline{\hat\phi}_{\frak d'\frak d} : E_{\frak d}\vert_{U_{\frak d'\frak d}} \to E_{\frak d'}
$$
is an embedding of vector bundle over an embedding of orbifold
$$
\underline\phi_{\frak d'\frak d} : U_{\frak d'\frak d} \to U_{\frak d'}
$$
such that
\begin{enumerate}
\item The equality
$$
s_{\frak d'} \circ \underline\phi_{\frak d'\frak d}
= \underline{\hat\phi}_{\frak d'\frak d}\circ s_{\frak d}
$$
holds on $U_{\frak d'\frak d}$.
\item The equality
$$
\psi_{\frak d'} \circ \underline\phi_{\frak d'\frak d}
= \psi_{\frak d}
$$
holds on $U_{\frak d'\frak d} \cap s_{\frak d}^{-1}(0)$.
\end{enumerate}
\item
The restriction of $ds_{\frak d'}$ to the normal direction induces an isomorphism
\begin{equation}\label{tangent***}
N_{U_{\frak d'\frak d}}U_{\frak d'} \cong
\frac{\hat{\phi}_{\frak d'\frak d}^*E_{\frak d'}}{E_\frak d\vert_{U_{\frak d'\frak d}}}
\end{equation}
as vector bundles on $U_{\frak d'\frak d}\cap s_{\frak d}^{-1}(0)$.
\item
If
$p \in \psi_{\frak d}(U_{\frak d'\frak d} \cap s_{\frak d}^{-1}(0)) \subset \psi_{\frak d'}(U_{\frak d'}
\cap s_{\frak d'}^{-1}(0))$,
$$
\underline\phi_{\frak d'\frak d}\circ \underline\phi_{\frak d p} = \underline\phi_{\frak d' p}, \quad
\underline{\hat\phi}_{\frak d'\frak d}\circ \underline{\hat\phi}_{\frak d p} = \underline{\hat\phi}_{\frak d' p}
$$
on $\underline{\phi}_{\frak d p}^{-1}(U_{\frak d' \frak d}) \cap U_{\frak d p} \cap U_{{\frak d}'p}$
and $E_p\vert_{\underline{\phi}_{\frak d p}^{-1}(U_{\frak d' \frak d}) \cap U_{\frak d p} \cap U_{{\frak d}'p}}$, respectively.
(We write $\underline{\phi}_{\frak d p}$ etc. instead of
$\underline\phi_{* p}$ etc. for the structure maps of $U_{\frak d}$.
Namely we replace $*$ by $\frak d$.)
\item
The space $U(D)$ is Hausdorff. The map
${\Pi}_{\frak d} :  U_{\frak d} \to U(D)$ is a topological embedding.
We have
\begin{equation}\label{Dddd}
U_{\frak d'\frak d} = {\Pi}_{\frak d}^{-1}{\Pi}_{\frak d'}(U_{\frak d'})
\end{equation}
and
\begin{equation}\label{Icompati}
{\Pi}_{\frak d'} \circ \underline{\phi}_{\frak d'\frak d} = {\Pi}_{\frak d}
\end{equation}
on $U_{\frak d'\frak d}$.
Moreover
$$
U(D) = \bigcup_{\frak d \in D} I_{\frak d}(U_{\frak d}).
$$
We call $U(D)$ the {\it total space} of our  mixed orbifold neighborhood.
\item
We define a subset $s_D^{-1}(0)$ of $U(D)$ by $s_D^{-1}(0)
= \bigcup_{\frak d \in D}{\Pi}_{\frak d}(s_\frak d^{-1}(0))$.
We define $\psi_D : s_D^{-1}(0) \to X$ such that $\psi_D\circ {\Pi}_{\frak d} = \psi_{\frak d}$
on $s_\frak d^{-1}(0) \subset U_{\frak d}$.
(This is well-defined by (7).)
We require that
$$
\psi_D : s_D^{-1}(0) \to X
$$
is a topological embedding onto a neighborhood of $X(D)$ in $X$.
\end{enumerate}
\end{defn}

Note that (\ref{Kislargerthan}) implies
\begin{equation}
X(D) \subset \bigcup_{\frak d \in D}
\psi_{\frak d} \left(s_{\frak d}^{-1}(0)\right).
\end{equation}
We also have the following:
\begin{lem}\label{mixednbdrestrict}
If $U(D)'$ is an open subset of $U(D)$ such that
$$
U(D)' \supset \psi_D^{-1}(X(D)) \cap \bigcup_{\frak d \in D} {\Pi}_{\frak d}(s_{\frak d}^{-1}(0)) .
$$
Then there exists a mixed orbifold neighborhood of $X(D)$ such that
its total space is the above $U(D)'$.
\end{lem}
\begin{proof}
We put $U'_{\frak d} = U_{\frak d} \cap {\Pi}_{\frak d}^{-1}(U(D)')$,
$U'_{\frak d'\frak d} = {\Pi}_{\frak d}^{-1}({\Pi}_{\frak d'}(U_{\frak d'}) \cap U(D)')$.
We define $E'_{\frak d}$ and various maps by restricting ones of $U(D)$.
It is straightforward to check that they satisfy the required
properties (1)-(11) of Definition \ref{mixednbd}.
\end{proof}
Another way to shrink $U(D)$ is as follows.
\begin{lem}\label{mixednbdrestrict2}
Let $
(\mathcal K_{\frak d} \subset U_{\frak d}, E_{\frak d},s_{\frak d}, \psi_{\frak d})
$
be a mixed orbifold neighborhood of $X(D)$.
Suppose we are given $U'_{\frak d} \subset U_{\frak d}$ for each $\frak d$
such that
$$
\psi_{\frak d}(s_{\frak d}^{-1}(0) \cap U'_{\frak d}) \supset \mathcal K_{\frak d}.
$$
Then  $U'_{\frak d}$ and the restriction of the other data define
a mixed orbifold neighborhood of $X(D)$.
\end{lem}
\begin{proof}
The proof is obvious.
\end{proof}
\par
The goal of the second part of the proof of Theorem \ref{goodcoordinateexists}
is to prove the following:
\begin{prop}\label{existmixed}
For any ideal $D$ there exists a mixed orbifold neighborhood of $X(D)$.
\end{prop}
\begin{proof}
The proof is by an induction on $\# D$.
If $\# D = 1$ then $D = \{\frak d\}$ with $\frak d$ maximal in $\frak D$.
We note that $X(\frak d)$ is compact when $\frak d$ is maximal.
We put $\mathcal K_{\frak d} = X(\frak d)$ that is compact.
We use Proposition \ref{purecover} to obtain a
pure orbifold neighborhood $U_{\frak d}$ of $\mathcal K_{\frak d} = X(\frak d)$.
The proposition is proved in this case.
\par
Suppose we have proved the proposition for all $D'$ with $\# D' < \#D$.
We will prove it for $D$. Let $\frak d_0$ be an element of $D$ that is
minimal. We put $D' = D \setminus \{\frak d_0\}$, which is
a sub-ideal of $D$ with $\# D' < \#D$. Therefore, by the induction hypothesis,
we have a mixed orbifold neighborhood of $X(D')$.
We denote it by $U^{(1)}_{\frak d}$, $\mathcal K_{\frak d}$,
$\underline{\phi}_{\frak d'\frak d}$ etc. (Here $\frak d,\frak d' \in D'$.)
\par
Let $\mathcal K_{\frak d_0}$ be a compact subset of $X(\frak d_0)$
such that
\begin{equation}
 \bigcup_{\frak d \in D'}\psi^{(1)}_{\frak d}
\left(U^{(1)}_{\frak d} \cap (s_{\frak d}^{(1)})^{-1}(0)\right)
\supset
\overline{(X(\frak d_0) \setminus \mathcal K_{\frak d_0})}.
\end{equation}
We apply Proposition \ref{purecover}  to  $\mathcal K_{\frak d_0}$ to obtain
$U^{(1)}_{\frak d_0}$.

The main part of the proof is to glue $U^{(1)}_{\frak d_0}$ with $X(D')$
to obtain the required mixed orbifold neighborhood of $X(D)$.
The construction is similar to the proof of Lemma \ref{ofdgluemainsublam}.

\begin{rem} We would like to remark that one outstanding difference between
the gluing maps $\underline \phi_{ij}$ for the pure orbifold neighborhoods and
$\underline\phi_{\frak d'\frak d}$ for the mixed orbifold neighborhoods is that the
former is a diffeomorphism while the latter is only an embedding.
\end{rem}

Because of this, we would like to repeat the detail of the gluing process.
\par
Let $U^{(2)}(D')$  be a relatively compact open subset of $U^{(1)}(D')$
satisfying
\begin{equation}\label{U2cond1}
U^{(2)}(D') \supset \psi_{D'}^{-1}(X(D')).
\end{equation}
We may choose it sufficiently close to $U^{(1)}(D')$ such that
\begin{equation}\label{U2cond2}
\bigcup_{\frak d \in D'}\psi^{(1)}_{\frak d}
\left(U^{(2)}_{\frak d} \cap (s_{\frak d}^{(1)})^{-1}(0)\right)
\supset
\overline{(X(\frak d_0) \setminus \mathcal K_{\frak d_0})}.
\end{equation}
Here $U^{(2)}_{\frak d} $ is obtained from $U^{(2)}(D')$ as in the proof of Lemma \ref{mixednbdrestrict}.
\par
Let $U^{(2)}_{\frak d_0} \subset U^{(1)}_{\frak d_0}$ be a relatively compact open subset
such that
\begin{equation}\label{U2cond3}
(\psi_{\frak d_0}^{(1)})^{-1}(\mathcal K_{\frak d_0}) \subset U^{(2)}_{\frak d_0}.
\end{equation}
We remark (\ref{U2cond2}) and (\ref{U2cond3}) imply
\begin{equation}\label{U2cond4}
\psi^{(1)}_{\frak d_0} \left((s_{\frak d_0}^{(1)})^{-1}(0)\cap U^{(2)}_{\frak d_0}\right)
\cup
\psi^{(1)}_{D'} \left((s_{D'}^{(1)})^{-1}(0)\cap U^{(2)}(D')\right)
\supset X(\frak d_0).
\end{equation}
\par
We put
\begin{equation}\label{Ldefn}
\aligned
\mathcal L_{\frak d_0} = X(\frak d_0) &\cap
\psi^{(1)}_{\frak d_0}\left ((s_{\frak d_0}^{(1)})^{-1}(0)
\cap U^{(2)}_{\frak d_0}\right)\\
&\cap\psi^{(1)}_{D'} \left((s_{D'}^{(1)})^{-1}(0)\cap U^{(2)}(D')\right).
\endaligned
\end{equation}
\par
For each $q \in \overline{\mathcal L}_{\frak d_0}$ we take an open neighborhood $U^{(01)}_q$ in
$U_q$ that satisfies the following conditions.
\begin{conds}
\begin{enumerate}
\item
If $\frak d > \frak d_0$ and $q \in \overline{\psi^{(1)}_{\frak d}((s_{\frak d}^{(1)})^{-1}(0)\cap U_{\frak d}^{(2)})}\cap
\overline{\mathcal L}_{\frak d_0}$,
then $U^{(01)}_q \subset U^{(1)}_{\frak d q}$.
Here $U^{(1)}_{\frak d q}$ is as in Definition \ref{pureofdnbd} (5) for the pure orbifold neighborhood $U^{(1)}_{\frak d}$.
\item
$U^{(01)}_q \subset U^{(1)}_{\frak d_0 q}$, where  $U^{(1)}_{\frak d_0 q}$ is as in Definition \ref{pureofdnbd} (5) for the pure orbifold neighborhood $U^{(1)}_{\frak d_0}$.
\end{enumerate}
\end{conds}
For a finite subset $\frak I = \{q_i \mid i = 1,\dots, I\} \subset \overline{\mathcal L}_{\frak d_0}$ and
each given $\frak d \in D'$, we write
\begin{equation}\label{eq:Ifrakd}
\frak I_{\frak d} = \frak I \cap  \overline{\psi^{(1)}_{\frak d}((s_{\frak d}^{(1)})^{-1}(0)\cap U_{\frak d}^{(2)})}
\cap \overline{\mathcal L}_{\frak d_0}.
\end{equation}
We also denote
$$
\mathcal U^{(21)}_{\frak d}:= \psi^{(1)}_{\frak d}((s_{\frak d}^{(1)})^{-1}(0)\cap U_{\frak d}^{(2)})
\subset \mathcal U^{(1)}_{\frak d}.
$$
By the compactness of $\overline{\mathcal L}_{\frak d_0}$, it follows that
there exists a finite subset $\frak J$ of $\overline{\mathcal L}_{\frak d_0}$
such that the following condition holds.
\begin{conds}\label{cond415}
For any $\frak d \in D'$ we have:
$$
\bigcup_{q_i \in \frak I_{\frak d}}
\mathcal U_{q_i}^{(01)} \supset
\overline{\mathcal U^{(21)}_{\frak d}}
\cap \overline{\mathcal L}_{\frak d_0}.
$$
Here $\mathcal U_{q_i}^{(01)} = \psi_{q_i}^{(1)}((s_{q_i}^{(1)})^{-1}(0)\cap U_{q_i}^{{(01)}})$.
\end{conds}
Since
$$
\overline{\mathcal L}_{\frak d_0} \subset \bigcup_{\frak d \in D'}
\overline{\mathcal U^{(21)}_{\frak d}},
$$
Condition \ref{cond415} implies
\begin{equation}
\bigcup_{q_i \in \frak J_{\frak d}} \mathcal U_{q_i}^{(01)} \supset \overline{\mathcal L}_{\frak d_0}.
\end{equation}
We may assume that $\{\mathcal U_{q_i}^{(01)}\}$ satisfies Condition \ref{cond46}.
\par
We next take a relatively compact open subset $U_{q_i}^{(01)-}$ of $U_{q_i}^{{(01)}}$ such that
the following holds.
\begin{conds}\label{cond416}
For any $\frak d \in D'$, we have:
$$
\bigcup_{q_i \in \overline{\mathcal U^{(21)}_{\frak d}
\cap \overline{\mathcal L}_{\frak d_0}}}
\mathcal U_{q_i}^{(01)-} \supset
\overline{\mathcal U^{(21)}_{\frak d}}
\cap \overline{\mathcal L}_{\frak d_0}.
$$
Here $\mathcal U_{q_i}^{(01)-} = \psi^{(1)}_{q_i}((s^{(1)}_{q_i})^{-1}(0)\cap U_{q_i}^{(01)-})$.
\end{conds}
We also assume that $\{\mathcal U_{q_i}^{(01)-}\}$ satisfies Condition \ref{cond46}.
\par
For $r \in \overline{\mathcal L}_{\frak d_0}$ we take an open neighborhood $U_r^{0}$ of $o_r$ in $U_r$ with the
following properties.
\begin{conds}\label{cond417}
\begin{enumerate}
\item
If $\frak d >\frak d_0$ and $r \in \overline{\mathcal U^{(21)}_{\frak d}}\cap \overline{\mathcal L}_{\frak d_0}$,
then $U_r^0 \subset U_{\frak d r}^{(1)}$.
\item
$U_r^0 \subset U_{\frak d_0 r}^{(1)}$.
\item
If $\underline{\phi}^{(1)}_{\frak d_0 r}(U_r^0) \cap
\overline{\underline{\phi}^{(1)}_{\frak d_0 q_i}(U_{q_i}^{(01)-})} \ne \emptyset$
then $U_r^0 \subset U_{q_i r}^{(1)} \cap (\underline{\phi}_{q_i r}^{(1)})^{-1}(U^{(01)}_{q_i})$.
\item If
$\underline{\phi}^{(1)}_{\frak d r}(U_r^0) \cap
\overline{\underline{\phi}^{(1)}_{\frak d q_i}(U_{q_i}^{(01)-})} \ne \emptyset$
then $U^0_r \subset U^{(1)}_{q_i r} \cap ( \underline{\phi}_{q_i r}^{(1)})^{-1}(U_{q_i}^{(01)})$.
\end{enumerate}
\end{conds}
The existence of such $U^0_r$ is proved in the same way as Lemma \ref{existsUr0}.
We choose a subset
$$
\frak J = \{r_j \mid r_j \in\overline{\mathcal L}_{\frak d_0}, j=1, \cdots, J \}
\subset \overline{\mathcal L}_{\frak d_0}
$$
such that the following holds for each $\frak d \in D'$
\begin{equation}
\bigcup_{r_j \in \frak J_{\frak d}} \mathcal U_{r_j}^0
\supset
\overline{\mathcal U^{(21)}_{\frak d}} \cap \overline{\mathcal L}_{\frak d_0},
\end{equation}
where
$$
\frak J_{\frak d} =\frak J \cap  \overline{\psi^{(1)}_{\frak d}((s_{\frak d}^{(1)})^{-1}(0)\cap U_{\frak d}^{(2)})}\cap \overline{\mathcal L}_{\frak d_0}.
$$
We now put
\begin{equation}
U_{\frak d\frak d_0}^{(1)}
=\bigcup_{r_ j \in \frak J_{\frak d}}
\bigcup_{q_i \in \frak J_{\frak d}}
\underline{\phi}^{(1)}_{\frak d_0 r_j}(U_{r_j}^0) \cap \underline{\phi}^{(1)}_{\frak d_0 q_i}(U_{q_i}^{(01)-}).
\end{equation}
This is an open subset of $U_{\frak d_0}^{(1)}$.
Since $U_{\frak d_0}^{(1)}$ is an orbifold, its open subset $U_{\frak d\frak d_0}^{(1)}$ is also an orbifold.
We remark that
\begin{equation}\label{Liscontained1}
\psi^{(1)}_{\frak d_0} (U_{\frak d\frak d_0}^{(1)} \cap (s^{(1)}_{\frak d_0})^{-1}(0))
\supset \overline{\mathcal L}_{\frak d_0} \cap \mathcal U^{(21)}_{\frak d}.
\end{equation}

The following lemma is an analog to Lemma \ref{ofdgluemainsublam}. As we mentioned before,
one difference is that the gluing map $\underline\phi^{(1)}_{\frak d \frak d_0}$
is not an open embedding.

\begin{lem}\label{lem418}
There exists an embedding of orbifolds
$\underline\phi^{(1)}_{\frak d \frak d_0} : U_{\frak d\frak d_0}^{(1)}
\to U_{\frak d}^{(1)}$ with the following properties.
\begin{enumerate}
\item
If $x = \underline{\phi}^{(1)}_{\frak d_0 r_j}(\tilde x_j)$,
then
\begin{equation}\label{440}
\underline\phi^{(1)}_{\frak d \frak d_0}(x) = \underline{\phi}^{(1)}_{\frak d r_j}(\tilde x_j).
\end{equation}
\item
There exists an embedding of vector bundles
$$
\underline{\hat\phi}^{(1)}_{\frak d \frak d_0}
: E_{\frak d_0}^{(1)}\vert_{U_{\frak d\frak d_0}^{(1)}} \to E_{\frak d}^{(1)}
$$
that covers $\underline\phi^{(1)}_{\frak d \frak d_0}$.
\item
If $\frak d > \frak d'  > \frak d_0$ then
$$
\underline{\phi}^{(1)}_{\frak d \frak d_0} =
\underline{\phi}^{(1)}_{\frak d \frak d'} \circ \underline{\phi}^{(1)}_{\frak d' \frak d_0}
$$
on
$
(\underline{\phi}^{(1)}_{\frak d' \frak d_0})^{-1}(U_{\frak d \frak d'}^{(1)})
\cap U_{\frak d \frak d_0}^{(1)}
$
and
$$
\underline{\hat\phi}^{(1)}_{\frak d \frak d_0} =
\underline{\hat\phi}^{(1)}_{\frak d \frak d'} \circ \underline{\hat\phi}^{(1)}_{\frak d' \frak d_0}
$$
on
$
E^{(1)}_{\frak d_0}\vert_{(\underline{\phi}^{(1)}_{\frak d' \frak d_0})^{-1}(U_{\frak d \frak d'}^{(1)})
\cap U_{\frak d \frak d_0}^{(1)}
}
$.
\item
We have
$$
s_{\frak d}^{(1)} \circ \underline{\phi}^{(1)}_{\frak d \frak d_0}
=
\underline{\hat\phi}^{(1)}_{\frak d \frak d_0} \circ s^{(1)}_{\frak d_0}
$$
on $U_{\frak d \frak d_0}^{(1)}$.
\item
We have
$$
\psi_{\frak d}^{(1)} \circ \underline{\phi}^{(1)}_{\frak d \frak d_0}
=
\psi_{\frak d_0}^{(1)}
$$
on $U_{\frak d \frak d_0}^{(1)} \cap (\overline s_{\frak d_0}^{(1)})^{-1}(0)$.
\item
The restriction of $ds_{\frak d}^{(1)}$ to the normal direction induces an isomorphism
\begin{equation}\label{tangent***}
N_{U_{\frak d\frak d_0}}U_{\frak d} \cong \frac{E_{\frak d}}
{\hat{\phi}_{\frak d\frak d_0}(E_{\frak d_0}\vert_{U_{\frak d\frak d_0}})}
\end{equation}
as vector bundles on $U_{\frak d\frak d_0}\cap s_{\frak d_0}^{-1}(0)$.
\end{enumerate}
\end{lem}
\begin{proof}
The proof is similar to the proof of Lemma \ref{ofdgluemainsublam}.
\par
Note that the right hand side of (\ref{440}) is well defined because of Condition \ref{cond417} (1).
We first show that the right hand side of (\ref{440}) is independent of $j$.
Suppose
$$
x = \underline{\phi}^{(1)}_{\frak d_0 r_j}(\tilde x_j)= \underline{\phi}^{(1)}_{\frak d_0 r_{j'}}(\tilde x_{j'})
\in \underline{\phi}^{(1)}_{\frak d_0 q_i}(U_{q_i}^{(01)-}).
$$
Then by Condition \ref{cond417} (3) we have $\tilde x_j \in U^{(1)}_{q_i r_j}$, $\tilde x_{j'}
\in U^{(1)}_{q_i r_{j'}}$
and $\underline{\phi}^{(1)}_{q_i r_j}(\tilde x_j) \in U_{q_i}^{(01)}$,
$\underline{\phi}^{(1)}_{q_i r_{j'}}(\tilde x_{j'}) \in U_{q_i}^{(01)}$ .
Since
$$
\underline{\phi}^{(1)}_{\frak d_0 q_i}(\underline{\phi}^{(1)}_{q_i r_j}(\tilde x_j))
= x
=
\underline{\phi}^{(1)}_{\frak d_0 q_i}(\underline{\phi}^{(1)}_{q_i r_{j'}}(\tilde x_{j'})),
$$
it follows that
$$
\underline{\phi}^{(1)}_{q_i r_j}(\tilde x_j)
=
\underline{\phi}^{(1)}_{q_i r_{j'}}(\tilde x_{j'}).
$$
Therefore
$$
\underline{\phi}^{(1)}_{\frak d r_j}(\tilde x_j)
=\underline{\phi}^{(1)}_{\frak d q_i}(\underline{\phi}^{(1)}_{q_i r_j}(\tilde x_j))
=
\underline{\phi}^{(1)}_{\frak d q_i}(\underline{\phi}^{(1)}_{\frak d r_{j'}}(\tilde x_{j'}))
=
\underline{\phi}^{(1)}_{\frak d r_{j'}}(\tilde x_{j'})
$$
as required. We remark that $\underline{\phi}^{(1)}_{\frak d_0 r_j}$ is an open embedding of orbifolds.
Therefore
$\underline\phi^{(1)}_{\frak d \frak d_0}$ defined by (\ref{440}) is an embedding of orbifolds.
\par
The proof of (2) is similar. Then the proofs of (3)-(6) are straightforward.
\end{proof}
The pure orbifold neighborhoods $U^{(1)}_{\frak d}$ ($\frak d > \frak d_0$) and  $U^{(1)}_{\frak d_0}$
together with $\underline\phi_{\frak d\frak d_0}$ etc. satisfy
the properties required in Definition \ref{mixednbd}
except the following two points.
\begin{enumerate}
\item[(A)]
Hausdorff-ness of the space $U(D)$ that is required in Definition \ref{mixednbd} (8).
\item[(B)]
Injectivity of the map $\psi_D$ that is  required in Definition \ref{mixednbd} (9).
\end{enumerate}
\par
In the rest of the proof  we shrink $U^{(1)}_{\frak d}$, $U^{(1)}_{\frak d_0}$
again so that (A), (B) above are satisfied.
We shrink in such a way appearing either in Lemma \ref{mixednbdrestrict} or in Lemma \ref{mixednbdrestrict2}.
Therefore the other properties required in Definition \ref{mixednbd} than (A) (B) hold
after shrinking.
\par
We put
\begin{equation}\label{44344}
U_{\frak d \frak d_0}^{(2)} = (\underline\phi^{(1)}_{\frak d \frak d_0})^{-1}
(U_{\frak d}^{(2)}) \cap U_{\frak d_0}^{(2)}.
\end{equation}
Let $\underline\phi^{(2)}_{\frak d \frak d_0}$ and
$\underline{\hat\phi}^{(2)}_{\frak d \frak d_0}$ be the restrictions of $\underline\phi^{(1)}_{\frak d \frak d_0}$
and $\underline{\hat\phi}^{(1)}_{\frak d \frak d_0}$
to $U_{\frak d \frak d_0}^{(2)}$  and $E_{\frak d_0}\vert_{U_{\frak d \frak d_0}^{(2)}}$.
\par
The inclusion (\ref{Liscontained1}) and the definition of $U_{\frak d\frak d_0}^{(2)}$ imply
\begin{equation}\label{Liscontained2}
\psi^{(2)}_{\frak d_0} (U_{\frak d\frak d_0}^{(2)} \cap (s^{(2)}_{\frak d_0})^{-1}(0))
\supset \overline{\mathcal L}_{\frak d_0} \cap \psi^{(2)}_{\frak d}((s_{\frak d}^{(2)})^{-1}(0)
\cap U_{\frak d}^{(2)}) \cap \psi_{\frak d_0}^{(2)}
((s_{\frak d_0}^{(2)})^{-1}(0) \cap U_{\frak d_0}^{(2)}).
\end{equation}
\begin{lem}\label{lem1}
There exist $U_{\frak d_0}^{(3)} \subset U_{\frak d_0}^{(2)}$
and $U_{\frak d}^{(3)} \subset U_{\frak d}^{(2)}$ such that the following holds.
We define $U_{\frak d \frak d_0}^{(3)}$ and $U_{\frak d \frak d'}^{(3)}$ for $\frak d,\frak d' \in D'$
with $\frak d > \frak d' >\frak d_0$ in the same way as (\ref{44344}).
Various bundles maps sections are defined by restrictions in an obvious way.
Then we have:
\begin{enumerate}
\item
(\ref{U2cond1})-(\ref{U2cond4}) hold when we replace ${*}^{(1)}$ by ${*}^{(3)}$.
(Here and hereafter $*$ is anything such as $U_{\frak d_0}$ etc.)
\item
The conclusion of Lemma \ref{lem418} holds when we replace ${*}^{(1)}$ by ${*}^{(3)}$.
\item
(\ref{Liscontained2}) holds when we replace ${*}^{(2)}$ by ${*}^{(3)}$.
\item
\begin{equation}\label{overremove1}
\psi_{\frak d_0}(s_{\frak d_0}^{-1}(0)  \cap U_{\frak d_0}^{(3)})
\cap
\psi_{\frak d}(s_{\frak d}^{-1}(0)  \cap U^{(3)}_{\frak d})
\subseteq
\psi_{\frak d_0}
(s_{\frak d_0}^{-1}(0)  \cap U_{\frak d\frak d_0}^{(3)})
\end{equation}
for $\frak d > \frak d_0$
\end{enumerate}
\end{lem}
The opposite inclusion of \eqref{overremove1} is obvious.
We note that (4) above implies that Property (B) (injectivity required in Definition \ref{mixednbd} (9))
holds for $U_{\frak d_0}^{(3)}$, $U_{\frak d}^{(3)}$.
\begin{proof}
By (\ref{Liscontained2}) and (\ref{Ldefn}) we have
$$
\aligned
&\psi_{\frak d}^{(2)}(s_{\frak d}^{-1}(0)  \cap U^{(2)}_{\frak d})
\cap
\psi_{\frak d_0}^{(2)}(s_{\frak d_0}^{-1}(0)  \cap U^{(2)}_{\frak d_0})
\cap
X(D) \\
&=\psi_{\frak d}^{(2)}(s_{\frak d}^{-1}(0)  \cap U^{(2)}_{\frak d})
\cap
\psi_{\frak d_0}^{(2)}(s_{\frak d_0}^{-1}(0)  \cap U^{(2)}_{\frak d_0})
\cap
X(\frak d_0) \\
&\subset
\overline{\mathcal L}_{\frak d_0} \cap \psi^{(2)}_{\frak d}((s_{\frak d}^{(2)})^{-1}(0)\cap U_{\frak d}^{(2)})
\subset
\psi^{(2)}_{\frak d_0} (U_{\frak d\frak d_0}^{(2)} \cap (s^{(2)}_{\frak d_0})^{-1}(0)).
\endaligned$$
Therefore we can choose open sets $V_{\frak d_0}, V_{D'} \subset X$ such that:
\begin{enumerate}
\item $V_{D'} \supset X(D')$.
\item
 $V_{D'} \supset \overline{(X(\frak d_0) \setminus \mathcal K_{\frak d_0})}$.
\item $\mathcal K_{\frak d_0} \subset V_{\frak d_0}$.
\item
$V_{D'} \cup V_{\frak d_0} \supset X(D)$.
\item
$V_{D'} \cap V_{\frak d_0}  \cap \psi^{(2)}_{\frak d}((s_{\frak d}^{(2)})^{-1}(0)\cap U_{\frak d}^{(2)})\subset \psi^{(2)}_{\frak d_0} (U_{\frak d\frak d_0}^{(2)} \cap (s^{(2)}_{\frak d_0})^{-1}(0)) $.
\end{enumerate}
We take $U_{\frak d_0}^{(3)} \subset U_{\frak d_0}^{(2)}$
and $U(D')^{(3)} \subset U(D')^{(2)}$,
(where $U(D')^{(2)} = \bigcup_{\frak d \in D'} \Pi_{\frak d}(U_{\frak d}^{(2)})$)
such that
$$
\aligned
&U_{\frak d_0}^{(3)} \cap (s_{\frak d_0}^{(2)})^{-1}(0)
=(\psi_{\frak d_0}^{(2)})^{-1}( V_{\frak d_0}),\\
&
U_{D'}^{(3)} \cap (s_{D'}^{(2)})^{-1}(0)
=(\psi_{D'}^{(2)})^{-1}( V_{D'}).
\endaligned
$$
We use it to obtain $U_{\frak d}^{(3)}$.
Then Formula (\ref{Liscontained2}) after replacing ${*}^{(2)}$  by ${*}^{(3)}$
determines $U_{\frak d\frak d_0}^{(3)}$.
The open sets $U_{\frak d\frak d'}^{(3)}$ for $\frak d,\frak d' \in D'$
$(\frak d > \frak d')$ are defined similarly. The bundles, maps, sections are defined by restriction
in an obvious way.
Then the conclusion of Lemma \ref{lem418} holds.
\par
Condition (1)-(4) above imply (\ref{U2cond1})-(\ref{U2cond4}),
respectively.
Condition (5) implies (\ref{overremove1}).
\end{proof}
\begin{rem}
A similar formula
\begin{equation}\label{overremove2}
\psi_{\frak d_1}(s_{\frak d_1}^{-1}(0)  \cap U^{(3)}_{\frak d_1})
\cap
\psi_{\frak d_2}(s_{\frak d_2}^{-1}(0)  \cap U^{(3)}_{\frak d_2})
\subseteq
\psi_{\frak d_1}
(s_{\frak d_1}^{-1}(0)  \cap U^{(3)}_{\frak d_2\frak d_1})
\end{equation}
for $\frak d_2 > \frak d_1 > \frak d_0$
is a consequence of Definition \ref{mixednbd} (9), applied to $D'$
and so is a part of induction hypothesis.
(More precisely (\ref{overremove2}) with $U^{(3)}_{\frak d_1}$ etc.
replaced by $U^{(1)}_{\frak d_1}$ etc. is the consequence of the
induction hypothesis. Then (\ref{overremove2}) follows easily from definition.)
\end{rem}
It remains to shrink so that (A) (Hausdorff-ness) holds.
The way we shrink here is similar to the construction of pure orbifold
neighborhood.  Note we are gluing many spaces $U_{\frak d}$.
We will reduce the problem to the gluing of two spaces
$U_{\frak d_0}$ and $U(D')$.
For this purpose we need to modify so that the maps $\underline\phi_{\frak d\frak d_0}$
can be glued to give a map $U_{\frak d_0} \to U(D')$.
The detail follows.
\par
We shrink again each of $U^{(3)}_{\frak d}$ with $\frak d > \frak d_0$
to obtain $U^{(4)}_{\frak d}$.
We take it so that $U^{(4)}_{\frak d}$ is relatively compact in $U^{(3)}_{\frak d}$
and $U^{(4)}(D') = \bigcup_{\frak d\in D'}\Pi_{\frak d}(U^{(4)}_{\frak d})$
still carries the structure of mixed orbifold neighborhood.
We also shrink $U^{(3)}_{\frak d_0}$ to a relatively compact subset
$U^{(4)}_{\frak d_0}$.
The domain of the coordinate change is defined by:
$$
U^{(4)}_{\frak d_2\frak d_1}
= \Pi^{-1}_{\frak d_1}(\Pi_{\frak d_1}(U^{(4)}_{\frak d_1})\cap
\Pi_{\frak d_2}(U^{(4)}_{\frak d_2}))
=
\underline\phi_{\frak d_2\frak d_1}^{-1}(U^{(4)}_{\frak d_2})
\cap U^{(4)}_{\frak d_1}
$$
and
$$
U^{(4)}_{\frak d\frak d_0}
= (\underline\phi_{\frak d\frak d_0}^{(1)})^{-1}(U_{\frak d}^{(4)})
\cap U^{(4)}_{\frak d_0}.
$$
Bundles and bundle maps etc. are obtained by restriction in an obvious way.
The conclusion of Lemma \ref{lem1} holds with $*^{(3)}$
replaced by $*^{(4)}$.
\par
For each point $x \in \psi_{\frak d_0}^{-1}(\mathcal K_{\frak d_0})$
we take $U_x$ a neighborhood of $x$ in $U^{(4)}_{\frak d_0}$
with the following property.
\begin{proper}\label{propertyforUxX}
If $\frak d_2 > \frak d_1 > \frak d_0$ and
$$
U_x \cap U^{(4)}_{\frak d_1 \frak d_0} \cap
U^{(4)}_{\frak d_2 \frak d_0} \ne \emptyset
$$
then
$$
\underline\phi^{(3)}_{\frak d_1\frak d_0}(U_x)
\subseteq
U^{(3)}_{\frak d_2 \frak d_1}.
$$
\end{proper}
Such a choice is possible because of the next lemma.
\begin{lem}
$$
(\psi_{\frak d_0}^{(3)})^{-1}(\mathcal K_{\frak d_0})
\cap  \overline{(U^{(4)}_{\frak d_1 \frak d_0} \cap
U^{(4)}_{\frak d_2 \frak d_0})}
\subset
(\underline\phi^{(3)}_{\frak d_1\frak d_0})^{-1}
(U^{(3)}_{\frak d_2\frak d_1}).
$$
\end{lem}
\begin{proof}
Note that
$
\frak K_{\frak d_0} = \psi_{\frak d_0}^{(3)}((s_{\frak d_0}^{(3)})^{-1}(0)
\cap \overline{U}_{\frak d_0}^{(4)})
$
is a compact subset of
$\psi_{\frak d_0}^{(3)}((s_{\frak d_0}^{(3)})^{-1}(0)  \cap U_{\frak d_0}^{(3)})$
and
$
\frak K_{\frak d_i} =  \psi_{\frak d_i}^{(3)}((s_{\frak d_i}^{(3)})^{-1}(0)
 \cap \overline{U}_{\frak d_i}^{(4)})
$
is a compact subset of
$\psi_{\frak d_i}^{(3)}((s_{\frak d_i}^{(3)})^{-1}(0)  \cap U^{(3)}_{\frak d_i})$
for $i=1,2$.
(Here $\overline{U}_{\frak d_0}^{(4)}$ is the closure of
${U}_{\frak d_0}^{(4)}$ in ${U}_{\frak d_0}^{(3)}$.)
Using
(\ref{overremove1}), (\ref{overremove2}),
we find
$$
\aligned
\frak K_0 \cap \frak K_i
&\subset
(\psi_{\frak d_0}^{(3)})(s_{\frak d_0}^{-1}(0)  \cap U_{\frak d_0}^{(3)})
\cap
\psi_{\frak d_i}(s_{\frak d_i}^{-1}(0)  \cap U_{\frak d_i})\\
&\subseteq
(\psi_{\frak d_0}^{(3)})
(s_{\frak d_0}^{-1}(0)  \cap U_{\frak d_i\frak d_0}^{(3)}),
\endaligned
$$
and
$$
\frak K_1 \cap \frak K_2 \subset
\psi^{(3)}_{\frak d_1}(s_{\frak d_1}^{-1}(0)  \cap U^{(3)}_{\frak d_1})
\cap
\psi^{(3)}_{\frak d_2}(s_{\frak d_2}^{-1}(0)  \cap U^{(3)}_{\frak d_2})
\subseteq
\psi^{(3)}_{\frak d_1}
((s^{(3)}_{\frak d_1})^{-1}(0)  \cap U^{(3)}_{\frak d_2\frak d_1}
).
$$
Now the lemma follows from:
$$
\psi^{(3)}_{\frak d_0}
\left(
(\psi_{\frak d_0}^{(3)})^{-1}(\mathcal K_{\frak d_0})
\cap  \overline{(U^{(4)}_{\frak d_1 \frak d_0} \cap
U^{(4)}_{\frak d_2 \frak d_0})}
\right)
\subset
\frak K_0 \cap \frak K_1 \cap \frak K_2.
$$
\end{proof}
We cover $\psi_{\frak d_0}^{-1}(\mathcal K_{\frak d_0})$ by
a finitely many sets $U_{x_i}$ ($i=1,\dots,I$) among such
$U_x$'s and let $U^{(5)}_{\frak d_0}$ be the  union
$\bigcup_{i\in I}U_{x_i}$ of them.
\begin{lem}
\begin{equation}\label{Mcformula-}
U_{\frak d_1 \frak d_0}^{(5)} \cap  U_{\frak d_2 \frak d_0}^{(5)}
\subset
(\underline\phi^{(5)}_{\frak d_1 \frak d_0})^{-1}(U_{\frak d_2 \frak d_1}^{(4)}).
\end{equation}
Here
$$
U_{\frak d \frak d_0}^{(5)}
= U_{\frak d \frak d_0}^{(4)} \cap U^{(5)}_{\frak d_0}
$$
and ${\underline\phi}^{(5)}_{\frak d \frak d_0}$ is the
restriction of ${\underline\phi}^{(3)}_{\frak d \frak d_0}$
to $U_{\frak d \frak d_0}^{(5)}$.
\end{lem}
\begin{proof}
Suppose $y \in U_{\frak d_1 \frak d_0}^{(5)} \cap  U_{\frak d_2 \frak d_0}^{(5)}$.
Then we have
$
{\underline\phi}^{(3)}_{\frak d_2\frak d_0}(y) \in U_{\frak d_2}^{(4)}
$
and
$
{\underline\phi}^{(3)}_{\frak d_1\frak d_0}(y)\in U_{\frak d_1}^{(4)}.
$
There exists $x_i$ ($i\in I$) with $y \in U_{x_i}$.
Since
$$
y \in U_{x_i} \cap U^{(5)}_{\frak d_1 \frak d_0} \cap
U^{(5)}_{\frak d_2 \frak d_0},
$$
Property \ref{propertyforUxX} implies:
$$
{\underline\phi}_{\frak d_1\frak d_0}^{(3)}(y)\in
{\underline\phi}_{\frak d_1\frak d_0}^{(3)}(U_{x_i})
\subseteq
U_{\frak d_2 \frak d_1}^{(3)}.
$$
By Lemma \ref{lem418},
$$
{\underline\phi}_{\frak d_2\frak d_1}{\underline\phi}_{\frak d_1\frak d_0}^{(3)}(y)
= {\underline\phi}_{\frak d_2\frak d_0}^{(3)}(y).
$$
Therefore
$$
{\underline\phi}^{(3)}_{\frak d_1\frak d_0}(y)
\in ({\underline\phi}^{(1)}_{\frak d_2 \frak d_1})^{-1}(U_{\frak d_2}^{(4)})\cap U_{\frak d_1}^{(4)}
=
U_{\frak d_2 \frak d_1}^{(4)}.
$$
Thus
$$
y \in
({\underline\phi}^{(5)}_{\frak d_1 \frak d_0})^{-1}(U_{\frak d_2 \frak d_1}^{(4)})
$$
as required.
\end{proof}
Let $U^{(6)}_{\frak d}$ be a relatively compact
subset of  $U^{(4)}_{\frak d}$ for $\frak d \in D'$ and
$U^{(6)}_{\frak d_0}$ a relatively compact
subset of  $U^{(5)}_{\frak d_0}$.
We may choose them so that Lemma \ref{lem1} (1)-(4)
holds when we replace $*^{(3)}$ and $*^{(2)}$ by $*^{(6)}$.
\par
We define
\begin{equation}
U^{(6)}(D)
= \bigcup_{\frak d \in D} U^{(6)}_{\frak d}/\sim.
\end{equation}
Here $\sim$ is defined as follows.
$x \sim y$ if and only if one of the following holds:
\begin{enumerate}
\item $x = y$.
\item $x \in U^{(6)}_{\frak d'}$, $y \in U^{(6)}_{\frak d'\frak d} \cap
U^{(6)}_{\frak d}$, $x = \underline{\phi}^{(6)}_{\frak d'\frak d}(y)$.
\item $y \in U^{(6)}_{\frak d'}$, $x \in U^{(6)}_{\frak d'\frak d} \cap
U^{(6)}_{\frak d}$, $y = \underline{\phi}^{(6)}_{\frak d'\frak d}(x)$.
\end{enumerate}
We define $U^{(6)}(D')$ in the same way.
\par\medskip
We define a map ${\Pi}_{\frak d} : U^{(6)}_{\frak d} \to U^{(6)}(D)$ by sending an element to its equivalence class.
We remark that we have a continuous map
$$
\underline{\phi}^{(6)}_{D'\frak d_0} : U^{(6)}_{D'\frak d_0} \to U^{(6)}(D')
$$
from
$
U^{(6)}_{D'\frak d_0} = \bigcup_{\frak d \in D'} {\Pi}^{(6)}_{\frak d}(U^{(6)}_{\frak d\frak d_0})
$
such that
$$
\underline{\phi}^{(6)}_{D'\frak d_0} = {\Pi}^{(6)}_{\frak d}\circ \underline{\phi}^{(6)}_{\frak d\frak d_0}
$$
holds on $U^{(6)}_{\frak d\frak d_0}$.
This is a consequence of Lemma \ref{lem418} (3)
and \eqref{Mcformula-}.
\begin{rem}\label{thankmac}
The authors thank D. McDuff who pointed out that the inclusion (\ref{Mcformula-})
is necessary to show such  $U^{(6)}_{D'\frak d_0}$ exists, during
our discussion at google group Kuranishi.
\end{rem}
\par
Now we are in the last step to achieve Hausdorff-ness.
We use a similar trick as in the last part of the proof of Proposition \ref{purecover} to
modify $U^{(6)}_{\frak d}$ etc. as follows.
Note that $U^{(6)}(D)$ can also be  written as
$$
U^{(6)}(D)
= (U^{(6)}(D') \cup U^{(6)}_{\frak d_0})/\sim,
$$
where $x \sim y$ if and only if one of the following holds.
\begin{enumerate}
\item $x = y$, or
\item $x \in U^{(6)}_{D'}$, $y \in U^{(6)}_{D'\frak d_0} \subset U^{(6)}_{\frak d_0}$,
$x = \underline{\phi}^{(6)}_{D'\frak d_0}(y)$, or
\item $y \in U^{(6)}_{D'}$, $x \in U^{(6)}_{D'\frak d_0} \subset U^{(6)}_{\frak d_0}$, $y
= \underline{\phi}^{(6)}_{D'\frak d_0}(x)$.
\end{enumerate}
\par\medskip
We also note that $U^{(6)}(D')$ is already Hausdorff (with respect to the quotient topology) by induction hypothesis.
\begin{rem}
In fact the obvious map $U^{(6)}(D') \to U(D')$ is injective and continuous.
(It may not be a topological embedding however. See Remark \ref{rem520}.)
\end{rem}
Now the rest of the construction is similar to the one in Lemmas \ref{lemma49} - \ref{XcapU}.
We take a relatively compact subset $U^{(7)}(D')$ of $U^{(6)}(D')$ such that
\begin{equation}\label{formula440}
U^{(7)}(D') \supset X(D')
\end{equation}
and a relatively compact subset $U^{(7)}_{\frak d_0}$ of $U^{(6)}_{\frak d_0}$ such that
\begin{equation}\label{formula441}
U^{(7)}_{\frak d_0} \supset (\psi^{(6)}_{\frak d_0})^{-1}(\mathcal K_{\frak d_0}).
\end{equation}
We take $W_{D'\frak d_0} \subset U^{(6)}_{D'\frak d_0}$ such that
$$
W_{D'\frak d_0} \supset (s_{\frak d_0}^{(6)})^{-1}(0) \cap \overline{(U^{(7)}_{\frak d_0} \cap
(\underline{\phi}^{(6)}_{D\frak d_0})^{-1}(U^{(7)}(D')))}.
$$
The existence of such $W_{D'\frak d_0}$ can be proved in the same way as Lemma \ref{lemma49}.
In the same way as Lemma \ref{lem410}, we can find $U^{(8)}_{\frak d_0}$ such that
together with $U^{(8)}(D') = U^{(7)}(D')$ it satisfies
\begin{equation}
\overline U^{(8)}_{\frak d_0} \cap \overline{(\underline{\phi}^{(6)}_{D'\frak d_0})^{-1}(U^{(8)}(D'))}
\subset W_{D'\frak d_0}.
\end{equation}
Moreover (\ref{formula440}), (\ref{formula441}) hold with $*^{(7)}$ replaced by $*^{(8)}$.
\par
We put $C_1 = \overline U^{(8)}_{\frak d_0}$, $C_2 = \overline U^{(8)}(D')$.
We put $C_{21} = (\underline{\phi}^{(6)}_{D'\frak d_0})^{-1}(C_2)$. We then define
$C = (C_1 \cup C_2)/\sim$ where $\sim$ is defined by using the restriction of $\underline{\phi}^{(6)}_{D'\frak d_0}$
to $C_{21}$, as before. (We put quotient topology on it.)
We can prove that $C_{21}$ is compact in the same way as Lemma \ref{C21compact}.
Therefore $C$ is Hausdorff by Proposition \ref{metrizable}.
\par
Let $\Pi_i : C_i \to C$ be the obvious map.
We put
$$
U(D) = \Pi_1(U^{(8)}_{\frak d_0}) \cup \Pi_2(U^{(8)}(D')) \subset C
$$
and will use a topology induced from the topology of $C$ on it.
We define
$$
U_{\frak d_0} = U^{(8)}_{\frak d_0}, \qquad
U_{\frak d} = \Pi_{\frak d}^{-1}(U_\frak d^{(8)}) \subset U^{(6)}_{\frak d}.
$$
They are orbifolds.
We obtain bundles, sections, maps, coordinate changes, on them by restriction in an obvious way.

The proof of Proposition \ref{existmixed} is complete.
\end{proof}
\begin{lem}\label{linearcond}
We may choose $U(D)$ so that the following holds in addition.
Let $\frak d_k > \frak d_0$.
If
$$
\bigcap_{k=1}^K
{\Pi}_{\frak d_k}(U_{\frak d_k}) \cap {\Pi}_{\frak d_0}(U_{\frak d_0})  \ne \emptyset
$$
then
$$
\bigcap_{k=1}^K
{\Pi}_{\frak d_k}(U_{\frak d_k}\cap s_{\frak d_k}^{-1}(0)) \cap {\Pi}_{\frak d_0}(U_{\frak d_0}\cap s_{\frak d_0}^{-1}(0))  \ne \emptyset.
$$
\end{lem}
\begin{proof}
We will modify $U(D)$ so that it satisfies this additional condition by induction on $\#D$.
\par
The inductive step is as follows.
We take $\frak d_0 \in D$ that is minimal in $D$.
We put $D' = D \setminus \{\frak d_0\}$.
\par
We modify $U_{\frak d_0}$ so that the conclusion of the lemma holds by
induction on $K$. We assume the conclusion of the lemma holds for
$K \le K_0-1$.
We consider the case of $K_0$.
Let
$$
\frak C = \{\{\frak d_1,\dots,\frak d_{K_0}\} \mid (\ref{sempty}), \text{$\frak d_i$ are all different.}\}
$$
\begin{equation}\label{sempty}
\bigcap_{k=1}^{{K_0}}
{\Pi}_{\frak d_k}(U_{\frak d_k}\cap s_{\frak d_k}^{-1}(0)) \cap {\Pi}_{\frak d_0}(U_{\frak d_0}\cap s_{\frak d_0}^{-1}(0)) =\emptyset.
\end{equation}
We shrink $U_{\frak d}$ a bit so that we may assume
\begin{equation}\label{semptycl}
\bigcap_{k=1}^{{K_0}}
\overline{{\Pi}_{\frak d_k}(U_{\frak d_k}\cap s_{\frak d_k}^{-1}(0)}) \cap \overline{{\Pi}_{\frak d_0}(U_{\frak d_0}\cap s_{\frak d_0}^{-1}(0)}) =\emptyset.
\end{equation}
for $\{\frak d_1,\dots,\frak d_{K_0}\} \in \frak C$.
\par
We replace $U_{\frak d_0}$ by
\begin{equation}\label{Ud_0prime}
U'_{\frak d_0} = U_{\frak d_0} \setminus \bigcup_{\{\frak d_1,\dots,\frak d_{K_0}\} \in \frak C} \bigcap_{k=1}^{K_0}\overline U_{\frak d_k\frak d_0}.
\end{equation}
We will prove $U_{\frak d}$ ($\frak d > \frak d_0$)
together with $U'_{\frak d_0}$ satisfies the required property for $K \le K_0$.
\par
We first consider the case $K\le K_0 -1$.
Suppose
$$
\bigcap_{k=1}^K
{\Pi}_{\frak d_k}(U_{\frak d_k}) \cap {\Pi}_{\frak d_0}(U'_{\frak d_0})  \ne \emptyset.
$$
Then
$$
\bigcap_{k=1}^K
{\Pi}_{\frak d_k}(U_{\frak d_k}) \cap {\Pi}_{\frak d_0}(U_{\frak d_0})  \ne \emptyset.
$$
Then by induction hypothesis we have
$$
\bigcap_{k=1}^K
{\Pi}_{\frak d_k}(U_{\frak d_k}\cap s_{\frak d_k}^{-1}(0)) \cap {\Pi}_{\frak d_0}(U_{\frak d_0}\cap s_{\frak d_0}^{-1}(0))  \ne \emptyset.
$$
We note that
\begin{equation}\label{UprimeandUons-0}
U'_{\frak d_0} \cap s^{-1}_{\frak d_0}(0) =U_{\frak d_0} \cap s^{-1}_{\frak d_0}(0)
\end{equation}
by (\ref{semptycl}), (\ref{Ud_0prime}). Therefore
$$
\bigcap_{k=1}^K
{\Pi}_{\frak d_k}(U_{\frak d_k}\cap s_{\frak d_k}^{-1}(0)) \cap {\Pi}_{\frak d_0}(U'_{\frak d_0}\cap s_{\frak d_0}^{-1}(0))  \ne \emptyset
$$
as required.
\par
We next consider the case $K=K_0$.
Suppose
$$
\bigcap_{k=1}^{K_0}
{\Pi}_{\frak d_k}(U_{\frak d_k}\cap s_{\frak d_k}^{-1}(0)) \cap {\Pi}_{\frak d_0}(U'_{\frak d_0}\cap s_{\frak d_0}^{-1}(0))  = \emptyset.
$$
Then by (\ref{UprimeandUons-0}) we have
$$
\bigcap_{k=1}^{K_0}
{\Pi}_{\frak d_k}(U_{\frak d_k}\cap s_{\frak d_k}^{-1}(0)) \cap {\Pi}_{\frak d_0}(U_{\frak d_0}\cap s_{\frak d_0}^{-1}(0))  = \emptyset.
$$
Namely $\{ \frak d_1,\dots,\frak d_{K_0}\} \in C$.
Therefore
$$
\bigcap_{k=1}^{K_0}
{\Pi}_{\frak d_k}(U_{\frak d_k}) \cap {\Pi}_{\frak d_0}(U'_{\frak d_0})  = \emptyset
$$
by (\ref{Ud_0prime}). The proof of the inductive step is complete.
\end{proof}

Now we are in the position to complete the proof of Theorem \ref{goodcoordinateexists}.
We apply Proposition \ref{existmixed}  to obtain a mixed orbifold neighborhood of $X(\frak D) = X$.
We put $\frak P = \frak D \subset \Z_{> 0}$.
The order is $\le$.
For $\frak d \in \frak D = \frak P$, we have
$U_{\frak d}$, $E_{\frak d}$, $s_{\frak d}$, $\psi_{\frak d}$ by Definition \ref{mixednbd} (2)(3).
Let us check Definition \ref{goodcoordinatesystem} (1)-(9).
\par
 Definition \ref{goodcoordinatesystem} (1)-(4) follows from Definition \ref{mixednbd} (2)(3).
Definition \ref{goodcoordinatesystem} (5)(6) follows from  Definition \ref{mixednbd} (4) - (8).
Definition \ref{goodcoordinatesystem} (7) follows from Definition \ref{mixednbd} (9).
Definition \ref{goodcoordinatesystem} (8)  is obvious since $\frak P \subset \Z_{> 0}$.
Condition  \ref{Joyce} in Definition \ref{goodcoordinatesystem} (9) follows from Definition \ref{mixednbd} (8) (9).
Conditions \ref{plusalpha} and  \ref{plusalpha2} in Definition \ref{goodcoordinatesystem} (9) follow from Lemma \ref{linearcond}
and (\ref{Dddd}).
Condition \ref{proper} follows from Definition \ref{mixednbd} (8)
especially Hausdorff-ness of $U(X)$.
\par
The proof of Theorem \ref{goodcoordinateexists} is now complete.
\end{proof}

\section{Appendix: a lemma on general topology}\label{gentoplem}

In this appendix we prove Proposition \ref{metrizable}.
We assume Assumption \ref{Kassumption}.

We first prove the following.
\begin{lem}\label{hausdorfflema}
$K(\frak P)$ is Hausdorff.
\end{lem}
\begin{rem}
Note $\{(x,y) \in (\coprod_{\frak p} K_{\frak p}) \times (\coprod_{\frak p} K_{\frak p}) \mid x\sim y\}$
is a closed subset. However this does not immediately imply that $K(\frak P)$ is Hausdorff.
(This is because $\Pi_{\frak p}$ is not an open mapping.)
\end{rem}
\begin{proof}
The proof is by induction on $\#\frak P$. Denote by $\pi: \coprod_{\frak p} K_{\frak p} \to K(\frak P)$ the projection.
\par
The case $\#\frak P = 1$ is trivial.
Suppose $\frak P = \{\frak p,\frak q\}$ and $\frak p>\frak q$.
We remark that $\Pi_{\frak q} : K_{\frak q} \to K(\frak P)$ and $\Pi_{\frak p} : K_{\frak p} \to K(\frak P)$
are both closed mappings.
\par
Let $p \ne q \in K(\frak P)$. We need to find open neighborhoods $A_p, \, A_q \subset K(\frak P)$
such that $A_p \cap A_q =\emptyset$.
There are 4 cases to consider.
We put $q = [x]$, $p =[y]$
\begin{enumerate}
\item $x, \, y\in K_{\frak q} \setminus K_{\frak p\frak q}$ or
$x, \, y\in K_{\frak p} \setminus \underline{\phi}_{\frak p\frak q}(K_{\frak p\frak q})$.
\item $x \in K_{\frak q}\setminus K_{\frak p\frak q}$,
$y\in K_{\frak p} \setminus \underline{\phi}_{\frak p\frak q}(K_{\frak p\frak q})$. Or the same with $x$ and $y$ exchanged.
\item $x \in K_{\frak q} \setminus K_{\frak p\frak q}$, $y \in
K_{\frak p\frak q}$. There are 3 similar cases where $x$, $y$ are exchanged and/or $\frak p$, $\frak q$  are
exchanged.
\item $x, \, y\in K_{\frak p\frak q}$.
\end{enumerate}
Case (1).
Suppose $x, \, y\in K_{\frak q} \setminus K_{\frak p\frak q}$.
Choose disjoint open subsets $A'_x, A'_y$ of $K_{\frak q}$ such that $x \in A'_x, y \in A'_y$.
Then
$A_x = \Pi_{\frak q}(A'_x \setminus K_{\frak p\frak q})$ and
$A_y = \Pi_{\frak q}(A'_y \setminus K_{\frak p\frak q})$ have required properties.
\par
\noindent
Case (2).
$A_x = \Pi_{\frak q}(K_{\frak q} \setminus K_{\frak p\frak q})$
and $A_y = \Pi_{\frak p}(K_{\frak p} \setminus \underline{\phi}_{\frak p\frak q}(K_{\frak p\frak q}))$ have required properties.
\par
\noindent
Case (3).
Note $K_{\frak q}$ is normal since it is Hausdorff and compact.
Therefore there exist disjoint open subsets $A'_x, A'_y$ of $K_{\frak q}$ such that
$x \in A'_x, \,\, K_{\frak p\frak q} \subset A'_y$.
Let $A''_y$ be an open subset of $K_{\frak p}$ containing $\underline{\phi}_{\frak p\frak q}(K_{\frak p\frak q})$.
Then $A_x = \Pi_{\frak q}(A'_x)$,
$A_y = \Pi_{\frak q}(A'_y) \cup \Pi_{\frak p}(A''_y)$ have required properties.
\par
\noindent
Case (4).
We take disjoint open subsets $A'_x$ and $A'_y$ of $K_{\frak q}$ such that
$x \in  A'_x$ and $y \in  A'_y$.
Since $K_{\frak q}$ is normal we may assume $\overline A'_x \cap \overline A'_y = \emptyset$.
\par
We take open subsets $A''_x$ and $A''_y$ of $K_{\frak p}$ such that
$A''_x \cap \underline{\phi}_{\frak p\frak q}(K_{\frak p\frak q}) = \underline{\phi}_{\frak p\frak q}(A'_x \cap K_{\frak p\frak q})$
and
$A''_y \cap \underline{\phi}_{\frak p\frak q}(K_{\frak p\frak q}) = \underline{\phi}_{\frak p\frak q}(A'_y \cap K_{\frak p\frak q})$.
Such $A''_x$ exists since  $\underline{\phi}_{\frak p\frak q}(A'_x \cap K_{\frak p\frak q})$
is an open subset of $\underline{\phi}_{\frak p\frak q}(K_{\frak p\frak q})$ with respect to the subspace topology.
However $A''_x \cap A''_y$ may be nonempty.
\par
Since $K_{\frak p}$ is normal and $\overline A'_x \cap \overline A'_y = \emptyset$,
we can find disjoint open subsets $A'''_x$ and $A'''_y$ of $K_{\frak p}$
such that
$\underline{\phi}_{\frak p\frak q}(\overline A'_x \cap K_{\frak p\frak q}) \subset A'''_x$,
$\underline{\phi}_{\frak p\frak q}(\overline A'_y \cap K_{\frak p\frak q}) \subset A'''_y$.
\par
Now
$
A_x = \Pi_{\frak q}(A'_x) \cup \Pi_{\frak p}(A''_x\cap A'''_x)
$,
$
A_y = \Pi_{\frak q}(A'_y) \cup \Pi_{\frak p}(A''_y\cap A'''_y)
$
have required properties.
\par
This proves the lemma for the case $\# \frak P = 2$.
\par\medskip
Now suppose the lemma hold for the case $\#\frak P < n$ for $n \in \N$
and prove the case of $\#\frak P  = n$.
Let $\frak p_0 \in \frak P$ be an element which is minimal
with respect to the partial order.
We put $\frak P' = \frak P \setminus \{\frak p_0\}$.
We obtain $K(\frak P')$. By the induction hypothesis, it is Hausdorff.
We put
\begin{equation}\label{closedcover}
K_{\frak P'\frak p_0} = \bigcup_{\frak p \in \frak P'} K_{\frak p\frak p_0} \subset K_{\frak p_0}.
\end{equation}
Let $\Pi'_{\frak p} : K_{\frak p} \to K(\frak P')$
($\frak p \in \frak P'$) be the map which sends an element to its equivalence class.
For $x \in K_{\frak p\frak p_0}$ we put
\begin{equation}\label{kiregirevarsho}
\underline\phi_{\frak P'\frak p_0}(x) = \Pi'_{\frak p}(\underline\phi_{\frak p\frak p_0}(x)).
\end{equation}
It is easy to see that (\ref{kiregirevarsho}) induces a map
\begin{equation}
\underline\phi_{\frak P'\frak p_0} : K_{\frak P'\frak p_0} \to K(\frak P').
\end{equation}
\begin{sublem}
$\underline\phi_{\frak P'\frak p_0}$ is continuous.
\end{sublem}
\begin{proof}
The restriction of $\underline\phi_{\frak P'\frak p_0}$ to each $K_{\frak p\frak p_0}$
is continuous by the definition of quotient topology.
Moreover (\ref{closedcover}) is a covering by closed sets.
The sublemma follows.
\end{proof}
Since $K_{\frak P'\frak p_0}$ is compact and $K(\frak P')$ is Hausdorff
by induction hypothesis, it follows that $\underline\phi_{\frak P'\frak p_0}$
is a topological embedding.
Thus $K_{\frak p_0}$, $K(\frak P')$, $K_{\frak P'\frak p_0}$, $\underline\phi_{\frak P'\frak p_0}$
satisfy the assumption of Lemma \ref{hausdorfflema}, where $\frak P = \{\frak P',\frak p_0\}$.
($\frak p_0 < \frak P'$.)
We then obtain a Hausdorff space, which we write $K'(\frak P)$.
The proof of Lemma \ref{hausdorfflema} is completed by the next sublemma.
\end{proof}
\begin{sublem}
$K'(\frak P)$ is homeomorphic to $K(\frak P)$.
\end{sublem}
\begin{proof} We have obvious projection maps
\begin{equation}\label{quotient1}
\pi'' : K_{\frak p_0}  \sqcup K(\frak P') \to K'(\frak P),
\end{equation}
\begin{equation}\label{quotient2}
\pi' : \coprod_{\frak p \in \frak P'} K_{\frak p} \to K(\frak P'),
\end{equation}
\begin{equation}\label{quotient3}
\pi : \coprod_{\frak p \in \frak P} K_{\frak p} \to K(\frak P).
\end{equation}
The topology of the right hand sides are quotient topology
with respect to these maps.
Composing (\ref{quotient1}) and (\ref{quotient2}) we obtain
\begin{equation}
\pi'' \circ (\text{id} \sqcup \pi' ) : \coprod_{\frak p \in \frak P} K_{\frak p} \to K'(\frak P).
\end{equation}
The topology of $K'(\frak P)$ is the quotient topology of this map.
\par
Using the fact that $\sim$ is an equivalence relation, it is easy to see
the following.
If $x,y \in \coprod_{\frak p \in \frak P} K_{\frak p}$, then
$\pi(x) = \pi(y)$ if and only if $\pi'' \circ (\text{id} \sqcup \pi' )(x)
= \pi'' \circ (\text{id} \sqcup \pi' )(y)$.
This implies the sublemma.
\end{proof}

We now recall that a family of subsets $\{U_i \mid i \in I\}$ of a topological space $X$ containing $x \in X$
is said to be a \emph{neighborhood basis} of $x$ if
\begin{enumerate}
\item each $U_i$ contains an open neighborhood of $x$,
\item for each open set $U$ containing $x$ there exists $i$ such that $U_i \subset U$.
\end{enumerate}
A family of open subsets $\{U_i \mid i \in I\}$ of a topological space $X$ is said to be a
basis of the open sets if for each $x$ the set $\{U_i \mid x \in U_i\}$
is a neighborhood basis of $x$. A topological space is said to satisfy the second axiom of countability
if there exists a countable basis of open subsets $\{U_i \mid i \in I\}$.

We next prove the following.
\begin{lem}\label{lem2kasan}
If $K_{\frak p}$ satisfies the second axiom of countability in addition, then
$K(\frak P)$ also satisfies the second axiom of countability.
\end{lem}
\begin{proof}
For each $\frak p$, we take a countable set $\frak U_{\frak p}
=\{U_{\frak p,i} \subset K_{\frak p} \mid i \in I_{\frak p}\}$
which is a basis of the open sets of $K_{\frak p}$. We may assume $\emptyset \in \frak U_{\frak p}$
and each $U_{\frak p,i}$ is open.
\par
For each $\vec i = (i_{\frak p})_{\frak p \in \frak P}$ ($i_{\frak p} \in I_{\frak p}$) we define
$U(\vec i)$ to be the interior of the set
\begin{equation}\label{U+basis1}
U^+(\vec i) := \bigcup_{\frak p \in \frak P}\Pi_{\frak p}(U_{\frak p,i_{\frak p}}).
\end{equation}
This is a countable family of open subsets of $K(\frak P)$.
We will prove that this family is a basis of open sets of $K(\frak P)$.
\par
Let $q \in K(\frak P)$, we put
\begin{equation}\label{defPxsss}
\frak P(q)
=
\{\frak p \in \frak P \mid q=[x], \, x \in K_{\frak p}\}.
\end{equation}
Here and hereafter we identify $K_{\frak p}$ to its image in $K(\frak P)$.
Note since $K(\frak P)$ is Hausdorff and $K_{\frak p}$ is compact,
the natural inclusion map $K_{\frak p} \to \coprod_{\frak p \in\frak P} K_{\frak p}$
induces a topological embedding $K_{\frak p} \to K(\frak P)$.
\par
For $\frak p \in \frak P(x)$, we have $x_{\frak p} \in K_{\frak p}$ with $[x_\frak p] = q$.
We take a countable neighborhood basis of $x_{\frak p}$.
We put $$
I_{\frak p}(q) = \{i \in I_{\frak p} \mid q = [x], \text{for some } x \in U_{\frak p,i}\}
$$
For each $\vec i = (i_{\frak p}) \in \prod_{\frak p\in \frak P(q)} I_{\frak p}(q)$,
we set
\begin{equation}\label{U+basis12}
U^+(\vec i) =\bigcup_{\frak p \in \frak P(q)} \Pi_{\frak p}(U_{\frak p,i_{\frak p}})
\subset K(\frak P).
\end{equation}
We claim that the collection $\{U^+(\vec i)  \mid \vec i \in \prod_{\frak p\in \frak P(q)} I_{\frak p}(q)\}$
is a neighborhood basis of $q$ in $K(\frak P)$ for any $q$.
The claim follows from Sublemmas \ref{9sublem1},\, \ref{9sublem2}.
\par
\begin{sublem}\label{9sublem1} The subset
$U^+(\vec i)$ is a neighborhood of $q$ in $K(\frak P)$.
\end{sublem}
\begin{proof}
For $\frak p \in \frak P(q)$ the set $K_{\frak p} \setminus U_{\frak p,i_{\frak p}}$
is a closed subset of $K_{\frak p}$ and so is compact.
Therefore $\Pi_{\frak p}(K_{\frak p} \setminus U_{\frak p,i_{\frak p}})$ is compact and so is closed.
\par
If $\frak p \notin \frak P(q)$ then we consider  $\Pi_{\frak p}(K_{\frak p})$ which is closed.
\par
Now we put
$$
K = \bigcup_{\frak p \in \frak P(q)}\Pi_{\frak p}(K_{\frak p} \setminus U_{\frak p,i_{\frak p}})
\cup
\bigcup_{\frak p \notin \frak P(q)}\Pi_{\frak p}(K_{\frak p}).
$$
This is a finite union of closed sets and so is closed.
It is easy to see that
$$
q \in K(\frak P) \setminus K \subset U^+(\vec i).
$$
\end{proof}
\begin{sublem} \label{9sublem2} The collection $\{U^+(\vec i)\}$ satisfies
the properties (2)  of the neighborhood basis above.
\end{sublem}
\begin{proof}
We remark that the map $K_{\frak p} \to K(\frak P)$ is a topological embedding.
Therefore $U\cap K_{\frak p}$ is an open set of $K_{\frak p}$.
Therefore for each $\frak p \in \frak P(q)$,
the set  $U\cap K_{\frak p}$ is a neighborhood of $x_{\frak p}$
in $K_{\frak p}$. By the definition of neighborhood basis in $K_\frak p$,
there exists $i_{\frak p}$ such that $U_{\frak p,i_{\frak p}} \subset U\cap K_{\frak p}$.
We put $\vec i = (i_{\frak p})$. Then
$U^+(\vec i) \subset U$ as required.
\end{proof}
We remark that $U^+(\vec i)$ in (\ref{U+basis12}) is a special case of $U^+(\vec i)$ in
(\ref{U+basis1}).
(We take $U_{\frak p,i_{\frak p}} = \emptyset$ for $\frak p \notin \frak P(x)$.)
The family $U(\vec i)$ is a countable basis of open sets of $K(\frak P)$.
The lemma is proved.
\end{proof}
\begin{lem}\label{loccomp}
If each $K_{\frak p}$ is locally compact, then
$K(\frak P)$ is locally compact.
\end{lem}
\begin{proof}
Let $x \in K(\frak P)$. We define $\frak P(x)$ by (\ref{defPxsss}).
For each $\frak p \in \frak P(x)$, we take a neighborhood basis of
$\{U_{\frak p,i} \mid I_{\frak p}\}$ of $x$ in $K_{\frak p}$ such that
$U_{\frak p,i}$ are all compact.
\par
For each $\vec i = (i_{\frak p}) \in \prod_{\frak p\in \frak P(x)} I_{\frak p}(x)$
we define $U^+(\vec i)$ by (\ref{U+basis12}).
They form a  neighborhood basis of $x$ in $K(\frak P)$ by Sublemmas \ref{9sublem1},
\ref{9sublem2}.
Since $U^+(\vec i)$ are all compact, the lemma follows.
\end{proof}
Combining Lemmas \ref{hausdorfflema}, \ref{lem2kasan}, \ref{loccomp}
and  a celebrated result by Urysohn
we obtain Proposition \ref{metrizable}.

\section{Orbifold via coordinate system}
\label{ofd}

In this section we review orbifold, its embedding and a bundle on it.
Let $X$ be a paracompact Hausdorff space.

\begin{defn}[orbifold]\label{ofddefn}
\begin{enumerate}
\item
An orbifold chart of $X$ at $p\in X$ is $(V_p,\Gamma_p,\psi_p)$
such that $V_p$ is a manifold on which a finite
group $\Gamma_p$ acts effectively, such that $o_p \in V_p$
is fixed by all the elements of $\Gamma_p$.
$\psi_p$ is a homeomorphism from a quotient space
$U_p = V_p/\Gamma_p$ onto a neighborhood of $p$ in $X$
such that $\psi_p(o_p) = p$.
\item
Let $(V_p,\Gamma_p,\psi_p)$ be as above and $q \in \psi_p(U_p)$.
A coordinate change is $(V_{pq},\phi_{pq},h_{pq})$ such that
$h_{pq} : \Gamma_q \to \Gamma_p$ is a group homomorphism,
$V_{pq} \subseteq V_q$ is a $\Gamma_q$ invariant open neighborhood of $o_q$ in
$V_q$ and $\phi_{pq} : V_{pq} \to V_p$ is an $h_{pq}$ equivariant
smooth open embedding of manifolds. We assume that
they satisfy
\begin{equation}\label{91form}
\psi_p \circ \underline\phi_{pq} = \psi_q
\end{equation}
where $\underline\phi_{pq}  : V_{pq}/\Gamma_q \to V_{p}/\Gamma_p$ is induced by $\phi_{pq}$.
\end{enumerate}
We call
$(\{ (V_p,\Gamma_p,\psi_p) \},\{(V_{pq},\phi_{pq},h_{pq})\})$
an {\it orbifold structure} on $X$.
\end{defn}
We can use (\ref{91form}) to show
$
\underline{\phi}_{pq} \circ \underline{\phi}_{qr}
= \underline{\phi}_{pr}
$
on
$
\underline{\phi}_{qr}^{-1}(U_{pq}) \cap U_{pr}.
$
We can use this fact to show that Definition \ref{ofddefn}
is equivalent to Definition \ref{ofddef1}.

\begin{rem}
We assumed the action of $\Gamma_p$ is effective.
We always do so in this article.
To emphasize this point we say sometimes
effective orbifold instead of orbifold.
\par
Sometimes effectivity of $\Gamma_p$-action is not assumed
in the definition of orbifold.
In such case definition of (uneffective) orbifold
becomes rather complicated if we use coordinate chart.
See \cite[Section 4]{fooo:overZ}.
The notion of morphisms between uneffective
orbifolds is also harder to define if we use
the language of coordinate chart.
(See Example \ref{notembedding} below.)
\end{rem}

\begin{defn}[embedding of orbifold]\label{emborf}
Let $X,Y$ be orbifolds and $F : X \to Y$ a continuous map.
$F$ is said to be an {\it embedding} of an orbifold if $F$ is a
topological embedding
and the following conditions are satisfied for each $q \in X$ and $p=F(q) \in Y$.
\begin{enumerate}
\item There exists an open subset $V^X_{pq} \subset V^X_q$ of the
chart of $q$ that is $\Gamma^X_q$ equivariant and containing $o_q$.
\item
There exists a smooth embedding
$F_{q} : V^X_{pq} \to V^Y_p$  of manifolds.
\item There exist a  group isomorphism $h^F_{pq} : \Gamma_q^X \to \Gamma^Y_p$ such that $F_{q}$ is $h^F_{pq}$ equivariant.
\item
The map $h^F_{pq}$ restricts to an isomorphism
$(\Gamma^X_q)_x \to (\Gamma^Y_p)_{\phi_{pq}(x)}$ for any
$x \in V^X_{pq}$.
\item $F_{q}$ induces the map $F\vert_{\psi_q(V^X_{pq}/\Gamma^X_q)} : \psi_q(V^X_{pq}/\Gamma^X_q)\to
\psi_p(V^Y_p/\Gamma^Y_p) \subset Y$.
\end{enumerate}
\end{defn}

\begin{rem}
\begin{enumerate}
\item
Note that in the above definition $F_q$ and $h^F_{pq}$ are
required to {\it exist} for $F$ to be an embedding of an orbifold,
but they are {\it not}   a part of the data which defines
an embedding of an orbifold.
In other words, two embeddings between orbifolds are equal if they coincide
set theoretically. (Namely they coincide as maps between sets.)
\item
In general, we need to be very careful to define the notion of morphisms between orbifolds.
In fact a natural framework to define the category of orbifolds is that of 2 category.
In other words, the correct notion to define is not two morphisms being the same
but being equivalent.
\item
On the other hand, as long as we use only effective orbifolds and
embeddings in the above sense, we do not need to use 2 category.
In fact, in such case orbifolds and embeddings between them
can be studied in a similar way as the case of  manifolds
and embeddings between them.
This simplifies the discussion significantly, especially when implementing
the Kuranishi structure in applications. Because of this,
we use only embeddings of orbifolds and no other maps between them
for the purpose of the study of Kuranishi structures.
\end{enumerate}
\end{rem}
It is easy to see that the (set theoretical) composition of embeddings of
orbifolds is an embedding of an orbifold.
\begin{rem}\label{notembedding}
Let $\Z_2$ act on $\R^2$ by $(x,y) \mapsto (-x,-y)$.
The quotient space $\R^2/\Z_2$ has a natural structure of orbifold.
We regard one point $\{p\}$ as a manifold and hence is an orbifold.
The map which sends $p$ to the equivalence class of
$(0,0)$ is a continuous map and is
a topological embedding  : $\{p\} \to \R^2/\Z_2$.
But it is not an embedding of orbifold in our sense.
In fact the isotropy group is not isomorphic.
\end{rem}
\begin{defn}\label{diffeoofd}
Two embeddings of orbifolds are said to be the {\it same} if they coincide set theoretically.
\par
An identity map is an embedding of orbifold.
\par
An embedding of an orbifold is said to be a {\it diffeomorphism} if it has an inverse (as set theoretical map) and if the inverse is an embedding of an orbifold.
\par
Suppose we have two orbifold structures on the same space $X$.
We say they are the {\it same} if the identity is a diffeomorphism
between two orbifold structures.
\par
An open set of an orbifold has an obvious orbifold structure.
\par
We say an embedding of an orbifold to be an
{\it open embedding} if it gives a diffeomorphism onto an open set
of the target.
\end{defn}
\begin{rem}
We remark that the definition of two orbifold structures being
the same we gave above is a straightforward generalization
of the well-know definition in the case of manifold structures.
We need to be careful if we try to generalize this definition to the case of Kuranishi structure.
\end{rem}

\begin{rem}\label{localchart}
Sometimes orbifolds are defined in a slightly different way,
as follows.
Let $X$ be a paracompact Hausdorff space and
$\bigcup_i \mathcal U_i = X$ be a locally finite cover.
We consider homeomorphisms
$$
\psi_i :U_i =  V_i/\Gamma_i \to \mathcal U_i,
$$
where $V_i$ is a smooth manifold and $\Gamma_i$
is a finite group with smooth and effective action on $V_i$.
We say that $(\mathcal U_i,\psi_i,V_i,\Gamma_i)$
defines an orbifold structure on $X$ if
the maps
\begin{equation}\label{psiijtrans}
\psi_{ji} = \psi^{-1}_j \circ \psi_i : U_{ji} = \psi_i^{-1}(\mathcal U_j)
\to U_j
\end{equation}
are open embeddings of orbifolds in the sense of
Definition \ref{diffeoofd} for each $i,j$.
\par
It is easy to show that this definition is equivalent to Definition \ref{ofddefn}.
\end{rem}
Various notions appearing in the theory of manifold, such as Riemannian metric,
differential form, integration etc. can be generalized to the
case of orbifold in a straightforward way.
\begin{defn}[orbibundle]\label{orbibundledef}
Let $X$ be a paracompact Hausdorff space.
Suppose we are given an orbifold structure
$(\{ (V_p,\Gamma_p,\psi_p) \},\{(V_{pq},\phi_{pq},h_{pq})\})$
on $X$ in the sense of
Definition \ref{ofddefn}.
\par
A vector bundle on $X$ (we call it an {\it orbibundle} sometimes also)
is an orbifold $\mathcal E$ together with a continuous map
$\pi : \mathcal E \to X$ and $E_{p}, \widehat{\phi}_{pq},
\widehat{\psi}_p$ for
each $p$ or $p,q$ such that the following holds.
\begin{enumerate}
\item
$E_p \to V_p$ is a $\Gamma_p$ equivariant
smooth vector bundle on a manifold $V_p$.
\item
$\widehat{\phi}_{pq} : E_q\vert_{V_{pq}} \to E_p$
is an $h_{pq}$ equivariant bundle map over ${\phi}_{pq} $,
that is a fiberwise isomorphism.
\item
$\widehat{\psi}_p : E_p/\Gamma_p \to \pi^{-1}(\mathcal U_p)$
is a diffeomorphism of orbifolds such that
$$
\begin{CD}
E_p/\Gamma_p @>{\widehat\psi_p}>> \pi^{-1}(\mathcal U_p) \\
@VVV  @VV{\pi}V \\
U_p @>>{\psi_p}>  \mathcal U_p
\end{CD}
$$
and
$$
\begin{CD}
(E_q\vert_{V_{pq}})/\Gamma_q @>{\widehat\psi_q}>> \pi^{-1}(\mathcal U_{pq}) \\
@V{\underline{\widehat{\phi}}_{pq}}VV  @VVV \\
E_p/\Gamma_p @>>{\widehat\psi_p}>  \pi^{-1}(\mathcal U_{p})
\end{CD}
$$
commute.
Here $\underline{\widehat{\phi}}_{pq}$ is induced by
${\widehat{\phi}}_{pq}$ and $\mathcal U_{pq}
= \psi_q(V_{pq}/\Gamma_q)$.
\item
The rank of the vector bundle $E_p$ is
independent of $p$.\footnote{This condition is automatic if $X$ is connected.}
We call it the {\it rank} of our vector bundle $\mathcal E$.
\end{enumerate}
More precisely we say $((\mathcal E,\pi),\{(E_{p},
\widehat{\psi}_p)\},\{\widehat{\phi}_{pq}\})$
is a vector bundle. Sometimes we say
$\mathcal E$ is a {\it vector bundle} etc. by an abuse of notation.
\end{defn}
\begin{exm}
If $X$ is an orbifold, its tangent bundle $TX$ is defined as an orbibundle on $X$
in an obvious way.
\end{exm}
We can define (Whitney) sum and tensor product of
orbibundles in an obvious way.

\begin{defn}[Embedding of vector bundle]\label{defn:embedding} Let $((\mathcal E^X,\pi),\{(E^X_{p},
\widehat{\psi}^X_p)\},\{\widehat{\phi}^X_{pq}\})$ and
$((\mathcal E^Y,\pi),\{(E^Y_{p},
\widehat{\psi}^Y_p)\},\{\widehat{\phi}^Y_{pq}\})$
be
vector bundles over orbifolds $X$ and $Y$,
respectively.
Let $F : X \to Y$ be an embedding of an orbifold.
Then an {\it embedding of a  vector bundle} $\hat F :
\EE_X \to \EE_Y$  over $F$ is an embedding of an orbifold
such that the following holds in addition.
\par
Let $q\in X$ and $p = F(q)$.
Let $V^X_{pq}$, $F_{q} : V^X_{pq} \to V^Y_p$,
$h^F_{pq} : \Gamma_q^X \to \Gamma^Y_p$
be as in Definition \ref{emborf}.
Then there exists an $h^F_{pq}$ equivariant
embedding of vector bundles
$$
\hat F_{pq} : E^X_q\vert_{V^X_{pq}} \to E^Y_p
$$
such that the following diagram commutes.
$$
\begin{CD}
E^X_q\vert_{V^X_{pq}/\Gamma^X_q} @>{\widehat\psi^X_q}>> \pi^{-1}(\mathcal U^X_{pq}) \\
@V{\underline{\hat F}_{pq}}VV  @VV{\widehat{F}}V \\
E^Y_p/\Gamma_p^Y @>>{\widehat\psi^Y_p}>  \pi^{-1}(\mathcal U^Y_{p})
\end{CD}
$$
Here $\underline{\hat F}_{pq}$ is induced from
${\hat F}_{pq}$.
\end{defn}
\begin{defn}
\begin{enumerate}
\item
An embedding of a vector bundle $\hat F :
\EE_X \to \EE_Y$ is said to be an  {\it open embedding}
of a vector bundle if it is an open embedding
as a map between orbifolds.
\item
An {\it isomorphism} between vector bundles is an
embedding
of a vector bundle that is a diffeomorphism
between orbifolds.
\item
Two vector bundles on a given orbifold $X$ is said to be
{\it isomorphic} if there exists an isomorphism which covers
the identity map between them.
\item
If $\mathcal E_a$ $a=1,2$ are orbibundles over the same orbifold $X$
and $\hat F :  \mathcal E_1 \to \mathcal E_2$ is an embedding
of an orbibundle over the identity map, then we say
$\mathcal E_1$ is a {\it subbundle} of $\mathcal E_2$.
\item
If $\mathcal E_1$ is a subbundle of $\mathcal E_2$,
we can define the quotient bundle $\mathcal E_2/\mathcal E_2$
in an obvious way.
\end{enumerate}
\end{defn}
It is easy to see that the (set theoretical) composition of
embeddings of orbibundle is an embedding of orbibundle.
\begin{rem}
If two orbifold structures on $X$ are the same,
the notion of vector bundles on it is the same in the following sense.
The isomorphism classes of vector bundles
on one structure corresponds
to the isomorphism classes of vector bundles
on the other structure by a canonical bijection.
The proof is easy and is left to the reader.
\end{rem}
\begin{rem}
Suppose a structure of orbifold on $X$ is given
by $(\mathcal U_i,\psi_i,V_i,\Gamma_i)$
as in Remark \ref{localchart}.
Then we can define a vector bundle on it
by $(E_i,\hat\psi_{ij})$ as follows.
\begin{enumerate}
\item
$E_i \to V_i$ is a $\Gamma_i$ equivariant
vector bundle.
\item
Let
$\psi_{ji} = \psi^{-1}_j \circ \psi_i : U_{ji} = \psi_i^{-1}(\mathcal U_j)
\to U_j$
be as in (\ref{psiijtrans}).
Then
$$
\hat\psi_{ij} :
E_i/\Gamma_i\vert_{U_{ji}} \to E_j/\Gamma_j
$$
be an open embedding of vector bundle in the above sense.
\end{enumerate}
It is easy to see that this definition is equivalent to
Definition \ref{orbibundledef}.
\end{rem}
\begin{lem}[Induced bundle]
Let $F : X \to Y$ be an embedding of an orbifold
and $\mathcal E^Y$ be a vector bundle on $Y$.
Then there exists a vector bundle $\mathcal E^X$ on $X$
with the same rank as $\mathcal E^Y$
and an embedding of an orbibundle $\hat F : \mathcal E^X
\to \mathcal E^Y$ which covers $F$.
\par
$(\mathcal E^X,\hat F)$ is unique in the following sense.
If $(\mathcal E_a^X,\hat F_a)$, $a=1,2$ are two such
choices, then there exists an isomorphism
of vector bundles $\hat I : \mathcal E_1^X \to \mathcal E_2^X$
such that it covers identity and satisfies
$
\hat F_2 \circ \hat I  = \hat F_1.
$
\end{lem}
The proof is easy and is left to the reader.
We call $\mathcal E^X$ the {\it induced bundle} and write
$F^*\mathcal E^Y$.
In the case $X$ is an open subset of $Y$,
$\mathcal E^X = F^* \mathcal E^Y$ is called the {\it restriction}
of $\mathcal E^Y$ to $X$.
\begin{rem}
When we consider a map between orbifolds which is not an
embedding, the induced bundle may not be well-defined.
\end{rem}
\begin{exm}
If $F : X \to Y$ is an embedding of an orbifold,
it induces an embedding of orbibundles
$dF : TX \to TY$ between their tangent bundles in an obvious way.
Therefore $TX$ is a subbundle of the pull back bundle $F^*TY$.
\end{exm}
\begin{defn}
If $F : X \to Y$ is an embedding of an orbifold, the normal bundle $N_XY$
is by definition the quotient bundle $F^*TY/TX$.
\end{defn}
\begin{rem}
We can prove tubular neighborhood theorem in an obvious way.
Namely $N_XY$ as an orbifold is
diffeomorphic to an open neighborhood of $F(X)$ in $Y$.
The proof is similar to the proof of tubular neighborhood theorem
for manifolds.
\end{rem}
Finally we mention two easy lemmas about gluing orbifolds
or orbibundles by a diffeomorphism or an isomorphism.
\par
Let $X$, $Y$ be orbifolds
and $U_{YX} \subset X$ be an open set.
Suppose $F : U_{YX} \to Y$
is an open embedding of orbifolds.
\par
We define an equivalence relation $\sim $ on
the disjoint union
$X \sqcup Y$ by
the following.
We say $x\sim y$ if and only if
\begin{enumerate}
\item
$x = y$, or
\item
$x \in U_{XY}$ and $y \in Y$ and $y = F(x)$,
\item
$y \in U_{XY}$ and $x \in Y$ and $x = F(y)$.
\end{enumerate}
Let $X \cup_F Y$ be the set of equivalence classes equipped
with quotient topology.
\begin{lem}\label{gluelemmaord}
If $X \cup_F Y$ is Hausdorff, it has a structure of an orbifold such that
the maps
\begin{equation}\label{XYembedtounion}
X \to X \cup_F Y, \qquad
Y \to X\cup_F Y
\end{equation}
which send $x \in X$ (resp. $y\in Y$) to its equivalence class
are open embeddings of orbifolds.
\end{lem}
The proof is easy and is left to the reader.

\begin{rem}
One needs to be very careful when trying to generalize
Lemma \ref{gluelemmaord} to
uneffective orbifolds.
\end{rem}

\begin{lem}
In the situation of Lemma \ref{gluelemmaord}
we assume in addition that we are given two orbibundles $\mathcal E^X$ over $X$ and
$\mathcal E^Y$ over $Y$ with the same rank.
Suppose furthermore that $F$ is covered by an open embedding
$\hat F : \mathcal E^X\vert_{U_{YX}} \to \mathcal E^Y$
of an orbibundle. Here $\mathcal E^X\vert_{U_{YX}}$ is the
restriction of an orbibundle.
\par
Then there exists a structure of orbibundle
on $\mathcal E^X \cup_{\hat F} \mathcal E^Y$ together
with embeddings of orbibundles
$$
\mathcal E^X \to \mathcal E^X \cup_{\hat F} \mathcal E^Y, \qquad
\mathcal E^Y \to \mathcal E^X \cup_{\hat F} \mathcal E^Y
$$
which cover the maps in (\ref{XYembedtounion}).
\end{lem}
The proof is easy and is left for the reader.
\par\newpage
\part{Construction of the Kuranishi structure 1: Gluing analysis in the
simple case}
\label{secsimple}

In this part we give detailed proof of the gluing analysis (or stretching the neck)
of the pseudo-holomorphic curve with node, and also decay estimate of the exponential
order with respect to the length of the neck.
This analysis implies the smoothness of the Kuranishi chart at infinity.
In Part \ref{secsimple} we discuss a simple case where
we have two pseudo-holomorphic maps from stable and irreducible curves
joined at one point.
The general case will be discussed in the next Part.

\section{Setting}
We will describe the general case in Part \ref{generalcase}. To simplify the notation and clarify the main analytic point of the proof
we prove the case where we glue holomorphic maps from two stable bordered Riemann surfaces
to $(X,L)$ in this section.
\par
Let $\Sigma_i$ be a bordered Riemann surface with one end.
($i=1,2$.)
We assume that there are compact subsets $K_i \subset \Sigma_i$ such that
$\Sigma_i \setminus K_i$ are half infinite cylinders.
For $T>0$, we put coordinates $[-5T, \infty) \times [0,1]$, resp.
$(-\infty, 5T] \times [0,1]$ on $\Sigma_i \setminus K_i$, $i=1,2$,
respectively.
We identify their ends as follows.
\begin{equation}
\aligned
\Sigma_1 &= K_1 \cup ((-5T,\infty)\times [0,1]), \\
\Sigma_2 &= ((-\infty,5T)\times [0,1]) \cup K_2.
\endaligned
\end{equation}
Here $K_i$ are compact and $\pm \infty$ are the ends.
We put
\begin{equation}
\Sigma_T = K_1 \cup ((-5T,5T)\times [0,1]) \cup  K_2.
\end{equation}
We use $\tau$ for the coordinate of the factors $(-5T,\infty)$,
$(-\infty,5T)$, or $(-5T,5T)$ and $t$ for the coordinate of the second factor
$[0,1]$.
\par
Let $X$ be a symplectic manifold with compatible (or tame) almost complex structure
and $L$ be its Lagrangian submanifold.
\par
Let
$$
u_i : (\Sigma_i,\partial \Sigma_i)  \to  (X,L), \qquad i=1,2
$$
be pseudo-holomorphic maps of finite energy.
Then, by the removable singularity theorem that is now standard, we have asymptotic value
\begin{equation}\label{tauinf}
\lim_{\tau\to \infty} u_1(\tau,t) \in L
\end{equation}
and
\begin{equation}\label{tauminf}
\lim_{\tau\to -\infty} u_2(\tau,t) \in L.
\end{equation}
The limits (\ref{tauinf}) and
(\ref{tauminf}) are independent of $t$.
\par
We assume that the limit (\ref{tauinf}) coincides with
(\ref{tauminf}) and denote it by $p_0 \in L$.
\par
We fix a coordinate of $X$ and of $L$ in a neighborhood of $p_0$.
So a trivialization of the tangent bundle $TX$ and $TL$ in a neighborhood of
$p_0$ is fixed.
Hereafter we assume the following:
\begin{equation}
\text{\rm Diam}(u_1([-5T,\infty)\times [0,1])) \le \epsilon_1,
\qquad
\text{\rm Diam}(u_2((-\infty,5T]\times [0,1])) \le \epsilon_1.
\end{equation}
\par
The maps $u_i$ determine  homology classes $\beta_i = [u_i] \in H_2(X,L)$.
\par
We take $K_i^{\rm obst}$ a compact subset of the interior of $K_i$ and take
\begin{equation}\label{Eitake}
E_i \subset \Gamma(K_i^{\rm obst};u_i^*TX\otimes \Lambda^{0,1})
\end{equation}
a finite dimensional linear subspace consisting of smooth sections supported in
$K_i^{\rm obst}$.
\par
For simplicity  we also fix a complex structure of the
source $\Sigma_i$. The version where it can move will be
discussed in Part \ref{generalcase}.
We also assume that $\Sigma_i$ equipped with marked points $\vec z_i$ is stable.
The process to add marked points to stabilize it will be discussed in Part \ref{generalcase} also.
Let
\begin{equation}\label{lineeq}
D_{u_i}\overline{\partial}
:
L^2_{m+1,\delta}((\Sigma_i,\partial \Sigma_i);u_i^*TX,u_i^*TL) \to
L^2_{m,\delta}(\Sigma_i;u_i^{*}TX \otimes \Lambda^{01})
\end{equation}
be the linearization of the Cauchy-Riemann equation.
Here we define the weighted Sobolev space we use as follows.
\begin{defn}(\cite[Section 7.1.3]{fooo:book1})\footnote{In \cite{fooo:book1} $L^p_1$ space
is used in stead of $L^2_{m}$ space.}
Let $L^2_{m+1,loc}((\Sigma_i,\partial \Sigma_i);u_i^*TX;u_i^*TL)$ be the set of the
sections $s$ of $u_i^*TX$ which is locally of $L^2_{m+1}$-class, (Namely its differential up to order $m+1$
is of $L^2$ class. Here $m$ is sufficiently large, say larger than $10$.)
We also assume $s(z) \in u_i^*TL$ for $z \in \partial \Sigma_i$.
\par
The weighted Sobolev space
$L^2_{m+1,\delta}((\Sigma_i,\partial \Sigma_i);u_i^*TX,u_i^*TL)$
is the set of all pairs $(s,v)$ of elements $s$ of $L^2_{m+1,loc}((\Sigma_i,\partial \Sigma_i);u_i^*TX;u_i^*TL)$
and $v \in T_{p_0}L$, (here  $p_0 \in L$ is the point (\ref{tauinf}) or (\ref{tauminf})) such that
\begin{equation}\label{weight1}
\sum_{k=0}^{m+1} \int_{\Sigma_i \setminus K_i}  e^{\delta\vert \tau \pm 5T\vert}\vert \nabla^k(s - \text{\rm Pal}(v))\vert^2 < \infty,
\end{equation}
where $\text{\rm Pal} : T_{p_0}X \to T_{u_i(\tau,t)}X$ is defined by the trivialization we fixed right after (\ref{tauminf}).
(Here $\pm$ is $+$ for $i=1$ and $-$ for $i=2$.)
The norm is defined as the sum of (\ref{weight1}), the norm of $v$ and the $L^2_{m+1}$ norm of $s$ on
$K_i$. (See (\ref{normformjula}).)
\par
$L^2_{m,\delta}(\Sigma_i;u_i^{*}TX \otimes \Lambda^{01})$
is defined similarly without boundary condition and with out $v$.
(See (\ref{normformjula52}).)
\end{defn}
When we define $D_{u_i}\overline{\partial}$ we forget $v$ component and use $s$ only.
\begin{rem}
The positive number $\delta$ is chosen as follows. (\ref{tauinf}) and a standard estimate imply that there exists
$\delta_1 > 0$ such that
\begin{equation}\label{approestu}
\left\vert\frac{d}{d\tau}u_i\right\vert_{C^k}(\tau,t) < C_ke^{-\delta_1 \vert \tau\vert}
\end{equation}
for any $k$. We choose $\delta$ smaller than $\delta_1/10$.
(\ref{approestu}) implies
$$
(D_{u_i}\overline\partial)(\text{\rm Pal}(v)) < C_ke^{-\delta_1 \vert \tau\vert/10}.
$$
Therefore (\ref{lineeq}) is defined and bounded.
\end{rem}
It is a standard fact that (\ref{lineeq}) is Fredholm.
\par
We work under the following assumption.
\begin{assump}\label{DuimodEi}
\begin{equation}\label{DuimodEi0}
D_{u_i}\overline{\partial} :
L^2_{m+1,\delta}((\Sigma_i,\partial \Sigma_i);u_i^*TX,u_i^*TL)
\to L^2_{m,\delta}(\Sigma_i;u_i^{*}TX \otimes \Lambda^{01})/E_i
\end{equation}
is surjective.
Moreover the following (\ref{Duievsurj}) holds.
Let
$
(D_{u_i}\overline{\partial})^{-1}(E_i)
$
be the kernel of (\ref{DuimodEi0}). We define
\begin{equation}\label{Duiev}
D{\rm ev}_{i,\infty} : L^2_{m+1,\delta}((\Sigma_i,\partial \Sigma_i);u_i^*TX,u_i^*TL)
\to  T_{p_0}L
\end{equation}
by
$$
D{\rm ev}_{i,\infty}(s,v) = v.
$$
Then
\begin{equation}\label{Duievsurj}
D{\rm ev}_{1,\infty} - D{\rm ev}_{2,\infty} : (D_{u_1}\overline{\partial})^{-1}(E_1) \oplus (D_{u_2}\overline{\partial})^{-1}(E_2)
\to T_{p_0}L
\end{equation}
is surjective.
\end{assump}
\par
Let us start stating the result.
Let
\begin{equation}\label{uprime}
u' : (\Sigma_T,\partial\Sigma_T) \to (X,L)
\end{equation}
be a smooth map.
We consider the following condition depending $\epsilon >0$.
\begin{conds}\label{nearbyuprime}
\begin{enumerate}
\item
$u'\vert_{K_i}$ is $\epsilon$-close to $u_i\vert_{K_i}$ in $C^1$ sense.
\item
The diameter of $u'([-5T,5T]\times [0,1])$ is smaller than $\epsilon$.
\end{enumerate}
\end{conds}
\par
We take $\epsilon_2$ sufficiently small compared to the `injectivity radius' of $X$
so that the next definition
makes sense.\footnote{More precisely, we assume that
$$
\{(x,y) \in X \times X \mid d(x,y) < \epsilon_2\} \subset {\rm E}(\{ (x,v) \in TX \mid \vert v\vert < \epsilon\}),
$$
where ${\rm E} : \{ (x,v) \in TX \mid \vert v\vert < \epsilon\} \to X$ is induced by an exponential map
of certain connection of $TX$. See (\ref{defE}).}
For $u'$ satisfying  Condition \ref{nearbyuprime} for $\epsilon <\epsilon_2$ :
$$
I_{u'} : E_i \to \Gamma(\Sigma_T;(u')^*TX \otimes \Lambda^{01})
$$
is the complex linear part of the parallel translation along
the short geodesic (between $u_i(z)$ and $u'(z)$. Here $z \in K_i^{\text{\rm obst}}$).
We put
\begin{equation}
E_i(u') =  I_{u'}(E_i).
\end{equation}
\par
The equation we study is
\begin{equation}\label{mainequation}
\overline\partial u' \equiv 0 , \quad \mod E_1(u') \oplus  E_2(u').
\end{equation}
\begin{rem}
In the actual construction of Kuranishi structure,
we take several $u_i$'s and take $E_i$'s for each of them. Then
in place of $E_1(u') \oplus  E_2(u')$ we take sum of
finitely many of them.
Here we simplify the notation. There are not so many differences  between the proof of
Theorem \ref{gluethm1} and the corresponding result in case we take several such $u_i$'s and  $E_i$'s.
See \cite[pages 4-5]{Fu1} and Section \ref{glueing}.
\end{rem}
Theorem \ref{gluethm1} describes all the solutions of (\ref{mainequation}).
To state this precisely we need a bit more notations.
\par
We consider the following condition for $u'_i :  (\Sigma_i,\partial\Sigma_i) \to (X,L)$.
\begin{conds}\label{uiconds}
\begin{enumerate}
\item
$u_i'\vert_{K_i}$ is $\epsilon$-close to $u_i\vert_{K_i}$ in $C^1$ sense.
\item
The diameter of $u_1'([-5T,\infty)\times [0,1])$, (resp.  $u_2'((-\infty,5T])\times [0,1])$)
is smaller than $\epsilon$.
\end{enumerate}
\end{conds}
Then we define
$$
I_{u'_i} : E_i \to \Gamma(\Sigma_i;(u_i')^*TX \otimes \Lambda^{01})
$$
by using the parallel transport in the same way as $I_{u'_T}$.
(This makes sense if $u'_i$ satisfies Condition \ref{uiconds} for $\epsilon < \epsilon_2$.)
We put
\begin{equation}
E_i(u'_i) =  I_{u'_i}(E_i).
\end{equation}
So we can define an equation
\begin{equation}\label{mainequationui}
\overline\partial u_i' \equiv 0 , \quad \mod E_i(u'_i).
\end{equation}
\begin{defn}
The set of solutions of equation (\ref{mainequationui}) with finite energy
and satisfying Condition \ref{uiconds} for $\epsilon = \epsilon_2$ is
denoted by $\mathcal M^{E_i}((\Sigma_i,\vec z_i);\beta_i)_{\epsilon_2}$.
Here $\beta_i$ is the homology class of $u_i$.
\end{defn}
\begin{rem}
In the usual story of pseudo-holomorphic curve, we identify  $u_i$ and $u_i'$ if there exists
a biholomorphic map $v : (\Sigma_i,\vec z_i)\to (\Sigma_i,\vec z_i)$ such that
$u'_i = u_i\circ v$.
In our situation where $\Sigma_i$ has no sphere or disk bubble and has nontrivial
boundary with at least one boundary marked points (that is $\tau =\pm\infty$), such $v$ is necessary the identity map.
Namely $\Sigma_i$ has no nontrivial automorphism.
\end{rem}
The surjectivity of (\ref{Duiev}), (\ref{Duievsurj}) and the implicit function theorem
imply that if $\epsilon_2$ is small then there exists a finite dimensional vector space $\tilde V_i$ and its neighborhood $V_i$ of $0$
such that
$$
\mathcal M^{E_i}((\Sigma_i,\vec z_i);\beta_i)_{\epsilon_2}
\cong
V_i.
$$
Since we assume that $\Sigma_i$ is nonsingular the group $\text{\rm Aut}((\Sigma_i,\vec z_i),u_i)$ is trivial.
(In the case when there is a sphere bubble, the automorphism group can be nontrivial. That case will
be discussed later.)
\par
For any $\rho_i \in V_i$ we denote by $u_i^{\rho_i} : (\Sigma_i,\partial\Sigma_i) \to (X,L)$
the corresponding solution of (\ref{mainequationui}).
\par
We have an evaluation map
$$
\text{\rm ev}_{i,\infty} : \mathcal M^{E_i}((\Sigma_i,\vec z_i);\beta_i)_{\epsilon_2} \to L
$$
that is smooth. Namely
$$
\text{\rm ev}_{i,\infty} (u'_i)  = \lim_{\tau\to \pm\infty}u'_i(\tau,t).
$$
(Here $\pm = +$ for $i=1$ and $-$ for $i=2$.)\footnote{This is a consequence
of the fact that $u_i$ is pseudo-holomorphic outside a compact set and has finite energy.}
We consider the fiber product:
\begin{equation}\label{fpmoduli}
\mathcal M^{E_1}((\Sigma_1,\vec z_1);\beta_1)_{\epsilon_2} \times_L \mathcal M^{E_2}((\Sigma_2,\vec z_2);\beta_2)_{\epsilon_2}.
\end{equation}
The surjectivity of (\ref{Duievsurj}) implies that this fiber product is transversal
so is
$$
V_1 \times_L V_2.
$$
And an element of $V_1 \times_L V_2$ is written as $\rho = (\rho_1,\rho_2)$.
\begin{defn}
Let $\beta =\beta_1 + \beta_2$.
We denote by $\mathcal M^{E_1+E_2}((\Sigma_T,\vec z);\beta)_{\epsilon}$
the set of solutions of (\ref{mainequation}) satisfying the Condition \ref{nearbyuprime}
with $\epsilon_2 = \epsilon$.
\end{defn}
\begin{thm}\label{gluethm1}
For each sufficiently small $\epsilon_3$ and sufficiently large $T$, there exist $\epsilon_1, \epsilon_2$ and a map
$$
\text{\rm Glu}_T :
\mathcal M^{E_1}((\Sigma_1,\vec z_1);\beta_1)_{\epsilon_2} \times_L \mathcal M^{E_2}((\Sigma_2,\vec z_2);\beta_2)_{\epsilon_2}
\to
\mathcal M^{E_1+E_2}((\Sigma_T,\vec z);\beta)_{\epsilon_1}
$$
that is a diffeomorphism to its image.
The image contains $\mathcal M^{E_1+E_2}((\Sigma_T,\vec z);\beta)_{\epsilon_3}$.
\end{thm}
The result about exponential decay estimate of this map is in Section \ref{subsecdecayT}.
(Theorem \ref{exdecayT}.)
\par\bigskip
\section{Proof of Theorem \ref{gluethm1} : 1 - Bump function and
weighted Sobolev norm}
\label{subsec12}
The proof of Theorem \ref{gluethm1}
was given in \cite[Section 7.1.3]{fooo:book1}.
The exponential decay estimate of the solution was
proved in \cite[Section A1.4]{fooo:book1} together with
a slightly modified version of the proof of Theorem \ref{gluethm1}.
Here we follow the proof of \cite[Section A1.4]{fooo:book1} and
give its more detail.
As mentioned there the origin of the proof is Donaldson's paper \cite{Don86I},
and its Bott-Morse version in \cite{Fuk96II}.
\par
We first introduce certain bump functions.
First let $\mathcal A_T \subset \Sigma_T$
and $\mathcal B_T \subset \Sigma_T$ be the domains
defined by
$$
\mathcal A_T = [-T-1,-T+1] \times [0,1],
\qquad
\mathcal B_T = [T-1,T+1] \times [0,1].
$$
We may regard
$\mathcal A_T,\mathcal B_T \subset \Sigma_i$.
The third domain is
$$
\mathcal X = [-1,1] \times [0,1] \subset \Sigma_T.
$$
We may also regard $\mathcal X \subset \Sigma_i$.
\par
Let $\chi_{\mathcal A}^{\leftarrow}$,
$\chi_{\mathcal A}^{\rightarrow}$ be smooth functions on $[-5T,5T]\times [0,1]$ such that
\begin{equation}
\chi_{\mathcal A}^{\leftarrow}(\tau,t)
= \begin{cases}
1   & \tau < -T-1 \\
0  & \tau > -T+1.
\end{cases}
\end{equation}
$$
\chi_{\mathcal A}^{\rightarrow}  = 1 - \chi_{\mathcal A}^{\leftarrow}.
$$
\par
We define
\begin{equation}
\chi_{\mathcal B}^{\leftarrow}(\tau,t)
= \begin{cases}
1   & \tau <  T-1 \\
0  & \tau >  T+1.
\end{cases}
\end{equation}
$$
\chi_{\mathcal B}^{\rightarrow}  = 1 - \chi_{\mathcal B}^{\leftarrow}.
$$
\par
We define
\begin{equation}
\chi_{\mathcal X}^{\leftarrow}(\tau,t)
= \begin{cases}
1   & \tau <  -1 \\
0  & \tau >  1.
\end{cases}
\end{equation}
$$
\chi_{\mathcal X}^{\rightarrow}  = 1 - \chi_{\mathcal X}^{\leftarrow}.
$$
We extend these functions to $\Sigma_T$ and $\Sigma_i$ ($i=1,2$) so that
they are locally constant outside $[-5T,5T]\times [0,1]$.
We denote them by the same symbol.
\par\medskip
We next introduce weighted Sobolev norms and their local versions for sections
on $\Sigma_T$ or $\Sigma_i$ as follows.
\par
We define $e_{i,\delta} : \Sigma_i \to [1,\infty)$ of $C^{\infty}$ class
as follows.
\begin{equation}\label{e1delta}
e_{1,\delta}(\tau,t)
\begin{cases}
=e^{\delta\vert \tau + 5T\vert} &\text{if $\tau > 1 -  5T$ }\\
=1   &\text{on  $K_1$} \\
\in [1,10] &\text{if $\tau < 1 -  5T$}
\end{cases}
\end{equation}
\begin{equation}\label{e2delta}
e_{2,\delta}(\tau,t)
\begin{cases}
=e^{\delta\vert \tau - 5T\vert} &\text{if $\tau <  5T-1$ }\\
=1   &\text{on  $K_2$} \\
\in [1,10] &\text{if $\tau >  5T-1$}
\end{cases}
\end{equation}
We also define $e_{T,\delta} : \Sigma_T \to [1,\infty)$ as follows:
\begin{equation}\label{e2delta}
e_{T,\delta}(\tau,t)
\begin{cases}
=e^{\delta\vert \tau - 5T\vert} &\text{if $1<\tau <  5T-1$ }\\
= e^{\delta\vert \tau + 5T\vert} &\text{if $-1>\tau > 1-5T$ }\\
=1   &\text{on  $K_1\cup K_2$} \\
\in [1,10] &\text{if $\vert\tau - 5T\vert < 1$ or $\vert\tau + 5T\vert < 1$} \\
\in [e^{5T\delta}/10,e^{5T\delta}] &\text{if $\vert\tau\vert < 1$}.
\end{cases}
\end{equation}
The weighted Sobolev norm we use for
$L^2_{m,\delta}(\Sigma_i;u_i^*TX\otimes \Lambda^{01})$
is
\begin{equation}
\Vert s\Vert^2_{L^2_{m,\delta}} = \sum_{k=0}^m \int_{\Sigma_i}
e_{i,\delta} \vert \nabla^k s\vert^2 \text{\rm vol}_{\Sigma_i}.
\end{equation}
For
$(s,v) \in L^2_{m+1,\delta}((\Sigma_i,\partial \Sigma_i);u_i^*TX,u_i^*TL)$
we define
\begin{equation}\label{normformjula}
\aligned
\Vert (s,v)\Vert^2_{L^2_{m+1,\delta}} = &\sum_{k=0}^{m+1} \int_{K_i}
\vert \nabla^k s\vert^2 \text{\rm vol}_{\Sigma_i}\\
&+
\sum_{k=0}^{m+1} \int_{\Sigma_i \setminus K_i}  e_{i,\delta}\vert \nabla^k(s - \text{\rm Pal}(v))\vert^2 \text{\rm vol}_{\Sigma_i}
+ \Vert v\Vert^2.
\endaligned
\end{equation}
We next define a weighted Sobolev norm for the sections on
$\Sigma_T$.
Let
$$
s \in L^2_{m+1}((\Sigma_T,\partial \Sigma_T);u^*TX,u^*TL).
$$
Since we take $m$ large, $s$ is continuous.
So $s(0,1/2) \in T_{u(0,1/2)}X$ is well defined.
There is a canonical trivialization of $TX$ in a neighborhood of $p_0$ that we fixed right
after (\ref{tauminf}). We use it to define $\text{\rm Pal}$ below.
We put
\begin{equation}\label{normformjula5}
\aligned
\Vert s\Vert^2_{L^2_{m+1,\delta}} = &\sum_{k=0}^{m+1} \int_{K_1}
\vert \nabla^k s\vert^2 \text{\rm vol}_{\Sigma_1} + \sum_{k=0}^{m+1} \int_{K_2}
\vert \nabla^k s\vert^2 \text{\rm vol}_{\Sigma_2}\\
&+
\sum_{k=0}^{m+1} \int_{[-5T,5T]\times [0,1]}  e_{T,\delta}\vert \nabla^k(s - \text{\rm Pal}(s(0,1/2)))\vert^2 \text{\rm vol}_{\Sigma_i}
\\&+ \Vert s(0,1/2)\Vert^2.
\endaligned
\end{equation}
For
$$
s \in L^2_{m}((\Sigma_T,\partial \Sigma_T);u^*TX\otimes \Lambda^{01})
$$
we define
\begin{equation}\label{normformjula52}
\Vert s\Vert^2_{L^2_{m,\delta}} = \sum_{k=0}^{m} \int_{\Sigma_T}
e_{T,\delta}\vert  \nabla^k s\vert^2 \text{\rm vol}_{\Sigma_1}.
\end{equation}
These norms were used in \cite[Section 7.1.3]{fooo:book1}.
\par
For a subset $W$ of $\Sigma_i$ or $\Sigma_T$ we define
$
\Vert s\Vert_{L^2_{m,\delta}(W\subset \Sigma_i)}
$,
$
\Vert s\Vert_{L^2_{m,\delta}(W\subset \Sigma_T)}
$
by restricting the domain of the integration (\ref{normformjula52}) or
(\ref{normformjula5}) to $W$.
\par
Let $(s_j,v_j) \in L^2_{m+1,\delta}((\Sigma_i,\partial \Sigma_i);u_i^*TX,u_i^*TL)$
for $j=1,2$. We define the inner product among them by:
\begin{equation}\label{innerprod1}
\aligned
\langle\!\langle (s_1,v_1),(s_2,v_2)\rangle\!\rangle_{L^2_{\delta}}
=
& \int_{\Sigma_i\setminus K_i}  (s_1-\text{\rm Pal}v_1,s_2-\text{\rm Pal}v_2)
\\
&+\int_{K_i} (s_1,s_2)
+
( v_1,v_2).
\endaligned
\end{equation}
We also use an exponential map. (The same map was used in \cite[pages 410-411]{fooo:book1}.)
We take a diffeomorphism
\begin{equation}\label{defE}
{\rm E} = ({\rm E} _1,{\rm E} _2): \{ (x,v) \in TX \mid \vert v\vert < \epsilon\} \to X \times
X
\end{equation}
to its image such that
$$
{\rm E} _1(x,v) = x, \quad \left.\frac{d {\rm E} _2(x,tv)}{dt}\right\vert_{t=0} = v
$$
and
$$
{\rm E} (x,v) \in L\times L, \qquad \text{for $x \in L$, $v \in T_xL$}.
$$
Furthermore we may take it so that
\begin{equation}\label{Einanbdofp0}
{\rm E} (x,v) = (x,x+v)
\end{equation}
on a neighborhood of $p_0$.
\par
To find such $\rm E$, we take a linear connection $\nabla$
(that may not be a Levi-Civita connection of a Riemannian metric)
of $TX$ such that $TL$ is parallel with respect to  $\nabla$.
We then use geodesic with respect to $\nabla$ to define an exponential
map. We then define $\rm E$ such that $t\mapsto  {\rm E} _2(x,tv)$ is a geodesic with initial direction $v$.
Note that we may take $\nabla$ so that in a neighborhood of $p_0$ it coincides with the standard trivial connection
with respect the coordinate we fixed. (\ref{Einanbdofp0}) follows.

\section{Proof of Theorem \ref{gluethm1} : 2 - Gluing by alternating method}
\label{alternatingmethod}

Let us start with
$$
u^{\rho} = (u_1^{\rho_1},u_2^{\rho_2})
\in \mathcal M^{E_1}((\Sigma_1,\vec z_1);\beta_1)_{\epsilon_2}
\times_L \mathcal M^{E_2}((\Sigma_2,\vec z_2);\beta_2)_{\epsilon_2}.
$$
Here $\rho_i \in V_i$ and the corresponding map $(\Sigma_i,\partial \Sigma_i) \to (X,L)$
is denoted by $u_i^{\rho_i}$.
Let $\rho = (\rho_1,\rho_2)$.
We put
$$
p^{\rho} = \lim_{\tau\to \infty} u_1^{\rho_1}(\tau,t) = \lim_{\tau\to -\infty} u_2^{\rho_2}(\tau,t).
$$
\par\medskip
\noindent{\bf Pregluing}:
\begin{defn}
We define
\begin{equation}
u_{T,(0)}^{\rho} =
\begin{cases}
\chi_{\mathcal B}^{\leftarrow} (u_1^{\rho_1} - p^{\rho}) +
\chi_{\mathcal A}^{\rightarrow} (u_2^{\rho_2} - p^{\rho}) + p^{\rho}
& \text{on $[-5T,5T] \times [0,1]$} \\
u_1^{\rho_1} & \text{on $K_1$} \\
u_2^{\rho_2} & \text{on $K_2$}.
\end{cases}
\end{equation}
\end{defn}
Note that we use the coordinate of the neighborhood of $p_0$ to define the sum in the first line.
\par\medskip
\noindent{\bf Step 0-3}:
\begin{lem}\label{lem18}
If $\delta < \delta_1/10$ then there exists $\frak e^{\rho} _{i,T,(0)} \in E_i$
such that
\begin{equation}\label{startingestimate}
\Vert\overline\partial u^{\rho} _{T,(0)} - \frak e^{\rho} _{1,T,(0)} - \frak e^{\rho} _{2,T,(0)}\Vert_{L_{m,\delta}^2} < C_{1,m}e^{-\delta T}.
\end{equation}
Moreover
\begin{equation}\label{frakeissmall}
\Vert \frak e^{\rho} _{i,T,(0)}\Vert_{L_{m}^2(K_i)} < \epsilon_{4,m}.
\end{equation}
Here $\epsilon_{4,m}$ is a positive number which we may choose arbitrarily small
by taking $V_i$ to be a sufficiently small neighborhood of zero in $\tilde V_i$.
\par
Moreover $\frak e^{\rho} _{i,T,(0)}$ is independent of $T$.
\end{lem}
\begin{proof}
We put
$$
\frak e_{i,T,(0)} = \overline\partial u^{\rho}_i \in E_i.
$$
Then by definition the support of
$\overline\partial u^{\rho} _{T,(0)} - \frak e^{\rho} _{1,T,(0)} - \frak e^{\rho} _{2,T,(0)}$
is in $[-5T,5T]\times [0,1]$.
Moreover it is estimated as (\ref{startingestimate}).
\end{proof}
\par\medskip
\noindent{\bf Step 0-4}:
\begin{defn}\label{deferfirst}
We put
$$
\aligned
{\rm Err}^{\rho}_{1,T,(0)}
&= \chi_{\mathcal X}^{\leftarrow} (\overline\partial u^{\rho} _{T,(0)} - \frak e^{\rho} _{1,T,(0)}), \\
{\rm Err}^{\rho}_{2,T,(0)}
&= \chi_{\mathcal X}^{\rightarrow} (\overline\partial u^{\rho} _{T,(0)} -  \frak e^{\rho} _{2,T,(0)}).
\endaligned$$
We regard them as elements of the weighted Sobolev spaces
$L^2_{m,\delta}((\Sigma_1,\partial \Sigma_1);(u_1^{\rho})^*TX\otimes \Lambda^{01})$
and
$L^2_{m,\delta}((\Sigma_2,\partial \Sigma_2);(u_2^{\rho})^*TX\otimes \Lambda^{01})$
respectively.
(We extend them by $0$ outside a compact set.)
\end{defn}
\par\medskip
\noindent{\bf Step 1-1}:
We first cut $u^{\rho}_{T,(0)}$ and extend to obtain maps $\hat u^{\rho}_{i,T,(0)} :
(\Sigma_i,\partial\Sigma_i) \to (X,L)$ $(i=1,2)$ as follows.
(This map is used to set the linearized operator (\ref{lineeqstep0}).)
\begin{equation}
\aligned
&\hat u^{\rho}_{1,T,(0)}(z) \\
&=
\begin{cases} \chi_{\mathcal B}^{\leftarrow}(\tau-T,t)   u^{\rho}_{T,(0)}(\tau,t) + \chi_{\mathcal B}^{\rightarrow}(\tau-T,t)p^{\rho}
&\text{if $z = (\tau,t) \in [-5T,5T] \times [0,1]$} \\
 u^{\rho}_{T,(0)}(z)
&\text{if $z \in K_1$} \\
p^{\rho}
&\text{if $z \in [5T,\infty)\times [0,1]$}.
\end{cases} \\
&\hat u^{\rho}_{2,T,(0)}(z) \\
&=
\begin{cases} \chi_{\mathcal A}^{\rightarrow}(\tau+T,t)   u^{\rho}_{T,(0)}(\tau,t) + \chi_{\mathcal A}^{\leftarrow}(\tau+T,t)p^{\rho}
&\text{if $z = (\tau,t) \in [-5T,5T] \times [0,1]$} \\
 u^{\rho}_{T,(0)}(z)
&\text{if $z \in K_2$} \\
p^{\rho}
&\text{if $z \in (-\infty,-5T]\times [0,1]$}.
\end{cases}
\endaligned
\end{equation}
Let

\begin{equation}\label{lineeqstep0}
\aligned
D_{\hat u^{\rho}_{i,T,(0)}}\overline{\partial}
:
L^2_{m+1,\delta}((\Sigma_i,\partial \Sigma_i);&(\hat u^{\rho}_{i,T,(0)})^*TX,(\hat u^{\rho}_{i,T,(0)})^*TL)
\\
&\to
L^2_{m,\delta}(\Sigma_i;(\hat u^{\rho}_{i,T,(0)})^{*}TX \otimes \Lambda^{01})
\endaligned
\end{equation}
be the linearization of the Cauchy-Riemann equation.
\begin{lem}\label{surjlineastep11}
We put $E_i = E_i(\hat u^{\rho}_{i,T,(0)})$.
We have
\begin{equation}\label{surj1step1}
\text{\rm Im}(D_{\hat u^{\rho}_{i,T,(0)}}\overline{\partial})  + E_i
=
L^2_{m,\delta}(\Sigma_i;(\hat u^{\rho}_{i,T,(0)})^{*}TX \otimes \Lambda^{01}).
\end{equation}
Moreover
\begin{equation}\label{Duievsurjstep1}
D{\rm ev}_{1,\infty} - D{\rm ev}_{2,\infty} : (D_{\hat u^{\rho}_{1,T,(0)}}\overline{\partial})^{-1}(E_1) \oplus (D_{\hat u^{\rho}_{2,T,(0)}}\overline{\partial})^{-1}(E_2)
\to T_{p^{\rho}}L
\end{equation}
is surjective.
\end{lem}
\begin{proof}
Since $\hat u^{\rho}_{i,T,(0)}$ is close to $u_i$ in exponential order, this is a consequence of
Assumption \ref{DuimodEi}.
\end{proof}
Note that $E_i(u'_i)$ actually depends on $u'_i$. So to obtain a linearized equation of (\ref{mainequation}) we
need to take into account of that effect.
Let $\Pi_{E_i(u_i')}$ be the projection to $E_i(u_i')$ with respect to the $L^2$ norm.
Namely
we put
\begin{equation}\label{form118}
\Pi_{E_i(u_i')}(A) = \sum_{a=1}^{\dim E_i}
\langle\!\langle A,\mathbf e_{i,a}(u'_i) \rangle\!\rangle_{L^2(K_1)}\mathbf e_{i,a}(u'_i),
\end{equation}
where $\mathbf e_{i,a}$, $a=1,\dots,\dim E_i(u_i')$ is an orthonormal basis of $E_i(u_i')$
which are supported in $K_i$.
\par
We put
\begin{equation}\label{DEidef}
(D_{u'_i}E_i)(A,v) = \frac{d}{ds}(\Pi_{E_i({\rm E} (u_i',sv))}(A))\vert_{s=0}.
\end{equation}
Here $v \in \Gamma((\Sigma_i,\partial \Sigma_i),(u'_i)^*TX,(u'_i)^*TL)$.
(Then ${\rm E} (u_i',sv)$ is a map $(\Sigma_i,\partial \Sigma_i) \to (X,L)$ defined in (\ref{defE}).)
\begin{rem}
We use an isomorphism
\begin{equation}\label{parallemap}
\Gamma(\Sigma_i ; {\rm E} (u_i',sv)^*TX\otimes \Lambda^{01})
\cong
\Gamma(\Sigma_i;(u'_i)^*TX \otimes \Lambda^{01})
\end{equation}
to define the right hand side of (\ref{DEidef}).
The map (\ref{parallemap}) is defined as follows.
Let $z \in \Sigma_i$. We have a path
$r \mapsto  {\rm E} (u_i'(z),rsv(z))$ joining $u'_i(z)$ to $ {\rm E} (u_i',sv)(z)$.
We use a connection $\nabla$ such that $TL$ is parallel to define
a parallel transport along this path. Its complex linear part
defines an isomorphism (\ref{parallemap}).
\par
We note that the same isomorphism (\ref{parallemap}) is used also to define $D_{u'_i}\overline \partial$.
Namely
$$
(D_{u'_i}\overline \partial)(v)
= \frac{d}{ds}(\overline\partial {\rm E} (u_i',sv))\vert_{s=0}
$$
where the right hand side is defined by using (\ref{parallemap}).
\end{rem}
\par
We put
$$
\Pi_{E_i(u_i')}^{\perp}(A) = A - \Pi_{E_i(u_i')}(A).
$$
\par
The equation (\ref{mainequationui}) is equivalent to the following
\begin{equation}\label{mainequationui2}
\Pi_{E_i(u_i')}^{\perp}\overline\partial u'_i
 = 0.
\end{equation}
We calculate the linearization
$$
\left.\frac{\partial}{\partial s} \Pi_{E_i({\rm E}(u_i',s V))}^{\perp} \overline\partial {\rm E}(u_i',s V))\right\vert_{s=0}
$$
to obtain the linearized equation:
\begin{equation}\label{linearized2221}
D_{u'_i}\overline\partial (V) - (D_{u'_i}E_i)(\overline\partial u'_i,V) \equiv 0
\mod E_i(u'_i).
\end{equation}
We note that
$$
\overline\partial \hat u^{\rho}_{i,T,(0)} - \frak e^{\rho} _{i,T,(0)}
$$
is exponentially small.  So we use the operator
\begin{equation}\label{144op}
V \mapsto D_{\hat u^{\rho}_{i,T,(0)}}\overline{\partial}(V) - (D_{\hat u^{\rho}_{i,T,(0)}}E_i)(\frak e^{\rho} _{i,T,(0)}, V)
\end{equation}
as an approximation of the linearization of (\ref{mainequationui2}).
\begin{lem}\label{lem112}
We put $E_i = E_i(\hat u^{\rho}_{i,T,(0)})$.
We have
\begin{equation}\label{surj1step1modififed}
\text{\rm Im}(D_{\hat u^{\rho}_{i,T,(0)}}\overline{\partial} - (D_{\hat u^{\rho}_{i,T,(0)}}E_i)(\frak e^{\rho} _{i,T,(0)}, \cdot))  + E_i
=
L^2_{m,\delta}(\Sigma_i;(\hat u^{\rho}_{i,T,(0)})^{*}TX \otimes \Lambda^{01}).
\end{equation}
Moreover
\begin{equation}\label{Duievsurjstep12}
\aligned
D{\rm ev}_{1,\infty} - &D{\rm ev}_{2,\infty} :
(D_{\hat u^{\rho}_{1,T,(0)}}\overline{\partial} - (D_{\hat u^{\rho}_{1,T,(0)}}E_1)(\frak e^{\rho} _{1,T,(0)}, \cdot))^{-1}(E_1) \\
&\oplus (D_{\hat u^{\rho}_{2,T,(0)}}\overline{\partial} -(D_{\hat u^{\rho}_{2,T,(0)}}E_2)(\frak e^{\rho} _{2,T,(0)}, \cdot))^{-1}(E_2)
\to T_{p^{\rho}}L
\endaligned
\end{equation}
is surjective.
\end{lem}
\begin{proof}
(\ref{frakeissmall}) implies that $ (D_{\hat u^{\rho}_{1,T,(0)}}E_1)(\frak e^{\rho} _{1,T,(0)}, \cdot)$ is small in operator norm.
The lemma follows from Lemma \ref{surjlineastep11}.
\end{proof}
\begin{rem}\label{independetmm}
Note that (\ref{frakeissmall}) is proved by taking $V_i$ in a small neighborhood of $0$ (in $\tilde V_i$)
with respect to the $C^m$ norm.
(Note $V_i \subset \mathcal M^{E_i}((\Sigma_i,\vec z_i);\beta_i)_{\epsilon_2}$ and $V_i$ consists of smooth maps.)
However we can take $V_i$ that is independent of $m$ and the conclusion of Lemma \ref{lem112} holds for $m$.
In fact the elliptic regularity implies that if the conclusion of Lemma \ref{lem112} holds for some $m$ then
it holds for all $m'>m$.
(The inequality (\ref{frakeissmall}) holds for that particular $m$ only. However this inequality is used to show
Lemma \ref{lem112} only.)
\end{rem}
We consider
\begin{equation}\label{gibker}
\aligned
&\text{\rm Ker}(D{\rm ev}_{1,\infty} - D{\rm ev}_{2,\infty})\\
\cap
&\left((D_{\hat u^{\rho}_{1,T,(0)}}\overline{\partial} - (D_{\hat u^{\rho}_{1,T,(0)}}E_1)(\frak e^{\rho} _{1,T,(0)}, \cdot)))^{-1}(E_1) \right.\\
&\quad\oplus
 \left.(D_{\hat u^{\rho}_{2,T,(0)}}\overline{\partial} - (D_{\hat u^{\rho}_{2,T,(0)}}E_2)(\frak e^{\rho} _{2,T,(0)}, \cdot))^{-1}(E_2)\right).
\endaligned\end{equation}
This is a finite dimensional subspace of
\begin{equation}\label{cocsissobolev}
\text{\rm Ker}(D{\rm ev}_{1,\infty} - D{\rm ev}_{2,\infty}) \cap
\bigoplus_{i=1}^2 L^2_{m+1,\delta}((\Sigma_i,\partial \Sigma_i);(\hat u^{\rho}_{i,T,(0)})^{*}TX,(\hat u^{\rho}_{i,T,(0)})^{*}TL)
\end{equation}
consisting of smooth sections.
\begin{defn}\label{defnfraH}
We denote by $\frak H(E_1,E_2)$ the intersection of the
$L^2$ orthogonal complement of (\ref{gibker}) with (\ref{cocsissobolev}).
Here the $L^2$ inner product is defined by (\ref{innerprod1}).
\end{defn}
\begin{defn}\label{def105}
We define $(V^{\rho}_{T,1,(1)},V^{\rho}_{T,2,(1)},\Delta p^{\rho}_{T,(1)})$ as follows.
\begin{equation}\label{144ffff}
\aligned
(D_{\hat u^{\rho}_{i,T,(0)}}\overline{\partial})(V^{\rho}_{T,i,(1)}) - &(D_{\hat u^{\rho}_{i,T,(0)}}E_i)(\frak e^{\rho} _{i,T,(0)},V^{\rho}_{T,i,(1)})
\\
&+ {\rm Err}^{\rho}_{i,T,(0)}
\in E_i(\hat u^{\rho}_{i,T,(0)}).
\endaligned
\end{equation}
\begin{equation}
D{\rm ev}_{\infty}(V^{\rho}_{T,1,(1)}) = D{\rm ev}_{-\infty}(V^{\rho}_{T,2,(1)}) = \Delta p^{\rho}_{T,(1)}.
\end{equation}
Moreover
$$
((V^{\rho}_{T,1,(1)},\Delta p^{\rho}_{T,(1)}),(V^{\rho}_{T,2,(1)},\Delta p^{\rho}_{T,(1)})) \in \frak H(E_1,E_2).
$$
\end{defn}
Lemma \ref{lem112} implies that such $(V^{\rho}_{T,1,(1)},V^{\rho}_{T,2,(1)},\Delta p^{\rho}_{T,(1)})$ exists and is unique.
\begin{lem}
If $\delta < \delta_1/10$, then
\begin{equation}\label{142form}
\Vert (V^{\rho}_{T,i,(1)},\Delta p^{\rho}_{T,(1)})\Vert_{L^2_{m+1,\delta}(\Sigma_i)} \le C_{2,m}e^{-\delta T},
\qquad \vert\Delta p^{\rho}_{T,(1)}\vert \le C_{2,m}e^{-\delta T}.
\end{equation}
\end{lem}
This is immediate from construction and the uniform boundedness of the right inverse of
$D_{\hat u^{\rho}_{i,T,(0)}}\overline{\partial} - (D_{\hat u^{\rho}_{i,T,(0)}}E_i)(\frak e^{\rho} _{i,T,(0)}, \cdot)$.
\par\medskip
\noindent{\bf Step 1-2}:
We use $(V^{\rho}_{T,1,(1)},V^{\rho}_{T,2,(1)},\Delta p^{\rho}_{T,(1)})$ to find an approximate solution $u^{\rho}_{T,(1)}$ of the next level.

\begin{defn}
We define $u^{\rho}_{T,(1)}(z) $ as follows.
(Here $\rm E$ is as in (\ref{defE}).)
\begin{enumerate}
\item If $z \in K_1$, we put
\begin{equation}
u^{\rho}_{T,(1)}(z)
=
{\rm E} (\hat u^{\rho}_{1,T,(0)}(z),V^{\rho}_{T,1,(1)}(z)).
\end{equation}
\item If $z \in K_2$, we put
\begin{equation}
u^{\rho}_{T,(1)}(z)
=
{\rm E} (\hat u^{\rho}_{2,T,(0)}(z),V^{\rho}_{T,2,(1)}(z)).
\end{equation}
\item
If $z  = (\tau,t) \in [-5T,5T]\times [0,1]$,
we put
\begin{equation}
\aligned
u^{\rho}_{T,(1)}(\tau,t) =
&\chi_{\mathcal B}^{\leftarrow}(\tau,t) (V^{\rho}_{T,1,(1)}(\tau,t) - \Delta p^{\rho}_{T,(1)})\\
&+\chi_{\mathcal A}^{\rightarrow}(\tau,t)(V^{\rho}_{T,2,(1)}(\tau,t)-\Delta p^{\rho}_{T,(1)})
+u^{\rho}_{T,(0)}(\tau,t) + \Delta p^{\rho}_{T,(1)}.
\endaligned
\end{equation}
\end{enumerate}
\end{defn}
We recall that $\hat u^{\rho}_{1,T,(0)}(z) = u^{\rho}_{T,(0)}(z)$ on $K_1$
and $\hat u^{\rho}_{2,T,(0)}(z) = u^{\rho}_{T,(0)}(z)$ on $K_2$.

\par\medskip
\noindent{\bf Step 1-3}:
Let $0 < \mu < 1$. We fix it throughout the proof.
\begin{lem}\label{mainestimatestep13}
There exists $\delta_2$ such that for any $\delta < \delta_2$, $T > T(\delta,m,\epsilon_{5,m})$
there exists $\frak e^{\rho} _{i,T,(1)} \in E_i$
with the following properties.
$$
\Vert\overline\partial u^{\rho} _{T,(1)} - (\frak e^{\rho} _{1,T,(0)} + \frak e^{\rho} _{1,T,(1)}) - (\frak e^{\rho} _{2,T,(0)} + \frak e^{\rho} _{2,T,(1)})\Vert_{L_{m,\delta}^2} < C_{1,m}\mu\epsilon_{5,m} e^{-\delta T}.
$$
(Here $C_{1,m}$ is the constant given in Lemma \ref{lem18}.)
Moreover
\begin{equation}\label{frakeissmall2}
\Vert \frak e^{\rho} _{i,T,(1)}\Vert_{L_{m}^2(K_i)} < C_{3,m}e^{-\delta T}.
\end{equation}
\end{lem}
\begin{proof}
The existence of $\frak e^{\rho} _{i,T,(1)}$ satisfying
$$
\Vert\overline\partial u^{\rho} _{T,(1)} - (\frak e^{\rho} _{1,T,(0)} + \frak e^{\rho} _{1,T,(1)}) - (\frak e^{\rho} _{2,T,(0)} + \frak e^{\rho} _{2,T,(1)})\Vert_{L_{m,\delta}^2(K_1\cup K_2\subset \Sigma_T)} < C_{1,m}\mu\epsilon_{5,m} e^{-\delta T}/10
$$
is a consequence of the fact that (\ref{linearized2221}) is the linearized equation of (\ref{mainequationui2}) and
the estimate (\ref{142form}). More explicitly we can prove it by a routine calculation as follows.
We first estimate on $K_1$. We have:
\begin{equation}\label{151ffff}
\aligned
\overline\partial ({\rm E} (\hat u^{\rho}_{1,T,(0)},V^{\rho}_{T,1,(1)}))& \\
=
\overline\partial ({\rm E} (\hat u^{\rho}_{1,T,(0)},0))
&+\int_0^1 \frac{\partial}{\partial s} \overline\partial ({\rm E} (\hat u^{\rho}_{1,T,(0)},sV^{\rho}_{T,1,(1)}))ds
\\
=\overline\partial ({\rm E} (\hat u^{\rho}_{1,T,(0)},0))
&+ (D_{\hat u^{\rho}_{1,T,(0)}}\overline\partial)(V^{\rho}_{T,1,(1)})
\\
&+\int_0^1 ds \int_0^s
\frac{\partial^2}{\partial r^2} \overline\partial ({\rm E} (\hat u^{\rho}_{1,T,(0)},rV^{\rho}_{T,1,(1)}))dr.
\endaligned\end{equation}
We remark
\begin{equation}\label{152ff}
\aligned
&\left\Vert\int_0^1 ds \int_0^s
\frac{\partial^2}{\partial r^2} \overline\partial ({\rm E} (\hat u^{\rho}_{1,T,(0)},rV^{\rho}_{T,1,(1)}))dr
\right\Vert_{L^2_m(K_1)}
\\
&\le C_{3,m} \Vert V^{\rho}_{T,1,(1)}\Vert_{L^2_{m+1,\delta}}^2
\le C_{4,m}e^{-2\delta T}.
\endaligned
\end{equation}
We have
\begin{equation}\label{155ff}
\aligned
\Pi_{E_1({\rm E} (\hat u^{\rho}_{1,T,(0)},V^{\rho}_{T,1,(1)}))}^{\perp}&\\
=
\Pi_{E_1(\hat u^{\rho}_{1,T,(0)})}^{\perp} &+
\int_0^1 \frac{\partial}{\partial s} \Pi_{E_i({\rm E} (\hat u^{\rho}_{1,T,(0)},sV^{\rho}_{T,1,(1)}))}^{\perp} ds
\\
=
\Pi_{E_1(\hat u^{\rho}_{1,T,(0)})}^{\perp} &-
(D_{\hat u^{\rho}_{1,T,(0)}}E_1)(\cdot,V^{\rho}_{T,1,(1)})\\
&+
\int_0^1 ds \int_0^s \frac{\partial^2}{\partial r^2} \Pi_{E_1({\rm E} (\hat u^{\rho}_{1,T,(0)},rV^{\rho}_{T,1,(1)}))}^{\perp} dr.
\endaligned\end{equation}
We can estimate the third term of the right hand side of (\ref{155ff}) in the same way as in (\ref{152ff}).
\par
On the other hand, (\ref{151ffff}) implies that
\begin{equation}\label{156ff}
\left\Vert\overline\partial ({\rm E} (\hat u^{\rho}_{1,T,(0)},V^{\rho}_{T,1,(1)}))
- \frak e^{\rho} _{1,T,(0)}\right\Vert_{L^2_{m}(K_1)}
\le C_{6,m}e^{-\delta T}.
\end{equation}
Therefore, using (\ref{155ff}) and (\ref{142form}),  we have
\begin{equation}
\aligned
&\bigg\Vert\Pi_{E_1({\rm E} (\hat u^{\rho}_{1,T,(0)},V^{\rho}_{T,1,(1)}))}^{\perp}\overline\partial ({\rm E} (\hat u^{\rho}_{1,T,(0)},V^{\rho}_{T,1,(1)}))
\\
&-
\Pi_{E_1(\hat u^{\rho}_{1,T,(0)},0)}^{\perp}\overline\partial ({\rm E} (\hat u^{\rho}_{1,T,(0)},V^{\rho}_{T,1,(1)}))\\
&-
\Pi_{E_1({\rm E} (\hat u^{\rho}_{1,T,(0)},V^{\rho}_{T,1,(1)}))}^{\perp}(\frak e^{\rho} _{1,T,(0)})
+
\Pi_{E_1(\hat u^{\rho}_{1,T,(0)},0)}^{\perp}(\frak e^{\rho} _{1,T,(0)})
\bigg\Vert_{L^2_{m}(K_1)}
\le C_{7,m}e^{-2\delta T}.
\endaligned\end{equation}
Therefore using (\ref{155ff}) we have:
\begin{equation}\label{160ff}
\aligned
&\Vert\Pi_{E_1({\rm E} (\hat u^{\rho}_{1,T,(0)},V^{\rho}_{T,1,(1)}))}^{\perp}\overline\partial ({\rm E} (\hat u^{\rho}_{1,T,(0)},V^{\rho}_{T,1,(1)}))
\\
&-
\Pi_{E_1(\hat u^{\rho}_{1,T,(0)},0)}^{\perp}\overline\partial ({\rm E} (\hat u^{\rho}_{1,T,(0)},V^{\rho}_{T,1,(1)}))\\
&+
(D_{\hat u^{\rho}_{1,T,(0)}}E_1)(\frak e^{\rho} _{1,T,(0)},V^{\rho}_{T,1,(1)})\Vert_{L^2_{m}(K_1)}
\le C_{8,m}e^{-2\delta T}.
\endaligned\end{equation}
By (\ref{144ffff}) and Definition \ref{deferfirst}, we have:
\begin{equation}\label{161ff}
\aligned
\overline\partial ({\rm E} (\hat u^{\rho}_{1,T,(0)},0))
&+ (D_{\hat u^{\rho}_{1,T,(0)}}\overline\partial)(V^{\rho}_{T,1,(1)})\\
&- (D_{\hat u^{\rho}_{1,T,(0)}}E_1)(\frak e^{\rho} _{1,T,(0)},V^{\rho}_{T,1,(1)})
\in  E_1(\hat u^{\rho}_{1,T,(0)})
\endaligned
\end{equation}
on $K_1$.
\par
(\ref{160ff}) and (\ref{161ff}) imply
\begin{equation}\label{162ff}
\aligned
&\Vert\Pi_{E_1({\rm E} (\hat u^{\rho}_{1,T,(0)},V^{\rho}_{T,1,(1)}))}^{\perp}\overline\partial ({\rm E} (\hat u^{\rho}_{1,T,(0)},V^{\rho}_{T,1,(1)}))
\\
&-
\Pi_{E_1(\hat u^{\rho}_{1,T,(0)},0)}^{\perp}\overline\partial ({\rm E} (\hat u^{\rho}_{1,T,(0)},V^{\rho}_{T,1,(1)}))\\
&+
\Pi_{E_1(\hat u^{\rho}_{1,T,(0)},0)}^{\perp}\overline\partial ({\rm E} (\hat u^{\rho}_{1,T,(0)},0))\\
&+
\Pi_{E_1(\hat u^{\rho}_{1,T,(0)},0)}^{\perp}(D_{\hat u^{\rho}_{1,T,(0)}}\overline\partial)(V^{\rho}_{T,1,(1)})\Vert_{L^2_{m}(K_1)}
\le C_{9,m}e^{-2\delta T}.
\endaligned\end{equation}
Combined with (\ref{151ffff}) and (\ref{152ff}), we have
\begin{equation}\label{eq164}
\aligned
&\Vert\Pi_{E_1({\rm E} (\hat u^{\rho}_{1,T,(0)},V^{\rho}_{T,1,(1)}))}^{\perp}(\overline\partial
({\rm E} (\hat u^{\rho}_{1,T,(0)},V^{\rho}_{T,1,(1)})))\Vert_{L^2_m(K_1)}\\
&\le C_{10,m}e^{-2\delta T}
\le C_{1,m}e^{-\delta T}\epsilon_{5,m}\mu/10,
\endaligned
\end{equation}
 for $T>T_m$ if we choose $T_m$ so that  $C_{10,m}e^{-\delta T_m} < C_{1,m}\epsilon_{5,m}\mu/10$.
\par
It follows from (\ref{156ff}) and (\ref{eq164}) that
$$
\Vert\Pi_{E_1({\rm E} (\hat u^{\rho}_{1,T,(0)},V^{\rho}_{T,1,(1)}))}(\overline\partial ({\rm E} (\hat u^{\rho}_{1,T,(0)},V^{\rho}_{T,1,(1)}))
- \frak e^{\rho} _{1,T,(0)}\Vert_{L^2_{m}(K_1)}
\le C_{11,m}e^{-\delta T}.
$$
Then (\ref{frakeissmall2}) follows,
by selecting
$$
\frak e^{\rho} _{1,T,(1)}
=
\Pi_{E_1({\rm E} (\hat u^{\rho}_{1,T,(0)},V^{\rho}_{T,1,(1)}))}(\overline\partial ({\rm E} (\hat u^{\rho}_{1,T,(0)},V^{\rho}_{T,1,(1)})
- \frak e^{\rho} _{1,T,(0)}).
$$
\par
The estimate on $K_2$ is the same.
\par
Let us estimate $\overline\partial u^{\rho} _{T,(1)}$ on $[-T+1,T-1]\times [0,1]$.
The inequality
$$
\Vert\overline\partial u^{\rho} _{T,(1)} \Vert_{L_{m,\delta}^2([-T+1,T-1]\times [0,1]\subset \Sigma_T)} < C_{1,m}\mu\epsilon_{5,m} e^{-\delta T}/10
$$
is also a consequence of the fact that (\ref{linearized2221}) is the linearized equation of (\ref{mainequationui2}) and
the estimate (\ref{142form}).
(Note the bump functions $\chi_{\mathcal B}^{\leftarrow}$ and $\chi_{\mathcal A}^{\rightarrow}$ are $\equiv 1$ there.)
On $\mathcal A_T$ we have
\begin{equation}\label{estimateatA1}
\overline\partial u^{\rho} _{T,(1)}
= \overline\partial(\chi_{\mathcal A}^{\rightarrow}(V^{\rho}_{T,2,(1)}-\Delta p^{\rho}_{T,(1)})+V^{\rho}_{T,1,(1)}
+u^{\rho}_{T,(0)}).
\end{equation}
Note
$$
\aligned
\Vert \overline\partial(\chi_{\mathcal A}^{\rightarrow}(V^{\rho}_{T,2,(1)}-\Delta p^{\rho}_{T,(1)})\Vert_{L^2_m(\mathcal A_T)}
&\le
C_{3,m}e^{-6T\delta}\Vert V^{\rho}_{T,2,(1)}-\Delta p^{\rho}_{T,(1)}\Vert_{L^2_{m+1,\delta}(\mathcal A_T \subset \Sigma_{2})}\\
&\le C_{12,m} e^{-7T\delta}.
\endaligned
$$
The first inequality follows from the fact the weight function $e_{2,\delta}$
is around $e^{6T\delta}$ on $\mathcal A_T$.
The second inequality follows from (\ref{142form}).
On the other hand the weight function $e_{T,\delta}$ is around $e^{4T\delta}$ at $\mathcal A_T$.\footnote{This drop of the weight is the
main part of the idea. It was used in \cite[page 414]{fooo:book1}. See \cite[Figure 7.1.6]{fooo:book1}.}
Therefore
\begin{equation}\label{2ff160}
\Vert \overline\partial(\chi_{\mathcal A}^{\rightarrow}(V^{\rho}_{T,2,(1)}-\Delta p^{\rho}_{T,(1)}))
\Vert_{L^2_{m,\delta}(\mathcal A_T\subset \Sigma_T)}
\le
C_{13,m} e^{-3T\delta}.
\end{equation}
Note
$$
{\rm Err}^{\rho}_{2,T,(0)}   = 0
$$
on $\mathcal A_T$. Using this in the same way as we did on $K_1$ we can show
\begin{equation}\label{2ff161}
\Vert\overline\partial(V^{\rho}_{T,1,(1)}
+u^{\rho}_{T,(0)})\Vert_{L^2_{m,\delta}(\mathcal A_T \subset \Sigma_T)}
\le C_{1,m}e^{-\delta T}\epsilon_{5,m}\mu/20
\end{equation}
for $T > T_m$.
Therefore by taking $T$ large we have
\begin{equation}\label{2ff162}
\Vert\overline\partial u^{\rho} _{T,(1)} \Vert_{L_{m,\delta}^2(\mathcal A_T\subset \Sigma_T)} < C_{1,m}\mu\epsilon_{5,m} e^{-\delta T}/10.
\end{equation}
(Note that the almost complex structure may not be integrable. So the almost complex structure may not be constant
with respect to the flat metric we are taking in the neighborhood of $p_0$. However
we can still deduce (\ref{2ff162}) from (\ref{2ff161}) and (\ref{2ff160}).)
\par
The estimate on $\mathcal B_T$ and on $([-5T,-T-1]\cup [T+1,5T]) \times [0,1]$ are similar.
The proof of Lemma \ref{mainestimatestep13} is complete.
\end{proof}
\par\medskip
\noindent{\bf Step 1-4}:
\begin{defn}
We put
$$
\aligned
{\rm Err}^{\rho}_{1,T,(1)}
&= \chi_{\mathcal X}^{\leftarrow} (\overline\partial u^{\rho} _{T,(1)} - (\frak e^{\rho} _{1,T,(0)}
+ \frak e^{\rho} _{1,T,(1)})) , \\
{\rm Err}^{\rho}_{2,T,(1)}
&= \chi_{\mathcal X}^{\rightarrow} (\overline\partial u^{\rho} _{T,(1)}  - (\frak e^{\rho} _{2,T,(0)} + \frak e^{\rho} _{2,T,(1)})).
\endaligned$$
We regard them as elements of the weighted Sobolev spaces
$L^2_{m,\delta}(\Sigma_1;(\hat u_{1,T,(1)}^{\rho})^*TX\otimes \Lambda^{01})$
and
$L^2_{m,\delta}(\Sigma_2;(\hat u_{2,T,(1)}^{\rho})^*TX\otimes \Lambda^{01})$
respectively.
(We extend them by $0$ outside a compact set.)
\end{defn}
We put $p^{\rho}_{(1)} = p^{\rho} + \Delta p^{\rho}_{T,(1)}$.
\par\medskip
We now come back to the Step 2-1 and continue.
In other words, we will prove the following by induction on $\kappa$.
\begin{eqnarray}
\left\Vert  (V^{\rho}_{T,i,(\kappa)},\Delta p^{\rho}_{T,(\kappa)})\right\Vert_{L^2_{m+1,\delta}(\Sigma_i)}
&<& C_{2,m}\mu^{\kappa-1}e^{-\delta T}, \label{form0182}
\\
\left\Vert \Delta p^{\rho}_{T,(\kappa)}\right\Vert
&<& C_{2,m}\mu^{\kappa-1}e^{-\delta T}, \label{form0183}
\\
\left\Vert  u^{\rho}_{T,(\kappa)} - u^{\rho}_{T,(0)}  \right\Vert_{L^2_{m+1,\delta}(\Sigma_{T})}
&<& C_{14,m}e^{-\delta T}, \label{form0184a}
\\
\left\Vert {\rm Err}^{\rho}_{i,T,(\kappa)} \right\Vert_{L^2_{m,\delta}(\Sigma_i)}
&<& C_{1,m}\epsilon_{5,m}\mu^{\kappa}e^{-\delta T}, \label{form0185}
\\
\left\Vert \frak e^{\rho} _{i,T,(\kappa)}\right\Vert_{L^2_{m}(K_i^{\text{\rm obst}})}
&<& C_{15,m}\mu^{\kappa-1}e^{-\delta T},
\quad \text{for $\kappa \ge 1$}. \label{form0186}
\end{eqnarray}
\begin{rem}
The left hand side of (\ref{form0184a}) is defined as follows.
We define $\frak u^{\rho}_{T,(\kappa)}$ by
$
u^{\rho}_{T,(\kappa)} = {\rm E}(u^{\rho}_{T,(\kappa-1)},\frak u^{\rho}_{T,(\kappa)}).
$
Then the left hand side of (\ref{form0184a}) is
$$
\Vert \frak u^{\rho}_{T,(\kappa)} \Vert_{L^2_{m+1,\delta}((\Sigma_T,\partial\Sigma_T);(u^{\rho}_{T,(\kappa-1)})^*TX,
(u^{\rho}_{T,(\kappa-1)})^*TL)}.
$$
\end{rem}
More precisely the claim we will prove is:
for any $\epsilon_{5,m}$ we can choose $T_m$ so that (\ref{form0182}) and (\ref{form0183}) imply (\ref{form0185}) and (\ref{form0186})
for given $T>T_m$,
and we can choose $\epsilon_{5,m}$ so that (\ref{form0185}) and (\ref{form0186})
for $\kappa$ implies (\ref{form0182}) and (\ref{form0183}) for $\kappa + 1$.
(It is easy to see that (\ref{form0182}) and (\ref{form0183}) imply (\ref{form0184a}).)
\par
Below we describe Steps $\kappa$-1,\dots,$\kappa$-4.
\par\medskip
\noindent{\bf Step $\kappa$-1}:

We first cut $u^{\rho}_{T,(\kappa-1)}$ and extend to obtain maps $\hat u^{\rho}_{i,T,(\kappa-1)} :
(\Sigma_i,\partial\Sigma_i) \to (X,L)$ $(i=1,2)$ as follows.
\begin{equation}
\aligned
&\hat u^{\rho}_{1,T,(\kappa-1)}(z) \\
&=
\begin{cases} \chi_{\mathcal B}^{\leftarrow}(\tau-T,t)   u^{\rho}_{T,(\kappa-1)}(\tau,t) + \chi_{\mathcal B}^{\rightarrow}(\tau-T,t)p^{\rho}_{(\kappa-1)}
&\text{if $z = (\tau,t) \in [-5T,5T] \times [0,1]$} \\
 u^{\rho}_{T,(\kappa-1)}(z)
&\text{if $z \in K_1$} \\
p^{\rho}_{T,(\kappa-1)}
&\text{if $z \in [5T,\infty)\times [0,1]$}.
\end{cases} \\
&\hat u^{\rho}_{2,T,(\kappa-1)}(z) \\
&=
\begin{cases} \chi_{\mathcal A}^{\rightarrow}(\tau+T,t)   u^{\rho}_{T,(\kappa-1)}(\tau,t) + \chi_{\mathcal A}^{\leftarrow}(\tau+T,t)p^{\rho}_{(\kappa-1)}
&\text{if $z = (\tau,t) \in [-5T,5T] \times [0,1]$} \\
 u^{\rho}_{T,(\kappa-1)}(z)
&\text{if $z \in K_2$} \\
p^{\rho}_{T,(\kappa-1)}
&\text{if $z \in (-\infty,-5T]\times [0,1]$}.
\end{cases}
\endaligned
\end{equation}
Let

\begin{equation}\label{lineeqstepkappa}
\aligned
D_{\hat u^{\rho}_{i,T,(\kappa-1)}}\overline{\partial}
:
L^2_{m+1,\delta}((\Sigma_i,\partial \Sigma_i);&(\hat u^{\rho}_{i,T,(\kappa-1)})^*TX,(\hat u^{\rho}_{i,T,(\kappa-1)})^*TL)
\\
&\to
L^2_{m,\delta}(\Sigma_i;(\hat u^{\rho}_{i,T,(\kappa-1)})^{*}TX \otimes \Lambda^{01}).
\endaligned
\end{equation}
\begin{lem}\label{surjlineastep11kappa}
We have
\begin{equation}\label{surj1step1}
\text{\rm Im}(D_{\hat u^{\rho}_{i,T,(\kappa-1)}}\overline{\partial})  + E_i
=
L^2_{m,\delta}(\Sigma_i;(\hat u^{\rho}_{i,T,(\kappa-1)})^{*}TX \otimes \Lambda^{01}).
\end{equation}
Moreover
\begin{equation}\label{Duievsurjstep1}
D{\rm ev}_{1,\infty} - D{\rm ev}_{2,\infty} : (D_{\hat u^{\rho}_{1,T,(0)}}\overline{\partial})^{-1}(E_1) \oplus (D_{\hat u^{\rho}_{2,T,(0)}}\overline{\partial})^{-1}(E_2)
\to T_{p^{\rho}_{T,(\kappa-1)}}L
\end{equation}
is surjective.
\end{lem}
\begin{proof}
Since $\hat u^{\rho}_{i,T,(\kappa-1)}$ is close to $u_i$ in exponential order, this is a consequence of
Assumption \ref{DuimodEi}.
\end{proof}
We denote
\begin{equation}
(\frak {se})^{\rho} _{i,T,(\kappa-1)} = \sum_{a=0}^{\kappa-1}\frak e^{\rho} _{i,T,(a)}.
\end{equation}
\begin{lem}\label{lem112kappa}
We have
\begin{equation}\label{surj1step1modififedkappa}
\aligned
&\text{\rm Im}(D_{\hat u^{\rho}_{i,T,(\kappa-1)}}\overline{\partial} -
(D_{\hat u^{\rho}_{i,T,(\kappa-1)}}E_i)((\frak {se})^{\rho} _{i,T,(\kappa-1)}, \cdot))  + E_i \\
&=
L^2_{m,\delta}(\Sigma_i;(\hat u^{\rho}_{i,T,(\kappa-1)})^{*}TX \otimes \Lambda^{01}).
\endaligned\end{equation}
Moreover
\begin{equation}\label{Duievsurjstep1kappa2}
\aligned
&D{\rm ev}_{1,\infty} - D{\rm ev}_{2,\infty} \\
: &(D_{\hat u^{\rho}_{1,T,(\kappa-1)}}\overline{\partial} -
(D_{\hat u^{\rho}_{1,T,(\kappa-1)}}E_1)((\frak {se})^{\rho} _{1,T,(\kappa-1)}, \cdot)))^{-1}(E_1)
\\
&\oplus (D_{\hat u^{\rho}_{2,T,(\kappa-1)}}\overline{\partial} -
(D_{\hat u^{\rho}_{2,T,(\kappa-1)}}E_2)((\frak {se})^{\rho} _{2,T,(\kappa-1)}, \cdot))^{-1}(E_2)
\to T_{p^{\rho}_{T,(\kappa-1)}}L
\endaligned
\end{equation}
is surjective.
\end{lem}
\begin{proof}
\begin{equation}\label{eq181}
\left\Vert \sum_{a=0}^{\kappa-1}\frak e^{\rho} _{i,T,(a)}\right\Vert_{L_{m}^2(K_i)} < \epsilon_{4,m} + C_{15,m}\frac{e^{-\delta T}}{1-\mu}.
\end{equation}
imply that $ (D_{\hat u^{\rho}_{1,T,(0)}}E_1)(\frak e^{\rho} _{1,T,(0)}, \cdot)$ is small in operator norm.
The lemma follows from Lemma \ref{surjlineastep11kappa}.
\end{proof}
Note that Remark \ref{independetmm} still applies to Lemma \ref{lem112kappa}.
\begin{defn}
We define $(V^{\rho}_{T,1,(\kappa)},V^{\rho}_{T,2,(\kappa)},\Delta p^{\rho}_{T,(\kappa)})$ as follows.
\begin{equation}\label{formula158}
\aligned
D_{\hat u^{\rho}_{i,T,(\kappa-1)}}(V^{\rho}_{T,i,(\kappa)}) - (D_{\hat u^{\rho}_{i,T,(\kappa-1)}}E_i)
&((\frak {se})^{\rho} _{i,T,(\kappa-1)}, V^{\rho}_{T,i,(\kappa)}) \\
&+ {\rm Err}^{\rho}_{i,T,(\kappa-1)}
\in E_i(\hat u^{\rho}_{i,T,(\kappa-1)}).
\endaligned
\end{equation}
\begin{equation}
D{\rm ev}_{1,\infty}(V^{\rho}_{T,1,(\kappa)}) = D{\rm ev}_{2,\infty}(V^{\rho}_{T,2,(\kappa)}) = \Delta p^{\rho}_{T,(\kappa)}.
\end{equation}
We also require
\begin{equation}\label{inHhanru}
((V^{\rho}_{T,1,(\kappa)},\Delta p^{\rho}_{T,(\kappa)}),(V^{\rho}_{T,2,(\kappa)},\Delta p^{\rho}_{T,(\kappa)}))
\in
\frak H(E_1,E_2).
\end{equation}
\end{defn}
Lemma \ref{lem112kappa} implies that such $(V^{\rho}_{T,1,(\kappa)},V^{\rho}_{T,2,(\kappa)},\Delta p^{\rho}_{T,(\kappa)})$ exists and is unique.
\begin{rem}
Note in (\ref{inHhanru}) we use the same space $\frak H(E_1,E_2)$ as in Definition \ref{def105}.
We may use the orthogonal complement of
$$
\text{\rm Ker}(D{\rm ev}_{1,\infty} - D{\rm ev}_{2,\infty}) \cap
\bigoplus_{i=1}^2 (D_{\hat u^{\rho}_{i,T,(\kappa-1)}}\overline{\partial}
- (D_{\hat u^{\rho}_{i,T,(\kappa-1)}}E_i)((\frak {se})^{\rho} _{i,T,(\kappa-1)}, \cdot))^{-1}(E_i)
$$
instead.  The reason why we use the same space as one in Definition \ref{def105}
here is that then a calculation we need to do for the exponential decay estimate of $T$ derivative
becomes a bit shorter.
Since $\hat u^{\rho}_{i,T,(\kappa)}$ is sufficiently close to $\hat u^{\rho}_{i,T,(0)}$, the unique existence of
$(V^{\rho}_{T,1,(\kappa)},V^{\rho}_{T,2,(\kappa)},\Delta p^{\rho}_{T,(\kappa)})$ satisfying (\ref{formula158}) - (\ref{inHhanru}) holds
by (\ref{eq181}).
\end{rem}
\begin{lem}\label{estimageVkappa}
If $\delta < \delta_1/10$ and $T>T(\delta,m)$, then
\begin{equation}\label{142form23}
\aligned
&\Vert (V_{T,i,(\kappa)},\Delta p^{\rho}_{T,(\kappa)})\Vert_{L^2_{m+1,\delta}(\Sigma_i)} \le C_{2,m}\mu^{\kappa-1}e^{-\delta T}, \\
&\vert\Delta p^{\rho}_{T,(\kappa)}\vert \le C_{2,m}\mu^{\kappa-1}e^{-\delta T}.
\endaligned
\end{equation}
\end{lem}
\begin{proof}
This follows from uniform boundedness of the inverse of (\ref{surj1step1modififedkappa}) together with
the $\kappa-1$ version of Lemma \ref{mainestimatestep13}.
(That is Lemma \ref{mainestimatestep13kappa}.)
\end{proof}
This lemma implies (\ref{form0182}) and (\ref{form0183}).
\par\medskip
\noindent{\bf Step $\kappa$-2}:
We use $(V^{\rho}_{T,1,(\kappa)},V^{\rho}_{T,2,(\kappa)},\Delta p^{\rho}_{T,(\kappa)})$ to find an approximate solution $u^{\rho}_{T,(\kappa)}$ of the next level.
\begin{defn}
We define $u^{\rho}_{T,(\kappa)}(z) $ as follows.
\begin{enumerate}
\item If $z \in K_1$, we put
\begin{equation}
u^{\rho}_{T,(\kappa)}(z)
=
{\rm E} (\hat u^{\rho}_{1,T,(\kappa-1)}(z),V^{\rho}_{T,1,(\kappa)}(z)).
\end{equation}
\item If $z \in K_2$, we put
\begin{equation}
u^{\rho}_{T,(\kappa)}(z)
=
{\rm E} (\hat u^{\rho}_{2,T,(\kappa-1)}(z),V^{\rho}_{T,2,(\kappa)}(z)).
\end{equation}
\item
If $z  = (\tau,t) \in [-5T,5T]\times [0,1]$,
we put
\begin{equation}
\aligned
u^{\rho}_{T,(\kappa)}(\tau,t) =
&\chi_{\mathcal B}^{\leftarrow}(\tau,t) (V^{\rho}_{T,1,(\kappa)}(\tau,t) - \Delta p^{\rho}_{T,(\kappa)})\\
&+\chi_{\mathcal A}^{\rightarrow}(\tau,t)(V^{\rho}_{T,2,(\kappa)}(\tau,t)-\Delta p^{\rho}_{T,(\kappa)})\\
&+u^{\rho}_{T,(\kappa-1)}(\tau,t)+\Delta p^{\rho}_{T,(\kappa)}.
\endaligned
\end{equation}
\end{enumerate}
\end{defn}
We note that $\hat u^{\rho}_{1,T,(\kappa-1)}(z) = u^{\rho}_{T,(\kappa-1)}(z)$ on $K_1$
and $\hat u^{\rho}_{2,T,(\kappa-1)}(z) = u^{\rho}_{T,(\kappa-1)}(z)$
on $K_2$.
\par
 (\ref{form0184a}) is immediate from the definition and (\ref{form0182}) and (\ref{form0183}), since $0<\mu<1$.
\par\medskip
\noindent{\bf Step $\kappa$-3}:
\begin{lem}\label{mainestimatestep13kappa}
For each $\epsilon_5>0$ we have the following.
If $\delta < \delta_2$ and $T > T(\delta,m,\epsilon_5)$, then there exists $\frak e^{\rho} _{i,T,(\kappa)} \in E_i$
such that
$$
\left\Vert\overline\partial u^{\rho} _{T,(\kappa)} - \sum_{a=0}^{\kappa}\frak e^{\rho} _{1,T,(a)} -\sum_{a=0}^{\kappa}\frak e^{\rho} _{2,T,(a)}
\right\Vert_{L_{m,\delta}^2} < C_{1,m}\mu^{\kappa}\epsilon_5 e^{-\delta T}.
$$
(Here $C_{1,m}$ is as in Lemma \ref{lem18}.)
Moreover
\begin{equation}\label{frakeissmall2kappa}
\Vert \frak e^{\rho} _{i,T,(\kappa)}\Vert_{L_{m}^2(K_i)} <C_{15,m}\mu^{\kappa-1}e^{-\delta T}.
\end{equation}
\end{lem}
\begin{proof}
The proof is  similar to the proof of Lemma \ref{mainestimatestep13} and proceed as follows.
We have:
\begin{equation}\label{2151ffff}
\aligned
\overline\partial ({\rm E} (\hat u^{\rho}_{1,T,(\kappa-1)},V^{\rho}_{T,1,(\kappa)}))& \\
=
\overline\partial ({\rm E} (\hat u^{\rho}_{1,T,(\kappa-1)},0))
&+\int_0^1 \frac{\partial}{\partial s} \overline\partial ({\rm E} (\hat u^{\rho}_{1,T,(\kappa-1)},sV^{\rho}_{T,1,(\kappa)}))ds
\\
=\overline\partial ({\rm E} (\hat u^{\rho}_{1,T,(\kappa-1)},0))
&+ (D_{\hat u^{\rho}_{1,T,(\kappa-1)}}\overline\partial)(V^{\rho}_{T,1,(\kappa)})
\\
&+\int_{0}^1 ds\int_0^s
\frac{\partial^2}{\partial r^2} \overline\partial ({\rm E} (\hat u^{\rho}_{1,T,(\kappa-1)},rV^{\rho}_{T,1,(\kappa)}))dr.
\endaligned\end{equation}
We remark
\begin{equation}\label{2152ff}
\aligned
&\left\Vert \int_{0}^1 ds\int_0^s
\frac{\partial^2}{\partial r^2} \overline\partial ({\rm E} (\hat u^{\rho}_{1,T,(\kappa-1)},rV^{\rho}_{T,1,(\kappa)}))dr
\right\Vert_{L^2_m(K_1)}
\\
&\le C_{4,m} \Vert V^{\rho}_{T,1,(\kappa)}\Vert_{L^2_{m+1,\delta}}^2
\le C_{5,m}e^{-2\delta T}\mu^{2(\kappa-1)}.
\endaligned
\end{equation}
We have
\begin{equation}\label{2155ff}
\aligned
\Pi_{E_1({\rm E} (\hat u^{\rho}_{1,T,(\kappa-1)},V^{\rho}_{T,1,(\kappa)}))}^{\perp}&\\
=
\Pi_{E_1(\hat u^{\rho}_{1,T,(\kappa-1)})}^{\perp} &+
\int_0^1 \frac{\partial}{\partial s} \Pi_{E_i({\rm E} (\hat u^{\rho}_{1,T,(\kappa-1)},sV^{\rho}_{T,1,(\kappa)}))}^{\perp} ds
\\
=
\Pi_{E_1(\hat u^{\rho}_{1,T,(\kappa-1)})}^{\perp} &-
(D_{\hat u^{\rho}_{1,T,(\kappa-1)}}E_1)(\cdot,V^{\rho}_{T,1,(\kappa)})\\
&+
\int_{0}^1 ds\int_0^s  \frac{\partial^2}{\partial r^2} \Pi_{E_1({\rm E} (\hat u^{\rho}_{1,T,(\kappa-1)},rV^{\rho}_{T,1,(\kappa)}))}^{\perp} dr.
\endaligned\end{equation}
We can estimate the third term of the right hand side of (\ref{2155ff}) in the same way as (\ref{2152ff}).
\par
On the other hand, (\ref{2151ffff}) implies that
\begin{equation}\label{2156ff}
\left\Vert\overline\partial ({\rm E} (\hat u^{\rho}_{1,T,(\kappa-1)},V^{\rho}_{T,1,(\kappa)}))
- \frak{se}^{\rho} _{1,T,(\kappa-1)}\right\Vert_{L^2_{m}(K_1)}
\le C_{6,m}e^{-\delta T}\mu^{\kappa-1}.
\end{equation}
Therefore
\begin{equation}
\aligned
&\Vert\Pi_{E_1({\rm E} (\hat u^{\rho}_{1,T,(\kappa-1)},V^{\rho}_{T,1,(\kappa)}))}^{\perp}
\overline\partial ({\rm E} (\hat u^{\rho}_{1,T,(\kappa-1)},V^{\rho}_{T,1,(\kappa)}))
\\
&-
\Pi_{E_1(\hat u^{\rho}_{1,T,(\kappa-1)},0)}^{\perp}\overline\partial ({\rm E} (\hat u^{\rho}_{1,T,(\kappa-1)},V^{\rho}_{T,1,(\kappa)}))\\
&+
(D_{\hat u^{\rho}_{1,T,(\kappa-1)}}E_1)(\frak{se}^{\rho} _{1,T,(\kappa-1)},V^{\rho}_{T,1,(\kappa)})\Vert_{L^2_{m}(K_1)}
\le C_{7,m}e^{-2\delta T}\mu^{\kappa-1}.
\endaligned\end{equation}
By (\ref{formula158}) we have:
\begin{equation}
\aligned
\overline\partial ({\rm E} (\hat u^{\rho}_{1,T,(\kappa-1)},0))
&+ (D_{\hat u^{\rho}_{1,T,(\kappa-1)}}\overline\partial)(V^{\rho}_{T,1,(\kappa)})\\
&- (D_{\hat u^{\rho}_{1,T,(\kappa-1)}}E_1)(\frak{se}^{\rho} _{1,T,(\kappa-1)},V^{\rho}_{T,1,(\kappa)})
\in  E_1(\hat u^{\rho}_{1,T,(\kappa-1)})
\endaligned
\end{equation}
on $K_1$.
\par
Summing up we have
\begin{equation}
\aligned
&\Vert\Pi_{E_1({\rm E} (\hat u^{\rho}_{1,T,(\kappa-1)},V^{\rho}_{T,1,(\kappa)}))}^{\perp}(\overline\partial
({\rm E} (\hat u^{\rho}_{1,T,(\kappa-1)},V^{\rho}_{T,1,(\kappa)})))\Vert_{L^2_m(K_1)}\\
&\le C_{10,m}e^{-2\delta T}\mu^{\kappa-1}
\le C_{1,m}e^{-\delta T}\epsilon_{5,m}\mu^{\kappa}/10
\endaligned
\end{equation}
for $T>T_m$.
\par
It follows from (\ref{2156ff}) that
$$
\Vert\Pi_{E_1({\rm E} (\hat u^{\rho}_{1,T,(\kappa-1)},V^{\rho}_{T,1,(\kappa)}))}(\overline\partial ({\rm E} (\hat u^{\rho}_{1,T,(\kappa-1)},V^{\rho}_{T,1,(\kappa)}))
- \frak{se}^{\rho} _{1,T,(\kappa-1)}\Vert_{L^2_{m}(K_1)}
\le C_{8,m}e^{-\delta T}\mu^{\kappa-1}.
$$
Then (\ref{frakeissmall2kappa}) follows by putting
$$
\aligned
\frak e^{\rho}_{1,T,(\kappa)}
&=
\Pi_{E_1({\rm E} (\hat u^{\rho}_{1,T,(\kappa-1)},V^{\rho}_{T,1,(\kappa)}))}(\overline\partial ({\rm E} (\hat u^{\rho}_{1,T,(\kappa-1)},V^{\rho}_{T,1,(\kappa)}))
- \frak{se}^{\rho} _{1,T,(\kappa-1)}\\
&\in E_1({\rm E} (\hat u^{\rho}_{1,T,(\kappa-1)},V^{\rho}_{T,1,(\kappa)}))
\cong E_1.
\endaligned
$$
\par
Let us estimate $\overline\partial u^{\rho} _{T,(\kappa)}$ on $[-T,T]\times [0,1]$.
The inequality
$$
\Vert\overline\partial u^{\rho} _{T,(\kappa)} \Vert_{L_{m,\delta}^2([-T,T]\times [0,1]\subset \Sigma_T)} < C_{1,m}
\mu^{\kappa}\epsilon_{5,m} e^{-\delta T}/10
$$
is also a consequence of the fact that (\ref{linearized2221}) is the linearized equation of (\ref{mainequationui2}) and
the estimate (\ref{142form23}).
(Note the bump functions $\chi_{\mathcal B}^{\leftarrow}$ and $\chi_{\mathcal A}^{\rightarrow}$ are $\equiv 1$ there.)
On $\mathcal A_T$ we have
\begin{equation}\label{2estimateatA1}
\overline\partial u^{\rho} _{T,(\kappa)}
= \overline\partial(\chi_{\mathcal A}^{\rightarrow}(V^{\rho}_{T,2,(\kappa)}-\Delta p^{\rho}_{T,(\kappa)})+V^{\rho}_{T,1,(\kappa)}
+u^{\rho}_{T,(\kappa-1)}).
\end{equation}
Note
$$
\aligned
\Vert \overline\partial(\chi_{\mathcal A}^{\rightarrow}(V^{\rho}_{T,2,(\kappa)}-\Delta p^{\rho}_{T,(\kappa)}))\Vert_{L^2_m(\mathcal A_T)}
&\le
C_{3,m}e^{-6T\delta}\Vert V^{\rho}_{T,2,(\kappa)}-\Delta p^{\rho}_{T,(\kappa)}\Vert_{L^2_{m+1,\delta}(\mathcal A_T \subset \Sigma_{2})}\\
&\le C_{12,m} e^{-7T\delta}\mu^{\kappa-1}.
\endaligned
$$
The first inequality follows from  the fact the weight function $e_{2,\delta}$
is around $e^{6T\delta}$ on $\mathcal A_T$. The second inequality follows from (\ref{142form23}).
On the other hand the weight function $e_{T,\delta}$ is around $e^{4T\delta}$ at $\mathcal A_T$.\footnote{This drop of the weight is the
main part of the idea. It was used in \cite[page 414]{fooo:book1}. See \cite[Figure 7.1.6]{fooo:book1}.}
Therefore
\begin{equation}\label{ff160}
\Vert \overline\partial(\chi_{\mathcal A}^{\rightarrow}(V^{\rho}_{T,2,(\kappa)}-\Delta p^{\rho}_{T,(\kappa)}))
\Vert_{L^2_{m,\delta}(\mathcal A_T\subset \Sigma_T)}
\le
C_{13,m} e^{-3T\delta}\mu^{\kappa-1}.
\end{equation}
Note
$$
{\rm Err}^{\rho}_{2,T,(\kappa-1)}   = 0
$$
on $\mathcal A_T$. Therefore in the same way as we did on $K_1$ we can show
\begin{equation}\label{ff161}
\Vert\overline\partial(V^{\rho}_{T,1,(\kappa)}
+u^{\rho}_{T,(\kappa-1)})\Vert_{L^2_{m,\delta}(\mathcal A_T \subset \Sigma_T)}
\le C_{1,m}e^{-\delta T}\epsilon_{5,m}\mu^{\kappa}/20
\end{equation}
for $T > T_m$.
Therefore by taking $T$ large we have
\begin{equation}\label{ff162}
\Vert\overline\partial u^{\rho} _{T,(\kappa)} \Vert_{L_{m,\delta}^2(\mathcal A_T\subset \Sigma_T)} < C_{1,m}\mu^{\kappa}\epsilon_{5,m} e^{-\delta T}/10.
\end{equation}
\par
The estimate on $\mathcal B_T$ and on $([-5T,-T-1]\cup [T+1,5T]) \times [0,1]$ are similar.
The proof of Lemma \ref{mainestimatestep13kappa} is complete.
\end{proof}
\par\medskip
\noindent{\bf Step $\kappa$-4}:
\begin{defn}
We put
$$
\aligned
{\rm Err}^{\rho}_{1,T,(\kappa)}
&= \chi_{\mathcal X}^{\leftarrow} \left(\overline\partial u^{\rho} _{T,(\kappa)} - \sum_{a=0}^{\kappa}\frak e^{\rho} _{1,T,(a)}  \right), \\
{\rm Err}^{\rho}_{2,T,(\kappa)}
&= \chi_{\mathcal X}^{\rightarrow} \left(\overline\partial u^{\rho} _{T,(\kappa)}
- \sum_{a=0}^{\kappa}\frak e^{\rho} _{2,T,(a)}\right).
\endaligned$$
We regard them as elements of the weighted Sobolev spaces
$L^2_{m,\delta}(\Sigma_1;(\hat u_{1,T,(\kappa)}^{\rho})^*TX\otimes \Lambda^{01})$
and
$L^2_{m,\delta}(\Sigma_2;(\hat u_{2,T,(\kappa)}^{\rho})^*TX\otimes \Lambda^{01})$
respectively.
(We extend them by $0$ outside a compact set.)
\end{defn}
We put $p^{\rho}_{(\kappa)} = p^{\rho}_{(\kappa-1)} + \Delta p^{\rho}_{T,(\kappa)}$.
\par
Lemma \ref{mainestimatestep13kappa} implies (\ref{form0185}) and (\ref{form0186}).
\par
\medskip
We have thus described all the induction steps.
For each fixed $m$ there exists $T_m$ such that if $T > T_m$ then
$$
\lim_{\kappa \to \infty} u^{\rho} _{T,(\kappa)}
$$
converges in $L_{m+1,\delta}^2$ sense to the solution of (\ref{mainequation}).
The limit is automatically of $C^{\infty}$ class by elliptic regularity.
We have thus constructed the map in Theorem \ref{gluethm1}.
We will prove its surjectivity and injectivity in Section \ref{surjinj} below.
Before doing so we prove an exponential decay estimate of its $T$ derivative.

\section{Exponential decay of $T$ derivatives}
\label{subsecdecayT}

We first state the result of this subsection.
We recall that for $T$ sufficiently large and $\rho = (\rho_1,\rho_2) \in V_1\times_L V_2$
we have defined $u^{\rho} _{T,(\kappa)}$. We denote its limit by
\begin{equation}\label{limitkappainf}
u^{\rho} _{T} = \lim_{\kappa \to \infty} u^{\rho} _{T,(\kappa)} : (\Sigma_T,\partial\Sigma_T) \to (X,L).
\end{equation}
The main result of this subsection
is an estimate of $T$ and $\rho$ derivatives of this map.
We prepare some notations to state the result.
\par
We change the coordinates of $\Sigma_i$ and $\Sigma_T$ as follows.
In the last section we put
$$
\Sigma_1 = K_1 \cup ([-5T,\infty) \times [0,1])
$$
and use $(\tau,t)$ for the coordinate of $[-5T,\infty) \times [0,1]$.
This identification depends on $T$. So we rewrite it to
$$
\Sigma_1 = K_1 \cup ([0,\infty) \times [0,1])
$$
and the coordinate for $[0,\infty) \times [0,1]$ is $(\tau',t)$ where
\begin{equation}
\tau' =\tau + 5T.
\end{equation}
Similarly we rewrite
$$
\Sigma_2 = ((-\infty,5T] \times [0,1]) \cup K_2
$$
to
$$
\Sigma_2 = ((-\infty,0] \times [0,1]) \cup K_2
$$
and use the coordinate $(\tau'',t)$ where
\begin{equation}
\tau'' =\tau - 5T.
\end{equation}
We may use either $(\tau',t)$ or $(\tau'',t)$ as the coordinate of
$
\Sigma_T \setminus (K_1\cup K_2)
$.
\par
Let $S$ be a positive number.
We have
$
K_i \subset \Sigma_T.
$
We put
\begin{equation}\label{1104}
\aligned
K_1^{+S} &= K_1 \cup ([0,S]\times [0,1]) \subset \Sigma_T, \\
K_2^{+S} &= ([-S,0]\times [0,1]) \cup K_2   \subset \Sigma_T.
\endaligned
\end{equation}
Here the inclusion $K_1 \cup ([0,S]\times [0,1]) \subset \Sigma_T$ is by using the coordinate
$\tau'$ and the inclusion $([-S,0]\times [0,1]) \cup K_2   \subset \Sigma_T$
is by using the coordinate $\tau''$.
\par
We may also regard $K_i^{+S} \subset \Sigma_i$.
Note that the spaces $K_i^{+S}$ are independent of $T$, as far as $10T>S$.

\par
We restrict the map $u^{\rho} _{T}$ to $K_i^{+S}$.
We thus obtain a map
$$
\text{\rm Glures}_{i,S} : [T_m,\infty) \times V_1 \times_L V_2 \to \text{\rm Map}_{L^2_{m+1}}((K_i^{+S},K_i^{+S}\cap\partial \Sigma_i),(X,L))
$$
by
\begin{equation}
\begin{cases}
\text{\rm Glures}_{1,S}(T,\rho)(x)  &=  u^{\rho} _{T}(x) \qquad x \in K_1\\
\text{\rm Glures}_{1,S}(T,\rho)(\tau',t) &= u^{\rho} _{T}(\tau',t) = u^{\rho}_T(\tau+5T,t)
\end{cases}
\end{equation}
\begin{equation}
\begin{cases}
\text{\rm Glures}_{2,S}(T,\rho)(x)  &=  u^{\rho} _{T}(x) \qquad x \in K_2\\
\text{\rm Glures}_{2,S}(T,\rho)(\tau'',t) &= u^{\rho} _{T}(\tau'',t) = u^{\rho}_T(\tau-5T,t)
\end{cases}
\end{equation}
Here $\text{\rm Map}_{L^2_{m+1}}((K_i^{+S},K_i^{+S}\cap\partial \Sigma_i),(X,L))$ is the space of maps of $L^2_{m+1}$ class
($m$ is sufficiently large, say $m>10$.) It has a structure of Hilbert manifold in an obvious way.
This Hilbert manifold is independent of $T$. So we can define $T$ derivative of a family of elements
of $\text{\rm Map}_{L^2_{m+1}}((K_i^{+S},K_i^{+S}\cap\partial \Sigma_i),(X,L))$ parametrized by $T$.
\begin{rem}\label{differentm}
The domain and the target of the map $\text{\rm Glures}_{i,S}$ depend on $m$.
However its image actually is in the set of smooth maps. Also
none of the constructions of $u^{\rho} _{T}$ depends on $m$.
(The proof of the convergence of (\ref{limitkappainf}) depends on $m$. So the number $T_m$ depends on $m$.)
Therefore the map $\text{\rm Glures}_{i,S}$ is {\it independent} of $m$ on the intersection of the domains.
Namely the map $\text{\rm Glures}_{i,S}$ constructed by using $L^2_{m_1}$ norm
coincides with  the map $\text{\rm Glures}_{i,S}$ constructed by using $L^2_{m_2}$ norm
on $[\max\{T_{m_1},T_{m_2}\},\infty) \times V_1\times_L V_2$.
\end{rem}

\begin{thm}\label{exdecayT}
For each $m$ and $S$ there exist $T(m), C_{16,m,S},
\delta > 0$ such that the following holds
for $T>T(m)$ and $n + \ell \le m - 10$ and $\ell > 0$.
\begin{equation}
\left\Vert \nabla_{\rho}^n \frac{d^{\ell}}{dT^{\ell}} \text{\rm Glures}_{i,S}\right\Vert_{L^2_{m+1-\ell}}
< C_{16,m,S}e^{-\delta T}.
\end{equation}
Here $\nabla_{\rho}^n$ is the $n$-th derivative in $\rho$ direction.
\end{thm}
\begin{rem}
Theorem \ref{exdecayT} is basically equivalent to \cite[Lemma A1.58]{fooo:book1}.
The proof below is basically the same as the one in \cite[page 776]{fooo:book1}.
We add some more detail.
\end{rem}
\begin{proof}

\par
The construction of $u_{T,(\kappa)}^{\rho}$ was by induction on $\kappa$.
We divide the inductive step of the construction of  $u_{T,(\kappa+1)}^{\rho}$
from  $u_{T,(\kappa)}^{\rho}$ into two.
\begin{enumerate}
\item[(Part A)]
Start from  $(V^{\rho}_{T,1,(\kappa)},V^{\rho}_{T,2,(\kappa)},\Delta p^{\rho}_{T,(\kappa)})$ and end with
${\rm Err}^{\rho}_{1,T,(\kappa)}$ and ${\rm Err}^{\rho}_{2,T,(\kappa)}$.
This is step $\kappa$-2,$\kappa$-3,$\kappa$-4.
\item[(Part B)]
Start from
${\rm Err}^{\rho}_{1,T,(\kappa)}$ and ${\rm Err}^{\rho}_{2,T,(\kappa)}$ and end with
$(V^{\rho}_{T,1,(\kappa+1)},V^{\rho}_{T,2,(\kappa+1)},\Delta p^{\rho}_{T,(\kappa+1)})$.
This is step $(\kappa+1)$-1.
\end{enumerate}
\par\medskip
We will prove the following inequality by induction on $\kappa$, under the assumption
$T > T(m)$, $\ell >0$, $n+\ell \le m-10$.
\begin{eqnarray}
\left\Vert \nabla_{\rho}^n \frac{\partial^{\ell}}{\partial T^{\ell}} (V^{\rho}_{T,i,(\kappa)},\Delta p^{\rho}_{T,(\kappa)})\right\Vert_{L^2_{m+1-\ell,\delta}(\Sigma_i)}
&<& C_{17,m}\mu^{\kappa-1}e^{-\delta T}, \label{form182}
\\
\left\Vert \nabla_{\rho}^n \frac{\partial^{\ell}}{\partial T^{\ell}} \Delta p^{\rho}_{T,(\kappa)}\right\Vert
&<& C_{17,m}\mu^{\kappa-1}e^{-\delta T}, \label{form183}
\\
\left\Vert \nabla_{\rho}^n \frac{\partial^{\ell}}{\partial T^{\ell}} u^{\rho}_{T,(\kappa)}\right\Vert_{L^2_{m+1-\ell,\delta}(K_i^{+5T+1})}
&<& C_{18,m}e^{-\delta T}, \label{form184}
\\
\left\Vert \nabla_{\rho}^n \frac{\partial^{\ell}}{\partial T^{\ell}}{\rm Err}^{\rho}_{i,T,(\kappa)} \right\Vert_{L^2_{m-\ell,\delta}(\Sigma_i)}
&<& C_{19,m}\epsilon_{6,m}\mu^{\kappa}e^{-\delta T}, \label{form185}
\\
\left\Vert \nabla_{\rho}^n \frac{\partial^{\ell}}{\partial T^{\ell}} \frak e^{\rho} _{i,T,(\kappa)}\right\Vert_{L^2_{m-\ell}(K_i^{\text{\rm obst}})}
&<& C_{19,m}
\mu^{\kappa-1}e^{-\delta T}. \label{form186}
\end{eqnarray}
More precisely, the claim we will prove is the following: For each $\epsilon_{6,m}$,
we can choose $T(m)$ so that (\ref{form182}) and (\ref{form183}) imply (\ref{form185}) and (\ref{form186})
for $T > T(m)$,
and we can choose $\epsilon_{6,m}$
so that (\ref{form185}) and (\ref{form186})
for $\kappa$ implies (\ref{form182}) and (\ref{form183}) for $\kappa + 1$.
(\ref{form184}) follows from (\ref{form182}) and (\ref{form183}).
\begin{rem}
We use $L^2_{m+1}$ norm on $K_i^{+5T+1}$ only in formula (\ref{form184}).
Note we use coordinate $(\tau',t)$ on $K_1^{+5T+1} \setminus K_1$,
and $(\tau'',t)$ on $K_2^{+5T+1} \setminus K_2$.
We remark also that $\Sigma_T = K_1^{+5T+1} \cup K_2^{+5T+1}$.
\end{rem}
\begin{rem}\label{rem136}
Note that $(V^{\rho}_{T,i,(\kappa)},\Delta p^{\rho}_{T,(\kappa)})$ appearing in (\ref{form182}) is an element of the weighted Sobolev space
$L^2_{m+1,\delta}((\Sigma_i,\partial \Sigma_i);(\hat u^{\rho}_{i,T,(\kappa-1)})^*TX,(\hat u^{\rho}_{i,T,(\kappa-1)})^*TL)$
that depends on $T$ and $\rho$. To make sense of $T$ and $\rho$ derivatives
we identify
$$
\aligned
&L^2_{m+1,\delta}((\Sigma_i,\partial \Sigma_i);(\hat u^{\rho}_{i,T,(\kappa-1)})^*TX,(\hat u^{\rho}_{i,T,(\kappa-1)})^*TL)\\
&\cong
L^2_{m+1,\delta}((\Sigma_i,\partial \Sigma_i);u_{i}^*TX,u_{i}^*TL)
\endaligned
$$
as follows. We find $V$ such that $\hat u^{\rho}_{i,T,(\kappa-1)} = {\rm E}(u_{i},V)$.
We use the parallel transport with respect to the path $r \mapsto {\rm E}(u_{i},rV)$ and
its complex linear part to define this isomorphism.
The same remark applies to (\ref{form185}) and (\ref{form186}).
\end{rem}
\begin{rem}
The square of the left hand side of (\ref{form182}), in case $i=1$,  is :
$$
\aligned
&\left\Vert \nabla_{\rho}^n \frac{\partial^{\ell}}{\partial T^{\ell}} V^{\rho}_{T,1,(\kappa)}\right\Vert^2_{L^2_{m+1-\ell}(K_1)}
\\
&+
\sum_{k=0}^{m+1-\ell} \int_{[0,\infty)\times [0,1]}  e_{1,T}(\tau,t)
\left\Vert\nabla_{\tau',t}^k\nabla_{\rho}^n \frac{\partial^{\ell}}{\partial T^{\ell}}
\big(
V^{\rho}_{T,i,(\kappa)} - \text{\rm Pal}(\Delta p^{\rho}_{T,(\kappa)})
\big)\right\Vert^2 d\tau'dt.
\endaligned
$$
Note that we apply Remark \ref{rem136} to define $T$ and $\rho$ derivatives in the above formula.
\par
The case $i=2$ is similar using $\tau''$ coordinate.
\end{rem}
\par\medskip
\noindent{\bf (Part A)} (See \cite[page 776 paragraph (A) and (B)]{fooo:book1}.)
\par
We assume (\ref{form182}) and (\ref{form183}).
\par
We find that
\begin{enumerate}
\item
\begin{equation}\label{1113formu}
{\rm Err}^{\rho}_{1,T,(\kappa)}(z)
=
\Pi_{E_1(\hat u^{\rho}_{1,T,(\kappa-1)})}^{\perp}
\overline \partial {\rm E} (\hat u^{\rho}_{1,T,(\kappa-1)}(z),V^{\rho}_{T,1,(\kappa)}(z))
\end{equation}
for $z \in K_1$.
\item
\begin{equation}\label{form1106}
\aligned
&{\rm Err}^{\rho}_{1,T,(\kappa)}(\tau') \\
= &(1-\chi(\tau'-5T))\overline\partial\big(\chi(\tau'-4T)(V^{\rho}_{T,2,(\kappa)}(\tau'-10T,t) - \Delta p^{\rho}_{T,(\kappa)})
\\
&\qquad\qquad\qquad\qquad+V^{\rho}_{T,1,(\kappa)}(\tau',t)
+u^{\rho}_{T,(\kappa-1)}(\tau',t)\big),
\endaligned
\end{equation}
for $(\tau',t) \in [0,\infty)\times [0,1]$. (Note $\tau' = \tau''+10T$
and the variable of $V^{\rho}_{T,2,(\kappa)}$ is $(\tau'',t)$.)
\end{enumerate}
Here $\chi : \R \to [0,1]$ is a smooth function such that
\begin{equation}\label{chichi}
\chi(\tau)
\begin{cases}
=0  & \tau < -1 \\
=1  & \tau > 1  \\
\in [0,1]  &\tau \in [-1,1].
\end{cases}
\end{equation}
\par
Note that in Formulas (\ref{form182})-(\ref{form186}) the Sobolev norms
in the left hand side are $L^2_{m+1-\ell,\delta}(\Sigma_i)$
etc. and are not $L^2_{m+1,\delta}(\Sigma_i)$ etc.
The origin of this loss of differentiability (in the sense of Sobolev space)
comes from the term  $V^{\rho}_{T,2,(\kappa)}(\tau'-10T)$.
In fact, we have
$$
\frac{\partial }{\partial T}V^{\rho}_{T_1,2,(\kappa)}(\tau'-10T)
= -10 \frac{\partial }{\partial \tau''}V^{\rho}_{T_1,2,(\kappa)}(\tau'-10T)
$$
for a fixed $T_1$.
Hence $\partial/\partial T$ is continuous as $L^2_{m+1} \to L^2_m$.
We remark in (\ref{form182}) for $i=2$ we use the coordinate $(\tau'',t)$
on $(-\infty,0]\times [0,1]$ to define $T$ derivative of $V_{T,2,(\kappa)}^{\rho}$.
\par
Taking this fact  into account the proof goes as follows.
\par
We can estimate $T$ and $\rho$ derivative of ${\rm Err}^{\rho}_{1,T,(\kappa)}$
on $K_1$ in the same way as the proof of
Lemma \ref{mainestimatestep13kappa}.
\begin{rem}\label{remark127}
The fact we use here is that the maps such as $(u,v) \mapsto {\rm E}(u,v)$,
$(u,v) \to \Pi^{\perp}_{E_i(u)}(v)$ are smooth maps from $L^2_{m+1,loc} \times L^2_{m+1,\delta}
\to L^2_{m+1,\delta}$  or
$L^2_{m+1,loc} \times L^2_{m,\delta}  \to L^2_{m,\delta}$ and $u \to \overline\partial u$ is a smooth map
$L^2_{m+1,\delta}  \to L^2_{m,\delta}$. (Since we assume $m$ sufficiently large this is a well-known fact.)
Moreover the map $T \mapsto u^{\rho}_{T,(\kappa-1)}$
and $T \mapsto V^{\rho}_{T,1,(\kappa)}$ are $C^{\ell}$ maps as
a map $[T(m),\infty) \to L^2_{m+1-\ell,\delta}$ with its differential estimated by
induction hypothesis (\ref{form184}) and (\ref{form182}).
\par
We note that $\rho \mapsto u^{\rho}_{T,(\kappa-1)}$ is smooth as a map
$V_1 \times_L V_2 \to L^2_{m+1,\delta}$.
\end{rem}
The estimates of $T$ and $\rho$ derivatives of (\ref{form1106}) are as follows.
\par
We first consider the domain $\tau' \in [4T+1,\infty)$.
There we have
\begin{equation}\label{form11010}
\aligned
{\rm Err}^{\rho}_{1,T,(\kappa)}(\tau',t)
= &(1-\chi(\tau'-5T))\overline\partial(V^{\rho}_{T,2,(\kappa)}(\tau'-10T,t)\\
&+V^{\rho}_{T,1,(\kappa)}(\tau',t)
+u^{\rho}_{T,(\kappa-1)}(\tau',t)- \Delta p^{\rho}_{T,(\kappa)}).
\endaligned
\end{equation}
By the same calculation as in the proof of Lemma \ref{mainestimatestep13kappa},
(\ref{form11010}) is equal to
$$
\aligned
(1-\chi(\tau'-5T))\int_0^1ds \int_0^s
\frac{\partial^2}{\partial r^2}& \overline\partial \big(r(V^{\rho}_{T,2,(\kappa)}(\tau'-10T) - \Delta p^{\rho}_{T,(\kappa)})\\
&+r(V^{\rho}_{T,1,(\kappa)}(\tau',t) - \Delta p^{\rho}_{T,(\kappa)})\\
&
+u^{\rho}_{T,(\kappa-1)}(\tau',t)
+ r\Delta p^{\rho}_{T,(\kappa)}\big) dr.
\endaligned
$$
(Note that we are away from the support of $E_i$.)\footnote{Note $\overline\partial$ is non-constant. So
$\overline\partial (r(V^{\rho}_{T,2,(\kappa)}(\tau'-10T) - \Delta p^{\rho}_{T,(\kappa)})
+r(V^{\rho}_{T,1,(\kappa)}(\tau',t)  - \Delta p^{\rho}_{T,(\kappa)})
+u^{\rho}_{T,(\kappa-1)}(\tau',t) + r\Delta p^{\rho}_{T,(\kappa)})$ is nonlinear on $r$.}
Using the fact that
$T \mapsto (V^{\rho}_{T,1,(\kappa)}(\tau',t)-\Delta p^{\rho}_{T,(\kappa)}) + (V^{\rho}_{T,2,(\kappa)}(\tau'-10T)-\Delta p^{\rho}_{T,(\kappa)})$
and $T \mapsto u^{\rho}_{T,(\kappa-1)}(\tau',t)$ are of $C^{\ell}$ class
as a map to $L^2_{m+1-\ell,\delta}$,
we can estimate it to obtain the required estimate (\ref{form185}) on
this part. We remark
$T \mapsto (V^{\rho}_{T,2,(\kappa-1)},\Delta p^{\rho}_{T,(\kappa-1)})$ is $C^{\ell}$ with exponential decay
estimate on $T$ derivatives as a map $[T(m),\infty) \to L^{2}_{m-\ell+1,\delta}$.
This follows from the induction hypothesis as follows.
\begin{equation}\label{1106}
\aligned
&\frac{\partial^{\ell}}{\partial T^{\ell}}
\left(\left.
V^{\rho}_{T,2,(\kappa)}(\tau' - 10 T)
\right)\right\vert_{T=T_1} \\
&
=
\sum_{\ell_1+\ell_2 = \ell} (-10)^{\ell_2}
\frac{\partial^{\ell_1}}{\partial T^{\ell_1}}\frac{\partial^{\ell_2}}{(\partial \tau'')^{\ell_2}}
V^{\rho}_{T,2,(\kappa)}(\tau' - 10 T_1).
\endaligned
\end{equation}
The $L^2_{m+1-\ell,\delta}$-norm of the right hand side can be estimated by
(\ref{form182}).
\par
We next consider $\tau' \in [0,4T+1]$.
There we have
\begin{equation}\label{form110622}
\aligned
{\rm Err}^{\rho}_{1,T,(\kappa)}(\tau',t)
= &\overline\partial(\chi(\tau'-4T)(V^{\rho}_{T,2,(\kappa)}(\tau'-10T)- \Delta p^{\rho}_{T,(\kappa)})
\\
&+V^{\rho}_{T,1,(\kappa)}(\tau',t)
+u^{\rho}_{T,(\kappa-1)}(\tau',t)).
\endaligned
\end{equation}
Note
$$
\overline\partial u^{\rho}_{T,(\kappa-1)}(\tau',t))
= {\rm Err}^{\rho}_{1,T,(\kappa-1)}(\tau',t),
$$
there.
Therefore we can calculate in the same way as the proof of Lemma \ref{mainestimatestep13kappa} to find
$$
\aligned
&\overline\partial(V^{\rho}_{T,1,(\kappa)}(\tau',t)
+u^{\rho}_{T,(\kappa-1)}(\tau',t))\\
&=\int_0^1ds \int_0^s
\frac{\partial^2}{\partial r^2} \overline\partial (
r(V^{\rho}_{T,1,(\kappa)}(\tau',t) - \Delta p^{\rho}_{T,(\kappa)})
+u^{\rho}_{T,(\kappa-1)}(\tau',t) + r\Delta p^{\rho}_{T,(\kappa)})dr.
\endaligned$$
We can again estimate the right hand side by using
 the fact that the maps
$T \mapsto (V^{\rho}_{T,1,(\kappa)}(\tau',t),\Delta p^{\rho}_{T,(\kappa)})$
and $T \mapsto u^{\rho}_{T,(\kappa-1)}(\tau',t)$ are of $C^{\ell}$ class
as a map to $L^2_{m+1-\ell,\delta}$ with estimate (\ref{form184}).
\par
Finally we observe that the ratio between weight function of $L_{m+1,\delta}^2(\Sigma_2)$
and of $L_{m+1,\delta}^2(\Sigma_T)$  is $e^{2T\delta}$ on $\tau = -T$ (that is $\tau' = 4T$).
We use this fact to estimate
$\overline\partial(\chi(\tau'-4T)(V^{\rho}_{T,2,(\kappa)}(\tau'-10T)- \Delta p^{\rho}_{T,(\kappa)}))$.
We thus obtain the required estimate
(\ref{form185}) for ${\rm Err}^{\rho}_{1,T,(\kappa)}$ on $\tau' \in [0,4T+1]$.
\par
We thus obtain an estimate for ${\rm Err}^{\rho}_{1,T,(\kappa)}(\tau',t)$.
\par
The estimate of derivatives of ${\rm Err}^{\rho}_{2,T,(\kappa)}(\tau',t)$ is similar.
Thus we have (\ref{form185}).
\par
We note that $\frak e^{\rho} _{i,T,(0)}$ is independent of $T$ as an element of $E_i$.
Among $ \frak e^{\rho} _{i,T,(\kappa)}$'s, the term $\frak e^{\rho} _{i,T,(0)}$ is the only
one that is not of exponential decay with respect to $T$. Once we note  this point
the rest of the proof of (\ref{form186}) is the same as the proof of Lemma \ref{mainestimatestep13kappa}.
\par
We finally prove (\ref{form184}). On $K_1$ we have
$$
u^{\rho}_{T,(\kappa)} = {\rm E} (u^{\rho}_{T,(\kappa-1)},V^{\rho}_{1,T,(\kappa)}).
$$
So using $\mu<1$, (\ref{form184}) follows from (\ref{form182}) on $K_1$.
\par
On $(\tau',t) \in [0,5T+1) \times [0,1]$ we have:
$$
\aligned
&u^{\rho}_{T,(\kappa)}(\tau',t) \\
&=
 V^{\rho}_{T,1,(\kappa)}(\tau',t)
+(1-\chi(\tau'-4T))(V^{\rho}_{T,2,(\kappa)}(\tau'-10T,t)- \Delta p^{\rho}_{T,(\kappa)})\\
&\qquad\qquad\qquad\qquad\qquad\qquad\qquad\qquad\qquad\qquad
+u^{\rho}_{T,(\kappa-1)}(\tau',t) \\
&=
\sum_{a=1}^{\kappa} V^{\rho}_{T,1,(a)}(\tau',t)
+ (1-\chi(\tau'-4T))\sum_{a=1}^{\kappa}(V^{\rho}_{T,2,(a)}(\tau'-10T,t)- \Delta p^{\rho}_{T,(a)})\\
&\qquad\qquad\qquad\qquad\qquad\qquad\qquad\qquad\qquad\qquad
+u^{\rho}_{T,(0)}(\tau',t).
\endaligned
$$
Then using a calculation similar to (\ref{1106}) we have (\ref{form182}) on $(\tau',t) \in [0,5T+1) \times [0,1]$.
\begin{rem}\label{Abremark}
In \cite{Abexotic} Abouzaid used $L^p_1$ norm for the maps $u$. He then proved that
the gluing map is continuous with respect to $T$ (that is $S$ in the notation of \cite{Abexotic})
but does not prove its differentiability with respect to $T$.
(Instead he used the technique to remove the part of the
moduli space with $T>T_0$.  See Subsection \ref{subsec342}.
This technique certainly works for the purpose of \cite{Abexotic}.)
In fact if we use $L^p_1$ norm instead of $L^2_m$ norm
then the left hand side of (\ref{form184}) becomes $L^p_{-1}$ norm
which is hard to use.
\par
Abouzaid mentioned in  \cite[Remark 5.1]{Abexotic} that this point is related
to the fact that quotients of Sobolev spaces by the diffeomorphisms in the
source are not naturally equipped with the structure of smooth Banach manifold.
Indeed in the situation when there is an automorphism on $\Sigma_2$,
for example $\Sigma_2$ is disk with one boundary marked point $(-\infty,t)$,
then the $T$ parameter is killed by a part of the automorphism.
So the shift of $V^{\rho}_{T,2,(\kappa)}$ by $T$ that appears in the
second term of (\ref{form1106}) will be equivalent to the action of
the automorphism group of $\Sigma_2$ in such a situation.
The shift of $T$ causes the loss of differentiability in the sense of Sobolev space
in the formulas (\ref{form182}) -(\ref{form186}).
However at the end of the day we can still get the differentiability of $C^{\infty}$ order and its
exponential decay by using various {\it weighted} Sobolev spaces with various $m$ simultaneously.
(See Remark \ref{differentm} also.)
\end{rem}
\par\medskip
\noindent{\bf (Part B)} (See \cite[page 776 the paragraph next to (B)]{fooo:book1}.)
\par
We assume  (\ref{form182})-(\ref{form186}) for $\kappa$ and will prove
(\ref{form182}) and (\ref{form183}) for $\kappa+1$.
This part is nontrivial only because  the construction here is global.
(Solving linear equation.) So we first review the set up of the function space that is independent of $T$.
\par
In Definition \ref{defnfraH}
we defined a function space
$\frak H(E_1,E_2)$,
that is a subspace of (\ref{cocsissobolev}).
Since  (\ref{cocsissobolev}) is still $T$ dependent we rewrite it a bit.
We consider $u_{i}^{\rho} : (\Sigma_i,\partial\Sigma_i) \to (X,L)$ that is $T$-{\it independent}.
\par
The maps
$\hat u^{\rho}_{i,T,(\kappa)}$ are
close to $u_{i}^{\rho}$.
(Namely the $C^0$ distance between them is smaller than
injectivity radius of $X$.)
We take a connection of $TX$ so that $L$ is
totally geodesic.
We use the complex linear part of the parallel
transport with respect to this connection,
to send
$$
\bigoplus_{i=1}^2 L^2_{m,\delta}((\Sigma_i,\partial \Sigma_i);(u_{i}^{\rho})^{*}TX,(u_{i}^{\rho})^{*}TL)
$$
to
$$
\bigoplus_{i=1}^2 L^2_{m,\delta}((\Sigma_i,\partial \Sigma_i);(\hat u^{\rho}_{i,T,(\kappa)})^{*}TX,(\hat u^{\rho}_{i,T,(\kappa)})^{*}TL).
$$
Note that
$\text{\rm Ker}(D{\rm ev}_{1,\infty} - D{\rm ev}_{2,\infty})$
is sent to
$\text{\rm Ker}(D{\rm ev}_{1,\infty} - D{\rm ev}_{2,\infty})$ by this map.
Therefore we obtain an isomorphism between
\begin{equation}\label{Tindependentfcs}
\text{\rm Ker}(D{\rm ev}_{1,\infty} - D{\rm ev}_{2,\infty}) \cap
\bigoplus_{i=1}^2 L^2_{m,\delta}((\Sigma_i,\partial \Sigma_i);(u_{i}^{\rho})^{*}TX,(u_{i}^{\rho})^{*}TL)
\end{equation}
and
\begin{equation}\label{cocsissobolevkappa}
\text{\rm Ker}(D{\rm ev}_{1,\infty} - D{\rm ev}_{2,\infty}) \cap
\bigoplus_{i=1}^2 L^2_{m,\delta}((\Sigma_i,\partial \Sigma_i);(\hat u^{\rho}_{i,T,(\kappa)})^{*}TX,(\hat u^{\rho}_{i,T,(\kappa)})^{*}TL).
\end{equation}
In case $\kappa =0$ we send
$\frak H(E_1,E_2)$ by this isomorphism  to obtain a subspace of (\ref{Tindependentfcs})
which we denote by $\frak H(E_1,E_2)$ by an abuse of notation.
We send it to the subspace of (\ref{cocsissobolevkappa})
and denote it by $\frak H(E_1,E_2;\kappa,T)$.
We thus have an isomorphism
$$
I_{1,\kappa,T} : \frak H(E_1,E_2) \to \frak H(E_1,E_2;\kappa,T).
$$
\par
We next use the parallel transport in the same way to find an isomorphism
$$
I_{2,\kappa,T} : L^2_{m,\delta}(\Sigma_i;(u_{i}^{\rho})^*TX \otimes \Lambda^{01}) \to
L^2_{m,\delta}(\Sigma_i;(\hat u^{\rho}_{i,T,(\kappa)})^{*}TX \otimes \Lambda^{01}).
$$
Thus the composition
$$
I_{2,\kappa,T} ^{-1}\circ \left(D_{\hat u^{\rho}_{i,T,(\kappa-1)}}\overline{\partial} -
(D_{\hat u^{\rho}_{i,T,(\kappa-1)}}E_i)((\frak {se})^{\rho} _{i,T,(\kappa-1)}, \cdot))\right)
\circ I_{1,\kappa,T}
$$
defines an operator
$$
D_{\kappa,T} :\frak H(E_1,E_2)
\to L^2_{m,\delta}(\Sigma_i;(u_{i}^{\rho})^*TX \otimes \Lambda^{01}).
$$
Here the domain and the target is independent of $T,\kappa$.
\begin{rem}
Note $D_{\hat u^{\rho}_{i,T,(\kappa-1)}}\overline{\partial} -
(D_{\hat u^{\rho}_{i,T,(\kappa-1)}}E_i)((\frak {se})^{\rho} _{i,T,(\kappa-1)},\cdot)$ is the differential operator
in (\ref{linearized2221}) and (\ref{144op}).
This differential operator gives the linearization of the right hand side
of (\ref{1113formu}).
\end{rem}
\par
We next eliminate $T,\kappa$ dependence of $E_i$.
We consider the finite dimensional subspace:
$$
E_i(\hat u^{\rho}_{i,T,(\kappa)})
\subset
L^2_{m,\delta}(\Sigma_i;(\hat u^{\rho}_{i,T,(\kappa)})^{*}TX \otimes \Lambda^{01}).
$$
Let us consider
$$
E_{i,(\kappa),T} = I_{2,\kappa,T} ^{-1}(E_i(\hat u^{\rho}_{i,T,(\kappa)}))
$$
that may depend on $T$.
However
$$
E_{i,(0)}  = I_{2,\kappa,T} ^{-1}(E_i(\hat u^{\rho}_{i,T,(0)}))
$$
is independent of $T$  since $\hat u^{\rho}_{i,T,(0)}= u_{i}^{\rho}$ on $K_i$.
Let $E_{i,(0)}^{\perp}$ be the $L^2$ orthogonal complement of $E_{i,(0)} $ in
$L^2_{m,\delta}(\Sigma_i;(\hat u^{\rho}_{i,T,(\kappa)})^{*}TX \otimes \Lambda^{01})$.
\par
We have
\begin{equation}\label{artificialsplitting}
E_{i,(\kappa),T} \oplus E_{i,(0)}^{\perp} = L^2_{m,\delta}(\Sigma_i;(u_{i}^{\rho})^*TX \otimes \Lambda^{01}).
\end{equation}
Therefore the inclusion induces an isomorphism
$$
E_{i,(0)}^{\perp} \cong L^2_{m,\delta}(\Sigma_i;(u_{i}^{\rho})^*TX \otimes \Lambda^{01})/E_{i,(\kappa),T}.
$$
We thus obtain
\begin{equation}
\overline D_{\kappa,T} :\frak H(E_1,E_2)
\to E_{i,(0)}^{\perp}.
\end{equation}
The induction hypothesis implies the following:
\begin{enumerate}
\item
There exist $C_{20,m}, C_{21,m} >0$ such that
\begin{equation}
C_{20,m} \Vert V\Vert_{L^2_{m+1,\delta}} \le \Vert\overline D_{0,T}(V)\Vert_{L^2_{m,\delta}} \le C_{21,m}  \Vert V\Vert_{L^2_{m+1,\delta}}.
\end{equation}
\item
\begin{equation}\label{estimateDover}
\Vert\overline D_{\kappa,T}(V) -  \overline D_{0,T}(V) \Vert_{L^2_{m,\delta}} \le C_{21,m}  e^{-\delta T}\Vert V\Vert_{L^2_{m+1,\delta}}.
\end{equation}
Moreover
\begin{equation}\label{11230}
\left\Vert\nabla_{\rho}^n \frac{\partial^{\ell}}{\partial T^{\ell}}  \overline D_{\kappa,T}(V)  \right\Vert_{L^2_{m-\ell,\delta}} \le C_{22,m}  e^{-\delta T}
\Vert V\Vert_{L^2_{m+1,\delta}}.
\end{equation}
\end{enumerate}
In fact, (\ref{11230}) follows from
\begin{eqnarray}
&\displaystyle\left\Vert \nabla_{\rho}^n \frac{\partial^{\ell}}{\partial T^{\ell}}\hat u^{\rho}_{i,T,(\kappa)}\right\Vert_{L^2_{m-\ell}(K_i)}
\le C_{23,m}e^{-\delta T},\label{1124}\\
&\displaystyle\left\Vert \nabla_{\rho}^n \frac{\partial^{\ell}}{\partial T^{\ell}}\hat u^{\rho}_{i,T,(\kappa)}\right\Vert_{L^2_{m-\ell}([S,S+1]\times [0,1])}
\le C_{23,m}e^{-\delta T}\label{1125}
\end{eqnarray}
for any $S\in [0,\infty)$.
(See also the Remark \ref{newremark11}.) Note that the weighted Sobolev norm
$
\Vert \nabla_{\rho}^n \frac{\partial^{\ell}}{\partial T^{\ell}}\hat u^{\rho}_{i,T,(\kappa)}\Vert_{L^2_{m-\ell,\delta}(\Sigma_i)}
$
can be large because
$$
\frac{\partial}{\partial T} \chi_{\mathcal B}^{\leftarrow}(\tau-T,t) u^{\rho}_{T,(\kappa-1)}
$$
is only estimated by $e^{-3\delta T}$ on the support of $\chi_{\mathcal B}^{\leftarrow}(\tau-T,t)$
but the weight $e_{1,\delta}$ is roughly $e^{7T\delta}$ on the support of $\chi_{\mathcal B}^{\leftarrow}(\tau-T,t)$.
However this does not cause any problem to prove (\ref{11230}). In fact the operator $\overline D_{\kappa,T}$
is a differential operator whose coefficient depends on $\hat u^{\rho}_{i,T,(\kappa)}$.
So to estimate the operator norm of its derivatives with respect to the {\it weighted} Sobolev norm,
we only need to estimate the local Sobolev norm without weight of  $\hat u^{\rho}_{i,T,(\kappa)}$,
that is provided by (\ref{1124}) and (\ref{1125}).
\par
We note that $\overline D_{0,T}$ is independent of $T$. So we write $\overline D_0$.
Now we have:
\begin{equation}
\aligned
\overline D_{\kappa,T}^{-1}
&= \left( ( 1 + (\overline D_{\kappa,T} -  \overline D_0)\overline D_0^{-1})\overline D_0\right)^{-1} \\
& = \overline D_0^{-1}\sum_{k=0}^{\infty} (-1)^k ((\overline D_{\kappa,T} -  \overline D_0))\overline D_0^{-1})^k.
\endaligned
\end{equation}
Therefore
\begin{equation}\label{derivativeDest}
\left\Vert \nabla_{\rho}^n \frac{\partial^{\ell}}{\partial T^{\ell}} \overline D_{\kappa,T}^{-1}(W)\right\Vert_{L^2_{m+1-\ell,\delta}}
\le C_{24,m} e^{-\delta}\Vert W\Vert_{L^2_{m,\delta}}
\end{equation}
for $\ell > 0$ and $\ell + n \le m$. (Here we assume $W$ is $T$ independent.)
Since
$$
(V^{\rho}_{T,1,(\kappa+1)},V^{\rho}_{T,2,(\kappa+1)},\Delta p^{\rho}_{T,(\kappa+1)})
=
(I_{1,\kappa,T}\circ \overline D_{\kappa,T}^{-1} \circ I_{2,\kappa,T}^{-1})
({\rm Err}^{\rho}_{1,T,(\kappa)},{\rm Err}^{\rho}_{2,T,(\kappa)}),
$$
(\ref{form185}) and (\ref{derivativeDest}) imply (\ref{form182}) and (\ref{form183}) for $\kappa +1$.
\par
The proof of Theorem \ref{exdecayT} is now complete.
\end{proof}
\begin{rem}\label{newremark11}
Let us add a few more explanation about the proof of
(\ref{estimateDover}) and
(\ref{11230}). Especially the relation between two operators $D_{\kappa,T}$ and $\overline D_{\kappa,T}$.
We consider the direct sum decomposition
\begin{equation}\label{artificialsplitting2}
E_{i,(\kappa),T} \oplus E_{i,(0)}^{\perp} = L^2_{m,\delta}(\Sigma_i;(u_{i}^{\rho})^*TX \otimes \Lambda^{01}).
\end{equation}
Note that this is not an orthogonal decomposition.
We take an isomorphism
$$
B_{i,(\kappa),T} : L^2_{m,\delta}(\Sigma_i;(u_{i}^{\rho})^*TX \otimes \Lambda^{01})
\to L^2_{m,\delta}(\Sigma_i;(u_{i}^{\rho})^*TX \otimes \Lambda^{01})
$$
such that according to the orthogonal decomposition
\begin{equation}\label{bettersplitting2}
E_{i,(0)} \oplus E_{i,(0)}^{\perp} = L^2_{m,\delta}(\Sigma_i;(u_{i}^{\rho})^*TX \otimes \Lambda^{01}).
\end{equation}
The restriction
$
B_{i,(\kappa),T}\vert_{E_{i,(0)}^{\perp}}
$
is the identity map and the restriction $B_{i,(\kappa),T}\vert_{E_{i,(0)}}$ is the canonical isomorphism
$$
A_{i,(\kappa),T}  : E_{i,(0)} \to E_{i,(0)}
$$
given by the parallel transportation. Namely we put
$$
B_{i,(\kappa),T} = A_{i,(\kappa),T} \circ \Pi_{E_{i,(0)}} + \Pi_{E_{i,(0)}^{\perp}}.
$$
It is easy to prove
\begin{equation}\label{estimateDoverB}
\Vert\overline B_{i,(\kappa),T}(V)  -  V \Vert_{L^2_{m,\delta}} \le C_{25,m}  e^{-\delta T}\Vert V\Vert_{L^2_{m,\delta}}.
\end{equation}
Moreover
\begin{equation}\label{11230B}
\left\Vert\nabla_{\rho}^n \frac{\partial^{\ell}}{\partial T^{\ell}}  B_{i,(\kappa),T}(V)   \right\Vert_{L^2_{m-\ell,\delta}} \le C_{26,m}  e^{-\delta T}
\Vert V\Vert_{L^2_{m,\delta}}.
\end{equation}
Note that
$$
C_{i,(\kappa),T} = \Pi_{E_{i,(0)}^{\perp}}\circ B_{i,(\kappa),T}^{-1}
$$
is the projection to the second factor in (\ref{artificialsplitting2}) and hence
\begin{equation}\label{relDandDover}
\overline D_{\kappa,T} = \Pi_{E_{i,(0)}^{\perp}}\circ B_{i,(\kappa),T}^{-1}\circ D_{\kappa,T}.
\end{equation}
We can use (\ref{estimateDoverB}), (\ref{11230B}) and (\ref{relDandDover}) to prove
(\ref{estimateDover}),
(\ref{11230}).
\end{rem}

\section{Surjectivity and injectivity of the gluing map}
\label{surjinj}

In this subsection we prove surjectivity and injectivity of the map $\text{\rm Glu}_T$ in Theorem \ref{gluethm1} and
complete the proof of Theorem \ref{gluethm1}.\footnote{
Here surjectivity means the second half of the statement of Theorem \ref{gluethm1}, that is
`The image contains $\mathcal M^{E_1+E_2}((\Sigma_T,\vec z);\beta)_{\epsilon_3}$.'}
The proof goes along the line of \cite{Don83}. (See also \cite{freedUhlen}.)
The surjectivity proof is written in  \cite[Section 14]{FOn} and
injectivity is proved in the same way.
(\cite[Section 14]{FOn} studies the case of pseudo-holomorphic curve without boundary.
It however can be adapted easily to the bordered case as we mentioned in
\cite[page 417 lines 21-26]{fooo:book1}.)
Here we explain the argument in our situation in more detail.
\par
We begin with the following a priori estimate.
\begin{prop}$($\cite[Lemma 11.2]{FOn}$)$\label{neckaprioridecay}
There exist $\epsilon_3, C_{25,m}, \delta_2 > 0$ such that
if $u : (\Sigma_T,\partial\Sigma_T) \to (X,L)$ is an element of
$\mathcal M^{E_1+E_2}((\Sigma_T,\vec z);\beta)_{\epsilon}$
for $0 <\epsilon<\epsilon_3$ then we have
\begin{equation}
\left\Vert \frac{\partial u}{\partial \tau}\right\Vert_{C^m([\tau-1,\tau+1] \times [0,1])} \le C_{27,m} e^{-\delta_2 (5T - \vert\tau\vert)}.
\end{equation}
\end{prop}
The proof is the same as \cite[Lemma 11.2]{FOn} that is proved  in \cite[Section 14]{FOn} and so is omitted.
\par
We also have the following:
\begin{lem}\label{gluedissmoothindex}
$\mathcal M^{E_1+E_2}((\Sigma_T,\vec z);\beta)_{\epsilon}$
is a smooth manifold of dimension $\dim V_1 + \dim V_2 - \dim L$.
\end{lem}
This is a consequence of the implicit function theorem and the index sum formula.
\par
\medskip
\begin{proof}[Proof of surjectivity]
During this proof we take $m$ sufficiently large and fix it.
We will fix $\epsilon$ and $T_0$ during the proof and assume $T>T_0$.
(They are chosen so that the discussion below works.)
Let $u :  (\Sigma_T,\partial\Sigma_T) \to (X,L)$ be an element of $\mathcal M^{E_1+E_2}((\Sigma_T,\vec z);\beta)_{\epsilon}$.
The purpose here is to show that $u$ is in the image of $\text{\rm Glu}_T$.
We define $u_i' : (\Sigma_i,\partial\Sigma_i) \to (X,L)$ as follows.
We put $p^u_0 = u(0,0) \in L$.
\begin{equation}
\aligned
&u'_{1}(z) \\
&=
\begin{cases} \chi_{\mathcal B}^{\leftarrow}(\tau-T,t)   u(\tau,t) + \chi_{\mathcal B}^{\rightarrow}(\tau-T,t)p_0^u
&\text{if $z = (\tau,t) \in [-5T,5T] \times [0,1]$} \\
 u(z)
&\text{if $z \in K_1$} \\
p_0^u
&\text{if $z \in [5T,\infty)\times [0,1]$}.
\end{cases} \\
&u'_2(z) \\
&=
\begin{cases} \chi_{\mathcal A}^{\rightarrow}(\tau+T,t)   u(\tau,t) + \chi_{\mathcal A}^{\leftarrow}(\tau+T,t)p_0^u
&\text{if $z = (\tau,t) \in [-5T,5T] \times [0,1]$} \\
 u(z)
&\text{if $z \in K_2$} \\
p_0^u
&\text{if $z \in (-\infty,-5T]\times [0,1]$}.
\end{cases}
\endaligned
\end{equation}
Proposition \ref{neckaprioridecay} implies
\begin{equation}
\Vert
\Pi_{E_i(u'_i)}\overline\partial u'_{i}
\Vert_{L^2_{m,\delta}(\Sigma_i)}
\le C_{28,m} e^{-\delta T}.
\end{equation}
Here we take $\delta < \delta_2/10$.
On the other hand, by assumption and elliptic regularity we have
\begin{equation}
\Vert
u'_i - u_i
\Vert_{L^2_{m+1,\delta}(\Sigma_i)}
\le C_{29,m} \epsilon.
\end{equation}
Therefore by an implicit function theorem we have the following:
\begin{lem}
There exists $\rho_i \in V_i$ such that
\begin{equation}\label{1121}
\Vert
u'_i - u^{\rho_i}_i
\Vert_{L^2_{m+1,\delta}(\Sigma_i)}
\le C_{30,m} e^{-\delta T},
\end{equation}
$\rho = (\rho_1,\rho_2) \in V_1 \times_L V_2$,
and
\begin{equation}
\vert\rho_i\vert \le C_{31,m}\epsilon.
\end{equation}
\end{lem}
(Note when $\rho_i = 0$, $u_i^{\rho_i} = u_i$.)
\par
By (\ref{1121}) we have
\begin{equation}\label{1123}
\Vert
u - u^{\rho}_T
\Vert_{L^2_{m+1,\delta}(\Sigma_T)}
\le C_{32,m} e^{-\delta T}.
\end{equation}
Here $u^{\rho}_T =\text{\rm Glu}_T(\rho)$.
\par
We take $V \in \Gamma((\Sigma_T,\partial \Sigma_T);(u^{\rho}_T)^*TX;(u^{\rho}_T)^*TL)$ so that
$$
u(z) = {\rm E} (u^{\rho}_T(z),V(z)).
$$
We define $u^s : (\Sigma_T,\partial \Sigma_T) \to (X,L)$ by
\begin{equation}
u^s(z) ={\rm E} (u^{\rho}_T(z),sV(z)).
\end{equation}
(\ref{1123}) implies
\begin{equation}
\Vert
\Pi^{\perp}_{(E_1+E_2)(u^s)}\overline\partial u^s
\Vert_{L^2_{m,\delta}(\Sigma_T)}
\le C_{33,m} e^{-\delta T}
\end{equation}
and
\begin{equation}
\left\Vert
\frac{\partial}{\partial s} u^s
\right\Vert_{L^2_{m+1,\delta}(K_i^{+S})}
\le C_{34,m} e^{-\delta T}
\end{equation}
for each $s \in [0,1]$.
\begin{lem}
If $T$ is sufficiently large, then there exists $\hat u^s :  (\Sigma_T,\partial \Sigma_T) \to (X,L)$ $(s \in [0,1])$
with the following properties.
\begin{enumerate}
\item
$$
\overline\partial \hat u^s \equiv 0 \mod (E_1+E_2)(\hat u^s).
$$
\item
\begin{equation}\label{sderovatove}
\left\Vert
\frac{\partial}{\partial s} \hat u^s
\right\Vert_{L^2_{m+1,\delta}(K_i^{+S})}
\le 2C_{35,m} e^{-\delta T}.
\end{equation}
\item $\hat u^s = u^s$ for $s=0,1$.
\end{enumerate}
\end{lem}
\begin{proof}
Run the alternating method described in Subsection \ref{alternatingmethod}
in one parameter family version.
Since $u^s$ is already a solution for $s=0,1$, it does not change.
\end{proof}
\begin{lem}\label{immersion}
The map
$\text{\rm Glu}_T : V_1 \times_L V_2 \to \mathcal M^{E_1+E_2}((\Sigma_T,\vec z);\beta)_{\epsilon}$
is an immersion if $T$ is sufficiently large.
\end{lem}
\begin{proof}
We consider the composition of $\text{\rm Glu}_T$ with
$$
\mathcal M^{E_1+E_2}((\Sigma_T,\vec z);\beta)_{\epsilon}
\to  {L^2_{m+1}}((K_i^{+S},K_i^{+S}\cap\partial \Sigma_i),(X,L))
$$
defined by restriction.
In the case $T = \infty$ this composition is obtained by restriction of maps.
By unique continuation, this is certainly an immersion for $T=\infty$.
Then Theorem \ref{exdecayT} implies that it is an immersion for sufficiently large $T$.
\end{proof}
Now we will prove that
$$
A = \{s \in [0,1] \mid \hat u^s \in \text{image of $\text{\rm Glu}_T$}\}
$$
is open and closed.
Lemma \ref{gluedissmoothindex} implies that $\mathcal M^{E_1+E_2}((\Sigma_T,\vec z);\beta)_{\epsilon}$
is a smooth manifold and has the same dimension as $V_1 \times_L V_2$.
Therefore Lemma \ref{immersion} implies that $A$ is open.
The closedness of $A$ follows from (\ref{sderovatove}).
\par
Note $0 \in A$. Therefore $1\in A$. Namely $u$ is in the image of $\text{\rm Glu}_T$ as required.
\end{proof}
\begin{proof}[Proof of injectivity]
Let $\rho^j = (\rho_1^j,\rho_2^j) \in V_1\times_L V_2$ for $j=0,1$. We assume
\begin{equation}
\text{\rm Glu}_T(\rho^0) = \text{\rm Glu}_T(\rho^1)
\end{equation}
and
\begin{equation}
\Vert \rho_i^j\Vert < \epsilon.
\end{equation}
We will prove that $\rho^0 = \rho^1$ if $T$ is sufficiently large and $\epsilon$ is
sufficiently small.
We may assume that $V_1\times_L V_2$ is connected and simply connected.
Then, we have a path $s \mapsto \rho^s =(\rho^s_1,\rho^s_2) \in V_1 \times_L V_2$
such that
\begin{enumerate}
\item $\rho^s = \rho^j$ for $s=1$, $j=0,1$.
\item
$$
\left\Vert \frac{\partial}{\partial s} \rho^s \right\Vert \le \Phi_1(\epsilon)
$$
where $\lim_{\epsilon \to 0}\Phi_1(\epsilon) = 0$.
\end{enumerate}
We define $V(s) \in  \Gamma((\Sigma_T,\partial \Sigma_T);(u^{\rho^0}_T)^*TX;(u^{\rho^0}_T)^*TL)$
such that
$$
u^{\rho^s}_T(z) = {\rm E} (u^{\rho^0}_T(z),V(s)(z)).
$$
(By (2) $u^{\rho^s}_T(z)$ is $C^0$-close to $u^{\rho^0}_T(z)$, as $\epsilon \to 0$. Therefore
there exists such a unique $V(s)$ if $\epsilon$ is small.)
Note $V(1) = V(0)$ since $u^{\rho^1} = u^{\rho^0}$.
Therefore for $w \in D^2 = \{w \in \C \mid \vert w\vert \le 1\}$ there exists $V(w)$ such that
\begin{enumerate}
\item
$V(s) = V(w)$ if $w = e^{2\pi\sqrt{-1}s}$.
\item
We put $w = x + \sqrt{-1}y$.
\begin{equation}
\left\Vert
\frac{\partial}{\partial x} V(w)
\right\Vert _{L^2_{m+1,\delta}(\Sigma_T)}
+
\left\Vert
\frac{\partial}{\partial y} V(w)
\right\Vert _{L^2_{m+1,\delta}(\Sigma_T)}
\le \Phi_2(\epsilon)
\end{equation}
where $\lim_{\epsilon \to 0}\Phi_2(\epsilon) = 0$.
\end{enumerate}
We put $u^w(z) = {\rm E}(u^{\rho^0}_T(z),V(w)(z))$.
\begin{lem}
If $T$ is sufficiently large and $\epsilon$ is sufficiently small
then there exists $\hat u^w :  (\Sigma_T,\partial \Sigma_T) \to (X,L)$ $(s \in [0,1])$
with the following properties.
\begin{enumerate}
\item
$$
\overline\partial \hat u^w \equiv 0 \mod (E_1+E_2)(\hat u^w).
$$
\item
\begin{equation}\label{sderovatove131}
\left\Vert
\frac{\partial}{\partial x} \hat u^w
\right\Vert_{L^2_{m+1,\delta}(K_i^{+S})}
+
\left\Vert
\frac{\partial}{\partial y} \hat u^w
\right\Vert_{L^2_{m+1,\delta}(K_i^{+S})}
\le \Phi_3(\epsilon)
\end{equation}
with $\lim_{\epsilon \to 0}\Phi_3(\epsilon) = 0$.
\item $\hat u^w = u^w$ for $w\in \partial D^2$.
\end{enumerate}
\end{lem}
\begin{proof}
Run the alternating method described in Subsection \ref{alternatingmethod}
in two parameter family version.
\end{proof}
\begin{lem}\label{liftinglemma}
If $T$ is sufficiently large and $\epsilon$ is sufficiently small, there exists
a smooth map $F : D^2 \to V_1\times_L V_2$ such that
\begin{enumerate}
\item
$\text{\rm Glu}_T(F(w)) = \hat u^{w}$.
\item If $s \in [0,1]$ then we have:
$$
F(e^{2\pi\sqrt{-1}s}) =\rho^s.
$$
\end{enumerate}
\end{lem}
\begin{proof}
Note that $\rho \mapsto \text{\rm Glu}_T(\rho)$ is a local diffeomorphism.
So we can apply the proof of homotopy lifting property as follows.
Let $D_{r}^2
= \{z \in \C \mid \vert z - (r-1) \vert \le r\}$.
We put
$$
A = \{r \in [0,1] \mid \text{$\exists$ $F : D^2_r \to V_1\times_L V_2$
satisfying (1) above and $F(-1) = \rho^{1/2}$}\}.
$$
Since $\text{\rm Glu}_T(\rho)$ is a local diffeomorphism,
$A$ is open. We can use (\ref{sderovatove131}) to show
closedness of $A$.
Since $0\in A$, it follows that $1\in A$. The proof of
Lemma \ref{liftinglemma} is complete.
\end{proof}
The proof of Theorem \ref{gluethm1} is now complete.
\end{proof}
\par
\newpage
\part{Construction of the Kuranishi structure 2: Construction in the
general case}\label{generalcase}

\section{Graph associated to a stable map}
\label{Graph}

We first recall the definition of the moduli space of (bordered) stable maps of genus zero.
\begin{defn}
Let $\beta \in H_2(X,L;\Z)$ and
$k,\ell \ge 0$.
The {\it  compactified moduli space
of pseudo-holomorphic disks with $k+1$ boundary marked points and $\ell$ interior marked points with boundary condition
given by $L$} that we denote by
$\mathcal M_{k+1,\ell}(\beta)$ is the set of
equivalence classes of $((\Sigma,\vec z,\vec z^{\text{\rm int}}),u)$, where:
\begin{enumerate}
\item
$\Sigma$ is a bordered semi-stable curve of genus zero with one boundary component $\partial\Sigma$.
\item
$u : (\Sigma,\partial\Sigma) \to (X,L)$ is a pseudo-holomorphic map of homology class $\beta$.
\item
$\vec z = (z_0,\dots,z_k)$ are boundary marked points.
None of them are singular points and they are all distinct.
We assume that they respect the cyclic order of $\partial\Sigma$.
\item
$\vec z^{\text{\rm int}} = (z^{\text{\rm int}}_1,\dots,z^{\text{\rm int}}_{\ell})$ are interior marked points of $\Sigma$.
None of them are singular points and they are all distinct.
\end{enumerate}
\par
We say $((\Sigma,\vec z,\vec z^{\text{\rm int}}),u)$ is {\it  equivalent} to
$((\Sigma',\vec z',\vec z^{\text{\rm int} \prime}),u')$ if there exists a biholomorphic map
$v : \Sigma' \to \Sigma$  such that $u\circ v = u'$ and $v(z'_i) = z_i$,   $v(z^{\text{\rm int} \prime}_i) = z^{\text{\rm int}}_i$.
\end{defn}
\begin{defn}
Let $\alpha \in H_2(X;\Z)$ and
$\ell \ge 0$.
The {\it compactified moduli space
of pseudo-holomorphic sphere with $\ell$ (interior) marked points} that we denote by
$\mathcal M_{\ell}^{\rm cl}(\alpha)$ is the set of the
equivalence classes of $((\Sigma,\vec z^{\text{\rm int}}),u)$, where:
\begin{enumerate}
\item
$\Sigma$ is a semi-stable curve of genus zero without  boundary.
\item
$u : \Sigma \to X$ is a pseudo-holomorphic map of homology class $\alpha$.
\item
$\vec z^{\text{\rm int}} = (z^{\text{\rm int}}_1,\dots,z^{\text{\rm int}}_{\ell})$ are marked points of $\Sigma$.
None of them are singular points and they are all distinct.
\end{enumerate}
\par
We say $((\Sigma,\vec z^{\text{\rm int}}),u)$ is {\it equivalent} to
$((\Sigma',\vec z^{\text{\rm int} \prime}),u')$ if there exists a biholomorphic map
$v : \Sigma' \to \Sigma$  such that $u\circ v = u'$ and $v(z^{\text{\rm int} \prime}_i) = z^{\text{\rm int}}_i$.
\end{defn}
\par
The topology of $\mathcal M_{\ell}^{\rm cl}(\alpha)$ is defined in \cite[Definition 10.3]{FOn}
and the topology of $\mathcal M_{k+1,\ell}(\beta)$ is defined in \cite[Definition 7.1.42]{fooo:book1}.
(See Definition \ref{convdefn}.)
\par
It is proved in \cite[Theorem 11.1 and Lemma 10.4]{FOn} that
$\mathcal M_{\ell}^{\rm cl}(\alpha)$
is compact and Hausdorff.
$\mathcal M_{k+1,\ell}(\beta)$  is also compact and Hausdorff. See \cite[Theorem 7.1.43]{fooo:book1} and the
references therein.
\par
We refer \cite[Section 2.1]{fooo:book1} for the moduli space $\mathcal M_{k+1,\ell}(\beta)$.
See also \cite{Liu}.
\par
We consider the case when $X$ is a point and denote the moduli space of that case by $\mathcal M_{k+1,\ell}$.
We call it Deligne-Mumford moduli space. (This is a slight abuse of notation since
Deligne-Mumford studied the case when there is no boundary.)
We define $\mathcal M_{\ell}^{\text{\rm cl}}$ in the same way.
\par
\begin{thm}\label{existsKura}
$\mathcal M_{\ell}^{\rm cl}(\alpha)$ has a Kuranishi structure (without boundary) and
$\mathcal M_{k+1,\ell}(\beta)$ has a Kuranishi structure with corners.
\end{thm}
\begin{rem}
\begin{enumerate}
\item
Theorem \ref{existsKura} in case of $\mathcal M_{\ell}^{\rm cl}(\alpha)$ is a special case of
\cite[Theorem 7.10]{FOn}.
In the case of $\mathcal M_{k+1,\ell}(\beta)$, Theorem \ref{existsKura}
is \cite[Theorem 2.1.29]{fooo:book1}.
\item
In the case of $\mathcal M_{k+1,\ell}(\beta)$
we need to describe the way how various moduli spaces with different $k$, $\ell$, $\beta$ are related
along their boundaries and corners, for the application. See \cite[Proposition 7.1.2]{fooo:book1}
for the precise statement on this point.
It is easy to see that the proof we will give in this note implies that version.
\end{enumerate}
\end{rem}
Below we give a detailed proof of Theorem \ref{existsKura}. The proof is
based on the proof in \cite{FOn}. The smoothness of coordinate at infinity
is useful especially in the case of $\mathcal M_{k+1,\ell}(\beta)$. On that point
we follow the method of \cite[Section 7.2 and Appendix A1.4]{fooo:book1}.
\begin{rem}
We discuss the case of genus zero here. We can handle the case of
moduli space of pseudo-holomorphic curves with or without boundary and of arbitrary genus and with arbitrary number
of boundary components, in the same way.  The case of multi Lagrangian submanifolds in pairwise clean intersection
can be also handled in the same way.
To slightly simplify the notation we restrict ourselves to the case of disks, that is
mainly used in our book \cite{fooo:book1} and spheres, that is asked in the google group `Kuranishi' explicitly.
In fact {\it no} new idea is required for generalization to higher genus etc. as far as the construction of Kuranishi structure concerns.
\end{rem}
In a way similar to \cite[Section 8]{FOn}, we stratify $\mathcal M_{k+1,\ell}(\beta)$ as follows.
For each element $\frak p = [(\Sigma,\vec z,\vec z^{\text{\rm int}}),u]$ of $\mathcal M_{k+1,\ell}(\beta)$
we associate $\mathcal G = \mathcal G_{\frak p}$, a graph with some extra data, as follows.
\par
A vertex $\rm v$ of $\mathcal G$ corresponds to $\Sigma_{\rm v}$ an irreducible component of $\Sigma$.
(It is either a disk or a sphere.)
We put data $\beta_{\rm v} = [u\vert _{\Sigma_{\rm v}}]$ that is either an element of $H_2(X,L;\Z)$ or
an element of $H_2(X;\Z)$.
\par
To each singular point $z$ of $\Sigma$ we associate an edge ${\rm e}_z$ of  $\mathcal G$.
The edge ${\rm e}_z$ joins two vertices ${\rm v}_1,{\rm v}_2$ such that $z \in \Sigma_{{\rm v}_i}$.
Note $z$ can be either boundary or interior singular points.
We also denote by $z_{\rm e}$ the singular point of $\Sigma$ corresponding to the edge $\rm e$.
\par
For each vertex ${\rm v}$ we also include the data which marked points
are contained in $\Sigma_{\rm v}$.
\begin{defn}
We call a graph $\mathcal G$ equipped with some other data described above, the
{\it  combinatorial type} of $\frak p = [(\Sigma,\vec z,\vec z^{\text{\rm int}}),u]$.
We denote by $\mathcal M_{k+1,\ell}(\beta;\mathcal G)$ the set of
$\frak p$ with combinatorial type $\mathcal G$.
\par
We write $\overset{\circ}{\mathcal M}_{k+1,\ell}(\beta)$ the stratum $\mathcal M_{k+1,\ell}(\beta;{\rm pt})$,
where ${\rm pt}$ is a graph without edge.\footnote{$\overset{\circ}{\mathcal M}_{k+1,\ell}(\beta)$ is
slightly smaller than the `interior' of ${\mathcal M}_{k+1,\ell}(\beta)$. Namely
elements of $\overset{\circ}{\mathcal M}_{k+1,\ell}(\beta)$ do not contain any disk or sphere bubble.
On the other hand, elements of the interior of ${\mathcal M}_{k+1,\ell}(\beta)$ may contain sphere bubble.}
\par
We say that $\mathcal G$ is {\it stable} if corresponding
pseudo-holomorphic curve is stable.
We say that $\mathcal G$ is {\it source stable} if
the marked bordered curve obtained by forgetting the map
is stable.
\end{defn}
Let $\mathcal G$ and $\mathcal G'$ be combinatorial types.
We say $\mathcal G \succ \mathcal G'$ if $\mathcal G'$ is obtained
from $\mathcal G$ by iterating the following process finitely many times.
\par
Take an edge ${\rm e}$ of $\mathcal G$. We shrink $\rm e$ and identify two
vertices ${\rm v}_1$, ${\rm v}_2$ contained in $\rm e$.
Let ${\rm v}$ be the vertex identified to ${\rm v}_1$, ${\rm v}_2$.
We put $\beta_{\rm v} =\beta_{{\rm v}_1} + \beta_{{\rm v}_2}$.
The marked points assigned to ${\rm v}_1$ or ${\rm v}_2$ will be
assigned to ${\rm v}$.
\begin{lem}
If
$$
\overline{\mathcal M_{k+1,\ell}(\beta;\mathcal G)}
\cap
\mathcal M_{k+1,\ell}(\beta;\mathcal G')
\ne \emptyset,
$$
then $\mathcal G \succ \mathcal G'$.
\end{lem}
The proof is easy so omitted.
\par
Sometimes we add the following data to $\mathcal G$.
\begin{enumerate}
\item
Orientation to each of the edge.
We call that $\mathcal G$ is {\it oriented} in case we include this data.\footnote{Actually in our case
of genus $0$ with at least one marked point there is a canonical way to
orient the edges as follows. We remove $z_{\rm e}$ from $\Sigma$. Then there is a
component which contains the $0$-th boundary marked point (or first interior marked point if $\partial \Sigma = \emptyset$).
If $\rm v$ is a vertex
contained in $\rm e$ we orient $\rm e$ so that $\rm v$ is inward if and only if the
corresponding irreducible component is in the connected component of $\Sigma$ minus boundary marked points
that contains $0$-th boundary marked point.}
\item
The length $T_{\rm e} \in \R_{>0}$ to each of the edges ${\rm e}$.
\end{enumerate}
We say an edge $\rm e$ is an {\it outgoing edge} of its vertex $\rm v$
and {\it incoming edge} of its vertex ${\rm v}'$ if the orientation of $\rm e$ is
goes from $\rm v$ to ${\rm v}'$.
By an abuse of terminology we say
$\rm v$ is an {\it incoming vertex} (resp. outgoing vertex) of the $\rm e$ if
$\rm e$ is an {\it incoming edge} (resp. outgoing edge) of $\rm v$.
\footnote{This might be different from the usual meaning of the English word
incoming and outgoing.}
\par
We use the following notation.
\par\medskip
$C^0_{\mathrm d}(\mathcal G) = $ the set of the vertices that correspond to a disk component.
\par
$C^0_{\mathrm  s}(\mathcal G) = $ the set of the vertices that correspond to a sphere component.
\par
$C^0(\mathcal G) = C^0_{\mathrm  d}(\mathcal G) \cup C^0_{\mathrm  s}(\mathcal G)$.
\par
$C^1_{\mathrm  o}(\mathcal G) = $ the set of the edges that correspond to a boundary singular point.
\par
$C^1_{\mathrm  c}(\mathcal G) = $ the set of the edges that correspond to an interior singular point.
\par
$C^1(\mathcal G) = C^1_{\mathrm  o}(\mathcal G) \cup C^1_{\mathrm  c}(\mathcal G)$.
\par\medskip
Here d,s,o,c indicate disk, sphere, open (string), closed (string), respectively.
\par
We define moduli space of marked stable maps from genus zero curve {\it without} boundary in the same way.
We denote it by $\mathcal M_{\ell}^{\text{\rm cl}}(\alpha)$ where $\alpha \in H_2(X;\Z)$.
($\ell$ is the number of (interior) marked points.)
In the same way we can associate a combinatorial type to it that is a graph $\mathcal G$.
In this case there is no $C^0_{\mathrm d}(\mathcal G)$ or $C^1_{\mathrm  o}(\mathcal G)$.
We define $\mathcal M_{\ell}^{\text{\rm cl}}(\alpha;\mathcal G)$,
$\overset{\circ}{\mathcal M}_{\ell}^{\text{\rm cl}}(\alpha)$,
in the same way.
\par
Let us introduce some more notations.
Let $\frak p \in \mathcal M_{k+1,\ell}(\beta)$. We put
$$
\frak p = (\frak x,u) = ((\Sigma,\vec z,\vec z^{\rm int}),u).
$$
Then we sometimes write $\frak x = \frak x_{\frak p}$, $\Sigma = \Sigma_{\frak p} =\Sigma_{\frak x}$,
$\vec z = \vec z_{\frak p} =\vec z_{\frak x}$, $\vec z^{\rm int} = \vec z^{\rm int}_{\frak p} =\vec z^{\rm int}_{\frak x}$.
We also write $u = u_{\frak p}$.
We use a similar notation in case  $\frak p \in \mathcal M_{\ell}^{\text{\rm cl}}(\alpha)$.
\begin{defn}
We put
\begin{equation}
\aligned
\Gamma_{\frak p} = \{ v : \Sigma_{\frak p} \to \Sigma_{\frak p}
\mid &\text{$v$ is a biholomorphic map},\,\, v(z_{\frak p,i}) =z_{\frak p,i}, \\
&v(z_{\frak p,i}^{\text{\rm int}}) =z_{\frak p,i}^{\text{\rm int}},
u_{\frak p} \circ v = u_{\frak p}.\}
\endaligned
\end{equation}
\begin{equation}\label{2132}
\aligned
\Gamma^+_{\frak p} = \{ v : \Sigma_{\frak p} \to \Sigma_{\frak p}
\mid &\text{$v$ is a biholomorphic map},\,\, v(z_{\frak p,i}) =z_{\frak p,i}, \\
&\exists \sigma\in \frak S_{\ell} \,\,v(z_{\frak p,i}^{\text{\rm int}}) =z_{\frak p,\sigma(i)}^{\text{\rm int}}, \,\,
u_{\frak p} \circ v = u_{\frak p}.\}
\endaligned
\end{equation}
Here $\frak S_{\ell}$ is the group of permutations of $\{1,\dots,\ell\}$.
\par
The assignment $v \mapsto \sigma$ defines a group homomorphism
\begin{equation}\label{permrep}
\Gamma^+_{\frak p} \to \frak S_{\ell}.
\end{equation}
When $\frak H$ is a subgroup of $\frak S_{\ell}$ we denote by
$\Gamma^{\frak H}_{\frak p}$ its inverse image by (\ref{permrep}).
We denote
$$
\mathcal M_{k+1,\ell}(\beta;\frak H)= \mathcal M_{k+1,\ell}(\beta)/\frak H,
$$
where $\frak H$ acts by permutation of the interior marked points.
\par
In case $X$ is a point we write $\mathcal M_{k+1,\ell}(\frak H)$ and define the groups $\Gamma^{\frak H}_{\frak x}$,
$\Gamma^{+}_{\frak x}$ for an element
$\frak x \in \mathcal M_{k+1,\ell}$.
Note that in our case of genus zero with at least one boundary marked point,
the group $\Gamma_{\frak x}$ is trivial.
(However this fact is never used in this article.)
\par
We define a similar notion in the case of $\mathcal M_{\ell}^{\text{\rm cl}}$ etc.
\end{defn}
\par\medskip
\section{Coordinate around the singular point}
\label{coordinateinf}

Let us assume that $\mathcal G$ is an oriented combinatorial type that is source stable
and $\frak H$ is a subgroup of $\frak S_{\ell}$.
Let $\frak x = [\Sigma,\vec z,\vec z^{\text{\rm int}}] \in \mathcal M_{k+1,\ell}(\frak H)$ with combinatorial type $\mathcal G$.
It is well-known that $\mathcal M_{k+1,\ell}(\frak H)$ is an effective orbifold with boundary and corners with its local
model $\frak V({\frak x})/\Gamma^{\frak H}_{\frak x}$. Let us describe this neighborhood in more detail below.
\par
For each ${\rm v} \in C^0_{\mathrm d}(\mathcal G)$, the element $\frak x$ determines a marked disk
$\frak x_{\rm v} \in \overset{\circ}{\mathcal M}_{k_{\rm v}+1,\ell_{\rm v}}$.
Here $k_{\rm v}$ is the sum of the number of edges $\in C^1_{\mathrm  o}(\mathcal G)$ containing ${\rm v}$ and the number of boundary marked
points assigned to ${\rm v}$. $\ell_{\rm v}$ is the sum of the number of edges $\in C^1_{\mathrm  c}(\mathcal G)$ containing
${\rm v}$ and the number of interior marked
points assigned to ${\rm v}$. (In other words the singular points of $\Sigma$ that is contained in  $\Sigma_{\rm v}$ is regarded as a
marked point of $\frak x_{\rm v}$.)
\par
For each ${\rm v} \in C^0_{\mathrm s}(\mathcal G)$, the element $\frak x$ determines a marked sphere
$\frak x_{\rm v} \in \overset{\circ}{\mathcal M}_{\ell_{\rm v}}^{\text{\rm cl}}$ in the same way.
\par
Let $\frak V(\frak x_{\rm v})/\Gamma^{\frak H}_{\frak x_{\rm v}}$  be the neighborhood of $\frak x_{\rm v}$ in ${\mathcal M}
_{k_{\rm v}+1,\ell_{\rm v}}({\frak H})$
or in ${\mathcal M}_{\ell_{\rm v}}^{\text{\rm cl}}(\frak H)$, respectively,
according to whether ${\rm v} \in C^0_{\rm d}(\mathcal G)$ or ${\rm v} \in C^0_{\rm s}(\mathcal G)$.
The group $\Gamma^{\frak H}_{\frak x}$ acts on the product $\prod  \frak V(\frak x_{\rm v})$. The quotient
$$
\frak V(\frak x;\mathcal G)/\Gamma^{\frak H}_{\frak x}
= \left(\prod_{{\rm v}\in C^0(\mathcal G)}  \frak V(\frak x_{\rm v})\right) / \Gamma^{\frak H}_{\frak x}
$$
is a neighborhood of $\frak x$ in $\mathcal M_{k+1,\ell}(\mathcal G;\frak H)$.
\par
A neighborhood of $\frak x$ in $\mathcal M_{k+1,\ell}(\frak H)$ is identified with
\begin{equation}\label{nbhdofstratum}
\left(  \frak V(\frak x;\mathcal G) \times \left(\prod_{{\rm e}\in C^1_{\mathrm  o}(\mathcal G)} (T_{{\rm e},0},\infty]\right)
\times \left(\prod_{{\rm e}\in C^1_{\mathrm  c}(\mathcal G)} ((T_{{\rm e},0},\infty] \times S^1)/\sim\right)\right)/\Gamma^{\frak H}_{\frak x}.
\end{equation}
\begin{rem}\label{rem:161}
The equivalence relation $\sim$ in (\ref{nbhdofstratum}) is defined as follows.
$(T,\theta) \sim (T',\theta')$ if $(T,\theta) = (T',\theta')$ or $T=T'=\infty$.
\par
The action of $\Gamma^{\frak H}_{\frak x}$ on
$$
\left(\prod_{{\rm e}\in C^1_{\mathrm  o}(\mathcal G)} (T_{{\rm e},0},\infty]\right)
\times \left(\prod_{{\rm e}\in C^1_{\mathrm  c}(\mathcal G)} ((T_{{\rm e},0},\infty] \times S^1)/\sim\right)
$$
is by exchanging the factors associated to the edges $\rm e$ and by  rotation of the $S^1$ factors.
(See the proof of Lemma \ref{Phisiequv}.)
\end{rem}
We will define a map from (\ref{nbhdofstratum}) to $\mathcal M_{k+1,\ell}({\frak H})$.
(See Definition \ref{def214}.)
We need to fix a coordinate of $\Sigma$ around each of the singular point for this purpose.
For the sake of consistency with the analytic construction in Section \ref{secsimple}, we use
cylindrical coordinate.
\begin{defn}\label{coordinatainfdef}
Let
\begin{equation}\label{fibrationsigma}
\pi : \frak M_{\frak x_{\rm v}} \to \frak V(\frak x_{\rm v})
\end{equation}
be a fiber bundle whose fiber is a two dimensional manifold
together with fiberwise complex structure.
This fiber bundle is the universal family in the sense of (2) below.
We call (\ref{fibrationsigma}) with extra data described below
a
{\it universal family with coordinate at infinity}
if the following conditions are satisfied.
\begin{enumerate}
\item
$\frak M_{\frak x_{\rm v}}$ has a fiberwise biholomorphic  $\Gamma^+_{\frak x_{\rm v}}$ action and
$\pi$ is $\Gamma^+_{\frak x_{\rm v}}$ equivariant.
\item
For $\frak y \in  \frak V(\frak x_{\rm v})$
the fiber $\pi^{-1}(\frak y)$ is biholomorphic to $\Sigma_{\frak y}$ minus marked points
corresponding to the singular points of $\frak y$.
\item
As a part of the data we fix a closed subset  $\frak K_{\frak x_{\rm v}} \subset \frak M_{\frak x_{\rm v}}$ such that
$\pi : \frak K_{\frak x_{\rm v}} \to  \frak V(\frak x_{\rm v})$ is proper.
\item
We consider the direct product
\begin{equation}\label{endproductstri}
\aligned
 \frak V(\frak x_{\rm v})
\times
&\bigcup_{{\rm e}\in C^1_{\mathrm o}(\mathcal G)
\atop \text{${\rm e}$ is an outgoing edge of ${\rm v}$}} (0,\infty) \times [0,1]
\\
&\cup
\bigcup_{{\rm e}\in C^1_{\mathrm o}(\mathcal G)
\atop \text{${\rm e}$ is an incoming edge of ${\rm v}$}} (-\infty,0) \times [0,1]
\\
&\cup
\bigcup_{{\rm e}\in C^1_{\mathrm c}(\mathcal G)
\atop \text{${\rm e}$ is an outgoing edge of ${\rm v}$}} (0,\infty) \times S^1
\\
&\cup
\bigcup_{{\rm e}\in C^1_{\mathrm c}(\mathcal G)
\atop \text{${\rm e}$ is an incoming edge of ${\rm v}$}} (-\infty,0) \times S^1.
\endaligned
\end{equation}
(Here and hereafter the symbols $\cup$ and $\bigcup$ in (\ref{endproductstri}) are the {\it disjoint} union.) \par
As a part of the data we fix a diffeomorphism between $\frak M_{\frak x_{\rm v}}\setminus \frak K_{\frak x_{\rm v}}$ and (\ref{endproductstri})
that commutes with the projection to $ \frak V(\frak x_{\rm v})$ and is a fiberwise biholomorphic map.
Moreover the diffeomorphism sends each end
corresponding to a singular point $z_{\rm e}$ to the end in (\ref{endproductstri})
corresponding to the edge ${\rm e}$.
\item
The diffeomorphism in (4) extends to a fiber preserving diffeomorphism
$$
\frak M_{\frak x_{\rm v}} \cong   \frak V(\frak x_{\rm v}) \times (\Sigma_{\frak x_{\rm v}} \setminus \{\text{singular points}\}).
$$
This diffeomorphism sends each of the interior or boundary marked points of the fiber of
$\frak y$ to the corresponding marked point of
$\{\frak y\} \times \Sigma_{\frak x_{\rm v}}$.
However, this diffeomorphism does {\it not} preserve fiberwise complex structure.
As a part of the data we fix this extension of diffeomorphism.
\item
The action of an element of $\Gamma^+_{\frak x_{\rm v}}$ on (\ref{endproductstri}) is given by exchanging the factors associated to the edges $\rm e$ and by  rotation of the $S^1$ factors.
\end{enumerate}
Hereafter we sometimes call a {\it coordinate at infinity} in place of a universal family with coordinate at infinity.
\end{defn}
\begin{exm}
Let $\frak x_{\text{\rm v}}$ be $S^2$ with $\ell + 2 $ marked points
$$
z_0=0, z_1 = \infty, z_2=1, \dots, z_{\ell+1} = e^{2\pi\sqrt{-1}(\ell-1)/\ell}.
$$
Let $\frak H \subset \frak S_{\ell+2}$ be the subgroup $\frak S_{\ell}$ consisting of elements
that fix $z_0, z_1$.
We assume that $z_0$ and $z_1$ correspond to singular points of $\frak x$.
It is easy to see that $\Gamma_{\frak x}^{\frak H} = \Z_{\ell}$.
Then $\Sigma_{\frak x_{\rm v}} \setminus \{z_0,z_1\} = \R \times S^1$
and the action of $\Gamma_{\frak x}^{\frak H}$ is given by rotation of the $S^1$ factors.
\end{exm}
\begin{defn}\label{defn288}
Suppose we are given a coordinate at infinity for each of $\frak x_{\rm v}$
where $\frak x_{\rm v}$ corresponds to an irreducible component of $\frak x$. We say that they
are {\it invariant under the $\Gamma^+_{\frak x}$-action} if the following holds.
\par
We define a fiber bundle
\begin{equation}\label{2149}
\pi : \underset{{\rm v}\in C^0(\mathcal G)}{\bigodot}\frak M_{\frak x_{\rm v}} \to \prod_{{\rm v}\in C^0(\mathcal G)}\frak V(\frak x_{\rm v})
\end{equation}
as follows.
We take projections
$\prod_{{\rm v}\in C^0(\mathcal G)}\frak V(\frak x_{\rm v}) \to \frak V(\frak x_{\rm v})$
and pull back the bundle (\ref{fibrationsigma}) by this projection.
We thus obtain a fiber bundle over $\prod_{{\rm v}\in C^0(\mathcal G)}\frak V(\frak x_{\rm v})$.
(\ref{2149}) is the disjoint union of those bundles over ${\rm v}\in C^0(\mathcal G)$.
In other words the fiber of
(\ref{2149}) at $(\frak y_{\rm v} : {\rm v} \in C^0(\mathcal G))$ is a disjoint union
of $\frak y_{\rm v}$'s.
\par
The fiber bundle (\ref{2149}) has a $\Gamma^+_{\frak x_{\rm v}}$-action.
We consider its restriction to
\begin{equation}\label{2endproductstri}
\pi : \underset{{\rm v}\in C^0(\mathcal G)}{\bigodot}(\frak M_{\frak x_{\rm v}}\setminus \frak K_{\frak x_{\rm v}})
\to \prod_{{\rm v}\in C^0(\mathcal G)}\frak V(\frak x_{\rm v}).
\end{equation}
The group $\Gamma^+_{\frak x}$ acts on the sum of the second factors of (\ref{endproductstri})
by exchanging the factors associated to the edges $\rm e$ and by  rotation of the $S^1$ factors.
We require that (\ref{2endproductstri}) is invariant under this action.
\par
Moreover we assume that the diffeomorphisms in Definition \ref{coordinatainfdef} (4)(5)
are $\Gamma^+_{\frak x}$ equivariant.
\end{defn}
Now we fix a coordinate at infinity for each of $\frak x_{\rm v}$ that is invariant under the
$\Gamma^{\frak H}_{\frak x}$ action.
We will use it to define a map from (\ref{nbhdofstratum}) to a neighborhood of $\frak x$ in $\mathcal M_{k+1,\ell}({\frak H})$ as follows.
Let  $(\frak y_{\rm v} : {\rm v} \in C^0(\mathcal G))$ and $\frak y_{\rm v} \in \frak V(\frak x_{\rm v})$.
Take a representative  $\Sigma_{\frak y_{\rm v}}$ of $\frak y_{\rm v}$.
We put $K_{\frak y_{\rm v}} = \Sigma_{\frak y_{\rm v}}\cap \frak K_{\frak x_{\rm v}}$.
The coordinate at infinity defines a biholomorphic map between
$\bigcup_{ {\rm v} \in C^0(\mathcal G))}\Sigma_{\frak y_{\rm v}} \setminus K_{\rm v}$ and
\begin{equation}\label{endidentify}
\aligned
&\bigcup_{{\rm e}\in C^1_{\mathrm o}(\mathcal G)
\atop \text{${\rm e}$ is an outgoing edge of ${\rm v}$}} (0,\infty) \times [0,1]
\\
&\cup
\bigcup_{{\rm e}\in C^1_{\mathrm o}(\mathcal G)
\atop \text{${\rm e}$ is an incoming edge of ${\rm v}$}} (-\infty,0) \times [0,1]
\\
&\cup
\bigcup_{{\rm e}\in C^1_{\mathrm c}(\mathcal G)
\atop \text{${\rm e}$ is an outgoing edge of ${\rm v}$}} (0,\infty) \times S^1
\\
&\cup
\bigcup_{{\rm e}\in C^1_{\mathrm c}(\mathcal G)
\atop \text{${\rm e}$ is an incoming edge of ${\rm v}$}} (-\infty,0) \times S^1.
\endaligned
\end{equation}
We write the coordinate of each summand of (\ref{endidentify})
by $(\tau'_{\rm e},t_{\rm e})$, $(\tau''_{\rm e},t_{\rm e})$, $(\tau'_{\rm e},t'_{\rm e})$, $(\tau''_{\rm e},t''_{\rm e})$ respectively.
(Here we identify $S^1 = \R/\Z$ so $t_{\rm e} \in [0,1]$ or $t'_{\rm e}, t''_{\rm e} \in \R/\Z$.)
\par
Now, let $((T_{\rm e};{\rm e}\in C^1_{\mathrm  o}(\mathcal G)) ,((T_{\rm e},\theta_{\rm e});{\rm e}\in C^1_{\mathrm  c}(\mathcal G))$
be an element of
\begin{equation}\label{infnbfparam}
 \left(\prod_{{\rm e}\in C^1_{\mathrm  o}(\mathcal G)} (T_{{\rm e},0},\infty]\right)
\times \left(\prod_{{\rm e}\in C^1_{\mathrm  c}(\mathcal G)} ((T_{{\rm e},0},\infty] \times S^1)/\sim\right).
\end{equation}
(Here $\theta_{\rm e} \in \R/\Z$.)
\begin{defn}\label{def29}
We denote the right hand side of (\ref{infnbfparam}) by
$(\vec T^{\rm o}_0,\infty] \times ((\vec T^{\rm c}_0,\infty] \times \vec S^1)$.
\end{defn}
\par
We first consider the case $T_{\rm e} \ne \infty$.
We define $\tau_{\rm e}$ for ${\rm e} \in C^1(\mathcal G)$ and $t_{\rm e}$ for ${\rm e}\in C^1_{\mathrm c}(\mathcal G)$
as follows.
\begin{eqnarray}
\tau_e &=& \tau'_{\rm e} - 5T_{\rm e}  = \tau''_{\rm e} + 5T_{\rm e}, \label{cctau1}\\
t_{\rm e} &=& t'_{\rm e} = t''_{\rm e} - \theta_{\rm e}.\label{ccttt1}
\end{eqnarray}
We note that (\ref{cctau1}), (\ref{ccttt1}) are consistent with the  notation of Section \ref{subsecdecayT}.
We consider
\begin{equation}\label{neckopen}
[-5T_{\rm e},5T_{\rm e}] \times [0,1]
\end{equation}
for each ${\rm e} \in  C^1_{\mathrm  o}(\mathcal G)$ with coordinate $(\tau_{\rm e},t_{\rm e})$
and
\begin{equation}\label{neckclosed}
[-5T_{\rm e},5T_{\rm e}] \times S^1
\end{equation}
for each ${\rm e} \in  C^1_{\mathrm  c}(\mathcal G)$ with coordinate $(\tau_{\rm e},t_{\rm e})$.
\par
We now consider the union
\begin{equation}\label{summandtoglued}
\aligned
\bigcup_{{\rm v} \in C^0(\mathcal G)}K_{\frak y_{\rm v}}
&\cup \bigcup_{{\rm e}\in C^1_{\mathrm o}(\mathcal G)} [-5T_{\rm e},5T_{\rm e}] \times [0,1]
\\
&\cup
\bigcup_{{\rm e}\in C^1_{\mathrm c}(\mathcal G)} [-5T_{\rm e},5T_{\rm e}] \times S^1.
\endaligned
\end{equation}
(\ref{cctau1}) and (\ref{ccttt1}) describe the way how we glue various summands in
(\ref{summandtoglued}) to obtain a bordered Riemann surface, that is nonsingular
in our case where $T_e \ne \infty$.
\begin{defn}\label{def214}
We denote by ${\overline\Phi}((\frak y_{\rm v};{\rm v} \in C^0(\mathcal G)),
(T_{\rm e};{\rm e}\in C^1_{\mathrm  o}(\mathcal G)),(T_{\rm e},\theta_{\rm e});{\rm e}\in C^1_{\mathrm  c}(\mathcal G))$
the element of $\mathcal M_{k+1,\ell}$ represented by the above bordered Riemann surface.
\par
Hereafter we write ${\frak y}= (\frak y_{\rm v};{\rm v} \in C^0(\mathcal G))$,
$\vec T^{\mathrm  o} = (T_{\rm e};{\rm e}\in C^1_{\mathrm  o}(\mathcal G))$,
$\vec T^{\mathrm c} = (T_{\rm e};{\rm e}\in C^1_{\mathrm  c}(\mathcal G))$, and
$\vec \theta = (\theta_{\rm e};{\rm e}\in C^1_{\mathrm  c}(\mathcal G))$.
We put $\vec T = (\vec T^{\mathrm  o},\vec T^{\mathrm  c})$.
We denote ${\overline{\Phi}}({\frak y},\vec T^{\mathrm  o},(\vec T^{\mathrm c},\vec{\theta}))
= {\overline{\Phi}}({\frak y},\vec T,\vec \theta)
\in \mathcal M_{k+1,\ell}$.
\par\medskip
We next consider the case when some $T_{\rm e} =\infty$.
We define a graph $\mathcal G'$ as follows :
We shrink all the edges ${\rm e}$ of $\mathcal G$ with $T_{\rm e} \ne \infty$.
Various data we associate to $\mathcal G'$ are induced by
the one associated to $\mathcal G$ in an obvious way.
The element
${\overline{\Phi}}({\frak y},\vec T^{\mathrm  o},(\vec T^{\mathrm  c},\vec{\theta}))$
is contained in $\mathcal M_{k+1,\ell}(\mathcal G')$.
Namely we glue  (\ref{summandtoglued}) to obtain a (noncompact) bordered Riemann surface $\Sigma'$.
Then we add a finite number of points (each corresponds to the edges with infinite length)
to obtain (singular) stable bordered curve ${\overline{\Phi}}({\frak y},\vec T^{\mathrm  o},(\vec T^{\mathrm  c},\vec{\theta}))$ such that
${\overline{\Phi}}({\frak y},\vec T^{\mathrm  o},(\vec T^{\mathrm  c},\vec{\theta}))$ minus singular points is $\Sigma'$.
\par
Thus we have defined
$$
{\overline\Phi} : \prod_{{\rm v}\in C^0(\mathcal G)}\frak V(\frak x_{\rm v}) \times
(\vec T^{\rm o}_0,\infty] \times ((\vec T^{\rm c}_0,\infty] \times \vec S^1)
\to
\mathcal M_{k+1,\ell}.
$$
\end{defn}
We define some terminology below.
\begin{defn}\label{defcoreandneck}
We call $K_{\frak y_{\rm v}}$ as in (\ref{summandtoglued}) a component of the {\it core}  of $\frak y$ or of ${\overline{\Phi}}({\frak y},\vec T^{\mathrm  o},(\vec T^{\mathrm  c},\vec{\theta}))$.
Each of the connected component of the second or third term of  (\ref{summandtoglued}) is
called a component of the {\it neck region}.
In case $T_{\rm}$ is infinity, there is a domain identified with
$([0,\infty) \cup (-\infty,0]) \times [0,1]$ or with $([0,\infty) \cup (-\infty,0]) \times S^1$ corresponding
to it. We call it also a component of the neck region. The union of all the components of the core and the neck region is ${\overline{\Phi}}({\frak y},\vec T^{\mathrm  o},(\vec T^{\mathrm  c},\vec{\theta}))$
minus singular points.
\end{defn}
\begin{rem}\label{rem216}
Note that $\mathcal M_{k+1,\ell}$ has an $\frak S_{\ell}$ action by permutation of the
interior marked points. A local chart of $\mathcal M_{k+1,\ell}$ at $\frak x$ is of the form $\frak V/\Gamma_{\frak x}$,
and a local chart of $\mathcal M_{k+1,\ell}/\frak S_{\ell}$ at $[\frak x]$ is of the form $\frak V/\Gamma^+_{\frak x}$.
\end{rem}
\begin{lem}\label{Phisiequv}
The map ${\overline\Phi}$ is $\Gamma^{+}_{\frak x}$ equivariant.
\end{lem}
\begin{proof}
We first define a $\Gamma^{+}_{\frak x}$ action on (\ref{infnbfparam}).
Note an element of $\Gamma^{+}_{\frak x}$ acts on the graph $\mathcal G$ in an obvious way.
So it determines the way how to exchange the factors of (\ref{infnbfparam}).
The rotation part of the action is defined as follows.
By Definition \ref{coordinatainfdef} (6) we can determine the rotation
of the $t_{\rm e}$ coordinate induced by an element of $\Gamma^{+}_{\frak x}$.
Therefore by (\ref{ccttt1}) the action on $\theta_{\rm e}$ coordinate is determined.
\par
Once we defined $\Gamma^{+}_{\frak x}$ action on (\ref{infnbfparam}) the
equivariance of the map  ${\overline{\Phi}}$  is immediate from definition.
\end{proof}
Note that the space (\ref{nbhdofstratum}) has a stratification.
(This stratification is induced by the stratification of $(0,\infty]$ that consists of
$(0,\infty)$ and $\{\infty\}$.
The map ${\overline{\Phi}}$ respects this stratification and stratification of $\mathcal M_{k+1,\ell}$ by $\{\mathcal M_{k+1,\ell}(\mathcal G)\}$.
Moreover ${\overline{\Phi}}$ is continuous and strata-wise smooth.
We do not discuss the smooth structure of (\ref{nbhdofstratum})  yet. (See Section  \ref{chart}.)
\par
We remark that the map ${\overline{\Phi}}$ {\it depends} on the choice of coordinate at infinity.
The next result describes how ${\overline{\Phi}}$ depends on the choice of coordinate at infinity.
\par
Let
\begin{equation}
\overline{\Phi}_1 :
\prod_{{\rm v}\in C^0(\mathcal G)}\frak V^{(1)}(\frak x_{\rm v}) \times
(\vec T^{\rm o}_0,\infty] \times ((\vec T^{\rm c}_0,\infty] \times \vec S^1)
\to
\mathcal M_{k+1,\ell}
\end{equation}
be the map in Definition \ref{def214}.
Suppose
$$
\frak Y_0 = \overline{\Phi}_1(\frak y_0,\vec T_{\frak Y_0},\vec \theta_{\frak Y_0})
$$
and $\mathcal G_{\frak Y_0}$ is the combinatorial type of $\frak Y_0$.
Note $\mathcal G_{\frak Y_0}$ is obtained from $\mathcal G_{\frak x}$ by shrinking
several edges. Therefore we may regard
$$
C^1(\mathcal G_{\frak Y_0}) \subseteq C^1(\mathcal G_{\frak x}).
$$
Namely we can canonically identify ${\rm e} \in C^1(\mathcal G_{\frak x})$ with an element of
${\rm e} \in C^1(\mathcal G_{\frak Y_0})$ if
$T_{\frak Y_0,{\rm e}} = \infty$.
\par
We take a coordinate at infinity of $\frak Y_0$. By Definition \ref{def214} it determines an embedding
\begin{equation}
\overline{\Phi}_2 :
\prod_{{\rm v}\in C^0(\mathcal G_{\frak Y_0})}\frak V^{(2)}(\frak Y_{0,\rm v}) \times
(\vec T^{\rm o}_1,\infty] \times ((\vec T^{\rm c}_1,\infty] \times \vec S^1)
\to
\mathcal M_{k+1,\ell}.
\end{equation}
Here an element of $(\vec T^{\rm o}_1,\infty] \times ((\vec T^{\rm c}_1,\infty] \times \vec S^1)$
is  $((T_{\rm e};{\rm e} \in C^1_{\rm o}(\mathcal G_{\frak Y_0}),((T_{\rm e},\theta_{\rm e}),
{\rm e} \in C^1_{\rm c}(\mathcal G_{\frak Y_0}))$.
\par
We put
\begin{equation}
\overline{\Phi}_{12} = \overline{\Phi}_1^{-1} \circ \overline{\Phi}_2.
\end{equation}
We next define $\Psi_{12}$.
Let $(\frak z_{\rm v}) \in \prod_{{\rm v}\in C^0(\mathcal G_{\frak Y_0})}\frak V^{(2)}(\frak Y_{0,\rm v})$.
We denote $\vec{\infty} \in (\vec T^{\rm o}_1,\infty] \times ((\vec T^{\rm c}_1,\infty] \times \vec S^1)$
to be the point whose components are all $\infty$.
Then
$\overline{\Phi}_2((\frak z_{\rm v}),\vec\infty)$ has the same combinatorial type
$\mathcal G_{\frak Y_0}$ as $\frak Y_0$.
We define $\Psi_{12}^{\frak y}((\frak z_{\rm v}))\in \prod_{{\rm v}\in C^0(\mathcal G)}\frak V^{(1)}(\frak x_{\rm v})$
and $\vec T',\vec\theta'$
by
$$
\overline{\Phi}_1^{-1}(\overline{\Phi}_2((\frak z_{\rm v}),\vec\infty))
=
(\Psi_{12}^{\frak y}((\frak z_{\rm v})),\vec T',\vec\theta').
$$
We note that $T'_{\rm e} = \infty$ if ${\rm e} \in C^1(\mathcal G_{\frak Y_0}) \subset C^1(\mathcal G_{\frak x})$.
Then we put
\begin{equation}\label{2167}
\Psi_{12}((\frak z_{\rm v}),\vec T,\vec\theta)
=
({\Psi}^{\frak y}_{12}((\frak z_{\rm v})),\vec T'',\vec\theta'')
\end{equation}
where
$$
T''_{\rm e} =
\begin{cases}
T_{\rm e}  &\text{if ${\rm e} \in C^1(\mathcal G_{\frak Y_0})$} \\
T'_{\rm e}      &\text{if ${\rm e} \in C^1(\mathcal G_{\frak x}) \setminus C^1(\mathcal G_{\frak Y_0})$},
\end{cases}
$$
$$
\theta''_{\rm e} =
\begin{cases}
\theta_{\rm e}  &\text{if ${\rm e} \in C_{\rm c}^1(\mathcal G_{\frak Y_0})$} \\
\theta'_{\rm e}      &\text{if ${\rm e} \in C_{\rm c}^1(\mathcal G_{\frak x}) \setminus C_{\rm c}^1(\mathcal G_{\frak Y_0})$}.
\end{cases}
$$
\begin{rem}
If $\frak Y_0$ has the same combinatorial type as $\frak x$ then $\Psi_{12}$ is the identity map.
Note that even in the case $\frak Y_0=\frak x$ the map $\overline{\Phi}_{12}$ may not be the identity
map since $\overline{\Phi}_j$ depends on the choice of coordinate at infinity.
\end{rem}
\par
Let $k_{T,{\rm e}} =0,1,\dots$, $k_{\theta,{\rm e}} = 0,1,2,\dots$
and define
$$
\frac{\partial^{\vert \vec k_{T}\vert}}{\partial T^{\vec k_{T}}}
= \prod_{{\rm e}\in C^1(\mathcal G)}\frac{\partial^{ k_{T,{\rm e}}}}{\partial T_{\rm e}^{k_{T,{\rm e}}}}.
$$
We define
$\frac{\partial^{\vert \vec k_{\theta}\vert}}{\partial T^{\vec k_{\theta}}}$ in the same way.
We put
$$
\vec k_{T}\cdot \vec T = \sum_{{\rm e}\in C^1(\mathcal G)} k_{T,{\rm e}}T_{\rm e},
\quad
\vec k_{\theta}\cdot \vec T^{\rm c} = \sum_{{\rm e}\in C^1_{\rm c}(\mathcal G)} k_{\theta,{\rm e}}T_{\rm e}.
$$
\begin{prop}\label{changeinfcoorprop}
In the above situation we have the following inequality
for any compact subset $\frak V_0(\frak x,\mathcal G)$ of $\frak V(\frak x,\mathcal G)$ :
\begin{equation}
\left\Vert
\frac{\partial^{\vert \vec k_{T}\vert}}{\partial T^{\vec k_{T}}} \frac{\partial^{\vert \vec k_{\theta}\vert}}{\partial \theta^{\vec k_{\theta}}}  ({\overline{\Phi}}_{12}
- {\Psi}_{12})
\right\Vert_{C^k} \le C_{1,k} e^{-\delta' (\vec k_{T}\cdot \vec T+\vec k_{\theta}\cdot \vec T^{\rm c})},
\label{2144}
\end{equation}
for $\vert \vec k_T\vert, \vert \vec k_{\theta}\vert \le k$
with  $\vert \vec k_T\vert + \vert \vec k_{\theta}\vert \ne 0$,
where the left hand sides are $C^k$ norm (as maps on $\frak y$)
and $\delta' > 0$ depends only on $\delta$ and $k$.
\end{prop}
\begin{rem}
The estimates in Proposition \ref{changeinfcoorprop} holds strata-wise. Namely
in the situation where some of $T_{\rm e}$ is infinity, we only consider $\vec k_{T}, \vec k_{\theta}$ such that $k_{T,\rm e}=k_{\theta,\rm e} = 0$
for the edges $\rm e$ with $T_{\rm e} = \infty$.
\end{rem}
\begin{rem}\label{metricfiber}
During the proof of Proposition \ref{changeinfcoorprop} and also
during various discussions in later sections,
we need metrics of the source and the target to define various norms etc.
For this purpose we take a Riemannian metric on $X$ and also a family of metrics
of the fibers of (\ref{fibrationsigma}) such that outside $K_{{\rm v}}$ it
coincides with the standard flat metric (via coordinates $\tau$ and $t$).
We include it in the data of universal family with coordinate at infinity.
Since we use it only to fix norm etc. it is not an important part of that data.
\end{rem}
Proposition \ref{changeinfcoorprop} is a generalization of \cite[Lemma A1.59]{fooo:book1} and will be used for the same purpose
later to derive the exponential decay estimate of the coordinate change of our Kuranishi structure.
We suspect  Proposition \ref{changeinfcoorprop} is not new. However for completeness sake the proof will be given later in Subsection \ref{proposss}.
\begin{rem}
In case $\frak Y_0 = \frak x$,
Proposition \ref{changeinfcoorprop} implies that
there exists $\vec{\Delta T} : \frak V_0(\frak x,\mathcal G) \to \R^{\# C^1(\mathcal G)}$,
$\vec{\Delta \theta} : \frak V_0(\frak x,\mathcal G) \to (S^1)^{\# C^1_{\text{\rm c}}(\mathcal G)}$
such that $T$ component (resp. $\theta$ component) of ${\overline{\Phi}}_{21}$  goes to $\vec T + \vec{\Delta T}$
(resp. $\vec {\theta}_e + \vec{\Delta \theta}$) in an exponential order as $T$ goes to infinity.
(\ref{2144}) implies that ${\frak y}$ component of ${\overline{\Phi}}_{21}$ goes to $\frak y$ in exponential
order  as $T$ goes to infinity.
\end{rem}
Proposition \ref{changeinfcoorprop} describes the coordinate change (change of the
parametrization) of the moduli space.
A coordinate at infinity determines a parametrization of the (bordered) curve itself, since
it includes the trivialization of the fiber bundle (\ref{fibrationsigma}).
Proposition \ref{reparaexpest} below describes the way how it changes when
we change the coordinate at infinity.
\par
Let  $\overline{\Phi}_{12} = \overline{\Phi}_{1}^{-1} \circ \overline{\Phi}_{2}$ be as in
Proposition \ref{changeinfcoorprop} and
let $(\frak y_j,\vec T_j,\vec \theta_j)$ ($j=1,2$) be in the domain of  $ \overline{\Phi}_{j}$.
We assume
\begin{equation}\label{2255}
(\frak y_1,\vec T_1,\vec \theta_1) = \overline{\Phi}_{12} (\frak y_2,\vec T_2,\vec \theta_2).
\end{equation}
Let $\Sigma_{(\frak y_j,\vec T_j,\vec \theta_j)}$ be a curve
representing $\overline{\Phi}_{j}(\frak y_j,\vec T_j,\vec \theta_j)$. It comes with coordinate at infinity.
By (\ref{2255}) and stability, there exists a {\it unique} isomorphism
\begin{equation}\label{2256}
\frak v_{(\frak y_2,\vec T_2,\vec \theta_2)} : \Sigma_{(\frak y_2,\vec T_2,\vec \theta_2)}
\to \Sigma_{(\frak y_1,\vec T_1,\vec \theta_1)}
\end{equation}
of marked curves.
\par
Let $K_{\rm v}^{(j)}$ be the core of $\Sigma_{(\frak y_j,\vec T_j,\vec \theta_j)}$.
We take a compact subset $K_{{\rm v},0}^{(2)} \subset K_{\rm v}^{(2)}$ such that
\begin{equation}
\frak v_{(\frak y_2,\vec T_2,\vec \theta_2)} (K_{{\rm v},0}^{(2)} ) \subset  K_{\rm v}^{(1)}
\end{equation}
for sufficiently large $\vec T_1$.
Note that the sets $K_{\rm v}^{(1)}$ and $K_{{\rm v},0}^{(2)}$ are independent of
$(\frak y_2,\vec T_2,\vec \theta_2)$.
Let
$$
C^k(K_{{\rm v},0}^{(2)},K_{\rm v}^{(1)})
$$
be the space of $C^k$ maps with $C^k$ topology.
The restriction of $\frak v_{(\frak y_2,\vec T_2,\vec \theta_2)}$ to
$K_{{\rm v},0}^{(2)}$ defines an element of it that we denote by
$$
\text{\rm Res}(\frak v_{(\frak y_2,\vec T_2,\vec \theta_2)}) \in C^k(K_{{\rm v},0}^{(2)},K_{\rm v}^{(1)}).
$$
\begin{prop}\label{reparaexpest}
There exist $C_{2,k}$, $T_k$ such that
for each ${\rm e}_0 \in C^1_{\rm c}(\mathcal G_{\frak y_2})$ we have
\begin{equation}
\aligned
&\left\Vert
\nabla_{\frak y_2}^n \frac{\partial^{\vert \vec k_{T}\vert}}{\partial T_2^{\vec k_{T}}}\frac{\partial^{\vert \vec k_{\theta}\vert}}{\partial \theta_2^{\vec k_{\theta}}}
\frac{\partial}{\partial T_{2,{\rm e}_0}}
\text{\rm Res}(\frak v_{(\frak y_2,\vec T_2,\vec \theta_2)})
\right\Vert_{C^k}
< C_{2,k}e^{-\delta_2  T_{2,{\rm e}_0}},
\\
&\left\Vert
\nabla_{\frak y_2}^n \frac{\partial^{\vert \vec k_{T}\vert}}{\partial T_2^{\vec k_{T}}}\frac{\partial^{\vert \vec k_{\theta}\vert}}{\partial \theta_2^{\vec k_{\theta}}}
\frac{\partial}{\partial \theta_{2,{\rm e}_0}}
\text{\rm Res}(\frak v_{(\frak y_2,\vec T_2,\vec \theta_2)})
\right\Vert_{C^k}
< C_{2,k}e^{-\delta_2  T_{2,{\rm e}_0}},
\endaligned
\end{equation}
if
each of $T_{2,{\rm e}}$ is greater than $T_k$ and
$\vert{\vec k_{T}}\vert + \vert{\vec k_{\theta}}\vert +n \le k$.
Here $\vec T_2 = (T_{2,{\rm e}}; {\rm e}\in  C^1(\mathcal G_{\frak y_2}))$,
$\vec \theta_2 = (\theta_{2,{\rm e}}; {\rm e}\in  C^1_{\rm c}(\mathcal G_{\frak y_2}))$.
\par
The first inequality also holds for ${\rm e}_0 \in C^1_{\rm o}(\mathcal G_{\frak y_2})$.
\end{prop}
We note that when all the numbers $T_{2,{\rm e}}$ are $\infty$,
$\overline{\Phi}_2(\frak y_2,\vec T_2,\vec \theta_2)$ has the same combinatorial type as $\frak Y_0$.
(Note $\overline{\Phi}_2$ gives a coordinate of the Deligne-Mumford moduli space
in a neighborhood of $\frak Y_0$.)
Then, integrating on $T_{2,{\rm e}}$, Proposition \ref{reparaexpest} implies:
\begin{cor}\label{corestimatecoochange}
\begin{equation}
\aligned
&\left\Vert
\nabla_{\frak y_2}^n \frac{\partial^{\vert \vec k_{T}\vert}}{\partial T_2^{\vec k_{T}}}\frac{\partial^{\vert \vec k_{\theta}\vert}}{\partial \theta_2^{\vec k_{\theta}}}
(\text{\rm Res}(\frak v_{(\frak y_2,\vec T_2,\vec \theta_2)}) - \text{\rm Res}(\frak v_{(\frak y_2,\vec \infty)})
\right\Vert_{C^k}
< C_{3,k}e^{-\delta_2  T_{2,{\rm min}}},
\\
&\left\Vert
\nabla_{\frak y_2}^n \frac{\partial^{\vert \vec k_{T}\vert}}{\partial T_2^{\vec k_{T}}}\frac{\partial^{\vert \vec k_{\theta}\vert}}{\partial \theta_2^{\vec k_{\theta}}}
(\text{\rm Res}(\frak v_{(\frak y_2,\vec T_2,\vec \theta_2)}) -\text{\rm Res}(\frak v_{(\frak y_2,\vec \infty)})
\right\Vert_{C^k}
< C_{3,k}e^{-\delta_2  T_{2,{\rm min}}},
\endaligned
\end{equation}
if $T_{2,{\rm e}} \ge T_{2,{\rm min}} > T_k$ for all $\rm e$  and $\vert{\vec k_{T}}\vert + \vert{\vec k_{\theta}}\vert +n \le k$.
Here $T_{2,{\rm min}}=\min (T_{2,{\rm e}} ; {\rm e} \in C^1(\mathcal G_{\frak y_2}))$.
\end{cor}
\par\medskip
In later subsections we also use a parametrized version of Propositions \ref{changeinfcoorprop} and \ref{reparaexpest},
which we discuss now.
\par
Let $Q$ be a finite dimensional manifold. Suppose we have a fiber bundle
\begin{equation}\label{fibrationsigmafami}
\pi : \tilde{\frak M}^{(2)}_{\frak x_{\rm v}} \to
Q_{\rm v} \times  \frak V(\frak x_{\rm v})
\end{equation}
that is a universal family (\ref{fibrationsigma}) when we restrict it to each of
$\{\xi\} \times \frak V(\frak x_{\rm v})$ for $\xi_{\rm v} \in Q_{\rm v}$.
We put
$$
Q = \prod_{{\rm v} \in C^0(\mathcal G)} Q_{\rm v}.
$$
\begin{defn}\label{def:Qfamily}
A $Q$-{\it parametrized family of coordinates at infinity} is a fiber bundle
(\ref{fibrationsigmafami}) and its trivialization so that
for each $\xi = (\xi_{\rm v})$ the restriction to $\{\xi_{\rm v}\} \times\frak V(\frak x_{\rm v}) $
gives a coordinate at infinity in the sense of Definition \ref{coordinatainfdef}.
\end{defn}
Suppose a $Q$-parametrized family of coordinate at infinity in the above sense is given.
Then we can perform the construction we already described
for each $\xi$ and  obtain a map
\begin{equation}
\overline{\Phi}_{2} :
Q \times \prod_{{\rm v}\in C^0(\mathcal G)}\frak V(\frak x_{\rm v}) \times
(\vec T^{\rm o}_0,\infty] \times ((\vec T^{\rm c}_0,\infty] \times \vec S^1)
\to \mathcal M_{k+1,\ell}.
\end{equation}
Note that for each $\xi \in Q$ it gives a diffeomorphism to a neighborhood of
$\frak x$ in $\mathcal M_{k+1,\ell}$.
\par
Suppose we have a (unparametrized) coordinate at infinity
that is a fiber bundle
$$
\pi : {\frak M}^{(1)}_{\frak x_{\rm v}} \to \frak V(\frak x_{\rm v})
$$
equipped with trivialization.
It induces an embedding
$$
\overline{\Phi}_{1} :
\prod_{{\rm v}\in C^0(\mathcal G)}\frak V(\frak x_{\rm v}) \times
(\vec T^{\rm o}_0,\infty] \times ((\vec T^{\rm c}_0,\infty] \times \vec S^1)
\to \mathcal M_{k+1,\ell}.
$$
They induce a map
\begin{equation}\label{coodinatechange12}
\aligned
\overline{\Phi}_{12}:
&Q \times \prod_{{\rm v}\in C^0(\mathcal G)}\frak V(\frak x_{\rm v}) \times
(\vec T^{\rm o \prime}_0,\infty] \times ((\vec T^{\rm c \prime}_0,\infty] \times \vec S^1)\\
&\to
\prod_{{\rm v}\in C^0(\mathcal G)}\frak V(\frak x_{\rm v}) \times
(\vec T^{\rm o}_0,\infty] \times ((\vec T^{\rm c}_0,\infty] \times \vec S^1)
\endaligned
\end{equation}
by the formula:
$$
\overline{\Phi}_{1}(\overline{\Phi}_{12}(\xi,\frak y,\vec T,\vec\theta))
= \overline{\Phi}_{2}(\xi,\frak y,\vec T,\vec\theta).
$$
Here $\vec T^{\rm o \prime}_0$ and $\vec T^{\rm c \prime}_0$
are sufficiently large compared with $\vec T^{\rm o}_0$ and $\vec T^{\rm c }_0$.
\par
Moreover we have a family of biholomorphic maps:
\begin{equation}\label{mapvbra}
\frak v_{(\xi,\frak y,\vec T,\vec\theta)} : \Sigma_{\vec T,\vec\theta}^{\frak y,\xi,(2)}
\to \Sigma_{\vec T',\vec \theta'}^{\frak y',(1)}.
\end{equation}
Here $(\frak y',\vec T',\vec\theta') = \overline{\Phi}_{12}(\xi,\frak y,\vec T,\vec \theta)$
and $\Sigma_{\vec T',\vec \theta'}^{\frak y',(1)}$, $\Sigma_{\vec T,\vec\theta}^{\frak y,\xi,(2)}$
 are marked bordered curves
representing $\overline{\Phi}_{1}(\overline{\Phi}_{12}(\xi,\rho,\vec T,\vec\theta))$ and
$\overline{\Phi}_{2}(\xi,\rho,\vec T,\vec\theta)$, respectively.
\begin{lem}\label{changeinfcoorproppara}
We  have $C_{4,k}$, $C_{5,k}$ such that:
\begin{equation}\label{2336}
\left\Vert
\nabla_{\xi}^{\vec k_{\xi}}\frac{\partial^{\vert \vec k_{T}\vert}}{\partial T^{\vec k_{T}}} \frac{\partial^{\vert \vec k_{\theta}\vert}}{\partial \theta^{\vec k_{\theta}}}
({\overline{\Phi}}_{12}(\xi,\frak y,\vec T,\vec \theta)
- \Psi_{12}(\frak y,\vec T,\vec \theta))
\right\Vert_{C^k} \le C_{4,k} e^{-\delta (\vec k_{T}\cdot \vec T+\vec k_{\theta}\cdot \vec T^{\rm c})}
\end{equation}
for $\vert \vec k_{\xi}\vert$, $\vert\vec k_T\vert, \vert \vec k_{\theta}\vert \le k$,
if each of $T_{{\rm e}}$ is greater than $T_k$.
The left hand sides are $C^k$ norm (as functions on $\frak y$).
Moreover for each ${\rm e}_0 \in C^1_{\rm c}(\mathcal G_{\frak y_2})$ we have
\begin{equation}\label{2337}
\aligned
&\left\Vert
\nabla_{\xi}^{\vec k_{\xi}}\nabla_{\frak y}^n \frac{\partial^{\vert \vec k_{T}\vert}}{\partial T^{\vec k_{T}}}\frac{\partial^{\vert \vec k_{\theta}\vert}}{\partial \theta^{\vec k_{\theta}}}
\frac{\partial}{\partial T_{{\rm e}_0}}
\text{\rm Res}(\frak v_{(\xi,\frak y,\vec T,\vec \theta)})
\right\Vert_{C^k}
< C_{5,k}e^{-\delta_2  T_{{\rm e}_0}},
\\
&\left\Vert
\nabla_{\xi}^{\vec k_{\xi}}\nabla_{\frak y}^n \frac{\partial^{\vert \vec k_{T}\vert}}{\partial T^{\vec k_{T}}}\frac{\partial^{\vert \vec k_{\theta}\vert}}{\partial \theta^{\vec k_{\theta}}}
\frac{\partial}{\partial \theta_{{\rm e}_0}}
\text{\rm Res}(\frak v_{(\xi,\frak y,\vec T,\vec \theta)})
\right\Vert_{C^k}
< C_{5,k}e^{-\delta_2  T_{{\rm e}_0}},
\endaligned
\end{equation}
if each of $T_{{\rm e}}$ is greater than $T_k$ and
$\vert \vec k_{\xi}\vert + \vert{\vec k_{T}}\vert + \vert{\vec k_{\theta}}\vert +n \le k$.
\par
The first inequality of (\ref{2337}) also holds for ${\rm e}_0 \in C^1_{\rm o}(\mathcal G_{\frak y_2})$.
\end{lem}
Note that (\ref{2336}), (\ref{2337}) are parametrized versions of
Propositions \ref{changeinfcoorprop}, \ref{reparaexpest}, respectively.
For the proof, see Section \ref{proposss}.

\par\medskip
\section{Stabilization of the source by adding marked points and obstruction bundles}
\label{stabilization}

Let $((\Sigma,\vec z,\vec z^{\text{\rm int}}),u) = (\frak x,u) \in \mathcal M_{k+1,\ell}(\beta;\mathcal G)$.
We assume that $\mathcal G$ is stable but is not source stable.
In Section \ref{secsimple} we assumed that the source is stable.
In order to carry out  analytic detail similar to the one in Section \ref{secsimple} in the general case,
we stabilize the source by adding marked points.
In other words, we use the method of \cite[appendix]{FOn} for this purpose.\footnote{
16 years of experience shows that the method of \cite[appendix]{FOn}
is easier to use in various applications than the method of \cite[Section 13]{FOn}.}
\begin{rem}\label{rem161}
We note that the method of \cite[appendix]{FOn}
had been  used earlier   in various places by many people.
A nonexhausting list of it is \cite[Proposition 7.11, Theorem 9.1]{Wo},
\cite[appendix]{FOn}, \cite[begining of Section 3 and the proof of Lemma 3.1]{LiTi98},
\cite[page 395]{Si}, \cite[page 424]{fooo:book1}, \cite[Section 4.3]{fooo:toricmir}.
See also \cite[(3.9)]{Rua99}.
\end{rem}
We recall:
\begin{defn}
An irreducible component $\frak x_{\rm v} = (\Sigma_{\rm v},\vec z_{\rm v},\vec z_{\rm v}^{\text{\rm int}})$ of $\frak x$ is
said to be {\it unstable},
if and only if one of the following holds:
\begin{enumerate}
\item
$\frak x_{\rm v} \in \mathcal M_{k_{\rm v}+1,\ell_{\rm v}}$ and
$k_{\rm v} + 1+ 2\ell_{\rm v} < 3$.
\item
$\frak x_{\rm v} \in \mathcal M^{\rm cl}_{\ell_{\rm v}}$ and
$\ell_{\rm v} < 3$.
\end{enumerate}
\end{defn}
There is at least one boundary marked point in case $\frak x_{\rm v}$ is a disk
($\frak x \in \mathcal M_{k+1,\ell}$ and $k+1 > 0$), and at least one interior marked point in case $\frak x_{\rm v}$ is a sphere.
(This is because it should be attached to a disk or to a sphere.) Note we assume $\ell \ge 1$ in case of $\mathcal M_{\ell}^{\rm cl}$.)
Therefore there are three cases where $\frak x_{\rm v}$ is unstable:
\begin{enumerate}
\item[(a)] $\frak x_{\rm v}$ is a disk. $\frak x_{\rm v} \in \mathcal M_{k_{\rm v}+1,\ell_{\rm v}}$
and $k_{\rm v} = 0$ or $1$. $\ell_{\rm v} = 0$.
\item[(b)]
$\frak x_{\rm v}$ is a sphere.
$\frak x_{\rm v} \in \mathcal M^{\rm cl}_{\ell_{\rm v}}$ and
$\ell_{\rm v} =2$.
\item [(c)]
$\frak x_{\rm v}$ is a sphere.
$\frak x_{\rm v} \in \mathcal M^{\rm cl}_{\ell_{\rm v}}$ and
$\ell_{\rm v} =1$.
\end{enumerate}
\begin{rem}
In the case of higher genus there are some other kinds of irreducible components that are unstable.
For example, $T^2$ without marked points is unstable. We can handle them in the same way.
If we consider also $\mathcal M_{0}^{\rm cl}(\alpha)$, then $\mathcal M_{0}^{\rm cl}$ also appears.
\end{rem}
\begin{defn}(\cite[Section 13 p989 and appendix p1047]{FOn})
A {\it minimal stabilization} is a choice of additional interior marked points, where we put one interior
marked point $w_{\rm v}$ of $\Sigma_{\rm v}$
for each $\frak x_{\rm v}$ satisfying (a) or (b) above and two interior marked points
$w_{{\rm v},1}$, $w_{{\rm v},2}$ for each $\frak x_{\rm v}$ satisfying (c) above,
so that the following holds.
\begin{enumerate}
\item
$w_{\rm v} \notin \vec z_{\rm v}^{\text{\rm int}}$. $w_{{\rm v},1}, w_{{\rm v},2} \notin \vec z_{\rm v}^{\text{\rm int}}$.
They are not singular.
\item
$u$ is an immersion at $w_{\rm v}$, $w_{{\rm v},1}$, $w_{{\rm v},2}$.
\item
Let $v \in \Gamma_{(\frak x,u)}^+$ such that $v\Sigma_{\rm v} = \Sigma_{\rm v'}$. Suppose $\frak x_{\rm v}$ satisfies (a) or (b) above.
Then $v w_{\rm v} = v' w_{\rm v'}$ for some $v' \in   \Gamma_{(\frak x_{\rm v'},u)}^+$.
Suppose $\frak x_{\rm v}$ satisfies (c) above. Then there exists $v' \in   \Gamma_{(\frak x_{\rm v'},u)}^+$
such that  $v w_{{\rm v},i} = v' w_{{\rm v}',i}$ for $i=1,2$.
\item
$w_{{\rm v},1}\ne v'w_{{\rm v},2}$ for any $v' \in   \Gamma_{(\frak x_{\rm v},u)}^+$.
\end{enumerate}
\end{defn}
(We add three marked points in the case of $\mathcal M_{0}^{\rm cl}$.)
\begin{defn}\label{symstabili}
A {\it symmetric stabilization} is a choice of additional marked points
$\vec w = (w_1,\dots,w_{\ell'}) \in \text{\rm Int}\,\Sigma$, such that:
\begin{enumerate}
\item
$\vec w \cap \vec z^{\text{\rm int}} = \emptyset$.
\item
$w_i \ne w_j$ for $i\ne j$.
\item
$u$ is an immersion at each $w_i$.
\item
$(\Sigma, \vec z, \vec w \cup \vec z^{\text{\rm int}})$ is stable.
\item
For each $v \in \Gamma_{(\frak x,u)}^+$ there exists $\sigma_{v} \in \frak S_{\ell'}$, such that
$$
v(w_i) = w_{\sigma_{v}(i)}.
$$
\end{enumerate}
\end{defn}
We note that a minimal stabilization induces a symmetric stabilization.
Namely we take
$$
\aligned
&\{v w_{\rm v} \mid v \in \Gamma_{(\frak x_{\rm v},u)}^+, \,\,
\text{$\frak x_{\rm v}$ satisfies (a) or (b)}\}\\
&\cup
\{v w_{{\rm v},i} \mid v \in \Gamma_{(\frak x_{\rm v},u)}^+, i=1,2, \,\,
\text{$\frak x_{\rm v}$ satisfies (c)}\}.
\endaligned
$$
Since the notion of symmetric stabilization is  more general, we use
symmetric stabilization in this note.
Symmetric stabilization was used in \cite{fooo:toricmir}.
\par
We write
$$
\frak x\cup \vec w = (\Sigma, \vec z, \vec z^{\text{\rm int}}\cup \vec w )
$$
when $\frak x = (\Sigma, \vec z, \vec z^{\text{\rm int}})$.
\begin{rem}
In our genus zero case, Definition \ref{symstabili} (4) implies that the automorphism
group of $(\Sigma, \vec z, \vec z^{\text{\rm int}}\cup \vec w)$ is trivial.\footnote{In the case of higher genus, we may include the triviality of the automorphism
as a part of the definition of the symmetric stabilization. If we do so then (\ref{sigmahomotoS}) is still an injective homomorphism.}
So we can define an {\it injective} homomorphism
\begin{equation}\label{sigmahomotoS}
\sigma : \Gamma_{(\frak x,u)} \to \frak S_{\ell'}
\end{equation}
by
$$
v(w_i) = w_{\sigma(i)}.
$$
(Here $\frak S_{\ell'}$ is the symmetric group of order $\ell'!$.)
We denote by $\frak H_{(\frak x,u)}$ the image of (\ref{sigmahomotoS}).
In a similar way we obtain an injective homomorphism
\begin{equation}
\sigma : \Gamma^+_{(\frak x,u)} \to \frak S_{\ell}\times \frak S_{\ell'}.
\end{equation}
We denote its image by $\frak H^+_{(\frak x,u)}$.
\end{rem}
\par
We use the notion of symmetric stabilization of $\frak x \in \mathcal M_{k+1,\ell}(\beta;\mathcal G)$ to define the notion of
obstruction bundle data as follows.

\begin{defn}\label{obbundeldata}
An {\it obstruction bundle data $\frak E_{\frak p}$ centered at}
$$
\frak p = (\frak x,u) = ((\Sigma,\vec z,\vec z^{\text{\rm int}}),u) \in \mathcal M_{k+1,\ell}(\beta;\mathcal G)
$$
is the data satisfying the conditions described below.
\begin{enumerate}
\item
A symmetric stabilization $\vec w = (w_1,\dots,w_{\ell'})$ of $(\frak x,u)$.
We denote by $\mathcal G_{\vec w \cup \frak x}$ the combinatorial type of $\vec w \cup \frak x$.
\item
A neighborhood $\frak V(\frak x_{\rm v} \cup \vec w_{\rm v})$ of $
\frak x_{\rm v} \cup \vec w_{\rm v} = (\Sigma_{\frak x_{\rm v}},\vec z_{\frak x_{\rm v},}\vec z^{\rm int}_{\rm v} \cup \vec w_{\rm v})$
in $\overset{\circ}{\mathcal M}
_{k_{\rm v}+1,\ell_{\rm v}+\ell_{\rm v}'}$ or $\overset{\circ}{\mathcal M}^{\rm cl}
_{\ell_{\rm v}+\ell_{\rm v}'}$.
Here $\frak x_{\rm v} \in \overset{\circ}{\mathcal M}
_{k_{\rm v}+1,\ell_{\rm v}}$ or $\in \overset{\circ}{\mathcal M}^{\rm cl}
_{\ell_{\rm v}+\ell_{\rm v}'}$ is an irreducible component of $\frak x$ and $\vec w_{\rm v}$ is a part of
$\vec w$ that is contained in this irreducible component.
\item
A universal family with coordinate at infinity of  $\frak x_{\rm v} \cup \vec w_{\rm v}$ defined on
$\frak V(\frak x_{\rm v} \cup \vec w_{\rm v})$. (We use the notation of
Definition \ref{coordinatainfdef}.)
We assume that it is invariant under the $\Gamma_{(\frak x\cup \vec w,u)}^{\frak H^+_{(\frak x,u)}}$ action in the sense we will explain later.
\item
A compact subset $K^{\rm obst}_{\rm v}$ such that $K^{\rm obst}_{\rm v} \times \frak V(\frak x_{\rm v}\cup \vec w_{\rm v})$  is contained in
$\frak K_{\frak x_{\rm v}}$, which is defined in
Definition \ref{coordinatainfdef} (3). We assume that they are  $\Gamma_{(\frak x\cup \vec w,u)}^{\frak H^+_{(\frak x,u)}}$ invariant
in the sense we will explain later.
We call $K^{\rm obst}_{\rm v}$ the {\it support of the obstruction bundle}.
\item
A $\frak y \in \frak V(\frak x_{\rm v} \cup \vec w_{\rm v})$-parametrized smooth family
of finite dimensional complex linear subspaces $E_{\frak p,{\rm v}}(\frak y,u)$
of
$$
\Gamma_0(\text{\rm Int}\,K^{\rm obst}_{\rm v}; u^*TX \otimes \Lambda^{01}).
$$
Here $\Gamma_0$ denotes the set of the smooth sections with compact support on the domain $\Sigma_{\frak y_{\rm v}}$ induced   by
$\frak y_{\rm v} \in \frak V(\frak x_{\rm v} \cup \vec w_{\rm v})$.
We regard $u : \Sigma_{\frak x_{\rm v}} \to X$ also as a map from $\Sigma_{\frak y_{\rm v}}$
by using the smooth trivialization of the universal family given as a part of
Definition \ref{coordinatainfdef} (5).
\par
We assume that $\bigoplus_{{\rm v} \in C^0(\mathcal G)} E_{{\rm v}}$ is invariant under the $\Gamma_{(\frak x\cup \vec w,u)}^{\frak H^+_{\frak p}}$ action
in the sense we will explain later.
\item
For each $\text{\rm v} \in  C^0_{\rm d}(\mathcal G_{\frak p})$ and $\frak y_{\rm v} \in \frak V(\frak x_{\rm v}\cup \vec w_{\rm v})$
the differential operator
\begin{equation}\label{ducomponents0}
\aligned
\overline D_{u} \overline\partial : &L^2_{m+1,\delta}((\Sigma_{\frak y_{\rm v}},\partial \Sigma_{\frak y_{\rm v}});
u^*TX, u^*TL) \\
&\to
L^2_{m,\delta}(\Sigma_{\frak y_{\rm v}}; u^*TX \otimes \Lambda^{01})/E_{\frak p,{\rm v}}(\frak y,u)
\endaligned
\end{equation}
is surjective.
(We define the above weighted Sobolev spaces in the same way as in Section \ref{subsec12}.
See Section \ref{glueing} for the precise definition in the general case.)
\par
If $\text{\rm v} \in  C^0_{\rm s}(\mathcal G_{\frak p})$ and $\frak y_{\rm v} \in \frak V(\frak x_{\rm v}\cup \vec w_{\rm v})$,
the differential operator
\begin{equation}\label{ducomponentssphere}
\aligned
\overline D_{u} \overline\partial : &L^2_{m+1,\delta}(\Sigma_{\frak y_{\rm v}};
u^*TX) \\
&\to
L^2_{m,\delta}(\Sigma_{\frak y_{\rm v}}; u^*TX \otimes \Lambda^{01})/E_{\frak p,{\rm v}}(\frak y,u)
\endaligned
\end{equation}
is surjective.
\item
The kernels of (\ref{ducomponents0}) and (\ref{ducomponentssphere}) satisfy a transversality property for evaluation maps
that is as described in Condition \ref{transeval}.
\item
For each $w_i \in \Sigma_{\rm v}$ we take a codimension $2$ submanifold $\mathcal D_{i}$ of $X$
such that $u(w_i) \in \mathcal D_i$ and
$$
u_*T_{w_i}\Sigma_{\rm v} + T_{u(w_i)}\mathcal D_i = T_{w_i}X.
$$
Moreover $\{\mathcal D_i\}$ is invariant under the $\Gamma_{\frak p}^+$ action in the following sense.
Let $v \in \Gamma_{\frak p}^+$ and $v(w_i) =w_{\sigma(i)}$ then
\begin{equation}\label{Dequivalent}
\mathcal D_i = \mathcal D_{\sigma(i)}.
\end{equation}
(Note $u(w_i) = u(w_{\sigma(i)})$ since $u\circ v = u$.)
\end{enumerate}
\end{defn}
\begin{conds}\label{transeval}
Suppose a vertex $\text{\rm v} \in C^0_{\rm d}(\mathcal G_{\frak p})$ is contained in an edge ${\rm e} \in C^1_{\rm o}(\mathcal G_{\frak p})$.
Let $z_{\rm e}$ be a singular point of $\Sigma_{\frak x}$ corresponding to the edge ${\rm e} \in C^1_{\rm o}(\mathcal G_{\frak p})$.
We define
\begin{equation}\label{evaluflagcomp}
\text{\rm ev}_{\rm v,e} :
L^2_{m+1,\delta}((\Sigma_{\frak y_{\rm v}},\partial \Sigma_{\frak y_{\rm v}});
u^*TX, u^*TL)
\to T_{u(z_{\rm e})}L
\end{equation}
by
$
s \mapsto \pm s(z_{\rm e})
$
where we take $+$ if $\text{\rm v}$ is an outgoing vertex of $\rm e$ and we take $-$
if $\text{\rm v}$ is an incoming vertex of $\rm e$.
If $\text{\rm v} \in C^0_{\rm d}(\mathcal G_{\frak p})$ and ${\rm e} \in C^1_{\rm c}(\mathcal G_{\frak p})$, then we define
\begin{equation}\label{evaluflagcomp20}
\text{\rm ev}_{\rm v,e} :
L^2_{m+1,\delta}((\Sigma_{\frak y_{\rm v}},\partial\Sigma_{\frak y_{\rm v}});
u^*TX,u^*TL)
\to T_{u(z_{\rm e})}X
\end{equation}
by the same formula.
In a similar way we define
\begin{equation}\label{evaluflagcompcl}
\text{\rm ev}_{\rm v,e} :
L^2_{m+1,\delta}(\Sigma_{\frak y_{\rm v}};
u^*TX) \to T_{u(z_{\rm e})}X,
\end{equation}
if ${\rm e} \in C^1_{\rm c}(\mathcal G_{\frak p})$ and $\text{\rm v} \in C^0_{\rm s}(\mathcal G_{\frak p})$
is its vertex.
\par
Combining all of (\ref{evaluflagcomp}), (\ref{evaluflagcomp20}), (\ref{evaluflagcompcl}) we obtain a map:
\begin{equation}
\aligned
\text{\rm ev}_{\mathcal G_{\frak p}} :
&\bigoplus_{\text{\rm v} \in C^0_{\rm d}(\mathcal G_{\frak p})}
L^2_{m+1,\delta}((\Sigma_{\frak y_{\rm v}},\partial \Sigma_{\frak y_{\rm v}});
u^*TX, u^*TL)
\\
&\oplus \bigoplus_{\text{\rm v} \in C^0_{\rm s}(\mathcal G_{\frak p})}
L^2_{m+1,\delta}(\Sigma_{\frak y_{\rm v}};
u^*TX)
\\
&\to
\bigoplus_{{\rm e} \in C^1_{\rm o}(\mathcal G_{\frak p})} T_{u(z_{\rm e})}L
\oplus
\bigoplus_{{\rm e} \in C^1_{\rm c}(\mathcal G_{\frak p})} T_{u(z_{\rm e})}X.
\endaligned
\end{equation}
The condition we require is
that the restriction of $\text{\rm ev}_{\mathcal G_{\frak p}} $ to
$$
\bigoplus_{\text{\rm v} \in C^0(\mathcal G_{\frak p})}  \text{\rm Ker}\overline D_{u_{\rm v}} \overline\partial
$$
is surjective.
\end{conds}
\begin{rem}
In \cite{fooo:book1} we used Kuranishi structures on $\mathcal M_{k+1,\ell}(\beta)$ so that
the evaluation maps $\text{\rm ev} : \mathcal M_{k+1,\ell}(\beta) \to L^{k+1} \times X^{\ell}$
are weakly submersive.
To construct Kuranishi structures satisfying this additional property,
we need to require an additional assumption to the obstruction bundle data. Namely we need to assume that the evaluation
maps at the marked points
$$\aligned
{\rm ev} :
&\bigoplus_{\text{\rm v} \in C^0_{\rm d}(\mathcal G_{\frak p})}
L^2_{m+1,\delta}((\Sigma_{\frak y_{\rm v}},\partial \Sigma_{\frak y_{\rm v}});
u^*TX, u^*TL)
\to
\prod_{i=0}^{k} T_{u(z_i)}L \times \prod_{i=1}^{\ell} T_{u(z_i^{\rm int})}X
\endaligned$$
are also surjective. But
we do not include it in the definition here since there are cases we do not assume it.
\end{rem}
We next explain the precise meaning of invariance under the action in (3), (4), (5).
The invariance in (3) is defined in Definition \ref{defn288}.
The $\Gamma_{(\frak x\cup \vec w,u)}^{\frak H^+_{(\frak x,u)}}$ action
on
$\frak K_{\frak x_{\rm v}}$ is induced by its action.
(See  Definition \ref{defn288}.)
So we require (the totality of) $K^{\rm obst}_{\rm v}$ is invariant under this action in (4).
To make sense of (5) we define a $\Gamma_{(\frak x\cup \vec w,u)}^{\frak H^+_{(\frak x,u)}}$ action on
\begin{equation}\label{sumobstructionspace}
\bigoplus_{{\rm v} \in C^0(\mathcal G)}\Gamma_0(\text{\rm Int}\,K^{\rm obst}_{\rm v}; u^*TX \otimes \Lambda^{01}).
\end{equation}
If $v \in \Gamma_{(\frak x\cup \vec w,u)}^{\frak H^+_{(\frak x,u)}}$ then
$v\Sigma_{\rm v} = \Sigma_{\rm v'}$ for some ${\rm v}'$ and
$K^{\rm obst}_{{\rm v}'} = vK^{\rm obst}_{{\rm v}}$ by (4). Moreover $u\circ v = u$ holds on $\Sigma_{\rm v}$.
Therefore we obtain
$$
v_* :
\Gamma_0(\text{\rm Int}\,K^{\rm obst}_{\rm v}; u^*TX \otimes \Lambda^{01})
\cong
\Gamma_0(\text{\rm Int}\,K^{\rm obst}_{{\rm v}'}; u^*TX \otimes \Lambda^{01}).
$$
They induce a $\Gamma_{(\frak x\cup \vec w,u)}^{\frak H^+_{(\frak x,u)}}$ action on (\ref{sumobstructionspace}).
Note that this is the case of the action at $\vec w \cup \frak p = (\vec w \cup \frak x,u)$.
When we move to a nearby point $(\frak y,u)$, the situation becomes slightly different, since $v_*\frak y = \frak y$
holds no longer.
We have a smooth trivialization of the bundle (\ref{fibrationsigma}). (Definition \ref{coordinatainfdef} (5).)
Namely we are given a diffeomorphism
$$
v : K_{\rm v}(\frak y) \to K_{\rm v'}(\frak y)
$$
between the cores.
(Here we write $K_{\rm v}(\frak y)$ in place of $K_{\rm v}$ to include its complex structure.)
However this is not a biholomorphic map. On the other hand
$$
v : K_{\rm v}(\frak y) \to K_{\rm v'}(v_*\frak y)
$$
is a biholomorphic map by Definition \ref{coordinatainfdef} (1). Therefore we still obtain  a map
\begin{equation}\label{sentbyv}
\aligned
v_* :
&\Gamma_0(\text{\rm Int}\,K^{\rm obst}_{\rm v}(\frak y); u^*TX \otimes \Lambda^{01})\\
&\cong
\Gamma_0(\text{\rm Int}\,K^{\rm obst}_{{\rm v}'}(v_*\frak y); (u\circ v^{-1})^*TX \otimes \Lambda^{01}).
\endaligned
\end{equation}
Definition \ref{obbundeldata} (5) means
$$
v_*\left(E_{\frak p,{\rm v}}(\frak y,u)\right) = E_{\frak p,{\rm v}'}
(v_*\frak y,u\circ v^{-1})
=
E_{\frak p,{\rm v}'}
(v_*\frak y,u)
$$
where the map $v_*$ appearing at the beginning of the formula is the map (\ref{sentbyv}).
\begin{rem}\label{rem23}
The condition (8), especially $u(w_i) \in \mathcal D_i$, is assumed only for $\frak p$ and $\vec w$. For the
general point $\frak V(\frak y_{\rm v} \cup \vec w_{\rm v})$ this condition is not assumed at this stage.
We put this condition only at later step (Section \ref{cutting}. See also Definition \ref{transconst}.) and only to the solutions of the equation.
\end{rem}
\begin{lem}\label{existobbundledata}
For each $\frak p$ there exists an obstruction bundle data $\frak E_{\frak p}$ centered at $\frak p$.
\end{lem}
\begin{proof}
Existence of symmetric stabilization is obvious.
We can find $E_{\frak p,\text{\rm v}}(\frak p \cup \vec w_{\frak p})$ for
${\rm v} \in C^0(\mathcal G_{\frak p \cup \vec w_{\frak p}})$  satisfying (7), (8) by the unique continuation
properties of the linearization of the Cauchy-Riemann equation.
We can make them  $\Gamma_{\frak p \cup \vec w_{\frak p}}^{\frak H^+_{\frak p}}$ invariant
by taking the union of the images of actions.
Then we extend them to a small neighborhood of $\frak p \cup \vec w_{\frak p}$ in a way such that (7), (8) are satisfied.
We  make them $\Gamma_{\frak p \cup \vec w_{\frak p}}^{\frak H^+_{\frak p}}$ invariant by taking average as follows.
Let $\frak y = (\frak y_{\rm v})$ such that $\frak y_{\rm v} \in \frak V(\frak x_{\rm v} \cup \vec w_{\rm v})$.
Using the trivialization of the bundle  (\ref{fibrationsigma})
we can define
$$
\frak I'_{\frak y} : \bigoplus_{{\rm v} \in C^0(\mathcal G_{\frak p \cup \vec w_{\frak p}})} E_{\frak p,{\rm v}}
\to \bigoplus_{{\rm v} \in C^0(\mathcal G)} \Gamma(\Sigma_{\frak y,{\rm v}}; u^*TX\otimes \Lambda^{01}).
$$
Note for $v \in \Gamma_{\frak p \cup \vec w_{\frak p}}^{\frak H^+_{\frak p}}$  the equality
$
v_* \circ \frak I'_{\frak y} = \frak I'_{v\frak y}\circ v_*
$
may not be satisfied. However since
$
v_* \circ \frak I'_{\frak p} = \frak I'_{\frak p}\circ v_*
$
we  may assume
$$
\Vert v_* \circ \frak I'_{\frak y} - \frak I'_{v\frak y}\circ v_* \Vert
$$
is small by taking $\frak V(\frak x_{\rm v} \cup \vec w_{\rm v})$ small.
Therefore
$$
\frak I_{\frak y} =
\frac{1}{\# \Gamma_{\frak p \cup \vec w_{\frak p}}}\sum_{v \in  \Gamma_{\frak p \cup \vec w_{\frak p}}^{\frak H^+_{\frak p}}}
(v^{-1})_* \circ  \frak I'_{v\frak y} \circ v_*
$$
is injective and close to $\frak I'_{\frak p}$.
We hence obtain the required  $E_{\frak p}(\frak y)$ by
$$
E_{\frak p}(\frak y) = \text{\rm Im}\frak I_{\frak y} .
$$
The existence of the codimension $2$ submanifolds $\mathcal D_i$ is obvious.
\end{proof}
The obstruction bundle data determines
$$
E_{\frak p}(\frak y,u) = \bigoplus_{{\rm v} \in C^0(\mathcal G_{\frak p \cup \vec w_{\frak p}})}E_{\frak p,{\rm v}}(\frak y,u)
\subset
\bigoplus_{{\rm v} \in C^0(\mathcal G_{\frak p \cup \vec w_{\frak p}})}L^2_{m,\delta}(\Sigma_{\frak y_{\rm v}}; u^*TX \otimes \Lambda^{01})
$$
for $\frak y \in \frak V(\frak x\cup \vec w)$.
This subspace plays the role of (a part of) the obstruction bundle of the Kuranishi structure we will construct.
To define our equation and thickened moduli space we need to extend the family of linear subspaces
$E_{\frak p}(\cdot)$
so that we  associate $E_{\frak p}(\frak q)$ to  an object $\frak q$ which is `close' to $\frak p$.
We will define this close-ness below.
(This is a generalization of Condition \ref{nearbyuprime}.)
\par
We use the map
$$
\overline{\Phi} : \prod_{\rm v \in C^{0}(\mathcal G_{\frak p \cup \vec w_{\frak p}})}\frak V(\frak x_{\rm v} \cup \vec w_{\rm v}) \times (\vec T^{\rm o}_0,\infty] \times ((\vec T^{\rm c}_0,\infty] \times \vec S^1)
\to \mathcal M_{k+1,\ell+\ell'}.
$$
(See Definition \ref{def29}.)
Let $\frak Y = \overline{\Phi}(\frak y,\vec T,\vec \theta)$ be an element of $\mathcal M_{k+1,\ell+\ell'}$ that is represented by
$(\Sigma_{\frak Y},\vec z_{\frak Y}, \vec z^{\text{\rm int}}_{\frak Y} \cup \vec w_{\frak Y})$.
By construction (\ref{summandtoglued}) we have
$$
\aligned
\Sigma_{\frak Y} = \bigcup_{{\rm v} \in C^0(\mathcal G_{\frak p \cup \vec w_{\frak p}})} K^{\frak Y}_{{\rm v} }
&\cup \bigcup_{{\rm e}\in C^1_{\mathrm o}(\mathcal G_{\frak p \cup \vec w_{\frak p}})} [-5T_{\rm e},5T_{\rm e}] \times [0,1]
\\
&\cup
\bigcup_{{\rm e}\in C^1_{\mathrm c}(\mathcal G_{\frak p \cup \vec w_{\frak p}})} [-5T_{\rm e},5T_{\rm e}] \times S^1.
\endaligned
$$
We called the second and the third summand the neck region.
In case $T_{\rm e} = \infty$ the product of the union of two half lines and $[0,1]$ or $S^1$ is
also called the neck region.
See  Definition \ref{defcoreandneck}.
\begin{defn}\label{epsiloncloseto}
Let $u' : (\Sigma_{\frak Y} ,\partial\Sigma_{\frak Y}) \to (X,L)$ be a smooth map in homology class $\beta$.
We say that $(\Sigma_{\frak Y},u')$ is {\it $\epsilon$-close to $\frak p$}
with respect to the given obstruction bundle data if the following holds.
\begin{enumerate}
\item
Since $\frak Y = \overline{\Phi}(\frak y,\vec T,\vec \theta)$
the core $K^{\frak Y}_{{\rm v} } \subset \Sigma_{\frak Y}$ is identified with $K^{\frak y}_{\rm v} \subset \Sigma_{\frak y}$.
We require
\begin{equation}\label{c1close}
\vert u - u'\vert_{C^{10}(K^{\frak Y}_{{\rm v} })} < \epsilon
\end{equation}
for each $\rm v$. (We regard $u$ as a map from $\Sigma_{\frak y}$ by using the
smooth trivialization of the universal family given as a part of
Definition \ref{coordinatainfdef} (4).)
\item
The map $u'$ is holomorphic on each of the neck region.
\item
The diameter of the $u'$ image of each of the connected component of the neck region is
smaller than $\epsilon$.
\item $T_{\rm e} > \epsilon^{-1}$ for each $\rm e$.
\end{enumerate}
\end{defn}
\begin{rem}
We use metrics of the source and of $X$ to define the left hand side of (\ref{c1close}).
See Remark \ref{metricfiber}.
\end{rem}
\begin{rem}
We note that Definition \ref{epsiloncloseto} is not a definition of
topology on certain set.
In fact, `$(\Sigma_{\frak Y},u')$ is close to $\frak p$' is defined only when $\frak p$ is an element of $\mathcal M_{k+1,\ell}(\beta)$,
but $(\Sigma_{\frak Y},u')$ may not be an element of $\mathcal M_{k+1,\ell}(\beta)$.
\par
Even in case $(\Sigma_{\frak Y},u') \in \mathcal M_{k+1,\ell}(\beta)$,
the fact that $(\Sigma_{\frak Y},u')$ is $\epsilon$-close to $\frak p$ does not
imply that
$\frak p$ is $\epsilon$-close $(\Sigma_{\frak Y},u')$.
In fact, if  $(\Sigma_{\frak Y},u')$ is $\epsilon$-close to $\frak p$  then
$\mathcal G_{\frak p} \succ \mathcal G_{\frak Y}$.
\par
On the other hand, we have the following.
If $(\Sigma_{\frak Y},u') \in \mathcal M_{k+1,\ell}(\beta)$
and is $\epsilon_1$-close to $\frak p$
and if $(\Sigma_{\frak Y'},u'')$ is $\epsilon_2$-close to
$(\Sigma_{\frak Y},u')$, then
$(\Sigma_{\frak Y'},u'')$ is $\epsilon_1 + o(\epsilon_2)$-close to $\frak p$.
(Here $\lim_{\epsilon_2\to 0} o(\epsilon_2) = 0$.)
\end{rem}
Let $\frak Y = \overline{\Phi}(\frak y,\vec T,\vec \theta)$ and $u' : (\Sigma_{\frak Y} ,\partial\Sigma_{\frak Y}) \to (X,L)$ be a smooth map in
homology class $\beta$ such that $(\Sigma_{\frak Y},u')$ is $\epsilon$-close to $\frak p$. We assume that $\epsilon$ is smaller than the injectivity radius of $X$.
Let ${\rm v} \in C^0(\mathcal G)$.
\begin{defn}\label{Emovevvv}
Suppose that we are given an obstruction bundle data $\frak E_{\frak p}$
centered at $\frak p$.
We define a map
\begin{equation}\label{Ivpdefn}
I^{{\rm v},\frak p}_{(\frak y,u),(\frak Y,u')} : E_{\frak p,\rm v}(\frak y,u) \to \Gamma_0(\text{\rm Int}\,K^{\rm obst}_{\rm v}; (u')^*TX \otimes \Lambda^{01})
\end{equation}
by using the complex linear part of the parallel transport along
the path of the form $t \mapsto {\rm E}(u(z),tv)$, where  $ {\rm E}(u(z),v) = u'(z)$.
(Note this is a short geodesic joining $u(z)$ and $u'(z)$ with respect to the connection which we used to define $\rm E$.)
Here we identify
$$
K^{\rm obst}_{\rm v} \subset K_{\rm v} \subset \Sigma_{\frak y}, \quad
K^{\rm obst}_{\rm v} \subset K_{\rm v} \subset \Sigma_{\frak Y}.
$$
\par
We write the image of (\ref{Ivpdefn}) by  $E_{\frak p,\rm v}(\frak Y,u')$.
\end{defn}
The map $I^{{\rm v},\frak p}_{(\frak y,u),(\frak Y,u')}$ is  $\Gamma_{(\frak x\cup \vec w,u)}^{\frak H^+_{(\frak x,u)}}$ invariant in the sense of
Lemma \ref{ptequiv} below.
Note  we have an injective homomorphism
$\Gamma_{(\frak x\cup \vec w,u)}^{\frak H^+_{(\frak x,u)}} \to \frak S_{\ell} \times \frak S_{\ell'}$
such that the $\Gamma_{(\frak x\cup \vec w,u)}^{\frak H^+_{(\frak x,u)}}$ action on the elements of
$\frak V(\frak x \cup \vec w)$ is identified with the permutation of the $\ell$ marked points in $\frak x$ and $\ell'$
marked points $\vec w$. (See (\ref{2132}).)
For $v \in \Gamma_{(\frak x\cup \vec w,u)}^{\frak H^+_{(\frak x,u)}}$ we define $v_*\frak Y$ by
permuting the marked points of $\frak Y$ in the same way.
If $(\frak Y,u')$ is $\epsilon$-close to $\frak p$ then $(v_*\frak Y,u')$ is $\epsilon$-close to $\frak p$.
Let $\rm v'$ be the vertex which is mapped from $\rm v$ by $v$ with respect to the $\Gamma_{(\frak x\cup \vec w,u)}^{\frak H^+_{(\frak x,u)}}$
action of $\mathcal G$. (See the discussion about Definition \ref{obbundeldata} (5)  we gave right above Remark \ref{rem23}.)
We remark that
$v_*\frak Y = \overline{\Phi}(v_*\frak y,v_*\vec T,v_*\vec \theta)$.
By using diffeomorphism in Definition \ref{def29}, we have a map
$v : \frak y \to v_*\frak y$.
Note there exists a map (diffeomorphism) $v : \Sigma_{\frak y} \to \Sigma_{\frak y}$ that permutes the marked points in the required way. However
this map is not holomorphic in general. It becomes biholomorphic
as  a map $v : \Sigma_{\frak y} \to \Sigma_{v_*\frak y}$.
\begin{lem}\label{ptequiv}
The following diagram commutes.
\begin{equation}
\CD
E_{\frak p,\rm v}(\frak y,u)
@>{I^{{\rm v},\frak p}_{(\frak y,u),(\frak Y,u')}}>> \Gamma_0(\text{\rm Int}\,K^{\rm obst}_{\rm v}(\frak Y); (u')^*TX \otimes \Lambda^{01})
\\
 @V{v_*}VV  @V V{v_*}V \\
E_{\frak p,\rm v'}(v_*\frak y,u\circ v^{-1})
@>{I^{{\rm v}',\frak p}_{(v_*\frak y,u),(v_*\frak Y,u'\circ v^{-1})}}>> \Gamma_0(\text{\rm Int}\,K^{\rm obst}_{\rm v'}(v_*\frak Y); (u'\circ v^{-1})^*TX \otimes \Lambda^{01})
\endCD
\label{A1.16}\end{equation}
Here we define $\rm v'$ by $v(K_{\rm v}) = K_{\rm v'}$.
\end{lem}
\begin{proof} The lemma follows from the fact that parallel transport etc. is
independent of the enumeration of the marked points.
(Note the left vertical arrow is well-defined by Definition \ref{obbundeldata} (5).)
\end{proof}
\begin{cor}\label{vindependenEcor}
$$
v_*\left(\bigoplus_{\rm v \in C^0(\mathcal G)}E_{\frak p,\rm v}(\frak Y,u')\right)
=
\bigoplus_{\rm v \in C^0(\mathcal G)}E_{\frak p,\rm v}(v_*\frak Y,u'\circ v^{-1}).
$$
\end{cor}
This is a consequence of Lemma \ref{ptequiv} and Definition \ref{obbundeldata} (5).
\par\medskip
We next show that the Fredholm regularity  (Definition \ref{obbundeldata} (6)) and
evaluation map transversality  (Definition \ref{obbundeldata} (7)) are preserved when
we take $(\frak Y,u')$ that is $\epsilon$-close to $\frak p$.
(See  Proposition \ref{linearMV}.)
To state them precisely we need some preparation.
\par
Let $\frak Y = \overline{\Phi}(\frak y,\vec T,\vec \theta)$ be an element of $\mathcal M_{k+1,\ell+\ell'}$ that is represented by
$(\Sigma_{\frak Y},\vec z_{\frak Y}, \vec z^{\text{\rm int}}_{\frak Y} \cup \vec w_{\frak Y})$.
We denote by $\mathcal G_{\frak Y}$ the combinatorial type of $\frak Y$.
(Here $\mathcal G_{\frak y}$ is the combinatorial type of $\frak y$ and
$\mathcal G_{\frak Y}$ is obtained from $\mathcal G_{\frak y}$
by shrinking the edges $\rm e$ such that $T_{\rm e} \ne \infty$.)
Let $\rm v \in C^0_{\rm d}(\mathcal G_{\frak Y})$. We have a differential operator
\begin{equation}\label{ducomponents20}
\aligned
D_{u',\rm v} \overline\partial : L^2_{m+1,\delta}((\Sigma_{\frak Y_{\rm v}},\partial \Sigma_{\frak Y_{\rm v}});&
(u')^*TX, (u')^*TL) \\
&\to
L^2_{m,\delta}(\Sigma_{\frak Y_{\rm v}}; (u')^*TX \otimes \Lambda^{01}).
\endaligned
\end{equation}
In case $\rm v \in C^0_{\rm s}(\mathcal G_{\frak Y})$ we have
\begin{equation}\label{ducomponents20c}
D_{u',\rm v} \overline\partial : L^2_{m+1,\delta}(\Sigma_{\frak Y_{\rm v}}; (u')^*TX)
\to
L^2_{m,\delta}(\Sigma_{\frak Y_{\rm v}}; (u')^*TX \otimes \Lambda^{01}).
\end{equation}
\begin{defn}\label{fredreg}
We say $(\frak Y,u')$ is {\it Fredholm regular} with respect to the obstruction bundle data $\frak E_{\frak p}$
if the sum of the image of (\ref{ducomponents20}) and $E_{\frak p,\rm v}(\frak Y,u')$ is
$L^2_{m,\delta}(\Sigma_{\frak Y_{\rm v}}; (u')^*TX \otimes \Lambda^{01})$ and
if the sum of the image of  (\ref{ducomponents20c}) and $E_{\frak p,\rm v}(\frak Y,u')$ is
$L^2_{m,\delta}(\Sigma_{\frak Y_{\rm v}}; (u')^*TX \otimes \Lambda^{01})$.
\end{defn}
Using this terminology, Definition \ref{obbundeldata} (6) means that
$(\frak x,u)$ is Fredholm regular with respect to the obstruction bundle data $\frak E_{\frak p}$.
\par
We next define the notion of evaluation map transversality.
\begin{defn}\label{maptrans}
A {\it flag} of $\mathcal G$ is a pair $(\mathrm v,\mathrm e)$ of edges $\rm e$ and its vertex $\rm v$.
Suppose $\mathcal G$ is oriented.  We say a flag $(\mathrm v,\mathrm e)$ is {\it incoming} if $\rm e$ is an
incoming edge. Otherwise it is said {\it outgoing}.
We denote by $z_{\rm e}$ the singular point corresponding to an edge $\rm e$.
\end{defn}
For each flag $(\rm v,\rm e)$ of $\mathcal G_{\frak Y}$, we define
\begin{equation}\label{evaluflagcomp2}
\text{\rm ev}_{\rm v,e} :
L^2_{m+1,\delta}((\Sigma_{\frak Y_{\rm v}},\partial \Sigma_{\frak Y_{\rm v}});
(u')^*TX, (u')^*TL)
\to T_{u'(z_{\rm e})}L,
\end{equation}
if ${\rm v} \in C^0_{\rm d}(\mathcal G_{\frak Y})$, $\rm e \in C^1_{\rm o}(\mathcal G_{\frak Y})$ in the same way as (\ref{evaluflagcomp}),
\begin{equation}\label{evaluflagcomp23}
\text{\rm ev}_{\rm v,e} :
L^2_{m+1,\delta}((\Sigma_{\frak Y_{\rm v}},\partial \Sigma_{\frak Y_{\rm v}}) ;
(u')^*TX, (u')^*TL)
\to T_{u'(z_{\rm e})}X,
\end{equation}
if ${\rm v} \in C^0_{\rm d}(\mathcal G_{\frak Y})$,
$\rm e \in C^1_{\rm c}(\mathcal G_{\frak Y})$
in the same way as (\ref{evaluflagcomp20}),
and
\begin{equation}\label{evaluflagcompcl2}
\text{\rm ev}_{\rm v,e} :
L^2_{m+1,\delta}(\Sigma_{\frak Y_{\rm v}};
(u')^*TX)
\to T_{u'(z_{\rm e})}X,
\end{equation}
if ${\rm e} \in C^1_{\rm c}(\mathcal G_{\frak Y})$ in the same way as (\ref{evaluflagcompcl}).
\par
Combining them we obtain
\begin{equation}\label{evaluationmap2diff}
\aligned
\text{\rm ev}_{\mathcal G_{\frak Y}} :
&\bigoplus_{\text{\rm v} \in C^0_{\rm d}(\mathcal G_{\frak Y})}
L^2_{m+1,\delta}((\Sigma_{\frak Y_{\rm v}},\partial \Sigma_{\frak Y_{\rm v}});
(u')^*TX, (u')^*TL)
\\
&\bigoplus_{\text{\rm v} \in C^0_{\rm s}(\mathcal G_{\frak Y})}
L^2_{m+1,\delta}(\Sigma_{\frak Y_{\rm v}};
(u')^*TX)
\\
&\to
\bigoplus_{{\rm e} \in C^1_{\rm o}(\mathcal G_{\frak Y})} T_{u'(z_{\rm e})}L
\oplus
\bigoplus_{{\rm e} \in C^1_{\rm c}(\mathcal G_{\frak Y})} T_{u'(z_{\rm e})}X.
\endaligned
\end{equation}
\begin{defn}\label{def:1720}
Suppose $(\frak Y,u')$ is Fredholm regular with respect to the obstruction bundle data $\frak E_{\frak p}$.
We say that
$(\frak Y,u')$ is {\it evaluation map transversal} with respect to the obstruction bundle data $\frak E_{\frak p}$
if the restriction of (\ref{evaluationmap2diff}) to the direct sum of the kernels of (\ref{evaluflagcomp2}), (\ref{evaluflagcomp23}) and of
(\ref{evaluflagcompcl2}) is surjective.
\end{defn}
Using this terminology, Definition \ref{obbundeldata} (7) means that
$(\frak x,u)$ is evaluation map transversal with respect to the obstruction bundle data $\frak E_{\frak p}$.
\par\medskip
Proposition \ref{linearMV} below says that Fredholm regularity and evaluation map transversality are preserved if
$(\frak Y,u')$ is sufficiently close to $\frak p$.
To state it we need to note the following point.
\par
When we define $\epsilon$-close-ness, we put the condition that the image of each
connected component of the neck region has diameter $< \epsilon$.
But we did not assume a similar condition for $\frak p$ and $\frak E_{\frak p}$ itself. So
in case when this condition is not satisfied for $\frak p$, there can not exist any object that is
$\epsilon$-close to $\frak p$. Especially $\frak p$ itself is not $\epsilon$-close to $\frak p$.
\par
However,
we can always modify the core $K_{\rm v}$ so that $\frak p$ itself becomes $\epsilon$-close to $\frak p$ as follows.
We take a positive number $R_{(\rm v,\rm e)}$ for each flag of $\mathcal G$ and write $\vec R$ the totality of
such $R_{(\rm v,\rm e)}$.
We put
\begin{equation}\label{extendedcore}
\aligned
K^{+\vec R}_{\rm v} =
K_{\rm v}
&\cup \bigcup_{{\rm e}\in C^1_{\mathrm o}(\mathcal G)
\atop \text{$(\rm v,{\rm e})$ is an outgoing flag}} (0,R_{(\rm v,\rm e)}] \times [0,1]
\\
&\cup
\bigcup_{{\rm e}\in C^1_{\mathrm o}(\mathcal G)
\atop \text{$(\rm v,{\rm e})$ is an incoming flag}} [-R_{(\rm v,\rm e)},0) \times [0,1]
\\
&\cup
\bigcup_{{\rm e}\in C^1_{\mathrm c}(\mathcal G)
\atop \text{$(\rm v,{\rm e})$ is an outgoing flag} }(0,R_{(\rm v,\rm e)}] \times S^1
\\
&\cup
\bigcup_{{\rm e}\in C^1_{\mathrm c}(\mathcal G)
\atop \text{$(\rm v,{\rm e})$ is an incoming flag}} [-R_{(\rm v,\rm e)},0) \times S^1.
\endaligned
\end{equation}
\begin{defn}
We can define an obstruction bundle data $\frak E_{\frak p}$ centered at $\frak p$ using $K^{+\vec R}_{\rm v}$ in place of $K_{\rm v}$.
We call it the obstruction bundle data obtained by {\it extending the core} and write
$\frak E_{\frak p}^{+ \vec R}$.
We call (\ref{extendedcore}) the {\it extended core}. (In case we need to specify $\vec R$ we call it the
$\vec R$-extended core.)
(\ref{extendedcore}) is a generalization of (\ref{1104}).
\end{defn}
\begin{prop}\label{linearMV}
Let $\frak p \in \mathcal M_{k+1,\ell}(\beta)$ and $\frak E_{\frak p}$ be an obstruction bundle data centered at $\frak p$.
Then there exist $\epsilon > 0$ and $\vec R$ with the following properties.
\begin{enumerate}
\item
If $(\frak Y,u')$ is $\epsilon$-close to $\frak p$ with respect to $\frak E_{\frak p}^{+ \vec R}$,
then $(\frak Y,u')$ is Fredholm regular with respect to $\frak E_{\frak p}^{+ \vec R}$.
\item
If $(\frak Y,u')$ is $\epsilon$-close to $\frak p$ with respect to $\frak E_{\frak p}^{+ \vec R}$,
then $(\frak Y,u')$ is evaluation map transversal with respect to $\frak E_{\frak p}^{+ \vec R}$.
\item
$\frak p$ is $\epsilon$-close to $\frak p$ with respect to $\frak E_{\frak p}^{+ \vec R}$.
\end{enumerate}
\end{prop}
\begin{proof}
By using the fact that the diameter of the $u'$ image of the connected component  of the neck region is small,
we can prove an exponential decay estimate of $u'$ on the neck region. This is an analogue of Lemma
\ref{neckaprioridecay} and its proof is the same as the proof of
\cite[Lemma 11.2]{FOn}.
Then the rest of the proof of (1),(2) is a version of the proof of Mayer-Vietoris principle of Mrowka \cite{Mr}.
See \cite[Proposition 7.1.27]{fooo:book1} or \cite[Lemma 8.5]{Fuk96II}.
(3) is obvious.
\end{proof}
So far we have discussed the case of bordered genus zero curve.
The case of genus zero curve without boundary is the same so we do not repeat it.
\footnote{Higher genus case is also the same.}
\par\medskip
\section{The differential equation and thickened moduli space}
\label{settin2}

To construct a Kuranishi neighborhood of each point in
our moduli space $\mathcal M_{k+1,\ell}(\beta)$ or $\mathcal M_{\ell}^{\rm cl}(\alpha)$,
we need to assign an obstruction bundle to each point of it.
To do so we follow the way we had written in \cite[end of the page 1003]{FOn} and
\cite[end of the page 423-middle of page 424]{fooo:book1}.
The outline of the argument is as follows.
For each $\frak p \in \mathcal M_{k+1,\ell}(\beta)$ we take an obstruction bundle data
$\frak E_{\frak p}$.
We then consider a closed neighborhood $\frak W_{\frak p}$ of $\frak p$ in  $\mathcal M_{k+1,\ell}(\beta)$
so  that its elements together with certain marked points added is
$\epsilon_{\frak p}$-close to $\frak p$ with respect to $\frak E_{\frak p}^{+\vec R}$.
Here we choose $\epsilon_{\frak p}$ and $\frak E_{\frak p}^{+\vec R}$ so that
Proposition  \ref{linearMV} holds.
We next take a finite number of $\frak p_c \in \mathcal M_{k+1,\ell}(\beta)$ such that
$$
\bigcup_{c} {\rm Int}\,\frak W_{\frak p_c} = \mathcal M_{k+1,\ell}(\beta).
$$
For $\frak p \in \mathcal M_{k+1,\ell}(\beta)$, we collect all $E_{\frak p_c}$
such that $\frak p_c$ satisfies $\frak p \in \frak W_{\frak p_c}$. The sum will be the obstruction bundle $\mathcal E_{\frak p}$
at $\frak p$.
Now we will describe this process in more detail below.
\par
We first define the subset $\frak W_{\frak p}$ in more detail.
We note that in Definition \ref{epsiloncloseto}, we need $\ell + \ell_{\frak p}$
interior marked points to define its $\epsilon$-close-ness to
an element $\frak p \in \mathcal M_{k+1,\ell}(\beta)$.
(Here $\ell_{\frak p}$ is the number of marked points we add
as a part of the  obstruction bundle data
$\frak E_{\frak p}$.)
We start with describing the process of forgetting those $\ell_{\frak p}$ marked points.
\par
\begin{defn}\label{transconst}
We consider the situation of Definition \ref{epsiloncloseto}.
Let $\frak Y = \Phi(\frak y,\vec T,\vec \theta)$ and let $u' : (\Sigma_{\frak Y} ,\partial\Sigma_{\frak Y}) \to (X,L)$ be a smooth map in the homology class $\beta$ that is $\epsilon$-close to $\frak p$.
We say $(\frak Y,u')$ satisfies the {\it transversal constraint} if for each $w_i \in \vec w$ we have
\begin{equation}
u'(w_i) \in \mathcal D_{\frak p,i}.
\end{equation}
\end{defn}
Let us explain the notation appearing in the above definition.
We have $\vec w_{\frak p}$, the additional marked points on $\Sigma_{\frak p}$ as a part of the obstruction bundle data
$\frak E_{\frak p}$. The element $\frak y$ is in a neighborhood $\frak V(\frak x_{\frak p}\cup \vec w_{\frak p})$.
(This neighborhood $\frak V(\frak x_{\frak p}\cup \vec w_{\frak p})$  is also a part of the date $\frak E_{\frak p}$.)
$(\vec T,\vec \theta)$ is as in Definition \ref{def214}.
Thus $\frak Y = {\overline\Phi}(\frak y,\vec T,\vec \theta)$ is a bordered genus zero curve with $k+1$
boundary and $\ell + \ell_{\frak p}$ interior marked points. ($\ell_{\frak p}$ is the number of points in $\vec w_{\frak p}$.)
We denote by $w_{\frak p,i}$ the $(\ell + i)$-th interior marked point. (It is $i$-th among the additional marked points.)
For each $i=1,\dots,\ell_{\frak p}$,  we took $\mathcal D_{\frak p,i}$ that is transversal to $u_{\frak p}(\Sigma_{\frak p})$ at $u_{\frak p}(w_i)$
as a part of the data $\frak E_{\frak p}$.
\par
\begin{lem}\label{transpermutelem}
For each $\frak p \in \mathcal M_{k+1,\ell}(\beta)$ and an obstruction bundle data $\frak E_{\frak p}$
centered at $\frak p$ there exists $\epsilon_{\frak p}$ such that the following holds.
\par
Let $\frak q = (\frak x_{\frak q},u_{\frak q}) \in \mathcal M_{k+1,\ell}(\beta)$.
We consider the set of symmetric marking $\vec w'_{\frak p}$ of $\frak x_{\frak q}$ with $\#\vec w'_{\frak p} = \ell_{\frak p}$,
such that the following holds.
\begin{enumerate}
\item
There exists $\frak y\in \frak V(\frak x_{\frak p}\cup \vec w_{\frak p})$ and
$(\vec T,\vec \theta) \in (\vec T^{\rm o}_0,\infty] \times ((\vec T^{\rm c}_0,\infty] \times \vec S^1)$
such that $\frak x_{\frak q}\cup\vec w'_{\frak p} = \overline{\Phi}(\frak y,\vec T,\vec \theta)$.
\item
$(\frak x_{\frak q}\cup \vec w'_{\frak p}, u_{\frak q})$ is $\epsilon_{\frak p}$-close to
$\frak p$.
\item
$(\frak x_{\frak q}\cup\vec w'_{\frak p}, u_{\frak q})$ satisfies the transversal constraint.
\end{enumerate}
Then the set of such $\vec w'_{\frak p}$ consists of a single $\Gamma_{\frak p}$ orbit
if it is nonempty. Here we regard $\Gamma_{\frak p} \subset \frak S_{\ell'}$
by (\ref{sigmahomotoS}) and $\Gamma_{\frak p}$ acts on the set of  $\vec w'_{\frak p}$'s by permutation.
\end{lem}
The proof of Lemma \ref{transpermutelem} is not difficult. We however postpone its proof to
Section \ref{cutting} where the transversal constraint is studied more systematically.
\par
We are now ready to provide the definition of $\frak W_{\frak p} \subset \mathcal M_{k+1,\ell}(\beta)$.
\par
First for each $\frak p \in  \mathcal M_{k+1,\ell}(\beta)$ we take and fix an
obstruction bundle data $\frak E_{\frak p}$.
Let $\vec w_{\frak p}$ be the additional marked points we take as a part of $\frak E_{\frak p}$.
We take $\epsilon_{\frak p}$ so that Proposition \ref{linearMV} and Lemma \ref{transpermutelem} hold.
Moreover we may change $\frak E_{\frak p}$ if necessary so that Proposition \ref{linearMV}
holds for $\frak E_{\frak p}^{+ \vec R} = \frak E_{\frak p}$.
\begin{defn}\label{openW+++}
$\frak W^+(\frak p)$ is the set of all $\frak q \in  \mathcal M_{k+1,\ell}(\beta)$
such that the set of $\vec w'_{\frak p}$ satisfying (1)-(3) of Lemma \ref{transpermutelem} is nonempty.
The constant $\epsilon_{\frak p}$ (which is often denoted by $\epsilon_{\frak p_c}$ or $\epsilon_c$)
is determined later. (See
Lemma \ref{nbdregmaineq} (Remark \ref{Lem258}), Proposition \ref{forgetstillstable}, Lemma \ref{setisopen},
Lemma \ref{injectivitypp3}, Sublemma \ref{lub296lem} (Remark \ref{rem298}.)
See also 2 lines above Definition \ref{obbundle1}.)
We note that $\frak W^+(\frak p)$ is open, as we will see in Subsection \ref{cutting}. See Remark \ref{rem:2111}.
\par
We choose a compact subset $\frak W_{\frak p} \subset \frak W^+(\frak p)$  that is a neighborhood of
$\frak p$.
We take $\frak W^0_{\frak p}$ that is a compact subset of ${\rm Int}\,\,\frak W_{\frak p}$
and is a neighborhood of
$\frak p$.
\par
We take a finite set $\{\frak p_{c} \mid c \in \frak C\} \subset \mathcal M_{k+1,\ell}(\beta)$ such that
\begin{equation}\label{mpccoverM}
\bigcup_{c \in \frak C} \text{\rm Int}\,\, \frak W^0_{\frak p_c} = \mathcal M_{k+1,\ell}(\beta).
\end{equation}
\end{defn}
We fix this set $\{\frak p_{c} \mid c \in \frak C\}$ in the rest of the
construction of the Kuranishi structure.
From now on none of the obstruction bundle data at $\frak p$ for
$\frak p \notin \frak C$ is used in this note.
\begin{defn}\label{defn243}
For $\frak p \in \mathcal M_{k+1,\ell}(\beta)$, we define
$$
\frak C(\frak p) = \{ c \in \frak C \mid \frak p \in \frak W_{\frak p_c}\}.
$$
We also choose additional marked points $\vec w_{c}^{\frak p}$
of $\frak x_{\frak p}$ for each $c\in \frak C(\frak p)$ such that
\begin{enumerate}
\item There exist $\frak y\in \frak V(\frak x_{\frak p_c}\cup \vec w_{\frak p_c})$ and
$(\vec T,\vec \theta) \in (\vec T^{\rm o}_0,\infty] \times ((\vec T^{\rm c}_0,\infty] \times \vec S^1)$
such that $\frak x_{\frak p}\cup\vec w_c^{\frak p} = \overline{\Phi}(\frak y,\vec T,\vec \theta)$.
\item
$(\frak x_{\frak p}\cup \vec w_c^{\frak p},u_{\frak p})$ is $\epsilon_{\frak p_c}$-close to
$\frak p_c$.
\item
$(\frak x_{\frak p}\cup \vec w_c^{\frak p},u_{\frak p})$ satisfies the transversal constraint.
\end{enumerate}
\end{defn}
\begin{lem}
For each $\frak p$ there exists a neighborhood $U$ of it so that if
$\frak q \in U$ then
$$
\frak C(\frak q) \subseteq \frak C(\frak p).
$$
\end{lem}
\begin{proof}
The lemma follows from the fact that $ \frak W_{\frak p_c}$ is closed.
\end{proof}
\par
We next define an obstruction bundle $\mathcal E_{\frak p}$ for each $\frak p = (\frak x_{\frak p},u_{\frak p})
\in \mathcal M_{k+1,\ell}(\beta)$.
Take $c \in  \frak C(\frak p)$. Let $\vec w^{\frak p}_c$ be as in Definition \ref{defn243}.
By Definition \ref{Emovevvv},
the map
\begin{equation}\label{Ivpdefn22}
I^{{\rm v},\frak p_c}_{(\frak y_c,u_c),\frak p\cup \vec w^{\frak p}_c)} : E_{\frak p_c,\rm v}(\frak y_c,u_c)
\to \Gamma_0(\text{\rm Int}\,K^{\rm obst}_{\rm v}; u_{\frak p}^*TX \otimes \Lambda^{01})
\end{equation}
is defined. Here $\frak x_{\frak p} \cup \vec w^{\frak p}_c = \overline{\Phi}(\frak y_c,\vec T_c,\vec \theta_c)$ and
$\frak y_c \in \frak V(\frak x_{\frak p_c}\cup \vec w_{\frak p_c})$.
Note
$
K_{\rm v}^{\rm obst} \subset K_{\rm v} \subset \frak x_{\frak p_c}.
$
We have also $K_{\rm v}^{\rm obst} \subset \frak x_{\frak p}$
since $\vec w^{\frak p}_c \cup \frak x_{\frak p} = \overline{\Phi}(\frak y_c,\vec T_c,\vec \theta_c)$.
\begin{lem}
The image
$
E_{c}(\frak p)
$ of (\ref{Ivpdefn22}) depends only on
$\frak p \in  \frak W^+_{\frak p_c}$ and  is independent of the choices of $\vec w_c^{\frak p}$
satisfying Definition \ref{defn243} (1)-(3).
\end{lem}
\begin{proof}
This is a consequence of Corollary \ref{vindependenEcor} and Lemma \ref{transpermutelem}.
\end{proof}
\begin{defn}\label{def:187}
We define
\begin{equation}\label{kuraequation}
\mathcal E_{\frak C(\frak p)}(\frak p)
= \sum_{c \in \frak C(\frak p)} E_{c}(\frak p).
\end{equation}
For $\frak A \subset \frak C(\frak p)$ we put
\begin{equation}
\mathcal E_{\frak A} (\frak p)
= \sum_{c \in \frak A} E_{c}(\frak p).
\end{equation}
\end{defn}
The defining equation of the thickened moduli space at $\frak p$ is
$$
\overline\partial u_{\frak p} \equiv 0  \mod \mathcal E_{\frak C(\frak p)}({\frak p}).
$$
We need to extend the subspace $\mathcal E_{\frak C(\frak p)}({\frak p})$ to a family of
subspaces parametrized by a neighborhood of $\frak p$.
Before doing so we need the following.
\begin{lem}\label{transbetweenEs}
By perturbing $E_{\frak p_c}$ (that is a part of the obstruction bundle data $\frak E_{\frak p_c}$) we may assume that
$$
E_{c}(\frak p) \cap E_{c'}(\frak p) = \{0\},
$$
if $c,c' \in \frak C(\frak p)$ and $c\ne c'$.
\end{lem}
\begin{proof}
The proof will be written in Section \ref{Lemma49}.
\end{proof}
Now we start extending the equation (\ref{kuraequation}) to an element $\frak q$ in a  `neighborhood' of $\frak p$.
We do not yet assume that $\frak q$ satisfies the transversal constraint (Definition \ref{transconst}).
So to define $E_{c}(\frak q)$ we need to include $\vec w'_c$ for all $c \in \frak C(\frak p)$ as marked points of $\frak q$.
We also take more marked points $\vec w_{\frak p}$ to stabilize $\frak p$ and take corresponding additional
marked points $\vec w'_{\frak p}$ on $\Sigma_{\frak q}$.
The marked points $\vec w_{\frak p}$ are used to fix the coordinate to perform the gluing construction in section \ref{glueing}.
$\vec w'_c$ is used to define the map (\ref{Ivpdefn22}).
Thus they have different roles.
\par
A technical point to take care of is the following.
We may assume that the $\ell_c$ components of $\vec w_c^{\frak p}$ are mutually different, for each $c$.
(This is because $\ell_c$ components of $\vec w_{\frak p_c}$ are mutually different.)
However there is no obvious way to arrange so that $\vec w_c^{\frak p} \cap \vec w_{c'}^{\frak p} = \emptyset$
for $c \ne c'$.
Note, in the usual stable map compactification, at the point where two or more marked points become coincide,
we put the `phantom bubble' so that they become  different points on this bubbled component.
For our purpose, the proof becomes simpler when we do {\it not} put a phantom bubble in case one of the components of $\vec w_{c}^{\frak p}$ coincides with
one of the components of $\vec w_{c'}^{\frak p}$ for $c \ne c'$.
Taking these points into account we  define $\mathcal M_{k+1,(\ell,\ell_{\frak p},(\ell_c))}(\beta,\frak p)_{\epsilon_{0},\vec T_{0}}$ below.
\par
We first review the situation we are working in and prepare some notations.
Let $\frak p \in \mathcal M_{k+1,\ell}(\beta)$.
We defined $\frak C(\frak p)$ in Definition \ref{defn243}.
For $c \in \frak C(\frak p)$ we fixed an obstruction bundle data $\frak E_{\frak p_c}$ centered at $\frak p_c$.
Additional marked points $\vec w_{\frak p_c}$ is a part of the data $\frak E_{\frak p_c}$. We put $\ell_c = \#\vec w_c$.
We also put $\epsilon_c = \epsilon_{\frak p_c}$ where the right hand side is as in Lemma \ref{transpermutelem}.
As mentioned before we take $\frak E_{\frak p_c}$ so that
Proposition \ref{linearMV}
holds for $\frak E_{\frak p_c}^{+ \vec R} = \frak E_{\frak p_c}$.
\begin{defn}\label{stabdata}
A {\it stabilization data} at $\frak p$ is the data as follows.
\begin{enumerate}
\item
A symmetric stabilization $\vec w_{\frak p} = (w_{\frak p,1},\dots,w_{\frak p,\ell_{\frak p}})$ of $\frak p$.
Let $\ell_{\frak p} = \#\vec w_{\frak p}$.
\item
For each $w_{\frak p,i}$ ($i=1,\dots,\ell_{\frak p}$), we take and fix $\mathcal D_{\frak p,i}$
such that it is a codimension two submanifold of $X$ and is transversal to $u_{\frak p}$ at
$u_{\frak p}(w_{\frak p,i})$. We also assume $u_{\frak p}(w_{\frak p,i}) \in \mathcal D_{\frak p,i}$.
\item
We assume that $\{\mathcal D_{\frak p,i}\mid i=1,\dots,\ell_{\frak p}\}$ is invariant under the
$\Gamma_{\frak p}$ action in the same sense as in Definition \ref{obbundeldata} (8) (\ref{Dequivalent}).
\item A coordinate at infinity of $\frak p \cup \vec w_{\frak p}$.
\item $\vec w_{\frak p} \cap \vec w_c^{\frak p} = \emptyset$ for any $c \in \frak C(\frak p)$.
\item
Let $K^{\rm obst}_{{\rm v},c}$ be the support of the obstruction bundle as in Definition \ref{obbundeldata} (4).
(Here ${\rm v} \in C^0(\mathcal G_{\frak p_c})$.)
Since $\frak x_{\frak p} = \overline{\Phi}(\frak y,\vec T,\vec \theta)$ we may regard
$K^{\rm obst}_{{\rm v},c} \subset \Sigma_{\frak p}$. We require
$$
K^{\rm obst}_{{\rm v},c} \subset \bigcup_{{\rm v}' \in C^0(\mathcal G_{\frak p})} {\rm Int} K_{{\rm v}'}.
$$
Here the right hand side is the core of the coordinate at infinity given by item (4)
Definition \ref{stabdata}.
\end{enumerate}
\end{defn}
A stabilization data at $\frak p$ is similarly defined as the obstruction bundle data centered at $\frak p$.
But it does not include $K_{\rm v}^{\rm obst}$ or $E_{\frak p,\rm v}$.
The stabilization data at $\frak p$ has no relation to the
obstruction bundle data at $\frak p$.\footnote{In case $\frak p = \frak p_c$ we have
both stabilization data and obstruction bundle data at $\frak p$.
The notation $\vec w_{\frak p}$ is used for both structures.
They may not be coincide. We use the same symbol for both since this can not cause any confusion and
the case $\frak p = \frak p_c$ does not play a role in our discussion.}
\par
We fix a metric on all the Deligne-Mumford moduli spaces.
Let $\frak V_{\epsilon_0}(\frak p \cup \vec w_{\frak p})$ be the $\epsilon_0$-neighborhood of
$\frak p \cup \vec w_{\frak p}$ in $\mathcal M_{k+1,\ell+\ell_{\frak p}}(\mathcal G_{(\frak p \cup \vec w_{\frak p})})$ where
$\mathcal G_{(\frak p \cup \vec w_{\frak p})}$ is the combinatorial type of $\frak p \cup \vec w_{\frak p}$.
\begin{defn}
[{\rm Definition of $
{\frak U}_{k+1,(\ell;\ell_{\frak p},(\ell_c))}(\beta,\frak p;\frak B)_{\epsilon_{0},\vec T_{0}}
$}]
\label{defn251}
We fix a stabilization data at $\frak p$ and an obstruction bundle data
centered at $\frak p_{c}$ for each $c \in \frak C(\frak p)$.
Let $\frak B \subset \frak C(\frak p)$.
For each $c \in \frak G(\frak p)$ we chose
$\vec w^{\frak p}_c$  in Definition \ref{defn243}.
\par
For $\epsilon_0 > 0$ and $\vec T_{0} = (\vec T^{\rm o}_0,\vec T^{\rm c}_0) =  (T_{{\rm e},0} : {\rm e} \in C^1(\mathcal G_{\frak p}))$
we consider the set of all $(\frak Y,u',(\vec w'_c ; c\in \frak B))$ such that
the following holds for some $\vec R$.
\begin{enumerate}
\item There exist $\frak y \in \frak V_{\epsilon_0}(\frak p \cup \vec w_{\frak p})$,
$(\vec T,\vec \theta) \in (\vec T^{\rm o}_0,\infty] \times ((\vec T^{\rm c}_0,\infty] \times \vec S^1)$
such that
$$
\frak Y = \overline{\Phi}(\frak y,\vec T,\vec \theta) \in \mathcal M_{k+1,\ell+\ell_{\frak p}}.
$$
\item
$u'$ is $\epsilon_0$-close to $u_{\frak p}$ on the extended core $K_{\rm v}^{+\vec R}$ of $\Sigma_{\frak p}$ in $C^{10}$-topology.
We use the coordinate at infinity of $\frak p \cup \vec w_{\frak p}$ that is included in the
stabilization data at $\frak p$,  to define
this $C^{10}$ close-ness.
\item
Moreover we assume that the diameter of the $u'$ image of each neck region of $\Sigma_{\frak Y}$
is smaller than $\epsilon_0$.
We assume furthermore that $u'$ is pseudo-holomorphic in the neck regions.
(The neck region here is the complement of the union of the extended cores  $K_{\rm v}^{+\vec R}$.)
\item
We write $\frak Y = \frak Y_0 \cup \vec w'_{\frak p}(\frak Y)$ where $\vec w'_{\frak p}(\frak Y)$
are $\ell_{\frak p}$ marked points that correspond to $\vec w_{\frak p}$.
We assume that $(\frak Y_0\cup \vec w'_c,u')$ is
$\epsilon_{0}$-close to $\frak p \cup \vec w^{\frak p}_c$ in the sense of Definition \ref{epsiloncloseto} after
extending the core of $\frak p \cup \vec w_c^{\frak p}$ by $\vec R$.
\end{enumerate}
\par
We say that $(\frak Y^{(1)},u^{\prime (1)},(\vec w^{\prime (1)}_c ; c\in \frak B))$ is
{\it weakly equivalent} to $(\frak Y^{(2)},u^{\prime (2)},(\vec w^{\prime (2)}_c ; c\in \frak B))$
if there exists a bi-holomorphic map $v : \frak Y^{(1)}\to \frak Y^{(2)}$ such that
\begin{enumerate}
\item[(a)]
$u^{\prime (1)} = u^{\prime (2)}\circ v$.
\item[(b)]
$v(w^{\prime (1)}_{c,i}) = w^{\prime (2)}_{c,\sigma_c(i)}$, where $\sigma_c \in \frak S_{\ell_c}$.
\item[(c)]
$v$ sends the $i$-th boundary marked point of $\frak Y^{(1)}$ to the
$i$-th boundary marked point of $\frak Y^{(2)}$.
$v$ sends 1-st,\dots,$\ell$-th interior marked points of $\frak Y^{(1)}$ to the
corresponding interior marked points of $\frak Y^{(2)}$.
$v$ sends $\ell+1$,\dots,$\ell+k$,\dots$\ell+\ell_{\frak p}$-th interior marked points of
$\frak Y^{(1)}$ to the $\ell+\sigma(1)$,\dots,$\ell+\sigma(k)$,\dots$\ell+\sigma(\ell_{\frak p})$-th interior marked points of
$\frak Y^{(2)}$, where $\sigma \in \frak S_{\ell_{\frak p}}$.
\end{enumerate}
\par
We denote by $\overline{\frak U}_{k+1,(\ell;\ell_{\frak p},(\ell_c))}(\beta,\frak p;\frak B)_{\epsilon_{0},\vec T_{0}}$
the set of all weak equivalence classes of   $(\frak Y,u',(\vec w'_c ; c\in \frak B))$
satisfying (1)-(4) above. (Here we use the weak equivalence relation defined by (a), (b), (c).)
\par
We say that $(\frak Y^{(1)},u^{\prime (1)},(\vec w^{\prime (1)}_c ; c\in \frak B))$ is
{\it equivalent} to $(\frak Y^{(2)},u^{\prime (2)},(\vec w^{\prime (2)}_c ; c\in \frak B))$
when $\sigma = \sigma_c =$ identity is satisfied in (a)-(c) above in addition.
Let
$$
{\frak U}_{k+1,(\ell;\ell_{\frak p},(\ell_c))}(\beta,\frak p;\frak B)_{\epsilon_{0},\vec T_{0}}
$$
be the set of equivalence classes of this equivalence relation.
\end{defn}
\begin{lem}\label{ippaiequivalnce}
We may choose  $\epsilon_{0}$ sufficiently small so that the following holds.
Suppose $(\frak Y^{(1)},u^{\prime (1)},(\vec w^{\prime (1)}_c ; c\in \frak B))$ is
weakly equivalent to $(\frak Y^{(2)},u^{\prime (2)},(\vec w^{\prime (2)}_c ; c\in \frak B))$
in the above sense and
$
\frak Y^{(j)} = \overline{\Phi}(\frak y^{(j)} ,\vec T^{(j)} ,\vec \theta^{(j)} ) \in \mathcal M_{k+1,\ell+\ell_{\frak p}}.
$
Then we have
$$
(\frak y^{(2)} ,\vec T^{(2)} ,\vec \theta^{(2)} ) = v_* (\frak y^{(1)} ,\vec T^{(1)} ,\vec \theta^{(1)} )
$$
for some $v \in \Gamma_{\frak p} \subset \frak S_{\ell_{\frak p}}$.
\end{lem}
\begin{proof}
The proof is by contradiction.
Suppose there exists a sequence of positive numbers $\epsilon_{0,a} \to 0$ and $(u'_{(j),a},(\frak y^{(j),a} ,\vec T^{(j),a} ,\vec \theta^{(j),a}),
(\vec w^{\prime (j), a}_c ; c\in \frak B))$
for $j=1,2$ and $a=1,2,\dots$ such that:
\begin{enumerate}
\item
The object
$(u'_{(1),a},(\frak y^{(1),a} ,\vec T^{(1),a} ,\vec \theta^{(1),a}),(\vec w^{\prime (1), a}_c ; c\in \frak B))$ is
weakly equivalent to the object $(u'_{(2),a},(\frak y^{(2),a} ,\vec T^{(2),a} ,\vec \theta^{(2),a}),(\vec w^{\prime (2),a}_c ; c\in \frak B))$.
\item
$
\frak Y^{(j),a} = \overline{\Phi}(\frak y^{(j),a} ,\vec T^{(j),a} ,\vec \theta^{(j),a} ) \in \mathcal M_{k+1,\ell+\ell_{\frak p}}.
$
\item The objects $(u'_{(j),a},(\frak y^{(j),a} ,\vec T^{(j),a} ,\vec \theta^{(j),a}),(\vec w^{\prime (j),a}_c ; c\in \frak B))$ are representatives of
elements of $\frak U_{k+1,(\ell;\ell_{\frak p},(\ell_c))}(\beta,\frak p;\frak B)_{\epsilon_{0,a},\vec T_{0}}$.
\item
There is no $v \in \Gamma_{\frak p}$ satisfying
$
(\frak y^{(2),a} ,\vec T^{(2),a} ,\vec \theta^{(2),a}) = v_* (\frak y^{(1),a} ,\vec T^{(1),a} ,\vec \theta^{(1),a})
$.
\end{enumerate}
We will deduce contradiction. By assumption there exist $\vec R_a \to \infty$  and
biholomorphic maps $v_a : \frak Y^{(1),a} \to \frak Y^{(2),a}$ such that
\begin{enumerate}
\item[(I)]
$\vert u'_{(2),a} \circ v_a - u'_{(1),a}\vert_{C^{10}(K_{{\rm v}}^{+\vec R_a})} < \epsilon_{0,a}$.
\item[(II)]
The diameter of $u'_{(j),a}$ image of each connected component of the complement of the union
of the extended cores $K_{{\rm v}}^{+\vec R_a}$ is smaller than $\epsilon_{0,a}$.
\item[(III)]
$v_a(w^{\prime (1), a}_{c,i}) =  w^{\prime (2), a}_{c,\sigma_c(i)}$, where $\sigma_c \in \frak S_{\ell_c}$.
\item[(IV)]
$v_a$ sends the $i$-th boundary marked point of $\frak Y^{(1),a}$ to the
$i$-th boundary marked point of $\frak Y^{(2),a}$.
$v_a$ sends 1-st,\dots,$\ell$-th interior marked points of
$\frak Y^{(1),a}$ to the
corresponding interior marked points of $\frak Y^{(2),a}$.
$v_a$ sends $\ell+1$,\dots,$\ell+k$,\dots$\ell+\ell_{\frak p}$-th interior marked points of
$\frak Y^{(1),a}$ to the $\ell+\sigma_a(1)$,\dots,$\ell+\sigma_a(k)$,\dots$\ell+\sigma_a(\ell_{\frak p})$-th interior marked points of
$\frak Y^{(2),a}$, where $\sigma_a \in \frak S_{\ell_{\frak p}}$.
\end{enumerate}
By (I) and (II) we may take a subsequence (still denoted by the same symbol) such that
$v_a$ converges to a biholomorphic map $v : \Sigma_{\frak p} \to \Sigma_{\frak p}$
such that $u_{\frak p} \circ v = u_{\frak p}$.
Then (III) and (IV) imply that $v \in \Gamma_{\frak p}$.
\par
So  changing $\frak Y^{(2),a}$ by $v$ we may assume $v = {\rm identity}$.
Therefore $v_a$ converges to identity.
The stability then implies that $v_a$ is identity. This contradicts to (4).
\end{proof}
\begin{defn}\label{defEc}
Let $\frak q^+ = (\frak Y,u',(\vec w'_c ; c\in \frak B))
\in \frak U_{k+1,(\ell,\ell_{\frak p},(\ell_c))}(\beta;\frak p;\frak B)_{\epsilon_{0},\vec T_{0}}$.
We define
$$
E_c(\frak q^+)
\subset \bigoplus_{\rm v \in C^0(\mathcal G(\frak Y))} \Gamma_0(\text{\rm Int}\,K^{\rm obst}_{\rm v}; (u')^*TX \otimes \Lambda^{01})
$$
as follows, where $\mathcal G(\frak Y)$ is the combinatorial type of $(\frak Y,u')$.
We regard $K^{\rm obst}_{\rm v}$ as a subset of $\frak Y$.
We note that $\frak p \cup \vec w^{\frak p}_{c}$ is $\epsilon_{\frak p_c}$-close to $\frak p_{c} \cup \vec w_{\frak p_c}$
and $(\frak Y\cup \vec w'_c,u')$ is
$\epsilon_{0}$-close to $\frak p \cup \vec w^{\frak p}_{c}$
in the sense of Definition \ref{epsiloncloseto}.
Therefore we have
\begin{equation}\label{Ivpdefn2211}
I^{{\rm v},\frak p_c}_{(\frak y_c,u_c),(\frak Y\cup \vec w'_c,u')} : E_{\frak p_c,\rm v}(\frak y_c,u_c) \to \Gamma_0(\text{\rm Int}\,K^{\rm obst}_{\rm v};
(u')^*TX \otimes \Lambda^{01}).
\end{equation}
Here $\frak p_c = (\frak x_c,u_c)$ and $\frak Y\cup \vec w'_c = \overline{\Phi}(\frak y_c,\vec T,\vec \theta)$.
We regard $K^{\rm obst}_c$ as a subset of $\frak x_c$ also.
(Note that the core of $\frak Y$ is canonically identified with the core of $\frak y_c$.) Then
we define
\begin{equation}
E_c(\frak q^+) = \sum_{\rm v \in C^0(\mathcal G(\frak Y))}
I^{{\rm v},\frak p_c}_{(\frak y_c,u_c),(\frak Y\cup \vec w'_c,u')}(E_{\frak p_c,\rm v}(\frak y_c,u_c))
\end{equation}
and put
\begin{equation}\label{kuraequationplus}
\mathcal E_{\frak B}(\frak q^+)
= \sum_{c \in \frak B} E_{c}(\frak q^+).
\end{equation}
For $\frak A \subset \frak B$ we put
\begin{equation}
\mathcal E_{\frak A}(\frak q^+)
= \sum_{c \in \frak A} E_{c}(\frak q^+).
\end{equation}
\end{defn}
\begin{rem}
When we define $E_c(\frak q^+)$, we use the additional marked points
$\vec w'_c$ and $\vec w_{\frak p_c}$ that are assigned to $\frak p_c$.
So this subspace is taken in a way independent of $\frak p$. This is important
to prove that the coordinate change satisfies the cocycle condition later.
We explained this point in \cite[the last three lines in the answer to question 4]{Fu1}.
\end{rem}
The next lemma is a consequence of Lemmas \ref{ippaiequivalnce} and \ref{ptequiv}.
\begin{lem}
Suppose that $(\frak Y^{(1)},u^{\prime (1)},(\vec w^{\prime (1)}_c ; c\in \frak B))$ is
weakly equivalent to $(\frak Y^{(2)},u^{\prime (2)},(\vec w^{\prime (2)}_c ; c\in \frak B))$
and $v$ is as in Lemma \ref{ippaiequivalnce}.
We put $\frak q^{+ (j)} = (\frak Y^{(j)},u^{\prime (j)},(\vec w^{\prime (j)}_c ; c\in \frak B))$.
Then
$$
E_c(\frak q^{+ (2)}) = v_*E_c(\frak q^{+ (1)}).
$$
\end{lem}
Now we define:
\begin{defn}\label{defthickened}
The {\it thickened moduli space}
$\mathcal M_{k+1,(\ell,\ell_{\frak p},(\ell_c))}(\beta;\frak p;\frak A;\frak B)_{\epsilon_{0},\vec T_{0}}$
is the subset of $\frak U_{k+1,(\ell,\ell_{\frak p},(\ell_c))}(\beta;\frak p;\frak B)_{\epsilon_{0},\vec T_{0}}$
consisting of the equivalence classes of elements
$\frak q^+ = (\frak Y,u',(\vec w'_c ; c\in \frak B)) \in \frak U_{k+1,(\ell,\ell_{\frak p},(\ell_c))}(\beta;\frak p;\frak B)_{\epsilon_{0},\vec T_{0}}$
that satisfy
\begin{equation}\label{mainequationformulamod}
\overline\partial u' \equiv 0 \mod \mathcal E_{\frak A}(\frak q^+).
\end{equation}
In case $\frak A = \frak B$ we write $\mathcal M_{k+1,(\ell,\ell_{\frak p},(\ell_c))}(\beta;\frak p;\frak A)_{\epsilon_{0},\vec T_{0}}$.

\end{defn}
\begin{lem}\label{nbdregmaineq}
Assume $\frak A \ne \emptyset$.
We can choose $\epsilon_{0}$, $\epsilon_{\frak p_c}$
sufficiently small and $\vec T_{0}$ sufficiently large such that the following holds after
extending the core of $\frak p \cup \vec w_{\frak p}$.
\begin{enumerate}
\item
If $\frak q^+ = (\frak Y,u',(\vec w'_c ; c \in \frak B))$ is in  $\frak U_{k+1,(\ell,\ell_{\frak p},(\ell_c))}(\beta;\frak p;\frak B)_{\epsilon_{0},\vec T_{0}}$
then the equation (\ref{mainequationformulamod}) is Fredholm regular.
\item
If $\frak q^+ = (\frak Y,u',(\vec w'_c ; c\in \frak B))$ is in  $\frak U_{k+1,(\ell,\ell_{\frak p},(\ell_c))}(\beta;\frak p;\frak B)_{\epsilon_{0},\vec T_{0}}$
then $\frak q^+$  is evaluation map transversal.
\item
$\frak p \in \frak U_{k+1,(\ell,\ell_{\frak p},(\ell_c))}(\beta;\frak p;\frak B)_{\epsilon_{0},\vec T_{0}}$.
\end{enumerate}
\end{lem}
Here the definition of Fredholm regularity is the same as Definition \ref{fredreg} and the
definition of evaluation map transversality is the same as Definition \ref{maptrans}.
The proof of Lemma \ref{nbdregmaineq} is the same as that of Proposition \ref{linearMV}.
\begin{rem}\label{Lem258}
More precisely we first choose $\epsilon_{\frak p_c}$ so that Lemma \ref{nbdregmaineq} holds
for $\frak q^+ =\frak p \cup \vec w_{\frak p}$.
(The choice of $\epsilon_{\frak p_c}$ is done at the stage when
we take $\frak M^+(\frak p_c)$ in Definition \ref{transconst}.) Then we take $\epsilon_0$ small so that the Lemma \ref{nbdregmaineq} holds
for any element $\frak q^+$ of  $\frak U_{k+1,(\ell,\ell_{\frak p},(\ell_c))}(\beta;\frak p;\frak B)_{\epsilon_{0},\vec T_{0}}$.
\end{rem}
\par
\begin{cor}\label{smoothness00}
If $\epsilon_{0}$, $\epsilon_{\frak p_c}$ small then
$\mathcal M_{k+1,(\ell,\ell_{\frak p},(\ell_c))}(\beta;\frak p;\frak A;\frak B)_{\epsilon_{0},\vec T_{0}}$
has a structure of smooth manifold stratawise. The dimension of the top stratum is
$$
\dim \mathcal M_{k+1,\ell}(\beta) + 2\sum_{c\in \frak B}\ell_c + 2\ell_{\frak p}
+ \sum_{c \in \frak A}\dim_{\R} E_c.
$$
Here $\dim \mathcal M_{k+1,\ell}(\beta)$ is a virtual dimension that is given by
$$
\dim \mathcal M_{k+1,\ell}(\beta) = k +1 + 2\ell -3 + 2\mu(\beta).
$$
($\mu(\beta)$ is the Maslov index.)
The dimension of the stratum  $\mathcal M_{k+1,(\ell,\ell_{\frak p},(\ell_c))}(\beta;\frak p;\frak A;\frak B;\mathcal G)_{\epsilon_{0},\vec T_{0}}$
is
$$
\dim \mathcal M_{k+1,\ell}(\beta) + 2\sum_{c\in \frak B}\ell_c + 2\ell_{\frak p}
+ \sum_{c \in \frak A}\dim_{\R} E_c - 2\#C^1_{\rm c}(\mathcal G)
- \#C^1_{\rm o}(\mathcal G).
$$
$\Gamma_{\frak p}$ acts effectively on $\mathcal M_{k+1,(\ell,\ell_{\frak p},(\ell_c))}(\beta;\frak p;\frak A;\frak B)_{\epsilon_{0},\vec T_{0}}$.
\end{cor}
Corollary \ref{smoothness00} is an immediate consequence of Lemma  \ref{nbdregmaineq}, implicit function
theorem and index calculation.
\begin{rem}
We can define the topology of $\mathcal M_{k+1,(\ell,\ell_{\frak p},(\ell_c))}(\beta;\frak p;\frak A;\frak B)_{\epsilon_{0},\vec T_{0}}$
in the same way as the topology of $\mathcal M_{k+1,\ell}(\beta)$.
We omit it here and will define the topology of $\mathcal M_{k+1,(\ell,\ell_{\frak p},(\ell_c))}(\beta;\frak p;\frak A;\frak B)_{\epsilon_{0},\vec T_{0}}$
in the next subsection. (Definition \ref{topthickmoduli}.)
\end{rem}
So far we have described the case of $\mathcal M_{k+1,\ell}(\beta)$.  The case of $\mathcal M^{\rm cl}_{\ell}(\alpha)$
is similar with obvious modification.
\par\medskip
\section{Gluing analysis in the general case}
\label{glueing}

The purpose of this section is to generalize Theorems \ref{gluethm1} and \ref{exdecayT}
to the case of the thickened moduli space  $\mathcal M_{k+1,(\ell,\ell_{\frak p},(\ell_c))}(\beta;\frak p;\frak A)_{\epsilon_{0},\vec T_{0}}$
we defined in the last subsection.
Actually this generalization is straightforward.
\par
We first state the result.
Let $\mathcal G_{\frak p}$ be the combinatorial type of $\frak p$.
We first consider the stratum
$\mathcal M_{k+1,(\ell,\ell_{\frak p},(\ell_c))}(\beta;\frak p;\frak A;\mathcal G_{\frak p})_{\epsilon_{0}}$.
We did not include $\vec T_{0}$ in the notation since this parameter does not play a role in our stratum.
(Note $T_{{\rm e},0}$ is the gluing parameter. We do not perform gluing to obtain an element in the
same stratum as $\frak p$.)
We write
\begin{equation}\label{2193}
V_{k+1,(\ell,\ell_{\frak p},(\ell_c))}(\beta;\frak p;\frak A;\frak B;\epsilon_0) = \mathcal M_{k+1,(\ell,\ell_{\frak p},(\ell_c))}(\beta;\frak p;\frak A;\frak B;\mathcal G_{\frak p})_{\epsilon_{0}}.
\end{equation}
This space in this subsection plays the role of $V_1 \times_L V_2$ in Theorem \ref{gluethm1}.
In case $\frak B = \frak A$, we put
$$
V_{k+1,(\ell,\ell_{\frak p},(\ell_c))}(\beta;\frak p;\frak A;\epsilon_0) :=
V_{k+1,(\ell,\ell_{\frak p},(\ell_c))}(\beta;\frak p;\frak A;\frak A;\epsilon_0).
$$
\begin{lem}\label{lem:191}
$V_{k+1,(\ell,\ell_{\frak p},(\ell_c))}(\beta;\frak p;\frak A;\frak B;\epsilon_0)$ has a  structure of smooth manifold.
\end{lem}
\begin{proof}
This is a special case of Corollary \ref{smoothness00} and
is a consequence of Lemma \ref{transpermutelem} (2) and (3).
We give a proof for completeness.
\par
Let $c \in \frak C({\frak p})$.
Since $\frak p \prec \frak p_c$, there exists a map
$\pi : \mathcal G_{\frak p_c} \to \mathcal G_{\frak p}$.
For each ${\rm v}' \in C^0(\mathcal G_{\frak p_c})$
we obtain an element
$\frak p_{c,{\rm v}'} \in \mathcal M_{k_{{\rm v}'+1},\ell_{\rm v}'}(\beta_{{\rm v}'})$
and
$\frak p_{c,{\rm v}'} \cup \vec w_{c,{\rm v}'}\in \mathcal M_{k_{{\rm v}'}+1,\ell_{\rm v}+
\ell_{c,{\rm v}'}}(\beta_{{\rm v}'})$.
For ${\rm v} \in C^0_{\rm d}(\mathcal G_{\frak p})$
the union of $\frak p_{c,{\rm v}'}$ for all ${\rm v}' $ with
$\pi({\rm v}') = {\rm v}$ is an element
$\frak p_{c,{\rm v}} \in \mathcal M_{k_{{\rm v}+1},\ell_{\rm v}}(\beta_{{\rm v}})$.
Together with the union of $\vec w_{c,{\rm v}'}$'s
it gives
$\frak p_{c,{\rm v}} \cup \vec w_{c,{\rm v}}\in \mathcal M_{k_{{\rm v}+1},\ell_{\rm v}+
\ell_{c,{\rm v}}}(\beta_{{\rm v}})$.
The obstruction bundle data centered at $\frak p_c$
induces one centered at $\frak p_{c,{\rm v}}$ in an obvious way.
\par
Let $\frak p_{\rm v} \in \mathcal M_{k_{{\rm v}+1,\ell_{\rm v}}}(\beta_{\rm v})$
be an element obtained by restricting various data of $\frak p$ to the irreducible component  of
$\frak x_{\frak p}$ corresponding to the
vertex ${\rm v}$ in an obvious way.
We have additional marked points $\vec w_c^{\frak p_{\rm v}}$
by restricting $\vec w_c^{\frak p}$.
Then
$\frak p_{\rm v} \cup\vec w_c^{\frak p_{\rm v}}$ is
$\epsilon_c$ close to $\frak p_{c,{\rm v}} \cup \vec w_{c,{\rm v}}$.
\par
We have taken the additional marked points $\vec w_{\frak p}$ on $\frak p$.
Let $\vec w_{\frak p,{\rm v}}$ be a part of it that lies on
the irreducible component $\frak p_{\rm v}$
Then $\frak p_{\rm v} \cup \vec w_{\frak p,{\rm v}}
\in  \mathcal M_{k_{{\rm v}+1},\ell_{\rm v}+\ell_{\frak p,{\rm v}}}(\beta_{{\rm v}})$.
\par
Using $\frak p_{c,{\rm v}}, \vec w_{c,{\rm v}}, \frak p_{\rm v},
\vec w_{\frak p,{\rm v}}, \vec w_{c}^{\frak p_{\rm v}}$ etc., we define
$
\mathcal M_{k_{\rm v}+1,(\ell_{\rm v},\ell_{\frak p,{\rm v}},(\ell_{c,{\rm v}}))}
(\beta_{\rm v};\frak p_{\rm v};\frak A;\frak B;\text{\rm point})_{\epsilon_{0}}.
$
(Note that $\frak p_{\rm v}$ is irreducible.
So the corresponding graph is
trivial, that is the graph without edge.)
We note again that $\frak p_{\rm v}$ is irreducible and is source stable.
So the thickened moduli space
$
\mathcal M_{k_{\rm v}+1,(\ell_{\rm v},\ell_{\frak p,{\rm v}},(\ell_{c,{\rm v}}))}
(\beta_{\rm v};\frak p_{\rm v};\frak A;\frak B;\text{\rm point})_{\epsilon_{0}}
$
is the set parametrized by the solutions of the equations
$$
\overline\partial u' \equiv 0 \mod \mathcal E_{\frak B}(u')
$$
together with the complex structure of the source.
By
Lemma \ref{transpermutelem} (2) the linearized operator of this equation
is surjective. Therefore
$
\mathcal M_{k_{\rm v}+1,(\ell_{\rm v},\ell_{\frak p,{\rm v}},(\ell_{c,{\rm v}}))}
(\beta_{\rm v};\frak p_{\rm v};\frak A;\frak B;\text{\rm point})_{\epsilon_{0}}
$ is a smooth manifold
on a neighborhood of $(\frak p_{\rm v},\vec w_{\frak p,{\rm v}},(\vec w_c^{\frak p_{\rm v}}))$ for each ${\rm v} \in C^0_{\rm d}(\mathcal G_{\frak p})$.
(Note that we add marked points so that there is no
automorphism of elements of $\mathcal M_{k_{\rm v}+1,\ell_{\rm v}}(\beta_{\rm v})$.
So it is not only an orbifold but is also a manifold.)
The case $\rm v \in C_{\rm s}^0(\mathcal G_{\frak p})$ can be discussed in the same way
and obtain
$\mathcal M^{\rm cl}_{(\ell_{\rm v},\ell_{\frak p,{\rm v}},(\ell_{c,{\rm v}}))}
(\beta_{\rm v};\frak p_{\rm v};\frak A;\frak B;\text{\rm point})_{\epsilon_{0}}
$, that is also a smooth manifold.
\par
We take the product of them
for all $\rm v \in C^0(\mathcal G_{\frak p})$.
By taking evaluation maps we have
$$
\aligned
&\prod_{\rm v \in C_{\rm d}^0(\mathcal G_{\frak p})}\mathcal M_{k_{\rm v}+1,(\ell_{\rm v},\ell_{\frak p,{\rm v}},(\ell_{c,{\rm v}}))}
(\beta_{\rm v};\frak p_{\rm v};\frak A;\frak B;\text{\rm point})_{\epsilon_{0}}
\\
&\times
\prod_{\rm v \in C_{\rm s}^0(\mathcal G_{\frak p})}
\mathcal M^{\rm cl}_{(\ell_{\rm v},\ell_{\frak p,{\rm v}},(\ell_{c,{\rm v}}))}
(\beta_{\rm v};\frak p_{\rm v};\frak A;\frak B;\text{\rm point})_{\epsilon_{0}}
\\
&\to
\left(\prod_{\rm e\in  C_{\rm o}^1(\mathcal G_{\frak p})} L
\times
\prod_{\rm e\in  C_{\rm c}^1(\mathcal G_{\frak p})} X\right)^2.
\endaligned
$$
Lemma \ref{transpermutelem} (3) implies that this map is transversal to the diagonal set
$\prod_{\rm e\in  C_{\rm o}^1(\mathcal G_{\frak p})} L
\times
\prod_{\rm e\in  C_{\rm c}^1(\mathcal G_{\frak p})} X
= L^{\# C_{\rm o}^1(\mathcal G_{\frak p})}\times X^{\#C_{\rm c}^1(\mathcal G_{\frak p})}$.
The inverse image of the diagonal set is $V_{k+1,(\ell,\ell_{\frak p},(\ell_c))}(\beta;\frak p;\frak A;\epsilon_0)$.
\end{proof}
\par
The gluing we will perform below defines a map
\begin{equation}
\aligned
\text{\rm Glu} : V_{k+1,(\ell,\ell_{\frak p},(\ell_c))}(\beta;\frak p;\frak A;\frak B;\epsilon_1)
&\times (\vec T^{\rm o}_0,\infty] \times ((\vec T^{\rm c}_0,\infty] \times \vec S^1)
\\
&\to
\mathcal M_{k+1,(\ell,\ell_{\frak p},(\ell_c))}(\beta;\frak p;\frak A;\frak B)_{\epsilon_{0},\vec T_{0}}.
\endaligned
\end{equation}\label{219433}
For a fixed $(\vec T,\vec \theta)$ we denote the restriction of $\text{\rm Glu}$
to $V_{k+1,(\ell,\ell_{\frak p},(\ell_c))}(\beta;\frak p;\frak A;\frak B;\epsilon_1) \times \{(\vec T,\vec \theta)\}$
by $\text{\rm Glu}_{(\vec T,\vec \theta)}$.
\begin{defn}
$\mathcal M_{k+1,(\ell,\ell_{\frak p},(\ell_c))}(\beta;\frak p;\frak A;\frak B;(\vec T,\vec \theta))_{\epsilon_{0},\vec T_{0}}$
is a subset of the space $\mathcal M_{k+1,(\ell,\ell_{\frak p},(\ell_c))}(\beta;\frak p;\frak A;\frak B)_{\epsilon_{0},\vec T_{0}}$
consisting of the equivalence classes of $(\frak Y,u')$ such that
$\frak Y ={\overline\Phi}(\frak y,\vec T,\vec \theta)$ where the combinatorial type of $\frak y$ is $\mathcal G_{\frak p}$.
In case $\frak A = \frak B$, we put
$$
\mathcal M_{k+1,(\ell,\ell_{\frak p},(\ell_c))}(\beta;\frak p;\frak A)_{\epsilon_{0},\vec T_{0}}
=
\mathcal M_{k+1,(\ell,\ell_{\frak p},(\ell_c))}(\beta;\frak p;\frak A;\frak A)_{\epsilon_{0},\vec T_{0}}.
$$
\end{defn}
\begin{thm}\label{gluethm3}
For each sufficiently small $\epsilon_3$,
and sufficiently large $\vec T$, there exist $\epsilon_2$, $\epsilon_4$ and a $\Gamma_{\frak p}^+$ equivariant map
\begin{equation}
\aligned
\text{\rm Glu}_{(\vec T,\vec \theta)} :
&V_{k+1,(\ell,\ell_{\frak p},(\ell_c))}(\beta;\frak p;\frak A;\frak B;\epsilon_4)\\
&\to
\mathcal M_{k+1,(\ell,\ell_{\frak p},(\ell_c))}(\beta;\frak p;\frak A;\frak B;(\vec T,\vec \theta))_{\epsilon_{2}}
\endaligned
\end{equation}
which is a diffeomorphism onto its image.
The image of  $\text{\rm Glu}_{(\vec T,\vec \theta)} $ contains the space $\mathcal M_{k+1,(\ell,\ell_{\frak p},(\ell_c))}(\beta;\frak p;\frak A;\frak B;(\vec T,\vec \theta))_{\epsilon_{3}}$.
\end{thm}
Here $\vec T$ being sufficiently large means that each of its component is sufficiently large.
Theorem \ref{gluethm3} is a generalization of Theorem \ref{gluethm1}.
\begin{defn}\label{topthickmoduli}
We define a topology on $\mathcal M_{k+1,(\ell,\ell_{\frak p},(\ell_c))}(\beta;\frak p;\frak A;\frak B;(\vec T,\vec \theta))_{\epsilon}$ for $\epsilon < \epsilon_3$ and $\vec T_0$ large so that
$\text{\rm Glu}$ is a homeomorphism to the image.
\par
It is easy to see that this topology coincides with the topology that is defined in the same way as the topology of $\mathcal M_{k+1,\ell}(\beta)$.
\end{defn}
To state a generalization of Theorem \ref{exdecayT}, that is the exponential decay estimate of $T$ derivatives,
we take $\vec R$ and the extended core $K_{\rm v}^{+\vec R}$ as in (\ref{extendedcore}).
By restriction we define a map
\begin{equation}\label{2195form}
\aligned
\mathcal M_{k+1,(\ell,\ell_{\frak p},(\ell_c))}&(\beta;\frak p;\frak A;\frak B;(\vec T,\vec \theta))_{\epsilon_{2},\vec T_{0}}
\\
&\to C^{\infty}((K_{\rm v}^{+\vec R},K_{\rm v}^{+\vec R}\cap \partial \Sigma_{\frak p,{\rm v}}),(X,L)).
\endaligned
\end{equation}
We compose it with $\text{\rm Glu}_{(\vec T,\vec \theta)}$ and obtain
$\text{\rm Glures}_{(\vec T^,\vec \theta),{\rm v},\vec R}$.
\begin{thm}\label{exdecayT33}
For each $m$ and $\vec R$ there exist $T(m)$, $C_{6,m,\vec R}$
and $\delta$ such that the following holds
for $T^{\rm o}_{\rm e}>T(m)$, $T^{\rm c}_{\rm e}>T(m)$ and $n + \vert{\vec k_{T}}\vert + \vert{\vec k_{\theta}}\vert \le m - 10$ and
$\vert{\vec k_{T}}\vert + \vert{\vec k_{\theta}}\vert > 0$.
\begin{equation}\label{est196}
\left\Vert \nabla_{\rho}^n \frac{\partial^{\vert \vec k_{T}\vert}}{\partial T^{\vec k_{T}}}\frac{\partial^{\vert \vec k_{\theta}\vert}}{\partial \theta^{\vec k_{\theta}}} \text{\rm Glures}_{(\vec T,\vec \theta), {\rm v},\vec R}\right\Vert_{L^2_{m+1-\vert{\vec k_{T}}\vert - \vert{\vec k_{\theta}}\vert}}
< C_{6,m,\vec R}e^{-\delta' (\vec k_{T}\cdot \vec T+\vec k_{\theta}\cdot \vec T^{\rm c})}.
\end{equation}
\par
Here $\nabla_{\rho}^n$ is the $n$-th derivative in $\rho \in V_{k+1,(\ell,\ell_{\frak p},(\ell_c))}(\beta;\frak p;\frak A;\frak B;\epsilon_2)$ direction
and $\delta'>0$ depends only on $\delta$ and $m$.
\end{thm}
The proofs of Theorems \ref{gluethm3} and \ref{exdecayT33} occupy the rest of this subsection.
We begin with introducing some notations.
Suppose that $(\frak x^{\rho,+},u^{\rho},(\vec w^{\rho}_{c}))$ is a representative of an element $\rho$ of  $V_{k+1,(\ell,\ell_{\frak p},(\ell_c))}(\beta;\frak p;\frak A;\epsilon_0)$.
We put $\Sigma_{\frak x^{\rho,+}} = \Sigma^{\rho}$.
Its marked points are denoted by $\vec z^{\rho}$, $\vec z^{{\rm int},\rho}$ and $\vec w^{\rho}_{\frak p}$,
$\vec w^{\rho}_{c}$.
Here $w$'s are additional marked points.
We divide each of the irreducible components $\Sigma^{\rho}_{\rm v}$ of $\Sigma^{\rho}$ as
\begin{equation}\label{sigmayv}
\aligned
K^{\rho}_{\rm v} &\cup \bigcup_{{\rm e}\in C^1_{\mathrm o}(\mathcal G)
\atop \text{${\rm e}$ is an outgoing edge of ${\rm v}$}} (0,\infty) \times [0,1]
\\
&\cup
\bigcup_{{\rm e}\in C^1_{\mathrm o}(\mathcal G)
\atop \text{${\rm e}$ is an incoming edge of ${\rm v}$}} (-\infty,0) \times [0,1]
\\
&\cup
\bigcup_{{\rm e}\in C^1_{\mathrm c}(\mathcal G)
\atop \text{${\rm e}$ is an outgoing edge of ${\rm v}$}} (0,\infty) \times S^1
\\
&\cup
\bigcup_{{\rm e}\in C^1_{\mathrm c}(\mathcal G)
\atop \text{${\rm e}$ is an incoming edge of ${\rm v}$}} (-\infty,0) \times S^1,
\endaligned
\end{equation}
where the coordinates of the 2-nd, 3-rd, 4-th, and 5-th summands are
$(\tau'_{\rm e},t_{\rm e})$, $(\tau''_{\rm e},t_{\rm e})$, $(\tau'_{\rm e},t_{\rm e}')$,
and $(\tau''_{\rm e},t_{\rm e}'')$, respectively.
Here $\tau'_{\rm e} \in (0,\infty)$, $\tau''_{\rm e} \in (-\infty,0)$.
\par
We call the end corresponding to $\rm e$ the {\it $\rm e$-th end}.
\par
We recall
\begin{eqnarray}\label{neckcoordinate2}
\tau_e &=& \tau'_{\rm e} - 5T_{\rm e}  = \tau''_{\rm e} + 5T_{\rm e}, \label{cctau12}\\
t_{\rm e} &=& t'_{\rm e} = t''_{\rm e} - \theta_{\rm e}.\label{ccttt12}
\end{eqnarray}
We put
$$
u^{\rho}_{\rm v} = u^{\rho}\vert_{K_{\rm v}},
\qquad
u^{\rho}_{\rm e} = u^{\rho}\vert_{\text{$\rm e$-th neck region}}.
$$
We denote by $\Sigma_{\frak Y} = \Sigma_{\vec T,\vec \theta}^{\rho}$ a representative of $\frak Y = \overline{\Phi}(\frak y,\vec T,\vec \theta)$.
The curve $\Sigma_{\vec T,\vec \theta}^{\rho}$ is a union
\begin{equation}\label{neddecomposit}
\aligned
\bigcup_{\rm v \in C^0(\mathcal G_{\frak p})}K^{\rho}_{\rm v}
 &\cup \bigcup_{{\rm e}\in C^1_{\mathrm o}(\mathcal G)} [-5T_{\rm e},5T_{\rm e}] \times [0,1] \\
 &\cup \bigcup_{{\rm e}\in C^1_{\mathrm c}(\mathcal G)} [-5T_{\rm e},5T_{\rm e}] \times S^1.
\endaligned
\end{equation}
The coordinates of the 2nd and 3rd terms are $\tau_{\rm e}$ and $t_{\rm e}$.
\par
We call $[-5T_{\rm e},5T_{\rm e}] \times [0,1]$ or $[-5T_{\rm e},5T_{\rm e}] \times S^1$ the
{\it $\rm e$-th neck}.
\par
In case $T_{\rm e} = \infty$, the curve $\Sigma_{\vec T,\vec \theta}^{\rho}$ contains $([0,\infty) \cup (-\infty,0]) \times [0,1]$ or
$([0,\infty) \cup (-\infty,0]) \times S^1$ corresponding to the $\rm e$-th edge.
We call   $([0,\infty) \times [0,1]$ (or $\times S^1$) the {\it outgoing $\rm e$-th end}
and  $(-\infty,0] \times [0,1]$ (or $S^1$) the {\it incoming $\rm e$-th end}.
\par
We call $K_{\rm v}$ the {\it $\rm v$-th core}.
\par
The restriction of $u^{\rho}$ to $K_{\rm v}$ is written as $u^{\rho}_{\rm v}$.
The restriction of $u^{\rho}$ to the $\rm e$-neck is written as $u^{\rho}_{\rm e}$.
\par
For each $\rm e$, let $\rm v_1$ and $\rm v_2$ be its incoming and outgoing vertices. We have
\begin{equation}\label{ryogawalimitcoin}
\lim_{\tau_{\rm e}\to -\infty} u_{\rm v_2}^{\rho}(\tau_{\rm e},t_{\rm e})
=
\lim_{\tau_{\rm e}\to \infty} u_{\rm v_1}^{\rho}(\tau_{\rm e},t_{\rm e}),
\end{equation}
and (\ref{ryogawalimitcoin}) is independent of $t_{\rm e}$.
We write this limit as $p^{\rho}_{\rm e}$.
We take a Darboux coordinate in a neighborhood of each $p^{\rho}_{\rm e}$ such that $L$ is flat in this
coordinate.  We choose the map $\rm E$ such that (\ref{Einanbdofp0}) holds in this neighborhood of  $p^{\rho}_{\rm e}$.
\par
For $\rm e \in C^{1}_{\rm o}(\mathcal G_{\frak p})$ with  $T_{\rm e} \ne \infty$, we define
\begin{equation}\label{DemanAetc}
\aligned
\mathcal A_{{\rm e},T} &= [-T_{\rm e}-1,-T_{\rm e}+1] \times [0,1] \subset [-5T_{\rm e},5T_{\rm e}] \times [0,1], \\
\mathcal B_{{\rm e},T} &= [T_{\rm e}-1,T_{\rm e}+1] \times [0,1] \subset [-5T_{\rm e},5T_{\rm e}] \times [0,1], \\
\mathcal X_{{\rm e},T} &= [-1,+1] \times [0,1] \subset [-5T_{\rm e},5T_{\rm e}] \times [0,1].
\endaligned
\end{equation}
In case $\rm e \in C^{1}_{\rm c}(\mathcal G_{\frak p})$, the sets
$\mathcal A_{{\rm e},T}$,  $\mathcal B_{{\rm e},T}$,   $\mathcal X_{{\rm e},T}$ are defined
in the same way as above replacing $[0,1]$
by $S^1$.
\par
If $\rm v$ is a vertex of $\rm e$ then $\mathcal A_{{\rm e},T}$, $\mathcal B_{{\rm e},T}$, $\mathcal X_{{\rm e},T}$
may be regarded as a subset of $\Sigma_{\rm v}^{\rho}$ also.
\par
Let $\chi_{{\rm e},\mathcal A}^{\leftarrow}$,
$\chi_{{\rm e},\mathcal A}^{\rightarrow}$ be smooth functions on $[-5T_{\rm e},5T_{\rm e}]\times [0,1]$
or $[-5T_{\rm e},5T_{\rm e}]\times S^1$ such that
\begin{equation}\label{eq201}
\chi_{{\rm e},\mathcal A}^{\leftarrow}(\tau_{\rm e},t_{\rm e})
= \begin{cases}
1   & \tau_{\rm e} < -T_{\rm e}-1 \\
0  & \tau_{\rm e} > -T_{\rm e}+1.
\end{cases}
\end{equation}
$$
\chi_{{\rm e},\mathcal A}^{\rightarrow}  = 1 - \chi_{{\rm e},\mathcal A}^{\leftarrow}.
$$
We define
\begin{equation}
\chi_{{\rm e},\mathcal B}^{\leftarrow}(\tau_{\rm e},t_{\rm e})
= \begin{cases}
1   & \tau_{\rm e} <  T_{\rm e}-1 \\
0  & \tau_{\rm e} >  T_{\rm e}+1.
\end{cases}
\end{equation}
$$
\chi_{{\rm e},\mathcal B}^{\rightarrow}  = 1 - \chi_{{\rm e},\mathcal B}^{\leftarrow}.
$$
We define
\begin{equation}
\chi_{{\rm e},\mathcal X}^{\leftarrow}(\tau_{\rm e},t_{\rm e})
= \begin{cases}
1   & \tau_{\rm e} <  -1 \\
0  & \tau_{\rm e} >  1.
\end{cases}
\end{equation}
$$
\chi_{{\rm e},\mathcal X}^{\rightarrow}  = 1 - \chi_{{\rm e},\mathcal X}^{\leftarrow}.
$$
We extend these functions to $\Sigma_{\vec T,\vec \theta}^{\rho}$ and $\Sigma_{\rm v}^{\rho}$ so that
they are locally constant on its core.
We denote them by the same symbol.
\par
We next introduce weighted Sobolev norms and their local versions for sections
on $\Sigma_{\rm v}^{\rho}$ as follows.
We define a smooth function $e_{\rm v,\delta} : \Sigma_{\rm v}^{\rho} \to [1,\infty)$
by
\begin{equation}\label{e1deltamulti}
e_{\rm v,\delta}(\tau_{\rm e},t_{\rm e})
\begin{cases}
=1   &\text{on  $K_{\rm v}$,} \\
=e^{\delta\vert \tau_{\rm e} + 5T_{\rm e}\vert} &\text{if $\tau_{\rm e} > 1 -  5T_{\rm e}$, and $\rm e$ is an outgoing edge of $\rm v$,}\\
\in [1,10] &\text{if $\tau_{\rm e} < 1 -  5T_{\rm e}$, and $\rm e$ is an outgoing edge of $\rm v$,}\\
=e^{\delta\vert \tau - 5T_{\rm e}\vert} &\text{if $\tau_{\rm e} <  5T_{\rm e}-1$, and $\rm e$ is an incoming edge of $\rm v$,}\\
\in [1,10] &\text{if $\tau_{\rm e} >  5T_{\rm e}-1$, and $\rm e$ is an incoming edge of $\rm v$.}
\end{cases}
\end{equation}
We also define a weight function $e_{\vec T,\delta} : \Sigma_{\vec T,\vec \theta}^{\rho} \to [1,\infty)$ as follows:
\begin{equation}\label{e2delta}
e_{\vec T,\delta}(\tau_{\rm e},t_{\rm e})
\begin{cases}
=e^{\delta\vert \tau_{\rm e} - 5T_{\rm e}\vert} &\text{if $1<\tau_{\rm e} <  5T_{\rm e}-1$,}\\
= e^{\delta\vert \tau + 5T_{\rm e}\vert} &\text{if $-1>\tau > 1-5T_{\rm e}$,}\\
=1   &\text{on  $K_{\rm v}$,} \\
\in [1,10] &\text{if $\vert\tau_{\rm e} - 5T_{\rm e}\vert < 1$ or $\vert\tau_{\rm e} + 5T_{\rm e}\vert < 1$,} \\
\in [e^{5T_{\rm e}\delta}/10,e^{5T_{\rm e}\delta}] &\text{if $\vert\tau_{\rm e}\vert < 1$}.
\end{cases}
\end{equation}
The weighted Sobolev norm we use for
$L^2_{m,\delta}(\Sigma^{\rho}_{\rm v};(u_{\rm v}^{\rho})^*TX\otimes \Lambda^{01})$
is given by
\begin{equation}
\Vert s\Vert^2_{L^2_{m,\delta}} = \sum_{k=0}^m \int_{\Sigma_{\rm v}^{\rho}}
e_{\rm v,\delta} \vert \nabla^k s\vert^2 \text{\rm vol}_{\Sigma_{\rm v}^{\rho}}.
\end{equation}
\begin{defn}\label{Sobolev263}
The Sobolev space
$L^2_{m+1,\delta}((\Sigma_{\rm v}^{\rho},\partial \Sigma_{\rm v}^{\rho});(u_{\rm v}^{\rho})^*TX,(u_{\rm v}^{\rho})^*TL)$
consists of elements $(s,\vec v)$ with the following properties.
\begin{enumerate}
\item $\vec v = (v_{\rm e})$ where $\rm e$ runs on the set of edges of $\rm v$  and
$v_{\rm e} \in T_{p^{\rho}_{\rm e}}(X)$ (in case $\rm e\in C^1_{\rm c}(\mathcal G)$) or $v_{\rm e} \in T_{p^{\rho}_{\rm e}}(L)$
(in case $\rm e\in C^1_{\rm o}(\mathcal G)$).
\item The following norm is finite.
\begin{equation}\label{normformjulamulti}
\aligned
\Vert (s,\vec v)\Vert^2_{L^2_{m+1,\delta}} = &\sum_{k=0}^{m+1} \int_{K_{\rm v}}
\vert \nabla^k s\vert^2 \text{\rm vol}_{\Sigma_i}+ \sum_{{\rm e} \text{: edges of $\rm v$}}\Vert v_{\rm e}\Vert^2\\
&+
\sum_{k=0}^{m+1}\sum_{{\rm e} \text{: edges of $\rm v$}} \int_{\text{e-th end}}  e_{\rm v,\delta}\vert \nabla^k(s - \text{\rm Pal}(v_{\rm e}))\vert^2 \text{\rm vol}_{\Sigma_{\rm v}^{\rho}}.
\endaligned
\end{equation}
\end{enumerate}
\end{defn}
\begin{defn}
We define
\begin{equation}\label{evaluationmkepot}
\aligned
D{\rm ev}_{\mathcal G_{\frak p}} :
&\bigoplus_{{\rm v}\in C^0_{\rm o}(\mathcal G_{\frak p})}
L^2_{m+1,\delta}((\Sigma_{{\rm v}}^{\rho},\partial \Sigma_{\rm v}^{\rho});(u_{\rm v}^{\rho})^*TX,(u_{\rm v}^{\rho})^*TL)\\
&\oplus \bigoplus_{{\rm v}\in C^0_{\rm c}(\mathcal G_{\frak p})}
L^2_{m+1,\delta}(\Sigma_{{\rm v}}^{\rho};(u_{\rm v}^{\rho})^*TX)
\\
&\to \bigoplus_{\rm e \in C^1_{\rm o}(\mathcal G_{\frak p})}T_{p^{\rho}_{\rm e}}L
\oplus \bigoplus_{\rm e \in C^1_{\rm c}(\mathcal G_{\frak p})}T_{p^{\rho}_{\rm e}}X
\endaligned\end{equation}
as in (\ref{evaluationmap2diff}).
\end{defn}
\begin{defn}\label{Lhattaato1}
We denote the kernel of (\ref{evaluationmkepot}) by
$$
L^2_{m+1,\delta}((\Sigma^{\rho},\partial \Sigma^{\rho});(u^{\rho})^*TX,(u^{\rho})^*TL).
$$
\end{defn}
\par
We next define  weighted Sobolev norms for the sections of various bundles on
$\Sigma_{\vec T,\vec \theta}^{\rho}$.
Let
$$
u' : (\Sigma_{\vec T,\vec \theta}^{\rho},\partial \Sigma_{\vec T,\vec \theta}^{\rho}) \to (X,L)
$$
be a smooth map of homology class $\beta$ that is pseudo-holomorphic in the neck region and has finite energy.
(We include the case when $u'$ is not pseudo-holomorphic in the neck region but satisfies the same
exponential decay estimate as the pseudo-holomorphic curve.)
We first consider the case when all $T_{\rm e} \ne \infty$.
In this case $\Sigma_{\vec T,\vec \theta}^{\rho}$ is compact.
We consider an element
$$
s \in L^2_{m+1}((\Sigma_{\vec T,\vec \theta}^{\rho},\partial \Sigma_{\vec T,\vec \theta}^{\rho});(u')^*TX,(u')^*TL).
$$
Since we take $m$ large, the section $s$ is continuous. We take a point $(0,1/2)_{\rm e}$ in the $\rm e$-th neck.
Since $s \in L^2_{m+1}$ its value $s((0,1/2)_{\rm e}) \in T_{u'((0,1/2)_{\rm e})}X$ is well-defined.
\par
We take a coordinate around $p_{\rm e}^{\rho}$ such that in case $\rm e \in C^1_{\rm o}(\mathcal G)$
our Lagrangian submanifold
$L$ is linear in this coordinate around $p_{\rm e}^{\rho}$.
We use this trivialization to find a canonical trivialization of $TX$ in a neighborhood of $p_{\rm e}^{\rho}$. We use
this trivialization to define $\text{\rm Pal}$ below.
We put
\begin{equation}\label{normformjula5multi}
\aligned
\Vert s\Vert^2_{L^2_{m+1,\delta}} = &\sum_{k=0}^{m+1}\sum_{\rm v} \int_{K_{\rm v}}
\vert \nabla^k s\vert^2 \text{\rm vol}_{\Sigma_{\rm v}^{\rho}}\\
&+
\sum_{k=0}^{m+1}\sum_{\rm e} \int_{\text{e-th neck}}  e_{\vec T,\delta}\vert \nabla^k(s - \text{\rm Pal}(s(0,1/2)_{\rm e}))\vert^2
dt_{\rm e}d\tau_{\rm e}
\\&+ \sum_{\rm e}\Vert s((0,1/2)_{\rm e}))\Vert^2.
\endaligned
\end{equation}
For a section
$
s \in L^2_{m}(\Sigma_{\vec T,\vec \theta}^{\rho};u^*TX\otimes \Lambda^{01})
$
we define
\begin{equation}\label{normformjula52multi}
\Vert s\Vert^2_{L^2_{m,\delta}} = \sum_{k=0}^{m} \int_{\Sigma_{\vec T,\vec \theta}^{\rho}}
e_{T,\delta}\vert  \nabla^k s\vert^2 \text{\rm vol}_{\Sigma_{\vec T,\vec \theta}^{\rho}}.
\end{equation}
\par
We next consider the case when some of the edges $\rm e$ have infinite length, namely $T_{\rm e} = \infty$.
Let $C^{1,\rm{inf}}_{\rm o}(\mathcal G_{\frak p},\vec T)$ (resp. $C^{1,\rm{inf}}_{\rm c}(\mathcal G_{\frak p},\vec T)$) be the set of elements $\rm e$ in $C^{1}_{\rm o}(\mathcal G_{\frak p})$ (resp. $C^{1}_{\rm c}(\mathcal G_{\frak p})$)
with $T_{\rm e} = \infty$ and let
$C^{1,\rm{fin}}_{\rm o}(\mathcal G_{\frak p},\vec T)$ (resp. $C^{1,\rm{fin}}_{\rm c}(\mathcal G_{\frak p},\vec T)$) be the set of elements
${\rm e} \in C^{1}_{\rm o}(\mathcal G_{\frak p})$ (resp. $C^{1}_{\rm c}(\mathcal G_{\frak p})$)
with $T_{\rm e} \ne \infty$.
Note the ends of $\Sigma_{\vec T,\vec \theta}^{\rho}$ correspond two to one to
$C^{1,\rm{inf}}_{\rm o}(\mathcal G_{\frak p},\vec T) \cup C^{1,\rm{inf}}_{\rm c}(\mathcal G_{\frak p},\vec T)$.
The ends that correspond to an element $\rm e$ of $C^{1,\rm{inf}}_{\rm o}(\mathcal G_{\frak p},\vec T)$ is
$([-5T_{\rm e},\infty) \times [0,1]) \cup (-\infty,5T_{\rm e}] \times [0,1])$ and
the ends that correspond to ${\rm e} \in C^{1,\rm{inf}}_{\rm c}(\mathcal G_{\frak p},\vec T)$ is
$([-5T_{\rm e},\infty) \times S^1) \cup (-\infty,5T_{\rm e}] \times S^1)$.
We have a weight function $e_{\rm v,\delta}(\tau_{\rm e},t_{\rm e})$ on it.
\begin{defn}\label{Lhattaato2}
An element of
$$
L^2_{m+1,\delta}((\Sigma_{\vec T,\vec \theta}^{\rho},\partial \Sigma_{\vec T,\vec \theta}^{\rho});(u')^*TX,(u')^*TL)
$$
is a pair $(s,\vec v)$ such that
\begin{enumerate}
\item $s$ is a section of
$(u')^*TX$ on $\Sigma_{\vec T,\vec \theta}^{\rho}$ minus singular points $z_{\rm e}$ with $T_{\rm e} = \infty$.
\item $s$ is locally of $L^2_{m+1}$ class.
\item On $\partial \Sigma_{\vec T,\vec \theta}^{\rho}$ the restriction of $s$ is in $(u')^*TL$.
\item $\vec v = (v_{\rm e})$ where ${\rm e}$ runs in $C^{1,{\rm inf}}({\mathcal G}_{\frak p},{\vec T})$ and
$v_{\rm e}$ is as in Definition \ref{Sobolev263} (1).
\item For each $\rm e$ with $T_{\rm e} = \infty$, the integral
\begin{equation}\label{intatinfedge}
\aligned
&\sum_{k=0}^{m+1}\int_{0}^{\infty} \int_{t_{\rm e}} e_{\rm v,\delta}(\tau_{\rm e},t_{\rm e})\vert\nabla^k(s(\tau_{\rm e},t_{\rm e})
 - {\rm Pal}(v_{\rm e}))\vert^2 d\tau_{\rm e} dt_{\rm e}\\
&+ \sum_{k=0}^{m+1}\int^{0}_{-\infty} \int_{t_{\rm e}} e_{\rm v,\delta}(\tau_{\rm e},t_{\rm e})\vert\nabla^k(s(\tau_{\rm e},t_{\rm e})
 - {\rm Pal}(v_{\rm e}))\vert^2 d\tau_{\rm e} dt_{\rm e}
  \endaligned
\end{equation}
is finite. (Here we integrate over ${t_{\rm e}\in [0,1]}$ (resp. ${t_{\rm e}\in S^1}$) if  ${\rm e }\in C_{\rm o}^{1,{\rm inf}}({\mathcal G}_{\frak p},\vec T)$
(resp.  ${\rm e} \in C_{\rm c}^{1,{\rm inf}}({\mathcal G}_{\frak p},\vec T)$).
\end{enumerate}
We define
\begin{equation}\label{2210}
\Vert (s,\vec v)\Vert^2_{L^2_{m+1,\delta}}
= (\ref{normformjula5multi}) + \sum_{{\rm e} \in C^{1,{\rm inf}}({\mathcal G}_{\frak p},\vec T)}(\ref{intatinfedge}) +
\sum_{{\rm e} \in C^{1,{\rm inf}}(\mathcal G_{\frak p},\vec T)}\Vert v_{{\rm e}}\Vert^2.
\end{equation}
An element of
$$
L^2_{m,\delta}(\Sigma_{\vec T,\vec \theta}^{\rho};(u')^*TX\otimes \Lambda^{01})
$$
is a section $s$ of the bundle $(u')^*TX\otimes \Lambda^{01}$ such that it is locally of  $L^2_{m}$-class and
\begin{equation}\label{intatinfedge3}
\aligned
&\sum_{k=0}^{m}\int_{0}^{\infty} \int_{t_{\rm e}} e_{\rm v,\delta}\vert\nabla^k s(\tau_{\rm e},t_{\rm e})\vert^2
 d\tau_{\rm e} dt_{\rm e}\\
&+ \sum_{k=0}^{m}\int^{0}_{-\infty} \int_{t_{\rm e}} e_{\rm v,\delta}\vert\nabla^k(s(\tau_{\rm e},t_{\rm e})
\vert^2 d\tau_{\rm e} dt_{\rm e}
  \endaligned
\end{equation}
is finite.
We define
\begin{equation}\label{2212}
\Vert s\Vert^2_{L^2_{m,\delta}}
= (\ref{normformjula52multi}) + \sum_{{\rm e} \in C^{1,{\rm inf}}({\mathcal G}_{\frak p},\vec T)}(\ref{intatinfedge3}).
\end{equation}
\end{defn}
\par
For a subset $W$ of $\Sigma_{\rm v}^{\rho}$ or $\Sigma^{\rho}_{\vec T,\vec \theta}$ we define
$
\Vert s\Vert_{L^2_{m,\delta}(W\subset \Sigma^{\rho}_{\rm v})}
$,
$
\Vert s\Vert_{L^2_{m,\delta}(W\subset \Sigma_{\vec T,\vec \theta}^{\rho})}
$
by restricting the domain of the integration (\ref{normformjula52multi}),
(\ref{normformjula5multi}), (\ref{2210}) or (\ref{2212}) to $W$.
\par
Let $(s_j,\vec v_j) \in L^2_{m+1,\delta}((\Sigma_{\rm v}^{\rho},\partial \Sigma_{\rm v}^{\rho});(u_{\rm v}^{\rho})^*TX,(u_{\rm v}^{\rho})^*TL)$
for $j=1,2$. We define an inner product among them by:
\begin{equation}\label{innerprod}
\aligned
&\langle\!\langle (s_1,\vec v_1),(s_2,\vec v_2)\rangle\!\rangle_{L^2_{\delta}}\\
=
& \sum_{{\rm e}\in C^1(\mathcal G_{\frak p})}\int_{\text{{\rm e}-th neck}} e_{\vec T,\delta} (s_1-\text{\rm Pal}(v_{1,{\rm e}}),
s_2-\text{\rm Pal}(v_{2,{\rm e}}))
\\
&+\sum_{{\rm v}\in C^0(\mathcal G_{\frak p})}\int_{K_{\rm v}} (s_1,s_2)
+\sum_{{\rm e}\in C^1(\mathcal G_{\frak p})}
(v_{1,\rm e},v_{2,\rm e}).
\endaligned
\end{equation}
\par
Now we start the gluing process.
Let us start with the maps
$$u_{\rm v}^{\rho} : (\Sigma_{\rm v}^{\rho},\partial \Sigma_{\rm v}^{\rho}) \to (X,L)$$
for each $\rm v$ so that
$(u_{\rm v}^{\rho};{\rm v} \in C^0(\mathcal G_{\frak p}))$ consists an element of $V_{k+1,(\ell,\ell_{\frak p},(\ell_c))}(\beta;\frak p;\frak A;\epsilon_2)$.
Let $(\vec T,\vec\theta) \in (\vec T^{\rm o}_0,\infty] \times ((\vec T^{\rm c}_0,\infty] \times \vec S^1)$.
For $\kappa = 0,1,2,\dots$, we will define a series of maps
\begin{eqnarray}
u_{\vec T,\vec \theta,(\kappa)}^{\rho} &:& (\Sigma_{\vec T,\vec \theta}^{\rho},\partial\Sigma_{\vec T,\vec \theta}^{\rho})
\to (X,L)
\\
\hat u_{{\rm v},\vec T,\vec \theta,(\kappa)}^{\rho} &:& (\Sigma_{\rm v}^{\rho},\partial \Sigma_{\rm v}^{\rho}) \to (X,L)
\end{eqnarray}
and elements
\begin{eqnarray}
\frak e^{\rho} _{c,\vec T,\vec \theta,(\kappa)} &\in& E_c = \bigoplus_{{\rm v}\in C^0(\mathcal G_{\frak p_{c}})} E_{c,\rm v}
\\
{\rm Err}^{\rho}_{{\rm v},\vec T,\vec \theta,(\kappa)}  &\in&
L^2_{m,\delta}(\Sigma_{\rm v}^{\rho};(\hat u_{{\rm v},\vec T,\vec \theta,(\kappa)}^{\rho})^*TX\otimes \Lambda^{01}).
\end{eqnarray}
Note $E_{c,\rm v} \subset \Gamma(K_{\rm v};u_{\frak p_c}^*TX \otimes \Lambda^{01})$ is a finite dimensional
space which we take as a part of the obstruction bundle data centered at $\frak p_c$.
\par
Moreover we will define
$V^{\rho}_{\vec T,\vec \theta,{\rm v},(\kappa)}$ for ${\rm v} \in C^0(\mathcal G_{\frak p})$ and $\Delta p^{\rho}_{{\rm e},\vec T,\vec \theta,(\kappa)}$
for ${\rm e} \in C^1(\mathcal G_{\frak p})$. The pair $((V^{\rho}_{\vec T,\vec \theta,{\rm v},(\kappa)}),(\Delta p^{\rho}_{{\rm e},\vec T,\vec \theta,(\kappa)}))$ is
an element of the weighted Sobolev space
$L^2_{m+1,\delta}((\Sigma^{\rho}_{\rm v},\partial \Sigma^{\rho}_{\rm v});
(\hat u_{{\rm v},\vec T,\vec \theta,(\kappa-1)}^{\rho})^*TX,(\hat u_{{\rm v},\vec T,\vec \theta,(\kappa-1)}^{\rho})^*TL)$.
\par
The construction of these objects is a straightforward generalization of the construction given by Section \ref{alternatingmethod} and proceed
by induction on $\kappa$ as follows.
\par\medskip
\noindent{\bf Pregluing}:
We first define an approximate solution $u_{\vec T,\vec \theta,(0)}^{\rho}$. For ${\rm e} \in C^1(\mathcal G_{\frak p})$ we denote by
${\rm v}_{\leftarrow}({\rm e})$  and ${\rm v}_{\rightarrow}({\rm e})$ its two vertices.
Here $\rm e$ is an outgoing edge of ${\rm v}_{\leftarrow}({\rm e})$ and is an incoming edge of
${\rm v}_{\rightarrow}({\rm e})$. We put:
\begin{equation}\label{22190}
u_{\vec T,\vec \theta,(0)}^{\rho} =
\begin{cases}
\chi_{{\rm e},\mathcal B}^{\leftarrow} (u_{{\rm v}_{\leftarrow}({\rm e})}^{\rho} - p^{\rho}_{\rm e}) +
\chi_{{\rm e},\mathcal A}^{\rightarrow} (u_{{\rm v}_{\rightarrow}({\rm e})}^{\rho} - p^{\rho}_{\rm e}) + p^{\rho}_{\rm e}
& \text{on the ${\rm e}$-th neck} \\
u_{\rm v}^{\rho} & \text{on $K_{\rm v}$} .
\end{cases}
\end{equation}
\par\medskip
\noindent{\bf Step 0-3}:
We next define
\begin{equation}
\sum_{c\in \frak A} \frak e^{\rho} _{c,\vec T,\vec \theta,(0)} = \overline\partial u_{\rm v}^{\rho},
\qquad \text{on $K_{\rm v}$}.
\end{equation}
Here we identify $E_c \cong E_c(u_{\rm v}^{\rho})$ on $K_{\rm v}$
by the parallel transport as we did in Definition \ref{defEc}.
See also Definition \ref{Emovevvv}.
Note that $\overline\partial u_{\rm v}^{\rho}$
is contained in $\oplus E_c$ since $(u_{\rm v}^{\rho};{\rm v} \in C^0(\mathcal G_{\frak p}))$
is an element of
 $V_{k+1,(\ell,\ell_{\frak p},(\ell_c))}(\beta;\frak p;\frak A;\epsilon_0)$.
\par
We put
\begin{equation}
\frak{se}^{\rho}_{\vec T,\vec \theta,(0)} :=  \sum_{c\in \frak A} \frak e^{\rho} _{c,\vec T,\vec \theta,(0)}.
\end{equation}
\par\medskip
\noindent{\bf Step 0-4}:
We next define
\begin{equation}\label{2222}
{\rm Err}^{\rho}_{{\rm v},\vec T,\vec \theta,(0)}  =
\begin{cases}
\chi_{{\rm e},\mathcal X}^{\leftarrow} \overline\partial u_{\vec T,\vec \theta,(0)}^{\rho}
& \text{on the ${\rm e}$-th neck if $\rm e$ is outgoing} \\
\chi_{{\rm e},\mathcal X}^{\rightarrow}  \overline\partial u_{\vec T,\vec \theta,(0)}^{\rho}
& \text{on the ${\rm e}$-th neck if $\rm e$ is incoming} \\
\overline\partial u_{\vec T,\vec \theta,(0)}^{\rho}  - \frak{se}^{\rho}_{\vec T,\vec \theta,(0)} & \text{on $K_{\rm v}$} .
\end{cases}
\end{equation}
See Remark \ref{rem270}.
\par\medskip
\noindent{\bf Step 1-1}:
We put
\begin{equation}
\aligned
&\hat u^{\rho}_{{\rm v},\vec T,\vec \theta,(0)}(z) \\
&=
\begin{cases} \chi_{{\rm e},\mathcal B}^{\leftarrow}(\tau_{\rm e}-T_{\rm e},t_{\rm e})
&\!\!\!\!\!\!u^{\rho}_{\vec T,\vec \theta,(0)}(\tau_{\rm e},t_{\rm e})
+ \chi_{{\rm e},\mathcal B}^{\rightarrow}(\tau_{\rm e}-T_{\rm e},t_{\rm e})p^{\rho}_{\rm e} \\
&\text{if $z = (\tau_{\rm e},t_{\rm e})$ is on the $\rm e$-th neck that is outgoing} \\
\chi_{{\rm e},\mathcal A}^{\rightarrow}(\tau_{\rm e}-T_{\rm e},t_{\rm e})
&\!\!\!\!\!\!u^{\rho}_{\vec T,\vec \theta,(0)}(\tau,t)
+ \chi_{{\rm e},\mathcal A}^{\leftarrow}(\tau_{\rm e}-T_{\rm e},t_{\rm e})p^{\rho}_{\rm e} \\
&\text{if $z = (\tau_{\rm e},t_{\rm e})$ is on the $\rm e$-th neck that is incoming} \\
 u^{\rho}_{{\rm v},\vec T,\vec \theta,(0)}(z)
&\text{if $z \in K_{\rm v}$.}
\end{cases}
\endaligned
\end{equation}
We denote the (covariant) linearization of the Cauchy-Riemann equation at this map $\hat u^{\rho}_{{\rm v},\vec T,\vec \theta,(0)}$  by
\begin{equation}\label{lineeqstep0vvv}
\aligned
D_{\hat u^{\rho}_{{\rm v},\vec T,\vec \theta,(0)}}\overline{\partial}
:
L^2_{m+1,\delta}((\Sigma_{\rm v},\partial \Sigma_{\rm v});&(\hat u^{\rho}_{{\rm v},\vec T,\vec \theta,(0)})^*TX,(\hat u^{\rho}_{{\rm v},\vec T,\vec \theta,(0)})^*TL)
\\
&\to
L^2_{m,\delta}(\Sigma_{\rm v};(\hat u^{\rho}_{{\rm v},\vec T,\vec \theta,(0)})^{*}TX \otimes \Lambda^{01}).
\endaligned
\end{equation}
\par
We next study the obstruction bundle $E_c$.
We recall that at $u^{\rho}_{\vec T,\vec \theta,(0)}$ the obstruction bundle $E_c(u^{\rho}_{\vec T,\vec \theta,(0)})$
was defined as follows. (See Definition \ref{defEc}.)
We use the added marked points $\vec w^{\rho}_c$ and consider $\Sigma_{\vec T,\vec \theta}^{\rho} \cup \vec w^{\rho}_c$.
Here, by abuse of notation, we include the $k+1$ boundary and $\ell$ interior
marked points in the notation $\Sigma_{\vec T,\vec \theta}^{\rho}$ .
(The additional marked points $\vec w_{\frak p}^{\rho}$ and $\vec w_{c}^{\rho}$ are {\it not} included.)
By assumption $\Sigma_{\vec T,\vec \theta}^{\rho} \cup \vec w^{\rho}_c$ is ($\epsilon_c$+$o(\epsilon_0)$)-close to
$\frak p_c$. Therefore the diffeomorphism between cores of $\Sigma_{\frak p_c}$ and of $\Sigma_{\vec T,\vec \theta}^{\rho}$
is determined, by the obstruction bundle data $\frak E_{\frak p_c}$. Using this diffeomorphism and the parallel transport we have
\begin{equation}\label{2224}
I^{{\rm v},\frak p_c}_{(\frak y_c,u_c),(\Sigma_{\vec T,\vec \theta}^{\rho} \cup \vec w^{\rho}_c,u^{\rho}_{\vec T,\vec \theta,(0)})} : E_{c,{\rm v}}(\frak y_c,u_c) \to \Gamma(K_{\rm v} ; (u^{\rho}_{\vec T,\vec \theta,(0)})^{*}TX \otimes \Lambda^{01}).
\end{equation}
The notation in (\ref{2224}) is as follows. There is a map $\pi : \mathcal G_{\frak p_c} \to\mathcal G_{\frak p}$
shrinking several edges. For ${\rm v} \in C^0(\mathcal G_{\frak p})$ we put
$$
E_{c,{\rm v}} = \bigoplus_{{\rm v}' \in C^0(\mathcal G_{\frak p_c}) \atop \pi({\rm v}') = {\rm v}}E_{c,{\rm v}'}
$$
where $E_{c,{\rm v}'}$ is the obstruction bundle that is included in the obstruction bundle data
$\frak C_{\frak p_c}$ at $\frak p_c$.
It determines
$
E_{c,{\rm v}}(\frak y_c,u_c) = \bigoplus_{{\rm v}' \in C^0(\mathcal G_{\frak p_c}) \atop \pi({\rm v}') = {\rm v}}E_{c,{\rm v}'}(\frak y_c,u_c)
$. Then (\ref{2224}) is defined by Definition \ref{Emovevvv}.
\begin{rem}\label{rem267}
In Definition \ref{stabdata} (6) we assumed
that the image of $K_{{\rm v},c}^{\rm obst}$ by the
diffeomorphism mentioned above is always contained in the core of $\Sigma_{\vec T,\vec \theta}^{\rho}$.
(Here $K_{{\rm v},c}^{\rm obst}$ is the support of $E_{c,\rm v}$.)
Note by the core we mean the core with respect to the coordinate at infinity that is
included as a part of the stabilization data at $\frak p$ here.
\end{rem}
The vector space $E_c(u^{\rho}_{\vec T,\vec \theta,(0)})$ is the sum over ${\rm v} \in C^0(\mathcal G_{\frak p})$ of the images of (\ref{2224}).
\par
We next consider the obstruction bundle at $\hat u^{\rho}_{{\rm v},\vec T,\vec \theta,(0)}$.
A technical point we need to take care of here is that the obstruction bundle we use is {\it not}
$E_c(\coprod_{{\rm v}\in C^0(\mathcal G_{\frak p})}\hat u^{\rho}_{{\rm v},\vec T,\vec \theta,(0)})$
but is slightly different from it.
Let $K_{{\rm v},c}^{\rm obst} \subset K_{\rm v} \subset \Sigma_{\vec T,\vec \theta}^{\rho}$ be the
image of the set $K_{{\rm v},c}^{\rm obst}$ by the above mentioned diffeomorphism that is induced by  the
stabilization data at $\frak p$.
We remark that we may regard $K_{\rm v}$ as a subset of $\Sigma_{\rm v}^{\rho}$ also by using the
stabilization data at $\frak p$.
Moreover on  $K_{\rm v}$ we have
$\hat u^{\rho}_{{\rm v},\vec T,\vec \theta,(0)}
= u^{\rho}_{\vec T,\vec \theta,(0)}$.
So we have
\begin{equation}\label{lem77maesiki}
\aligned
\text{Image of (\ref{2224})} \,\,\subset \,\,
&\Gamma(K_{\rm v} ; (u^{\rho}_{\vec T,\vec \theta,(0)})^{*}TX \otimes \Lambda^{01})\\
&=
\bigoplus_{{\rm v} \in C^0(\mathcal G_{\frak p})}\Gamma(K_{\rm v} ; (\hat u^{\rho}_{{\rm v},\vec T,\vec \theta,(0)})^{*}TX \otimes \Lambda^{01}).
\endaligned
\end{equation}
\begin{defn}\label{defn277}
We regard the left hand side of (\ref{lem77maesiki})
as a subspace of
$$
\bigoplus_{{\rm v} \in C^0(\mathcal G_{\frak p})}\Gamma(K_{\rm v} ; (\hat u^{\rho}_{{\rm v},\vec T,\vec \theta,(0)})^{*}TX \otimes \Lambda^{01})
$$
and denote it by
$$
\bigoplus_{{\rm v} \in C^0(\mathcal G_{\frak p})} E'_c(\hat u^{\rho}_{{\rm v},\vec T,\vec \theta,(0)})
\subset
\bigoplus_{{\rm v} \in C^0(\mathcal G_{\frak p})} L^2_{m,\delta}(\Sigma_{\rm v};(\hat u^{\rho}_{{\rm v},\vec T,\vec \theta,(0)})^{*}TX \otimes \Lambda^{01}).
$$
We also define
$$
\mathcal E'_{\frak p,{\rm v},\frak A}(\hat u^{\rho}_{{\rm v},\vec T,\vec \theta,(0)}) =
\bigoplus_{c\in \frak A}
E'_c(\hat u^{\rho}_{{\rm v},\vec T,\vec \theta,(0)}),
\quad
\mathcal E'_{\frak p,\frak A}(\hat u^{\rho}_{\vec T,\vec \theta,(0)}) =
\bigoplus_{{\rm v} \in C^0(\mathcal G_{\frak p})}\mathcal E'_{\frak p,{\rm v},\frak A}(\hat u^{\rho}_{{\rm v},\vec T,\vec \theta,(0)}).
$$
\end{defn}
\begin{rem}
The reason why $E'_c(\hat u^{\rho}_{{\rm v},\vec T,\vec \theta,(0)}) \ne E_c(\hat u^{\rho}_{{\rm v},\vec T,\vec \theta,(0)})$
is as follows.
The union of the domains of $\hat u^{\rho}_{{\rm v},\vec T,\vec \theta,(0)}$ over $\rm v$ is
$\Sigma_{\frak p}$. When we identify the core of $\Sigma_{\frak p}$ with the core of $\Sigma_{\vec T,\vec \theta}^{\rho}$,
we use the additional marked points $\vec w_{\frak p}$ included in the stabilization data at $\frak p$.
We now consider the two diffeomorphisms:
\begin{eqnarray}
K_{{\rm v},c}^{\rm obst}    &\longrightarrow&  \text{Core of $\Sigma_{\vec T,\vec \theta}^{\rho}$}    \longrightarrow
\text{Core of $\Sigma_{\frak p}$} \label{2225}\\
K_{{\rm v},c}^{\rm obst}    &\longrightarrow& \text{Core of $\Sigma_{\frak p}$}.\label{2226}
\end{eqnarray}
We note that the diffeomorphism of the second arrow of (\ref{2225}) is defined by using the additional marked points $\vec w_{\frak p}$.
The other arrows are defined by using the additional marked points $\vec w_{\frak p_c}$.
Therefore in general (\ref{2225}) $\ne$ (\ref{2226}).
The definition of  $E'_c(\hat u^{\rho}_{{\rm v},\vec T,\vec \theta,(0)})$ uses (\ref{2225})
and the definition of  $E_c(\hat u^{\rho}_{{\rm v},\vec T,\vec \theta,(0)})$ uses (\ref{2226}).
This phenomenon does not occur in the situation of Part \ref{secsimple}.
This is because we took $\frak p = \frak p_c$ in Part \ref{secsimple}.
\end{rem}
\begin{rem}\label{rem270}
In the situation of Part \ref{secsimple} we have
${\rm Err}^{\rho}_{{\rm v},\vec T,\vec \theta,(0)} = 0$ on the core $K_{\rm v}$.
However this is not the case in the current situation.
In fact, by definition we have
\begin{equation}\label{2227}
\sum_{\rm v \in C^0(\mathcal G_{\frak p})}{\rm Err}^{\rho}_{{\rm v},\vec T,\vec \theta,(0)}
=
\overline\partial u_{\vec T,\vec \theta,(0)}^{\rho}  - \frak{se}^{\rho}_{\vec T,\vec \theta,(0)}
\end{equation}
and
\begin{equation}\label{2228}
\frak{se}^{\rho}_{\vec T,\vec \theta,(0)} = \sum_{c\in \frak A} \frak e^{\rho} _{c,\vec T,\vec \theta,(0)} = \overline\partial u_{\rm v}^{\rho}
\end{equation}
on $K_{\rm v}$. Moreover $u_{\rm v}^{\rho} = u_{\vec T,\vec \theta,(0)}^{\rho}$ on $K_{\rm v}$.
However (\ref{2227}) is nonzero because  the way how we identify an element $\frak e^{\rho} _{c,\vec T,\vec \theta,(0)} \in E_c$
as a section on $K_{\rm v}$ are different between the case of $u_{\rm v}^{\rho}$ and of $u_{\vec T,\vec \theta,(0)}^{\rho}$.
Namely, in (\ref{2227}) we regard $\frak e^{\rho} _{c,\vec T,\vec \theta,(0)}$
(that is a part of $ \frak{se}^{\rho}_{\vec T,\vec \theta,(0)}$)
as an element of $E_c(u_{\vec T,\vec \theta,(0)}^{\rho})$.
In (\ref{2228}) we regard $\frak e^{\rho} _{c,\vec T,\vec \theta,(0)}$
as an element of $E_c(u_{\rm v}^{\rho})$.
\par
We identify $K_{\rm v} \subset \Sigma_{\vec T,\vec \theta}^{\rho}$ with
$K_{\rm v} \subset \Sigma^{\rho}_{\rm v}$ by using the stabilization data at $\frak p$.
Thus $\frak e^{\rho} _{c,\vec T,\vec \theta,(0)}$ in (\ref{2227})
is also regarded as an element of  $E'_c(u_{\rm v}^{\rho})$.
So ${\rm Err}^{\rho}_{{\rm v},\vec T,\vec \theta,(0)}$ is nonzero on $K_{\rm v}$
because of $E'_c(u_{\rm v}^{\rho}) \ne E_c(u_{\rm v}^{\rho})$.
But this difference is of exponentially small. Namely we have the next lemma.
\end{rem}
\begin{lem}\label{0estimatekkk}
Put $T_{\rm min} = \min\{T_{\rm e} \mid {\rm e} \in C^1(\mathcal G_{\frak p})\}$. Then there exists $T_m$ such that the following inequality holds
\begin{equation}\label{2251}
\left\Vert
\nabla_{\rho}^n \frac{\partial^{\vert \vec k_{T}\vert}}{\partial T^{\vec k_{T}}}\frac{\partial^{\vert \vec k_{\theta}\vert}}{\partial \theta^{\vec k_{\theta}}}
{\rm Err}^{\rho}_{{\rm v},\vec T,\vec \theta,(0)} \right\Vert_{L^2_{m-\vert{\vec k_{T}}\vert - \vert{\vec k_{\theta}}\vert-1,\delta}(\Sigma_{\rm v}^{\rho})}
< C_{7,m}e^{-\delta  T_{{\rm min}}}
\end{equation}
for $\vert{\vec k_{T}}\vert + \vert{\vec k_{\theta}}\vert \le m - 10$ and $T_{{\rm min}} > T_m$.
\end{lem}
The proof is given later right after the proof of Lemma \ref{288lam}.
\par
In Definition \ref{defn277} we defined $\mathcal E'_{\frak p,{\rm v},\frak A}(\cdot)$ for
$\cdot = \hat u^{\rho}_{{\rm v},\vec T,\vec \theta,(0)}$. We next extend it
to nearby maps.
Let $u'_{\rm v} : (\Sigma_{\rm v}^{\rho},\partial\Sigma_{\rm v}^{\rho}) \to (X,L)$ be a smooth map
which is sufficiently close to $\hat u^{\rho}_{{\rm v},\vec T,\vec \theta,(0)}$ in
$C^{10}$ sense on $K_{\rm v}$.
We define  $\mathcal E'_{\frak p,{\rm v}, \frak A}(u'_{\rm v})$ as follows.
We identify $K_{\rm v}$ with a subset of  $\Sigma_{\vec T,\vec \theta}^{\rho}$
by using the additional marked points $\vec w_{\frak p}^{\rho}$.
Take any $u'' : (\Sigma_{\vec T,\vec \theta}^{\rho},\partial \Sigma_{\vec T,\vec \theta}^{\rho}) \to (X,L)$
that coincides with $u'_{\rm v}$ on $K_{\rm v}$ and is enough close to
$u_{\frak p}$ so that
$E_{c}(u'')=\bigoplus_{{\rm v}'\in C^0(\mathcal G_{\frak p_c})}
E_{c,{\rm v}'}(u'')$ is defined.
We put
$$
E_{c,{\rm v}}(u'') = \bigoplus_{{\rm v}' \in C^0 (\mathcal G_{\frak p_c}) \atop \pi({\rm v}') = {\rm v}}E_{c,{\rm v}'}(u'').
$$
By definition, $E_{c,\rm v}(u'')$ is independent of $u''$ but depends only on $u'_{\rm v}$
and is in
$\Gamma(K_{\rm v} ; (u'_{\rm v})^{*}TX \otimes \Lambda^{01})$.
Again using the diffeomorphism which is defined by the marked points $\vec w_{\frak p}^{\rho}$
we identify this space as a subspace of
$\Gamma(\Sigma_{\rm v}^{\rho} ; (u'_{\rm v})^{*}TX \otimes \Lambda^{01})$.
That is by definition $E'_{c,\rm v}(u'_{\rm v})$.
(This is the case ${\rm v} \in C^0_{\rm d}(\mathcal G_{\frak p})$.
The case of  ${\rm v} \in C^0_{\rm s}(\mathcal G_{\frak p})$ is similar.)
We put
\begin{equation}\label{2230}
\mathcal E'_{\frak p,{\rm v},\frak A}(u'_{\rm v}) = \sum_{c \in \frak A}E'_{c,{\rm v}}(u'_{\rm v}),
\qquad
\mathcal E'_{\frak p,\frak A}(u') = \sum_{{\rm v} \in C^0(\mathcal G_{\frak p})} \mathcal E'_{\frak p,{\rm v},\frak A}(u'_{\rm v}).
\end{equation}
\par
Let
$$
\Pi_{\mathcal E'_{\frak p,\frak A}(u')} ~:~
\bigoplus_{{\rm v} \in C^0(\mathcal G_{\frak p})}L^2_{m,\delta}(\Sigma_{\rm v}^{\rho};(u'_{{\rm v}})^{*}TX \otimes \Lambda^{01})
\to
\mathcal E'_{\frak p,\frak A}(u')
$$
be the $L^2$-orthogonal projection.
We next define its derivation by an element
$$
v = (v_{\rm v}) \in \bigoplus_{{\rm v} \in C^0_{\rm d}(\mathcal G_{\frak p})}\Gamma((\Sigma^{\rho}_{\rm v},\partial\Sigma^{\rho}_{\rm v});(u_{\rm v}')^{*}TX,(u_{\rm v}')^{*}TL)
\oplus
\bigoplus_{{\rm v} \in C^0_{\rm s}(\mathcal G_{\frak p})}\Gamma(\Sigma^{\rho}_{\rm v};(u_{\rm v}')^{*}TX)
$$
by
\begin{equation}\label{DEidefvv}
(D_{u'_{\rm v}}\mathcal E'_{\frak p,\frak A})((A_{\rm v}),(v_{\rm v})) = \frac{d}{ds}(\Pi_{\mathcal E'_{\frak p,\frak A}({\rm E} (u_{\rm v}',sv_{\rm v}))}(A_{\rm v}))\vert_{s=0}
\end{equation}
as in (\ref{DEidef}),
where
$$
A_{\rm v} \in L^2_{m,\delta}(\Sigma_{\rm v}^{\rho};(u_{\rm v}')^{*}TX \otimes \Lambda^{01}).
$$
\par
We use the operator
\begin{equation}\label{144opvv}
V \mapsto D_{\hat u^{\rho}_{{\rm v},\vec T,\vec \theta,(0)}}\overline{\partial}(V) - (D_{\hat u^{\rho}_{{\rm v},\vec T,\vec \theta,(0)}}
\mathcal E'_{\frak p,\frak A})(\frak{se}^{\rho}_{\vec T,\vec \theta,(0)} , V)
\end{equation}
as the linearization of the Cauchy-Riemann equation modulo
$\mathcal E'_{\frak A}$.
\footnote{Here we consider $\mathcal E_{\frak A}$ and not $\mathcal E'_{\frak A}$.
Note we are studying the Cauchy-Riemann equation for $u^{\rho}_{\vec T,\vec\theta,(0)}$.
The obsutruction space $\mathcal E'_{\frak A}(\hat u^{\rho}_{{\rm v},\vec T,\vec \theta,(0)})$
is sent to $\mathcal E_{\frak A}(u^{\rho}_{\vec T,\vec\theta,(0)})$ by the identification using the stabilization data at $\frak p$.}
\par
We recall that
$$
L^2_{m+1,\delta}((\Sigma^{\rho},\partial \Sigma^{\rho});(\hat u^{\rho}_{\vec T,\vec \theta,(0)})^*TX,(\hat u^{\rho}_{\vec T,\vec \theta,(0)})^*TL)
$$
is the kernel of (\ref{evaluationmkepot}) for $\hat u^{\rho}_{\vec T,\vec \theta,(0)} = (\hat u^{\rho}_{{\rm v},\vec T,\vec \theta,(0)})
_{{\rm v}\in C^0(\mathcal G_{\frak p})}$.
The direct sum of (\ref{144opvv}) induces an operator on
$
L^2_{m+1,\delta}((\Sigma^{\rho},\partial \Sigma^{\rho});(\hat u^{\rho}_{\vec T,\vec \theta,(0)})^*TX,(\hat u^{\rho}_{\vec T,\vec \theta,(0)})^*TL)
$
by restriction.
\begin{lem}\label{28011}
The sum of the image of the direct sum of the operators (\ref{144opvv}) on
$$L^2_{m+1,\delta}((\Sigma^{\rho},\partial \Sigma^{\rho});(\hat u^{\rho}_{\vec T,\vec \theta,(0)})^*TX,(\hat u^{\rho}_{\vec T,\vec \theta,(0)})^*TL)$$ and the subspace $\mathcal E'_{\frak p,\frak A}(\hat u^{\rho}_{\vec T,\vec \theta,(0)})$
is
$$
\bigoplus_{{\rm v}\in C^0(\mathcal G_{\frak p})} L^2_{m,\delta}(\Sigma^{\rho}_{\rm v};(\hat u^{\rho}_{{\rm v},\vec T,\vec \theta,(0)})^{*}TX \otimes \Lambda^{01})
$$
if
$\vec T$ is sufficiently large.
\end{lem}
\begin{proof}
This is a consequence of Lemma \ref{nbdregmaineq}.
\end{proof}
Lemma \ref{28011} is a generalization of Lemma \ref{lem112}.
\begin{defn}
The $L^2$ orthogonal complement of
$$
\left(D_{\hat u^{\rho}_{\vec T,\vec \theta,(0)}}\overline{\partial} - (D_{\hat u^{\rho}_{\vec T,\vec \theta,(0)}}\mathcal E'_{\frak p,\frak A}))(\frak{se}^{\rho}_{\vec T,\vec \theta,(0)} , \cdot)\right)^{-1}(\mathcal E'_{\frak p,\frak A}(\hat u^{\rho}_{\vec T,\vec \theta,(0)}))
$$
in
$$
L^2_{m+1,\delta}((\Sigma^{\rho},\partial \Sigma^{\rho});(\hat u^{\rho}_{\vec T,\vec \theta,(0)})^*TX,(\hat u^{\rho}_{\vec T,\vec \theta,(0)})^*TL)
$$
is denoted by $\frak H(\rho,\vec T,\vec \theta)$.
\par
We take $\vec T = \vec\infty = (\infty,\dots,\infty)$
and write $\frak H(\rho) = \frak H(\rho,\vec\infty,\vec \theta_0)$.
Then the restriction of  (\ref{144opvv}) to $\frak H(\rho)$ induces an isomorphism
to
$$
\bigoplus_{{\rm v}\in C^0(\mathcal G_{\frak p})} L^2_{m,\delta}(\Sigma^{\rho}_{\rm v};(\hat u^{\rho}_{{\rm v},\vec T,\vec \theta,(0)})^{*}TX \otimes \Lambda^{01})/\mathcal E'_{\frak p,\frak A}(\hat u^{\rho}_{\vec T,\vec \theta,(0)})
$$
for sufficiently large $\vec T$.
\end{defn}
\begin{defn}
We define
$V^{\rho}_{\vec T,\vec \theta,{\rm v},(1)}$ for ${\rm v} \in C^0(\mathcal G_{\frak p})$
and $\Delta p^{\rho}_{{\rm e},\vec T,\vec \theta,(1)}$ for ${\rm e} \in C^1(\mathcal G_{\frak p})$
so that $((V^{\rho}_{\vec T,\vec \theta,{\rm v},(1)})_{\rm v},(\Delta p^{\rho}_{{\rm e},\vec T,\vec \theta,(1)})_{\rm e}) \in \frak H(\rho)$
is the unique element such that
\begin{equation}
\aligned
D_{\hat u^{\rho}_{{\rm v},\vec T,\vec \theta,(0)}}\overline{\partial}(V^{\rho}_{\vec T,\vec \theta,{\rm v},(1)})
&- (D_{\hat u^{\rho}_{{\rm v},\vec T,\vec \theta,(0)}}\mathcal E'_{\frak p,\frak A})(\frak{se}^{\rho}_{\vec T,\vec \theta,(0)} ,
V^{\rho}_{\vec T,\vec \theta,{\rm v},(1)})\\
&+ {\rm Err}^{\rho}_{{\rm v},\vec T,\vec \theta,(0)}
\in \mathcal E'_{\frak p,\frak A}(\hat u^{\rho}_{{\rm v},\vec T,\vec \theta,(0)})
\endaligned
\end{equation}
and
\begin{equation}
\lim_{\tau_{\rm e} \to \pm \infty} V^{\rho}_{\vec T,\vec \theta,{\rm v},(1)}(\tau_{\rm e},t_{\rm e}) = \Delta p^{\rho}_{{\rm e},\vec T,\vec \theta,(1)},
\end{equation}
where $\pm \infty= + \infty$ if $\rm e$ is outgoing and $=-\infty$ if $\rm e$ is incoming.
\end{defn}
\par\medskip
\noindent{\bf Step 1-2}:
\begin{defn}
We define $u_{\vec T,\vec \theta,(1)}^{\rho}(z)$ as follows.
(Here $\rm E$ is the map as in (\ref{defE}).)
\begin{enumerate}
\item If $z \in K_{\rm v}$, we put
\begin{equation}
u_{\vec T,\vec \theta,(1)}^{\rho}(z)
=
{\rm E} (u_{\vec T,\vec \theta,(0)}^{\rho}(z),V^{\rho}_{\vec T,\vec \theta,{\rm v},(1)}(z)).
\end{equation}
\item
If $z  = (\tau_{\rm e},t_{\rm e}) \in [-5T_{\rm e},5T_{\rm e}]\times [0,1]$ or $S^1$,
we put
\begin{equation}
\aligned
u_{\vec T,\vec \theta,(1)}^{\rho}(\tau_{\rm e},t_{\rm e}) =
&\chi_{{\rm v}_{\leftarrow}({\rm e}),\mathcal B}^{\leftarrow}(\tau_{\rm e},t_{\rm e}) (V^{\rho}_{\vec T,\vec \theta,{\rm v}_{\leftarrow}({\rm e}),(1)}
(\tau_{\rm e},t_{\rm e}) -
\Delta p^{\rho}_{{\rm e},\vec T,\vec \theta,(1)})\\
&+\chi_{{\rm v}_{\rightarrow}({\rm e}),\mathcal A}^{\rightarrow}(\tau_{\rm e},t_{\rm e})(V^{\rho}_{\vec T,\vec \theta,{\rm v}_{\rightarrow}({\rm e}),(1)}(\tau_{\rm e},t_{\rm e})-\Delta p^{\rho}_{{\rm e},\vec T,\vec \theta,(1)})\\
&+u_{\vec T,\vec \theta,(0)}^{\rho}(\tau_{\rm e},t_{\rm e}) +\Delta p^{\rho}_{{\rm e},\vec T,\vec \theta,(1)}.
\endaligned
\end{equation}
\end{enumerate}
\end{defn}
\par\medskip
\noindent{\bf Step 1-3}:
We define:
\begin{equation}
\frak e^{\rho} _{1,T,(1)}
=
\Pi_{\mathcal E_{\frak p,\frak A}({\rm E} (u_{\vec T,\vec \theta,(0)}^{\rho},V^{\rho}_{\vec T,\vec \theta,{\rm v},(1)})}(\overline\partial
{\rm E} (u_{\vec T,\vec \theta,(0)}^{\rho},V^{\rho}_{\vec T,\vec \theta,{\rm v},(1)}))
\end{equation}
and
\begin{equation}
\frak{se}^{\rho}_{\vec T,\vec \theta,(1)} =  \frak{e}^{\rho}_{\vec T,\vec \theta,(0)} + \frak e^{\rho} _{\vec T,\vec \theta,(1)}.
\end{equation}
\par\medskip
\noindent{\bf Step 1-4}:
We take $0 < \mu < 1$ and fix it throughout the proof of this subsection.
\begin{defn}
We put
\begin{equation}
{\rm Err}^{\rho}_{{\rm v},\vec T,\vec \theta,(1)}  =
\begin{cases}
\chi_{{\rm e},\mathcal X}^{\leftarrow} \overline\partial u_{\vec T,\vec \theta,(1)}^{\rho}
& \text{on the ${\rm e}$-th neck if $\rm e$ is outgoing} \\
\chi_{{\rm e},\mathcal X}^{\rightarrow}  \overline\partial u_{\vec T,\vec \theta,(1)}^{\rho}
& \text{on the ${\rm e}$-th neck if $\rm e$ is incoming} \\
\overline\partial u_{\vec T,\vec \theta,(0)}^{\rho}  - \frak{se}^{\rho}_{\vec T,\vec \theta,(1)} & \text{on $K_{\rm v}$} .
\end{cases}
\end{equation}
We extend them by $0$ outside a compact set and will regard them as elements of the function space
$L^2_{m,\delta}(\Sigma_{\rm v}^{\rho};(\hat u^{\rho}_{{\rm v},\vec T,\vec \theta,(1)})^{*}TX \otimes \Lambda^{01})$,
where $\hat u^{\rho}_{{\rm v},\vec T,\vec \theta,(1)}$ will be defined in the next step.
\end{defn}
We put
$p^{\rho}_{{\rm e},\vec T,\vec \theta,(1)} = p^{\rho}_{{\rm e},\vec T,\vec \theta,(0)} + \Delta p^{\rho}_{{\rm e},\vec T,\vec \theta,(1)}$.
\par
We now come back to Step 2-1 and continue inductively on $\kappa$.
\par
The main estimate of those objects are the next lemma.
We put $R_{(\rm v,e)} = 5T_{\rm e} + 1$ and $\vec R = (R_{(\rm v,e)})$.
\begin{prop}\label{expesgen1}
There exist $T_m, C_{8,m}, C_{9,m}, C_{10,m}, \epsilon_{5,m} > 0$ and $0<\mu<1$ such that
the following inequalities hold if $T_{\rm e}>T_m$ for all $\rm e$.
We put $\vec{T}=(T_e; e\in C^1(\mathcal G_{\frak p}))$ and
$T_{\rm min} = \min \{T_{\rm e} \mid {\rm e} \in C^1(\mathcal G_{\frak p})\}$.
\begin{eqnarray}
\left\Vert \left((V^{\rho}_{\vec T,\vec \theta,{\rm v},(\kappa)}),(\Delta p^{\rho}_{{\rm e},\vec T,\vec \theta,(\kappa)})\right)\right\Vert_{L^2_{m+1,\delta}(\Sigma^{\rho}_{\rm v})}
&<& C_{8,m}\mu^{\kappa-1}e^{-\delta T_{\rm min}}, \label{form0182vv}
\\
\left\Vert (\Delta p^{\rho}_{{\rm e},\vec T,\vec \theta,(\kappa)})\right\Vert
&<& C_{8,m}\mu^{\kappa-1}e^{-\delta T_{\rm min}}, \label{form0183vv}
\\
\left\Vert  u_{\vec T,\vec \theta,(\kappa)}^{\rho}- u_{\vec T,\vec \theta,(0)}^{\rho}  \right\Vert_{L^2_{m+1,\delta}
(K_{\rm v}^{+\vec R})}
&<& C_{9,m}e^{-\delta T_{\rm min}}, \label{form0184}
\\
\left\Vert{\rm Err}^{\rho}_{{\rm v},\vec T,\vec \theta,(\kappa)} \right\Vert_{L^2_{m,\delta}(\Sigma^{\rho}_{\rm v})}
&<& C_{10,m}\epsilon_{5,m}\mu^{\kappa}e^{-\delta T_{\rm min}}, \label{form0185vv}
\\
\left\Vert \frak e^{\rho} _{\vec T,\vec \theta,(\kappa)}\right\Vert_{L^2_{m}(K_{\rm v}^{\text{\rm obst}})}
&<& C_{10,m}\mu^{\kappa-1}e^{-\delta T_{\rm min}},
\label{form0186vv}
\end{eqnarray}
where we assume $\kappa \ge 1$ in (\ref{form0186vv}).
\end{prop}
\begin{proof}
The proof is the same as the discussion in Subsection \ref{alternatingmethod}
and so is omitted.\footnote{Actually we need some new argument for the case $\kappa = 0$ of
(\ref{form0185vv}). We will discuss it later during the proof of Lemma \ref{sublem278}.}
\end{proof}
(\ref{form0182vv}) implies that the limit of
$
u_{\vec T,\vec \theta,(\kappa)}^{\rho}
$
converges as $\kappa$ goes to $\infty$  after $C^k$ topology for each $k$ if $T_{\rm e} > T_{k+10}$ for all $\rm e$.
We define
\begin{equation}
{\rm Glu}_{\vec T,\vec \theta}(\rho) =  \lim_{\kappa\to \infty} u_{\vec T,\vec \theta,(\kappa)}^{\rho} =
u_{\vec T,\vec \theta}^{\rho}.
\end{equation}
(\ref{form0185vv}) and (\ref{form0186vv}) imply
$$
\overline\partial u_{\vec T,\vec \theta}^{\rho}
= \sum_{\kappa=0}^{\infty}\frak e^{\rho} _{\vec T,\vec \theta,(\kappa)}
\in \mathcal E_{\frak A}(\overline\partial u_{\vec T,\vec \theta}^{\rho}).
$$
Therefore
$$
u_{\vec T,\vec \theta}^{\rho} \in
\mathcal M_{k+1,(\ell,\ell_{\frak p},(\ell_c))}(\beta;\frak p;\frak A;(\vec T^{\rm o},\vec T^{\rm c},\vec \theta))_{\epsilon_{2},\vec T_0}.
$$
We thus have defined ${\rm Glu}_{\vec T,\vec \theta}$.
\par
We next prove Theorem \ref{exdecayT33}.
The main part of the proof is the next lemma.
\begin{prop}\label{expesgen2Tdev}
There exist $T_m, C_{11,m}, C_{12,m}, C_{13,m}, C_{14,m}, \epsilon_{2,m} > 0$ and $0<\mu<1$ such that
the following inequalities hold if $T_{\rm e}>T_m$ for all $\rm e$.
\par
Let ${\rm e}_0 \in C^1_{\rm o}(\mathcal G_{\frak p})$.
Then for each $\vec k_{T}$,  $\vec k_{\theta}$ we have
\begin{equation}
\aligned
&\left\Vert \nabla_{\rho}^n \frac{\partial^{\vert \vec k_{T}\vert}}{\partial T^{\vec k_{T}}}\frac{\partial^{\vert \vec k_{\theta}\vert}}{\partial \theta^{\vec k_{\theta}}}
\frac{\partial}{\partial T_{{\rm e}_0}}((V^{\rho}_{\vec T,\vec \theta,{\rm v},(\kappa)}),(\Delta p^{\rho}_{{\rm e},\vec T,\vec \theta,(\kappa)}))\right\Vert_{L^2_{m+1-\vert{\vec k_{T}}\vert - \vert{\vec k_{\theta}}\vert-1,\delta}(\Sigma^{\rho}_{\rm v})}\\
&< C_{11,m}\mu^{\kappa-1}e^{-\delta T_{{\rm e}_0}}, \label{form0182vv2}
\endaligned
\end{equation}
\begin{equation}
\displaystyle\left\Vert \nabla_{\rho}^n \frac{\partial^{\vert \vec k_{T}\vert}}{\partial T^{\vec k_{T}}}\frac{\partial^{\vert \vec k_{\theta}\vert}}{\partial \theta^{\vec k_{\theta}}}
\frac{\partial}{\partial T_{{\rm e}_0}} (\Delta p^{\rho}_{{\rm e},\vec T,\vec \theta,(\kappa)})\right\Vert
< C_{11,m}\mu^{\kappa-1}e^{-\delta T_{{\rm e}_0}}, \label{form0183v2v}
\end{equation}
\begin{equation}
\left\Vert \nabla_{\rho}^n \frac{\partial^{\vert \vec k_{T}\vert}}{\partial T^{\vec k_{T}}}\frac{\partial^{\vert \vec k_{\theta}\vert}}{\partial \theta^{\vec k_{\theta}}}
\frac{\partial}{\partial T_{{\rm e}_0}}  u_{\vec T,\vec \theta,(\kappa)}^{\rho} \right\Vert_{L^2_{m+1-\vert{\vec k_{T}}\vert - \vert{\vec k_{\theta}}\vert
-1,\delta}(K_{\rm v}^{+\vec R})}
< C_{12,m}e^{-\delta T_{{\rm e}_0}}, \label{form01842}
\end{equation}
\begin{equation}
\aligned
&\left\Vert
\nabla_{\rho}^n \frac{\partial^{\vert \vec k_{T}\vert}}{\partial T^{\vec k_{T}}}\frac{\partial^{\vert \vec k_{\theta}\vert}}{\partial \theta^{\vec k_{\theta}}}
\frac{\partial}{\partial T_{{\rm e}_0}}
{\rm Err}^{\rho}_{{\rm v},\vec T,\vec \theta,(\kappa)} \right\Vert_{L^2_{m-\vert{\vec k_{T}}\vert - \vert{\vec k_{\theta}}\vert-1,\delta}(\Sigma^{\rho}_{\rm v})}
\\
&< C_{13,m}\epsilon_{6,m}\mu^{\kappa}e^{-\delta  T_{{\rm e}_0}}, \label{form0185vv2}
\endaligned
\end{equation}
\begin{equation}
\left\Vert
\nabla_{\rho}^n\frac{\partial^{\vert \vec k_{T}\vert}}{\partial T^{\vec k_{T}}}\frac{\partial^{\vert \vec k_{\theta}\vert}}{\partial \theta^{\vec k_{\theta}}}
\frac{\partial}{\partial T_{{\rm e}_0}}
\frak e^{\rho} _{\vec T,\vec \theta,(\kappa)}\right\Vert_{L^2_{m-\vert{\vec k_{T}}\vert - \vert{\vec k_{\theta}}\vert-1}(K_{\rm v}^{\text{\rm obst}})}
< C_{14,m}\mu^{\kappa-1}e^{-\delta  T_{{\rm e}_0}},
\label{form0186vv2}
\end{equation}
for $\vert{\vec k_{T}}\vert + \vert{\vec k_{\theta}}\vert + n < m -11$.
\par
Let ${\rm e}_{0} \in C^1_{\rm c}(\mathcal G_{\frak p})$. Then the same inequalities as above hold if we replace
$\frac{\partial}{\partial T_{{\rm e}_0}}$ by $\frac{\partial}{\partial \theta_{{\rm e}_0}}$.
\end{prop}
\begin{proof}[Proposition \ref{expesgen2Tdev} $\Rightarrow$ Theorem \ref{exdecayT33}]
Note if $k_{{\rm e}_0} \ne 0$ or $\theta_{{\rm e}_0} \ne 0$ then
$$
\vec k_{T}\cdot \vec T+\vec k_{\theta}\cdot \vec T^{\rm c}
\le 2k \max \{T_{\rm e} \mid k_{T,{\rm e}} \ne 0, \,\text{or }
k_{\theta,{\rm e}} \ne 0 \}.
$$
It is then easy to see that Proposition \ref{expesgen2Tdev} implies Theorem \ref{exdecayT33} by putting $\delta' = \delta/2k$.
\end{proof}
\begin{proof}[Proof of Proposition \ref{expesgen2Tdev}]
The proof is mostly the same as the argument of Subsection \ref{subsecdecayT}.
The new part is the proof of the next lemma.
\begin{lem}\label{sublem278}
Let ${\rm e}_{0} \in C^1_{\rm c}(\mathcal G_{\frak p})$.
We have
\begin{equation}\label{2251}
\left\Vert
\nabla_{\rho}^n \frac{\partial^{\vert \vec k_{T}\vert}}{\partial T^{\vec k_{T}}}\frac{\partial^{\vert \vec k_{\theta}\vert}}{\partial \theta^{\vec k_{\theta}}}
\frac{\partial}{\partial T_{{\rm e}_0}}
{\rm Err}^{\rho}_{{\rm v},\vec T,\vec \theta,(0)} \right\Vert_{L^2_{m-\vert{\vec k_{T}}\vert - \vert{\vec k_{\theta}}\vert-1,\delta}(\Sigma^{\rho}_{\rm v})}
< C_{15,m}e^{-\delta  T_{{\rm e}_0}}
\end{equation}
and
\begin{equation}\label{2252}
\left\Vert
\nabla_{\rho}^n \frac{\partial^{\vert \vec k_{T}\vert}}{\partial T^{\vec k_{T}}}\frac{\partial^{\vert \vec k_{\theta}\vert}}{\partial \theta^{\vec k_{\theta}}}
\frac{\partial}{\partial \theta_{{\rm e}_0}}
{\rm Err}^{\rho}_{{\rm v},\vec T,\vec \theta,(0)} \right\Vert_{L^2_{m-\vert{\vec k_{T}}\vert - \vert{\vec k_{\theta}}\vert-1,\delta}(\Sigma^{\rho}_{\rm v})}
< C_{15,m}e^{-\delta  T_{{\rm e}_0}}.
\end{equation}
\end{lem}
\begin{proof}
We recall (\ref{2222}),
\begin{equation}
{\rm Err}^{\rho}_{{\rm v},\vec T,\vec \theta,(0)}  =
\begin{cases}
\chi_{{\rm e},\mathcal X}^{\leftarrow} \overline\partial u_{\vec T,\vec \theta,(0)}^{\rho}
& \text{on the ${\rm e}$-th neck if $\rm e$ is outgoing} \\
\chi_{{\rm e},\mathcal X}^{\rightarrow}  \overline\partial u_{\vec T,\vec \theta,(0)}^{\rho}
& \text{on the ${\rm e}$-th neck if $\rm e$ is incoming} \\
\overline\partial u_{\vec T,\vec \theta,(0)}^{\rho}  - \frak{se}^{\rho}_{\vec T,\vec \theta,(0)} & \text{on $K_{\rm v}$} .
\end{cases}
\end{equation}
We first estimate ${\rm Err}^{\rho}_{{\rm v},\vec T,\vec \theta,(0)}$ on the neck region.
Let ${\rm e} \in C^1_{\rm c}(\mathcal G_{\frak p})$ is an outgoing edge of $\rm v$.
Let ${\rm v}'$ be the other vertex of $\rm e$.
We have
\begin{equation}\label{neckerror1}
\aligned
&\hskip-3.4cm{\rm Err}^{\rho}_{{\rm v},\vec T,\vec \theta,(0)} (\tau'_{\rm e},t'_{\rm e}) \\
=
(1-\chi(\tau'_{\rm e}-5T_{\rm e}))\overline \partial
\bigg(
&p^{\rho}_{\rm e} + (1-\chi(\tau'_{\rm e}-6T_{\rm e}))
(u^{\rho}_{\rm v}(\tau'_{\rm e},t'_{\rm e}) - p^{\rho}_{\rm e})\\
&+ \chi(\tau'_{\rm e}-4T_{\rm e})(u^{\rho}_{{\rm v}'}(\tau'_{\rm e}-10T_{\rm e},t'_{\rm e}+\theta_{\rm e})
- p^{\rho}_{\rm e})\bigg).
\endaligned
\end{equation}
Note that we use the coordinates $(\tau'_{\rm e},t'_{\rm e})$ for $u^{\rho}_{\rm v}$ and
$(\tau''_{\rm e},t''_{\rm e})$ for $u^{\rho}_{{\rm v}'}$.
(See (\ref{cctau12}), (\ref{ccttt12}).)
The function $\chi$ is as in (\ref{chichi}).
\par
If ${\rm e}_{0} \ne {\rm e}$, then $\partial/\partial T_{{\rm e}_0}$ or $\partial/\partial \theta_{{\rm e}_0}$
of (\ref{neckerror1}) is zero.
\par
Let us study $\partial/\partial T_{{\rm e}}$ or $\partial/\partial \theta_{{\rm e}}$ of (\ref{neckerror1})
in case ${\rm e}_{0} = {\rm e}$.
We apply $\partial/\partial \theta_{{\rm e}}$ to the third line of (\ref{neckerror1})
to obtain
\begin{equation}\label{neckerror2}
\aligned
&(1-\chi(\tau'_{\rm e}-5T_{\rm e}))
\frac{\partial}{\partial \theta_{{\rm e}}}\overline \partial
\left(\chi(\tau'_{\rm e}-4T_{\rm e})u^{\rho}_{{\rm v}'}(\tau'_{\rm e}-10T_{\rm e},t'_{\rm e}+\theta_{\rm e})
\right)\\
&=
(1-\chi(\tau'_{\rm e}-5T_{\rm e}))\chi(\tau'_{\rm e}-4T_{\rm e})\overline \partial
\bigg(\frac{\partial}{\partial t'_{{\rm e}}}
u^{\rho}_{{\rm v}'}(\tau'_{\rm e}-10T_{\rm e},t'_{\rm e}+\theta_{\rm e})
\bigg).
\endaligned\end{equation}
Support of (\ref{neckerror2}) is in the domain
$
4T_{\rm e} - 1 \le \tau'_{\rm e}  \le 5T_{\rm e} + 1
$
that is
$
-6T_{\rm e}  - 1 \le \tau''_{\rm e}  \le -5T_{\rm e} + 1
$. There the $C^m$ norm of $u^{\rho}_{{\rm v}'}$ is estimated as
$$
\Vert u^{\rho}_{{\rm v}'}\Vert_{C^m([-6T_{\rm e}  - 1,-5T_{\rm e} + 1))} \le C_{11,m}\,e^{-5T_{\rm e} \delta_1}.
$$
On the other hand, the weight function $e_{\rm v,\delta}$ given in
\eqref{e1deltamulti} is estimated by $e^{5T_{\rm e}\delta}$ on the support.
(See (\ref{e1deltamulti}).)
Therefore this term has the required estimated.
(Note $\delta < \delta_1/10$.)
The other term or other case of the estimate on the neck region is similar.
\par\medskip
We next estimate  ${\rm Err}^{\rho}_{{\rm v},\vec T,\vec \theta,(0)}$ on the core.
As we explained in Remark \ref{rem270} this is nonzero because of the difference of the
parametrization of the core.
So to study it, we need to discuss the dependence of the parametrization of the core on the
coordinate at infinity. Proposition \ref{reparaexpest}, Corollary \ref{corestimatecoochange} and  Lemma \ref{changeinfcoorproppara}
give the estimate we need to study.
\par
We consider $\frak p_c$ and the obstruction bundle data $\frak E_{\frak p_c}$ there.
Let $\mathcal G_c$ be the combinatorial type of $\frak p_c$. Note
$\frak p \in \frak W_{\frak p_c}$ and  $(\frak x_{\frak p}\cup \vec w_c^{\frak p}, u_{\frak p})$ is $\epsilon_{\frak p_c}$-close to
$\frak p_c$.
Let  $\mathcal G(\frak p,c)$ be the combinatorial type of $(\frak x_{\frak p}\cup \vec w_c^{\frak p}, u_{\frak p})$.
By Definition \ref{epsiloncloseto} (1) we have
$
\mathcal G_c \succ \mathcal G(\frak p,c)
$.
Let
$$
\frak x_{\frak p}\cup \vec w_c^{\frak p} = {\overline\Phi}(\frak y_1,\vec T_{1},\vec \theta_{1}).
$$
Note that the singular point of $\frak p$ corresponds one to one to the edges $\rm e$ of $\frak y_{1}$ such that
$T_{1,{\rm e}} = \infty$.
\par
For each ${\rm v}' \in C^0(\mathcal G_{\frak p_c})$,
we denote the corresponding core of $\Sigma_{\frak p_c}$ by $K_{{\rm v}'}^c$.
We may also regard
$$
K_{{\rm v}'}^c \subset \Sigma_{\frak p}.
$$
Let $\pi : \mathcal G_{\frak p_c} \to \mathcal G_{\frak p}$ be a map
shrinking the edges $\rm e$ with $T_{\rm e} \ne \infty$.
We put ${\rm v} = \pi({\rm v}')$.
Then there exists $\vec R$ such that
\begin{equation}\label{inclusionisv}
K_{{\rm v}'}^c \subset K_{\rm v}^{+\vec R}.
\end{equation}
Here the right hand side is the core of the coordinate at infinity of $\frak p$,
that is included in the stabilization data of $\frak p$.
The inclusion (\ref{inclusionisv}) is
obtained from the map $\frak v_{\xi,\frak y,\vec T,\vec\theta}$
appearing in Lemma \ref{changeinfcoorproppara} as follows.
\par
We put
$$
\{{\rm v}(i)\mid i=1,\dots,n_{c,{\rm v}}\}
= \{{\rm v}' \in C^0(\mathcal G_{\frak p_c}) \mid \pi({\rm v}') = {\rm v}\}.
$$
We consider the union
$$
K^{c}_{{\rm v},0} = \bigcup_{i=1}^{n_{c,{\rm v}}}
K^{\rm obst}_{{\rm v}(i)} \subset \Sigma_{\frak p_c}.
$$
\par
We consider $\Sigma_{\vec T,\vec \theta}^{\rho}$ that is a domain of
$u_{\vec T,\vec \theta,(0)}^{\rho}$. The parameter $\rho$ includes both the marked points $\vec w^{\rho}_c$ and $\vec w^{\rho}_{\frak p}$.
By forgetting $\vec w^{\rho}_{\frak p}$ we have an embedding
$$
\frak v_{c,{\rm v}(i),\rho,\vec T,\vec\theta} : K^{c}_{{\rm v}(i)} \to \Sigma^{\rho}_{\vec T,\vec\theta}.
$$
(Here the parameter $\vec w^{\rho}_{\frak p}$
(that is a part of $\rho$) plays the role of the parameter
$\xi \in Q$ in Lemma \ref{changeinfcoorproppara}.)
\par
By forgetting $\vec w^{\rho}_{c}$ we have an embedding
$$
\frak v_{\frak p,{\rm v},\rho,\vec T,\vec\theta} : K_{\rm v} \to \Sigma^{\rho}_{\vec T,\vec\theta}.
$$
\par
We consider $K_{{\rm v}(i)}^{\rm obst} \subset K^c_{{\rm v}(i)}$ that is a compact set we fixed as a part of the
obstruction bundle data centered at $\frak p_{c}$.
By Remark \ref{rem267}, we may assume
$$
\frak v_{c,{\rm v}(i),\rho,\vec T,\vec\theta}(K_{{\rm v}(i)}^{\rm obst})
\subset \frak v_{\frak p,{\rm v},\rho,\vec T,\vec\theta}(K_{\rm v}).
$$
Therefore taking union over $i=1,\dots,n_{c,{\rm v}}$
we obtain
\begin{equation}\label{corechangeest}
\frak v_{(\frak p,c),{\rm v},\rho,\vec T,\vec\theta} :=
\frak v_{\frak p,{\rm v},\rho,\vec T,\vec\theta}^{-1}\circ
\left(\coprod_{i=1}^{n_{c,{\rm v}}} \frak v_{c,{\rm v}(i),\rho,\vec T,\vec\theta}
\right)
:
K^{c}_{{\rm v},0} \to K_{\rm v}.
\end{equation}
We denote this map by
$$
\text{\rm Res}(\frak v_{(\frak p,c),{\rm v},\rho,\vec T,\vec\theta}) \in C^m(K^{c}_{{\rm v},0},K_{\rm v}).
$$
We can estimate it by using Lemma \ref{changeinfcoorproppara} that is a family version of
Proposition \ref{reparaexpest} and
Corollary \ref{corestimatecoochange}. (See Lemma \ref{288lam} below.)
\par
We next describe the way how $\frak v_{(\frak p,c),{\rm v},\rho,\vec T,\vec\theta}$ and its
estimate are related to the estimate of ${\rm Err}^{\rho}_{{\rm v},\vec T,\vec \theta,(0)}$.
We first recall that
$$
\overline\partial u^{\rho}_{\rm v} \in \bigoplus_{c\in \frak A} E_{c,{\rm v}}
$$
by assumption. We denote by
$\frak e^{\rho} _{c,\vec T,\vec \theta,(0)}$ the sum of its $E_{c,{\rm v}}$ components over $\rm v$.
It is actually independent of $\vec T,\vec \theta$. So we write it
$\frak e^{\rho}_{c,(0)}$ here.
We remark that we identify
$$
E_{c,{\rm v}} \subset \Gamma_0(K_{\rm v};(u^{\rho}_{\rm v})^{*}TX\otimes \Lambda^{01})
$$
using the obstruction bundle data centered at $\frak p_c$.
Here $K_{\rm v} \subset \Sigma_{\frak y}$.
(Note that the combinatorial type of $\frak y$ is the same as $\frak p$.)
\par
In (\ref{22190}), we used $u^{\rho}_{\rm v}$ to obtain a map
$$
u_{\vec T,\vec \theta,(0)}^{\rho} :
(\Sigma_{\vec T,\vec \theta}^{\rho},\partial\Sigma_{\vec T,\vec \theta}^{\rho})
\to (X,L).
$$
Moreover $u^{\rho}_{\rm v} = u_{\vec T,\vec \theta,(0)}^{\rho}$
on $K_{\rm v}$.
However
$$
E_{c,{\rm v}}(u^{\rho}_{\rm v} ) \ne E_{c,{\rm v}}(u_{\vec T,\vec \theta,(0)}^{\rho}),
$$
as subsets of
$$
\Gamma(K_{\rm v};(u^{\rho}_{\rm v})^*TX \otimes \Lambda^{01})
=
\Gamma(K_{\rm v};(u_{\vec T,\vec \theta,(0)}^{\rho})^*TX \otimes \Lambda^{01}).
$$
In fact, $ E_{c,{\rm v}}(u_{\vec T,\vec \theta,(0)}^{\rho})$ is defined by the
diffeomorphism $\frak v_{(\frak p,c),{\rm v},\rho,\vec T,\vec\theta}$
and $E_{c,{\rm v}}(u^{\rho}_{\rm v} )$ is defined by the diffeomorphism
$\frak v_{(\frak p,c),{\rm v},\rho,\vec \infty}$.
\par
Therefore, by definition, ${\rm Err}^{\rho}_{{\rm v},\vec T,\vec \theta,(0)}$ on $K_{\rm v}$ is
\begin{equation}
\overline\partial u_{\vec T,\vec \theta,(0)}^{\rho} - \sum_c\frak e^{\rho}_{c,(0)}
= \sum_c \left(\frak e^{\rho,1}_{c,(0)} - \frak e^{\rho,2}_{c,(0)}\right),
\end{equation}
where $\frak e^{\rho,1}_{c,(0)} \in \bigoplus_{{\rm v} \in C^0(\mathcal G_{\frak p})}
E_{c,{\rm v}}(u_{\vec T,\vec \theta,(0)}^{\rho})$ and $\frak e^{\rho,2}_{c,(0)}
\in \bigoplus_{{\rm v} \in C^0(\mathcal G_{\frak p})}E_{c,{\rm v}}(u^{\rho}_{\rm v} )$ are defined as follows:
\begin{equation}
\aligned
\frak e^{\rho,1}_{c,(0)}(
\frak v_{(\frak p,c),{\rm v},\rho,\vec T,\vec\theta}(z))
&= \text{\rm Pal}_{u_{\frak p_c,{\rm v}}(z),u^{\rho}_{\rm v}(\frak v_{(\frak p,c),{\rm v},\rho,\vec T,\vec\theta}(z)}(\frak e^{\rho}_{c,(0)}), \\
\frak e^{\rho,2}_{c,(0)}(
\frak v_{(\frak p,c),{\rm v},\rho,\vec \infty}(z))
&= \text{\rm Pal}_{u_{\frak p_c,{\rm v}}(z),u^{\rho}_{\rm v}(\frak v_{(\frak p,c),{\rm v},\rho,\vec \infty}(z)}(\frak e^{\rho}_{c,(0)}).
\endaligned
\end{equation}
Thus Lemma \ref{288lam} below implies
$$
\left\Vert{\rm Err}^{\rho}_{{\rm v},\vec T,\vec \theta,(0)} \right\Vert_{L^2_{m,\delta}(K_{\rm v})}
< C_{8,m}\epsilon_{1,m}e^{-\delta T_{\rm min}}.
$$
This is the case  $\kappa = 0$ of
(\ref{form0185vv})   on $K_{\rm v}$.
\par
Proposition \ref{reparaexpest}  implies the estimate (\ref{2251}) and (\ref{2252})
on $K_{\rm v}$.
The proof of Lemma \ref{sublem278} is complete assuming Lemma \ref{288lam}.
\end{proof}
\begin{lem}\label{288lam}
There exist $C_{15,k}$, $T_k$ such that for each ${\rm e} \in C^1_{\rm c}(\mathcal G_{\frak p})$ we have:
\begin{equation}
\aligned
&\left\Vert
\nabla_{\frak y_2}^n \frac{\partial^{\vert \vec k_{T}\vert}}{\partial T_2^{\vec k_{T}}}\frac{\partial^{\vert \vec k_{\theta}\vert}}{\partial \theta_2^{\vec k_{\theta}}}
\frac{\partial}{\partial T_{2,{\rm e}_0}}
(\frak v_{(\frak p,c),{\rm v},\rho,\vec T,\vec\theta})
\right\Vert_{C^k}
< C_{15,k}e^{-\delta_2  T_{2,{\rm e}_0}},
\\
&\left\Vert
\nabla_{\frak y_2}^n \frac{\partial^{\vert \vec k_{T}\vert}}{\partial T_2^{\vec k_{T}}}\frac{\partial^{\vert \vec k_{\theta}\vert}}{\partial \theta_2^{\vec k_{\theta}}}
\frac{\partial}{\partial \theta_{2,{\rm e}_0}}
(\frak v_{(\frak p,c),{\rm v},\rho,\vec T,\vec\theta})
\right\Vert_{C^k}
< C_{15,k}e^{-\delta_2  T_{2,{\rm e}_0}},
\endaligned
\end{equation}
whenever $T_{2,{\rm e}}$ is greater than $T_k$ and
$\vert{\vec k_{T}}\vert + \vert{\vec k_{\theta}}\vert +n \le k$.
\par
The first inequality holds for ${\rm e} \in C^1_{\rm o}(\mathcal G_{\frak p})$ also.
\end{lem}
\begin{proof}
It suffices to prove the same estimate for
$\frak v_{\frak p,{\rm v},\rho,\vec T,\vec\theta}$ and $\frak v_{c,{\rm v}(i),\rho,\vec T,\vec\theta}$.
Note $\rho \in V_{k+1,(\ell,\ell_{\frak p},(\ell_c))}(\beta;\frak p;\frak A;\epsilon_0)$
contains various data.
We use only a part of such a data. We below recall the parameter space
which contains only the data we use below.
\par
Let
$\frak V(\frak x_{\frak p}\cup \vec w_{c}^{\frak p})$ be a neighborhood of
$\frak x_{\frak p}\cup \vec w_{c}^{\frak p}$ in the stratum of the
Deligne-Mumford moduli space that consists of elements of the
same combinatorial type as $\frak x_{\frak p}\cup \vec w_{c}^{\frak p}$.
We also take $\frak V(\frak x_{\frak p}\cup \vec w_{\frak p})$
and $\frak V(\frak x_{\frak p}\cup \vec w_{c}^{\frak p}\cup \vec w_{\frak p})$
that are neighborhoods in the stratum of the
Deligne-Mumford moduli space of $\frak x_{\frak p}\cup \vec w_{\frak p}$ and
$\frak x_{\frak p}\cup \vec w_{c}^{\frak p}\cup \vec w_{\frak p}$, respectively.
\par
We can take those three neighborhoods so that there exist $Q_1$ and $Q_2$ such that
\begin{equation}\label{2282iso}
Q_1 \times \frak V(\frak x_{\frak p}\cup \vec w_{c}^{\frak p})
\cong
\frak V(\frak x_{\frak p}\cup \vec w_{c}^{\frak p}\cup \vec w_{\frak p})
\cong
Q_2 \times \frak V(\frak x_{\frak p}\cup \vec w_{\frak p})
\end{equation}
and that the isomorphisms in (\ref{2282iso}) is compatible with the forgetful maps
$$
\frak V(\frak x_{\frak p}\cup \vec w_{c}^{\frak p}\cup \vec w_{\frak p}) \to \frak V(\frak x_{\frak p}\cup \vec w_{c}^{\frak p})
$$
and
$$
\frak V(\frak x_{\frak p}\cup \vec w_{c}^{\frak p}\cup \vec w_{\frak p}) \to
\frak V(\frak x_{\frak p}\cup \vec w_{\frak p}).
$$
We consider the universal family
$$
\frak M(\frak x_{\frak p}\cup \vec w_{c}^{\frak p}\cup \vec w_{\frak p}) \to
\frak V(\frak x_{\frak p}\cup \vec w_{c}^{\frak p}\cup \vec w_{\frak p}).
$$
Together with other data it gives a coordinate at infinity.
We take any of them.
\par
Using  (\ref{2282iso}), this coordinate at infinity of
$\frak x_{\frak p}\cup \vec w_{c}^{\frak p}\cup \vec w_{\frak p}$ induces
a $Q_1$-parametrized family of coordinates at infinity of
$\frak x_{\frak p}\cup \vec w_{c}^{\frak p}$ and
a $Q_2$-parametrized family of coordinates at infinity of
$\frak x_{\frak p}\cup \vec w_{\frak p}$.
(See Definition \ref{def:Qfamily} for the definition of a $Q$-parametrized family of coordinates at infinity.)
\par
Compared with the given coordinate at infinities of
$\frak x_{\frak p}\cup \vec w_{c}^{\frak p}$ and of $\frak x_{\frak p}\cup \vec w_{\frak p}$
we obtain the maps
$\frak v_{\frak p,{\rm v},\rho,\vec T,\vec\theta}$ and $\frak v_{c,{\rm v}(i),\rho,\vec T,\vec\theta}$.
Therefore Lemma \ref{288lam} follows from Lemma \ref{changeinfcoorproppara}.
\end{proof}
We thus have completed the first step of the induction to prove Proposition \ref{expesgen2Tdev}.
The other steps are similar to the proof of Theorem \ref{exdecayT}.
\par
When we study $T_{\rm e}$ and $\theta_{\rm e}$ derivatives
and prove Lemma \ref{expesgen2Tdev},
we again need to
estimate the $T_{\rm e}$ and $\theta_{\rm e}$ derivatives of the map
$$
E_{c} \to \Gamma_0(K_{\rm v};(u_{\vec T,\vec \theta,(\kappa)}^{\rho})^*TX \otimes \Lambda^{01}).
$$
This map is defined by using the diffeomorphism $\frak v_{(\frak p,c),{\rm v},\rho,\vec T,\vec\theta}$.
Therefore we can use Lemma \ref{288lam} in the same way as above to obtain
the required estimate.\footnote
{We remark that $E_c$ is a finite dimensional vector space consisting of smooth sections
with compact support. So estimating the effect of change of variables
of its element by $\frak v_{(\frak p,c),{\rm v},\rho,\vec T,\vec\theta}$ is easy using Lemma \ref{288lam}.}
\par
The proof of Proposition \ref{expesgen2Tdev} is complete.
\end{proof}
\begin{proof}[Proof of Lemma \ref{0estimatekkk}]
We can prove Lemma \ref{0estimatekkk} by integrating the inequality
in Lemma \ref{288lam}.
\end{proof}
Thus we have proved Theorem \ref{exdecayT33}.
\par
We can use it in the same way as in Section \ref{surjinj} to prove surjectivity and
injectivity of the map ${\rm Glu}_{\vec T,\vec \theta}$.
\par
To show that ${\rm Glu}_{\vec T,\vec \theta}$ is $\Gamma_{\frak p}^+$-equivariant, we only need to remark that
if $\frak p_c \in \frak C(\frak p)$ then $\Gamma_{\frak p}^+ \subseteq \Gamma_{\frak p_c}^+$.
(In fact all the constructions are equivariant.)
\par
The proof of Theorem \ref{gluethm3} is complete.
\qed
\begin{rem}
We close this subsection with another technical remark.
Theorems \ref{gluethm3} and \ref{exdecayT33} imply that
$$
\aligned
\text{\rm Glu} : V_{k+1,(\ell,\ell_{\frak p},(\ell_c))}(\beta;\frak p;\frak A;\frak B;\epsilon_1) &\times (\vec T^{\rm o}_0,\infty] \times ((\vec T^{\rm c}_0,\infty] \times \vec S^1)
\\
&\to
\mathcal M_{k+1,(\ell,\ell_{\frak p},(\ell_c))}(\beta;\frak p;\frak A;\frak B)_{\epsilon_{0},\vec T_{0}}
\endaligned
$$
is a strata-wise $C^m$ diffeomorphism  if
$T_{{\rm e},0}$ for all $\rm e$ is larger than a number {\it depending} on $m$.
Using Theorem \ref{exdecayT33} we can define
smooth structures on both sides so that the map becomes a $C^{m}$ diffeomorphism.
(See Section \ref{chart}. We will use $s_{\rm e}= T_{\rm e}^{-1}$ as a coordinate.)
\par
Note that the domain and the target of $\text{\rm Glu} $ have strata-wise
$C^{\infty}$ structure.\footnote{This is an easy consequence of implicit function theorem.}
However, the construction we gave does not show that $\text{\rm Glu}$
is of $C^{\infty}$-class.
This is not really an issue for our purpose of defining virtual fundamental chain or cycle.
Indeed, Kuranishi structure of $C^k$ class with sufficiently large $k$ is enough for such a purpose.
($C^1$-structure is enough.)
\par
On the other hand, as we will explain in Section \ref{toCinfty},  Theorems \ref{gluethm3} and \ref{exdecayT33} are
enough to prove the existence of Kuranishi structure of $C^{\infty}$ class.
Except in Section \ref{toCinfty}, we fix $m$ and will construct a Kuranishi structure of $C^m$ class.
For this purpose we choose $T_{{\rm e},0}$ so that it is larger than $T_{10m}$.
Therefore  our construction of $\text{\rm Glu}$ works on $L^2_{10m+1,\delta}$.
\end{rem}
\par\medskip
\section{Cutting down the solution space by transversals}
\label{cutting}

In Section \ref{glueing}, we described the thickened moduli space
$\mathcal M_{k+1,(\ell,\ell_{\frak p},(\ell_c))}(\beta;\frak p;\frak A;\frak B)_{\epsilon_{0},\vec T_{0}}$
by a gluing construction.
Its dimension is given by
$$\aligned
&\dim
\mathcal M_{k+1,(\ell,\ell_{\frak p},(\ell_c))}(\beta;\frak p;\frak A;\frak B)_{\epsilon_{0},\vec T_{0}}\\
&= \text{\rm virdim}\, \mathcal M_{k+1,\ell}(\beta)
+ \dim_{\R} \mathcal E_{\frak A} + (2\ell_{\frak p} + 2\sum_{c\in \frak B}\ell_c)\\
&=
k+1+2\ell + \mu(\beta) - 3
+ \dim_{\R} \mathcal E_{\frak A} + (2\ell_{\frak p} + 2\sum_{c\in \frak B}\ell_c).
\endaligned$$
\par
Note that the dimension of the Kuranishi neighborhood of $\frak p$
in $\mathcal M_{k+1,\ell}(\beta)$ must be
$\text{\rm virdim} \mathcal M_{k+1,\ell}(\beta) + \dim_{\R} \mathcal E_{\frak A}$.
Therefore
we need to cut down this moduli space
$\mathcal M_{k+1,(\ell,\ell_{\frak p},(\ell_c))}(\beta;\frak p;\frak A;\frak B)_{\epsilon_{0},\vec T_{0}}$
to obtain a Kuranishi neighborhood.
We do so by requiring the transversal constraint as in
Definition \ref{transconst}.
We will define it below in a slightly generalized form.
(For example, we define it for $(\frak x,u)$ such that $u$ is not necessarily
pseudo-holomorphic but satisfies the equation $\overline\partial u \equiv 0
\mod \mathcal E_{\frak A}(u)$ only.)
\par
Let $\frak p \in \mathcal M_{k+1,\ell}(\beta)$ and $\emptyset \ne \frak A \subseteq \frak B \subseteq \frak C(\frak p)$.
We consider a subset $\frak B^- \subseteq \frak B$ with $\frak A \subseteq \frak B^-$.
Let $\vec w_{\frak p} = (w_{\frak p,1},\dots,w_{\frak p,\ell_{\frak p}})$
be a symmetric stabilization of $\frak x_{\frak p}$ that is a part
of the stabilization data at $\frak p$.
Let $I\subset \{1,\dots,\ell_{\frak p}\}$ and we consider
$\vec w^{-}_{\frak p} = (w_{\frak p,i}; i \in I)$.
For simplicity of notation we put
$I= \{1,\dots,\ell^-_{\frak p}\}$.
We assume that $\vec w^-_{\frak p}$ is already a symmetric stabilization of $\frak x_{\frak p}$.
It induces a stabilization data at $\frak p$ in an obvious way.
We thus obtain
$\mathcal M_{k+1,(\ell,\ell^-_{\frak p},(\ell_c))}(\beta;\frak p;\frak A;\frak B^-)_{\epsilon_{0},\vec T_{0}}$.
\begin{defn}\label{def7289}
An element $(\frak Y,u',(\vec w'_c;c\in \frak B))$
of $\mathcal M_{k+1,(\ell,\ell_{\frak p},(\ell_c))}(\beta;\frak p;\frak A;\frak B)_{\epsilon_{0},\vec T_{0}}$ is said to satisfy the (partial)
{\it transversal constraint} for $\vec w_{\frak p} \setminus \vec w^-_{\frak p}$
and $\frak B \setminus \frak B^-$ if the following holds.
\begin{enumerate}
\item If $i > \ell^-_{\frak p}$ then
$
u'(w'_{\frak p,i}) \in \mathcal D_{\frak p,i}.
$
Here $w'_{\frak p,i}$, $i=1,\dots,\ell_{\frak p}$ denote the
$(\ell+1)$-th, \dots, $(\ell+\ell_{\frak p})$-th interior marked points of
$\frak Y$.
\item
If $c \in \frak B \setminus \frak B^-$ and $i=1,\dots,\ell_{c}$ then
$
u'(w'_{c,i}) \in \mathcal D_{c,i}$. Here
$\vec w'_c = (w'_{c,1},\dots,w'_{c,\ell_c})$.
\end{enumerate}
We denote by
$$
\mathcal M_{k+1,(\ell,\ell_{\frak p},(\ell_c))}(\beta;\frak p;\frak A;\frak B)
^{\vec w^-_{\frak p},\frak B^-}_{\epsilon_{0},\vec T_{0}}
$$
the set of all elements of the thickened moduli space $\mathcal M_{k+1,(\ell,\ell_{\frak p},(\ell_c))}(\beta;\frak p;\frak A;\frak B)_{\epsilon_{0},\vec T_{0}}$
satisfying  transversal constraint for $\vec w_{\frak p} \setminus \vec w^-_{\frak p}$
and $\frak B \setminus \frak B^-$.
\end{defn}
Our next goal is to show that $\mathcal M_{k+1,(\ell,\ell_{\frak p},(\ell_c))}(\beta;\frak p;\frak A;\frak B)
^{\vec w^-_{\frak p},\frak A^-}_{\epsilon_{0},\vec T_{0}}$
is homeomorphic to
$\mathcal M_{k+1,(\ell,\ell^-_{\frak p},(\ell_c))}(\beta;\frak p;\frak A;\frak B^-)_{\epsilon_{0},\vec T_{0}}$. (Proposition \ref{forgetstillstable}.)
To prove this we first define an appropriate forgetful map.
\begin{defn}\label{defnforget}
Let $(\frak Y,u',(\vec w'_c;c\in \frak B)) \in \mathcal M_{k+1,(\ell,\ell_{\frak p},(\ell_c))}(\beta;\frak p;\frak A;\frak B)_{\epsilon_{0},\vec T_{0}}$.
Note $\frak Y = \frak Y_0 \cup \vec w_{\frak p}$ and $\vec w_{\frak p}$ consists of $\ell_{\frak p}$ interior marked points.
We take only $\ell^-_{\frak p}$ of them and put $\vec w_{\frak p}^-$
and put $\frak Y^- = \frak Y_0 \cup \vec w_{\frak p}^-$.
We assume that  $\frak Y^-$ is stable and $\frak x_{\frak p} \cup \vec w_{\frak p}^-$ is also stable.
We also assume that $\Gamma_{\frak p}$ preserves  $\vec w_{\frak p}$ as a set.
We define the forgetful map by:
\begin{equation}\label{2265}
\frak{forget}_{\frak B,\frak B^-;\vec w_{\frak p},\vec w^-_{\frak p}}
(\frak Y,u',(\vec w'_c;c\in \frak B))
=
(\frak Y^-,u',(\vec w'_c;c\in \frak B^-)).
\end{equation}
\end{defn}
\begin{lem}\label{fraAforget1}
The map
$\frak{forget}_{\frak B,\frak B^-;\vec w_{\frak p},\vec w^-_{\frak p}}$
defines
$$
\mathcal M_{k+1,(\ell,\ell_{\frak p},(\ell_c))}(\beta;\frak p;\frak A;\frak B)_{\epsilon_{0},\vec T_{0}}
\to
\mathcal M_{k+1,(\ell,\ell^-_{\frak p},(\ell_c))}(\beta;\frak p;\frak A;\frak B^-)_{\epsilon_{0},\vec T_{0}}.
$$
This map is a continuous and strata-wise smooth submersion.  The fiber is
$2(\ell_{\frak p} - \ell^-_{\frak p}) + 2\sum_{c\in \frak B \setminus \frak B^-}\ell_c$
dimensional.
\end{lem}
\begin{proof}
We note that $\frak Y^-$ is still stable.
(This is because $\frak x_{\frak p} \cup \vec w^{-}_{\frak p}$ is stable.)
Therefore $\frak{forget}_{\frak B,\frak B^-;\vec w_{\frak p},\vec w^-_{\frak p}}$
preserves stratification.
Note we forget the position of the
$\ell_{\frak p} - \ell^-_{\frak p} + \sum_{c\in \frak B \setminus \frak B^-}\ell_c$
marked points. There is no constraint for those marked points
other than those coming from the condition that
$(\frak Y,u')$ is $\epsilon_0$-close to $(\frak x_{\frak p}\cup \vec w_{\frak p},u_{\frak p})$ and $(\frak Y_0 \cup \vec w'_c,u')$
are $\epsilon_0$-close to $\frak p \cup \vec w_c^{\frak p}$ for all $c \in \mathcal A$.
These are open conditions.
Therefore this map is a strata-wise smooth submersion and the fiber is
$2(\ell_{\frak p} - \ell^-_{\frak p}) + 2\sum_{c\in \frak B \setminus \frak B^-}\ell_c$
dimensional.
\end{proof}
\begin{prop}\label{forgetstillstable}
The following holds if $\epsilon_{0},\epsilon_{\frak p_{c}}$ are sufficiently small.
\begin{enumerate}
\item
The space
$\mathcal M_{k+1,(\ell,\ell_{\frak p},(\ell_c))}(\beta;\frak p;\frak A;\frak B)
^{\vec w^-_{\frak p},\frak B^-}_{\epsilon_{0},\vec T_0}$
is  a strata-wise smooth submanifold of our thickened moduli space
$\mathcal M_{k+1,(\ell,\ell_{\frak p},(\ell_c))}(\beta;\frak p;\frak A;\frak B)_{\epsilon_{0},\vec T_{0}}$
of codimension $2(\ell_{\frak p} - \ell^-_{\frak p}) + 2\sum_{c\in \frak B \setminus \frak B^-}\ell_c$.
\item
The restriction of $\frak{forget}_{\frak B,\frak B^-;\vec w_{\frak p},\vec w^-_{\frak p}}$ induces a homeomorphism
$$
\mathcal M_{k+1,(\ell,\ell_{\frak p},(\ell_c))}(\beta;\frak p;\frak A;\frak B)_{\epsilon_{0},\vec T_{0}}^{\vec w^-_{\frak p},\frak B^-}
\to
\mathcal M_{k+1,(\ell,\ell^-_{\frak p},(\ell_c))}(\beta;\frak p;\frak A;\frak B^-)_{\epsilon_{0},\vec T_{0}}
$$
that is a  strata-wise diffeomorphism.
\end{enumerate}
\end{prop}
\begin{rem}
Note that if $c \in \mathcal B$ then $\frak p \in \frak M_{\frak p_c}$
and $\epsilon_c$ is used to define $\frak M_{\frak p_c}$.
(See Definition \ref{openW+++}.)
\end{rem}
\begin{proof}
We consider the evaluation maps at the
$(\ell_{\frak p} - \ell^-_{\frak p}) + \sum_{c\in \frak B \setminus \frak B^-}\ell_c$
marked points that we forget by the map
$\frak{forget}_{\frak B,\frak B^-;\vec w_{\frak p},\vec w^-_{\frak p}}$.
It defines a continuous and strata-wise smooth map
\begin{equation}\label{evatforgottenmark}
\mathcal M_{k+1,(\ell,\ell_{\frak p},(\ell_c))}(\beta;\frak p;\frak A;\frak B)_{\epsilon_{0},\vec T_{0}}
\to X^{(\ell_{\frak p} - \ell^-_{\frak p}) + \sum_{c\in \frak B \setminus \frak B^-}\ell_c}.
\end{equation}
We consider the submanifold
\begin{equation}\label{transforgottenmark}
\prod_{i=\ell^-_{\frak p}+1}^{\ell_{\frak p}}\mathcal D_{\frak p,i}
\times
\prod_{c \in \frak B \setminus \frak B^-}\prod_{i=1}^{\ell_{c}}\mathcal D_{c,i}
\end{equation}
of the right hand side of (\ref{evatforgottenmark}). By Proposition \ref{linearMV} (2), the map
(\ref{evatforgottenmark}) is transversal to
(\ref{transforgottenmark}) at $\frak p$ if
$\epsilon_{\frak p_c}$ is sufficiently small.
Therefore we may assume (\ref{evatforgottenmark}) is transversal to
(\ref{transforgottenmark}) everywhere.
Since $\mathcal M_{k+1,(\ell,\ell_{\frak p},(\ell_c))}(\beta;\frak p;\frak A;\frak B)
^{\vec w^-_{\frak p},\frak B^-}_{\epsilon_{0},\vec T_{0}}$
is the inverse image of (\ref{transforgottenmark})
by the map (\ref{evatforgottenmark}), the statement
(1) follows.
\par
By choosing $\epsilon_0$ sufficiently small
we can ensure that the image under the map (\ref{evatforgottenmark}) of
each fiber of the map $\frak{forget}_{\frak B,\frak B^-;\vec w_{\frak p},\vec w^-_{\frak p}}$
intersects with the submanifold (\ref{transforgottenmark}) at one point.
Moreover by stability the elements of
$\mathcal M_{k+1,(\ell,\ell^-_{\frak p},(\ell_c))}(\beta;\frak p;\frak A;\frak B^-)_{\epsilon_{0},\vec T_{0}}
$
have no automorphism. The statement (2) follows.
\end{proof}
We next consider a similar but  a slightly different case of transversal
constraint. Namely:
\begin{defn}\label{defn9797}
An element $(\frak Y,u',(\vec w'_c;c\in \frak B))$
of $\mathcal M_{k+1,(\ell,\ell_{\frak p},(\ell_c))}(\beta;\frak p;\frak A;\frak B)_{\epsilon_{0},\vec T_{0}}$ is said to satisfy the
{\it transversal constraint at all additional marked points} if the following holds.
Let $w'_{\frak p,i}$, $i=1,\dots,\ell_{\frak p}$ denote the
$(\ell+1)$-th, \dots, $(\ell+\ell_{\frak p})$-th interior marked points of
$\frak Y$.
We put
$\vec w'_c = (w'_{c,1},\dots,w'_{c,\ell_c})$.
\begin{enumerate}
\item  For all $i = 1,\dots, \ell_{\frak p}$ we have
$
u'(w'_{\frak p,i}) \in \mathcal D_{\frak p,i}.
$
\item
For all $c \in \frak B$ and $i=1,\dots,\ell_{\frak c}$ we have
$
u'(w'_{c,i}) \in \mathcal D_{c,i}.
$
\end{enumerate}
\par
We denote by
$\mathcal M_{k+1,(\ell,\ell_{\frak p},(\ell_c))}(\beta;\frak p;\frak A;\frak B)
^{{\rm trans}}_{\epsilon_{0},\vec T_{0}}$
the set of all elements of the thickened moduli space $\mathcal M_{k+1,(\ell,\ell_{\frak p},(\ell_c))}(\beta;\frak p;\frak A;\frak B)_{\epsilon_{0},\vec T_{0}}$
satisfying  transversal constraint
at all additional marked points.
\end{defn}
\begin{lem}\label{transstratasmf}
The set
$\mathcal M_{k+1,(\ell,\ell_{\frak p},(\ell_c))}(\beta;\frak p;\frak A;\frak B)
^{{\rm trans}}_{\epsilon_{0},\vec T_{0}}$
is a closed subset of our space
$\mathcal M_{k+1,(\ell,\ell_{\frak p},(\ell_c))}(\beta;\frak p;\frak A;\frak B)_{\epsilon_{0},\vec T_{0}}$
and is a strata-wise smooth submanifold of codimension
$2\ell_{\frak p} + 2\sum_{c\in \frak B}\ell_c$.
\end{lem}
\begin{rem}
We note that the map $\rm Glu$ is a homeomorphism onto its image of the thickened moduli space
$\mathcal M_{k+1,(\ell,\ell_{\frak p},(\ell_c))}(\beta;\frak p;\frak A;\frak B)
^{{\rm trans}}_{\epsilon_{0},\vec T_{0}}$.
\end{rem}
\begin{proof}
By Proposition \ref{forgetstillstable} it suffices to consider the case $\frak A = \frak B$.
By the way similar to the proof of Proposition \ref{forgetstillstable}
we define
\begin{equation}\label{evatforgottenmark2}
\mathcal M_{k+1,(\ell,\ell_{\frak p},(\ell_c))}(\beta;\frak p;\frak A)_{\epsilon_{0},\vec T_{0}}
\to X^{\ell_{\frak p}  + \sum_{c\in \frak A}\ell_c}
\end{equation}
that is an evaluation map at all the added marked points.
Note at point $\frak p$, when we perturb the added marked points
$\vec w'_{\frak p}$ and $\vec w'_{\frak c}$ we still obtain an element of
the thickened moduli space.  This is because the map $u_{\frak p}$ is pseudo-holomorphic.
Therefore, the evaluation map
(\ref{evatforgottenmark2}) is transversal to
\begin{equation}\label{transforgottenmark2}
\prod_{i=1}^{\ell_{\frak p}}\mathcal D_{\frak p,i}
\times
\prod_{c \in \frak A }\prod_{i=1}^{\ell_{c}}\mathcal D_{c,i}
\end{equation}
at $\frak p$.
It implies that (\ref{evatforgottenmark2}) is transversal to
(\ref{transforgottenmark2}) everywhere on $\mathcal M_{k+1,(\ell,\ell_{\frak p},(\ell_c))}(\beta;\frak p;\frak A)_{\epsilon_{0},\vec T_{0}}$,
if  $\epsilon_0$ is small.
Since
$\mathcal M_{k+1,(\ell,\ell_{\frak p},(\ell_c))}(\beta;\frak p;\frak A)
^{{\rm trans}}_{\epsilon_{0},\vec T_{0}}$
is the inverse image of (\ref{transforgottenmark2})
by the map (\ref{evatforgottenmark2}), the lemma follows.
\end{proof}
\begin{defn}\label{zerosetofspreddf}
We denote by
$\mathcal M_{k+1,(\ell,\ell_{\frak p},(\ell_c))}(\beta;\frak p;\frak A)
^{{\rm trans}}_{\epsilon_{0},\vec T_{0}} \cap \frak s^{-1}(0)$
the set of all
$(\frak Y,u',(\vec w'_c;c\in \frak A))
\in \mathcal M_{k+1,(\ell,\ell_{\frak p},(\ell_c))}(\beta;\frak p;\frak A)
^{{\rm trans}}_{\epsilon_{0},\vec T_{0}}$
such that $u'$ is pseudo-holomorphic.
\end{defn}
Our space
$\mathcal M_{k+1,(\ell,\ell_{\frak p},(\ell_c))}(\beta;\frak p;\frak A)
^{{\rm trans}}_{\epsilon_{0},\vec T_{0}} \cap \frak s^{-1}(0)$
is a closed subset of the moduli space
$\mathcal M_{k+1,(\ell,\ell_{\frak p},(\ell_c))}(\beta;\frak p;\frak A)
^{{\rm trans}}_{\epsilon_{0},\vec T_{0}}$.
\par
By forgetting all the additional marked points we obtain a map
\begin{equation}\label{forget1}
\frak{forget} :
\mathcal M_{k+1,(\ell,\ell_{\frak p},(\ell_c))}(\beta;\frak p;\frak A)
^{{\rm trans}}_{\epsilon_{0},\vec T_{0}} \cap \frak s^{-1}(0)
\to \mathcal M_{k+1,\ell}(\beta).
\end{equation}
We recall that we have  injective homomorphisms
$$
\aligned
&\Gamma_{\frak p} \to \frak S_{\ell_{\frak p}} \times \prod_{c\in \frak A}\frak S_{\ell_c}, \\
&\Gamma^+_{\frak p} \to \frak S_{\ell} \times \frak S_{\ell_{\frak p}} \times \prod_{c\in \frak A}\frak S_{\ell_c}.
\endaligned
$$
The group $\Gamma^+_{\frak p}$ acts on
$\mathcal M_{k+1,(\ell,\ell_{\frak p},(\ell_c))}(\beta;\frak p;\frak A)_{\epsilon_{0},\vec T_{0}}$
as follows. We regard
$\Gamma^+_{\frak p} \subset \frak S_{\ell} \times \frak S_{\ell_{\frak p}} \times \prod_{c\in \frak A}\frak S_{\ell_c}$.
Then the action of $\Gamma_{\frak p}^+$ on $\mathcal M_{k+1,(\ell,\ell_{\frak p},(\ell_c))}(\beta;\frak p;\frak A)_{\epsilon_{0},\vec T_{0}}$ is by exchanging the interior marked points.
It is easy to see that
$\mathcal M_{k+1,(\ell,\ell_{\frak p},(\ell_c))}(\beta;\frak p;\frak A)
^{{\rm trans}}_{\epsilon_{0},\vec T_{0}}$
is invariant under this action.
Therefore (\ref{forget1}) induces a map
\begin{equation}\label{forget2}
\overline{\frak{forget}} :
\bigg(\mathcal M_{k+1,(\ell,\ell_{\frak p},(\ell_c))}(\beta;\frak p;\frak A)
^{{\rm trans}}_{\epsilon_{0},\vec T_{0}} \cap \frak s^{-1}(0)\bigg)/\Gamma_{\frak p}
\to \mathcal M_{k+1,\ell}(\beta).
\end{equation}
\begin{rem}
The map
(\ref{forget2}) induces a map
$$
\bigg(\mathcal M_{k+1,(\ell,\ell_{\frak p},(\ell_c))}(\beta;\frak p;\frak A)
^{{\rm trans}}_{\epsilon_{0},\vec T_{0}} \cap \frak s^{-1}(0)\bigg)/\Gamma_{\frak p}^+
\to \mathcal M_{k+1,\ell}(\beta)/\frak S_{\ell}.
$$
See Remark \ref{rem216}.
We can use this remark to construct an $\frak S_{\ell}$ invariant Kuranishi structure on
$\mathcal M_{k+1,\ell}(\beta)$.
\end{rem}
\begin{prop}\label{charthomeo}
The map
(\ref{forget2}) is a homeomorphism onto an open neighborhood of $\frak p$.
\end{prop}
\begin{proof}
The geometric intuition behind this proposition is clear.
We will give a detailed proof below for completeness sake.
We first review the definition of the topology of $\mathcal M_{k+1,\ell}(\beta)$ given in
\cite[Definition 10.2, 10.3]{FOn}, \cite[Definition 7.1.39, 7.1.42]{fooo:book1}.
\begin{defn}\label{convdefn}
Let
$\frak p_a = ((\Sigma_a,\vec z_a,\vec z^{\rm int}_a),u_a),
\frak p_{\infty} = ((\Sigma,\vec z,\vec z^{\rm int}),u)
\in \mathcal M_{k+1,\ell}(\beta)$.
We {\it assume} $(\Sigma_a,\vec z_a,\vec z^{\rm int}_a)$ and
$(\Sigma,\vec z,\vec z^{\rm int})$ are stable.
We say that a sequence $((\Sigma_a,\vec z_a,\vec z^{\rm int}_a),u_a)$
{\it stably converges} to
$((\Sigma,\vec z,\vec z^{\rm int}),u)$ and write
$$
\underset{a\to \infty}{\rm lims}\,\, \frak p_a = \frak p_{\infty}
$$
if the following holds.
\begin{enumerate}
\item
We assume $$\lim_{a\to \infty}(\Sigma_a,\vec z_a,\vec z^{\rm int}_a)
= (\Sigma,\vec z,\vec z^{\rm int})$$
in the Deligne-Mumford moduli space $\mathcal M_{k+1,\ell}$.
We take a coordinate at infinity of $(\Sigma,\vec z,\vec z^{\rm int})$.
It determines a diffeomorphism between cores of
$\Sigma_a$ and of $\Sigma$ for large $a$.
\item
For each $\epsilon$ we can extend the core appropriately
so that there exists $a_0$ such that (2),(3) hold
for $a > a_0$.
$$
\vert u_a - u\vert_{C^1(\text{Core})} < \epsilon.
$$
Here we regard $u_a$ and $u$ as maps from the core of $\Sigma_a$ and $\Sigma$ by the above
mentioned diffeomorphism.
\item
The diameter of the image of each connected component of the
 neck region by $u_a$ is smaller than $\epsilon$.
\end{enumerate}
\end{defn}
\begin{defn}
Let $\frak p_a = ((\Sigma_a,\vec z_a,\vec z^{\rm int}_a),u_a)$,
$\frak p_{\infty} = ((\Sigma,\vec z,\vec z^{\rm int}),u)
\in \mathcal M_{k+1,\ell}(\beta)$.
We say that $\frak p_a$ converges to $\frak p_{\infty}$ and write
$$
\lim_{a\to \infty}\frak p_a  = \frak p_{\infty}
$$
if there exist $\ell' \ge 0$ and
$\frak q_a = ((\Sigma_a,\vec z_a,\vec z^{\rm int}_a\cup \vec z^{+,\rm int}_a),u_a)$,
$\frak q_{\infty} = ((\Sigma,\vec z,\vec z^{\rm int}\cup \vec z^{+,\rm int}_{\infty}),u)
\in \mathcal M_{k+1,\ell+\ell'}(\beta)$
such that
\begin{equation}\label{Sconv2}
\underset{a\to \infty}{\rm lims}\,\, \frak q_a = \frak q_{\infty}
\end{equation}
and
\begin{equation}\label{conv2}
\frak{forget}_{(k+1;\ell+\ell'),(k+1;\ell)}(\frak q_a)= \frak p_a,
\quad
\frak{forget}_{(k+1;\ell+\ell'),(k+1;\ell)}(\frak q_{\infty})= \frak p_{\infty}.
\end{equation}
\end{defn}
Here
$$
\frak{forget}_{(k+1;\ell+\ell'),(k+1;\ell)} :
\mathcal M_{k+1,\ell+\ell'}(\beta)
\to \mathcal M_{k+1,\ell}(\beta)
$$
is a map forgetting $(\ell+1)$-st,\dots,$(\ell+\ell')$-st
(interior) marked points (and shrinking the irreducible components that become
unstable. See \cite[p 419]{fooo:book1}.)
\par
Now we prove the following:
\begin{lem}\label{setisopen}
If $\epsilon_0$, $\epsilon_{\frak p_c}$ are sufficiently small, then
the image of (\ref{forget2}) is an open subset of $\mathcal M_{k+1,\ell}(\beta)$.
\end{lem}
\begin{proof}
Let
$$
\frak p' \in \overline{\frak{forget}}
\bigg(\big(\mathcal M_{k+1,(\ell,\ell_{\frak p},(\ell_c))}(\beta;\frak p;\frak A)
^{{\rm trans}}_{\epsilon_{0},\vec T_{0}} \cap \frak s^{-1}(0)\big)/\Gamma_{\frak p}\bigg)
$$
and
$
\frak p_a \in \mathcal M_{k+1,\ell}(\beta)
$
such that
$
\lim_{a\to\infty} \frak p_a = \frak p'.
$
We will prove
$$
\frak p_a \in \overline{\frak{forget}}
\bigg(\big(\mathcal M_{k+1,(\ell,\ell_{\frak p},(\ell_c))}(\beta;\frak p;\frak A)
^{{\rm trans}}_{\epsilon_{0},\vec T_{0}} \cap \frak s^{-1}(0)\big)/\Gamma_{\frak p}\bigg)
$$
for all sufficiently large $a$.
\par
We put $\frak p' = (\frak Y_0,u')$ and
$$
(\frak Y_0\cup \vec w'_{\frak p},u',(\vec w'_c;c\in \frak A))
\in \mathcal M_{k+1,(\ell,\ell_{\frak p},(\ell_c))}(\beta;\frak p;\frak A)
^{{\rm trans}}_{\epsilon_{0},\vec T_{0}} \cap \frak s^{-1}(0).
$$
We also put $\frak p_a = (\frak x_{\frak p_a},u_{\frak p_a})$.
By Definition \ref{convdefn},
there exists $\frak q_a, \frak q_{\infty} \in \mathcal M_{k+1,\ell+\ell'}(\beta)$ such that
(\ref{Sconv2}) holds and
\begin{equation}\label{conv3}
\frak{forget}_{(k+1;\ell+\ell'),(k+1;\ell)}(\frak q_a)= \frak p_a,
\quad
\frak{forget}_{(k+1;\ell+\ell'),(k+1;\ell)}(\frak q_{\infty})= \frak p'.
\end{equation}
Let $\vec z^{+,\rm int}_a \subset \frak x_{\frak q_a}, \vec z^{+,\rm int}_{\infty} \subset \frak x_{\frak q_{\infty}}$
be the interior marked points that are not the marked points of $\frak p_a$ or of $\frak p'$.
By perturbing $\frak q_a$ and $\frak q_{\infty}$ a bit we may assume
\begin{equation}\label{notinDsss}
\aligned
u_{\frak q_a}(z^{+,\rm int}_{a,i})
&\notin
\bigcup_{i=1}^{\ell_{\frak p}} \mathcal D_{\frak p,i}
\cup
\bigcup_{c \in \frak A}\bigcup_{i=1}^{\ell_{c}} \mathcal D_{c,i},\\
u_{\frak q_{\infty}}(z^{+,\rm int}_{{\infty},i})
&\notin
\bigcup_{i=1}^{\ell_{\frak p}} \mathcal D_{\frak p,i}
\cup
\bigcup_{c \in \frak A}\bigcup_{i=1}^{\ell_{c}} \mathcal D_{c,i}.
\endaligned
\end{equation}
We consider the map $\Sigma_{\frak q_{\infty}} \to \Sigma_{\frak p'}$ that
shrinks the irreducible components which become unstable after forgetting $(\ell+1)$-th,\dots,
$(\ell + \ell')$-th marked points $\vec z^{+,\rm int}_{\infty}$ of $\frak x_{\frak q_{\infty}}$.
By (\ref{notinDsss}) none of the points $\vec w'_{\frak p}$, $\vec w'_c$ are
contained in the image of the irreducible components of $\Sigma_{\frak q_{\infty}}$ that we shrink.
Therefore $\vec w'_{\frak p}, \vec w'_c \subset \Sigma_{\frak p'}$ may be regarded as
points of $\Sigma_{\frak q_{\infty}}$.
\par
Then by extending the core if necessary we may assume that $\vec w'_{\frak p}, \vec w'_c$
are in the core of $\Sigma_{\frak q_{\infty}}$. Here we use the
coordinate at infinity that appears in the definition of $\underset{a\to \infty}{\rm lims}\,\, \frak q_a = \frak q_{\infty}$.
\par
We note that
$$
u_{\frak q_{\infty}}(w'_{\frak p,i}) \in \mathcal D_{\frak p,i},
\quad
u_{\frak q_{\infty}}(w'_{c,i}) \in \mathcal D_{c,i}.
$$
We also note that $u_{\frak q_a}$ converges to $u_{\frak q_{\infty}}$ in $C^1$-topology
on the core. Moreover  $u_{\frak q_{\infty}}$ is transversal to $\mathcal D_{\frak p,i}$ (resp. $\mathcal D_{c,i}$)
at $u_{\frak q_{\infty}}(w'_{\frak p,i})$ (resp. $u_{\frak q_{\infty}}(w'_{c,i})$).
Therefore, for sufficiently large $a$ there exist $w'_{a,\frak p,i}, w'_{a,c,i} \in \Sigma_{\frak q_a}$
with the following properties.
\begin{enumerate}
\item $u_{\frak q_a}(w'_{a,\frak p,i}) \in \mathcal D_{\frak p,i}$.
\item $u_{\frak q_a}(w'_{a,c,i}) \in \mathcal D_{c,i}$.
\item $\lim_{a \to \infty} w'_{a,\frak p,i} = w'_{\frak p,i}$.
\item $\lim_{a \to \infty} w'_{a,c,i} = w'_{c,i}$.
\end{enumerate}
Here in the statements (3) and (4) we use the identification of the core of $\Sigma_{\frak q_a}$ and of $\Sigma_{\frak q_{\infty}}$
induced by  the
coordinate at infinity that appears in the definition of $\underset{a\to \infty}{\rm lims}\,\, \frak q_a = \frak q_{\infty}$.
We send $w'_{a,\frak p,i}$ by the map $\Sigma_{\frak q_a} \to \Sigma_{\frak p_a}$ and denote it by the
same symbol. We thus obtain $\vec w'_{a,\frak p} \subset \Sigma_{\frak p_a}$.
The additional marked points $w'_{a,c,i}$ induce $\vec w'_{a,c} \subset \Sigma_{\frak p_a}$ in the same way.
\par
Using (1)-(4) above and the fact that $u_{\frak q_a}$ converges to $u_{\frak q_{\infty}}$ in $C^1$-topology
we can easily show that
$$
(\frak x_{\frak p_a} \cup \vec w'_{a,\frak p}, u_{\frak p_a}, (\vec w'_{a,c};c\in \frak A))
\in \mathcal M_{k+1,(\ell,\ell_{\frak p},(\ell_c))}(\beta;\frak p;\frak A)
^{{\rm trans}}_{\epsilon_{0},\vec T_{0}} \cap \frak s^{-1}(0)
$$
for sufficiently large $a$. Thus we have
$$
\aligned
\frak p_a &= \overline{\frak{forget}}((\frak x_{\frak p_a} \cup \vec w'_{a,\frak p}, u_{\frak p_a}, (\vec w'_{a,c};c\in \frak A))) \\
&\in
\overline{\frak{forget}}
\bigg(\big(\mathcal M_{k+1,(\ell,\ell_{\frak p},(\ell_c))}(\beta;\frak p;\frak A)
^{{\rm trans}}_{\epsilon_{0},\vec T_{0}} \cap \frak s^{-1}(0)\big)/\Gamma_{\frak p}\bigg)
\endaligned
$$
for sufficiently large $a$. The proof of Lemma \ref{setisopen} is complete.
\end{proof}
\begin{lem}\label{injectivitypp}
If $\epsilon_0$ is sufficiently small, then the map
(\ref{forget2}) is injective.
\end{lem}
\begin{proof}
The proof is by contradiction.
We assume that there exists $\epsilon_{0}^{(n)}$ with $\epsilon_{0}^{(n)} \to 0$
as $n\to \infty$, and
\begin{equation}\label{nearpisassump}
\aligned
&(\frak Y_{j;(n),0}\cup \vec w'_{j;(n),\frak p},u'_{j;(n)},(\vec w'_{j;(n),c};c\in \frak A))\\
&\in \mathcal M_{k+1,(\ell,\ell_{\frak p},(\ell_c))}(\beta;\frak p;\frak A)
^{{\rm trans}}_{\epsilon^{(n)}_{0},\vec T_{0}} \cap \frak s^{-1}(0)
\endaligned
\end{equation}
for $j=1,2$.
Here we extend the core of the
coordinate at infinity of $\frak p$ by $\vec R_{(n)} \to \infty$ to define the right hand side of (\ref{nearpisassump}).
We assume
\begin{equation}\label{Y1eqY2}
(\frak Y_{1;{(n)},0},u'_{1;{(n)}}) \sim  (\frak Y_{2;{(n)},0},u'_{2;{(n)}})
\end{equation}
in $\mathcal M_{k+1,\ell}(\beta)$ but
\begin{equation}\label{Y1isY2}
\aligned
&[(\frak Y_{1;{(n)},0}\cup \vec w'_{1;{(n)},\frak p},u'_{1;{(n)}},(\vec w'_{1;{(n)},c};c\in \frak A))]\\
&\ne
[(\frak Y_{2;{(n)},0}\cup \vec w'_{2;{(n)},\frak p},u'_{2;{(n)}},(\vec w'_{2;{(n)},c};c\in \frak A))]
\endaligned
\end{equation}
in $(\big(\mathcal M_{k+1,(\ell,\ell_{\frak p},(\ell_c))}(\beta;\frak p;\frak A)
^{{\rm trans}}_{\epsilon_{0}^{(n)},\vec T_{0}} \cap \frak s^{-1}(0)\big)/\Gamma_{\frak p}$.
We will deduce contradiction.
\par
The condition
(\ref{Y1eqY2}) implies that there exists $v_{(n)} : \Sigma_{\frak Y_{1;{(n)},0}} \to \Sigma_{\frak Y_{2;{(n)},0}}$
with the following properties.
\begin{enumerate}
\item $v_{(n)}$ is a biholomorphic map.
\item $u'_{2;{(n)}} \circ v_{(n)} = u'_{1;{(n)}}$.
\item
$v_{(n)}$ sends $k+1$ boundary marked points and $\ell$ interior marked points of
$\frak Y_{1;{(n)},0}$ to the corresponding marked points of $\frak Y_{2;{(n)},0}$.
\end{enumerate}
We take a coordinate at infinity associated to the stabilization data at $\frak p$.
Then (\ref{nearpisassump}) implies that the core of $\frak Y_{j;{(n)},0}$ ($j=1,2$) is
identified with the extended core $(K_{\rm v}^{\frak p})^{+\vec R_{(n)}}$ of $\frak p$. This identification may not preserve complex
structures but preserves the $k+1$ boundary and $\ell+\ell'$ interior marked points.
Therefore $v_{(n)}$ induces
$$
v_{(n)} : K_{0,{\rm v}}^{\frak p} \to (K_{\rm v}^{\frak p})^{+\vec R_{(n)}}
$$
where $K_{0,{\rm v}}^{\frak p}$ is a compact set such that
$w'_{1;(n),\frak p,i}, w'_{1;(n),c,i} \in K_{0,{\rm v}}^{\frak p}$.
(We may extend the core so that we can find such $K_{0,{\rm v}}^{\frak p}$.)
\par
We may take $\vec R_{(n)} \to \infty$ so that the $u'_{j;(n)}$ image of each of the connected components of the complement of
$(K_{\rm v}^{\frak p})^{+\vec R_{(n)}}$ has diameter $< \epsilon_{0}^{(n)}$.
\par
We consider the complex structure of $\Sigma_{\frak p}$ on $(K_{\rm v}^{\frak p})^{+\vec R_{(n)}}$ and denote it by
$j_{\frak p}$.
Then we have
\begin{equation}\label{almstjequiv}
\lim_{n\to \infty}\Vert (v_{(n)} )_* j_{\frak p} -  j_{\frak p}\Vert_{C^1\left((K_{\rm v}^{\frak p})^{+\vec R^-_{(n)}}\right)} = 0
\end{equation}
where  $\vec R^{-}_{(n)} \to \infty$ is chosen so that
$v_{(n)}((K_{\rm v}^{\frak p})^{+\vec R^-_{(n)}}) \subset (K_{\rm v}^{\frak p})^{+\vec R_{(n)}}$.
\par
On the  other hand by Property (4) above we have
\begin{equation}\label{almostucomp}
\lim_{n\to \infty}\Vert u\circ v_{(n)} - u \Vert_{C^1\left((K_{\rm v}^{\frak p})^{+\vec R^-_{(n)}}\right)}  =0.
\end{equation}
We use (\ref{almstjequiv}) and (\ref{almostucomp}) to prove the following.
\begin{sublem}
After taking a subsequence if necessary, there exists
$v' \in\Gamma_{\frak p}$ such that
$$
\lim_{n\to\infty}\Vert v_{(n)} - v' \Vert_{C^1((K^{\frak p}_{{\rm v}})^{+\vec R})}  = 0
$$
for any $\vec R$.
\end{sublem}
\begin{proof}
Since $v_{(n)}$ is biholomorphic with respect to a pair of complex structures converging
to $(j_{\frak p},j_{\frak p})$,  we can use
Gromov compactness to show that it converges in compact $C^{\infty}$ topology outside finitely many points
after taking a subsequence if necessary. Let $v'$ be the limit.
By the Property (2) above we have $u \circ v' = u$.
\par
On the irreducible component of $\frak x_{\frak p}$ where $u$ is not constant, we use $u \circ v' = u$ together with the
fact that $v_{(n)}$ is biholomorphic to show that there is no bubble on this component.
Namely $v_{(n)}$ converges everywhere on this component.
\par
The irreducible component of $\frak x_{\frak p}$  where $u$ is trivial is stable since $\frak p$ is stable.
We note that $v'$ preserves the marked points of $\frak p$.
It implies that $v'$ is not a constant map on this component.
Then using the fact that   $v_{(n)}$ is biholomorphic we can again
show that there is no bubble on this component.
\par
We thus proved that $v_{(n)}$ converges to $v'$ everywhere. It is then easy to see that $v' \in \Gamma_{\frak p}$.
\end{proof}
By replacing $(\frak Y_{2;(n),0}\cup \vec w'_{2;(n),\frak p},u'_2,(\vec w'_{2;(n),c};c\in \frak A))$
using the action of $v' \in \Gamma_{\frak p}$, we may assume that
\begin{equation}\label{closetoid}
\lim_{n\to \infty}\Vert v_{(n)} - {\rm identity} \Vert_{C^1(K^{\frak p}_{0,{\rm v}})} =0.
\end{equation}
Then, $u'_{1;(n)}(w'_{1;(n),\frak p,i}), u'_{2;(n)}(w'_{2;(n),\frak p,i}) \in \mathcal D_{\frak p,i}$  imply
\begin{equation}
v_{(n)}(w'_{1;(n),\frak p,i},) = w'_{2;(n),\frak p,i}.
\end{equation}
\par
We next take coordinate at infinity associated to the obstruction bundle data centered at $\frak p_c$.
Then we can think of the restriction $v_{(n)} : K^{\frak p_c}_{0,{\rm v}} \to K^{\frak p_c}_{{\rm v}}$,
which satisfies
\begin{equation}\label{closetoid2}
\lim_{n\to \infty}\Vert v_{(n)} - {\rm identity} \Vert_{C^1(K^{\frak p_c}_{{0,\rm v}})} = 0.
\end{equation}
(In fact, we may take $\vec R$ so that for each ${\rm v}\in C^0(\mathcal G_{\frak p_c})$ we have ${\rm v}'
\in C^0(\mathcal G_{\frak p})$ such that
$K_{\rm v}^{\frak p_c} \subset (K_{{\rm v}'}^{\frak p})^{+\vec R}$.)
\par
Then, $u'_{1;(n)}(w'_{1;(n),c,i}), u'_{2;(n)}(w'_{2;(n),c,i}) \in \mathcal D_{c,i}$  imply
\begin{equation}\label{closetoid299}
v_{(n)}(w'_{1;(n),c,i}) = w'_{2;(n),c,i}.
\end{equation}
Property (1),(2) and (\ref{closetoid2}), (\ref{closetoid299})  contradict to
(\ref{Y1isY2}). The proof of Lemma \ref{injectivitypp} is complete.
\end{proof}
\begin{lem}\label{injectivitypp3}
If $\epsilon_0$, $\epsilon_{\frak p_c}$ are sufficiently small, then
(\ref{forget2}) is a homeomorphism onto its image.
\end{lem}
\begin{proof}
It is easy to see that the map (\ref{forget2})  is continuous.
It is injective by Lemma \ref{injectivitypp}.
It suffices to show that the converse is continuous.
The proof of the continuity of the converse is similar to the proof of
Lemma \ref{setisopen}. We however repeat the detail of the proof
for completeness sake.
Let
$$
(\frak x_{\frak p_a} \cup \vec w'_{a,\frak p}, u_{\frak p_a}, (\vec w'_{a,c};c\in \frak A))
\in \mathcal M_{k+1,(\ell,\ell_{\frak p},(\ell_c))}(\beta;\frak p;\frak A)
^{{\rm trans}}_{\epsilon_{0},\vec T_{0}} \cap \frak s^{-1}(0)
$$
and
$$
(\frak x_{\frak p_{\infty}} \cup \vec w'_{{\infty},\frak p}, u_{\frak p_{\infty}}, (\vec w'_{{\infty},c};c\in \frak A))
\in \mathcal M_{k+1,(\ell,\ell_{\frak p},(\ell_c))}(\beta;\frak p;\frak A)
^{{\rm trans}}_{\epsilon_{0},\vec T_{0}} \cap \frak s^{-1}(0).
$$
We put $\frak p_{\infty} = (\frak x_{\frak p_{\infty}},u_{\frak p_{\infty}})$,
$\frak p_{a} = (\frak x_{\frak p_{a}},u_{\frak p_{a}})$
and assume
\begin{equation}
\lim_{a\to \infty} \frak p_{a} = \frak p_{\infty}
\end{equation}
in $\mathcal M_{k+1,\ell}(\beta)$.
\par
By Definition \ref{convdefn},
there exist $\frak q_a, \frak q_{\infty} \in \mathcal M_{k+1,\ell+\ell'}(\beta)$ such that
(\ref{Sconv2}) and
\begin{equation}\label{conv3}
\frak{forget}_{(k+1;\ell+\ell'),(k+1;\ell)}(\frak q_a)= \frak p_a,
\quad
\frak{forget}_{(k+1;\ell+\ell'),(k+1;\ell)}(\frak q_{\infty})= \frak p_{\infty}.
\end{equation}
Let $\vec z^{+,{\rm int}}_a \subset \frak x^{+,{\rm int}}_{\frak q_a}, \vec z^{+,{\rm int}}_{\infty} \subset \frak x_{\frak q_{\infty}}$
be the marked points of $\frak q_a, \frak q_{\infty}$ that  are not marked points of  $\frak p_a$ or of $\frak p_{\infty}$.
By perturbing $\frak q_a$ and $\frak q_{\infty}$ a bit we may assume
(\ref{notinDsss}).
\par
We consider the map $\Sigma_{\frak q_{\infty}} \to \Sigma_{\frak p_{\infty}}$ that
shrinks the components which become unstable after forgetting $(\ell+1)$-th,\dots,
$(\ell + \ell')$-th marked points $\vec z^{+,{\rm int}}_{\infty}$ of $\frak x_{\frak q_{\infty}}$.
By (\ref{notinDsss}) none of the points $\vec w'_{\infty,\frak p}$, $\vec w'_{\infty,c}$ are
contained in the image of the components of $\Sigma_{\frak q_{\infty}}$ that we shrink.
So $\vec w'_{\infty,\frak p}, \vec w'_{\infty,c} \subset \Sigma_{\frak p'}$ may be regarded as
points of $\Sigma_{\frak q_{\infty}}$.
\par
Then by extending the core if necessary we may regard that $\vec w'_{\infty,\frak p}, \vec w'_{\infty,c}$
are in the core of $\Sigma_{\frak q_{\infty}}$. Here we use the
coordinate at infinity that appears in the definition of $\underset{a\to \infty}{\rm lims}\,\, \frak q_a = \frak q_{\infty}$.
\par
We remark that
$
u_{\frak q_{\infty}}(w'_{\infty,c,i}) \in \mathcal D_{c,i}.
$
We also remark that $u_{\frak q_a}$ converges to $u_{\frak q_{\infty}}$ in $C^1$-topology
on the core. Moreover  $u_{\frak q_{\infty}}$ is transversal to $\mathcal D_{\frak p,i}$ (resp. $\mathcal D_{c,i}$)
at $u_{\frak q_{\infty}}(w'_{\infty,\frak p,i})$ (resp. $u_{\frak q_{\infty}}(w'_{\infty,c,i})$).
Therefore, for sufficiently large $a$ there exist $w''_{a,\frak p,i}, w''_{a,c,i} \in \Sigma_{\frak q_a}$
with the following properties.
\begin{enumerate}
\item $u_{\frak q_a}(w''_{a,\frak p,i}) \in \mathcal D_{\frak p,i}$.
\item $u_{\frak q_a}(w''_{a,c,i}) \in \mathcal D_{c,i}$.
\item $\lim_{a \to \infty} w''_{a,\frak p,i} = w'_{\infty,\frak p,i}$.
\item $\lim_{a \to \infty} w''_{a,c,i} = w'_{\infty,c,i}$.
\end{enumerate}
Here in (3)(4) we use the identification of the core of $\Sigma_{\frak q_a}$ and of $\Sigma_{\frak q_{\infty}}$
induced by  the
coordinate at infinity that appears in the definition of $\underset{a\to \infty}{\rm lims}\,\, \frak q_a = \frak q_{\infty}$.
We send $w''_{a,\frak p,i}$ by the map $\Sigma_{\frak q_a} \to \Sigma_{\frak p_a}$ and denote it by the
same symbol. We thus obtain $\vec w''_{a,\frak p} \subset \Sigma_{\frak p_a}$.
The additional marked points $w''_{a,c,i}$ induce $\vec w''_{a,c} \subset \Sigma_{\frak p_a}$ in the same way.
\begin{sublem}\label{lub296lem}
$w''_{a,\frak p,i} = w'_{a,\frak p,i}$ and $w''_{a,c,i} =  w'_{a,c,i}$
if $\epsilon_0$ and $\epsilon_{\frak p_c}$ are small and $a$ is large.
\end{sublem}
\begin{proof}
Note $(\frak x_{\frak p_a}\cup \vec w_{a,\frak p_a},u_{\frak p_a})$
and $(\frak x_{\frak p_{\infty}}\cup \vec w_{{\infty},\frak p_{\infty}},u_{\frak p_{\infty}})$
are both $\epsilon_0$-close to $(\frak x_{\frak p},\vec w_{\frak p},u_{\frak p})$.
Then we can choose $\epsilon_0$ small so that (3) above implies
$$
d(w'_{a,\frak p,i},w''_{a,\frak p,i}) \le 3\epsilon_0
$$
for sufficiently large $a$.
We can also show that
$$
d(w'_{a,c,i},w''_{a,c,i}) \le 3(o(\epsilon_0)+ \epsilon_{\frak p_c})
$$
in the same way. (Here $\lim_{\epsilon_0 \to 0}o(\epsilon_0) =0$.)
On the other hand we have $u_{\frak q_a}(w'_{a,\frak p,i}) \in \mathcal D_{\frak p,i}$,
$u_{\frak q_a}(w'_{a,c,i}) \in \mathcal D_{c,i}$.
They imply the sublemma.
\end{proof}
\begin{rem}\label{rem298}
In the last step we need to assume $\epsilon_{\frak p_c}$ small.
More precisely, when we take $\epsilon_{\frak p_c}$ at the stage of Definition \ref{openW+++} we require the
following.
\par
If
$
d(w'_{c,i},w''_{c,i}) \le  4\epsilon_{\frak p_c}
$,
$w'_{c,i},w''_{c,i} \in \Sigma_{\frak p}$ and $u_{\frak p}(w'_{c,i}) \in \mathcal D_{c,i}$,
$u_{\frak p}(w''_{c,i}) \in \mathcal D_{c,i}$, then $w'_{c,i} = w''_{c,i}$.
\par
We next choose $\epsilon_0$ so small  that the same statement holds for $\frak p_a$, with
$4\epsilon_{\frak p}$ replaced by $3\epsilon_{\frak p_c}$.
\end{rem}
Now (3)(4) above imply
$$
\lim_{a\to\infty}(\frak x_{\frak p_a} \cup \vec w'_{a,\frak p}, u_{\frak p_a}, (\vec w'_{a,c};c\in \frak A))
= (\frak x_{\frak p_{\infty}} \cup \vec w'_{{\infty},\frak p}, u_{\frak p_{\infty}}, (\vec w'_{{\infty},c};c\in \frak A))
$$
in $\mathcal M_{k+1,(\ell,\ell_{\frak p},(\ell_c))}(\beta;\frak p;\frak A)
^{{\rm trans}}_{\epsilon_{0},\vec T_{0}} \cap \frak s^{-1}(0)$
as required.
\end{proof}
The proof of Proposition \ref{charthomeo} is complete.
\end{proof}
\begin{proof}[Proof of Lemma \ref{transpermutelem}]
Lemma \ref{transpermutelem} is actually the same as Lemma \ref{injectivitypp}
except the following  point.
We remark that at the stage when we state Lemma \ref{transpermutelem} we did not prove Theorems \ref{gluethm3} and \ref{exdecayT33}.
In fact, to fix the obstruction bundle $E_c$ we used Lemma \ref{transpermutelem}.
However the argument here is not circular by the following reason.
\par
When we prove Lemma \ref{transpermutelem}, we take an obstruction bundle data centered at $\frak p$ only,
the same point as the one we start the gluing construction. We use the obstruction bundle induced
by this obstruction bundle data
to go through the gluing argument (proof of Theorems \ref{gluethm3} and \ref{exdecayT33}.)
We do not need the conclusion of Lemma \ref{transpermutelem} for the gluing argument.
Then we obtain $\text{\rm Glu}$. We use this map to go through the proof of Lemma \ref{injectivitypp} and prove
Lemma \ref{transpermutelem}.
\end{proof}
\begin{rem}\label{rem:2111}
In Definition \ref{openW+++} we mentioned that we prove open-ness of the set $\frak W^+(\frak p)$
in Subsection \ref{cutting}. Indeed it follows from Lemma \ref{setisopen}.
We remark that open-ness of $\frak W^+(\frak p)$ was used to define the set $\frak C(\frak p)$
and so was used in the proof of Theorems \ref{gluethm3} and \ref{exdecayT33}.
However the argument is not circular by the same reason as we explained
in the proof of Lemma \ref{transpermutelem} above.
\end{rem}
\par\medskip
\section{Construction of Kuranishi chart}
\label{chart}
In Lemma \ref{fraAforget1}, Proposition \ref{forgetstillstable}, Lemma \ref{transstratasmf},
{\it strata-wise} differentiable structures of the spaces
$\mathcal M_{k+1,(\ell,\ell_{\frak p},(\ell_c))}(\beta;\frak p;\frak A;\frak B)
^{\vec w^-_{\frak p},\frak B^-}_{\epsilon_{0},\vec T_{0}}$
and
$\mathcal M_{k+1,(\ell,\ell_{\frak p},(\ell_c))}(\beta;\frak p;\frak A)
^{{\rm trans}}_{\epsilon_{0},\vec T_{0}}$
or maps among them are discussed.
These spaces are actually differentiable manifolds with corners and the maps are differentiable maps
between them.
As we mentioned in \cite[page 771-773]{fooo:book1} this is a consequence of
the exponential decay estimate (Theorems \ref{exdecayT} and \ref{exdecayT33}).
We first discuss this point in detail here.
\par
Let $V_{k+1,(\ell,\ell_{\frak p},(\ell_c))}(\beta;\frak p;\frak A;\frak B;\epsilon_1)$ be as in
(\ref{2193}). We put
\begin{equation}
\aligned
&V_{k+1,(\ell,\ell_{\frak p},(\ell_c))}(\beta;\frak p;\frak A;\frak B;\epsilon_1)^{\vec w^-_{\frak p},\frak B^-}\\
&=
\mathcal M_{k+1,(\ell,\ell_{\frak p},(\ell_c))}(\beta;\frak p;\frak A)
^{\vec w^-_{\frak p},\frak B^-}_{\epsilon_{2},\vec T_{0}}
\cap
V_{k+1,(\ell,\ell_{\frak p},(\ell_c))}(\beta;\frak p;\frak A;\epsilon_1)
\endaligned
\end{equation}
\begin{equation}
\aligned
&V_{k+1,(\ell,\ell_{\frak p},(\ell_c))}(\beta;\frak p;\frak A;\epsilon_1)^{\rm trans}\\
&=
\mathcal M_{k+1,(\ell,\ell_{\frak p},(\ell_c))}(\beta;\frak p;\frak A)
^{\rm trans}_{\epsilon_{2},\vec T_{0}}
\cap
V_{k+1,(\ell,\ell_{\frak p},(\ell_c))}(\beta;\frak p;\frak A;\epsilon_1).
\endaligned
\end{equation}
(See Definitions \ref{def7289}, \ref{defn9797}.)
We note that the right hand side is independent of $\epsilon_2$
and $\vec T_0$ if $\epsilon_1$ is
sufficiently small.
By Proposition \ref{forgetstillstable} (1) and
Lemma \ref{transstratasmf},
$V_{k+1,(\ell,\ell_{\frak p},(\ell_c))}(\beta;\frak p;\frak A;\frak B;\epsilon_1)^{\vec w^-_{\frak p},\frak B^-}$ and
$V_{k+1,(\ell,\ell_{\frak p},(\ell_c))}(\beta;\frak p;\frak A;\epsilon_1)^{\rm trans}$ are $C^m$-submanifolds.
\par
The next proposition says that the thickened moduli spaces
$\mathcal M_{k+1,(\ell,\ell_{\frak p},(\ell_c))}(\beta;\frak p;\frak A;\frak B)
^{\vec w^-_{\frak p},\frak B^-}_{\epsilon_{0},\vec T_{0}}$
and
$\mathcal M_{k+1,(\ell,\ell_{\frak p},(\ell_c))}(\beta;\frak p;\frak A)
^{\rm trans}_{\epsilon_{0},\vec T_{0}}$
are graphs of
the maps ${\rm End}_{\vec w^-_{\frak p},\frak B^-}$ and
${\rm End}_{\vec w^-_{\frak p},\frak B^-}$,
which enjoy exponential decay estimate.
\begin{prop}\label{graphdecaythem}
There exist strata-wise $C^m$-maps
$$
\aligned
{\rm End}_{\vec w^-_{\frak p},\frak B^-} :
&V_{k+1,(\ell,\ell_{\frak p},(\ell_c))}(\beta;\frak p;\frak A;\frak B;\epsilon_1)^{\vec w^-_{\frak p},\frak B^-}
\times (\vec T^{\rm o}_0,\infty] \times ((\vec T^{\rm c}_0,\infty] \times \vec S^1)\\
&\to
V_{k+1,(\ell,\ell_{\frak p},(\ell_c))}(\beta;\frak p;\frak A;\frak B;\epsilon_1)
\endaligned
$$
and
$$
\aligned
{\rm End}_{\rm trans} :
V_{k+1,(\ell,\ell_{\frak p},(\ell_c))}(\beta;\frak p;\frak A;\epsilon_1)^{\rm trans}
&\times (\vec T^{\rm o}_0,\infty] \times ((\vec T^{\rm c}_0,\infty] \times \vec S^1)\\
&\to
V_{k+1,(\ell,\ell_{\frak p},(\ell_c))}(\beta;\frak p;\frak A;\epsilon_0)
\endaligned
$$
with the following properties.
\begin{enumerate}
\item
$\mathcal M_{k+1,(\ell,\ell_{\frak p},(\ell_c))}(\beta;\frak p;\frak A;\frak B)
^{\vec w^-_{\frak p},\frak B^-}_{\epsilon_{0},\vec T_{0}}$
is described by the map ${\rm End}_{\vec w^-_{\frak p},\frak B^-}$
as follows:
$$\aligned
&\mathcal M_{k+1,(\ell,\ell_{\frak p},(\ell_c))}(\beta;\frak p;\frak A;\frak B)
^{\vec w^-_{\frak p},\frak B^-}_{\epsilon_{0},\vec T_{0}}\\
&=
\bigg\{
{\rm Glu}({\rm End}_{\vec w^-_{\frak p},\frak B^-}(\frak q,(\vec T,\vec\theta)), \vec T,\vec\theta)
\\ &\qquad\mid
(\frak q,(\vec T,\vec\theta)) \in
V_{k+1,(\ell,\ell_{\frak p},(\ell_c))}(\beta;\frak p;\frak A;\frak B;\epsilon_1)^{\vec w^-_{\frak p},\frak B^-}
\times (\vec T^{\rm o}_0,\infty] \times ((\vec T^{\rm c}_0,\infty] \times \vec S^1)
\bigg\}.
\endaligned$$
We also have
$$\aligned
&\mathcal M_{k+1,(\ell,\ell_{\frak p},(\ell_c))}(\beta;\frak p;\frak A)
^{\rm trans}_{\epsilon_{0},\vec T_{0}}\\
&=
\bigg\{
{\rm Glu}({\rm End}_{\rm trans}(\frak q,(\vec T,\vec\theta)), \vec T,\vec\theta)
\\ &\qquad\mid
(\frak q,(\vec T,\vec\theta)) \in
V_{k+1,(\ell,\ell_{\frak p},(\ell_c))}(\beta;\frak p;\frak A;\epsilon_1)^{\rm trans}
\times (\vec T^{\rm o}_0,\infty] \times ((\vec T^{\rm c}_0,\infty] \times \vec S^1)
\bigg\}.
\endaligned$$
\item
The maps ${\rm End}_{\vec w^-_{\frak p},\frak B^-}$ and ${\rm End}_{\rm trans}$  enjoy the following exponential decay estimate.
\begin{equation}
\left\Vert \nabla_{\frak q}^n \frac{\partial^{\vert \vec k_{T}\vert}}{\partial T^{\vec k_{T}}}\frac{\partial^{\vert \vec k_{\theta}\vert}}{\partial \theta^{\vec k_{\theta}}} {\rm End}_{\vec w^-_{\frak p},\frak B^-}\right\Vert_{C^0}
< C_{16,m,\vec R}e^{-\delta' (\vec k_{T}\cdot \vec T+\vec k_{\theta}\cdot \vec T^{\rm c})}
\end{equation}
\begin{equation}
\left\Vert \nabla_{\frak q}^n \frac{\partial^{\vert \vec k_{T}\vert}}{\partial T^{\vec k_{T}}}\frac{\partial^{\vert \vec k_{\theta}\vert}}{\partial \theta^{\vec k_{\theta}}} {\rm End}_{\rm trans}\right\Vert_{C^0}
< C_{16,m,\vec R}e^{-\delta' (\vec k_{T}\cdot \vec T+\vec k_{\theta}\cdot \vec T^{\rm c})}
\end{equation}
if $n + \vert \vec k_{T}\vert+\vert \vec k_{\theta}\vert  \le m$.
Here $\nabla_{\frak q}^n$ is a derivation of the direction of the parameter space
$V_{k+1,(\ell,\ell_{\frak p},(\ell_c))}(\beta;\frak p;\frak A;\frak B;\epsilon_1)^{\vec w^-_{\frak p},\frak B^-}$ or of
the parameter space $V_{k+1,(\ell,\ell_{\frak p},(\ell_c))}(\beta;\frak p;\frak A;\epsilon_1)^{\rm trans}$.
\end{enumerate}
\end{prop}
\begin{proof}
We prove the estimate for the case of $\mathcal M_{k+1,(\ell,\ell_{\frak p},(\ell_c))}(\beta;\frak p;\frak A;\frak B)
^{\vec w^-_{\frak p},\frak B^-}_{\epsilon_{0},\vec T_{0}}$.
The other case is entirely similar.
\par
We consider the evaluation map (\ref{evatforgottenmark})
\begin{equation}\label{forgetmoikkai}
\mathcal M_{k+1,(\ell,\ell_{\frak p},(\ell_c))}(\beta;\frak p;\frak A;\frak B)_{\epsilon_{0},\vec T_{0}}
\to X^{(\ell_{\frak p} - \ell^-_{\frak p}) + \sum_{c\in \frak B \setminus \frak B^-}\ell_c}
\end{equation}
and compose it with (\ref{2193})
$$
\aligned
\text{\rm Glu} : V_{k+1,(\ell,\ell_{\frak p},(\ell_c))}(\beta;\frak p;\frak A;\frak B;\epsilon_1) &\times (\vec T^{\rm o}_0,\infty] \times ((\vec T^{\rm c}_0,\infty] \times \vec S^1)
\\
&\to
\mathcal M_{k+1,(\ell,\ell_{\frak p},(\ell_c))}(\beta;\frak p;\frak A;\frak B)_{\epsilon_{0},\vec T_{0}}
\endaligned
$$
to obtain
\begin{equation}\label{form314}
\aligned
{\rm ev}_{\vec w^-_{\frak p},\frak B^-} :
V_{k+1,(\ell,\ell_{\frak p},(\ell_c))}(\beta;\frak p;\frak A;\frak B;\epsilon_1) &\times (\vec T^{\rm o}_0,\infty] \times ((\vec T^{\rm c}_0,\infty]\times \vec S^1) \\
&\to X^{(\ell_{\frak p} - \ell^-_{\frak p}) + \sum_{c\in \frak B \setminus \frak B^-}\ell_c}.
\endaligned
\end{equation}
\begin{lem}\label{evalexdecay11}
The map ${\rm ev}_{\vec w^-_{\frak p},\frak B^-} $ enjoys the following exponential decay estimate.
\begin{equation}
\left\Vert \nabla_{\rho}^n \frac{\partial^{\vert \vec k_{T}\vert}}{\partial T^{\vec k_{T}}}\frac{\partial^{\vert \vec k_{\theta}\vert}}{\partial \theta^{\vec k_{\theta}}} {\rm ev}_{\vec w^-_{\frak p},\frak B^-}\right\Vert_{C^0}
< C_{17,m,\vec R}e^{-\delta' (\vec k_{T}\cdot \vec T+\vec k_{\theta}\cdot \vec T^{\rm c})},
\end{equation}
if $n + \vert \vec k_{T}\vert+\vert \vec k_{\theta}\vert  \le m$.
Here $\nabla_{\rho}^n$ is a derivation of the direction of the parameter space $V_{k+1,(\ell,\ell_{\frak p},(\ell_c))}(\beta;\frak p;\frak A;\frak B;\epsilon_1)^{\vec w^-_{\frak p},\frak B^-}$.
\end{lem}
\begin{proof}
We remark that (\ref{form314}) factors through
\begin{equation}\label{22942}
\aligned
\text{\rm Glures} : &V_{k+1,(\ell,\ell_{\frak p},(\ell_c))}(\beta;\frak p;\frak A;\frak B;\epsilon_1) \times (\vec T^{\rm o}_0,\infty] \times ((\vec T^{\rm c}_0,\infty] \times \vec S^1)\\
&\to \prod_{{\rm v}\in C^0(\mathcal G_{\frak p})}L^2_{m}((K_{\rm v}^{+\vec R},K_{\rm v}^{+\vec R}\cap \partial \Sigma_{\frak p,{\rm v}}),(X,L)).
\endaligned
\end{equation}
In fact we may take $\vec R$ so that all the marked points are in the extended core
$\bigcup_{{\rm v}\in C^0(\mathcal G_{\frak p})}K_{\rm v}^{+\vec R}$.
Therefore the lemma is an immediate consequence of Theorem \ref{exdecayT33}.
\end{proof}
By definition, we have:
$$
V_{k+1,(\ell,\ell_{\frak p},(\ell_c))}(\beta;\frak p;\frak A;\frak B;\epsilon_0)^{\vec w^-_{\frak p},\frak B^-}
=
{\rm ev}_{\vec w^-_{\frak p},\frak B^-}^{-1}
\bigg(
\prod_{i=\ell^-_{\frak p}+1}^{\ell_{\frak p}}\mathcal D_{\frak p,i}
\times
\prod_{c \in \frak B \setminus \frak B^-}\prod_{i=1}^{\ell_{c}}\mathcal D_{c,i}
\bigg).
$$
(See the proof of Proposition \ref{forgetstillstable}.)
Proposition \ref{graphdecaythem} is then a consequence of Lemma \ref{evalexdecay11} and the
implicit function theorem.
\end{proof}
We next change the coordinate of $(\vec T^{\rm o}_0,\infty] \times ((\vec T^{\rm c}_0,\infty] \times \vec S^1)$.
The original coordinates are
$((T_{{\rm e}}),(\theta_{\rm e})) \in (\vec T^{\rm o}_0,\infty] \times ((\vec T^{\rm c}_0,\infty] \times \vec S^1)$.
\begin{defn}
We define
\begin{equation}\label{2294}
\aligned
s_{{\rm e}} = \frac{1}{T_{{\rm e}}} \in \left[0,\frac{1}{T_{{\rm e},0}}\right), \qquad &\text{if ${\rm e} \in C^1_{{\rm o}}(\mathcal G_{\frak p})$},
\\
\frak z_{{\rm e}} = \frac{1}{T_{{\rm e}}}\exp(2\pi\sqrt{-1}\theta_{\rm e}) \in D^2\left(\frac{1}{T_{{\rm e},0}}\right), \qquad &\text{if ${\rm e} \in C^1_{{\rm c}}(\mathcal G_{\frak p})$}.
\endaligned
\end{equation}
We also put $s_{{\rm e}} = 0$ (resp. $\frak z_{{\rm e}} = 0$) if $T_{\rm e} = \infty$.
Here we put
$
 D^2(r) = \{ z \in \C \mid \vert z\vert < r\}.
$
\end{defn}
By this change of coordinates,
$(\vec T^{\rm o}_0,\infty] \times ((\vec T^{\rm c}_0,\infty] \times \vec S^1)$
is identified with
\begin{equation}\label{2295}
\prod_{{\rm e} \in C^1_{{\rm o}}(\mathcal G_{\frak p})}\left[0,\frac{1}{T_{{\rm e},0}}\right)
\times
\prod_{{\rm e} \in C^1_{{\rm o}}(\mathcal G_{\frak p})} D^2\left(\frac{1}{T_{{\rm e},0}}\right).
\end{equation}
\begin{defn}
We denote the right hand side of (\ref{2295}) as $[0,(\vec T^{\rm o}_0)^{-1}) \times D^2((\vec T^{\rm c}_0)^{-1})$.
\end{defn}
\begin{rem}\label{rcstratifiation}
The space
$[0,(\vec T^{\rm o}_0)^{-1}) \times D^2((\vec T^{\rm c}_0)^{-1})$ has a stratification
that is induced by the stratification
$$[0,1/T_{{\rm e},0})
= \{0\} \cup (0,1/T_{{\rm e},0})
$$
and
$$
D^2(1/T_{{\rm e},0}) = \{0\} \cup (D^2(1/T_{{\rm e},0}) \setminus \{0\} ).
$$
This stratification corresponds to the stratification of $(\vec T^{\rm o}_0,\infty] \times ((\vec T^{\rm c}_0,\infty] \times \vec S^1)$ that
we defined before, by the homeomorphism (\ref{2294}).
\end{rem}
We note that $[0,(\vec T^{\rm o}_0)^{-1}) \times D^2((\vec T^{\rm c}_0)^{-1})$
is a smooth manifold with corner. The above stratification is finer
than its stratification associated to the structure of manifold with corner.
\par
We then regard $\text{\rm Glu}$ as a map
\begin{equation}\label{2296}
\aligned
\text{\rm Glu}' : V_{k+1,(\ell,\ell_{\frak p},(\ell_c))}(\beta;\frak p;\frak A;\frak B;\epsilon_1)
&\times [0,(\vec T^{\rm o}_0)^{-1}) \times D^2((\vec T^{\rm c}_0)^{-1})\\
&\to \mathcal M_{k+1,(\ell,\ell_{\frak p},(\ell_c))}(\beta;\frak p;\frak A;\frak B)
_{\epsilon_{0},\vec T_{0}}.
\endaligned
\end{equation}
\begin{cor}
The inverse image
$$
(\text{\rm Glu}')^{-1}\left(
\mathcal M_{k+1,(\ell,\ell_{\frak p},(\ell_c))}(\beta;\frak p;\frak A;\frak B)
^{\vec w^-_{\frak p},\frak B^-}_{\epsilon_{0},\vec T_{0}}\right)
$$
is a $C^m$-submanifold of
$V_{k+1,(\ell,\ell_{\frak p},(\ell_c))}(\beta;\frak p;\frak A;\frak B;\epsilon_1)
\times [0,(\vec T^{\rm o}_0)^{-1}) \times D^2((\vec T^{\rm c}_0)^{-1})$.
It is transversal to the strata of the stratification mentioned in Remark \ref{rcstratifiation}.
\par
The same holds for
$\mathcal M_{k+1,(\ell,\ell_{\frak p},(\ell_c))}(\beta;\frak p;\frak A)
^{\rm trans}_{\epsilon_{0},\vec T_{0}}$.
\end{cor}
This is an immediate consequence of Proposition \ref{graphdecaythem}.
\begin{rem}
We can actually promote this $C^m$ structure to a $C^{\infty}$-structure as we will explain in Subsection \ref{toCinfty}.
The same remark applies to all the constructions of Subsections \ref{chart}-\ref{kstructure}.
\end{rem}
\begin{defn}\label{defVVVVV}
We put
$$
V_{k+1,\ell}((\beta;\frak p;\frak A);\epsilon_{0},\vec T_{0})
=
\mathcal M_{k+1,(\ell,\ell_{\frak p},(\ell_c))}(\beta;\frak p;\frak A)
^{\rm trans}_{\epsilon_{0},\vec T_{0}}
$$
and regard it as a $C^m$-manifold with corner so that $\text{\rm Glu}'$ is a $C^m$-diffeomorphism.
\end{defn}
\begin{lem}
The action of $\Gamma_{\frak p}$ on $V_{k+1,\ell}((\beta;\frak p;\frak A);\epsilon_{0},\vec T_{0})$ is of $C^m$-class.
\end{lem}
\begin{proof}
Note the $\Gamma_{\frak p}$-action on $(\vec T^{\rm o}_0,\infty] \times ((\vec T^{\rm c}_0,\infty] \times \vec S^1)$ is by
exchanging the factors associated to the edges $\rm e$ and by the rotation of the $S^1$ factors. Therefore it becomes a smooth action on
$ [0,(\vec T^{\rm o}_0)^{-1}) \times D^2((\vec T^{\rm c}_0)^{-1})$.
By construction $\text{\rm Glu}'$ is $\Gamma_{\frak p}$-equivariant. The lemma follows.
\end{proof}
The orbifold
$V_{k+1,\ell}((\beta;\frak p;\frak A);\epsilon_{0},\vec T_{0})/\Gamma_{\frak p}$ is a chart of the Kuranishi neighborhood of
$\frak p$ which we define in this section.
Note we may assume that the action of $\Gamma_{\frak p}$ to $V_{k+1,\ell}((\beta;\frak p;\frak A);\epsilon_{0},\vec T_{0})$
is effective, by  increasing the obstruction bundle if necessary.
\par\smallskip
We next define an obstruction bundle.
Recall that we fixed a complex vector space $E_c$ for each $c \in \frak A$.
($E_c = \bigoplus_{{\rm v} \in C^0(\mathcal G_{\frak p_c})}E_{c,{\rm v}}$ and $E_{c,{\rm v}}$ is a subspace of
$
\Gamma_0(\text{\rm Int}\,K^{\rm obst}_{\rm v}; u_{\frak p_c}^*TX \otimes \Lambda^{01})
$.)
By Definition \ref{obbundeldata} (5), $E_c$ carries a $\Gamma_{\frak p_c}$ action.
It follows that $\Gamma_{\frak p} \subset \Gamma_{\frak p_c}$, because
$\frak p\cup \vec w^{\frak p}_c$ is $\epsilon_c$-close to $\frak p_c \cup \vec w_{\frak p_c}$.
Therefore we have a $\Gamma_{\frak p}$-action on
$$
E_{\frak A} = \bigoplus_{c\in \frak A} E_c.
$$
\begin{defn}\label{obbundle1}
The obstruction bundle of our Kuranishi chart is the bundle
\begin{equation}
\frac{\left( V_{k+1,\ell}((\beta;\frak p;\frak A);\epsilon_{0},\vec T_{0}) \times E_{\frak A}\right)}{\Gamma_{\frak p}}
\to
\frac{\left( V_{k+1,\ell}((\beta;\frak p;\frak A);\epsilon_{0},\vec T_{0}) \right)}{\Gamma_{\frak p}}.
\end{equation}
\end{defn}
We next define the  Kuranishi map, that is a section of the obstruction bundle.
Let $\frak q^+= (\frak x_{\frak q},u_{\frak q};(\vec w^{\frak q}_{c}; c\in \frak A)) \in \mathcal M_{k+1,(\ell,\ell_{\frak p},(\ell_c))}(\beta;\frak p;\frak A)
^{\rm trans}_{\epsilon_{0},\vec T_{0}}$. By definition we have
$$
\overline\partial u_{\frak q} \in \mathcal E_{\frak A}(\frak q^+).
$$
By Definition \ref{defEc} we have an isomorphism (\ref{Ivpdefn2211})
\begin{equation}
I^{{\rm v},\frak p_c}_{(\frak y_c,u_c),(\frak x_{\frak q}\cup \vec w^{\frak q}_{c},u_{\frak q})} : E_{\frak p_c,\rm v}(\frak y_c,u_c) \to \Gamma_0(\text{\rm Int}\,K^{\rm obst}_{\rm v};
(u_{\frak q})^*TX \otimes \Lambda^{01}).
\end{equation}
The direct sum of the right hand side over $c \in \frak A$ and ${\rm v} \in C^0(\mathcal G_{\frak p_c})$
is by definition $\mathcal E_{\frak A}(\frak q^+)$.
Sending the element $\overline\partial u_{\frak q}$ by the inverse of
$I^{{\rm v},\frak p_c}_{(\frak y_c,u_c),(\frak x_{\frak q}\cup \vec w^{\frak q}_{c},u_{\frak q})}$
we obtain an element
\begin{equation}\label{2299}
\bigoplus_{c \in \frak A
\atop
{\rm v} \in C^0(\mathcal G_{\frak p_c})}{I^{{\rm v},\frak p_c}_{(\frak y_c,u_c),(\frak x_{\frak q}\cup \vec w^{\frak q}_{c},u_{\frak q})}}^{-1}(\overline\partial u_{\frak q})
\in E_{\frak A}.
\end{equation}
\begin{defn}
We denote the element (\ref{2299}) by $\frak s(\frak q^+)$.
The section $\frak s$ is called the {\it Kuranishi map}.
\end{defn}
\begin{lem}
The section
$\frak s$ defined above is a section of $C^m$-class of the obstruction bundle in Definition \ref{obbundle1} and is
$\Gamma_{\frak p}$-equivariant.
\end{lem}
\begin{proof}
The $\Gamma_{\frak p}$-equivariance is immediate from its construction.
\par
To prove that $\frak s$ is of $C^m$-class,
we first remark that $\frak s$ is extended to the thickened moduli space
$\mathcal M_{k+1,(\ell,\ell_{\frak p},(\ell_c))}(\beta;\frak p;\frak A)_{\epsilon_{0},\vec T_{0}}$
by the same formula.
We consider the composition of $\frak q^+ \mapsto \frak s(\frak q^+)$ with
the map $\text{\rm Glu}'$ (\ref{2296}). Since $K_{\rm v}^{\rm obst}$ lies in the core
this composition factors through \text{\rm Glures} (\ref{22942}).
(Here we identify  $ (\vec T^{\rm o}_0,\infty] \times ((\vec T^{\rm c}_0,\infty] \times \vec S^1)$
with $[0,(\vec T^{\rm o}_0)^{-1}) \times D^2((\vec T^{\rm c}_0)^{-1})$.)
Therefore by Theorem \ref{exdecayT33} we have
\begin{equation}
\left\Vert \nabla_{\rho}^n \frac{\partial^{\vert \vec k_{T}\vert}}{\partial T^{\vec k_{T}}}\frac{\partial^{\vert \vec k_{\theta}\vert}}{\partial \theta^{\vec k_{\theta}}}(\frak s\circ\text{\rm Glu}) \right\Vert_{C^0}
< C_{18,m,\vec R}e^{-\delta' (\vec k_{T}\cdot \vec T+\vec k_{\theta}\cdot \vec T^{\rm c})},
\end{equation}
if $n + \vert \vec k_{T}\vert+\vert \vec k_{\theta}\vert  \le m$.
Therefore $\frak s$ is of $C^{m}$-class.
\end{proof}
We note that the zero set of the section $\frak s$ coincides with the set
$$\mathcal M_{k+1,(\ell,\ell_{\frak p},(\ell_c))}(\beta;\frak p;\frak A)
^{{\rm trans}}_{\epsilon_{0},\vec T_{0}} \cap \frak s^{-1}(0)$$
which we defined in Definition \ref{zerosetofspreddf}.
\begin{defn}\label{defpsi}
We define a local parametrization map
$$
\psi :
\frac{\frak s^{-1}(0)}{\Gamma_{\frak p}} \to \mathcal M_{k+1,\ell}(\beta)
$$
to be the map (\ref{forget2}).
\end{defn}
Proposition \ref{charthomeo} implies that $\psi$ is a homeomorphism to an open neighborhood of $\frak p$.
\par
In summary we have proved the following:
\begin{prop}\label{chartprop}
Let $\frak p \in \mathcal M_{k+1,\ell}(\beta)$.
We take a stabilization data at $\frak p$ and $\frak A \subset \frak C(\frak p)$. ($\frak A \ne \emptyset$.)
Then there exists a Kuranishi neighborhood of $\mathcal M_{k+1,\ell}(\beta)$ at $\frak p$. Namely :
\begin{enumerate}
\item An (effective) orbifold
$V_{k+1,\ell}((\beta;\frak p;\frak A);\epsilon_{0},\vec T_{0})/\Gamma_{\frak p}$.
\item
A vector bundle
$$
\frac{\left( V_{k+1,\ell}((\beta;\frak p;\frak A);\epsilon_{0},\vec T_{0}) \times E_{\frak A}\right)}{\Gamma_{\frak p}}
\to
\frac{\left( V_{k+1,\ell}((\beta;\frak p;\frak A);\epsilon_{0},\vec T_{0}) \right)}{\Gamma_{\frak p}}
$$
on it.
\item
Its section $\frak s$ of $C^m$-class.
\item
A homeomorphism
$$
\psi :
\frac{\frak s^{-1}(0)}{\Gamma_{\frak p}} \to \mathcal M_{k+1,\ell}(\beta)
$$
onto an open neighborhood of $\frak p$ in $\mathcal M_{k+1,\ell}(\beta)$.
\end{enumerate}
\end{prop}
Before closing this subsection, we prove that the evaluation maps
on $\mathcal M_{k+1,\ell}(\beta)$ are extended to our Kuranishi
neighborhood as $C^m$-maps.
\par
We consider the map
\begin{equation}\label{evaluationext}
{\rm ev} :
V_{k+1,\ell}((\beta;\frak p;\frak A);\epsilon_{0},\vec T_{0})
=\mathcal M_{k+1,(\ell,\ell_{\frak p},(\ell_c))}(\beta;\frak p;\frak A)^{\rm trans}_{\epsilon_{0},\vec T_{0}} \to L^{k+1} \times X^{\ell}
\end{equation}
that is the evaluation map at the $0$-th,\dots,$k$-th boundary marked points and
$1$st - $\ell$-th interior marked points.
\begin{lem}
The map
(\ref{evaluationext}) is a $C^m$-map and is $\Gamma_{\frak p}$-equivariant.
\end{lem}
\begin{proof}
We first remark that (\ref{evaluationext}) extends to
 $\mathcal M_{k+1,(\ell,\ell_{\frak p},(\ell_c))}(\beta;\frak p;\frak A)_{\epsilon_{0},\vec T_{0}}$.
Its composition with $\text{\rm Glu}$ factors through \text{\rm Glures} (\ref{22942}).
Therefore by Theorem \ref{exdecayT33} we have
\begin{equation}
\left\Vert \nabla_{\rho}^n \frac{\partial^{\vert \vec k_{T}\vert}}{\partial T^{\vec k_{T}}}\frac{\partial^{\vert \vec k_{\theta}\vert}}{\partial \theta^{\vec k_{\theta}}}({\rm ev}\circ\text{\rm Glu}) \right\Vert_{C^0}
< C_{19,m,\vec R}e^{-\delta' (\vec k_{T}\cdot \vec T+\vec k_{\theta}\cdot \vec T^{\rm c})},
\end{equation}
if $n + \vert \vec k_{T}\vert+\vert \vec k_{\theta}\vert  \le m$.
Therefore ${\rm ev}$ is of $C^{m}$-class. $\Gamma_{\frak p}$ equivariance is immediate from definition.
\end{proof}
\begin{rem}
Proposition \ref{chartprop} holds and can be proved when we replace
$\mathcal M_{k+1,\ell}(\beta)$ by  $\mathcal M_{\ell}^{\rm cl}(\alpha)$. The proof is the same.
\end{rem}
\par\medskip

\section{Coordinate change - I: Change of the stabilization and of the coordinate at infinity}
\label{changeofmarking}

In this subsection and the next, we define coordinate change between Kuranishi neighborhoods we constructed in the last
subsection and prove a version of compatibility of the coordinate changes.
In Subsection \ref{kstructure} we will adjust the sizes of the Kuranishi neighborhoods and of the domains of the coordinate changes
so that they literally satisfy the definition of the Kuranishi structure.
\par
We begin with recalling the facts we have proved so far.
We take a finite set $\{\frak p_c \mid c \in \frak C\} \subset \mathcal M_{k+1,\ell}(\beta)$ and
fix an obstruction bundle data $\frak E_{\frak p_c}$ centered at each $\frak p_c$.
\par
Let $\frak w_{\frak p}$  be a stabilization data at $\frak p\in \mathcal M_{k+1,\ell}(\beta)$.
The stabilization data $\frak w_{\frak p}$ consists of
the following:
\begin{enumerate}
\item
 The additional
marked points $\vec w_{\frak p}$ of $\frak x_{\frak p}$.
\item
The codimension 2 submanifolds $\mathcal D_{\frak p,i}$.
\item
A coordinate at infinity of $\frak x_{\frak p} \cup \vec w_{\frak p}$.
\end{enumerate}
By an abuse of notation we denote the coordinate at infinity also by $\frak w_{\frak p}$ from now on.
Let $\ell_{\frak p} = \# \vec w_{\frak p}$ and $\frak A \subset \frak C(\frak p)$.
We always assume that $\frak A \ne \emptyset$.
\par
By taking a sufficiently small $\epsilon_0$ and sufficiently large $\vec T_{0}$, we obtained a
Kuranishi neighborhood at $\frak p$ by Proposition \ref{chartprop}.
The Kuranishi neighborhood is $V_{k+1,\ell}((\beta;\frak p;\frak A);\epsilon_{0},\vec T_{0})/\Gamma_{\frak p}$.
This Kuranishi neighborhood depends on $\epsilon_0, \vec{T}_0$ as well as $\frak w_{\frak p}$.
During the construction of the coordinate change, we need to shrink this Kuranishi neighborhood several
times. We use a pair of positive numbers $(\frak o,\mathcal T)$ to specify the size as follows.
We consider
\begin{equation}\label{gluemapsss}
\aligned
\text{\rm Glu} : V_{k+1,(\ell,\ell_{\frak p},(\ell_c))}(\beta;\frak p;\frak A;\epsilon_1) &\times
(\vec T^{\rm o}_0,\infty] \times ((\vec T^{\rm c}_0,\infty] \times \vec S^1)
\\
&\to
\mathcal M^{{\frak w}_{\frak p}}_{k+1,(\ell,\ell_{\frak p},(\ell_c))}(\beta;\frak p;\frak A)_{\epsilon_{0},\vec T_{0}}.
\endaligned
\end{equation}
\begin{rem}
Here and hereafter we include the symbol ${\frak w}_{\frak p}$ in the notation of the
thickened moduli space, to show the stabilization data at $\frak p$
that we use to define it.
In fact the dependence of the thickened moduli space on the stabilization data is an important
point to study in this subsection.
\end{rem}
$V_{k+1,(\ell,\ell_{\frak p},(\ell_c))}(\beta;\frak p;\frak A;\epsilon_1)$ is a smooth manifold.
We fix a metric on it. Let
\begin{equation}\label{gluemapsssdom}
B_{\frak o}^{{\frak w}_{\frak p}}(\frak p;V_{k+1,(\ell,\ell_{\frak p},(\ell_c))}(\beta;\frak p;\frak A;\epsilon_1))
\end{equation}
be the $\frak o$ neighborhood of $\frak p$ in this space.
We put $T_{{\rm e},0} = \mathcal T$ for all $\rm e$ and denote it by $\vec{\mathcal T}$.
Since this space is independent of $\epsilon_1$ if $\frak o$ is sufficiently small
compared to $\epsilon_1$ we omit $\epsilon_1$ from the notation.
We consider
\begin{equation}\label{gluemapsssdom2}
B^{{\frak w}_{\frak p}}_{\frak o}(\frak p;V_{k+1,(\ell,\ell_{\frak p},(\ell_c))}(\beta;\frak p;\frak A))\times (\vec{\mathcal T},\infty]
\times ((\vec{\mathcal T},\infty] \times \vec S^1).
\end{equation}

\begin{defn}
We say that $(\frak o,\mathcal T)$ is $\frak w_{\frak p}$ {\it admissible} if
the domain of the map (\ref{gluemapsss}) includes (\ref{gluemapsssdom2}).
We say it is admissible if it is clear which stabilization data we take.
\par
We say $(\frak o,\mathcal T)>(\frak o',\mathcal T')$ if $\frak o>\frak o'$ and $1/\mathcal T > 1/\mathcal T'$.
\end{defn}
\begin{defn}
We denote by
$
V(\frak p,\frak w_{\frak p};(\frak o,\mathcal T);\frak A)
$
the intersection of the image of the set (\ref{gluemapsssdom2}) by the map
(\ref{gluemapsss}) and
$\mathcal M^{{\frak w}_{\frak p}}_{k+1,(\ell,\ell_{\frak p},(\ell_c))}(\beta;\frak p;\frak A)_{\epsilon_{0},\vec T_{0}}^{\rm trans}$.
\par
The restrictions of the obstruction bundle, Kuranishi map, and the map $\psi$ to
$
V(\frak p,\frak w_{\frak p};(\frak o,\mathcal T);\frak A)
$
are written as $\mathcal E_{\frak p,\frak w_{\frak p};(\frak o,\mathcal T);\frak A}$.
and $\frak s_{\frak p,\frak w_{\frak p};(\frak o,\mathcal T);\frak A}$,
$\psi_{\frak p,\frak w_{\frak p};(\frak o,\mathcal T);\frak A}$,
respectively.
\par
They define a Kuranishi neighborhood.
Sometimes we denote
by $
V(\frak p,\frak w_{\frak p};(\frak o,\mathcal T);\frak A)
$ this Kuranishi neighborhood, by an abuse of notation.
\end{defn}
The main result of this subsection is the following.
\begin{prop}\label{prop2117}
Let $\frak w^{(j)}_{\frak p}$, ($j=1,2$) be stabilization data at $\frak p$
and $\frak A \supseteq \frak A^{(1)} \supseteq \frak A^{(2)} \ne \emptyset$.
Suppose $(\frak o^{(1)},\mathcal T^{(1)})$ is $\frak w^{(1)}_{\frak p}$ admissible.
\par
Then there exists $(\frak o^{(2)}_0,\mathcal T^{(2)}_0)$ such that if
$(\frak o^{(2)},\mathcal T^{(2)}) < (\frak o^{(2)}_0,\mathcal T^{(2)}_0)$
then $(\frak o^{(2)},\mathcal T^{(2)})$ is $\frak w_{\frak p}^{(2)}$ admissible and
we have a coordinate change from
$
V(\frak p,\frak w^{(2)}_{\frak p};(\frak o^{(2)},\mathcal T^{(2)});\frak A^{(2)})
$
to
$
V(\frak p,\frak w^{(1)}_{\frak p};(\frak o^{(1)},\mathcal T^{(1)});\frak A^{(1)})
$.
Namely there exists $(\phi_{12},\widehat{\phi}_{12})$ with the
following properties.
\begin{enumerate}
\item
$$
\phi_{12} : V(\frak p,\frak w^{(2)}_{\frak p};(\frak o^{(2)},\mathcal T^{(2)});\frak A^{(2)})
\to V(\frak p,\frak w^{(1)}_{\frak p};(\frak o^{(1)},\mathcal T^{(1)});\frak A^{(1)})
$$
is a $\Gamma_{\frak p}$-equivariant $C^m$ embedding.
\item
$$
\widehat{\phi}_{12} :
\mathcal E_{\frak p,\frak w^{(2)}_{\frak p};(\frak o^{(2)},\mathcal T^{(2)});\frak A^{(2)}}
\to
\mathcal E_{\frak p,\frak w^{(1)}_{\frak p};(\frak o^{(1)},\mathcal T^{(1)});\frak A^{(1)}}
$$
is a $\Gamma_{\frak p}$-equivariant embedding of vector bundles of $C^m$-class
that covers $\phi_{12}$.
\item
The next equality holds.
$$
\frak s_{\frak p,\frak w^{(1)}_{\frak p};(\frak o^{(1)},\mathcal T^{(1)});\frak A^{(1)}}
\circ
\phi_{12}
=
\widehat{\phi}_{12}
\circ
\frak s_{\frak p,\frak w^{(2)}_{\frak p};(\frak o^{(2)},\mathcal T^{(2)});\frak A^{(2)}}.
$$
\item
The next equality holds on
$\frak s_{\frak p,\frak w^{(2)}_{\frak p};(\frak o^{(2)},\mathcal T^{(2)});\frak A^{(2)}}^{-1}(0)$.
$$
\tilde\psi_{\frak p,\frak w^{(1)}_{\frak p};(\frak o^{(1)},\mathcal T^{(1)});\frak A^{(1)}}
\circ
\phi_{12}
=
\tilde\psi_{\frak p,\frak w^{(2)}_{\frak p};(\frak o^{(2)},\mathcal T^{(2)});\frak A^{(2)}}.
$$
Here
$\tilde\psi_{\frak p,\frak w^{(1)}_{\frak p};(\frak o^{(1)},\mathcal T^{(1)});\frak A^{(1)}}$
is the composition of
$\psi_{\frak p,\frak w^{(1)}_{\frak p};(\frak o^{(1)},\mathcal T^{(1)});\frak A^{(1)}}$
and the projection map
$$
V(\frak p,\frak w^{(1)}_{\frak p};(\frak o^{(1)},\mathcal T^{(1)});\frak A^{(1)})
\to
V(\frak p,\frak w^{(1)}_{\frak p};(\frak o^{(1)},\mathcal T^{(1)});\frak A^{(1)})
/
\Gamma_{\frak p}.
$$
The definition of
$\tilde\psi_{\frak p,\frak w^{(2)}_{\frak p};(\frak o^{(2)},\mathcal T^{(2)});\frak A^{(2)}}$
is similar.
\item
Let $\frak q^{(2)} \in V(\frak p,\frak w^{(2)}_{\frak p};(\frak o^{(2)},\mathcal T^{(1)});\frak A^{(2)})$
and $\frak q^{(1)} = \phi_{12}(\frak p)$.
Then the derivative of $\frak s_{\frak p,\frak w^{(2)}_{\frak p};(\frak o^{(2)},\mathcal T^{(2)});\frak A^{(2)}}$
induces an isomorphism
$$
\frac
{T_{\frak q^{(1)}}V(\frak p,\frak w^{(1)}_{\frak p};(\frak o^{(1)},\mathcal T^{(1)});\frak A^{(1)})}
{T_{\frak q^{(2)}}V(\frak p,\frak w^{(2)}_{\frak p};(\frak o^{(2)},\mathcal T^{(2)});\frak A^{(2)})}
\cong
\frac
{\left(\mathcal E_{\frak p,\frak w^{(1)}_{\frak p};(\frak o^{(1)},\mathcal T^{(1)});\frak A^{(1)}}\right)_{\frak q^{(1)}}}
{\left(\mathcal E_{\frak p,\frak w^{(2)}_{\frak p};(\frak o^{(2)},\mathcal T^{(2)});\frak A^{(2}}\right)_{\frak q^{(2)}}}.
$$
\end{enumerate}
\end{prop}
\begin{proof}
We divide the proof into several cases.
\par\medskip
\noindent{\bf Case 1}: The case $\vec w_{\frak p}^{(1)} = \vec w_{\frak p}^{(2)}$,
$\mathcal D^{(1)}_{\frak p,i} = \mathcal D^{(2)}_{\frak p,i}$ and $\frak A^{(1)} = \frak A^{(2)}$.
\par
This is the case when only the coordinate at infinity $\frak w_{\frak p}^{(1)}$ is different from $\frak w_{\frak p}^{(2)}$.
A part of the data of the coordinate at infinity is a fiber bundle
(\ref{fibrationsigma}) that is:
\begin{equation}\label{fibrationsigma2}
\pi : \frak M^{(j)}_{(\frak x_{\frak p} \cup \vec w_{\frak p})_{\rm v}}
\to \frak V^{(j)}((\frak x_{\frak p} \cup \vec w_{\frak p})_{\rm v})
\end{equation}
where $\frak V^{(j)}((\frak x_{\frak p} \cup \vec w_{\frak p})_{\rm v})$
is a neighborhood of $(\frak x_{\frak p} \cup \vec w_{\frak p})_{\rm v}$
in the Deligne-Mumford moduli space
$\mathcal M_{k_{\rm v}+1,\ell_{\rm v}}$ or
$\mathcal M_{\ell_{\rm v}}^{\rm cl}$.
(${\rm v} \in C^0(\mathcal G_{\frak x_{\frak p} \cup \vec w_{\frak p}})$.)
We choose $\frak V^{(2)-}((\frak x_{\frak p} \cup \vec w_{\frak p})_{\rm v})
\subset \frak V^{(j)}((\frak x_{\frak p}\cup \vec w_{\frak p})_{\rm v})$
an open neighborhood of $(\frak x_{\frak p}\cup \vec w_{\frak p})_{\rm v}$ so that
\begin{equation}
\frak V^{(2)-}((\frak x_{\frak p} \cup \vec w_{\frak p})_{\rm v})\subset
\frak V^{(1)}((\frak x_{\frak p} \cup \vec w_{\frak p})_{\rm v}).
\end{equation}
We put
$\frak M^{(2) - }_{(\frak x_{\frak p}\cup \vec w_{\frak p})_{\rm v}}
=
\pi^{-1}(\frak V^{(2)-}((\frak x_{\frak p} \cup \vec w_{\frak p})_{\rm v}))$.
Then there exists a unique bundle map
$$
\Phi_{12} : \frak M^{(2) - }_{(\frak x_{\frak p}\cup \vec w_{\frak p})_{\rm v}}
\to \frak M^{(1)}_{(\frak x_{\frak p}\cup \vec w_{\frak p})_{\rm v}}
$$
that preserves the marked points and is a fiberwise biholomorphic map.
This is because of the stability.
By extending the core of ${\frak w}^{(2)}_{\frak p}$
we may assume
\begin{equation}\label{extendcore2banme}
\Phi_{12}(\frak K^{(2) - }_{(\frak x_{\frak p}\cup \vec w_{\frak p})_{\rm v}} )
\supset
(\frak K^{(1)}_{(\frak x_{\frak p}\cup \vec w_{\frak p})_{\rm v}})
\cap
\pi^{-1}(\frak V^{(2)-}((\frak x_{\frak p} \cup \vec w_{\frak p})_{\rm v})).
\end{equation}
\begin{lem}\label{lemma21118}
Let $\epsilon_0$ and $\mathcal T^{(1)}$ be given, then
there exist $\epsilon'_0$, $\mathcal T^{(2)}$ such that
\begin{equation}\label{lem2118formula}
\mathcal M^{\frak w_{\frak p}^{(2) -}}_{k+1,(\ell,\ell_{\frak p},(\ell_c))}(\beta;\frak p;\frak A)
_{\epsilon'_{0},\vec{\mathcal T}^{(2)}}
\subset
\mathcal M^{\frak w_{\frak p}^{(1)}}_{k+1,(\ell,\ell_{\frak p},(\ell_c))}(\beta;\frak p;\frak A)
_{\epsilon_{0},\vec{\mathcal T}^{(1)}}.
\end{equation}
\end{lem}
Here we define $\frak w_{\frak p}^{(2)-}$ from $\frak w_{\frak p}^{(2)}$ by shrinking
$\frak V^{(2)}((\frak x_{\frak p} \cup \vec w_{\frak p})_{\rm v})$
to
$\frak V^{(2)-}((\frak x_{\frak p} \cup \vec w_{\frak p})_{\rm v})$
and extending the core so that (\ref{extendcore2banme}) is satisfied and use it
to define the left hand side.
\begin{proof}
Since the equation (\ref{mainequationformulamod}) is independent of the
stabilization data at $\frak p$,
it suffices to show
$$
\frak U_{k+1,(\ell,\ell_{\frak p},(\ell_c))}^{\frak w_{\frak p}^{(2)-}}(\beta;\frak p)_{\epsilon'_{0},\vec{\mathcal T}^{(2)}}
\subseteq
\frak U^{\frak w_{\frak p}^{(1)}}_{k+1,(\ell,\ell_{\frak p},(\ell_c))}(\beta;\frak p)_{\epsilon_{0},\vec{\mathcal T}^{(1)}}.
$$
Here the meaning of the symbol `$(2)-$' and `$(1)$' is similar to (\ref{lem2118formula}).
\par
An element of
$
\frak U_{k+1,(\ell,\ell_{\frak p},(\ell_c))}^{\frak w_{\frak p}^{(2)-}}(\beta;\frak p)_{\epsilon'_{0},\vec{\mathcal T}^{(2)}}
$ is
$(\frak Y_0 \cup \vec w'_{\frak p},u',(\vec w'_c))$.
Let us check that it satisfies (1)-(4) of Definition \ref{defn251} applied to
$\frak U^{\frak w_{\frak p}^{(1)}}_{k+1,(\ell,\ell_{\frak p},(\ell_c))}(\beta;\frak p)_{\epsilon_{0},\vec{\mathcal T}^{(1)}}$.
\par
(1) is obvious. (2) follows from (\ref{extendcore2banme}). (4) is also obvious.
\par
We will prove (3).
We note that $\frak p$ is $\epsilon_0$ close to $\frak p$ itself by our choice.
So the  diameter of the $u_{\frak p}$ image of each connected component of
the neck region (with respect to ${\frak w}^{(1)}$) is smaller than $\epsilon_0$.
We take $\epsilon'_0$ so that
the diameter of the $u_{\frak p}$ image of each connected component of
the neck region (with respect to ${\frak w}^{(1)}$) is smaller than $\epsilon_0 - 2\epsilon'_0$.
Now since the $C^0$ distance between $u'$ and $u_{\frak p}$
on the core of ${\frak w}^{(2)}$ is small than $\epsilon'_0$,
$$\aligned
&u'\big(
\text{e-th neck with respect to $\frak w^{(1)}_{\frak p}$}
\big)\\
&\subset
\text{$\epsilon'_0$ neighborhood of
$u_{\frak p}\big(
\text{e-th neck with respect to $\frak w^{(2)}_{\frak p}$}
\big)$.}
\endaligned$$
(3) follows.
\end{proof}
Using the fact that $\mathcal D^{(1)}_{\frak p,i} = \mathcal D^{(2)}_{\frak p,i}$,
Lemma \ref{lemma21118} implies
\begin{equation}\label{lem21110formula}
\mathcal M^{\frak w_{\frak p}^{(2)-}}_{k+1,(\ell,\ell_{\frak p},(\ell_c))}(\beta;\frak p;\frak A)
^{\rm trans}_{\epsilon'_{0},\vec{\mathcal T}^{(2)}}
\subset
\mathcal M^{\frak w_{\frak p}^{(1)}}_{k+1,(\ell,\ell_{\frak p},(\ell_c))}(\beta;\frak p;\frak A)
^{\rm trans}_{\epsilon_{0},\vec{\mathcal T}^{(1)}}.
\end{equation}
Let
\begin{equation}\label{gluemapsss1}
\aligned
\text{\rm Glu}^{(1)} :
&B^{\frak w_{\frak p}^{(1)}}_{\frak o^{(1)}}(\frak p;V_{k+1,(\ell,\ell_{\frak p},(\ell_c))}(\beta;\frak p;\frak A)) \\
&\times (\vec{\mathcal T}^{(1)},\infty]
\times ((\vec{\mathcal T}^{(1)},\infty] \times \vec S^1)
\to
\mathcal M^{\frak w_{\frak p}^{(1)}}_{k+1,(\ell,\ell_{\frak p},(\ell_c))}(\beta;\frak p;\frak A)_{\epsilon_{0},\vec{\mathcal T}^{(1)}}
\endaligned
\end{equation}
and
$$
\aligned
\text{\rm Glu}^{(2)-} :
&B^{\frak w_{\frak p}^{(2)-}}_{\frak o^{(2)}}(\frak p;V_{k+1,(\ell,\ell_{\frak p},(\ell_c))}(\beta;\frak p;\frak A)) \\
&\times (\vec{\mathcal T}^{(2)},\infty]
\times ((\vec{\mathcal T}^{(2)},\infty] \times \vec S^1)
\to
\mathcal M^{\frak w_{\frak p}^{(2)-}}_{k+1,(\ell,\ell_{\frak p},(\ell_c))}(\beta;\frak p;\frak A)_{\epsilon_{0},\vec{\mathcal T}^{(2)}}
\endaligned
$$
be appropriate restrictions of (\ref{gluemapsss}).
Its image is an open neighborhood of $\frak p \cup \vec w_{\frak p}$.
Therefore there exists $(\frak o^{(2)}_0, \mathcal T^{(2)}_0)$
such that for any $(\frak o^{(2)},\mathcal T^{(2)}) < (\frak o^{(2)}_0,\mathcal T^{(2)}_0)$ we have
\begin{equation}\label{gluemapsss2}
\aligned
&\text{\rm Glu}^{(2)-}\big(B^{\frak w_{\frak p}^{(2)-}}_{\frak o^{(2)}}(\frak p;V_{k+1,(\ell,\ell_{\frak p},(\ell_c))}(\beta;\frak p;\frak A)) \times (\vec{\mathcal T}^{(2)},\infty] \times ((\vec{\mathcal T}^{(2)},\infty] \times \vec S^1)\big)
\\
&\subset
\text{\rm Glu}^{(1)}\big(B^{\frak w_{\frak p}^{(1)}}_{\frak o^{(1)}}(\frak p;V_{k+1,(\ell,\ell_{\frak p},(\ell_c))}(\beta;\frak p;\frak A)) \times (\vec{\mathcal T}^{(1)},\infty] \times ((\vec{\mathcal T}^{(1)},\infty] \times \vec S^1)\big).
\endaligned
\end{equation}
This in turn implies
$$
V(\frak p,\frak w^{(2)}_{\frak p};(\frak o^{(2)},\mathcal T^{(2)});\frak A)
\subset V(\frak p,\frak w^{(1)}_{\frak p};(\frak o^{(1)},\mathcal T^{(1)});\frak A).
$$
Let $\phi_{12}$ be this natural inclusion.
\begin{lem}\label{2120lem}
 $\phi_{12}$  is a $C^m$-map.
\end{lem}
\begin{proof}
Let
$$
\aligned
&\hat V(\frak p,\frak w^{(2)}_{\frak p};(\frak o^{(2)},\mathcal T^{(2)});\frak A)\\
&\subset
B^{\frak w_{\frak p}^{(2)-}}_{\frak o^{(2)}}(\frak p;V_{k+1,(\ell,\ell_{\frak p},(\ell_c))}(\beta;\frak p;\frak A)) \times (\vec{\mathcal T}^{(2)},\infty] \times ((\vec{\mathcal T}^{(2)},\infty] \times \vec S^1)
\endaligned$$
be the inverse image of
$V(\frak p,\frak w^{(2)-}_{\frak p};(\frak o^{(2)},\mathcal T^{(2)});\frak A)$ by $\text{\rm Glu}^{(2)-}$
and let
$$
\aligned
&\hat V(\frak p,\frak w^{(1)}_{\frak p};(\frak o^{(1)},\mathcal T^{(1)});\frak A)\\
&\subset
B^{\frak w_{\frak p}^{(1)}}_{\frak o^{(1)}}(\frak p;V_{k+1,(\ell,\ell_{\frak p},(\ell_c))}(\beta;\frak p;\frak A)) \times (\vec{\mathcal T}^{(1)},\infty] \times ((\vec{\mathcal T}^{(1)},\infty] \times \vec S^1)
\endaligned$$
be the inverse image of $V(\frak p,\frak w^{(1)}_{\frak p};(\frak o^{(1)},\mathcal T^{(1)});\frak A)$
by  $\text{\rm Glu}^{(1)}$.
\par
We consider the maps
$$\aligned
&B^{\frak w_{\frak p}^{(1)}}_{\frak o^{(1)}}(\frak p;V_{k+1,(\ell,\ell_{\frak p},(\ell_c))}(\beta;\frak p;\frak A))
\to \prod_{{\rm v}\in C^0(\mathcal G_{\frak p})}\frak V^{(1)}((\frak x_{\frak p} \cup \vec w_{\frak p})_{\rm v})
\\
&B^{\frak w_{\frak p}^{(2)-}}_{\frak o^{(2)}}(\frak p;V_{k+1,(\ell,\ell_{\frak p},(\ell_c))}(\beta;\frak p;\frak A))
\to
\prod_{{\rm v}\in C^0(\mathcal G_{\frak p})}\frak V^{(2)-}((\frak x_{\frak p} \cup \vec w_{\frak p})_{\rm v})
\endaligned$$
that forget the maps.  (Namely it sends $(\frak y,u')$ to $\frak y$.)
\par
We then define a map
\begin{equation}\label{evaluandTfac}
\aligned
&\frak F^{(1)} : \hat V(\frak p,\frak w^{(1)}_{\frak p};(\frak o^{(1)},\mathcal T^{(1)});\frak A)
\\\to
&\prod_{{\rm v}\in C^0(\mathcal G_{\frak p})} C^m((K_{{\rm v},(1)}^{+\vec R_{(1)}},K_{{\rm v},(1)}^{+\vec R_{(1)}}\cap
\partial \Sigma_{{\rm v},(1)}),(X,L))\\
&\times \prod_{{\rm v}\in C^0(\mathcal G_{\frak p})}\frak V^{(1)}((\frak x_{\frak p} \cup \vec w_{\frak p})_{\rm v})\\
&\times (\vec{\mathcal T}^{(1)},\infty] \times ((\vec{\mathcal T}^{(1)},\infty] \times \vec S^1).
\endaligned
\end{equation}
Here the first factor is induced  by the map
$$
\mathcal M^{\frak w_{\frak p}^{(1)}}_{k+1,(\ell,\ell_{\frak p},(\ell_c))}(\beta;\frak p;\frak A)_{\epsilon_{0},\vec T_{0}}
\to L^2_{m+10}((K_{{\rm v},(1)}^{+\vec R_{(1)}},K_{{\rm v},(1)}^{+\vec R_{(1)}}\cap \partial \Sigma_{{\rm v},(1)}),(X,L))
$$
that is the map $\text{\rm Glu}^{(1)}$
followed by the restriction of the domain  to the
core $K_{{\rm v},(1)}^{+\vec R_{(1)}}$. (See (\ref{2195form}).)
(We put the symbol $(1)$ in $K_{{\rm v},(1)}^{+\vec R_{(1)}})$ to clarify that this core is induced by
$\frak w_{\frak p}^{(1)}$.)
We chose $T_{{\rm e},0}$ so that the gluing construction works for $L^2_{10m+1}$.
(See the end of Section \ref{glueing}.)
The second and the third factors are the obvious projections.
The map $\frak F^{(1)}$ is a $C^m$
embedding of the $C^m$ manifold $\hat V(\frak p,\frak w^{(2)}_{\frak p};(\frak o^{(2)},\mathcal T^{(2)});\frak A)$,
with corners.
\par
We also consider a similar embedding
\begin{equation}\label{evaluandTfac2}
\aligned
&\frak F^{(2)} : \hat V(\frak p,\frak w^{(2)}_{\frak p};(\frak o^{(2)},\mathcal T^{(2)});\frak A)
\\
\to
&\prod_{{\rm v}\in C^0(\mathcal G_{\frak p})} C^{2m}((K_{{\rm v},(2)}^{+\vec R_{(2)}}),K_{{\rm v},(2)}^{+\vec R_{(2)}}
\cap \partial \Sigma_{{\rm v},(2)}),(X,L))\\
&\times \prod_{{\rm v}\in C^0(\mathcal G_{\frak p})}\frak V^{(2)}((\frak x_{\frak p} \cup \vec w_{\frak p})_{\rm v})\\
&\times (\vec{\mathcal T}^{(2)},\infty] \times ((\vec{\mathcal T}^{(2)},\infty] \times \vec S^1).
\endaligned
\end{equation}
\par
We denote by $\frak X(1,m)$ the right hand side of (\ref{evaluandTfac}) and
by $\frak X(2,2m)$ the right hand side of (\ref{evaluandTfac2}).
\par
We next study the change of parametrization of  the core.
Let us use the notation in Proposition \ref{reparaexpest}.
For $(\rho,\vec T,\vec \theta) \in \prod_{{\rm v} \in C^0(\mathcal G_{\frak p})}\frak V^{(2)}((\frak x_{\frak p} \cup \vec w_{\frak p})_{\rm v})
\times (\vec{\mathcal T}^{(2)},\infty] \times ((\vec{\mathcal T}^{(2)},\infty] \times \vec S^1)$
we have a map
$$
\frak v_{\rho,\vec T,\vec \theta} : \Sigma_{\vec T,\vec \theta}^{(2)} \to \Sigma_{\overline{\Phi}_{12}({\rho,\vec T,\vec \theta})}^{(1)}.
$$
The source $\Sigma_{\vec T,\vec \theta}^{(2)}$  is obtained  using the coordinate at infinity
$\frak w_{\frak p}^{(2)}$ and the target $\Sigma_{\overline{\Phi}_{12}({\rho,\vec T,\vec \theta})}^{(1)}$
is obtained  using the coordinate at infinity
$\frak w_{\frak p}^{(1)}$.
We may assume that
$$
\frak v_{\rho,\vec T,\vec \theta}(K_{{\rm v},(2)}^{+\vec R_{(2)}})
\subset K_{{\rm v},(1)}^{+\vec R_{(1)}}.
$$
We then define a map
$$
\frak H_{12} : \frak X(2,2m) \to \frak X(1,m)
$$
by the formula
\begin{equation}\label{defnH12}
\frak H_{12}(u,(\rho,\vec T,\vec \theta)) = (u\circ \frak v_{(\rho,\vec T,\vec \theta)},\overline{\Phi}_{12}(\rho,\vec T,\vec \theta)).
\end{equation}
\begin{sublem}\label{sublem1}
$\frak H_{12}$ is a $C^m$-map.
\end{sublem}
\begin{proof}
By Proposition \ref{changeinfcoorprop}, the map $\overline{\Phi}_{12}$ is a $C^m$ diffeomorphism.
Therefore the second and the third factors of $\frak H_{12}$ is a $C^m$-map.
The first factor is of $C^m$-class because of Proposition \ref{reparaexpest} and a well-known fact that the map
$
C^{m}(M_1,M_2) \times C^{2m}(M_2,M_3)  \to C^{m}(M_1,M_3)
$
given by
$
(v,u) \mapsto u\circ v
$
is a $C^m$ map.
\end{proof}
On the other hand we have:
\begin{sublem}\label{sublem2}
$$
\frak H_{12} \circ \frak F^{(2)} = \frak F^{(1)} \circ \phi_{12}.
$$
\end{sublem}
This is immediate from the construction.
\par
Since $\frak F^{(2)}$ and $\frak F^{(1)}$ are both $C^m$ embeddings,
Sublemmas \ref{sublem1} and \ref{sublem2} imply  Lemma \ref{2120lem}.
\end{proof}
The map $\phi_{12}$ is obviously $\Gamma_{\frak p}$ equivariant.
We then define
$$
\aligned
\widehat{\phi}_{12} = \phi_{12} \times {\rm identity} :
&V(\frak p,\frak w^{(2)}_{\frak p};(\frak o^{(2)},\mathcal T^{(2)});\frak A) \times \bigoplus_{c\in \frak A} E_c\\
&\subset V(\frak p,\frak w^{(1)}_{\frak p};(\frak o^{(1)},\mathcal T^{(1)});\frak A)\times \bigoplus_{c\in \frak A} E_c.
\endaligned
$$
Conditions (2)-(5) are trivial to verify.
It  also follows that the maps obtained are $\Gamma_{\frak p}$-equivariant.
(During the proof of Proposition \ref{prop2117}, the $\Gamma_{\frak p}$-equivariance is always trivial
to prove. So we do not mention it any more.)
\par\medskip
\noindent{\bf Case 2}: The case ${\frak w}_{\frak p}^{(1)} = {\frak w}_{\frak p}^{(2)}$ and $\frak A^{(1)} \ne \frak A^{(2)}$.
\par
Assume that $\frak B \supseteq \frak A^{(1)} \supset \frak A^{(2)}$ ($\frak B\subseteq \frak C(\frak p)$).
If we regard
$$
V(\frak p,\frak w_{\frak p};(\frak o^{(1)},\mathcal T^{(1)});\frak A^{(1)})
\subset
\mathcal M^{{\frak w}_{\frak p}^{(1)}}_{k+1,(\ell,\ell_{\frak p},(\ell_c))}(\beta;\frak p;\frak A^{(1)};\frak B)
^{\rm trans}_{\epsilon_{0},\vec{\mathcal T}^{(1)}},
$$
then we may also regard
$$
V(\frak p,\frak w_{\frak p};(\frak o^{(1)},\mathcal T^{(1)});\frak A^{(2)})
\subset
\mathcal M^{{\frak w}_{\frak p}^{(1)}}_{k+1,(\ell,\ell_{\frak p},(\ell_c))}(\beta;\frak p;\frak A^{(1)};\frak B)
^{\rm trans}_{\epsilon_{0},\vec{\mathcal T}^{(1)}}.
$$
Moreover
\begin{equation}\label{2315emb}
V(\frak p,\frak w_{\frak p};(\frak o^{(1)},\mathcal T^{(1)});\frak A^{(2)})
\subset
V(\frak p,\frak w_{\frak p};(\frak o^{(1)},\mathcal T^{(1)});\frak A^{(1)}).
\end{equation}
We can show that (\ref{2315emb}) is a $C^m$-map in the same way as the proof of Lemma
\ref{2120lem}.
(Actually the proof is easier since there is no coordinate change of the source
and so $\frak H_{12}$  is the identity map in the situation of Case 2.)
\par
Furthermore an element $(\frak Y,u',(\vec w'_{a,c};c\in \frak B))$ of
$V(\frak p,\frak w_{\frak p};(\frak o^{(1)},\mathcal T^{(1)});\frak A^{(1)})$ is in
$V(\frak p,\frak w_{\frak p};(\frak o^{(1)},\mathcal T^{(1)});\frak A^{(2)})$ if and only if
\begin{equation}\label{2315+1}
\frak s_{\frak p,\frak w_{\frak p};(\frak o^{(1)},\mathcal T^{(1)});\frak A^{(1)}}(\frak Y,u',(\vec w'_{a,c};c\in \frak B))
= \overline\partial u'
\in \mathcal E_{\frak p,\frak w_{\frak p};(\frak o^{(1)},\mathcal T^{(1)});\frak A^{(2)}}.
\end{equation}
We put $\frak q^+ = (\frak Y,u',(\vec w'_{a,c};c\in \frak B))$.
By Lemmas \ref{nbdregmaineq} and \ref{transstratasmf}, $d_{\frak q^+}\frak s$
induces an isomorphism:
$$
\frac
{T_{\frak q^+}V(\frak p,\frak w_{\frak p};(\frak o^{(1)},\mathcal T^{(1)});\frak A^{(1)})}
{T_{\frak q^+}V(\frak p,\frak w_{\frak p};(\frak o^{(1)},\mathcal T^{(1)});\frak A^{(2)})}
\cong
\frac
{\left(\mathcal E_{\frak p,\frak w_{\frak p};(\frak o^{(1)},\mathcal T^{(1)});\frak A^{(1)}}\right)_{\frak q^{+}}}
{\left(\mathcal E_{\frak p,\frak w_{\frak p};(\frak o^{(1)},\mathcal T^{(1)});\frak A^{(2)}}\right)_{\frak q^{+}}}.
$$
We have thus obtained a coordinate change in this case.
\par\medskip
The other two cases are as follows.
\par\smallskip
\noindent{\bf Case 3}: The case $\vec{w}_{\frak p}^{(1)} \subset \vec{w}_{\frak p}^{(2)}$ and $\frak A^{(1)} = \frak A^{(2)}$.
The stabilization data  ${\frak w}_{\frak p}^{(1)}$ is induced from ${\frak w}_{\frak p}^{(2)}$.
\par\smallskip
\noindent{\bf Case 4}: The case $\vec{w}_{\frak p}^{(1)} \supset \vec{w}_{\frak p}^{(2)}$ and $\frak A^{(1)} = \frak A^{(2)}$.
The stabilization data  ${\frak w}_{\frak p}^{(2)}$ is induced from ${\frak w}_{\frak p}^{(1)}$.
\par\smallskip
Let  us explain the notion that `stabilization data  ${\frak w}_{\frak p}^{(1)}$ is induced from ${\frak w}_{\frak p}^{(2)}$.'
Suppose $\vec{w}_{\frak p}^{(1)} \subset \vec{w}_{\frak p}^{(2)}$.
Let
\begin{equation}\label{2149quote}
\pi : \underset{{\rm v}\in C^0(\mathcal G_{\frak p})}{{\bigodot}}\frak M^{(2)}_{(\frak x_{\frak p}\cup \vec w^{(2)}_{\frak p})_{\rm v}}
\to \prod_{{\rm v}\in C^0(\mathcal G_{\frak p})}\frak V^{(2)}((\frak x_{\frak p}\cup \vec w^{(2)}_{\frak p})_{\rm v})
\end{equation}
be the fiber bundle (\ref{2149}) that is a part of
the data included in ${\frak w}_{\frak p}^{(2)}$.
Here $\frak V^{(2)}((\frak x_{\frak p}\cup \vec w^{(2)}_{\frak p})_{\rm v})$ is an open neighborhood of
$(\frak x_{\frak p}\cup w^{(2)}_{\frak p})_{\rm v}$ in $\mathcal M_{k_{\rm v}+1,\ell_{\rm v} + \ell^{(2)}_{\rm v}}$
or in $\mathcal M^{\rm cl}_{\ell_{\rm v} + \ell^{(2)}_{\rm v}}$.
(They are contained in the top stratum of the Deligne-Mumford moduli spaces.)
\par
Forgetful map of the marked points in $\vec{w}_{\frak p}^{(2)} \setminus \vec{w}_{\frak p}^{(1)}$ induces a
map
$$
\frak{forget}_{\rm v} :
\mathcal M_{k_{\rm v}+1,\ell_{\rm v} + \ell^{(2)}_{\rm v}} \to \mathcal M_{k_{\rm v}+1,\ell_{\rm v} + \ell^{(1)}_{\rm v}}
$$
etc.
We put
$$
\frak{forget}_{\rm v}(\frak V^{(2)}((\frak x_{\frak p}\cup \vec w^{(2)}_{\frak p})_{\rm v}))
=
\frak V^{(1),+}((\frak x_{\frak p}\cup \vec w^{(1)}_{\frak p})_{\rm v}).
$$
We take $\frak V^{(1)}((\frak x_{\frak p}\cup \vec w^{(1)}_{\frak p})_{\rm v}) \subset \frak V^{(1),+}((\frak x_{\frak p}\cup \vec w^{(1)}_{\frak p})_{\rm v})$
that is a neighborhood of $(\frak x_{\frak p}\cup \vec w^{(1)}_{\frak p})_{\rm v}$
such that there exists a section
\begin{equation}\label{sect}
\frak{sect}_{\rm v} : \frak V^{(1)}((\frak x_{\frak p}\cup \vec w^{(1)}_{\frak p})_{\rm v})
\to \frak{forget}(\frak V^{(2)}((\frak x_{\frak p}\cup \vec w^{(2)}_{\frak p})_{\rm v})).
\end{equation}
Then we can pull back (\ref{2149quote}) by $\frak{sect} = (\frak{sect}_{\rm v})$ to obtain a fiber bundle
\begin{equation}\label{2149quote2}
\pi : \underset{{\rm v}\in C^0(\mathcal G_{\frak p})}{\bigodot}\frak M^{(1)}_{(\frak x_{\frak p}\cup \vec w^{(1)}_{\frak p})_{\rm v}}
\to \prod_{{\rm v}\in C^0(\mathcal G_{\frak p})}\frak V^{(1)}((\frak x_{\frak p}\cup \vec w^{(1)}_{\frak p})_{\rm v}).
\end{equation}
Moreover we can pull back a trivialization of the fiber bundle (\ref{2149quote}) to one of the fiber bundle (\ref{2149quote2}).
Thus we obtain a coordinate at infinity of $(\frak x_{\frak p}\cup w^{(1)}_{\frak p})_{\rm v}$.
\begin{defn}
We call the coordinate at infinity obtained as above the {\it coordinate at infinity induced from ${\frak w}_{\frak p}^{(2)}$}.
\par
We also take codimension $2$ submanifolds $\mathcal D^{(1)}_{\frak p,i}$ that are included as a part of
the stabilization data ${\frak w}_{\frak p}^{(1)}$,
so that $\mathcal D^{(1)}_{\frak p,i} = \mathcal D^{(2)}_{\frak p,i}$ for $i=1,\dots,\#\vec w^{(1)}_{\frak p}$.
We thus have obtained a stabilization date   ${\frak w}_{\frak p}^{(1)}$.
We call it  {\it the stabilization data induced from ${\frak w}_{\frak p}^{(2)}$}.
\end{defn}
We now construct a coordinate change of the Kuranishi structures in Case 3.
In Definition \ref{defnforget} we defined a forgetful map
$$
\aligned
\frak{forget}_{\frak B,\frak B^-;\vec w_{\frak p},\vec w^-_{\frak p}} &:
\mathcal M^{\frak w_{\frak p}}_{k+1,(\ell,\ell_{\frak p},(\ell_c))}(\beta;\frak p;\frak A;\frak B)_{\epsilon'_{0},\vec{\mathcal T}^{(2)}}\\
&\to
\mathcal M^{\frak w_{\frak p}^{-}}_{k+1,(\ell,\ell^-_{\frak p},(\ell_c))}(\beta;\frak p;\frak A;\frak B^-)_{\epsilon_{0},\vec{\mathcal T}^{(1)}}.
\endaligned
$$
Here we shrink the base space of (\ref{2149quote}) so that this map
is well-defined.
We need to extend the core of the domain and replace $\epsilon_0$ by $\epsilon'_0$ in the same way
as in Lemma \ref{2120lem}. We then obtain a stabilization data, which we denote by $\frak m^{(2)-}_{\frak p}$.
\par
Taking $\vec w_{\frak p} = \vec w_{\frak p}^{(2)-}$ and $\vec w^-_{\frak p} = \vec w^{(1)}_{\frak p}$
and $\frak B^- = \frak B$ we have
$$
\aligned
\frak{forget}_{\frak B,\frak B;\vec w_{\frak p}^{(2)-},\vec w_{\frak p}^{(1)}} &:
\mathcal M^{{\frak w}_{\frak p}^{(2)-}}_{k+1,(\ell,\ell_{\frak p}^{(2)},(\ell_c))}(\beta;\frak p;\frak A;\frak B)_{\epsilon'_{0},\vec{\mathcal T}^{(2)}}\\
&\to
\mathcal M^{{\frak w}_{\frak p}^{(1)}}_{k+1,(\ell,\ell^{(1)}_{\frak p},(\ell_c))}(\beta;\frak p;\frak A;\frak B)_{\epsilon_{0},\vec{\mathcal T}^{(1)}}.
\endaligned
$$
It induces a map
\begin{equation}\label{2318}
\aligned
\mathcal M^{{\frak w}_{\frak p}^{(2)}}_{k+1,(\ell,\ell_{\frak p}^{(2)},(\ell_c))}(\beta;\frak p;\frak A;\frak B)^{\vec w_{\frak p}^{(2)}\setminus \vec w_{\frak p}^{(1)}}_{\epsilon'_{0},\vec{\mathcal T}^{(2)}}\to
\mathcal M^{{\frak w}_{\frak p}^{(1)}}_{k+1,(\ell,\ell^{(1)}_{\frak p},(\ell_c))}(\beta;\frak p;\frak A;\frak B)_{\epsilon_{0},\vec{\mathcal T}^{(1)}}
\endaligned
\end{equation}
which is a strata-wise differentiable open embedding by Proposition \ref{forgetstillstable}.
We denote the map (\ref{2318}) by $\tilde\phi_{12}$.
\begin{lem}\label{tildephecm}
$\tilde\phi_{12}$ is of $C^m$-class in a neighborhood of $\frak p \cup \vec w_{\frak p}^{(2)}$.
\end{lem}
\begin{proof}
The proof is similar to the proof of Lemma \ref{2120lem}.
We use Lemma \ref{changeinfcoorproppara}
which is a parametrized version of Propositions \ref{changeinfcoorprop} and \ref{reparaexpest}.
Let
$$
\frak x_{\frak p}\cup \vec w_{\frak p}^{(2)} = \tilde{\frak x} \in
\prod_{{\rm v}\in C^0(\mathcal G_{\frak p})}\frak V^{(2)}((\frak x_{\frak p}\cup \vec w^{(2)}_{\frak p})_{\rm v})
$$
and $\frak{forget}(\tilde{\frak x}) = \frak x = \frak x_{\frak p}\cup \vec w_{\frak p}^{(1)}$.
Let $\frak V^{(1)-}((\frak x_{\frak p}\cup \vec w^{(1)}_{\frak p})_{\rm v})$
be a neighborhood of $\frak p$.
\par
Let $\frak{sect}_{(1),{\rm v}}$ be the section we chose in (\ref{sect}).
It gives a stabilization data $\frak w_{\frak p}^{(1)}$.
We take
$$
\frak{sect}_{(2),{\rm v}} :
Q_{\rm v} \times  \frak V^{(1)-}((\frak x_{\frak p}\cup \vec w^{(1)}_{\frak p})_{\rm v})
\to \frak V^{(2)}((\frak x_{\frak p}\cup \vec w^{(2)}_{\frak p})_{\rm v})
$$
such that the following condition is satisfied.
\begin{conds}
\begin{enumerate}
\item
$
\frak{forget}(\frak{sect}_{(2),{\rm v}}(\xi,\frak y_{\rm v})) = \frak y_{\rm v}.
$
\item
$\frak{sect}_{(2),{\rm v}}$ is a diffeomorphism onto an open neighborhood of $\tilde{\frak x}_{\rm v}$.
\end{enumerate}
\end{conds}
Pulling back  $\frak w^{(2)}_{\frak p}$ by $\frak{sect}_{(2)}$ we have a $Q = \prod Q_{\rm v}$-parametrized family
of stabilization data, which we call $\tilde{\frak w}^{(2)}_{\frak p}$.
We denote the image of $\frak{sect}_{(2),\rm v}$ by
$ \frak V^{(2)-}((\frak x_{\frak p}\cup \vec w^{(2)}_{\frak p})_{\rm v})$.
\par
We use $\frak w_{\frak p}^{(1)}$ in the same was as in the proof of
Lemma \ref{2120lem} to obtain
\begin{equation}\label{evaluandTfacsaa}
\aligned
&\frak F^{(1)} : \hat V^-(\frak p,\frak w^{(1)}_{\frak p};(\frak o^{(1)},\mathcal T^{(1)});\frak A)
\\
\to
&\prod_{{\rm v}\in C^0(\mathcal G_{\frak p})} C^{m}((K_{{\rm v},(1)}^{+\vec R_{(1)}},K_{{\rm v},(1)}^{+\vec R_{(1)}}
\cap \partial \Sigma_{{\rm v},(1)}),(X,L))\\
&\times \prod_{{\rm v}\in C^0(\mathcal G_{\frak p})}\frak V^{(1)-}((\frak x_{\frak p} \cup \vec w_{\frak p})_{\rm v})\\
&\times ((\vec{\mathcal T}^{(1)},\infty] \times ((\vec{\mathcal T}^{(1)},\infty] \times \vec S^1).
\endaligned
\end{equation}
(Here we put $-$ in  $ \hat V^-(\frak p,\frak w^{(1)}_{\frak p};(\frak o^{(1)},\mathcal T^{(1)});\frak A)$ to clarify that this space uses $\frak V^{(1)-}((\frak x_{\frak p} \cup \vec w_{\frak p})_{\rm v})$.)
\par
We use $\frak w_{\frak p}^{(2)}$ to obtain
\begin{equation}\label{evaluandTfac2paaa}
\aligned
&\frak F^{(2)} : \hat V^-(\frak p,\frak w^{(2)}_{\frak p};(\frak o^{(2)},\mathcal T^{(2)});\frak A)
\\
\to
&\prod_{{\rm v}\in C^0(\mathcal G_{\frak p})} C^{2m}((K_{{\rm v},(2)}^{+\vec R_{(2)}}),K_{{\rm v},(2)}^{+\vec R_{(2)}}\cap \partial \Sigma_{{\rm v},(2)}),(X,L))\\
&\times \prod_{{\rm v}\in C^0(\mathcal G_{\frak p})}\frak V^{(2)-}((\frak x_{\frak p} \cup \vec w_{\frak p})_{\rm v})\\
&\times ((\vec{\mathcal T}^{(2)},\infty] \times ((\vec{\mathcal T}^{(2)},\infty] \times \vec S^1).
\endaligned
\end{equation}
Let $\frak X(1,m)$, $\frak X(2,2m)$ be the spaces
in the right hand side of (\ref{evaluandTfacsaa}), (\ref{evaluandTfac2paaa})
respectively.
\par
We apply Lemma \ref{changeinfcoorproppara}
to the family of coordinates at infinity $\tilde{\frak w}^{(2)}_{\frak p}$
and the coordinate at infinity ${\frak w}^{(1)}_{\frak p}$.
It gives  estimates of the map $\overline{\Phi}_{12}$ defined in
(\ref{coodinatechange12})
and $\frak v_{(\xi,\rho,\vec T,\vec\theta)}$
as in (\ref{mapvbra}).
\par
We define $
\frak H_{12} : \frak X(2,2m) \to \frak X(1,m)
$ by
\begin{equation}\label{defnH12p}
\frak H_{12}(u,\frak{sect}_{(2)}(\xi,\rho),(\vec T,\vec \theta)) = (u\circ \frak v_{(\xi,\rho,\vec T,\vec \theta)},\overline{\Phi}_{12}(\xi,\rho,\vec T,\vec \theta)).
\end{equation}
By construction we have
\begin{equation}
\frak H_{12} \circ \frak F^{(2)} = \frak F^{(1)} \circ \tilde{\phi}_{12}.
\end{equation}
Lemma \ref{changeinfcoorproppara}  implies that
$\frak H_{12}$ is a $C^m$-map.
Moreover $\frak F^{(1)}$ and $\frak F^{(2)}$ are $C^m$-embeddings.
Therefore $\tilde\phi_{12}$ is a $C^m$-map on
$ \hat V^-(\frak p,\frak w^{(2)}_{\frak p};(\frak o^{(2)},\mathcal T^{(2)});\frak A)$.
The proof of Lemma \ref{tildephecm} is complete.
\end{proof}

We go back to the construction of coordinate change in Case 3.
By requiring the transversal constraint at all the marked points,
$\tilde\phi_{12}$ induces a required
coordinate change $\phi_{12}$.
Since $\frak A^{(1)} = \frak A^{(2)}$, it is easy to find the bundle map $\hat\phi_{12}$
that has the required properties.

\begin{rem}
Note that the map (\ref{2318}) and the coordinate change $\phi_{12}$ we obtain are independent of the choice of the section of (\ref{sect}).
But $\phi_{12}$ depends on the codimension 2 submanifolds we take, since the process to take $\rm trans$ depends on them.
We use the coordinate at infinity
(or the map $\frak{sect}_{\rm v}$ of (\ref{sect})) only to {\it prove} that $\phi_{12}$ is of $C^m$-class.
\end{rem}
Using the fact that the map (\ref{2318}) is a local diffeomorphism the construction of the
coordinate change in Case 4 is an inverse of one in Case 3.
\par\medskip
We have thus constructed the coordinate change in the 4 cases above.
The general case can be constructed by a composition of them.
\par
Let us be given $({\frak w}_{\frak p}^{(1)},\frak A^{(1)})$ and $({\frak w}_{\frak p}^{(2)},\frak A^{(2)})$.
We say that the pair $(({\frak w}_{\frak p}^{(1)},\frak A^{(1)}),({\frak w}_{\frak p}^{(2)},\frak A^{(2)}))$
is of Type 1,2,3,4, if we can apply Case 1,2,3,4, respectively.
We say the coordinate change obtained {\it the coordinate change of Type} 1,2,3,4, respectively.
\begin{lem}\label{lem2121}
For given $({\frak w}_{\frak p}^{(1)},\frak A^{(1)})$ and $({\frak w}_{\frak p}^{(6)},\frak A^{(6)})$
with $\vec w^{(1)}\cap \vec w^{(6)} = \emptyset$,
there exist $({\frak w}_{\frak p}^{(j)},\frak A^{(j)})$ for $j=2,\dots,5$ such that:
\par
The pair $(({\frak w}_{\frak p}^{(1)},\frak A^{(1)}),({\frak w}_{\frak p}^{(2)},\frak A^{(2)}))$ is of type $2$,
\par
The pair $(({\frak w}_{\frak p}^{(2)},\frak A^{(2)}),({\frak w}_{\frak p}^{(3)},\frak A^{(3)}))$ is of type $1$,
\par
The pair $(({\frak w}_{\frak p}^{(3)},\frak A^{(3)}),({\frak w}_{\frak p}^{(4)},\frak A^{(4)}))$ is of type $3$,
\par
The pair $(({\frak w}_{\frak p}^{(4)},\frak A^{(4)}),({\frak w}_{\frak p}^{(5)},\frak A^{(5)}))$ is of type $4$,
\par
The pair $(({\frak w}_{\frak p}^{(5)},\frak A^{(5)}),({\frak w}_{\frak p}^{(6)},\frak A^{(6)}))$ is of type $1$.
\end{lem}
\begin{proof}
We put $({\frak w}_{\frak p}^{(2)},\frak A^{(2)}) = ({\frak w}_{\frak p}^{(1)},\frak A^{(6)})$
and
$\frak A^{(j)} = \frak A^{(6)}$ for all $j=2,\dots,6$.
\par
Let $\vec w^{(4)}_{\frak p} = \vec w^{(1)}_{\frak p} \cup \vec w^{(6)}_{\frak p}$.
(Note this is a disjoint union by assumption.)
We take (any) coordinate at infinity for $\frak x_{\frak p} \cup \vec w^{(4)}
_{\frak p}$.
The codimension 2 submanifolds are determined from the data given in $\frak w^{(1)}_{\frak p}$ and  $\frak  w^{(6)}_{\frak p}$.
We thus defined $({\frak w}_{\frak p}^{(4)},\frak A^{(4)})$.
\par
We take the  coordinates at infinity that is induced from  ${\frak w}_{\frak p}^{(4)}$
so that the set of additional marked points are $\vec w^{(1)}_{\frak p}$ and $\vec w^{(6)}_{\frak p}$. We thus obtain $({\frak w}_{\frak p}^{(3)},\frak A^{(3)})$
and $({\frak w}_{\frak p}^{(5)},\frak A^{(5)})$, respectively.
It is easy to see that they have required properties.
\end{proof}
\begin{rem}
We need the hypothesis $\vec w^{(1)}_{\frak p}\cap \vec w^{(6)}_{\frak p} = \emptyset$ in Lemma \ref{lem2121}. Otherwise
it might happen that $w^{(1)}_{\frak p,i} = w^{(6)}_{\frak p,j}$ but
$\mathcal D^{(1)}_{\frak p,i} \ne \mathcal D^{(6)}_{\frak p,j}$.
\end{rem}
By Lemma \ref{lem2121} we can define a coordinate change for the pairs  $({\frak w}_{\frak p}^{(1)},\frak A^{(1)})$ and $({\frak w}_{\frak p}^{(2)},\frak A^{(2)})$ as the composition of 5 coordinate changes.
We have thus  constructed the required coordinate change
$$
\phi_{12} : V(\frak p,\frak w^{(2)}_{\frak p};(\frak o^{(2)},\mathcal T^{(2)});\frak A^{(2)})
\to V(\frak p,\frak w^{(1)}_{\frak p};(\frak o^{(1)},\mathcal T^{(1)});\frak A^{(1)})
$$
in case $\vec w^{(1)}_{\frak p} \cap \vec w^{(2)}_{\frak p} = \emptyset$.
\par
In  general cases we take ${\frak w}_{\frak p}^{(0)}$ such that
$\vec w^{(1)}_{\frak p} \cap \vec w^{(0)}_{\frak p} = \vec w^{(2)}_{\frak p} \cap \vec w^{(0)}_{\frak p} = \emptyset$
and put
$$
\phi_{12} = \phi_{10} \circ \phi_{02}.
$$
The proof of Proposition \ref{prop2117} is complete.
\end{proof}
We remark that in the proof of Lemma \ref{lem2121} we made a choice of coordinate at infinity
of  $\frak x_{\frak p} \cup \vec w^{(4)}_{\frak p}$.
We also take ${\frak w}^{(0)}_{\frak p}$ at the last step of the proof of Proposition \ref{prop2117}.
However the resulting coordinate change is
independent of these choices if we shrink the domain.
Namely we have the following Lemma \ref{lem123pre}.
We put
\begin{equation}\label{constructedUUU}
U(\frak p,\frak w_{\frak p};(\frak o,\mathcal T);\frak A)
=
V(\frak p,\frak w^{(1)}_{\frak p};(\frak o^{(1)},\mathcal T^{(1)});\frak A^{(1)})/\Gamma_{\frak p}.
\end{equation}
This is an orbifold.
\begin{rem}
We may replace the obstruction bundle by a bigger one and
may assume that the orbifold (\ref{constructedUUU}) is effective.
This is because the action of $\Gamma_{\frak p}$ on
$\Gamma(\Sigma_{\frak p};u_{\frak p}^*TX\otimes \Lambda^{01})$
is effective and it is still effective if we restrict to the space of
smooth sections
supported on a compact subset of a core.
\footnote{The same applies to the case of
higher genus. The only exception is the case $X$ is a point.}
In this article we always assume the effectivity of orbifolds.
\end{rem}
The map
$$
\phi_{12} : V(\frak p,\frak w^{(2)}_{\frak p};(\frak o^{(2)},\mathcal T^{(2)});\frak A^{(2)})
\to V(\frak p,\frak w^{(1)}_{\frak p};(\frak o^{(1)},\mathcal T^{(1)});\frak A^{(1)})
$$
induces
$$
\underline\phi_{12} : U(\frak p,\frak w^{(2)}_{\frak p};(\frak o^{(2)},\mathcal T^{(2)});\frak A^{(2)})
\to U(\frak p,\frak w^{(1)}_{\frak p};(\frak o^{(1)},\mathcal T^{(1)});\frak A^{(1)})
$$
which is an embedding of orbifold.
The embedding of vector bundle $\widehat{\phi}_{12}$ induces
$\underline{\widehat{\phi}}_{12}$ that is an embedding of orbibundles.
\begin{lem}\label{lem123pre}
We use the notation in Proposition \ref{prop2117}.
If two different choices of
$(\frak o^{(2),j}_0,\mathcal T^{(2),j}_0)$ $(j=1,2)$
and
$(\phi_{12}^j,\widehat{\phi}_{12}^j)$ $(j=1,2)$ are made,
then there exists
$(\frak o^{(3)},\mathcal T^{(3)})$ such that
$(\frak o^{(3)},\mathcal T^{(3)}) < (\frak o^{(2),j}_0,\mathcal T^{(2),j}_0)$ $(j=1,2)$
and
$$
(\underline\phi_{12}^1,\widehat{\underline\phi}_{12}^1) = (\underline\phi_{12}^2,\widehat{\underline\phi}_{12}^2)
$$
on
$
U(\frak p,\frak w^{(2)}_{\frak p};(\frak o^{(3)},\mathcal T^{(3)});\frak A^{(2)})
$.
\end{lem}
\begin{proof}
We first prove the next lemma.
\begin{lem}\label{lem123}
Let $\vec w_{\frak p}^{(1)} \subset \vec w_{\frak p}^{(2)}$.
Let  ${\frak w}_{\frak p}^{(i,j)}$ $i=1,2$, $j=1,2$ be the stabilization data at $\frak p$ such that
the additional marked points associated to ${\frak w}_{\frak p}^{(i,j)}$
is $\vec w_{\frak p}^{(j)}$.
\par
We assume that $(({\frak w}_{\frak p}^{(i,1)},\frak A),({\frak w}_{\frak p}^{(i,2)},\frak A))$
is type $3$.\footnote{Namely we assume that the coordinate at infinity of
${\frak w}_{\frak p}^{(i,1)}$ is induced by that of ${\frak w}_{\frak p}^{(i,2)}$ and
the submanifolds we assigned in
Definition \ref{stabdata} (2)  coincide each other when they are
assigned to the same marked points.}
\par
Let $\phi_{(i,j);(i',j')}$ be the coordinate change from the coordinate associated with
${\frak w}_{\frak p}^{(i',j')}$ to one associated with ${\frak w}_{\frak p}^{(i,j)}$.
Then we have
\begin{equation}\label{comb4sq}
\underline\phi_{(1,1);(1,2)}\circ \underline\phi_{(1,2);(2,2)} = \underline\phi_{(1,1);(2,1)}\circ \underline\phi_{(2,1);(2,2)}
\end{equation}
on a small neighborhood of $\frak p$ in the Kuranishi chart associated with ${\frak w}_{\frak p}^{(2,2)}$.
The same equality holds for $\widehat{\underline\phi}_{(i,j);(i',j')}$.
\par
The same conclusion holds when $\vec w_{\frak p}^{(2)} \subset \vec w_{\frak p}^{(1)}$ and replace `type $3$' by `type $4$'.
\end{lem}
\begin{rem}
We remark that the difference between $\vec w_{\frak p}^{(1)}$ and $\vec w_{\frak p}^{(2)}$ is
coordinate at infinity.
\end{rem}
\begin{proof}
This lemma as well as several other lemmas that appear later, is a consequence of the
following general observation.
\par
We consider an open subset $\mathcal U \subset \mathcal M_{k+1,\ell+\ell'}$ of the Deligne-Mumford moduli space.
Let
$$
\pi : \frak M(\mathcal U) \to \mathcal U
$$
be the restriction of
the universal family to $\mathcal U$.
Suppose we have a {\it topological space} $\Xi$ consisting of (an appropriate equivalence classes of)
pairs $(\frak x,u')$ where
$\frak x \in \mathcal U$ and $u' : \pi^{-1}(\frak x) \to X$ is a smooth map.
\footnote{Here equivalence relation is defined by an appropriate reparametrization of the source
by a biholomorphic map.}
Here we emphasis that we regard $\Xi$ as a topological space and do not need to use any
other structure such as smooth structure.
\par
Suppose $(V_i,\Gamma_i,E_i,\frak s_i,\psi_i)$ is a Kuranishi neighborhood at $\frak p$.
We assume that the coordinate change $\underline\phi_{ji}$ is defined as follows:
Suppose that there exists a homeomorphism $\Phi_i :  V_i/\Gamma_i \to \Xi$ onto an open neighborhood of $\frak x$ with
$\frak x = \Phi_{i}(\frak p)$ for all of $i$ and
$$
\underline\phi_{ji} = \Phi_{j}^{-1} \circ \Phi_{i}
$$
holds on a neighborhood of $\frak p$.
Then we have
$$
\underline\phi_{12}\circ \underline\phi_{23} = \underline\phi_{13}
$$
on a neighborhood of $\frak p$.
This observation is obvious.
\begin{rem}
Note it is important here that we only need to check set theoretical equality.
This is because our orbifolds are always effective orbifolds and we consider only
embeddings as maps between them.
Therefore we do not need to study orbifold structure or smooth structure to
prove compatibility of the coordinate changes etc.
\end{rem}
\begin{rem}\label{rem124}
Later we will use a slightly more general case.
Namely we consider the case when there are $V_{i,j}$ and $\Phi_{i,j} :  V_{i,j}/\Gamma_{i,j} \to \Xi$ for $(i,j) = (1,1),\dots,(1,m)$
and $(i,j) = (2,1),\dots,(2,n)$. We assume $V_{1,1} = V_{2,1}$ and $V_{1,m} = V_{2,n}$.
Suppose $\frak x = \Phi_{i,j}(\frak p)$ is independent of $i,j$ and $\Phi_{i,j}$
is a homeomorphism onto a neighborhood of $\frak x$.
We put:
$
\underline\phi_{(i,j)(i,j+1)} = \Phi_{i,j}^{-1} \circ \Phi_{i,j+1}.
$
Then we have
$$
\underline\phi_{(1,1)(1,2)}\circ\dots\circ \underline\phi_{(1,m-1)(1,m)}
=
\underline\phi_{(2,1)(2,2)}\circ\dots\circ \underline\phi_{(2,n-1)(2,n)}
$$
on a neighborhood of $\frak p$.
This is again obvious.
\end{rem}
Now we apply the observation above to the situation of Lemma \ref{lem123}.
The role of $\Xi$ is taken by
$$
\mathcal M^{{\frak w}_{\frak p}^{(2,2)}}_{k+1,(\ell,\ell_{\frak p}^{(2)},(\ell_c))}(\beta;\frak p;\frak A;\frak B)
^{\rm trans}_{\epsilon_{0},\vec{\mathcal T}^{(2)}}.
$$
We note that this set depends on the coordinate at infinity.
However Lemma \ref{lemma21118} implies that it is independent of
the coordinate at infinity on a neighborhood of $\frak p$.
We have thus proved (\ref{comb4sq}).
\par
Note the bundle maps  $\widehat{\underline\phi}_{(i,j);(i',j')}$ are
nothing but the identity maps on the fiber in our situation.
The proof of Lemma \ref{lem123} is complete.
\end{proof}
Lemma \ref{lem123pre} for the case  $\vec w^{(1)}_{\frak p}\cap \vec w^{(2)}_{\frak p} = \emptyset$ is
immediate from Lemma \ref{lem123}.
\par
Let us prove the general case. We need to prove the independence of the coordinate change of the
choice of ${\frak w}_{\frak p}^{(0)}$.
Let ${\frak w}_{\frak p}^{(0,1)}$, ${\frak w}_{\frak p}^{(0,2)}$ be two such choices.
Namely we assume
$\vec w^{(1)}_{\frak p}\cap \vec w^{(0,i)}_{\frak p} =\vec w^{(2)}_{\frak p}\cap \vec w^{(0,i)}_{\frak p} = \emptyset$
for $i=1,2$.
We first assume $\vec w^{(0,1)}_{\frak p} \cap \vec w^{(0,2)}_{\frak p} = \emptyset$ in addition.
We put $\vec w^{(0)}_{\frak p} = \vec w^{(0,1)}_{\frak p} \cup \vec w^{(0,2)}_{\frak p}$.
We take a stabilization data ${\frak w}_{\frak p}^{(0)}$ so that the codimension $2$ submanifolds are
induced by ${\frak w}_{\frak p}^{(0,i)}$.
Then,
$\underline\phi_{(0,i),0}$
are composition of coordinate change of type 3 and of type 1
and
$
\underline\phi_{0,(0,i)}
$
are composition of coordinate change of type 4 and of type 1.
Therefore from the first part of the proof we have
$$\aligned
\underline\phi_{1(0,1)}\circ \underline\phi_{(0,1)2}
&= \underline\phi_{1(0,1)}\circ \underline\phi_{(0,1)0}\circ \underline\phi_{0(0,1)}\circ\underline\phi_{(0,1)2}\\
&=\underline\phi_{10}\circ \underline\phi_{02}
=\underline\phi_{1(0,2)}\circ \underline\phi_{(0,2)2}
\endaligned$$
as required.\footnote{Here $\phi_{1(0,1)}$ is the coordinate change from the
Kuranishi chart associated with $\vec w^{(0,1)}_{\frak p}$ to the one
associated with $\vec w^{1}_{\frak p}$.
The notation of other coordinate changes are similar.}
\par
To remove the condition $\vec w^{(0,1)}_{\frak p} \cap \vec w^{(0,2)}_{\frak p} = \emptyset$
it suffices to remark that there exists ${\vec w}^{(0,3)}_{\frak p}$
such that
$\vec w^{(1)}_{\frak p}\cap \vec w^{(0,3)}_{\frak p} =\vec w^{(2)}_{\frak p}\cap \vec w^{(0,3)}_{\frak p} = \emptyset$
and
$\vec w^{(0,1)}_{\frak p} \cap \vec w^{(0,3)}_{\frak p} =\vec w^{(0,2)}_{\frak p} \cap \vec w^{(0,3)}_{\frak p} = \emptyset$.
The proof of Lemma \ref{lem123pre} is complete.
\end{proof}
Now we prove the compatibility of the coordinate transformations stated
in Proposition \ref{prop2117}.
\begin{lem}\label{lem125}
Let $({\frak w}_{\frak p}^{(j)},\frak A^{(j)})$ be a pair of stabilization data at $\frak p$ and $\frak A^{(j)} \subset \frak C(\frak p)$,
for $j=1,2,3$.
Suppose $\frak A^{(1)} \supseteq \frak A^{(2)} \supseteq \frak A^{(3)} \ne \emptyset$
and let $(\frak o^{(1)},\mathcal T^{(1)})$  be admissible for $({\frak w}_{\frak p}^{(1)},\frak A^{(1)})$.
\par
By Proposition \ref{prop2117} we have admissible $(\frak o^{(2)},\mathcal T^{(2)})$ and $(\frak o^{(3)},\mathcal T^{(3)})$
such that the coordinate change
$$(\phi_{1j},\hat{\phi}_{1j}) :
V(\frak p,\frak w^{(j)}_{\frak p};(\frak o,\mathcal T);\frak A^{(j)}) \to
V(\frak p,\frak w^{(1)}_{\frak p};(\frak o^{(1)},\mathcal T^{(1)});\frak A^{(2)})
$$
exists if $(\frak o^{(j)},\mathcal T^{(j)}) > (\frak o,\mathcal T)$. (Here $j=2,3$).
\par
By Proposition \ref{prop2117} there exists admissible $(\frak o^{(4)},\mathcal T^{(4)})$ such that
a coordinate change
$$(\phi_{23},\hat{\phi}_{23})
:
V(\frak p,\frak w^{(3)}_{\frak p};(\frak o,\mathcal T);\frak A^{(3)}) \to
V(\frak p,\frak w^{(2)}_{\frak p};(\frak o^{(2)},\mathcal T^{(2)});\frak A^{(2)})
$$
exists if $(\frak o^{(4)},\mathcal T^{(4)}) > (\frak o,\mathcal T)$.
\par
Now there exists $(\frak o^{(5)},\mathcal T^{(5)})$ with
$(\frak o^{(5)},\mathcal T^{(5)}) < (\frak o^{(j)},\mathcal T^{(j)})$
$(j=3,4)$ such that
we have
\begin{equation}
(\underline\phi_{13},\hat{\underline\phi}_{13}) = (\underline\phi_{12},\hat{\underline\phi}_{12})\circ(\underline\phi_{23},\hat{\underline\phi}_{23})
\end{equation}
on $U(\frak p,\frak w^{(3)}_{\frak p};(\frak o^{(5)},\mathcal T^{(5)});\frak A^{(3)})$.
\end{lem}
\begin{proof}
We first prove the case when $\vec w_{\frak p}^{(1)}$, $\vec w_{\frak p}^{(2)}$, $\vec w_{\frak p}^{(3)}$ are mutually disjoint.
\par
We note that we may assume $\frak A^{(1)} = \frak A^{(2)} = \frak A^{(3)} $.
In fact the coordinate change of type 2 (that is the coordinate change which replaces $\frak A$ by its subset $\frak A^-$),
is defined by inclusion of the domains so that $\frak A^-$ is obtained from $\frak A$ by the equation (\ref{2315+1}).
This process commutes with other types of coordinate changes.
So we assume $\frak A^{(1)} = \frak A^{(2)} = \frak A^{(3)} = \frak A$.
\par
We also note that the composition of two coordinate changes of type $j$ (for $j=1,\dots,4$) is again
a coordinate change of type $j$.
\par
Now using Lemma \ref{lem2121}, we can find
$\frak w_{\frak p}^{(i,j)}$ $i=1,2,3$, $j=2,\dots,6$ such that
$(\frak w_{\frak p}^{(i,j)},\frak w_{\frak p}^{(i,j+1)})$ is as in the conclusion of Lemma  \ref{lem2121}
and
$$
{\frak w}_{\frak p}^{(1,2)} = {\frak w}_{\frak p}^{(3,2)} = {\frak w}_{\frak p}^{(1)},
\quad
{\frak w}_{\frak p}^{(1,6)} = {\frak w}_{\frak p}^{(2,2)} = {\frak w}_{\frak p}^{(2)}, \quad
{\frak w}_{\frak p}^{(2,6)} = {\frak w}_{\frak p}^{(3,6)} = {\frak w}_{\frak p}^{(3)}.
$$
Then
$$
\aligned
\underline\phi_{12} &= \underline\phi_{(1,2)(1,3)}\circ\underline\phi_{(1,3)(1,4)}\circ\underline\phi_{(1,4)(1,5)}\circ\underline\phi_{(1,5)(1,6)}, \\
\underline\phi_{23} &= \underline\phi_{(2,2),(2,3)}\circ\underline\phi_{(2,3)(2,4)}\circ\underline\phi_{(2,4)(2,5)}\circ\underline\phi_{(2,5)(1,6)}, \\
\underline\phi_{13} &= \underline\phi_{(3,2)(3,3)}\circ\underline\phi_{(3,3)(1,4)}\circ\underline\phi_{(3,4)(3,5)}\circ\underline\phi_{(3,5)(3,6)}.
\endaligned
$$
Therefore we can apply the general observation mentioned
in the course of the proof of
Lemma \ref{lem123} in the form of Remark \ref{rem124} to prove Lemma \ref{lem125}
in our case.
\par
In fact we can take $\Xi$ as follows.
We consider $\vec w_{\frak p}^{(i,4)}$ for $i=1,2,3$ and put
$\vec w_{\frak p} = \vec w_{\frak p}^{(1,4)} \cup \vec w_{\frak p}^{(2,4)} \cup \vec w_{\frak p}^{(3,4)}$.
We take (any) coordinate at infinity of $\frak x_{\frak p} \cup \vec w_{\frak p}$.
We take the codimension 2 submanifolds $\mathcal D_{\frak p,i}$
(that is a part of the data  ${\frak w}_{\frak p}$) so that they coincide with those
taken for ${\frak w}^{(i)}_{\frak p}$, $i=1,2,3$. (Note we use the assumption that
$\vec w_{\frak p}^{(1)}$, $\vec w_{\frak p}^{(2)}$, $\vec w_{\frak p}^{(3)}$ are mutually disjoint here.)
We have thus defined the stabilization data $\frak w_{\frak p}$.
Then
$$
\Xi =
\mathcal M^{{\frak w}_{\frak p}}_{k+1,(\ell,\ell_{\frak p}^{(+)},(\ell_c))}(\beta;\frak p;\frak A;\frak B)
^{\rm trans}_{\epsilon_{0},\vec{\mathcal T}_1},
$$
where $\ell_{\frak p}^{(+)} = \#\vec w_{\frak p}$.
\par
We finally remove the condition that $\vec w_{\frak p}^{(1)}$, $\vec w_{\frak p}^{(2)}$, $\vec w_{\frak p}^{(3)}$ are mutually disjoint.
We take $\vec{w}_{\frak p}^{(4)}, \vec{w}_{\frak p}^{(5)}$ such that
$$
\vec w_{\frak p}^{(i)} \cap \vec{w}_{\frak p}^{(4)} = \emptyset = \vec w_{\frak p}^{(i)} \cap \vec{w}_{\frak p}^{(5)}
$$
for $i=1,2,3$ and $\vec w_{\frak p}^{(4)} \cap \vec{w}_{\frak p}^{(5)} = \emptyset$.
We also take codimension two transversal submanifolds $\mathcal D_i$
for each of those additional marked points.
We have thus obtained the stabilization data $\frak w_{\frak p}^{(4)}$,
$\frak w_{\frak p}^{(5)}$. Then we have
$$
\underline\phi_{12}\circ \underline\phi_{23} = \underline\phi_{15}\circ \underline\phi_{52}\circ\underline\phi_{24}\circ \underline\phi_{43}
=   \underline\phi_{15}\circ \underline\phi_{54}\circ \underline\phi_{43} = \underline\phi_{15}\circ \underline\phi_{53}
= \underline\phi_{13}.
$$
Here the first and the last equalities are the definitions. The second and the third equalities
follow from the case of Lemma \ref{lem125} which we already proved. The proof of
Lemma \ref{lem125} is complete.
\end{proof}

\section{Coordinate change - II: Coordinate change among different strata}
\label{differentstratum}

In this subsection we construct coordinate changes between the
Kuranishi charts we constructed in Proposition \ref{chartprop} for the general case.
Let $\frak p(1) \in \mathcal M_{k+1,\ell}(\beta)$.
We take a stabilization data $\frak w_{\frak p(1)}$ at $\frak p(1)$ and
$\frak A^{(1)} \subseteq \frak C(\frak p(1))$.
We use them to define Kuranishi neighborhood
$V(\frak p(1),\frak w_{\frak p(1)};(\frak o^{(1)},\mathcal T^{(1)});\frak A^{(1)})$
given in Definition \ref{defpsi}.
Let
\begin{equation}\label{322eq}
\psi_{\frak p(1),\frak w_{\frak p(1)};(\frak o^{(1)},\mathcal T^{(1)});\frak A^{(1)}} :
\frak s_{\frak p(1),\frak w_{\frak p(1)};(\frak o^{(1)},\mathcal T^{(1)});\frak A^{(1)}}^{-1}(0)
/\Gamma_{\frak p(1)} \to \mathcal M_{k+1,\ell}(\beta)
\end{equation}
be the map in Proposition \ref{chartprop}.
We assume that $\frak p(2)$ is contained  in its image.
\par
We will define the notion of {\it induced stabilization data} at $\frak p(2)$.
We recall that the stabilization data $\frak w_{\frak p(1)}$
includes the fiber bundle (\ref{2149})
\begin{equation}\label{21492}
\pi : \underset{{\rm v}\in C^0(\mathcal G_{\frak p(1)})}{\bigodot}\frak M^{(1)}((\frak x_{\frak p(1)}\cup \vec w_{\frak p(1)})_{\rm v})\to \prod_{{\rm v}\in C^0(\mathcal G_{\frak p(1)})}\frak V^{(1)}((\frak x_{\frak p(1)}\cup \vec w_{\frak p(1)})_{\rm v}).
\end{equation}
Here
$\frak V^{(1)}((\frak x_{\frak p(1)}\cup \vec w_{\frak p(1)})_{\rm v})$
is a neighborhood of
$(\frak x_{\frak p(1)}\cup \vec w_{\frak p(1)})_{\rm v}$
in the Deligne-Mumford moduli space
$\mathcal M_{k_{\rm v}+1,\ell_{\rm v}+\ell_{\frak p(1),\rm v}}$.
The product in the right hand side of (\ref{21492}) is identified with a
neighborhood of $\frak x_{\frak p(1)} \cup \vec w_{\frak p(1)}$
in the stratum $\mathcal M_{k+1,\ell+\ell_{\frak p(1)}}(\mathcal G_{\frak p(1)\cup \vec w_{\frak p(1)}})$ of the
Deligne-Mumford moduli space
$\mathcal M_{k+1,\ell+\ell_{\frak p(1)}}$.
We denote this neighborhood by
$\frak V(\frak x_{\frak p(1)}\cup \vec w_{\frak p(1)})$.
\begin{conds}
We consider a symmetric stabilization $\vec w_{\frak p(2)}$
on $\frak x_{\frak p(2)}$, an element $\sigma_0\in
\frak V(\frak x_{\frak p(1)}\cup \vec w_{\frak p(1)})$  and
$(\vec S^{\rm o}_{0},(\vec S^{\rm c}_{0},\vec \theta_0)))
\in (\vec{\mathcal T}^{(1)},\infty] \times ((\vec{\mathcal T}^{(1)},\infty] \times \vec S^1)$
that satisfy the following two conditions.
\begin{enumerate}\label{condindtransdata}
\item
$
\frak x_{\frak p(2)}\cup \vec w_{\frak p(2)}
= \overline{\Phi}(\sigma_0;\vec S^{\rm o}_{0},(\vec S^{\rm c}_{0},\vec \theta_0)).
$
\item
$\frak p(2) \cup \vec w_{\frak p(2)}$
satisfies the transversal constraint at all marked points.
Namely for each $i=1,\dots,\ell_{\frak p(1)}$ we have
$$
u_{\frak p(2)}(w_{\frak p(2),i}) \in \mathcal D_{\frak p(1),i}.
$$
Here $\mathcal D_{\frak p(1),i}$ is a codimension 2 submanifold
included in the stabilization data $\frak w_{\frak p(1)}$.
(We remark $\#\vec w_{\frak p(2)} = \# \vec w_{\frak p(1)} = \ell_{\frak p(1)}$.)
\end{enumerate}
\end{conds}
An element of $\Gamma_{\frak p(1)}$ is regarded as an
element of the permutation group $\frak S_{\ell_{\frak p(1)}}$.
So it transforms $\vec w_{\frak p(2)}$ by permutation.
The group $\Gamma_{\frak p(1)}$ acts also
on the set of pairs $(\sigma_0;\vec S^{\rm o}_{0},(\vec S^{\rm c}_{0},\vec \theta_0))$. We then have the following:
\begin{lem}\label{dependcompsst}
The set of triples $(\vec w_{\frak p(2)},\sigma_0;\vec S^{\rm o}_{0},(\vec S^{\rm c}_{0},\vec \theta_0))$
satisfying Condition \ref{condindtransdata}
consists of a single $\Gamma_{\frak p(1)}$-orbit.
\end{lem}
\begin{proof}
This is an immediate consequence of
Proposition \ref{charthomeo}.
\end{proof}
We continue the construction of the induced stabilization data
at $\frak p(2)$.
Let $\mathcal G_{\frak p(2)\cup \vec w_{\frak p(2)}}$
be the combinatorial type of $\frak p(2)\cup \vec w_{\frak p(2)}$.
In general it is different from the combinatorial
type  $\mathcal G_{\frak p(1)\cup \vec w_{\frak p(1)}}$
of $\frak p(1)\cup \vec w_{\frak p(1)}$.
In fact the graph $\mathcal G_{\frak p(2)\cup \vec w_{\frak p(2)}}$
is obtained from the graph $\mathcal G_{\frak p(1)\cup \vec w_{\frak p(1)}}$
by shrinking all the edges $\rm e$ such that $S_{0,{\rm e}} \ne \infty$.
We denote by $C^{1,{\rm fin}}(\mathcal G_{\frak p(1)\cup \vec w_{\frak p(1)}})$
the set of edges $\rm e$ with $S_{0,{\rm e}} \ne \infty$.
We have
\begin{equation}
C^{1}(\mathcal G_{\frak p(1)\cup \vec w_{\frak p(1)}})
=
C^{1,{\rm fin}}(\mathcal G_{\frak p(1)\cup \vec w_{\frak p(1)}})
\sqcup
C^{1}(\mathcal G_{\frak p(2)\cup \vec w_{\frak p(2)}}).
\end{equation}
Here the right hand side is the {\it disjoint} union.
Choose $\Delta S \in \R_{>0}$ that is sufficiently smaller than $S_{0,{\rm e}}$.
(We may take for example $\Delta S = 1$.)
\par
Let $\frak V^{(2)}(\frak x_{\frak p(2)}\cup \vec w_{\frak p(2)})$
be a neighborhood of
$\frak x_{\frak p(2)}\cup \vec w_{\frak p(2)}$
in the stratum
$\mathcal M_{k+1,\ell+\ell_{\frak p(1)}}(\mathcal G_{\frak p(2)\cup \vec w_{\frak p(2)}})$
of the Deligne-Mumford moduli space $\mathcal M_{k+1,\ell+\ell_{\frak p(1)}}$.
We can take them so that there exists an identification
\begin{equation}\label{2324}
\aligned
\frak V^{(2)}(\frak x_{\frak p(2)}\cup \vec w_{\frak p(2)})
=&\,\frak V^{(1)}(\frak x_{\frak p(1)}\cup \vec w_{\frak p(1)}) \\
&\times \prod_{{\rm e} \in C^{1,{\rm fin}}_{\rm o}(\mathcal G_{\frak p(1)\cup \vec w_{\frak p(1)}})} ((S_{0,\rm e}-\Delta S,S_{0,\rm e}+\Delta S)
\times [0,1])\\
&\times
 \prod_{{\rm e} \in C^{1,{\rm fin}}_{\rm c}(\mathcal G_{\frak p(1)\cup \vec w_{\frak p(1)}})} ((S_{0,\rm e}-\Delta S,S_{0,\rm e}+\Delta S)\times S^1).
\endaligned
\end{equation}
\par
Let $\overline{\rm v}$ be a vertex of $\mathcal G_{\frak p(2)\cup \vec w_{\frak p(2)}}$.
We take the subgraph
$\mathcal G_{\frak p(1)\cup \vec w_{\frak p(1)},\overline{\rm v}}$
of the graph
$\mathcal G_{\frak p(1)\cup \vec w_{\frak p(1)}}$
as follows.
There exists a map
$\mathcal G_{\frak p(1)\cup \vec w_{\frak p(1)}}
\to \mathcal G_{\frak p(2)\cup \vec w_{\frak p(2)}}$
that shrinks the edges $\rm e$ with $S_{0,\rm e} \ne \infty$.
An edge ${\rm e} \in C^1(\mathcal G_{\frak p(1)\cup \vec w_{\frak p(1)}})$
is an edge of $\mathcal G_{\frak p(1)\cup \vec w_{\frak p(1)},\overline{\rm v}}$
if it goes to the point $\overline{\rm v}$ by this map,
or it goes to the edge containing $\overline{\rm v}$ by this map.
Then we have
\begin{equation}\label{232400}
\aligned
&\frak V^{(2)}((\frak x_{\frak p(2)}\cup \vec w_{\frak p(2)})_{\overline{\rm v}}) \\
=&\prod_{{\rm v} \in C^0(\mathcal G_{\frak p(1)\cup \vec w_{\frak p(1)},\overline{\rm v}})}\frak V^{(1)}((\frak x_{\frak p(1)}\cup \vec w_{\frak p(1)})_{\rm v}) \\
&\times \prod_{{\rm e} \in C^{1,{\rm fin}}_{\rm o}(\mathcal G_{\frak p(1)\cup \vec w_{\frak p(1)},\overline{\rm v}})}
((S_{0,\rm e}-\Delta S,S_{0,\rm e}+\Delta S)\times [0,1])\\
&\times
 \prod_{{\rm e} \in C^{1,{\rm fin}}_{\rm c}(\mathcal G_{\frak p(1)\cup \vec w_{\frak p(1)},\overline{\rm v}})} ((S_{0,\rm e}-\Delta S,S_{0,\rm e}+\Delta S)\times S^1).
\endaligned
\end{equation}
The universal family over the Deligne-Mumford moduli space
restricts to a fiber bundle
\begin{equation}\label{2326}
\pi : \frak M^{(2)}((\frak x_{\frak p(2)}\cup \vec w_{\frak p(2)})_{\overline{\rm v}})\to
\frak V^{(2)}((\frak x_{\frak p(2)}\cup \vec w_{\frak p(2)})_{\overline{\rm v}}).
\end{equation}
The fiber at $(\sigma;\vec S^{\rm o},(\vec S^{\rm c},\vec \theta))$
of this bundle,
which we denote by $\Sigma_{(\sigma;\vec S^{\rm o},(\vec S^{\rm c},\vec \theta))}$, is the union of the following three types of
2 dimensional manifolds.
\begin{enumerate}
\item[(I)]
For each ${\rm v} \in C^0(\mathcal G_{\frak p(2)\cup \vec w_{\frak p(2)}})$
we consider the core $K_{\rm v}^{\sigma_{\rm v}}$ that is contained in
$\Sigma_{\sigma_{\rm v}}$.
(Here $\sigma_{\rm v} \in \frak V^{(1)}((\frak x_{\frak p(1)}\cup \vec w_{\frak p(1)})_{\rm v})$ is a component of $\sigma$ and $\Sigma_{\sigma_{\rm v}}$ is a
Riemann surface corresponding to this element ${\sigma_{\rm v}}$.)
\item[(II)]
If ${\rm e}\in C^1_{\mathrm o}(\mathcal G_{\frak p(2)\cup \vec w_{\frak p(2)}})$,
$S_{0,{\rm e}} = \infty$ and ${\rm e}$ goes to an outgoing edges of $\overline{\rm v}$,
we have  $[0,\infty) \times [0,1]$.
\par
If ${\rm e}\in C^1_{\mathrm o}(\mathcal G_{\frak p(2)\cup \vec w_{\frak p(2)}})$,
$S_{0,{\rm e}} = \infty$ and ${\rm e}$ goes to an incoming edge of $\overline{\rm v}$,
we have  $(-\infty,0] \times [0,1]$.
\par
If ${\rm e}\in C^1_{\mathrm c}(\mathcal G_{\frak p(2)\cup \vec w_{\frak p(2)}})$,
$S_{0,{\rm e}} = \infty$ and ${\rm e}$ goes to an outgoing edge of $\overline{\rm v}$,
we have  $[0,\infty) \times S^1$.
\par
If ${\rm e}\in C^1_{\mathrm c}(\mathcal G_{\frak p(2)\cup \vec w_{\frak p(2)}})$,
$S_{0,{\rm e}} = \infty$ and ${\rm e}$ goes to an incoming edge of $\overline{\rm v}$,
we have  $(-\infty,0] \times S^1$.
\item[(III)]
If ${\rm e}\in C^1_{\mathrm o}(\mathcal G_{\frak p(2)\cup \vec w_{\frak p(2)}})$,
$S_{0,{\rm e}} \ne \infty$,
we have $[-5S_{{\rm e}},5S_{{\rm e}}] \times [0,1]$.
If ${\rm e}\in C^1_{\mathrm c}(\mathcal G_{\frak p(2)\cup \vec w_{\frak p(2)}})$,
$S_{0,{\rm e}} \ne \infty$,
we have $[-5S_{{\rm e}},5S_{{\rm e}}] \times S^1$.
\end{enumerate}
\begin{defn}
The core $K_{\overline{\rm v}}$ of $\Sigma_{(\sigma;\vec S^{\rm o},(\vec S^{\rm c},\vec \theta))}$ is the union of the subsets of type I or type III.
\end{defn}
On the complement of the core, the fiber bundle
(\ref{2326}) has a
trivialization, that is given by the identification of the subsets of type II
with the standard set mentioned there. This trivialization preserves complex structures.
\par
This trivialization extends to the subsets of type I. In fact, such an extension is a part of the data
included in the coordinate at infinity of $\frak w_{\frak p(1)}$.
Note that this extension of trivialization does not respect the fiberwise complex structure.
\par
Note, however, that this trivialization does {\it not} extend to the trivialization
of the fiber bundle (\ref{2326}) if there exists an edge
${\rm e}\in C^1_{\mathrm c}(\mathcal G_{\frak p(2)\cup \vec w_{\frak p(2)}})$
with
$S_{0,{\rm e}} \ne \infty$.
In fact, there exists an $S^1$ factor
in (\ref{232400}) that corresponds to such an edge $\rm e$ and our
fiber bundle has nontrivial monodromy around it,
that is the Dehn twist at the domain $[-5S_{0,{\rm e}},5S_{0,{\rm e}}] \times S^1$.
\par
Therefore to find a coordinate at infinity that satisfies Definition
\ref{coordinatainfdef} (5) we need to restrict the domain.
We take a sufficiently small $\Delta\theta$
(for example $\Delta\theta = 1/10$) and put
\begin{equation}\label{232400aa}
\aligned
&\frak V((\frak x_{\frak p(2)}\cup \vec w_{\frak p(2)})_{\overline{\rm v}}) \\
=&\prod_{{\rm v} \in C^0(\mathcal G_{\frak p(1)\cup \vec w_{\frak p(1)},\overline{\rm v}})}\frak V((\frak x_{\frak p(1)}\cup \vec w_{\frak p(1)})_{\rm v}) \\
&\times \prod_{{\rm e} \in C^{1,{\rm fin}}_{\rm o}(\mathcal G_{\frak p(1)\cup \vec w_{\frak p(1)},\overline{\rm v}})}
((S_{0,\rm e}-\Delta S,S_{0,\rm e}+\Delta S)\times [0,1])\\
&\times
 \prod_{{\rm e} \in C^{1,{\rm fin}}_{\rm c}(\mathcal G_{\frak p(1)\cup \vec w_{\frak p(1)},\overline{\rm v}})} ((S_{0,\rm e}-\Delta S,S_{0,\rm e}+\Delta S)\times
(\theta_{0,{\rm e}}-\Delta{\theta},\theta_{0,{\rm e}}+\Delta{\theta})).
\endaligned
\end{equation}
(Note $\frak x_{\frak p(2)}\cup \vec w_{\frak p(2)}
= \overline{\Phi}(\sigma_0;\vec S^{\rm o}_{0},(\vec S^{\rm c}_{0},\vec \theta_0))$ and $\theta_{0,{\rm e}}$ is a component of
$\vec \theta_0)$.)
\par
We consider the fiber bundle
\begin{equation}\label{23262}
\pi : \frak M((\frak x_{\frak p(2)}\cup \vec w_{\frak p(2)})_{\overline{\rm v}})\to
\frak V((\frak x_{\frak p(2)}\cup \vec w_{\frak p(2)})_{\overline{\rm v}})
\end{equation}
in place of (\ref{2326}).
\par
Now we can extend the trivialization of the fiber bundle
defined in the complement of
the core, to the trivialization that is defined everywhere. (But it does not preserve the
complex structures.)
We  have thus defined a coordinate at infinity of $\frak p(2)$.
\par
We take the codimension 2 submanifolds $\mathcal D_{\frak p(1),i}$
that is a part of $\frak w_{\frak p(1)}$
and put
$$
\mathcal D_{\frak p(2),i}
=
\mathcal D_{\frak p(1),i}.
$$
\begin{defn}
The stabilization data at $\frak p(2)$ that is obtained as above is called the
{\it stabilization data induced by $\frak w_{\frak p(1)}$}.
\end{defn}
\begin{rem}
There are more than one
ways of extending the trivialization of the fiber bundle that is given
on the part of type I and type II to the whole space.
However the way to do so is
determined if we take the following two families of diffeomorphisms.
\begin{enumerate}
\item A family of diffeomorphisms from the rectangles
$[-5S_{{\rm e}},5S_{{\rm e}}] \times [0,1]$ to
$[0,1]\times [0,1]$ so that they are obvious isometries
in a neighborhood of $\partial[-5S_{{\rm e}},5S_{{\rm e}}] \times [0,1]$.
Here the parameter is $S_{{\rm e}} \in (S_{0,\rm e}-\Delta S,S_{0,\rm e}+\Delta S)$.
\item
A family of diffeomorphisms from the annuli
$[-5S_{{\rm e}},5S_{{\rm e}}] \times S^1$ to
$[0,1]\times S^1$ so that
they are obvious isometries
in a neighborhood of $\{-5S_{{\rm e}}\} \times S^1$
and is the rotation by $\theta_{\rm e}$ in a
neighborhood  $\{5S_{{\rm e}}\} \times S^1$.
Here the parameter is $S_{{\rm e}} \in (S_{0,\rm e}-\Delta S,S_{0,\rm e}+\Delta S)$
and $\theta_{\rm e} \in (\theta_{0,{\rm e}}-\Delta{\theta},\theta_{0,{\rm e}}+\Delta{\theta})$.
\end{enumerate}
Such families of diffeomorphisms obviously exist. We can take one
and use it whenever we define the induced coordinate at infinity.
In that sense the notion of induced coordinate at infinity
and of induced stabilization data is well-defined.
(Namely it can be taken independent of $\frak p(1)$ for example.)
\end{rem}
\par
In Section \ref{coordinateinf}, we discussed how the parametrization changes when we change the
coordinate at infinity. There we defined a map $\Psi_{12}$.
(See (\ref{2167}).) The following is obvious from definition. We use the notation in Propositions \ref{changeinfcoorprop} and \ref{reparaexpest}.
\begin{lem}\label{lempsi12}
If we take the induced core on $\frak Y_0$ then
$\overline{\Phi}_{12} = \Psi_{12}$.
Moreover $\frak v_{\frak y_2,\vec T_2,\vec\theta_2}$ is the identity map on the core $K_{\rm v}$.
\end{lem}
The first main result of this subsection is the following.
\begin{prop}\label{prop21333}
Let $\frak p(1) \in \mathcal M_{k+1,\ell}(\beta)$ and  take a stabilization data $\frak w_{\frak p(1)}$
at $\frak p(1)$ and admissible $(\frak o^{(1)}, \mathcal T^{(1)})$. Let $\frak p(2)$ be in the image of (\ref{322eq}).
We take the induced stabilization data $\frak w_{\frak p(2)}$.
Let $\frak A \subseteq \frak C(\frak p(2)) \subseteq \frak C(\frak p(1))$.
\par
Then there exists an admissible $(\frak o_0^{(2)}, \mathcal T_0^{(2)})$
such that if $(\frak o^{(2)}, \mathcal T^{(2)}) < (\frak o_0^{(2)}, \mathcal T_0^{(2)})$
there exists a coordinate change
$$
(\phi_{12},\hat\phi_{12}) :
V(\frak p(2),\frak w_{\frak p(2)};(\frak o^{(2)},\mathcal T^{(2)});\frak A)\to
V(\frak p(1),\frak w_{\frak p(1)};(\frak o^{(1)},\mathcal T^{(1)});\frak A).
$$
\end{prop}
\begin{proof}
We have maps
\begin{equation}\label{gluemapsss1a}
\aligned
\text{\rm Glu}^{(1)} :&B_{\frak o^{(1)}}^{\frak w_{\frak p(1)}}(\frak p;V_{k+1,(\ell,\ell_{\frak p},(\ell_c))}(\beta;\frak p(1);\frak A)) \times (\vec{\mathcal T}^{(1)},\infty] \times ((\vec{\mathcal T}^{(1)},\infty] \times \vec S^1)
\\
&\to
\mathcal M^{\frak w_{\frak p(1)}}_{k+1,(\ell,\ell_{\frak p(1)},(\ell_c))}(\beta;\frak p(1);\frak A)_{\epsilon_{0,1},\vec{\mathcal T}^{(1)}}
\endaligned
\end{equation}
and
\begin{equation}\label{gluemapsss1a2}
\aligned
\text{\rm Glu}^{(2)} :&B_{\frak o^{(2)}}^{\frak w_{\frak p(2)}}(\frak p;V_{k+1,(\ell,\ell_{\frak p},(\ell_c))}(\beta;\frak p(2);\frak A)) \times (\vec{\mathcal T}^{(2)},\infty] \times ((\vec{\mathcal T}^{(2)},\infty] \times \vec S^1)
\\
&\to
\mathcal M^{\frak w_{\frak p(2)}}_{k+1,(\ell,\ell_{\frak p(1)},(\ell_c))}(\beta;\frak p(2);\frak A)_{\epsilon_{0,2},\vec{\mathcal T}^{(2)}}
\endaligned
\end{equation}
by  the gluing constructions at $\frak p(1)$ and at $\frak p(2)$ respectively.
(More precisely for a given $\epsilon_{0,2}$, the map (\ref{gluemapsss1a2}) is defined
by choosing $\frak o^{(2)}$ small and $\mathcal T^{(2)}$ large.)
\par
By the assumption and
Proposition \ref{charthomeo}, there exists $\vec w_{\frak p(2),c}$ such that
$$
(\frak p(2)\cup  \vec w_{\frak p(2)},(\vec w_{\frak p(2),c}))
\in \mathcal M^{\frak w_{\frak p(1)}}_{k+1,(\ell,\ell_{\frak p(1)},(\ell_c))}(\beta;\frak p(1);\frak A)_{\epsilon_{0,1},\vec{\mathcal T}^{(1)}}^{\rm trans}.
$$
We observe
$$
(\frak p(2)\cup  \vec w_{\frak p(2)},(\vec w_{\frak p(2),c}))
\in \mathcal M^{\frak w_{\frak p(2)}}_{k+1,(\ell,\ell_{\frak p(2)},(\ell_c))}(\beta;\frak p(1);\frak A)_{\epsilon_{0,2},\vec{\mathcal T}^{(2)}}
$$
and the image of (\ref{gluemapsss1a2}) defines a neighborhood basis when
we move  $\epsilon_{0,2}$.
Therefore by taking $\epsilon_{0,2}$ small and $\mathcal T^{(2)}$ large, we may assume that
\begin{equation}\label{thikenincl232}
\aligned
&\mathcal M^{\frak w_{\frak p(2)}}_{k+1,(\ell,\ell_{\frak p(1)},(\ell_c))}(\beta;\frak p(2);\frak A)_{\epsilon_{0,2},\vec {\mathcal T^{(2)}}}\\
&\subset
\mathcal M^{\frak w_{\frak p(1)}}_{k+1,(\ell,\ell_{\frak p(1)},(\ell_c))}(\beta;\frak p(1);\frak A)_{\epsilon_{0,1},\vec{\mathcal T}^{(1)}}
\endaligned\end{equation}
and this is an open embedding.
By construction, the element of the thickened moduli space
$\mathcal M^{\frak w_{\frak p(2)}}_{k+1,(\ell,\ell_{\frak p(2)},(\ell_c))}(\beta;\frak p(2);\frak A)_{\epsilon_{0,2},\vec{\mathcal T}^{(2)}}$
satisfies the transversal constraint at all additional marked points with respect to $\frak w_{\frak p(1)}$
if and only if  the transversal constraint at all additional marked points with respect to $\frak w_{\frak p(2)}$ is satisfied.
Therefore
\begin{equation}\label{2332}
\aligned
&\mathcal M^{\frak w_{\frak p(2)}}_{k+1,(\ell,\ell_{\frak p(1)},(\ell_c))}(\beta;\frak p(2);\frak A)_{\epsilon_{0,2},\vec{\mathcal T}^{(2)}}^{\rm trans}\\
&\subset
\mathcal M^{\frak w_{\frak p(1)}}_{k+1,(\ell,\ell_{\frak p(1)},(\ell_c))}(\beta;\frak p(1);\frak A)_{\epsilon_{0,1},\vec{\mathcal T}^{(1)}}^{\rm trans}
\endaligned\end{equation}
and this is an open embedding.
We thus can define a continuous strata-wise $C^m$-map $\phi_{12}$ as the inclusion map.
It is an open embedding of $C^m$-class strata-wise.
\begin{lem}\label{coodinatedifferentp}
$\phi_{12}$ is of $C^m$-class.
\end{lem}
\begin{proof}
The proof is similar to that of Lemma \ref{2120lem}.
We repeat the detail for completeness.
Let
$\hat V(\frak p(j),\frak w_{\frak p(j)};(\frak o^{(j)},\mathcal T^{(j)});\frak A)$
be the inverse image of
$V(\frak p(j),\frak w_{\frak p(j)};(\frak o^{(j)},\mathcal T^{(j)});\frak A)$
by $\text{\rm Glu}^{(j)}$. (Here $j=1,2$). It suffices to show that
$$
\aligned
\tilde\phi_{12} = (\text{\rm Glu}^{(1)})^{-1}\circ \text{\rm Glu}^{(2)} :
&\hat V(\frak p(2),\frak w_{\frak p(2)};(\frak o^{(2)},\mathcal T^{(2)});\frak A) \\
&\to
\hat V(\frak p(1),\frak w_{\frak p(1)};(\frak o^{(1)},\mathcal T^{(1)});\frak A)
\endaligned
$$
is of $C^m$-class.
We obtain maps
\begin{equation}\label{evaluandTfacjj}
\aligned
\frak F^{(j)} :
\hat V&(\frak p(j),\frak w_{\frak p(j)};(\frak o^{(j)},\mathcal T^{(j)});\frak A)\\
\to
&\prod_{{\rm v}\in C^0(\mathcal G_{\frak p(1)})} C^m((K_{\rm v}^{+\vec R},K_{\rm v}^{+\vec R} \cap \partial\Sigma_{{\rm v},(1)}),(X,L))\\
&\times \prod_{{\rm v}\in C^0(\mathcal G_{\frak p(1)})}\frak V((\frak x_{\frak p(1)} \cup \vec w_{\frak p(1)})_{\rm v})\times (\vec{\mathcal T}^{(j)},\infty] \times (\vec{\mathcal T}^{(j)},\infty] \times \vec S^1)
\endaligned
\end{equation}
in the same way as (\ref{evaluandTfac}) for $j=1,2$.
We remark here that we take the graph $\mathcal G_{\frak p(1)}$ for the case $j=2$ also.
By applying Theorem \ref{exdecayT33} we find that (\ref{evaluandTfacjj}) is an
$C^m$-embedding for $j=1$.
\par
We will prove that (\ref{evaluandTfacjj}) is a $C^m$-embedding for $j=2$ also.
It follows from Theorem \ref{exdecayT33} applied to the gluing at $\frak p(2)$ that
$\frak F^{(2)}$ is of $C^m$-class.
We put $\frak F^{(2)} = (\frak F^{(2)}_1,\frak F^{(2)}_2)$. Here $\frak F^{(2)}_1$ (resp.$\frak F^{(2)}_2$) is a map to the factor in the
second line (resp. third line).
It suffices to show that $\frak F^{(2)}_1$  is a $C^m$-embedding on each of the fiber of $\frak F^{(2)}_2$.
Note that the factors of the third line parametrize the complex structure of the source.
The fact that $\frak F^{(2)}_1$ is an embedding on the fiber of $T_{\rm e} = \infty$ follows from
Theorem \ref{exdecayT33} applied to the gluing at $\frak p(1)$.
Then we apply Theorem \ref{exdecayT33} to the gluing at $\frak p(2)$ to show that  $\frak F^{(2)}_1$ is an embedding on the fiber
of $\frak F^{(2)}_2$ if $\mathcal T^{(2)}$ is sufficiently large.
\par
Now using the obvious fact that $\frak F^{(1)}\circ \tilde\phi_{12} = \frak F^{(2)}$, we conclude that
$\tilde\phi_{12}$ is a $C^m$-embedding.
\end{proof}
\begin{rem}
Contrary to the case of the proof of Lemma \ref{2120lem},
we do have $\frak F^{(1)}\circ \tilde\phi_{12} = \frak F^{(2)}$.
This is because we are using the coordinate at infinity $\frak w_{\frak p(2)}$ that is
induced from $\frak w_{\frak p(1)}$ and so the parametrization of the core is the same.
\end{rem}
We thus have defined $\phi_{12}$. We define $\hat\phi_{12} = \phi_{12} \times {\rm identity}$.
It is easy to see that  $\phi_{12}$ is $\Gamma_{\frak p(2)}$-equivariant.
Other properties are also easy to prove. The proof of Proposition \ref{prop21333} is now complete.
\end{proof}
\begin{rem}
In Lemma \ref{dependcompsst} we proved that the two choices of $\vec w_{(2)}$ are transformed each other
under the $\Gamma_{\frak p(1)}$ action.
More precisely we have the following.
The action of $\Gamma_{\frak p(1)}$ is given by the permutation of the marked points $\vec w_{(2)}$.
If $\gamma \in \Gamma_{\frak p(2)}$ the permutation of $\vec w_{(2)}$ gives an equivalent element.
Namely there exists a biholomorphic map $\frak x_{\frak p(2)} \cup \vec w_{(2)} \to \frak x_{\frak p(2)} \cup \gamma\vec w_{(2)}$.
\par
In case $\gamma \notin \Gamma_{\frak p(2)}$,
$\frak x_{\frak p(2)} \cup \vec w_{(2)}$ is not biholomorphic to $\frak x_{\frak p(2)} \cup \gamma\vec w_{(2)}$.
Each of the choice $\vec w_{(2)}$ and $\gamma\vec w_{(2)}$
induces a stabilization data at $\frak p(2)$, which we write $\frak w_{(2)}$ and $\gamma\frak w_{(2)}$
respectively.
They define the coordinate changes.
We remark that there is a canonical diffeomorphism
$$
\mathcal M^{\frak w_{\frak p(2)}}_{k+1,(\ell,\ell_{\frak p(1)},(\ell_c))}(\beta;\frak p(2);\frak A)^{\rm trans}_{\epsilon_{0,2},\vec{\mathcal T}^{(2)}}
\cong
\mathcal M^{\gamma\frak w_{\frak p(2)}}_{k+1,(\ell,\ell_{\frak p(1)},(\ell_c))}(\beta;\frak p(2);\frak A)^{\rm trans}_{\epsilon_{0,2},\vec{\mathcal T}^{(2)}}
$$
by permutation of the marked points.
Namely we have
$$
\gamma :
V(\frak p(2),\frak w_{\frak p(2)};(\frak o^{(2)},\mathcal T^{(2)});\frak A)
\to
V(\frak p(2),\gamma\frak w_{\frak p(2)};(\frak o^{(2)},\mathcal T^{(2)});\frak A).
$$
On the other hand $\gamma \in \Gamma_{\frak p(1)}$ acts on
$V(\frak p(1),\frak w_{\frak p(1)};(\frak o^{(1)},\mathcal T^{(1)});\frak A)$.
Since our construction is $\Gamma_{\frak p(1)}$ equivariant we have
$$
\gamma \circ \phi_{12} = \phi_{12} \circ \gamma.
$$
Here $\phi_{12}$ in the left hand side uses $\frak w_{\frak p(2)}$ and
$\phi_{12}$ in the right hand side uses $\gamma\frak w_{\frak p(2)}$.
This is the same  as the case of coordinate change of the charts of orbifolds.
\par
We remark that this phenomenon finally causes an appearance of $\gamma^{\alpha}_{pqr}$
in Definition \ref{Definition A1.5}.
\end{rem}
Combined with the result of the last subsection, Proposition \ref{prop21333} implies the following.
\begin{cor}\label{cchanges}
Let $\frak p(1) \in \mathcal M_{k+1,\ell}(\beta)$. We take a stabilization data $\frak w_{\frak p(1)}$
at $\frak p(1)$ and admissible $(\frak o^{(1)}, \mathcal T^{(1)})$.
\par
Let $\frak p(2)$ be in the image of (\ref{322eq}).
We take a stabilization data $\frak w_{\frak p(2)}$ at $\frak p(2)$.
\par
Then there exists an admissible $(\frak o_0^{(2)}, \mathcal T_0^{(2)})$
such that the following holds for $\frak A^{(j)} \subseteq \frak C(\frak p(j))$ with
$\frak A^{(2)} \subseteq \frak A^{(1)}$.
\par
For any $(\frak o^{(2)}, \mathcal T^{(2)}) < (\frak o_0^{(2)}, \mathcal T_0^{(2)})$ then
there exists a coordinate change
$(\phi_{12},\hat\phi_{12})$ from the Kuranishi chart
$V(\frak p(2),\frak w_{\frak p(2)};(\frak o^{(2)},\mathcal T^{(2)});\frak A^{(2)})$
to the Kuranishi chart
$V(\frak p(1),\frak w_{\frak p(1)};(\frak o^{(1)},\mathcal T^{(1)});\frak A^{(1)})$.
\end{cor}
\begin{proof}
Let $\frak w_{\frak p(2)}'$ be the stabilization data at $\frak p(2)$ induced by
$\frak w_{\frak p(1)}$.
Then the required  coordinate change
is obtained by composing the three coordinate changes associated to the pairs,
$((\frak w_{\frak p(1)},\frak A^{(1)}),(\frak w_{\frak p(1)},\frak A^{(2)}))$,
$((\frak w_{\frak p(1)},\frak A^{(2)}),(\frak w_{\frak p(2)}',\frak A^{(2)}))$,
$((\frak w_{\frak p(2)}',\frak A^{(2)}),(\frak w_{\frak p(2)},\frak A^{(2)}))$.
They are obtained by Proposition \ref{prop2117},
Proposition \ref{prop21333}, Proposition \ref{prop2117}, respectively.
\end{proof}
\begin{rem}
By construction the coordinate change given in Corollary \ref{cchanges} is
independent of the choices involved in
the definition, in a neighborhood of $\frak p(2)$.
\end{rem}
We next prove the compatibility of the coordinate changes in
Corollary \ref{cchanges}.
\begin{prop}\label{compaticoochamain}
Let $\frak p(1) \in \mathcal M_{k+1,\ell}(\beta)$. We take a stabilization data $\frak w_{\frak p(1)}$
at $\frak p(1)$ and admissible  $(\frak o^{(1)}, \mathcal T^{(1)})$.
\par
Let $\frak p(2)$ be in the image of (\ref{322eq}).
We take a stabilization data $\frak w_{\frak p(2)}$ at $\frak p(2)$.
Let $(\frak o_0^{(2)}, \mathcal T_0^{(2)})$ be as in Corollary \ref{cchanges}.
\par
Then there exists $
\epsilon_{7} = \epsilon_7(\frak p(1),\frak w_{\frak p(1)},\frak p(2),\frak w_{\frak p(2)})$
with the following properties for each $(\frak o^{(2)}, \mathcal T^{(2)}) < (\frak o_0^{(2)}, \mathcal T_0^{(2)})$.
\par
Let $\frak p(3) \in \mathcal M_{k+1,\ell}(\beta)$.
We assume $d(\frak p(2),\frak p(3)) < \epsilon_{7}$.\footnote{$d$ here is any metric on $\mathcal M_{k+1,\ell}(\beta)$.}
Then for any stabilization data $\frak m_{\frak p(3)}$ at $\frak p(3)$,
there exists admissible $(\frak o_0^{(3)}, \mathcal T_0^{(3)})$ such that
if $(\frak o^{(3)}, \mathcal T^{(3)}) < (\frak o_0^{(3)}, \mathcal T_0^{(3)})$ and
$\frak A^{(j)} \subseteq \frak C(\frak p(j))$ ($j=1,2,3$) with
$\frak A^{(1)} \supseteq \frak A^{(2)}\supseteq \frak A^{(3)}$, then
we have the following.
\begin{enumerate}
\item There exists a coordinate change
$$(\phi_{23},\hat\phi_{23})
:
V(\frak p(3),\frak w_{\frak p(3)};(\frak o^{(3)},\mathcal T^{(3)});\frak A^{(3)})\to
V(\frak p(2),\frak w_{\frak p(2)};(\frak o^{(2)},\mathcal T^{(2)});\frak A^{(2)})
$$
as in Corollary \ref{cchanges}.
\item
There exists a coordinate change
$$
(\phi_{13},\hat\phi_{13})
:
V(\frak p(3),\frak w_{\frak p(3)};(\frak o^{(3)},\mathcal T^{(3)});\frak A^{(3)})
\to
V(\frak p(1),\frak w_{\frak p(1)};(\frak o^{(1)},\mathcal T^{(1)});\frak A^{(1)})
$$
as in Corollary \ref{cchanges}.
\item
We have
$$
(\underline\phi_{13},\hat{\underline\phi}_{13})  = (\underline\phi_{12},\hat{\underline\phi}_{12}) \circ (\underline\phi_{23},\hat{\underline\phi}_{23}).
$$
Here
$$
(\phi_{12},\hat\phi_{12})
:
V(\frak p(2),\frak w_{\frak p(2)};(\frak o^{(2)},\mathcal T^{(2)});\frak A^{(2)})
\to
V(\frak p(1),\frak w_{\frak p(1)};(\frak o^{(1)},\mathcal T^{(1)});\frak A^{(1)})
$$
is the coordinate change in Corollary \ref{cchanges} and
$$
(\underline\phi_{12},\hat{\underline\phi}_{12})
:
U(\frak p(2),\frak w_{\frak p(2)};(\frak o^{(2)},\mathcal T^{(2)});\frak A^{(2)})
\to
U(\frak p(1),\frak w_{\frak p(1)};(\frak o^{(1)},\mathcal T^{(1)});\frak A^{(1)})
$$
is induced by it.
\end{enumerate}
\end{prop}
\begin{proof}
By the same reason as in the case of Proposition \ref{prop21333}
we may assume $\frak A^{(1)} = \frak A^{(2)} = \frak A^{(3)} = \frak A$.
So we will assume it throughout the proof.
\par
We first prove the following.
\begin{lem}\label{lem2144}
Let $\frak w_{\frak p(2)}^{(1)}$ be the stabilization data at $\frak p(2)$ induced by $\frak w_{\frak p(1)}$ and
$\frak w_{\frak p(3)}^{(1)}$ the stabilization data at $\frak p(3)$ induced by $\frak w_{\frak p(2)}^{(1)}$.
Then $\frak w_{\frak p(3)}^{(1)}$ is the stabilization data induced by $\frak w_{\frak p(1)}$.
\end{lem}
The proof is obvious.
\begin{lem}\label{lem21409}
Let $\frak w_{\frak p(1)}^{(1)} = \frak w_{\frak p(1)}$ and $\frak w_{\frak p(2)}^{(1)}$, $\frak w_{\frak p(3)}^{(1)}$ be as in
Lemma \ref{lem2144}.
We denote by $(\phi_{ij},\hat\phi_{ij})$ (for $1\le i < j\le 3$) the coordinate changes induced by the pair
$(\frak w_{\frak p(i)}^{(1)},\frak w_{\frak p(j)}^{(1)})$.
Then we have
\begin{equation}\label{2140formula}
(\phi_{12},\hat\phi_{12})\circ (\phi_{23},\hat\phi_{23}) = (\phi_{13},\hat\phi_{13})
\end{equation}
in a neighborhood of $\frak p(3)$.
\end{lem}
\begin{proof}
We can choose $\epsilon_{0,j}$ $(j=1,2,3)$ such that
$$
\aligned
&\mathcal M^{\frak w^{(1)}_{\frak p(3)}}_{k+1,(\ell,\ell_{\frak p(1)},(\ell_c))}(\beta;\frak p(3);\frak A)_{\epsilon_{0,3},\vec{\mathcal T}^{(3)}}^{\rm trans}\\
&\subset
\mathcal M^{\frak w^{(1)}_{\frak p(2)}}_{k+1,(\ell,\ell_{\frak p(1)},(\ell_c))}(\beta;\frak p(2);\frak A)_{\epsilon_{0,2},\vec {\mathcal T}^{(2)}}^{\rm trans}\\
&\subset
\mathcal M^{\frak w^{(1)}_{\frak p(1)}}_{k+1,(\ell,\ell_{\frak p(1)},(\ell_c))}(\beta;\frak p(1);\frak A)_{\epsilon_{0,1},\vec{\mathcal T}^{(1)}}^{\rm trans}.
\endaligned
$$
The maps (\ref{2140formula}) are all induced by this inclusion in a neighborhood of $\frak p(3)$.
Hence the lemma.
\end{proof}
The proof of the next lemma is the main part of the proof of Proposition \ref{compaticoochamain}.
\begin{lem}\label{2141}
Let $\frak p(2) \in \mathcal M_{k+1,\ell}(\beta)$ and let
$\frak w_{\frak p(2)}^{(1)}$, $\frak w_{\frak p(2)}^{(2)}$ be two stabilization data at $\frak p(2)$.
Suppose a $\frak w_{\frak p(2)}^{(1)}$-admissible $(\frak o^{(21)}, \mathcal T^{(21)})$ is given.
Take an $\frak w_{\frak p(2)}^{(2)}$-admissible $(\frak o_0^{(22)}, \mathcal T_0^{(22)})$ such that if  $(\frak o^{(22)}, \mathcal T^{(22)}) < (\frak o_0^{(22)}, \mathcal T_0^{(22)})$
then there exists a coordinate change
$$
(\phi_{(21)(22)},\hat\phi_{(21)(22)}) :
V(\frak p(2),\frak w^{(2)}_{\frak p(2)};(\frak o^{(22)},\mathcal T^{(22)});\frak A)\to
V(\frak p(2),\frak w^{(1)}_{\frak p(2)};(\frak o^{(21)},\mathcal T^{(21)});\frak A)
$$
as in Proposition \ref{prop2117}.
\par
Then there exists
$\epsilon_{8} = \epsilon_8(\frak p(2),\frak w_{\frak p(2)}^{(1)},\frak w_{\frak p(2)}^{(2)},(\frak o^{(21)},
\mathcal T^{(21)}),(\frak o^{(22)}, \mathcal T^{(22)}))$
such that if $\frak p(3)  \in \mathcal M_{k+1,\ell}(\beta)$, $d(\frak p(2),\frak p(3)) < \epsilon_{8}$ the following holds.
\begin{enumerate}
\item
There exists a stabilization data $\frak m_{\frak p(3)}^{(1)}$ at $\frak p(3)$ induced from $\frak w_{\frak p(2)}^{(1)}$
and a stabilization data  $\frak m_{\frak p(3)}^{(2)}$at $\frak p(3)$ induced from $\frak w_{\frak p(2)}^{(2)}$.
\item
There exists a $\frak w_{\frak p(3)}^{(1)}$-admissible $(\frak o_0^{(31)}, \mathcal T_0^{(31)})$ such that if
$(\frak o^{(31)}, \mathcal T^{(31)}) < (\frak o_0^{(31)}, \mathcal T_0^{(31)})$ then the coordinate change
$$
(\phi_{(21)(31)},\hat\phi_{(21)(31)}) :
V(\frak p(3),\frak w^{(1)}_{\frak p(3)};(\frak o^{(31)},\mathcal T^{(31)});\frak A)
\to
V(\frak p(2),\frak w^{(1)}_{\frak p(2)};(\frak o^{(21)},\mathcal T^{(21)});\frak A)
$$
as in Proposition \ref{prop21333} exists.
\item
There exists a $\frak w_{\frak p(3)}^{(2)}$-admissible $(\frak o_0^{(32)}, \mathcal T_0^{(32)})$ such that if
$(\frak o^{(32)}, \mathcal T^{(32)}) < (\frak o_0^{(32)}, \mathcal T_0^{(32)})$ then the
coordinate change
$$
(\phi_{(22)(32)},\hat\phi_{(22)(32)}) :
V(\frak p(3),\frak w^{(2)}_{\frak p(3)};(\frak o^{(32)},\mathcal T^{(32)});\frak A)
\to
V(\frak p(2),\frak w^{(2)}_{\frak p(2)};(\frak o^{(22)},\mathcal T^{(22)});\frak A)
$$
as in Proposition \ref{prop21333} exists.
\item
There exists a $\frak w_{\frak p(3)}^{(2)}$-admissible $(\frak o_0^{(32)\prime}, \mathcal T_0^{(32)\prime})$ such that if
$(\frak o^{(32)\prime}, \mathcal T^{(32)\prime}) < (\frak o_0^{(32)\prime}, \mathcal T_0^{(32)\prime})$ then the
coordinate change
$$
(\phi_{(31)(32)},\hat\phi_{(31)(32)}) :
V(\frak p(3),\frak w^{(2)}_{\frak p(3)};(\frak o^{(32)\prime},\mathcal T^{(32)\prime});\frak A)
\to
V(\frak p(3),\frak w^{(1)}_{\frak p(3)};(\frak o^{(31)},\mathcal T^{(31)});\frak A)
$$
as in Proposition \ref{prop2117} exists.
\item
Suppose $(\frak o^{(32)\prime\prime}, \mathcal T^{(32)\prime\prime}) <  (\frak o_0^{(32)\prime}, \mathcal T_0^{(32)\prime})$
and
$(\frak o^{(32)\prime\prime}, \mathcal T^{(32)\prime\prime}) <  (\frak o_0^{(32)}, \mathcal T_0^{(32)})$.
Then we have
\begin{equation}
\aligned
&(\underline\phi_{(21)(22)},\hat{\underline\phi}_{(21)(22)})\circ ({\underline\phi}_{(22)(32)},\hat{\underline\phi}_{(22)(32)})\\
&=
({\underline\phi}_{(21)(31)},\hat{\underline\phi}_{(21)(31)})\circ ({\underline\phi}_{(31)(32)},\hat{\underline\phi}_{(31)(32)})
\endaligned\end{equation}
on
$U(\frak p(3),\frak w^{(2)}_{\frak p(3)};(\frak o^{(32)\prime\prime},\mathcal T^{(32)\prime\prime});\frak A)$.
\end{enumerate}
\end{lem}
\begin{rem}
The statement
(1) above was proved at the beginning of this subsection.
The statements (2) and (3) above were proved by Proposition \ref{prop21333}.
The statement (4) above was proved by Proposition \ref{prop2117}.
So only the statement (5) is new in Lemma \ref{2141}.
\end{rem}
\begin{proof}[Lemma \ref{2141} $\Rightarrow$ Proposition \ref{compaticoochamain}]
Let $\frak w_{\frak p(2)}^{(1)}$ be the stabilization data at $\frak p(2)$ induced by
$\frak w_{\frak p(1)}$.
\par
We apply Lemma \ref{2141} to $\frak w_{\frak p(2)}^{(1)}$
and $\frak w_{\frak p(2)}^{(2)} = \frak w_{\frak p(2)}$.
We then obtain $\epsilon_{8}$.
This $\epsilon_{8}$ is $\epsilon_7$ in Proposition \ref{compaticoochamain}.
Suppose
$\frak p(3)  \in \mathcal M_{k+1,\ell}(\beta)$, $d(\frak p(2),\frak p(3)) < \epsilon_{8}$.
We obtain $\frak m_{\frak p(3)}^{(1)},  \frak m_{\frak p(3)}^{(2)}$ from Lemma \ref{2141} (1).
\par
Using the pair of stabilization data $(\frak w_{\frak p(1)},\frak w_{\frak p(2)}^{(1)})$
we obtain the coordinate change $(\phi_{1(21)},\hat\phi_{1(21)})$
by Proposition \ref{prop21333}.
\par
Using the pair of stabilization data $(\frak m_{\frak p(3)}^{(2)},\frak m_{\frak p(3)})$
we obtain the coordinate change $(\phi_{(32)3},\hat\phi_{(32)3})$
by Proposition \ref{prop2117}.
\par
Now by using Lemma \ref{2141} (5) we have
\begin{equation}\label{eq2337}
\aligned
&({\underline\phi}_{1(21)},\hat{\underline\phi}_{1(21)}) \circ
({\underline\phi}_{(21)(22)},\hat\phi_{(21)(22)})\circ
({\underline\phi}_{(22)(32)},\hat{\underline\phi}_{(22)(32)})\circ
({\underline\phi}_{(32)3},\hat{\underline\phi}_{(32)3})\\
&=
({\underline\phi}_{1(21)},\hat{\underline\phi}_{1(21)}) \circ
({\underline\phi}_{(21)(31)},\hat{\underline\phi}_{(21)(31)})\circ
({\underline\phi}_{(31)(32)},\hat{\underline\phi}_{(31)(32)})\circ
({\underline\phi}_{(32)3},\hat{\underline\phi}_{(32)3})
\endaligned
\end{equation}
in a neighborhood of $[\frak p(3)]$.
\par
By definition of
$(\phi_{12},\hat{\phi}_{12})$, $(\phi_{23},\hat{\phi}_{23})$
given in the proof of Corollary \ref{cchanges}, we have
\begin{equation}
({\underline\phi}_{1(21)},\hat{\underline\phi}_{1(21)}) \circ
({\underline\phi}_{(21)(22)},\hat{\underline\phi}_{(21)(22)})
= ({\underline\phi}_{12},\hat{\underline\phi}_{12})
\end{equation}
and
\begin{equation}
({\underline\phi}_{(22)(32)},\hat{\underline\phi}_{(22)(32)})\circ
({\underline\phi}_{(32)3},\hat{\underline\phi}_{(32)3})
=
({\underline\phi}_{23},\hat{\underline\phi}_{23}).
\end{equation}
On the other hand, by Lemma \ref{lem21409}
$(\phi_{1(21)},\hat\phi_{1(21)}) \circ
(\phi_{(21)(31)},\hat\phi_{(21)(31)})$
is the coordinate change given by Proposition \ref{prop21333}.
By Lemma \ref{lem125},
$(\phi_{(31)(32)},\hat\phi_{(31)(32)})\circ
(\phi_{(32)3},\hat\phi_{(32)3})$
is the coordinate change given by  Proposition \ref{prop2117}.
Therefore,
by the definition given in the proof of Corollary \ref{cchanges},
\begin{equation}\label{eq2340}
\aligned
({\underline\phi}_{13},\hat{\underline\phi}_{13}) =
({\underline\phi}_{1(21)},\hat{\underline\phi}_{1(21)}) &\circ
({\underline\phi}_{(21)(31)},\hat{\underline\phi}_{(21)(31)})\\
&\circ
({\underline\phi}_{(31)(32)},\hat{\underline\phi}_{(31)(32)})\circ
({\underline\phi}_{(32)3},\hat{\underline\phi}_{(32)3}).
\endaligned
\end{equation}
Proposition \ref{compaticoochamain} follows from
(\ref{eq2337})-(\ref{eq2340}).
\end{proof}
\begin{proof}[Proof of Lemma \ref{2141}]
By definition, the coordinate change
$(\phi_{(21)(22)},\hat\phi_{(21)(22)})$ is a composition of finitely many
coordinate changes that are one of the types 1,3,4. (The notion of coordinate changes of
type 1,3,4 is defined right before Lemma \ref{lem2121}.)
Therefore it suffices to prove the lemma in the case when
$(\phi_{(21)(22)},\hat\phi_{(21)(22)})$ is one of types 1,3,4.
We prove each of those cases below.
\par\medskip
\noindent{\bf Case 1}:
$(\phi_{(21)(22)},\hat\phi_{(21)(22)})$ is of type 1.
\par
We use the notation in the proof of Proposition \ref{prop2117}
with $\frak p$ being replaced by $\frak p(2)$ or $\frak p(3)$.
\par
By Lemma \ref{lemma21118} we have
\begin{equation}\label{lem2118formula23}
\mathcal M^{\frak w_{\frak p(2)}^{(2) -}}_{k+1,(\ell,\ell_{\frak p},(\ell_c))}(\beta;\frak p(2);\frak A)
_{\epsilon'_{0,2},\vec{\mathcal T}^{(2) \prime}}
\subset
\mathcal M^{\frak w_{\frak p(2)}^{(1)}}_{k+1,(\ell,\ell_{\frak p},(\ell_c))}(\beta;\frak p(2);\frak A)
_{\epsilon_{0,2},\vec{\mathcal T}^{(1)}}.
\end{equation}
(Here we replace $\epsilon_{0}, \epsilon'_{0}$ in (\ref{lem2118formula}) by
$\epsilon_{0,2}, \epsilon'_{0,2}$. We also put $\ell_{\frak p} = \ell_{\frak p(2)}
= \ell_{\frak p(3)}$.)
Also by Lemma \ref{lemma21118} we have
\begin{equation}\label{lem2118formula233}
\mathcal M^{\frak w_{\frak p(3)}^{(2) -}}_{k+1,(\ell,\ell_{\frak p},(\ell_c))}(\beta;\frak p(3);\frak A)
_{\epsilon'_{0,3},\vec{\mathcal T}^{(3) \prime}}
\subset
\mathcal M^{\frak w_{\frak p(3)}^{(1)}}_{k+1,(\ell,\ell_{\frak p},(\ell_c))}(\beta;\frak p(3);\frak A)
_{\epsilon_{0,3},\vec{\mathcal T}^{(3)}}.
\end{equation}
(Here we replace $\epsilon_{0}, \epsilon'_{0}$ in (\ref{lem2118formula}) by
$\epsilon_{0,3}, \epsilon'_{0,3}$.)
\par
By the definition of type 1 we use the same codimension 2 submanifolds
to put the transversal constraint. Therefore we have
\begin{equation}\label{lem2118formula23tra}
\mathcal M^{\frak w_{\frak p(2)}^{(2) -}}_{k+1,(\ell,\ell_{\frak p},(\ell_c))}(\beta;\frak p(2);\frak A)
_{\epsilon'_{0,2},\vec{\mathcal T}^{(2) \prime}}^{\rm trans}
\subset
\mathcal M^{\frak w_{\frak p(2)}^{(1)}}_{k+1,(\ell,\ell_{\frak p},(\ell_c))}(\beta;\frak p(2);\frak A)
_{\epsilon_{0,2},\vec{\mathcal T}^{(2)}}^{\rm trans}
\end{equation}
and
\begin{equation}\label{lem2118formula233tra}
\mathcal M^{\frak w_{\frak p(3)}^{(2) -}}_{k+1,(\ell,\ell_{\frak p},(\ell_c))}(\beta;\frak p(3);\frak A)
_{\epsilon'_{0,3},\vec{\mathcal T}^{(3) \prime}}^{\rm trans}
\subset
\mathcal M^{\frak w_{\frak p(3)}^{(1)}}_{k+1,(\ell,\ell_{\frak p},(\ell_c))}(\beta;\frak p(3);\frak A)
_{\epsilon_{0,3},\vec {\mathcal T}^{(3)}}^{\rm trans}.
\end{equation}
On the other hand, by (\ref{2332}) we have
\begin{equation}\label{2332ss}
\mathcal M^{\frak w_{\frak p(3)}^{(1)}}_{k+1,(\ell,\ell_{\frak p},(\ell_c))}(\beta;\frak p(3);\frak A)_{\epsilon_{0,3},\vec {\mathcal T}^{(3)}}^{\rm trans}
\subset
\mathcal M^{\frak w_{\frak p(2)}^{(1)}}_{k+1,(\ell,\ell_{\frak p},(\ell_c))}(\beta;\frak p(2);\frak A)_{\epsilon_{0,2},\vec {\mathcal T}^{(2)}}^{\rm trans}.
\end{equation}
Note that the stabilization data $\frak w_{\frak p(2)}^{(2) -}$ and $\frak w_{\frak p(3)}^{(2) -}$ appearing in
(\ref{lem2118formula23tra}) and (\ref{lem2118formula233tra}) are obtained by extending the core of
the coordinate at infinity included in $\frak w_{\frak p(2)}^{(2)}$ and $\frak w_{\frak p(3)}^{(2)}$, respectively.
Therefore by further extending the core we may assume that $\frak w_{\frak p(3)}^{(2) -}$
is induced from $\frak w_{\frak p(2)}^{(2) -}$.
Therefore again by  (\ref{2332}) we have
\begin{equation}\label{2332s2s}
\mathcal M^{\frak w_{\frak p(3)}^{(2)-}}_{k+1,(\ell,\ell_{\frak p},(\ell_c))}(\beta;\frak p(3);\frak A)_{\epsilon''_{0,3},\vec{\mathcal T}^{(3)\prime\prime}}^{\rm trans}
\subset
\mathcal M^{\frak w_{\frak p(2)}^{(2)-}}_{k+1,(\ell,\ell_{\frak p},(\ell_c))}(\beta;\frak p(2);\frak A)_{\epsilon'_{0,2},\vec {\mathcal T}^{(2) \prime}}^{\rm trans}.
\end{equation}
By definition, the coordinate changes $\phi_{(21)(22)}$, $\phi_{(31)(32)}$, $\phi_{(21)(31)}$,
$\phi_{(22)(32)}$ are the inclusion maps  (\ref{lem2118formula23tra}), (\ref{lem2118formula233tra}),
(\ref{2332ss}) and (\ref{2332s2s}) in neighborhoods of $\frak p(2)$,  $\frak p(3)$, $\frak p(3)$, $\frak p(3)$,
respectively. The lemma is proved in this case.
\par\medskip
\noindent{\bf Case 2}:
Void.
\par\medskip
\noindent{\bf Case 4}:
$(\phi_{(21)(22)},\hat\phi_{(21)(22)})$ is of type 4.
\par
We have $\vec w_{\frak p(2)}^{(1)} \supset \vec w_{\frak p(2)}^{(2)}$.
Therefore $\vec w_{\frak p(3)}^{(1)} \supset \vec w_{\frak p(3)}^{(2)}$.
It follows that $(\phi_{(31)(32)},\hat\phi_{(31)(32)})$ is also of type 4.
We have the following commutative diagram.
\begin{equation}\label{cd347}
\begin{CD}
\mathcal M^{\frak w_{\frak p(2)}^{(1)}}_{k+1,(\ell,\ell_{\frak p(2)}^{(1)},(\ell_c))}(\beta;\frak p(2);\frak A)
_{\epsilon_{0},\vec{\mathcal T}^{(1)}} @ > {\frak{forget}_{\frak A,\frak A;\vec w_{\frak p(2)}^{(1)},\vec w_{\frak p(2)}^{(2)}}} >>
\mathcal M^{\frak w_{\frak p(2)}^{(2)}}_{k+1,(\ell,\ell_{\frak p(2)}^{(2)},(\ell_c))}(\beta;\frak p(2);\frak A)
_{\epsilon_{0},\vec{\mathcal T}^{(2)}} \\
@ AA{\subset}A @ AA{\subset}A\\
\mathcal M^{\frak w_{\frak p(3)}^{(1)}}_{k+1,(\ell,\ell_{\frak p(3)}^{(1)},(\ell_c))}(\beta;\frak p(3);\frak A)
_{\epsilon'_{0},\vec{\mathcal T}^{(1) \prime}} @ > {\frak{forget}_{\frak A,\frak A;\vec w_{\frak p(3)}^{(1)},\vec w_{\frak p(3)}^{(2)}}} >>
\mathcal M^{\frak w_{\frak p(3)}^{(2)}}_{k+1,(\ell,\ell_{\frak p(3)}^{(2)},(\ell_c))}(\beta;\frak p(3);\frak A)
_{\epsilon'_{0},\vec{\mathcal T}^{(2) \prime}}
\end{CD}
\end{equation}
We note that we use the same codimension 2 submanifold to put transversal constraint.
Therefore (\ref{cd347}) induces:
\begin{equation}\label{cd348}
\begin{CD}
\mathcal M^{\frak w_{\frak p(2)}^{(1)}}_{k+1,(\ell,\ell_{\frak p(2)}^{(1)},(\ell_c))}(\beta;\frak p(2);\frak A)
_{\epsilon_{0},\vec{\mathcal T}^{(1)}}^{\rm trans} @ > {\frak{forget}_{\frak A,\frak A;\vec w_{\frak p(2)}^{(1)},\vec w_{\frak p(2)}^{(2)}}} >>
\mathcal M^{\frak w_{\frak p(2)}^{(2)}}_{k+1,(\ell,\ell_{\frak p(2)}^{(2)},(\ell_c))}(\beta;\frak p(2);\frak A)
_{\epsilon_{0},\vec{\mathcal T}^{(2)}}^{\rm trans} \\
@ AA{\subset}A @ AA{\subset}A\\
\mathcal M^{\frak w_{\frak p(3)}^{(1)}}_{k+1,(\ell,\ell_{\frak p(3)}^{(1)},(\ell_c))}(\beta;\frak p(3);\frak A)
_{\epsilon'_{0},\vec{\mathcal T}^{(1) \prime}}^{\rm trans} @ > {\frak{forget}_{\frak A,\frak A;\vec w_{\frak p(3)}^{(1)},\vec w_{\frak p(3)}^{(2)}}} >>
\mathcal M^{\frak w_{\frak p(3)}^{(2)}}_{k+1,(\ell,\ell_{\frak p(3)}^{(2)},(\ell_c))}(\beta;\frak p(3);\frak A)
_{\epsilon'_{0},\vec{\mathcal T}^{(2) \prime}}^{\rm trans}
\end{CD}
\end{equation}
The commutativity of (\ref{cd348}) is Lemma \ref{2141} in this case.
\par\medskip
\noindent{\bf Case 3}:
$(\phi_{(21)(22)},\hat\phi_{(21)(22)})$ is of type 3.
\par
We obtain the following commutative diagram in the same way.

\begin{equation}\label{cd349}
\begin{CD}
\mathcal M^{\frak w_{\frak p(2)}^{(1)}}_{k+1,(\ell,\ell_{\frak p(2)}^{(1)},(\ell_c))}(\beta;\frak p(2);\frak A)
_{\epsilon_{0},\vec{\mathcal T}^{(1)}}^{\rm trans} @ < {\frak{forget}_{\frak A,\frak A;\vec w_{\frak p(2)}^{(2)},\vec w_{\frak p(2)}^{(1)}}} <<
\mathcal M^{\frak w_{\frak p(2)}^{(2)}}_{k+1,(\ell,\ell_{\frak p(2)}^{(2)},(\ell_c))}(\beta;\frak p(2);\frak A)
_{\epsilon_{0},\vec{\mathcal T}^{(2)}}^{\rm trans} \\
@ AA{\subset}A @ AA{\subset}A\\
\mathcal M^{\frak w_{\frak p(3)}^{(1)}}_{k+1,(\ell,\ell_{\frak p(3)}^{(1)},(\ell_c))}(\beta;\frak p(3);\frak A)
_{\epsilon'_{0},\vec{\mathcal T}^{(1) \prime}}^{\rm trans} @ < {\frak{forget}_{\frak A,\frak A;\vec w_{\frak p(3)}^{(2)},\vec w_{\frak p(3)}^{(1)}}} <<
\mathcal M^{\frak w_{\frak p(3)}^{(2)}}_{k+1,(\ell,\ell_{\frak p(3)}^{(2)},(\ell_c))}(\beta;\frak p(3);\frak A)
_{\epsilon'_{0},\vec{\mathcal T}^{(2) \prime}}^{\rm trans}
\end{CD}
\end{equation}
All the above arrows are diffeomorphisms locally. This implies the lemma in this case.
The proof of Lemma \ref{2141} is complete.
\end{proof}
The proof of Proposition  \ref{compaticoochamain} is complete.
\end{proof}
\par\medskip
\section{Wrap-up of the construction of Kuranishi structure}
\label{kstructure}
In this subsection we complete the proof of Theorem \ref{existsKura}.
We will prove the case of $\mathcal M_{k+1,\ell}(\beta)$.
The case of $\mathcal M^{\rm cl}_{\ell}(\alpha)$ is the same.
\par
In this subsection we fix a stabilization data $\frak w_{\frak p}$ at $\frak p$ for each
$\frak p$ and
always use it. We also take $\frak A = \frak C(\frak p)$ unless otherwise
specified. So we omit them from the notation of Kuranishi chart.
We write $\frak d =(\frak o,\mathcal T)$. Thus we write
$$
(V(\frak p;\frak d),\mathcal E_{(\frak p;\frak d)},\frak s_{(\frak p;\frak d)},
\psi_{(\frak p;\frak d)})
$$
to denote our Kuranishi neighborhood.
\par
For simplicity of notation we denote by $\tilde\psi_{(\frak p;\frak d)}$ the composition
of $\psi_{(\frak p;\frak d)}$ and the projection
$\frak s_{(\frak p;\frak d)}^{-1}(0) \to \frak s_{(\frak p;\frak d)}^{-1}(0)/\Gamma_{\frak p}$.
\par
The next lemma is the main technical lemma we use for the construction.
\begin{lem}\label{lem2143}
There exist finite subsets $\frak P_j = \{\frak p(j,i) \mid i=1,\dots,N_j\} \subset \mathcal M_{k+1,\ell}(\beta)$
for $j=1,2,3$ and admissible $\frak d(j,1,i) > \frak d(j,2,i)$ for $j=1,2,3$, $i=1,\dots,N_j$
such that they satisfy the following properties.
\begin{enumerate}
\item
If $j=1,2,3$ then
$$
\bigcup_{i=1}^{N_j} \tilde\psi_{(\frak p(j,i);\frak d(j,2,i))}
(\frak s_{(\frak p(j,i);\frak d(j,2,i))}^{-1}(0))
= \mathcal M_{k+1,\ell}(\beta).
$$
\item
The following holds for $j > j'$. If
$$
\frak p(j,i) \in \tilde\psi_{(\frak p(j',i');\frak d(j',2,i'))}
(\frak s_{(\frak p(j',i');\frak d(j',2,i'))}^{-1}(0)),
$$
then there exists a coordinate change
$$
\phi_{(j',i'),(j,i)}:
V(\frak p(j,i);\frak d(j,1,i))
\to
V(\frak p(j',i');\frak d(j',1,i'))
$$
as in Corollary \ref{cchanges}.
\item
Let $j=1$ or $2$, $i_1,\dots,i_m \in \{1,\dots,N_{j+1}\}$.
Suppose
$$
\bigcap_{n=1}^{m} \tilde\psi_{(\frak p(j+1,i_n);\frak d(j+1,1,i_n))}
(\frak s_{(\frak p(j+1,i_n);\frak d(j+1,1,i_n))}^{-1}(0))\\
\ne \emptyset,
$$
then there exists $i$ independent of $n$ such that
$$
\frak p(j+1,i_n) \in \tilde\psi_{(\frak p(j,i);\frak d(j,2,i))}
(\frak s_{(\frak p(j,i);\frak d(j,2,i))}^{-1}(0))
$$
for any $n = 1,\dots,m$.
\item
Let $i_j \in \{1,\dots,N_j\}$.
If
$$
\frak p(3,i_3) \in \tilde\psi_{(\frak p(2,i_2);\frak d(2,2,i_2))}
(\frak s_{(\frak p(2,i_2);\frak d(2,2,i_2))}^{-1}(0))
$$
and
$$
\frak p(2,i_2) \in \tilde\psi_{(\frak p(1,i_1);\frak d(1,2,i_1))}
(\frak s_{(\frak p(1,i_1);\frak d(1,2,i_1))}^{-1}(0)),
$$
then there exists a coordinate change
$$
\phi_{(1,i_1),(3,i_3)}:
V(\frak p(3,i_3);\frak d(3,1,i_3))
\to
V(\frak p(1,i_1);\frak d(1,1,i_1))
$$
as in Corollary \ref{cchanges}.
Moreover we have
\begin{equation}\label{2451}
{\underline\phi}_{(1,i_1),(2,i_2)}\circ {\underline\phi}_{(2,i_2),(3,i_3)}
={\underline\phi}_{(1,i_1),(3,i_3)}
\end{equation}
everywhere on $U(\frak p(3,i_3);\frak d(3,1,i_3))
= V(\frak p(3,i_3);\frak d(3,1,i_3))/\Gamma_{\frak p(3,i_3)}$.
\end{enumerate}
\end{lem}
\begin{proof}
For each $\frak p \in \mathcal M_{k+1,\ell}(\beta)$, we take admissible
$\frak d(\frak p,1;1) >\frak d(\frak p,1;2) >\frak d(\frak p,1;3)$.
Then we have $\frak P_1 = \{\frak p(1,i) \mid i=1,\dots,N_1\}$
such that
\begin{equation}\label{2351}
\bigcup_{i=1}^{N_1} \tilde\psi_{(\frak p(1,i);\frak d(\frak p(1,i),1;3))}
(\frak s_{(\frak p(1,i);\frak d(\frak p(1,i),1;3))}^{-1}(0))
= \mathcal M_{k+1,\ell}(\beta).
\end{equation}
We put
$\frak d(1,1,i) = \frak d(\frak p(1,i),1;1)$, $\frak d(1,2,i) = \frak d(\frak p(1,i),1;2)$.
Then, since
\begin{equation}\label{2351ato}
\aligned
&\tilde\psi_{(\frak p(1,i);\frak d(\frak p(1,i),1;3))}
(\frak s_{(\frak p(1,i);\frak d(\frak p(1,i),1;3)))}^{-1}(0))\\
&\subset
\tilde\psi_{(\frak p(1,i);\frak d(1,2,i))}
(\frak s_{(\frak p(1,i);\frak d(1,2,i))}^{-1}(0)),
\endaligned
\end{equation}
Lemma \ref{lem2143} (1) hods for $j=1.$
\par\smallskip
For each $\frak p \in \mathcal M_{k+1,\ell}(\beta)$
we take an admissible $\frak d(\frak p,2;1)$ so that the following conditions hold.
\begin{conds}\label{conds144}
\begin{enumerate}
\item[(a)]
If
$
\frak p \in \tilde\psi_{(\frak p(1,i);\frak d(1,2,i))}
(\frak s_{(\frak p(1,i);\frak d(1,2,i))}^{-1}(0)),
$
then there exists a coordinate change
$$
\phi_{(1,i),(2,\frak p)}:
V(\frak p;\frak d(\frak p,2;1))
\to
V(\frak p(1,i);\frak d(1,1,i))
$$
as in Corollary \ref{cchanges}.
\item[(b)]
If
$$
\tilde\psi_{(\frak p;\frak d(\frak p,2;1))}
(\frak s_{(\frak p;\frak d(\frak p,2;1))}^{-1}(0))\cap
\tilde\psi_{(\frak p(1,i);\frak d(\frak p(1,i),1;3))}
(\frak s_{(\frak p(1,i);\frak d(\frak p(1,i),1;3))}^{-1}(0))
\ne \emptyset
$$
then
$$
\tilde\psi_{(\frak p;\frak d(\frak p,2;1))}
(\frak s_{(\frak p;\frak d(\frak p,2;1))}^{-1}(0))
\subseteq
\tilde\psi_{(\frak p(1,i);\frak d(\frak p(1,i),1;2))}
(\frak s_{(\frak p(1,i);\frak d(\frak p(1,i),1;2))}^{-1}(0)).
$$
\item[(c)]
Let $\epsilon_{9}(\frak p)$ be the positive number we define below.
If an element $\frak q \in \mathcal M_{k+1,\ell}(\beta)$
satisfies
$
\frak q \in \tilde{\psi}_{(\frak p;\frak d(\frak p,2;1))}
(\frak s_{(\frak p;\frak d(\frak p,2;1))}^{-1}(0)),
$
then $d(\frak p,\frak q) < \epsilon_{9}(\frak p)$.
\end{enumerate}
\end{conds}
Here $\epsilon_{9}(\frak p)$ is defined as follows.
For each $i=1,\dots,N_1$ we put
$\frak p(1) = \frak p(1,i)$, $\frak p(2) = \frak p$
and apply Proposition  \ref{compaticoochamain}.
We then obtain $\epsilon_{7}(i,\frak p)$.
We define
$$
\epsilon_{9}(\frak p) = \min \{ \epsilon_{7}(i,\frak p)
\mid i=1,\dots,N_1\}.
$$
The existence of such $\frak d(\frak p,2;1)$ is obvious.
Furthermore for each $\frak p \in \mathcal M_{k+1,\ell}(\beta)$,
we take $\frak d(\frak p,2;2),\frak d(\frak p,2;3)$
such that
$\frak d(\frak p,2;1) >\frak d(\frak p,2;2) >\frak d(\frak p,2;3)$.
Then we have $\frak P_2 = \{\frak p(2,i) \mid i=1,\dots,N_2\}$
such that
\begin{equation}\label{2352}
\bigcup_{i=1}^{N_2} \tilde\psi_{(\frak p(2,i);\frak d(\frak p(2,i),2;3))}
(\frak s_{(\frak p(2,i);\frak d(\frak p(2,i),2;3))}^{-1}(0))
= \mathcal M_{k+1,\ell}(\beta).
\end{equation}
We put
$\frak d(2,1,i) = \frak d(\frak p(2,i),2;1)$, $\frak d(2,2,i) = \frak d(\frak p(2,i),2;2)$.
Then (\ref{2352}) and $\frak d(\frak p,2;2) >\frak d(\frak p,2;3)$
imply Lemma \ref{lem2143} (1) for $j=2$.
Lemma \ref{lem2143} (2) for $(j,j')=(2,1)$ follows immediately from  Condition \ref{conds144}
(a).
\begin{sublem}\label{sublem2145}
Lemma \ref{lem2143} (3) holds for $j=1$.
\end{sublem}
\begin{proof}
Suppose
$$
\bigcap_{n=1}^m
\tilde\psi_{(\frak p(2,i_n);\frak d(2,1,i_n))}
(\frak s_{(\frak p(2,i_n);\frak d(2,1,i_n))}^{-1}(0))
\ne \emptyset.
$$
Then (\ref{2351}) implies that there exists $i$ such that
$$
\aligned
&\bigcap_{n=1}^m
\tilde\psi_{(\frak p(2,i_n);\frak d(2,1,i_n))}
(\frak s_{(\frak p(2,i_n);\frak d(2,1,i_n))}^{-1}(0))\\
&\cap
\tilde\psi_{(\frak p(1,i);\frak d(\frak p(1,i),1;3))}
(\frak s_{(\frak p(1,i);\frak d(\frak p(1,i),1;3))}^{-1}(0))
\ne \emptyset.
\endaligned
$$
Therefore Condition  \ref{conds144} (b) and (\ref{2351ato}) imply
$$
\tilde\psi_{(\frak p(2,i_n);\frak d(2,1,i_n))}
(\frak s_{(\frak p(2,i_n);\frak d(2,1,i_n))}^{-1}(0))\subset
\tilde\psi_{(\frak p(1,i);\frak d(1,2,i))}
(\frak s_{(\frak p(1,i);\frak d(1,2,i))}^{-1}(0))
$$
for any $n$.
In particular
$$
\frak p(2,i_n)  \in \tilde\psi_{(\frak p(1,i);\frak d(1,2,i))}
(\frak s_{(\frak p(1,i);\frak d(1,2,i))}^{-1}(0))
$$
as required.
\end{proof}
For each $\frak p \in \mathcal M_{k+1,\ell}(\beta)$
we take an admissible $\frak d(\frak p,3;1)$ so that the following conditions hold.
\begin{conds}\label{cond146}
\begin{enumerate}
\item[(a)]
If
$
\frak p \in \tilde\psi_{(\frak p(2,i);\frak d(2,2,i))}
(\frak s_{(\frak p(2,i);\frak d(2,2,i))}^{-1}(0)),
$
then there exists a coordinate change
$$
\phi_{(2,i),(3,\frak p)}:
V(\frak p;\frak d(\frak p,3;1))
\to
V(\frak p(2,i);\frak d(2,1,i))
$$
as in Corollary \ref{cchanges}.
\item[(b)]
If
$$
\tilde\psi_{(\frak p;\frak d(\frak p,3;1))}
(\frak s_{(\frak p;\frak d(\frak p,3;1))}^{-1}(0))\cap
\tilde\psi_{(\frak p(2,i);\frak d(\frak p(2,i),2;3))}
(\frak s_{(\frak p(2,i);\frak d(\frak p(2,i),2;3))}^{-1}(0))
\ne \emptyset,
$$
then
$$
\tilde\psi_{(\frak p;\frak d(\frak p,3;1))}
(\frak s_{(\frak p;\frak d(\frak p,3;1))}^{-1}(0))
\subseteq
\tilde\psi_{(\frak p(2,i);\frak d(\frak p(2,i),2;2))}
(\frak s_{(\frak p(2,i);\frak d(\frak p(2,i),2;2))}^{-1}(0)).
$$
\item[(c)]
Void
\item[(d)]
Let $(i_1,i_2)$ be an arbitrary pair of integers such that
$$
\aligned
\frak p &\in \tilde\psi_{(\frak p(2,i_2);\frak d(2,2,i_2))}
(\frak s_{(\frak p(2,i_2);\frak d(2,2,i_2))}^{-1}(0)),\\
\frak p(2,i_2) &\in \tilde\psi_{(\frak p(1,i_1);\frak d(1,2,i_1))}
(\frak s_{(\frak p(1,i_1);\frak d(1,2,i_1))}^{-1}(0)).
\endaligned
$$
Then
$$
\frak d(\frak p,3;1) < \frak d(i_1,i_2).
$$
Here the right hand side is defined below.
\item[(e)]
Under the same assumption as in (d),  there exists a coordinate change
$$
\phi_{(1,i_1),(3,\frak p)}:
V(\frak p;\frak d(\frak p,3;1))
\to
V(\frak p(1,i_1);\frak d(1,1,i_1))
$$
as in Corollary \ref{cchanges}.
\end{enumerate}
\end{conds}
The definition of $\frak d(i_1,i_2)$ is as follows.
We put $\frak p(1) = \frak p(1,i_1)$ and $\frak p(2) = \frak p(2,i_2)$
and $\frak p(3) = \frak p$.
We also put $(\frak o^{(1)}, \mathcal T^{(1)}) = \frak d(1,1,i_1)$,
$(\frak o^{(2)}, \mathcal T^{(2)}) = \frak d(2,1,i_2)$.
Using Condition \ref{conds144} (c)
we can apply Proposition \ref{compaticoochamain} to obtain
$(\frak o^{(3)}_0, \mathcal T^{(3)}_0)$, which we put
$\frak d(i_1,i_2)$.
\par
Existence of $\frak d(\frak p,3;1)$ is obvious.
Furthermore for each $\frak p\in \mathcal M_{k+1,\ell}(\beta)$, we take
$\frak d(\frak p,3;2)$ with
 $\frak d(\frak p,3;1) > \frak d(\frak p,3;2)$.
Then we have $\frak P_3 = \{\frak p(3,i) \mid i=1,\dots,N_3\}$
such that
\begin{equation}\label{23523}
\bigcup_{i=1}^{N_3} \tilde\psi_{(\frak p(3,i);\frak d(\frak p(3,i),3;2))}
(\frak s_{(\frak p(3,i);\frak d(\frak p(3,i),3;2)))}^{-1}(0))
= \mathcal M_{k+1,\ell}(\beta).
\end{equation}
We put
$\frak d(3,1,i) = \frak d(\frak p(3,i),3;1)$, $\frak d(3,2,i) = \frak d(\frak p(3,i),3;2)$.
\par
Now Lemma \ref{lem2143} (1) for $j=3$ follows from (\ref{23523}).
Lemma \ref{lem2143} (2) for $(j,j')=(3,2), (3,1)$ follows from Condition \ref{cond146} (a),(e).
The proof of Lemma \ref{lem2143} (3) for $j=2$ is the same as the proof of Sublemma \ref{sublem2145}.
\par
Finally Lemma \ref{lem2143} (4) is a consequence of Condition \ref{conds144} (c), Condition \ref{cond146} (d)(e),
and Proposition \ref{compaticoochamain}.
The proof of Lemma \ref{lem2143} is complete.
\end{proof}
\begin{proof}
[Proof of Theorem \ref{existsKura}]
We start the construction of a Kuranishi structure on $\mathcal M_{k+1,\ell}(\beta)$.
Let $\frak p \in \mathcal M_{k+1,\ell}(\beta)$.
There exists $i(\frak p) \in \{1,\dots,N_3\}$ such that
$$
\frak p \in \tilde\psi_{(\frak p(3,i(\frak p));\frak d(\frak p(3,i(\frak p)),3;2))}
(\frak s_{(\frak p(3,i(\frak p));\frak d(\frak p(3,i(\frak p)),3;2))}^{-1}(0)).
$$
We take any such $i(\frak p)$ and fix it.
Choose $\hat{\frak p} \in V(\frak p(3,i(\frak p));\frak d(3,1,i(\frak p)))$
such that
$$
\tilde\psi_{(\frak p(3,i(\frak p));\frak d(\frak p(3,i(\frak p)),3;2))}(\hat{\frak p}) =\frak p.
$$
We have an embedding
$$
\bigoplus_{c \in \frak C(\frak p)} \mathcal E_c \subset \mathcal E_{(\frak p(3,i(\frak p));\frak d(\frak p(3,i(\frak p)),3;2))}
$$
of vector bundles. (This is because $\frak C(\frak p) \subseteq \frak C(\frak p(3,i(\frak p)))$.)
We take a neighborhood $V_{\frak p}$ of $\hat{\frak p}$ in
the set
$$
W_{\frak p}^{(3)} = \bigg\{ \frak v \in V(\frak p(3,i(\frak p));\frak d(3,1,i(\frak p))) \mid
\frak s_{(\frak p(3,i(\frak p));\frak d(\frak p(3,i(\frak p)),3;2)))}(\frak v)
\in \bigoplus_{c \in \frak C(\frak p)} \mathcal E_c\bigg\}
$$
such that  $V_{\frak p}$ is $\Gamma_{\frak p}$ invariant. The sum
$\bigoplus_{c \in \frak C(\frak p)} \mathcal E_c$ defines a $\Gamma_{\frak p}$ equivariant vector bundle
on $V_{\frak p}$ that we denote by $E_{\frak p}$.
The restriction to $V_{\frak p}$ of the section
$\frak s_{(\frak p(3,i(\frak p));\frak d(\frak p(3,i(\frak p)),3;2))}$ and the map
$\tilde\psi_{(\frak p(3,i(\frak p));\frak d(\frak p(3,i(\frak p)),3;2))}$ (divided by $\Gamma_{\frak p(3,i(\frak p))}$)
is our $\frak s_{\frak p}$ and $\psi_{\frak p}$.
We can show easily that $(V_{\frak p},\Gamma_{\frak p},E_{\frak p},\frak s_{\frak p},\psi_{\frak p})$
is a Kuranishi chart of $\frak p$.
\par
We next define coordinate changes.
Let $\frak q \in \psi_{\frak p}(\frak s^{-1}_{\frak p}(0))$.
It implies $\frak C(\frak q) \subseteq \frak C(\frak p)$.
\par
We note that $i(\frak p)$ may be different from $i(\frak q)$.
On the other hand, we have
$$
\aligned
&\tilde\psi_{(\frak p(3,i(\frak p));\frak d(3,1,i(\frak p)))}
(\frak s_{(\frak p(3,i(\frak p));\frak d(3,1,i(\frak p)))}^{-1}(0))\\
&\cap
\tilde\psi_{(\frak p(3,i(\frak q));\frak d(3,1,i(\frak q)))}
(\frak s_{(\frak p(3,i_(\frak q));\frak d(3,1,i(\frak q)))}^{-1}(0))
\ne \emptyset.
\endaligned
$$
In fact, $\frak q$ is contained in the intersection.
Therefore by Lemma \ref{lem2143} (3),
there exists $i(\frak p,\frak q)$ such that
$$
\frak p(3,i(\frak p)), \frak p(3,i(\frak q))  \in \tilde\psi_{(\frak p(2,i(\frak p,\frak q));\frak d(2,2,i(\frak p,\frak q)))}
(\frak s_{(\frak p(2,i(\frak p,\frak q));\frak d(2,2,i(\frak p,\frak q)))}^{-1}(0)).
$$
Therefore by Lemma \ref{lem2143} (2),
we have coordinate changes:
$$
\phi_{(\frak p(2,i(\frak p,\frak q))(3,i(\frak p))} : V(\frak p(3,i(\frak p));\frak d(3,1,i(\frak p))) \to V(\frak p(2,i(\frak p,\frak q));\frak d(2,1,i(\frak p,\frak q)))
$$
and
$$
\phi_{(\frak p(2,i(\frak p,\frak q))(3,i(\frak q))} : V(\frak p(3,i(\frak q));\frak d(3,1,i(\frak q))) \to V(\frak p(2,i(\frak p,\frak q));\frak d(2,1,i(\frak p,\frak q))).
$$
We write them sometimes as $\phi_{(\frak p\frak q)\frak p}$, $\phi_{(\frak p\frak q)\frak q}$ for simplicity.
\par
By the compatibility of $\psi$ with coordinate changes,
$$
\frak q \in \phi_{(\frak p(2,i(\frak p,\frak q))(3,i(\frak p))}(V_{\frak p})
\cap  \phi_{(\frak p(2,i(\frak p,\frak q))(3,i(\frak q))} (V_{\frak q}).
$$
We consider
$$
W_{\frak q}^{(2)} =
\bigg\{ \frak v \in V(\frak p(2,i(\frak p,\frak q));\frak d(2,1,i(\frak p,\frak q)))
\,\,\bigg\vert\,\,
\frak s_{(\frak p(2,i(\frak p,\frak q));\frak d(2,1,i(\frak p,\frak q)))}(\frak v)
\in \bigoplus_{c \in \frak C(\frak q)} \mathcal E_c\bigg\}.
$$
Both
$
\phi_{(\frak p(2,i(\frak p,\frak q))(3,i(\frak p))}(V_{\frak p}) \cap W^{(2)}_{\frak q}
$
and
$\phi_{(\frak p(2,i(\frak p,\frak q))(3,i(\frak q))} (V_{\frak q})$
are open subsets of $W^{(2)}_{\frak q}$. This fact is proved by dimension counting and by  the fact that
$\phi_{(\frak p(2,i(\frak p,\frak q))(3,i(\frak p))}$
and
$\phi_{(\frak p(2,i(\frak p,\frak q))(3,i(\frak q))}$
are embeddings.
\par
We put
\begin{equation}
V_{\frak p\frak q}
= \phi_{(\frak p(2,i(\frak p,\frak q))(3,i(\frak q))}^{-1}
(\phi_{(\frak p(2,i(\frak p,\frak q))(3,i(\frak p))}(V_{\frak p}) \cap W_{\frak q} )
\end{equation}
and
\begin{equation}\label{eq23777}
\phi_{\frak p\frak q} = \phi_{(\frak p(2,i(\frak p,\frak q))(3,i(\frak p))}^{-1}
\circ \phi_{(\frak p(2,i(\frak p,\frak q))(3,i(\frak q))}.
\end{equation}
We can define $\hat{\phi}_{\frak p\frak q}$ by using
$\hat{\phi}_{(\frak p(2,i(\frak p,\frak q))(3,i(\frak p))}$ and
$\hat{\phi}_{(\frak p(2,i(\frak p,\frak q))(3,i(\frak q))}$.
We have thus constructed a coordinate change.
\par
We finally prove the compatibility of coordinate changes.
Let $\frak q \in \tilde\psi_{\frak p}(\frak s^{-1}_{\frak p}(0))$,
and $\frak r \in \tilde\psi_{\frak q}(\frak s^{-1}_{\frak q}(0))$.
We then obtain $i(\frak p,\frak q)$, $i(\frak p,\frak r)$, $i(\frak q,\frak r)$ as above.
\par
We note that
$$
\aligned
&\tilde\psi_{(\frak p(2,i(\frak p,\frak q));\frak d(2,2,i(\frak p,\frak q)))}
(\frak s_{(\frak p(2,i(\frak p,\frak q));\frak d(2,2,i(\frak p,\frak q)))}^{-1}(0))\\
&\supseteq
\tilde\psi_{(\frak p(3,i(\frak q));\frak d(3,1,i(\frak q)))}
(\frak s_{(\frak p(3,i(\frak q));\frak d(3,1,i(\frak q)))}^{-1}(0))\\
&\supseteq \tilde\psi_{\frak q}(\frak s_{\frak q}^{-1}(0))
\ni \frak r.
\endaligned
$$
Therefore
$$
\aligned
&\tilde\psi_{(\frak p(2,i(\frak p,\frak q));\frak d(2,2,i(\frak p,\frak q))}
(\frak s_{(\frak p(2,i(\frak p,\frak q));\frak d(2,2,i(\frak p,\frak q)))}^{-1}(0))\\
&\cap \tilde\psi_{(\frak p(2,i(\frak q,\frak r));\frak d(2,2,i(\frak q,\frak r)))}
(\frak s_{(\frak p(2,i(\frak p,\frak r));\frak d(2,2,i(\frak p,\frak r)))}^{-1}(0))\\
&\cap \tilde\psi_{(\frak p(2,i(\frak p,\frak r));\frak d(2,2,i(\frak p,\frak r)))}
(\frak s_{(\frak p(2,i(\frak p,\frak r));\frak d(2,2,i(\frak p,\frak r)))}^{-1}(0))
\endaligned
$$
is nonempty. Therefore Lemma \ref{lem2143} (2) and (3)
imply that there exists
$i(\frak p,\frak q,\frak r)$ such that we have coordinate changes:
$$
\aligned
{\underline\phi}_{(\frak p(1,i(\frak p,\frak q,\frak r))(2,i(\frak p,\frak q))} : &U(\frak p(2,i(\frak p,\frak q));\frak d(2,1,i(\frak p,\frak q)))\\
&\to  U(\frak p(1,i(\frak p,\frak q,\frak r));\frak d(1,1,i(\frak p,\frak q,\frak r)))\\
{\underline\phi}_{(\frak p(1,i(\frak p,\frak q,\frak r))(2,i(\frak q,\frak r))} : &U(\frak p(2,i(\frak q,\frak r));\frak d(2,1,i(\frak q,\frak r)))\\
&\to  U(\frak p(1,i(\frak p,\frak q,\frak r));\frak d(1,1,i(\frak p,\frak q,\frak r)))\\
{\underline\phi}_{(\frak p(1,i(\frak p,\frak q,\frak r))(2,i(\frak p,\frak r))} : &U(\frak p(2,i(\frak p,\frak r));\frak d(2,1,i(\frak p,\frak r)))\\
&\to  U(\frak p(1,i(\frak p,\frak q,\frak r));\frak d(1,1,i(\frak p,\frak q,\frak r))).
\endaligned
$$
\footnote{Here
$U(\frak p(2,i(\frak p,\frak q));\frak d(2,1,i(\frak p,\frak q)))
= V(\frak p(2,i(\frak p,\frak q));\frak d(2,1,i(\frak p,\frak q)))/
\Gamma_{\frak p(2,i(\frak p,\frak q))}$ etc..}
We write them as
${\underline\phi}_{(\frak p\frak q\frak r)(\frak p\frak q)}$, ${\underline\phi}_{(\frak p\frak q\frak r)(\frak q\frak r)}$,
${\underline\phi}_{(\frak p\frak q\frak r)(\frak p\frak r)}$.
By Lemma \ref{lem2143} (4) we obtain
$$
\aligned
{\underline\phi}_{(\frak p(1,i(\frak p,\frak q,\frak r))(3,i(\frak p))} : &U(\frak p(3,i(\frak p));\frak d(3,1,i(\frak p)))\\
&\to  U(\frak p(1,i(\frak p,\frak q,\frak r));\frak d(1,1,i(\frak p,\frak q,\frak r)))\\
{\underline\phi}_{(\frak p(1,i(\frak p,\frak q,\frak r))(3,i(\frak q))} : &U(\frak p(3,i(\frak q));\frak d(3,1,i(\frak q)))\\
&\to  U(\frak p(1,i(\frak p,\frak q,\frak r));\frak d(1,1,i(\frak p,\frak q,\frak r)))\\
{\underline\phi}_{(\frak p(1,i(\frak p,\frak q,\frak r))(3,i(\frak r))} : &U(\frak p(3,i(\frak r));\frak d(3,1,i(\frak r)))\\
&\to  U(\frak p(1,i(\frak p,\frak q,\frak r));\frak d(1,1,i(\frak p,\frak q,\frak r))).
\endaligned
$$
We write them as
${\underline\phi}_{(\frak p\frak q\frak r)\frak p}$,
${\underline\phi}_{(\frak p\frak q\frak r)\frak q}$, ${\underline\phi}_{(\frak p\frak q\frak r)\frak r}$.
\par
By Lemma \ref{lem2143}  (4) we have
$$
\aligned
&{\underline\phi}_{(\frak p\frak q\frak r)(\frak p\frak q)} \circ {\underline\phi}_{(\frak p\frak q)\frak p} = {\underline\phi}_{(\frak p\frak q\frak r)\frak p},
\qquad
{\underline\phi}_{(\frak p\frak q\frak r)(\frak p\frak q)} \circ {\underline\phi}_{(\frak p\frak q)\frak q} = {\underline\phi}_{(\frak p\frak q\frak r)\frak q}, \\
&\phi_{(\frak p\frak q\frak r)(\frak q\frak r)} \circ {\underline\phi}_{(\frak q\frak r)\frak q} = {\underline\phi}_{(\frak p\frak q\frak r)\frak q},
\qquad
{\underline\phi}_{(\frak p\frak q\frak r)(\frak q\frak r)} \circ {\underline\phi}_{(\frak q\frak r)\frak r} = {\underline\phi}_{(\frak p\frak q\frak r)\frak r}, \\
&{\underline\phi}_{(\frak p\frak q\frak r)(\frak p\frak r)} \circ {\underline\phi}_{(\frak p\frak r)\frak r} = {\underline\phi}_{(\frak p\frak q\frak r)\frak r},
\qquad
{\underline\phi}_{(\frak p\frak q\frak r)(\frak p\frak r)} \circ {\underline\phi}_{(\frak p\frak r)\frak p} = {\underline\phi}_{(\frak p\frak q\frak r)\frak p}.
\endaligned
$$
Now we calculate:
$$
\aligned
{\underline\phi}_{\frak p\frak q} \circ {\underline\phi}_{\frak q\frak r}
&= {\underline\phi}_{(\frak p\frak q)\frak p}^{-1}\circ  {\underline\phi}_{(\frak p\frak q)\frak q}\circ
{\underline\phi}_{(\frak q\frak r)\frak q}^{-1}\circ  \phi_{(\frak q\frak r)\frak r}\\
&= {\underline\phi}_{(\frak p\frak q\frak r)\frak p}^{-1}\circ {\underline\phi}_{(\frak p\frak q\frak r)(\frak p\frak q)} \circ {\underline\phi}_{(\frak p\frak q)\frak q}
\circ \phi_{(\frak q\frak r)\frak q}^{-1}\circ\phi_{(\frak p\frak q\frak r)(\frak q\frak r)}^{-1}\circ  {\underline\phi}_{(\frak p\frak q\frak r)\frak r}\\
&=\phi_{(\frak p\frak q\frak r)\frak p}^{-1}\circ\phi_{(\frak p\frak q\frak r)\frak r}\\
&=\phi_{(\frak p\frak r)\frak p}^{-1}
\circ{\underline\phi}_{(\frak p\frak q\frak r)(\frak p\frak r)}^{-1}\circ{\underline\phi}_{(\frak p\frak q\frak r)(\frak p\frak r)}\circ{\underline\phi}_{(\frak p\frak r)\frak r}\\
&= {\underline\phi}_{(\frak p\frak r)\frak p}^{-1}
\circ{\underline\phi}_{(\frak p\frak r)\frak r} = {\underline\phi}_{\frak p\frak r}.
\endaligned
$$
Note (\ref{2451}) holds everywhere on $U(\frak p(3,i_3);\frak d(3,1,i_3))$.
Therefore we can perform the above calculation everywhere on
${\underline\phi}_{\frak q\frak r}^{-1}(U_{\frak p\frak q}) \cap U_{\frak p\frak r}$. (The maps
appearing in the intermediate stage of the calculation are defined in larger domain.)
\par
The proof of the consistency of the bundle maps $\hat{{\underline\phi}}_{\frak p\frak q},$ $\hat{{\underline\phi}}_{\frak q\frak r}$, $\hat{{\underline\phi}}_{\frak p\frak r}$
is the same by using
$\hat{{\underline\phi}}_{(\frak p\frak q\frak r)\frak r}$ etc.
\par
The proof of Theorem \ref{existsKura} is now complete.
\end{proof}
\par\medskip

\section{Appendix: Proof of Proposition \ref{changeinfcoorprop}}
\label{proposss}
In this subsection we prove Propositions \ref{changeinfcoorprop}, \ref{reparaexpest} and Lemma \ref{changeinfcoorproppara}.
It seems likely that there are several different ways to prove them.
We prove the proposition by the alternating method similar to those in the proof of
Theorems \ref{gluethm1}, \ref{exdecayT}, \ref{gluethm3}, \ref{exdecayT33}.
\par
In view of Lemma \ref{lempsi12}, it suffices to prove them in the case $\frak x = \frak Y_0$.
So we assume it throughout this subsection.
\par
We start with describing the situation.
We consider the universal bundle (\ref{fibrationsigma}).
The base space $\frak V(\frak x_{\rm v})$ is a neighborhood of $\frak x_{\rm v}$ in the
Deligne-Mumford moduli space.
Suppose we have two coordinates at infinity, which we write $\frak w{(j)}$, $j=1,2$.
We denote the universal bundle (\ref{fibrationsigma}) over $\frak V(\frak x_{\rm v})$ that is a part of $\frak w{(j)}$
by
\begin{equation}\label{universal}
\pi^{(j)} : \frak M^{(j)}_{\frak x_{\rm v}} \to \frak V(\frak x_{\rm v}).
\end{equation}
Actually $\pi^{(1)} = \pi^{(2)}$ but we distinguish them.\footnote{To prove
Lemma \ref{changeinfcoorproppara}, we need to consider a parametrized family and
so the parameter $\xi$ should be added to many of the objects we define.
To simplify the notation we omit them.}
The fiber at the base point $\frak x_{\rm v}$ is written as $\Sigma^{(j)}_{\rm v}$ and the
fiber at $\rho_{\rm v} \in \frak V(\frak x_{\rm v})$ is written as $\Sigma^{\rho, (j)}_{\rm v}$.
\par
We have an isomorphism
\begin{equation}\label{universaliso}
\hat\varphi_{12} : \frak M^{(2)}_{\frak x_{\rm v}}  \to \frak M^{(1)}_{\frak x_{\rm v}}
\end{equation}
of fiber bundles that preserves fiberwise complex structures and marked points.
Such an isomorphism is unique since we assumed $\frak x_{\rm v}$ to be stable.
\par
By Definition \ref{coordinatainfdef} (5) we have a trivialization:
\begin{equation}\label{universaltrivia}
\varphi^{(j)}_{\rm v} : \Sigma^{(j)}_{\rm v} \times \frak V(\frak x_{\rm v}) \to \frak M^{(j)}_{\frak x_{\rm v}}.
\end{equation}
The map $\varphi^{(j)}_{\rm v}$ is a diffeomorphism of fiber bundles of $C^{\infty}$-class, and preserves the complex structure
on the neck (ends). Moreover it preserves $\Gamma_{\frak p}$-action and marked points.
\par
Let $\rho = (\rho_{\rm v})$.
The restriction of the composition $(\varphi^{(1)}_{\rm v})^{-1}  \circ \hat\varphi_{12}\circ \varphi^{(2)}_{\rm v}$
to the fiber at $\rho_{\rm v} \in \frak V(\frak x_{\rm v})$ becomes a diffeomorphism
\begin{equation}\label{3362}
u^{\rho}_{\rm v} : (\Sigma^{(2)}_{\rm v},j^{(2)}_{\rho}) \to (\Sigma^{(1)}_{\rm v},j^{(1)}_{\rho}).
\end{equation}
We note that $u^{\rho}_{\rm v}$ is a diffeomorphism and is biholomorphic in the neck region.
(Note that the complex structure of the neck region is fixed by the definition of
coordinate at infinity.)
It is also biholomorphic (everywhere) with respect to the family of complex structures, $j^{(1)}_{\rho}, j^{(2)}_{\rho}$
parametrized by $\rho$.
\par
The map (\ref{3362}) preserves the marked points and is $\Gamma_{\frak x}$-equivariant.
We also assume the image of the neck region by $u^{\rho}_{\rm v}$ is contained in the neck region.
(We can always assume so by extending the neck of the coordinate at infinity $\frak w(1)$ of the source.)
Hereafter we write $\Sigma^{\rho,(j)}_{\rm v} = (\Sigma^{(j)}_{\rm v},j^{(j)}_{\rho})$ in case we do not
need to write $j^{(j)}_{\rho}$ explicitly.
\begin{rem}
We fix a trivialization as a smooth fiber bundle since it is important to fix a parametrization to
study $\rho$ derivative of the $\rho$-parametrized family of maps from the fibers.
\end{rem}
\par
In (\ref{2256}) we introduced the map
\begin{equation}\nonumber
\frak v_{(\frak y_2,\vec T_2,\vec \theta_2)} : \Sigma_{(\frak y_2,\vec T_2,\vec \theta_2)}
\to \Sigma_{(\frak y_1,\vec T_1,\vec \theta_1)}.
\end{equation}
Here the marked bordered curves $\Sigma_{(\frak y_j,\vec T_j,\vec \theta_j)}$ ($j=1,2$)
are obtained by gluing $\Sigma^{(j)}_{\rm v}$ in a way parametrized by
$\frak y_j,\vec T_j,\vec \theta_j$.
The idea of the proof is to construct the map $\frak v_{(\frak y_2,\vec T_2,\vec \theta_2)}$
by gluing the maps $u^{\rho}_{\rm v}$ using the alternating method.
In this subsection we use the notation $u$, $\rho$ in place of $\frak v$, $\frak y$.
\par
We introduce several function spaces. Let
$$
\frak y = \rho = (\rho_{\rm v}) \in \prod_{{\rm v} \in C^{0}(\mathcal G_{\frak x})}\frak V(\frak x_{\rm v}).
$$
We write $\Sigma^{\rho, (j)}_{\vec T,\vec \theta}$ using the notation used in the gluing construction in Subsection \ref{glueing}.
\par
We use the decomposition (\ref{sigmayv}) and (\ref{neddecomposit})
with coordinate (\ref{neckcoordinate2}).
The domain (\ref{DemanAetc}) are also used.
We use the bump functions (\ref{eq201})-(\ref{e2delta}).
\par
On the function space
\begin{equation}\label{363form}
L^2_{m,\delta}(\Sigma^{\rho,(2)}_{\rm v};(u_{\rm v}^{\rho})^*T\Sigma^{\rho',(1)}_{\rm v}\otimes \Lambda^{01})
\end{equation}
we define the norm
\begin{equation}
\Vert s\Vert^2_{L^2_{m,\delta}} = \sum_{k=0}^m \int_{\Sigma_{\rm v}^{\rho}}
e_{\rm v,\delta} \vert \nabla^k s\vert^2 \text{\rm vol}_{\Sigma_{\rm v}^{\rho}}.
\end{equation}
We modify Definition \ref{Sobolev263} as follows.
\begin{defn}\label{Sobolev26322}
The Sobolev space
$$
L^2_{m+1,\delta}((\Sigma_{\rm v}^{\rho,(2)},\partial \Sigma_{\rm v}^{\rho,(2)});(u_{\rm v}^{\rho})^*T\Sigma_{\rm v}^{\rho',(1)},(u_{\rm v}^{\rho})^*T\partial \Sigma_{\rm v}^{\rho',(1)})
$$
consists of elements $(s,\vec v)$ with the following properties.
\begin{enumerate}
\item $\vec v = (v_{\rm e})$ where $\rm e$ runs on the set of edges of $\rm v$  and
$$
v_{\rm e} = c_1\frac{\partial}{\partial \tau_{\rm e}} + c_2\frac{\partial}{\partial t_{\rm e}}
$$ (in case $\rm e\in C^1_{\rm c}(\mathcal G)$) or
$$
v_{\rm e} = c\frac{\partial}{\partial \tau_{\rm e}}
$$
(in case $\rm e\in C^1_{\rm o}(\mathcal G)$).
Here $c,c_1,c_2 \in \R$.
\item The following norm is finite.
\begin{equation}\label{normformjulamulti}
\aligned
&\Vert (s,\vec v)\Vert^2_{L^2_{m+1,\delta}} \\= &\sum_{k=0}^{m+1} \int_{K_{\rm v}}
\vert \nabla^k s\vert^2 \text{\rm vol}_{\Sigma_i}+ \sum_{{\rm e} \text{: edges of $\rm v$}}\Vert v_{\rm e}\Vert^2\\
&+
\sum_{k=0}^{m+1}\sum_{{\rm e} \text{: edges of $\rm v$}} \int_{\text{e-th end}}  e_{\rm v,\delta}\vert \nabla^k(s - \text{\rm Pal}(v_{\rm e}))\vert^2 \text{\rm vol}_{\Sigma_{\rm v}^{\rho}}.
\endaligned
\end{equation}
Here $\text{\rm Pal}$ is defined by  the canonical trivialization of the tangent bundle
on the neck region.
\end{enumerate}
\end{defn}
In case ${\rm v}\in C^0_{\rm s}(\mathcal G_{\frak x})$ we use the function space $L^2_{m+1,\delta}(\Sigma_{{\rm v}}^{\rho,(2)};
(u_{\rm v}^{\rho})^*T\Sigma_{\rm v}^{\rho',(1)})$
in place of
$L^2_{m+1,\delta}((\Sigma_{\rm v}^{\rho,(2)},\partial \Sigma_{\rm v}^{\rho,(2)});(u_{\rm v}^{\rho})^*T\Sigma_{\rm v}^{\rho',(1)},(u_{\rm v}^{\rho})^*T\partial \Sigma_{\rm v}^{\rho',(1)})$.

We do not assume any condition similar to Definition \ref{Lhattaato1}
and put
\begin{equation}\label{consthomo0thcomplex}
\aligned
&L^2_{m+1,\delta}((\Sigma^{\rho,(2)},\partial \Sigma^{\rho,(2)})
;(u_{\rm v}^{\rho})^*T\Sigma^{\rho',(1)},(u^{\rho})^*T\partial \Sigma^{\rho',(1)})\\
=
&\bigoplus_{{\rm v}\in C^0_{\rm d}(\mathcal G_{\frak x})}
L^2_{m+1,\delta}((\Sigma_{{\rm v}}^{\rho,(2)},\partial \Sigma_{\rm v}^{\rho,(2)});(u_{\rm v}^{\rho})^*T\Sigma_{\rm v}^{\rho',(1)},
(u_{\rm v}^{\rho})^*T\partial \Sigma_{\rm v}^{\rho',(1)})\\
&\oplus \bigoplus_{{\rm v}\in C^0_{\rm s}(\mathcal G_{\frak x})}
L^2_{m+1,\delta}(\Sigma_{{\rm v}}^{\rho,(2)};(u_{\rm v}^{\rho})^*T\Sigma_{\rm v}^{\rho',(1)}).
\endaligned
\end{equation}
The sum of (\ref{363form}) over $\rm v$ is denoted by
$$
L^2_{m,\delta}(\Sigma^{\rho,(2)};(u^{\rho})^*T\Sigma^{\rho',(1)}\otimes \Lambda^{01}).
$$
\par
We next define weighted Sobolev norms for the sections of various bundles on
$\Sigma_{\vec T,\vec \theta}^{\rho,(2)}$.
Here $\Sigma_{\vec T,\vec \theta}^{\rho,(2)}$ was denoted by $\Sigma_{\vec T,\vec \theta}^{\rho}$
in Subsection \ref{glueing}.
Let
$$
u' : (\Sigma_{\vec T,\vec \theta}^{\rho,(2)},\partial \Sigma_{\vec T,\vec \theta}^{\rho,(2)}) \to
(\Sigma_{\vec T',\vec \theta'}^{\rho',(1)},\partial \Sigma_{\vec T',\vec \theta'}^{\rho',(1)})
$$
be a diffeomorphism that sends each neck region of the source to the corresponding
neck region of the target.
We first consider the case when all $T_{\rm e} \ne \infty$.
In this case $\Sigma_{\vec T,\vec \theta}^{\rho,(j)}$ is compact.
We consider an element
$$
s \in L^2_{m+1}((\Sigma_{\vec T,\vec \theta}^{\rho,(2)},\partial \Sigma_{\vec T,\vec \theta}^{\rho,(2)});(u')^*T\Sigma_{\vec T',\vec \theta'}^{\rho',(1)},
(u')^*T\partial \Sigma_{\vec T',\vec \theta'}^{\rho',(1)}).
$$
Since we take $m$ large the section $s$ is continuous. We take a point $(0,1/2)_{\rm e}$ in the $\rm e$-th neck.
So $s((0,1/2)_{\rm e}) \in T_{u'((0,1/2)_{\rm e})}\Sigma_{\vec T',\vec \theta'}^{\rho',(1)}$ is well-defined.
\par
We use a canonical trivialization of the tangent bundle in the
neck regions to define $\text{\rm Pal}$ below.
We put
\begin{equation}\label{3655multi}
\aligned
\Vert s\Vert^2_{L^2_{m+1,\delta}} = &\sum_{k=0}^{m+1}\sum_{\rm v} \int_{K_{\rm v}}
\vert \nabla^k s\vert^2 \text{\rm vol}_{\Sigma_{\rm v}^{\rho}}\\
&+
\sum_{k=0}^{m+1}\sum_{\rm e} \int_{\text{e-th neck}}  e_{\vec T,\delta}\vert \nabla^k(s - \text{\rm Pal}(s(0,1/2)_{\rm e}))\vert^2
dt_{\rm e}d\tau_{\rm e}
\\&+ \sum_{\rm e}\Vert s((0,1/2)_{\rm e}))\Vert^2.
\endaligned
\end{equation}
For a section
$
s \in L^2_{m}(\Sigma_{\vec T,\vec \theta}^{\rho,(2)};(u')^*T\Sigma_{\vec T',\vec \theta'}^{\rho',(1)}\otimes \Lambda^{01})
$
we define
\begin{equation}\label{366}
\Vert s\Vert^2_{L^2_{m,\delta}} = \sum_{k=0}^{m} \int_{\Sigma_T}
e_{T,\delta}\vert  \nabla^k s\vert^2 \text{\rm vol}_{\Sigma_T}.
\end{equation}
\par
We next consider the case when some of the edges $\rm e$ have infinite length, namely $T_{\rm e} = \infty$.
Let $C^{1,\rm{inf}}_{\rm o}(\mathcal G_{\frak x},\vec T)$ (resp. $C^{1,\rm{inf}}_{\rm c}(\mathcal G_{\frak x},\vec T)$) be the set of elements $\rm e$ in $C^{1}_{\rm o}(\mathcal G_{\frak x})$ (resp. $C^{1}_{\rm c}(\mathcal G_{\frak x})$)
with $T_{\rm e} = \infty$ and
$C^{1,\rm{fin}}_{\rm o}(\mathcal G_{\frak x},\vec T)$ (resp. $C^{1,\rm{fin}}_{\rm c}(\mathcal G_{\frak x},\vec T)$) be the set of elements $C^{1}_{\rm o}(\mathcal G_{\frak x})$ (resp. $C^{1}_{\rm c}(\mathcal G_{\frak x})$)
with $T_{\rm e} \ne \infty$.
Note the ends of $\Sigma_{\vec T,\vec \theta}^{\rho}$ correspond two to one to
$C^{1,\rm{inf}}_{\rm o}(\mathcal G_{\frak x},\vec T) \cup C^{1,\rm{inf}}_{\rm c}(\mathcal G_{\frak x},\vec T)$.
The ends that correspond to an element of $C^{1,\rm{inf}}_{\rm o}(\mathcal G_{\frak x},\vec T)$ is
$([-5T_{\rm e},\infty) \times [0,1]) \cup (-\infty,5T_{\rm e}] \times [0,1])$ and
the ends that correspond to $C^{1,\rm{inf}}_{\rm c}(\mathcal G_{\frak p},\vec T)$ is
$([-5T_{\rm e},\infty) \times S^1) \cup (-\infty,5T_{\rm e}] \times S^1)$.
We have a weight function $e_{\rm v,\delta}(\tau_{\rm e},t_{\rm e})$ on it.
\begin{defn}\label{Def3.32}
An element of
$$
L^2_{m+1,\delta}((\Sigma_{\vec T,\vec \theta}^{\rho,(2)},\partial \Sigma_{\vec T,\vec \theta}^{\rho,(2)});
(u')^*T\Sigma_{\vec T',\vec \theta'}^{\rho',(1)},(u')^*T\partial\Sigma_{\vec T',\vec \theta'}^{\rho',(1)})
$$
is a pair $(s,\vec v)$ such that:
\begin{enumerate}
\item $s$ is a section of
$(u')^*T\Sigma_{\vec T',\vec \theta'}^{\rho',(1)}$ on $\Sigma_{\vec T,\vec \theta}^{\rho,(2)}$ minus singular points $z_{\rm e}$
corresponding to the edges $\rm e$ with $T_{\rm e} = \infty$.
\item $s$ is locally of $L^2_{m+1}$ class.
\item On $\partial \Sigma_{\vec T,\vec \theta}^{\rho,(2)}$ the restriction of $s$ is in $(u')^*T\partial\Sigma_{\vec T',\vec \theta'}^{\rho',(1)}$.
\item $\vec v = (v_{\rm e})$ where ${\rm e}$ runs in $C^{1,{\rm inf}}({\mathcal G}_{\frak p},{\vec T})$
and $v_{\rm e}$ is as in Definition \ref{Sobolev26322} (1).
\item For each $\rm e$ with $T_{\rm e} = \infty$ the integral
\begin{equation}\label{intatinfedge356}
\aligned
&\sum_{k=0}^{m+1}\int_{0}^{\infty} \int_{t_{\rm e}} e_{\rm v,\delta}(\tau_{\rm e},t_{\rm e})\vert\nabla^k(s(\tau_{\rm e},t_{\rm e})
 - {\rm Pal}(v_{\rm e}))\vert^2 d\tau_{\rm e} dt_{\rm e}\\
&+ \sum_{k=0}^{m+1}\int^{0}_{-\infty} \int_{t_{\rm e}} e_{\rm v,\delta}(\tau_{\rm e},t_{\rm e})\vert\nabla^k(s(\tau_{\rm e},t_{\rm e})
 - {\rm Pal}(v_{\rm e}))\vert^2 d\tau_{\rm e} dt_{\rm e}
  \endaligned
\end{equation}
is finite. (Here we integrate over ${t_{\rm e}\in [0,1]}$ (resp. ${t_{\rm e}\in S^1}$) if  ${\rm e }\in C_{\rm o}^{1,{\rm inf}}({\mathcal G}_{\frak p},\vec T)$
(resp.  ${\rm e} \in C_{\rm c}^{1,{\rm inf}}({\mathcal G}_{\frak p},\vec T)$).
\item
The section $s$  vanishes at each marked points.
\end{enumerate}
We define
\begin{equation}\label{368}
\Vert (s,\vec v)\Vert^2_{L^2_{m+1,\delta}}
= (\ref{3655multi}) + \sum_{{\rm e} \in C^{1,{\rm inf}}({\mathcal G}_{\frak p},\vec T)}(\ref{intatinfedge356}) +
\sum_{{\rm e} \in C^{1,{\rm inf}}(\mathcal G_{\frak p},\vec T)}\Vert v_{{\rm e}}\Vert^2.
\end{equation}
An element of
$$
L^2_{m,\delta}(\Sigma_{\vec T,\vec \theta}^{\rho,(2)};(u')^*T\Sigma_{\vec T',\vec \theta'}^{\rho',(1)}\otimes \Lambda^{01})
$$
is a section $s$ of the bundle $(u')^*T\Sigma_{\vec T',\vec \theta'}^{\rho',(1)}\otimes \Lambda^{01}$ such that it is locally of  $L^2_{m}$-class and
\begin{equation}\label{369}
\aligned
&\sum_{k=0}^{m}\int_{0}^{\infty} \int_{t_{\rm e}} e_{\rm v,\delta}\vert\nabla^k s(\tau_{\rm e},t_{\rm e})\vert^2
 d\tau_{\rm e} dt_{\rm e}\\
&+ \sum_{k=0}^{m}\int^{0}_{-\infty} \int_{t_{\rm e}} e_{\rm v,\delta}\vert\nabla^k(s(\tau_{\rm e},t_{\rm e})
\vert^2 d\tau_{\rm e} dt_{\rm e}
  \endaligned
\end{equation}
is finite.
We define
\begin{equation}\label{370}
\Vert s\Vert^2_{L^2_{m,\delta}}
= (\ref{366}) + \sum_{{\rm e} \in C^{1,{\rm inf}}({\mathcal G}_{\frak p},\vec T)}(\ref{369}).
\end{equation}
\end{defn}
\par
For a subset $W$ of $\Sigma_{{\rm v}}^{\rho,(2)}$ or $\Sigma_{\vec T,\vec \theta}^{\rho,(2)}$ we define
$
\Vert s\Vert_{L^2_{m,\delta}(W\subset \Sigma_{{\rm v}}^{\rho,(2)})}
$,
$
\Vert s\Vert_{L^2_{m,\delta}(W\subset \Sigma_{\vec T,\vec \theta}^{\rho,(2)})}
$
by restricting the domain of the integration (\ref{3655multi}),
(\ref{366}), (\ref{368}) or (\ref{370}) to $W$.
\par
We consider maps $u_{\rm v}^{\rho} : (\Sigma_{\rm v}^{\rho,(2)},\partial \Sigma_{\rm v}^{\rho,(2)}) \to
(\Sigma_{\rm v}^{\rho,(1)},\partial \Sigma_{\rm v}^{\rho,(1)})$
in (\ref{3362}),
for all $\rm v$. We write $u^{\rho} = (u_{\rm v}^{\rho})$.
\par
We next define a vector space that corresponds to a fiber of the  `obstruction bundle' in our situation.
Let
$
u' : (\Sigma_{\rm v}^{\rho,(2)},\partial \Sigma^{\rho,(2)}_{\rm v}) \to
(\Sigma_{\rm v}^{\rho',(1)},\partial \Sigma_{\rm v}^{\rho',(1)})
$
be a diffeomorphism that sends each of the neck region of the source to the corresponding
neck region of the target.
We define
$$
E_{\rm v}^{\rho}(u') \subset \Gamma_0(K^{\rho,(2)}_{\rm v}, (u')^*T\Sigma^{\rho',(1)}_{\rm v}\otimes \Lambda^{01})
$$
as follows.
\par
We may identify $\frak V(\frak x_{\rm v})$ as an open subset of certain Euclidean space.
Let $\frak e_{\rm v} \in T_{\rho'_{\rm v}}\frak V(\frak x_{\rm v})$.
We define
\begin{equation}\label{fromdefto01}
\frak I^{\rho'}_{\rm v}(u',\frak e_{\rm v})
=
\left.\frac{d}{dt}(\overline\partial^{\rho,\rho' + t\frak e_{\rm v}}u_{\rm v}^{\rho})\right\vert_{t=0}.
\end{equation}
Here $\overline\partial^{\rho,\rho + t\frak e_{\rm v}}$ is the $\overline\partial$ operator
with respect to the complex structure $j^{(1)}_{\rho'+ t\frak e_{\rm v}}$ (on the target) and $j^{(2)}_{\rho}$
(on the source).
We thus obtain a map:
\begin{equation}\label{cokisdetan}
\frak I^{\rho'}_{\rm v}(u',\cdot) :   T_{\rho'_{\rm v}}\frak V(\frak x_{\rm v}) \to
L^2_{m,\delta}(\Sigma_{\rm v}^{\rho,(2)};(u')^*T\Sigma^{\rho',(1)}_{\rm v}\otimes \Lambda^{01}).
\end{equation}
Since the complex structure is independent of $\rho$ on the neck region,
the image of (\ref{cokisdetan}) is
contained in $\Gamma_0(K^{\rho,(2)}_{\rm v}, (u')^*T\Sigma^{\rho',(1)}_{\rm v}\otimes \Lambda^{01})$,
that is, the set of smooth sections supported on the interior of the core.
\begin{defn}
We denote by
$E^{\rho}_{\rm v}(u')$  the image of (\ref{cokisdetan}).
\end{defn}
We consider the linearization of the Cauchy-Riemann equation associated to the biholomorphic map $u'$ that is
\begin{equation}\label{CRatuuu}
\aligned
D_{u'}\overline \partial:
&L^2_{m+1,\delta}((\Sigma^{\rho,(2)}_{\rm v},\partial \Sigma^{\rho,(2)}_{\rm v})
;(u')^*T\Sigma_{\rm v}^{\rho',(1)},(u')^*T\partial \Sigma^{\rho',(1)}_{\rm v})\\
&\to L^2_{m,\delta}(\Sigma_{\rm v}^{\rho,(2)};(u')^*T\Sigma^{\rho',(1)}_{\rm v}\otimes \Lambda^{01}).
\endaligned\end{equation}
\begin{lem}\label{lem3511}
If $u'$ is sufficiently close to $u^{\rho}_{\rm v}$ then
the kernel of (\ref{CRatuuu}) is zero
and we have
\begin{equation}\label{cokistang}
{\rm Im}(D_{u'}\overline \partial) \oplus E^{\rho}_{\rm v}(u')
=
L^2_{m,\delta}(\Sigma_{\rm v}^{\rho,(2)};(u')^*T\Sigma^{\rho',(1)}_{\rm v}\otimes \Lambda^{01}).
\end{equation}
\end{lem}
\begin{proof}
We first consider the case $u' = u^{\rho}_{\rm v}$, that is a biholomorphic map. Then
the kernel is identified with the set of holomorphic vector fields on $\Sigma^{\rho,(2)}$ that vanish on the singular points and
marked points. Such a vector field is necessary zero by stability.
\par
By the standard result of deformation theory, the cokernel is identified with the deformation space of the complex structures,
since $u^{\rho}_{\rm v}$ is biholomorphic.
Therefore (\ref{cokistang}) holds.
\par
We then find that the conclusion holds if $u'$ is sufficiently close to $u^{\rho}_{\rm v}$
so that $D_{u'}\overline \partial$ is close to $D_{u^{\rho}_{\rm v}}\overline \partial$
in operator norm and $E^{\rho}_{\rm v}(u')$ is close to $E^{\rho}_{\rm v}(u^{\rho}_{\rm v})$, in the sense that we can choose their orthonormal basis
that are close to each other.
\end{proof}
\begin{rem}
`Sufficiently close' is a bit imprecise way to state the lemma.
In the case we apply the lemma, we can easily check  that the last part of the proof works.
\end{rem}
We next take a map
\begin{equation}\label{appendEdef}
{\rm E} : \{(z,v) \in T\Sigma^{(1)} \mid \vert v\vert \le \epsilon\} \to \Sigma^{(1)}
\end{equation}
such that
\begin{enumerate}
\item
${\rm E}(z,0) =  z$ and
$$
\left.\frac{d}{dt} {\rm E}(z,tv)\right\vert_{t=0} = v.
$$
\item
If $(z,v) \in T\partial \Sigma^{(1)}$ then  ${\rm E}(z,v) \in \partial \Sigma^{(1)}$.
\item ${\rm E}(z,v) =  z+v$ on the neck region.
\end{enumerate}
\par\medskip
Now we start the gluing construction.
Let $(\vec T,\vec\theta) \in (\vec T^{\rm o}_0,\infty] \times ((\vec T^{\rm c}_0,\infty] \times \vec S^1)$.
For $\kappa =0,1,2,\dots$,  we will define a series of maps
\begin{eqnarray}
u_{\vec T,\vec \theta,(\kappa)}^{\rho} &:& (\Sigma_{\vec T,\vec \theta}^{\rho,(2)},\partial\Sigma_{\vec T,\vec \theta}^{\rho,(2)})
\to (\Sigma_{\vec T^{(\kappa)},\vec \theta^{(\kappa)}}^{\rho_{(\kappa)},(1)},\partial\Sigma_{\vec T^{(\kappa)},\vec \theta^{(\kappa)}}^{\rho_{(\kappa)},(1)})
\\
\hat u_{{\rm v},\vec T,\vec \theta,(\kappa)}^{\rho} &:& (\Sigma_{\rm v}^{\rho,(2)},\partial \Sigma_{\rm v}^{\rho,(2)}) \to
(\Sigma_{\rm v}^{\rho_{(\kappa)},(1)},\partial \Sigma_{\rm v}^{\rho_{(\kappa)},(1)}),
\end{eqnarray}
(we will explain $\rho_{(\kappa)}$, $\vec T^{(\kappa)}$ and $\vec \theta^{(\kappa)}$ below) and elements
\begin{eqnarray}
\frak e^{\rho}_{{\rm v},\vec T,\vec \theta,(\kappa)} &\in& E_{\rm v}(\hat u_{{\rm v},\vec T,\vec \theta,(\kappa)}^{\rho})
\\
{\rm Err}^{\rho}_{{\rm v},\vec T,\vec \theta,(\kappa)}  &\in&
L^2_{m,\delta}(\Sigma_{\rm v}^{\rho,(2)};(\hat u_{{\rm v},\vec T,\vec \theta,(\kappa)}^{\rho})^*T\Sigma_{\rm v}^{\rho_{(\kappa)},(1)}\otimes \Lambda^{01}).
\end{eqnarray}
\par
Moreover we will define
$V^{\rho}_{\vec T,\vec \theta,{\rm v},(\kappa)}$ for ${\rm v} \in C^0(\mathcal G)$, $\Delta T^{\rho}_{\vec T,\vec \theta,(\kappa),{\rm v},{\rm e}}
\in \R$
for ${\rm e} \in C^1(\mathcal G)$ and
$\Delta \theta^{\rho}_{\vec T,\vec \theta,(\kappa),{\rm v},{\rm e}} \in \R$ for ${\rm e} \in C^1_{\rm c}(\mathcal G)$.
We put
$$
\aligned
&v^{\rho}_{\vec T,\vec \theta,(\kappa),{\rm v},{\rm e}} = \Delta T^{\rho}_{\vec T,\vec \theta,(\kappa),{\rm v},{\rm e}}\frac{\partial}{\partial \tau_{\rm e}},
\qquad \text{for ${\rm e} \in C^1_{\rm o}(\mathcal G)$},
\\
&v^{\rho}_{\vec T,\vec \theta,(\kappa),{\rm v},{\rm e}} = \Delta T^{\rho}_{\vec T,\vec \theta,(\kappa),{\rm v},{\rm e}}\frac{\partial}{\partial \tau_{\rm e}}
+ \Delta \theta^{\rho}_{\vec T,\vec \theta,(\kappa),{\rm v},{\rm e}}\frac{\partial}{\partial t_{\rm e}}
\qquad \text{for ${\rm e} \in C^1_{\rm c}(\mathcal G)$}.
\endaligned$$
The pair $((V^{\rho}_{\vec T,\vec \theta,{\rm v},(\kappa)}),(v^{\rho}_{\vec T,\vec \theta,
(\kappa),{\rm v},{\rm e}}))$
becomes an element of
$$
L^2_{m+1,\delta}((\Sigma_{\rm v}^{\rho,(2)},\partial \Sigma_{\rm v}^{\rho,(2)});(\hat u_{{\rm v},\vec T,\vec \theta,(\kappa-1)}^{\rho})^*T\Sigma_{\vec T^{(\kappa)},\vec \theta^{(\kappa)}}^{\rho_{(\kappa)},(1)},(\hat u_{{\rm v},\vec T,\vec \theta,(\kappa-1)}^{\rho})^*T\partial \Sigma_{\vec T^{(\kappa)},\vec \theta^{(\kappa)}}^{\rho_{(\kappa)},(1)}).
$$
\par
The vectors $\vec T^{(\kappa)}$ and $\vec \theta^{(\kappa)}$ are determined by
$\Delta T^{\rho}_{\vec T,\vec \theta,(1),{\rm v},{\rm e}}, \dots, \Delta T^{\rho}_{\vec T,\vec \theta,(\kappa-1),{\rm v},{\rm e}}$
and $\Delta \theta^{\rho}_{\vec T,\vec \theta,(1),{\rm v},{\rm e}}, \dots, \Delta \theta^{\rho}_{\vec T,\vec \theta,(\kappa-1),{\rm v},{\rm e}}$
as follows.
For each $\rm e$ let ${\rm v}_{\leftarrow}({\rm e})$ and ${\rm v}_{\rightarrow}({\rm e})$ be the vertices for which
$\rm e$ is outgoing (resp. incoming) edge. We put:
\begin{equation}\label{3407aa}
10T^{(\kappa)}_{\rm e} = 10T_{\rm e} - \sum_{a=0}^{\kappa}\Delta T^{\rho}_{\vec T,\vec \theta,(a),{\rm v}_{\leftarrow}({\rm e}),{\rm e}}
+ \sum_{a=0}^{\kappa}\Delta T^{\rho}_{\vec T,\vec \theta,(a),{\rm v}_{\rightarrow}({\rm e}),{\rm e}}
\end{equation}
\begin{equation}\label{3408aa}
\theta^{(\kappa)}_{\rm e} = \theta_{\rm e} + \sum_{a=0}^{\kappa}\Delta \theta^{\rho}_{\vec T,\vec \theta,(a),{\rm v}_{\leftarrow}({\rm e}),{\rm e}}
- \sum_{a=0}^{\kappa}\Delta \theta^{\rho}_{\vec T,\vec \theta,(a),{\rm v}_{\rightarrow}({\rm e}),{\rm e}}.
\end{equation}
\begin{rem}
As induction proceeds, we will modify the length of the neck region a bit
from $T_{\rm e}$ to $T^{(\kappa)}_{\rm e}$.
We also modify $\theta_{\rm e}$ (that is the parameter to tell
how much we twist the $S^1$
direction when we  glue the pieces to obtain our curve) to $\theta^{(\kappa)}_{\rm e}$.
\end{rem}
The elements $\rho_{(\kappa)} = (\rho_{{\rm v},(\kappa)})$ is defined from
$\frak e^{\rho}_{{\rm v},\vec T,\vec \theta,(\kappa)}$ inductively as follows.
\par
\begin{equation}\label{3383}
\frak I^{\rho_{(\kappa-1)}}_{\rm v}(\hat u_{{\rm v},\vec T,\vec \theta,(\kappa-1)}^{\rho},\rho_{{\rm v},(\kappa)}-\rho_{{\rm v},(\kappa-1)})
=
\frak e^{\rho}_{{\rm v},\vec T,\vec \theta,(\kappa)}.
\end{equation}
So $T^{(\kappa)}_{\rm e}, \theta^{(\kappa)}_{\rm e}$ and $\rho_{{\rm v},(\kappa)}$ {\it depend} on $\rho, \vec T,\vec \theta$.
\begin{rem}
The construction of these objects are very much similar to that of
Subsection \ref{glueing}.
Note that $(\Sigma^{(1)},\partial\Sigma^{(1)})$ plays the role of $(X,L)$ here.
(In fact $\partial\Sigma^{(1)}$ is a Lagrangian submanifold of $\Sigma^{(1)}$.)
However the construction here is different from one in Subsection \ref{glueing}
in the following two points.
\begin{enumerate}
\item We will  construct a map $u$ that not only  satisfies $\overline\partial u \equiv 0 \mod E^{\rho}_{\rm v}$
but is also a genuine holomorphic map. The linearized equation (\ref{CRatuuu}) is {\it not}
surjective. We will kill the cokernel by deforming the complex structure
of the target. Namely $\rho \ne \rho_{(\kappa)}$ in general.
\item
We do {\it not} require
$\Delta T^{\rho}_{\vec T,\vec \theta,(\kappa),{\rm v}_{\leftarrow}({\rm e}),{\rm e}}
= \Delta T^{\rho}_{\vec T,\vec \theta,(\kappa),{\rm v}_{\rightarrow}({\rm e}),{\rm e}}$ or
$\Delta \theta^{\rho}_{\vec \theta,\vec \theta,(\kappa),{\rm v}_{\leftarrow}({\rm e}),{\rm e}}
= \Delta \theta^{\rho}_{\vec \theta,\vec \theta,(\kappa),{\rm v}_{\rightarrow}({\rm e}),{\rm e}}$.
This condition corresponds to $D{\rm ev}_{\mathcal G_{\frak p}}(V,\Delta p) =0$
that we put in Definition \ref{Lhattaato1}. Here we did not put a similar condition in
(\ref{consthomo0thcomplex}).
Instead we deform the complex structure of the target again.
Namely $T^{(\kappa)}_{\rm e} \ne T_{\rm e}$, $\theta^{(\kappa)}_{\rm e} \ne \theta_{\rm e}$ in general.
\end{enumerate}
\end{rem}
Now we start the construction of the above objects
by induction on $\kappa$.
\par\medskip
\noindent{\bf Pregluing}:
Since $u^{\rho}_{\rm v} : \Sigma_{\rm v}^{\rho,(2)} \to \Sigma_{\rm v}^{\rho,(1)}$
is biholomorphic and sends the neck region to the corresponding neck region, there exists
$\Delta T^{\rho}_{\vec T,\vec \theta,(0),{\rm v},{\rm e}}
\in \R$
for ${\rm e} \in C^1(\mathcal G)$ and
$\Delta \theta^{\rho}_{\vec T,\vec \theta,(0),{\rm v},{\rm e}} \in \R$ for ${\rm e} \in C^1_{\rm c}(\mathcal G)$
such that
\begin{equation}
\vert
u^{\rho}_{\rm v}(\tau_{\rm e},t_{\rm e}) - (\tau_{\rm e}+\Delta T^{\rho}_{\vec T,\vec \theta,(0),{\rm v},{\rm e}},t_{\rm e}
+ \Delta \theta^{\rho}_{\vec T,\vec \theta,(0),{\rm v},{\rm e}})
\vert
\le
C_1e^{-\delta_1 \vert\tau_{\rm e}\vert}.
\end{equation}
Note that in case ${\rm e} \in C^1_{\rm o}(\mathcal G)$ we put
$\Delta \theta^{\rho}_{\vec T,\vec \theta,(0),{\rm v},{\rm e}} = 0$.
\par
We identify the $\rm e$-th neck region of $\Sigma_{\vec T^{(\kappa)},\vec \theta^{(\kappa)}}^{\rho_{(\kappa)},(2)}$
with
$$
[-5T_{\rm e} + \frak s\Delta T_{{\rm e},(\kappa)}^{\leftarrow}, 5T_{\rm e} + \frak s\Delta T_{{\rm e},(\kappa)} ^{\rightarrow} ]
\times [0,1] \,\,\text{\rm or $S^1$},
$$
where
$$
\aligned
\frak s\Delta T_{{\rm e},(\kappa)}^{\leftarrow}
&= \sum_{a=0}^{\kappa}\Delta T^{\rho}_{\vec T,\vec \theta,(a),{\rm v}_{\leftarrow}({\rm e}),{\rm e}},\\
\frak s\Delta T_{{\rm e},(\kappa)}^{\rightarrow}
&= \sum_{a=0}^{\kappa}\Delta T^{\rho}_{\vec T,\vec \theta,(a),{\rm v}_{\rightarrow}({\rm e}),{\rm e}}.
\endaligned
$$
We also denote
$$
\aligned
\frak s\Delta \theta_{{\rm e},(\kappa)} ^{\leftarrow}
&= \sum_{a=1}^{\kappa}\Delta \theta^{\rho}_{\vec T,\vec \theta,(a),{\rm v}_{\leftarrow}({\rm e}),{\rm e}}, \\
\frak s\Delta \theta_{{\rm e},(\kappa)} ^{\rightarrow}
&= \sum_{a=1}^{\kappa}\Delta \theta^{\rho}_{\vec T,\vec \theta,(a),{\rm v}_{\rightarrow}({\rm e}),{\rm e}}.
\endaligned
$$
We use the symbol $\tau_{\rm e}^{(\kappa)}$ as the coordinate of the first factor.
The symbol $t_{\rm e}^{(\kappa)}$ denotes the coordinate of the second factor that is given
by
$$
t_{\rm e}^{(\kappa)} = t_{\rm e} + \frak s\Delta \theta_{{\rm e},(\kappa)} ^{\leftarrow}
$$
in case ${\rm e} \in C^1_{\rm c}(\mathcal G_{\frak x})$. Here $t_{\rm e}$ is the
canonical coordinate of $S^1$. In case ${\rm e} \in C^1_{\rm o}(\mathcal G_{\frak x})$,
$t_{\rm e}^{(\kappa)} = t_{\rm e}$.
\par
We have
\begin{equation}\label{3394tt}
\tau_{\rm e}^{(\kappa)}
= \tau'_{\rm e} - 5T_{\rm e} + \frak s\Delta T_{{\rm e},(\kappa)}^{\leftarrow}
=\tau''_{\rm e} + 5T_{\rm e} + \frak s\Delta T_{{\rm e},(\kappa)}^{\rightarrow} .
\end{equation}
(Hence $\tau'_{\rm e} = \tau''_{\rm e} + 10T_{\rm e}
- \frak s\Delta T^{\leftarrow}_{{\rm e},(\kappa)}
+ \frak s\Delta T^{\rightarrow}_{{\rm e},(\kappa)}
= \tau''_{{\rm e},(\kappa)} + 10T_{{\rm e}}^{(\kappa)}$. See (\ref{3407aa}).)
\par
In case ${\rm e} \in C^1_{\rm c}(\mathcal G_{\frak x})$ we also have
\begin{equation}\label{3395tt}
t_{\rm e}^{(\kappa)}
= t'_{\rm e}  + \frak s\Delta \theta_{{\rm e},(\kappa)}^{\leftarrow}
=t''_{\rm e} -\theta_{\rm e} +
\frak s\Delta \theta_{{\rm e},(\kappa)}^{\rightarrow}.
\end{equation}
(Hence $t'_{\rm e} = t''_{\rm e} - \theta_{\rm e}^{(\kappa)}$. See (\ref{3408aa}).)
\par
We define the map $
{\rm id}^{\rho,\vec T,\vec \theta}_{{\rm e},(\kappa)}
$
from the $\rm e$-th neck of $\Sigma_{\vec T,\vec \theta}^{\rho,(2)}$ to
the $\rm e$-th neck of $\Sigma_{\vec T^{(\kappa)},\vec \theta^{(\kappa)}}^{\rho_{(\kappa)},(1)}$
by
\begin{equation}
{\rm id}^{\rho,\vec T,\vec \theta}_{{\rm e},(\kappa)} : (\tau_{\rm e},t_{\rm e}) \mapsto
(\tau_{\rm e}^{(\kappa)} ,t_{\rm e}^{(\kappa)} ) = (\tau_{\rm e},t_{\rm e}).
\end{equation}
\par
We now put:
\begin{equation}\label{2219aa}
u_{\vec T,\vec \theta,(0)}^{\rho} =
\begin{cases}
\chi_{{\rm e},\mathcal B}^{\leftarrow} (u_{{\rm v}_{\leftarrow}({\rm e})}^{\rho} - {\rm id}^{\rho,\vec T,\vec \theta}_{{\rm e},(0)}) +
\chi_{{\rm e},\mathcal A}^{\rightarrow} (u_{{\rm v}_{\rightarrow}({\rm e})}^{\rho} - {\rm id}^{\rho,\vec T,\vec \theta}_{{\rm e},(0)})
+ {\rm id}^{\rho,\vec T,\vec \theta}_{{\rm e},(0)}
& \text{on the ${\rm e}$-th neck} \\
u_{\rm v}^{\rho} & \text{on $K_{\rm v}$}.
\end{cases}
\end{equation}
\par\medskip
\noindent{\bf Step 0-4}:
We next define
\begin{equation}\label{222233}
{\rm Err}^{\rho}_{{\rm v},\vec T,\vec \theta,(0)}  =
\begin{cases}
\chi_{{\rm e},\mathcal X}^{\leftarrow} \overline\partial u_{\vec T,\vec \theta,(0)}^{\rho}
& \text{on the ${\rm e}$-th neck if $\rm e$ is outgoing} \\
\chi_{{\rm e},\mathcal X}^{\rightarrow}  \overline\partial u_{\vec T,\vec \theta,(0)}^{\rho}
& \text{on the ${\rm e}$-th neck if $\rm e$ is incoming} \\
0 & \text{on $K_{\rm v}$} .
\end{cases}
\end{equation}
\par\medskip
\noindent{\bf Step 1-1}:
Let ${\rm id}_{{\rm v},{\rm e}}$ be the identity map from the neck region of $\Sigma^{(2)}_{\rm v}$
to the neck region of $\Sigma^{(1)}_{\rm v}$.  (It does not coincide with $u^{\rho}_{\rm v}$ there.)
We set:
\begin{equation}
\Delta^{{\rm v}_{\leftarrow}({\rm e}),{\rm e}}_{\vec T,\vec \theta,(0)} = (\frak s\Delta T_{{\rm e},(0)}^{\leftarrow},
\frak s\Delta \theta_{{\rm e},(0)}^{\leftarrow}),
\quad
\Delta^{{\rm v}_{\rightarrow}({\rm e}),{\rm e}}_{\vec T,\vec \theta,(0)} = (\frak s\Delta T_{{\rm e},(0)}^{\rightarrow},
\frak s\Delta \theta_{{\rm e},(0)}^{\rightarrow}).
\end{equation}
(In case ${\rm e} \in C^1_{\rm o}(\mathcal G_{\frak x})$ we set
$\frak s\Delta \theta_{{\rm e},(0)}^{\leftarrow} = \frak s\Delta \theta_{{\rm e},(0)}^{\rightarrow}=0$.)
We then define
\begin{equation}
{\rm id}_{{\rm v},{\rm e}}^{\vec T,\vec \theta,(0)}
= {\rm id}_{{\rm v},{\rm e}} + \Delta^{{\rm v},{\rm e}}_{\vec T,\vec \theta,(0)}.
\end{equation}
Now, we put
\begin{equation}
\aligned
&\hat u^{\rho}_{{\rm v},\vec T,\vec \theta,(0)}(z) \\
&=
\begin{cases} \chi_{{\rm e},\mathcal B}^{\leftarrow}(\tau_{\rm e}-T_{\rm e},t_{\rm e})
&\!\!\!\!\!\!u^{\rho}_{\vec T,\vec \theta,(0)}(\tau_{\rm e},t_{\rm e})
+ \chi_{{\rm e},\mathcal B}^{\rightarrow}(\tau_{\rm e}-T_{\rm e},t_{\rm e}){\rm id}_{{\rm v},{\rm e}}^{\vec T,\vec \theta,(0)} \\
&\text{if $z = (\tau_{\rm e},t_{\rm e})$ is on the $\rm e$-th neck that is outgoing} \\
\chi_{{\rm e},\mathcal A}^{\rightarrow}(\tau_{\rm e}-T_{\rm e},t_{\rm e})
&\!\!\!\!\!\!u^{\rho}_{\vec T,\vec \theta,(0)}(\tau,t)
+ \chi_{{\rm e},\mathcal A}^{\leftarrow}(\tau_{\rm e}-T_{\rm e},t_{\rm e}){\rm id}_{{\rm v},{\rm e}}^{\vec T,\vec \theta,(0)}\\
&\text{if $z = (\tau_{\rm e},t_{\rm e})$ is on the $\rm e$-th neck that is incoming} \\
 u^{\rho}_{{\rm v},\vec T,\vec \theta,(0)}(z)
&\text{if $z \in K_{\rm v}$.}
\end{cases}
\endaligned
\end{equation}
\begin{defn}
We define
$V^{\rho}_{\vec T,\vec \theta,{\rm v},(1)}$ for ${\rm v} \in C^0(\mathcal G_{\frak p})$
and real numbers $\Delta T^{\rho}_{\vec T,\vec \theta,(1),{\rm v}_{\leftarrow}({\rm e}),{\rm e}}$,
$\Delta T^{\rho}_{\vec T,\vec \theta,(1),{\rm v}_{\rightarrow}({\rm e}),{\rm e}}$ for ${\rm e} \in C^1(\mathcal G_{\frak p})$
and
$\Delta \theta^{\rho}_{\vec T,\vec \theta,(1),{\rm v}_{\leftarrow}({\rm e}),{\rm e}}$,
$\Delta \theta^{\rho}_{\vec T,\vec \theta,(1),{\rm v}_{\rightarrow}({\rm e}),{\rm e}}$ for ${\rm e} \in C^1_{\rm c}(\mathcal G_{\frak p})$
so that the following conditions are satisfied.
\begin{equation}
D_{\hat u^{\rho}_{{\rm v},\vec T,\vec \theta,(0)}}\overline{\partial}(V^{\rho}_{\vec T,\vec \theta,{\rm v},(1)})
- {\rm Err}^{\rho}_{{\rm v},\vec T,\vec \theta,(0)}\\
=
\frak e^{\rho}_{{\rm v},\vec T,\vec \theta,(0)} \in E_{\rm v}(\hat u_{{\rm v},\vec T,\vec \theta,(0)}^{\rho})
\end{equation}
and
\begin{equation}
\aligned
&\lim_{\tau_{\rm e} \to \infty} \left(V^{\rho}_{\vec T,\vec \theta,{\rm v}_{\leftarrow}({\rm e}),(1)}(\tau_{\rm e},t_{\rm e}) -
\Delta T^{\rho}_{\vec T,\vec \theta,(1),{\rm v}_{\leftarrow}({\rm e}),{\rm e}}\frac{\partial}{\partial \tau_{\rm e}}\right) = 0,\\
&\lim_{\tau_{\rm e} \to -\infty} \left(V^{\rho}_{\vec T,\vec \theta,{\rm v}_{\rightarrow}({\rm e}),(1)}(\tau_{\rm e},t_{\rm e}) -
\Delta T^{\rho}_{\vec T,\vec \theta,(1),{\rm v}_{\rightarrow}({\rm e}),{\rm e}}\frac{\partial}{\partial \tau_{\rm e}}\right) = 0,
\endaligned
\end{equation}
if ${\rm e} \in C_{\rm o}^1(\mathcal G_{\frak p})$,
\begin{equation}
\aligned
&\lim_{\tau_{\rm e} \to \infty} \left(V^{\rho}_{\vec T,\vec \theta,{\rm v}_{\leftarrow}({\rm e}),(1)}(\tau_{\rm e},t_{\rm e}) -
\Delta T^{\rho}_{\vec T,\vec \theta,(1),{\rm v}_{\leftarrow}({\rm e}),{\rm e}}\frac{\partial}{\partial \tau_{\rm e}}
-
\Delta \theta^{\rho}_{\vec T,\vec \theta,(1),{\rm v}_{\leftarrow}({\rm e}),{\rm e}}\frac{\partial}{\partial t_{\rm e}}\right) = 0,\\
&\lim_{\tau_{\rm e} \to -\infty} \left(V^{\rho}_{\vec T,\vec \theta,{\rm v}_{\rightarrow}({\rm e}),(1)}(\tau_{\rm e},t_{\rm e}) -
\Delta T^{\rho}_{\vec T,\vec \theta,(1),{\rm v}_{\rightarrow}({\rm e}),{\rm e}}\frac{\partial}{\partial \tau_{\rm e}}
-
\Delta \theta^{\rho}_{\vec T,\vec \theta,(1),{\rm v}_{\rightarrow}({\rm e}),{\rm e}}\frac{\partial}{\partial t_{\rm e}}\right) = 0,
\endaligned
\end{equation}
if ${\rm e} \in C^1_{\rm c}(\mathcal G_{\frak p})$.
\end{defn}
The unique existence of such objects is a consequence of Lemma \ref{lem3511}.
\par
We define $\rho_{(1)}$ by (\ref{3383}).
\par\medskip
\noindent{\bf Step 1-2}:
\begin{defn}
We define $u_{\vec T,\vec \theta,(1)}^{\rho}(z)$ as follows.
(Here $\rm E$ is as in (\ref{appendEdef}).)
\begin{enumerate}
\item If $z \in K_{\rm v}$ we put
\begin{equation}\label{3419}
u_{\vec T,\vec \theta,(1)}^{\rho}(z)
=
{\rm E} (u_{\vec T,\vec \theta,(0)}^{\rho},V^{\rho}_{\vec T,\vec \theta,{\rm v},(1)}(z)).
\end{equation}
\item
If $z  = (\tau_{\rm e},t_{\rm e}) \in [-5T_{\rm e},5T_{\rm e}]\times [0,1]$ or $S^1$,
we put
\begin{equation}\label{3420}
\aligned
&u_{\vec T,\vec \theta,(1)}^{\rho}(\tau_{\rm e},t_{\rm e}) \\
&=\chi_{{\rm v}_{\leftarrow}({\rm e}),\mathcal B}^{\leftarrow}(\tau_{\rm e},t_{\rm e}) (V^{\rho}_{\vec T,\vec \theta,{\rm v}_{\leftarrow}({\rm e}),(1)}
(\tau_{\rm e},t_{\rm e}) -
 (\Delta T_{{\rm e},(1)}^{\leftarrow},
\Delta \theta_{{\rm e},(1)}^{\leftarrow}))\\
&+\chi_{{\rm v}_{\rightarrow}({\rm e}),\mathcal A}^{\rightarrow}(\tau_{\rm e},t_{\rm e})(V^{\rho}_{\vec T,\vec \theta,{\rm v}_{\rightarrow}({\rm e}),(1)}(\tau_{\rm e},t_{\rm e})
-
(\Delta T_{{\rm e},(1)}^{\rightarrow},
\Delta \theta_{{\rm e},(1)}^{\rightarrow}))\\
&+u_{\vec T,\vec \theta,(0)}^{\rho}(\tau_{\rm e},t_{\rm e}).
\endaligned
\end{equation}
\end{enumerate}
Here we use the coordinate $(\tau_{\rm e}^{(1)},t_{\rm e}^{(1)})$
given in (\ref{3394tt}) and (\ref{3395tt}) for the {\it target}.
\end{defn}
\par
We remark that $\tau^{(0)}_{\rm e} = \tau^{(1)}_{\rm e} - \Delta T^{\leftarrow}_{{\rm e},(1)}$.
Therefore, in a neighborhood of $\{-5T_{\rm e}\} \times [0,1]\times S^1$,
(\ref{3419}) and (\ref{3420}) are consistent.
\par\medskip
\noindent{\bf Step 1-3}:
We recall that $\rho_{{\rm v},(1)}$ is defined by
\begin{equation}\label{338322}
\frak I^{\rho_{(0)}}_{\rm v}(\hat u_{{\rm v},\vec T,\vec \theta,(0)}^{\rho},\rho_{{\rm v},(1)}-\rho_{{\rm v},(0)})
=
\frak e^{\rho}_{{\rm v},\vec T,\vec \theta,(1)}.
\end{equation}
(Note $\rho_{{\rm v},(0)} = 0$.)
\par\medskip
\noindent{\bf Step 1-4}:
\begin{defn}
We put
\begin{equation}
{\rm Err}^{\rho}_{{\rm v},\vec T,\vec \theta,(1)}  =
\begin{cases}
\chi_{{\rm e},\mathcal X}^{\leftarrow} \overline\partial u_{\vec T,\vec \theta,(1)}^{\rho}
& \text{on ${\rm e}$-th neck if $\rm e$ is outgoing} \\
\chi_{{\rm e},\mathcal X}^{\rightarrow}  \overline\partial u_{\vec T,\vec \theta,(1)}^{\rho}
& \text{on ${\rm e}$-th neck if $\rm e$ is incoming} \\
\overline\partial u_{\vec T,\vec \theta,(1)}^{\rho}   & \text{on $K_{\rm v}$} .
\end{cases}
\end{equation}
We extend them by $0$ outside a compact set and will regard them as elements of the function space
$L^2_{m,\delta}(\Sigma^{\rho,(2)}_{\rm v};(\hat u^{\rho}_{{\rm v},\vec T,\vec \theta,(1)})^{*}T\Sigma_{\vec T^{(1)},\vec \theta^{(1)}}^{\rho_{(1)},(1)} \otimes \Lambda^{01})$,
where $\hat u^{\rho}_{{\rm v},\vec T,\vec \theta,(1)}$ will be defined in the next step.
\end{defn}
We thus come back to Step 2-1 and continue.
We obtain the following estimate by induction on $\kappa$.
We put $R_{\rm e} = 5T_{\rm e} + 1$.
\begin{lem}\label{expesgen1sec3}
There exist $T_m, C_{2,m},\dots, C_{8,m}, \epsilon_{1,m} > 0$ and $0<\mu<1$ such that
the following inequalities hold if $T_{\rm e}>T_m$ for all $\rm e$.
We put $T_{\rm min} = \min\{ T_{\rm e}\mid {\rm e} \in C^1(\mathcal G_{\frak p})\}$.
\begin{eqnarray}
\left\Vert
\left((V^{\rho}_{\vec T,\vec \theta,{\rm v},(\kappa)}),(v^{\rho}_{\vec T,\vec \theta,{\rm v},{\rm e},(\kappa)})\right)\right\Vert_{L^2_{m+1,\delta}
(\Sigma^{\rho,{(2)}}_{\rm v})}
&<& C_{2,m}\mu^{\kappa-1}e^{-\delta T_{\rm min}}, \label{form0182vv3}
\\
\left\Vert (v^{\rho}_{\vec T,\vec \theta,{\rm v},{\rm e},(\kappa)})\right\Vert
&<& C_{3,m}\mu^{\kappa-1}e^{-\delta T_{\rm min}}, \label{form0183vv3}
\\
\left\Vert  u_{\vec T,\vec \theta,(\kappa)}^{\rho}- u_{\vec T,\vec \theta,(0)}^{\rho}  \right\Vert_{L^2_{m+1,\delta}((K_{\rm v}^{(2)})^{+\vec R})}
&<& C_{4,m}e^{-\delta T_{\rm min}}, \label{form0184vv3}
\\
\left\Vert{\rm Err}^{\rho}_{{\rm v},\vec T,\vec \theta,(\kappa)} \right\Vert_{L^2_{m,\delta}(\Sigma^{\rho,{(2)}}_{\rm v})}
&<& C_{5,m}\epsilon_{1,m}\mu^{\kappa}e^{-\delta T_{\rm min}}, \label{form0185vv3}
\\
\left\Vert \frak e^{\rho} _{\vec T,\vec \theta,(\kappa)}\right\Vert_{L^2_{m}((K_{\rm v}^{(2)})^{+\vec R})}
&<& C_{6,m}\mu^{\kappa-1}e^{-\delta T_{\rm min}},
\label{form0186vv3vv3}
\\
\left\Vert \Delta T^{\rho}_{\vec T,\vec \theta,(\kappa),{\rm v},{\rm e}}
\right\Vert
&<& C_{7,m}\mu^{\kappa-1}e^{-\delta T_{\rm min}}, \label{Tconverges3}\\
\left\Vert
\Delta \theta^{\rho}_{\vec T,\vec \theta,(\kappa),{\rm v},{\rm e}}\right\Vert
&<& C_{8,m}\mu^{\kappa-1}e^{-\delta T_{\rm min}} \label{rhoconverges3}.
\end{eqnarray}
\end{lem}
The proof is the same as the proof of Proposition \ref{expesgen1} and so is omitted.
We note that (\ref{form0186vv3vv3}) and (\ref{3383}) imply
\begin{equation}\label{rhoconverges}
\left\Vert \rho_{(\kappa)} - \rho \right\Vert
< C_{9,m}\mu^{\kappa-1}e^{-\delta T_{\rm min}}.
\end{equation}
Therefore the limit
$$
\lim_{\kappa \to \infty} \rho_{(\kappa)} = \rho'(\rho,\vec T,\vec \theta)
$$
exists.
(\ref{Tconverges3}) and (\ref{rhoconverges3}) imply that
$$
\lim_{\kappa \to \infty} \frak s\Delta T^{\rho}_{\vec T,\vec \theta,(\kappa),{\rm v},{\rm e}} = \frak s\Delta T^{\rho}_{\vec T,\vec \theta,(\infty),{\rm v},{\rm e}}
$$
and
$$
\lim_{\kappa \to \infty} \frak s\Delta \theta^{\rho}_{\vec T,\vec \theta,(\kappa),{\rm v},{\rm e}}
= \frak s\Delta \theta^{\rho}_{\vec T,\vec \theta,(\infty),{\rm v},{\rm e}}
$$
converge.
We put
$$
\vec T'(\rho,\vec T,\vec \theta) = \vec T + \frak s\Delta \vec T^{\rho}_{\vec T,\vec \theta,(\infty)},
\quad
\vec \theta'(\rho,\vec T,\vec \theta) = \vec \theta + \frak s\Delta \vec \theta^{\rho}_{\vec T,\vec \theta,(\infty)}.
$$
Then
(\ref{form0184vv3}) implies that
$$
\lim_{\kappa \to \infty} u_{\vec T,\vec \theta,(\kappa)}^{\rho}
$$
converges to a map
$$
u_{\vec T,\vec \theta,(\infty)}^{\rho} :
(\Sigma_{\vec T,\vec \theta}^{\rho,(2)},\partial\Sigma_{\vec T,\vec \theta}^{\rho,(2)})
\to
(\Sigma_{\vec T'(\rho,\vec T,\vec \theta) ,\vec \theta'(\rho,\vec T,\vec \theta)}^{\rho'(\rho,\vec T,\vec \theta),(1)},
\partial\Sigma_{\vec T'(\rho,\vec T,\vec \theta) ,\vec \theta'(\rho,\vec T,\vec \theta)}^{\rho'(\rho,\vec T,\vec \theta),(1)})
$$
in $L^2_{m+1}$ topology.
(Note the union of $(K_{\rm v}^{(2)})^{+\vec R}$ for various
$\rm v$ covers $\Sigma_{\vec T,\vec \theta}^{\rho_{(\kappa)},(2)}$.)  The formula
(\ref{form0185vv3}) then implies that $u_{\vec T,\vec \theta,(\infty)}^{\rho}$ is a biholomorphic map.
\par
Therefore, using the notation in Proposition \ref{changeinfcoorprop} we have
\begin{equation}\label{form408}
{\overline{\Phi}}_{12}(\rho,\vec T,\vec\theta)
=
(\rho'(\rho,\vec T,\vec \theta),\vec T'(\rho,\vec T,\vec \theta') ,\vec \theta'(\rho,\vec T,\vec \theta)).
\end{equation}
Using the notation in Proposition \ref{reparaexpest} we have
\begin{equation}\label{3409}
\frak v_{(\rho,\vec T,\vec\theta)}
=
u_{\vec T,\vec \theta,(\infty)}^{\rho}.
\end{equation}
The $T_{\rm e}$ etc. derivative of the objects we constructed enjoy the following estimate.
\begin{lem}\label{expesgen2Tdevss}
There exist $T_m, C_{10,m}, \dots, C_{16,m}, \epsilon_{2,m} > 0$ and $0<\mu<1$ such that
the following inequalities hold if $T_{\rm e}>T_m$ for all $\rm e$.
\par
Let ${\rm e}_0 \in C^1(\mathcal G_{\frak p})$.
Then for each $\vec k_{T}$,  $\vec k_{\theta}$ we have
\begin{equation}
\aligned
&\left\Vert \nabla_{\rho}^n \frac{\partial^{\vert \vec k_{T}\vert}}{\partial T^{\vec k_{T}}}\frac{\partial^{\vert \vec k_{\theta}\vert}}{\partial \theta^{\vec k_{\theta}}}
\frac{\partial}{\partial T_{{\rm e}_0}}
\left((V^{\rho}_{\vec T,\vec \theta,{\rm v},(\kappa)}),(v^{\rho}_{\vec T,\vec \theta,{\rm v},{\rm e},(\kappa)})\right)
\right\Vert_{L^2_{m+1-\vert{\vec k_{T}}\vert - \vert{\vec k_{\theta}}\vert-n-1,\delta}(\Sigma_{\rm v}^{\rho,(2)})}\\
&< C_{10,m}\mu^{\kappa-1}e^{-\delta T_{{\rm e}_0}}, \label{form0182vv233}
\endaligned
\end{equation}
\begin{equation}
\displaystyle\left\Vert \nabla_{\rho}^n \frac{\partial^{\vert \vec k_{T}\vert}}{\partial T^{\vec k_{T}}}\frac{\partial^{\vert \vec k_{\theta}\vert}}{\partial \theta^{\vec k_{\theta}}}
\frac{\partial}{\partial T_{{\rm e}_0}}
(v^{\rho}_{\vec T,\vec \theta,{\rm v},{\rm e},(\kappa)})\right\Vert
< C_{11,m}\mu^{\kappa-1}e^{-\delta T_{{\rm e}_0}}, \label{form0183v2v33}
\end{equation}
\begin{equation}
\left\Vert \nabla_{\rho}^n \frac{\partial^{\vert \vec k_{T}\vert}}{\partial T^{\vec k_{T}}}\frac{\partial^{\vert \vec k_{\theta}\vert}}{\partial \theta^{\vec k_{\theta}}}
\frac{\partial}{\partial T_{{\rm e}_0}}  u_{\vec T,\vec \theta,(\kappa)}^{\rho} \right\Vert_{L^2_{m+1-\vert{\vec k_{T}}\vert - \vert{\vec k_{\theta}}\vert
-n-1,\delta}((K_{\rm v}^{(2))})^{+\vec R})}
< C_{12,m}e^{-\delta T_{{\rm e}_0}}, \label{form0184233}
\end{equation}
\begin{equation}
\aligned
&\left\Vert
\nabla_{\rho}^n \frac{\partial^{\vert \vec k_{T}\vert}}{\partial T^{\vec k_{T}}}\frac{\partial^{\vert \vec k_{\theta}\vert}}{\partial \theta^{\vec k_{\theta}}}
\frac{\partial}{\partial T_{{\rm e}_0}}
{\rm Err}^{\rho}_{{\rm v},\vec T,\vec \theta,(\kappa)} \right\Vert_{L^2_{m-\vert{\vec k_{T}}\vert - \vert{\vec k_{\theta}}\vert-n-1,\delta}
(\Sigma_{\rm v}^{\rho,{(2)}})}
\\
&< C_{13,m}\epsilon_{2,m}\mu^{\kappa}e^{-\delta  T_{{\rm e}_0}}, \label{form0185vv233}
\endaligned
\end{equation}
\begin{equation}
\left\Vert
\nabla_{\rho}^n\frac{\partial^{\vert \vec k_{T}\vert}}{\partial T^{\vec k_{T}}}\frac{\partial^{\vert \vec k_{\theta}\vert}}{\partial \theta^{\vec k_{\theta}}}
\frac{\partial}{\partial T_{{\rm e}_0}}
\frak e^{\rho} _{\vec T,\vec \theta,(\kappa)}\right\Vert_{L^2_{m-\vert{\vec k_{T}}\vert - \vert{\vec k_{\theta}}\vert-n-1}(K_{\rm v}^{(2)})}
< C_{14,m}\mu^{\kappa-1}e^{-\delta  T_{{\rm e}_0}},
\label{form0186vv233}
\end{equation}
\begin{equation}\label{415}
\left\Vert
\nabla_{\rho}^n \frac{\partial^{\vert \vec k_{T}\vert}}{\partial T^{\vec k_{T}}}\frac{\partial^{\vert \vec k_{\theta}\vert}}{\partial \theta^{\vec k_{\theta}}}
\frac{\partial}{\partial T_{{\rm e}_0}}
\Delta T^{\rho}_{\vec T,\vec \theta,(\kappa),{\rm v},{\rm e}}
\right\Vert
< C_{15,m}\mu^{\kappa-1}e^{-\delta T_{{\rm e}_0}},
\end{equation}
\begin{equation}\label{416}
\left\Vert
\nabla_{\rho}^n \frac{\partial^{\vert \vec k_{T}\vert}}{\partial T^{\vec k_{T}}}\frac{\partial^{\vert \vec k_{\theta}\vert}}{\partial \theta^{\vec k_{\theta}}}
\frac{\partial}{\partial T_{{\rm e}_0}}
\Delta \theta^{\rho}_{\vec T,\vec \theta,(\kappa),{\rm v},{\rm e}}\right\Vert
< C_{16,m}\mu^{\kappa-1}e^{-\delta T_{{\rm e}_0}}
\end{equation}
for $\vert{\vec k_{T}}\vert + \vert{\vec k_{\theta}}\vert + n < m -11$.
\par
Let ${\rm e}_{0} \in C^1_{\rm c}(\mathcal G_{\frak p})$. Then the same inequalities as above hold if we replace
$\frac{\partial}{\partial T_{{\rm e}_0}}$ by $\frac{\partial}{\partial \theta_{{\rm e}_0}}$.
\end{lem}
The proof is mostly the same as that of Proposition \ref{expesgen2Tdev}.
The difference is the following point only.
We remark that in (\ref{form0182vv233}),
(\ref{form0183v2v33}), (\ref{form0184233}),
(\ref{form0185vv233}) the norm is
$L^2_{m+1-\vert{\vec k_{T}}\vert - \vert{\vec k_{\theta}}\vert-n-1,\delta}$ norm.
On the other hand, in (\ref{form0182vv2}),
(\ref{form0183v2v}), (\ref{form01842}), (\ref{form0185vv2}), the norm was
$L^2_{m+1-\vert{\vec k_{T}}\vert - \vert{\vec k_{\theta}}\vert-1,\delta}$ norm.
The reason is as follows.
We remark that in our case
$$
T^{(\kappa)}_{\rm e} = T_{\rm e} - \frac{1}{10}\sum_{a=0}^{\kappa}\Delta T^{\rho}_{\vec T,\vec \theta,(a),{\rm v}_{\leftarrow}({\rm e}),{\rm e}}
+ \frac{1}{10}\sum_{a=0}^{\kappa}\Delta T^{\rho}_{\vec T,\vec \theta,(a),{\rm v}_{\rightarrow}({\rm e}),{\rm e}}
$$
is $\rho$ dependent. When we study $\rho$ derivative in the inductive steps,
we need to take $\rho$ derivative of
$$
\hat u^{\rho}_{{\rm v}',\vec T,\vec \theta,{(\kappa)}}(\tau'_{\rm e}-10T^{(\kappa)}_{\rm e},t'_{\rm e}+\theta^{(\kappa)}_{\rm e})
$$
etc..
Then there will be a term including $\tau''_{\rm e}$ or $t''_{\rm e}$ derivative of $\hat u^{\rho}_{{\rm v}',\vec T,\vec \theta,{(\kappa)}}$.
\par
Except this point  the proof of Lemma \ref{expesgen2Tdevss} is the same as the proof of Proposition \ref{expesgen2Tdev}
and so is omitted.
\begin{proof}[Proof of Proposition \ref{changeinfcoorprop}]
We note that (\ref{form0186vv2}) and (\ref{3383}) imply
\begin{equation}\label{rhoconvergesTTT}
\left\Vert \nabla_{\rho}^n \frac{\partial^{\vert \vec k_{T}\vert}}{\partial T^{\vec k_{T}}}\frac{\partial^{\vert \vec k_{\theta}\vert}}{\partial \theta^{\vec k_{\theta}}}
\frac{\partial}{\partial T_{{\rm e}_0}}
(\rho_{(\kappa)} - \rho) \right\Vert
< C_{17,m}\mu^{\kappa-1}e^{-\delta T_{\rm e}}
\end{equation}
and the same formula
with
$\frac{\partial}{\partial T_{{\rm e}_0}}$ replaced by $\frac{\partial}{\partial \theta_{{\rm e}_0}}$ if ${\rm e}_{0} \in C^1_{\rm c}(\mathcal G_{\frak p})$.
(\ref{form408}), (\ref{rhoconvergesTTT}), (\ref{415}) and (\ref{416}) imply (\ref{2144}).
\end{proof}
\begin{proof}[Proof of
Proposition \ref{reparaexpest}]
This is an immediate consequence of (\ref{3409}) and (\ref{form0184233}).
\end{proof}
\begin{proof}[Proof of Lemma \ref{changeinfcoorproppara}]
This is a parametrized version and the proof is the same as above.
\end{proof}

\section{Appendix: From $C^m$ structure to $C^{\infty}$ structure}
\label{toCinfty}

In this subsection we will prove that the Kuranishi structure
of $C^{m}$-class, which we obtained in Section \ref{generalcase},
is actually of $C^{\infty}$-class.
\par
We consider the embedding $\frak F^{(1)}$ (see the formula (\ref{evaluandTfac})) which we constructed in the proof of
Lemma \ref{2120lem}. Here we fix $m$.
\begin{lem}\label{lem3143}
The image of $\frak F^{(1)}$ is a $C^{\infty}$ submanifold.
\end{lem}
\begin{proof}
We first note several obvious facts.
Let $\frak M$ be a Banach manifold and $X \subset \frak M$ be a
subset. Then the statement that $X$ is a $C^{m'}$-submanifold of
finite dimension is well-defined. And the $C^{m'}$-structure of $X$
as a submanifold is unique if exists. Here $m'$ is one of $0,1,\dots,\infty$.
Moreover $X$ is a $C^{\infty}$-submanifold if and only if for each $p \in X$ and $m'$ there
exists a neighborhood  $U$ of $p$ such that $U \cap X$ is a submanifold of
$C^{m'}$-class.
\par
Now we prove the lemma.
Let $\frak q$ be in the image of $\frak F^{(1)}$ and take any $m'$.
Let $\frak w_{\frak p}$ be the stabilization data at $\frak p$
that we used to define $\frak F^{(1)}$.
We take the stabilization data $\frak w_{\frak q}$ on $\frak q$ that is induced by $\frak w_{\frak p}$.
We define ${\rm Glue}$ at $\frak q$ using the stabilization data $\frak w_{\frak q}$.
Then, as in the proof of Lemma \ref{coodinatedifferentp}, we obtain
\begin{equation}\label{evaluandTfacjj22}
\aligned
\frak F^{(2)} :
\hat V&(\frak q,\frak w_{\frak q};(\frak o',\mathcal T';\frak A))\\
\to
&\prod_{{\rm v}\in C^0(\mathcal G_{\frak q})} C^{m'}((K_{\rm v}^{+\vec R},K_{\rm v}^{+\vec R}
\cap \partial \Sigma_{\frak q,{\rm v}}),(X,L))\\
&\times \prod_{{\rm v}\in C^0(\mathcal G_{\frak p})}\frak V((\frak x_{\frak p} \cup \vec w_{\frak q)_{\rm v}})\times ((\vec{\mathcal T'},\infty] \times (\vec{\mathcal T}',\infty] \times \vec S^1).
\endaligned
\end{equation}
Let us denote the target of $\frak F^{(j)}$ by $\frak X(j)$.
The map $\frak F^{(2)}$ is a $C^{m'}$-embedding.
We define $\pi_{m,m'} : \frak X(2) \to \frak X(1)$ so that it is the identity map
for the second factor and the inclusion map
$$
C^{m'}((K_{\rm v}^{+\vec R},K_{\rm v}^{+\vec R}
\cap \partial \Sigma_{\frak q,{\rm v}}),(X,L))
\to
C^{m}((K_{\rm v}^{+\vec R},K_{\rm v}^{+\vec R}
\cap \partial \Sigma_{\frak q,{\rm v}}),(X,L))
$$
for the first factor. This map is of $C^{\infty}$ class.
We note that
$$
\pi_{m,m'} \circ \frak F^{(2)} = \frak F^{(1)} \circ \phi_{12},
$$
since we use the induced stabilization data for $\frak q$.
We already proved that $\varphi_{12}$ is a diffeomorphism of $C^m$-class
to an open subset.
Moreover $\frak F^{(2)}$ is an embedding of  $C^{m'}$-class.
Therefore a neighborhood of $\frak q$ of the image of $\frak F^{(1)}$
is a submanifold of $C^{m'}$-class.
The proof of Lemma \ref{lem3143} is complete.
\end{proof}
We define a $C^{\infty}$ structure of the Kuranishi neighborhood so that
$\frak F^{(1)}$ is a diffeomorphism to its image.
\begin{lem}
The coordinate change $\phi_{12}$ we defined is a diffeomorphism of $C^{\infty}$-class.
\end{lem}
\begin{proof}
We prove the case of $\phi_{12}$ in Lemma \ref{2120lem}.
We consider the following commutative diagram.
\begin{equation}\label{cd347333}
\begin{CD}
 \hat V(\frak p,\frak w^{(1)}_{\frak p};(\frak o^{(1)},\mathcal T^{(1)});\frak A)
_{\epsilon_{0},\vec{\mathcal T}^{(1)}} @ > {\frak F^{(1)}_{m'}} >>
\frak X^{(1)}_{m'} @>>{\pi_{m,m'}}> \frak X^{(1)}_{m}\\
@ AA{\subset}A @ AA\frak H_{12}A @ AA\frak H_{12}A\\
 \hat V(\frak p,\frak w^{(2)}_{\frak p};(\frak o^{(2)},\mathcal T^{(2)});\frak A)
_{\epsilon_{0},\vec{\mathcal T}^{(2)}}
@ > {\frak F^{(2)}_{2m'}}>>   \frak X^{(2)}_{2m'} @>>{\pi_{2m,2m'}}>
\frak X^{(2)}_{2m}
\end{CD}
\end{equation}
Here
\begin{equation}
\aligned
\frak X^{(2)}_{2m'} :
=&\prod_{{\rm v}\in C^0(\mathcal G_{\frak p})} C^{2m'}((K_{\rm v}^{+\vec R},K_{\rm v}^{+\vec R}\cap
\partial \Sigma_{\frak p,{\rm v}}),(X,L))\\
&\times \prod_{{\rm v}\in C^0(\mathcal G_{\frak p})}\frak V((\frak x_{\frak p} \cup \vec w_{\frak p})_{\rm v})\times ((\vec{\mathcal T}^{(2)},\infty] \times (\vec{\mathcal T}^{(2)},\infty] \times \vec S^1)
\endaligned
\end{equation}
is the space appearing in (\ref{evaluandTfac}), (\ref{evaluandTfac2}) and the map
$\frak F^{(2)}_{2m'}$ is defined as in (\ref{evaluandTfacjj}).
(We include $2m'$ in the notation to specify the function space we use.)
The space
$\frak X^{(1)}_{m'}$ and the map $\frak F^{(1)}_{m'}$ are
similarly defined.
The two maps $\frak H_{12}$ in the vertical arrow
are given by
$$
\frak H_{12}(u,(\rho,\vec T,\vec \theta))) = (u\circ \frak v_{(\rho,\vec T,\vec \theta)},\overline{\Phi}_{12}(\rho,\vec T,\vec \theta)).
$$
The maps in the horizontal lines are of $C^{\infty}$ class by definition.
The map $\frak H_{12}$ in the second vertical line is of $C^{m'}$ class by Sublemma \ref{sublem1}.
The map $\frak H_{12}$ in the third vertical line is one used in the proof of Lemma \ref{2120lem}.
Therefore $\phi_{12}$ is of $C^{m'}$-class at $\frak p$.
Note we can start at arbitrary point $\frak q$ in the image of $\frak F^{(2)}$
and prove that $\phi_{12}$ is of $C^{m'}$-class for any $m'$ at any point $\frak q$,
by using the proof of Lemma \ref{lem3143}.
This implies the lemma in the case of $\phi_{12}$ in Lemma \ref{2120lem}.
\par
In the other cases, the proof of the smoothness of the coordinate change is similar.
\end{proof}
We have thus  proved that the Kuranishi structure we obtained is of $C^{\infty}$-class.
\par\medskip
\section{Appendix: Proof of Lemma \ref{transbetweenEs}}
\label{Lemma49}
\begin{proof}[Proof of Lemma \ref{transbetweenEs}]
\begin{sublem}
There exists a finite dimensional smooth and compact family $\frak M$ of pairs $(\Sigma,u')$ such that
each element of $\mathcal M_{k+1,\ell}(\beta)$ appears as its member.
\end{sublem}
\begin{proof}
Run the gluing argument of Section \ref{glueing} at each point $\frak p \in \mathcal M_{k+1,\ell}(\beta)$
using an obstruction bundle data given at that point.
We then obtain a neighborhood of each $\frak p$ in a finite dimensional manifold.
We can take finitely many of them to cover  $\mathcal M_{k+1,\ell}(\beta)$ by compactness.
\end{proof}
We take a finite number of $\frak p_c$ so that (\ref{mpccoverM}) is satisfied.
For each $c$ and $N\in \Z_+$ we take $E_{c,N} \subset \bigoplus_{\rm v}\Gamma_0({\rm Int}\,K^{\rm obst}_{\rm v};
u_{\frak p_c}^*TX \otimes \Lambda^{01})$ that is isomorphic to the $N$ copies of $E_c$ as
a $\Gamma_{\frak p_c}$ vector space and $E_c \subset E_{c,N}$.
\par
We consider the space of $\Gamma_{\frak p_c}$-equivariant embeddings $\sigma_c : E_c \to E_{c,N}$ in the neighborhood of the original embedding.
Each $\sigma_c$ determines a perturbed $E_c$ which we write $E^{\sigma_c}_c$.
\par
The condition that $E^{\sigma_c}_c(\frak q) \cap E^{\sigma_{c'}}_{c'} (\frak q) \ne \{0\}$ for some $\frak q \in \frak M$ such that
$\frak q \cup \vec w'_{c}$ is
$\epsilon_{\frak p_c}$ close to $\frak p_c$ defines a subspace of the set of $(\sigma_c)_{c\in \frak C}$'s whose
codimension depends on the number of $c$'s, the dimension of $E_c$ and the dimension of $\frak M$ and $N$.
By taking $N$ huge, we may assume that such $(\sigma_c)_{c\in \frak C}$ is nowhere dense.
Namely the conclusion holds after perturbing $E_c$ by arbitrary small amount in  $E_{c,N}$.
\end{proof}

\par
\newpage
\part{$S^1$ equivariant Kuranishi structure and Floer's moduli space}
\label{S1equivariant}

In Part \ref{S1equivariant},
we explain Kuranishi structure on the space of
connecting orbits in the Floer theory for periodic Hamiltonian system
and use it to calculate the Floer homology of periodic Hamiltonian
system. This was done in \cite{FOn} but here we give more detailed proof than one in \cite{FOn}.
\par
Section \ref{abstractS1equiv} contains abstract theory of $S^1$-equivariant Kuranishi structure.
We define the notion of Kuranishi structure which admits  a locally free $S^1$ action
and its good coordinate system.
We explain how the construction in Part \ref{Part2} (that is basically the same as \cite{FOn} and \cite[Sectin A1]{fooo:book1})
can be modified so that all the constructions are $S^1$-equivariant.
\par
In Section \ref{Floersequation} we review the moduli space of
solutions of Floer's equation
(the Cauchy-Riemann equation perturbed by a Hamiltonian vector field).
\par
In Section \ref{construS1equiv}  we study the case of time independent Hamiltonian
and prove in detail that the Floer's  moduli space
has an $S^1$ equivariant Kuranishi structure in that case.
\par
In Section \ref{calcu}, we prove in detail that the
Floer homology of periodic Hamiltonian system is
isomorphic to the singular homology.
Namely it provides the detail of the proof of \cite[Theorem 20.5]{FOn2}.
\par\medskip
We also remark that, at the stage of the year 2012 (when this article is written), there are two proofs of isomorphism between
Floer homology of periodic Hamiltonian system and ordinary homology of a symplectic manifold $X$.
One is in \cite{FOn} and uses identification with Morse complex in the case Hamiltonian
is small and time independent.
This proof is the same as the one taken in this article.
(Namely the proof given in this article coincides with the one in
\cite{FOn} except some technical detail.)
The other uses Bott-Morse and de Rham theory and is in \cite[Section 26]{fooospectr}.
(Several other proofs are written in 1996 by Ruan \cite{Rua99}, Liu-Tian \cite{LiuTi98} also.)
This second proof has its origin in (the proof of) \cite[Theorem 1.2]{Fuk96II}.
\par
On the other hand, there is a third method using the Lagrangian Floer cohomology of the diagonal.
In this third method we do not need to study $S^1$ equivariant
Kuranishi structure at all. See Remark \ref{lagrangeproof}.

\section{Definition of $S^1$ equivariant Kuranishi structure, its good coordinate system and
perturbation}\label{abstractS1equiv}

We define an $S^1$ equivariant Kuranishi structure below.
A notion of $T^n$ equivariant Kuranishi structure (in a strong sense) is defined in
\cite[Definition B.4]{fooo:toric1}.
However that definition applies to the case when $T^n$ acts on the
target. The $S^1$ equivariant Kuranishi structure we use to study
time independent Hamiltonian is different therefrom since our $S^1$ action comes from the automorphism of the source.
Definition \ref{chertS1act} gives a definition in the current case.
Let $X$ be a Hausdorff metrizable space on which $S^1$ acts.
We assume that the isotropy group of every element is finite.

\begin{defn}\label{chertS1act}
Let
$(V_p, E_p, \Gamma_p, \psi_p, s_p)$ be a Kuranishi neighborhood of $p\in X$ as in
Definition \ref{Definition A1.1}, except replace the assumption $\gamma o_p = o_p$ for
$\gamma \in \Gamma_p$ by Condition (8) below.\footnote{$o_p \in V_p$ is a point such that
$\psi_p([o_p]) = p$.}
We define a {\it locally free $S^1$ action} on this chart as follows.
\begin{enumerate}
\item
There exists a group $G_p$ acting effectively on $V_p$ and $E_p$.
\item
$G_p \supset \Gamma_p$ and the  $\Gamma_p$  action extends to the $G_p$ action.
\item
The identity component $G_{p,0}$ of $G_p$ is isomorphic to $S^1$. We fix
an isomorphism $\frak h_p : S^1 \to G_{p,0}$.
\item
$G_p$ is generated by $\Gamma_p$ and $G_{p,0}$.
\item
$G_{p,0}$ commutes with  the action of $\Gamma_p$.
\item
The isotropy group at every point of $G_p$ action on $V_p$ is finite.
\item $s_p$ and $\psi_p$ are $G_p$ equivariant.
\item
The $G_{p,0}$ orbit of $o_p$ is invariant of $G_p$.
\end{enumerate}
\end{defn}
\begin{rem}
Note the Conditions (4), (5) imply that $G_p$ is isomorphic to the direct product
$\Gamma_p \times S^1$.
\end{rem}
The next example shows a  reason why we remove the assumption $\gamma o_p = o_p$.
\begin{exm}\label{zeifert}
We take $V_p = S^1 \times D^2$ and $\Gamma_p = \Z_2$ such that the
nontrivial element of $\Gamma_p$ acts by $(t,z) \mapsto (t+1/2,-z)$.
(Here $S^1 = \R/\Z$.) $G_{p,0} = S^1$ acts on $V_p$ by rotating the first factor $S^1$.
The action of $\Gamma_p$ is free. The quotient space $V_p/\Gamma_p$ is a manifold.
The induced $S^1$ action on $V_p/\Gamma_p$ is locally free but is not free.
The quotient space  $V_p/G_p$ is an orbifold $D^2/\Z_2$.
See Example \ref{zeifertgeo}.
\end{exm}
\begin{defn}\label{chengeS1act}
Let $(V_p, E_p, \Gamma_p, \psi_p, s_p)$ and
$(V_q, E_q, \Gamma_q, \psi_q, s_q)$ be Kuranishi neighborhoods of
$p \in X$ and $q \in \psi_p(s_p^{-1}(0)/\Gamma_p)$, respectively.
We assume that they carry locally free $S^1$ actions. ($G_p$ and $G_q$.)
Let a triple
$(\hat\phi_{pq},\phi_{pq},h_{pq})$ be a coordinate change in the sense of Definition \ref{Definition A1.3}.
We say it is {\it $S^1$ equivariant} if the following holds.
\begin{enumerate}
\item
$h_{pq}$ extends to a group homomorphism $G_q \to G_p$, which we denote by $\frak h_{pq}$.
\item
$V_{pq}$ is $G_q$ invariant.
\item
$\phi_{pq} : V_{pq} \to V_p$ is $\frak h_{pq}$-equivariant.
\item
$\hat\phi_{pq}$ is $\frak h_{pq}$-equivariant.
\item
$
\frak h_{pq}\circ \frak h_{q} = \frak h_{p}.
$
\end{enumerate}
\end{defn}
\begin{defn}\label{Definition A1.5S}
Let  $(V_p, E_p, \Gamma_p, \psi_p, s_p)$
and
$(\hat\phi_{pq},\phi_{pq},h_{pq})$
define a Kuranishi structure in the sense of Definition \ref{Definition A1.5} on $X$.
A {\it locally free $S^1$ action} on $X$ is assigned by a locally free $S^1$ action
in the sense of Definition \ref{chertS1act} on each chart $(V_p, E_p, \Gamma_p, \psi_p, s_p)$
so that the coordinate change is $S^1$ equivariant.
\end{defn}
\begin{rem}
Let $r \in \psi_q((V_{pq}\cap s_q^{-1}(0))/\Gamma_q)$,
$q \in \psi_p(s_p^{-1}(0)/\Gamma_p)$. There exists $\gamma_{pqr}^{\alpha} \in \Gamma_p$
for each connected component
$(\phi_{qr}^{-1}(V_{pq}) \cap V_{qr} \cap V_{pr})_\alpha$ of $\phi_{qr}^{-1}(V_{pq}) \cap V_{qr} \cap V_{pr}$ by
Definition \ref{Definition A1.5} (2). We automatically have
\begin{equation}\label{addiionalcond}
\frak h_{pq} \circ \frak h_{qr} = \gamma_{pqr}^{\alpha}
\cdot \frak h_{pr} \cdot (\gamma_{pqr}^{\alpha})^{-1},
\end{equation}
because $S^1$ lies in the center and this formula is already assumed for $\Gamma_r$.
\end{rem}

\begin{lemdef}\label{lemdef1}
If $X$ has a Kuranishi structure with a locally free $S^1$ action
then $X/S^1$ has an induced Kuranishi structure.
\end{lemdef}
\begin{proof}
Let $p\in X$. We take $o_p \in V_p$ and
choose a local transversal $\overline V_p$ to the $S^1$ orbit
$G_{p,0}o_p$.
We put
$$
\Gamma^+_p = \{ \gamma \in G_p \mid \gamma o_p = o_p\}.
$$
By Definition \ref{chertS1act} (6),
$\Gamma^+_p$ is a finite group.
We may choose $\overline V_p$ so that it is invariant under $\Gamma^+_p$.
We restrict $E_p$ to $\overline V_p$ to obtain $\overline E_p$.
The Kuranishi map $s_p$ induces $\overline s_p$.
\par
We may shrink our Kuranishi neighborhood and may assume
that
\begin{equation}\label{equiv1}
V_p = G_{p,0}\cdot\overline V_p
\end{equation}
for all $p$.
\par
For $x \in \overline V_p$, satisfying $s_p(x) = 0$,
we define  $\overline\psi_p(x)$ to be
the
equivalence class of $\psi_p(x)$ in $X/S^1$.
It is easy to see that
$(\overline V_p,\Gamma^+_p,\overline E_p,\overline s_p,\overline\psi_p)$ is a Kuranishi chart of $X/S^1$ at
$[p]$.
\par
Let $[q] \in \overline\psi_p(\tilde q)$, where $\tilde q \in \overline V_p$.
Choose $q \in X$ such that $q = \psi_p(\tilde q)$.
We have a coordinate transformation $(V_{pq},\phi_{pq},\hat\phi_{pq})$
and a group homomorphism $\frak h_{pq} : G_q\to G_p$
such that
$$
 \tilde q = g_{pq}(o_q) \cdot \phi_{pq}(o_q)
$$
holds for some $g_{pq}(o_q) \in G_{p,0}$. Moreover there exists  a smooth map
$$
g_{pq} : \overline V_{pq} \to G_{p,0}
$$
such that it coincides with $g_{pq}(o_q)$ at $o_q$ and
$$
g_{pq}(x) \cdot \phi_{pq}(x) \in \overline V_{p}.
$$
Here $\overline V_{pq}$ is a neighborhood of $o_q$ in $\overline V_q$.
We define
$$
\overline{\phi}_{pq}(x) =  g_{pq}(x) \cdot \phi_{pq}(x).
$$
We shrink $V_{pq}$ and may assume
\begin{equation}\label{equiV2}
V_{pq} = G_{p,0}\cdot\overline V_{pq}.
\end{equation}
\par
By definition
$$
\Gamma^+_{q} = \{ \gamma \in G_q \mid \gamma o_q = o_q\}.
$$
Using the fact that
$$
\{\gamma \in \Gamma_{p}\mid \gamma\cdot \phi_{pq}(o_q) = \phi_{pq}(o_q) \} = h_{pq}(\Gamma_{q})
$$
and $G_{p,0}$ is contained in the center, we find that
$$
\{ \gamma \in G_p \mid \gamma \overline{\phi}_{pq}(o_q) = \overline{\phi}_{pq}(o_q)\} = \frak h_{pq}(\Gamma^+_{q})
\subset \Gamma_p^+.
$$
We denote by $\overline h_{pq}$ the restriction of $\frak h_{pq}$ to $\Gamma^+_{q}$.
It is easy to see that
$\overline{\phi}_{pq}$ is $\overline h_{pq}$ equivariant.
We can lift $\overline{\phi}_{pq}$  to $\hat{\overline{\phi}}_{pq}$ using ${\phi}_{pq}$ and $G_p$ action on $E_p$.
\par
We have thus constructed a coordinate change of our Kuranishi structure on $X/S^1$.
It is straightforward to check the compatibility among the coordinate changes.
\end{proof}
We next define a good coordinate system.
We note that in Part \ref{Part2} we defined a chart of good coordinate system
as an orbifold that is not necessarily a global quotient.
So we define a notion of locally free $S^1$ action on orbifold.
\begin{defn}
Let $U$ be an orbifold on which $S^1$ acts effectively as a topological group.
We assume that the isotropy group of this $S^1$ action is always finite.
We say that the action is a {\it smooth action on orbifold}
if the following holds for each $p \in U$.
\par
There exists an $S^1$ equivariant neighborhood $U_p$ of $p$ in $U$ and
$V_p$ a manifold on which $G_p$ acts. $(V_p,\Gamma_p,\psi_p)$ is a
chart of $U$ as an orbifold.
The conditions (1)-(6)  in Definition \ref{chertS1act} hold and $\psi_p$ is $G_p$ equivariant.
Moreover the
$S^1$ action on $V_p/S^1 \subset U$ induced by $\frak h_p : S^1 \to G_{p,0}$ coincides with the
given $S^1$ action.
\end{defn}

Let $S^1$ act effectively on $X$ and assume that its isotropy group is finite.
\begin{defn}\label{goodcoordinatesystemS}
Suppose $X$ has a locally free $S^1$ equivariant Kuranishi structure.
An \emph{$S^1$ equivariant good coordinate system} on it is $(U_{\frak p}, E_{\frak p}, \psi_{\frak p}, s_{\frak p})$,
$(U_{\frak p \frak q},\hat\phi_{\frak p \frak q},\phi_{\frak  p\frak q})$
as in Definition \ref{goodcoordinatesystem}. We require furthermore the following in addition.
\begin{enumerate}
\item
There exists a smooth $S^1$ action on $U_{\frak p}$ and $E_{\frak p}$.
\item $\psi_{\frak p}$, $s_{\frak p}$ are $S^1$ equivariant.
\item $U_{\frak p \frak q}$ is $S^1$ invariant and $\hat\phi_{\frak p \frak q},\phi_{\frak  p\frak q}$ are $S^1$ equivariant.
\end{enumerate}
\end{defn}
Note the notion of $S^1$-equivariance of maps or subsets are defined set theoretically.
\begin{lem}
If $(U_{\frak p}, E_{\frak p}, \psi_{\frak p}, s_{\frak p})$,
$(U_{\frak p \frak q},\hat\phi_{\frak p \frak q},\phi_{\frak  p\frak q})$
is an $S^1$ equivariant good coordinate system
then it induces
a good coordinate system of $X/S^1$, that is
$(\overline U_{\frak p}, \overline E_{\frak p}, \overline \psi_{\frak p}, \overline s_{\frak p})$,
$(\overline U_{\frak p \frak q},\hat{\overline \phi}_{\frak p \frak q},\overline \phi_{\frak  p\frak q})$, where $\overline U_{\frak p} = U_{\frak p}/S^1$
etc..
\end{lem}
\begin{proof}
Apply the construction of Lemma-Definition \ref{lemdef1} locally.
\end{proof}
\begin{prop}\label{S1goodcoordinate}
For any locally free $S^1$ equivariant Kuranishi structure  we can find an $S^1$ equivariant good coordinate system.
\end{prop}
\begin{proof}
The proof uses the construction of good coordinate system in Section
\ref{sec:existenceofGCS}.
We defined and used the notion of pure and mixed orbifold neighborhood there.
We constructed them for Kuranishi structure.
We will use pure and mixed orbifold neighborhood of the
Kuranishi structure on $X/S^1$ and extend them to ones on $X$.
The detail follows.
\par
We stratify
$\overline X = X/S^1 = \bigcup_{\frak d}\overline X(\frak d)$ where $[p]\in \overline X(\frak d)$
if $\dim \overline U_{[p]} = \frak d$.
So $X(\frak d+1)/S^1 = \overline X(\frak d)$.
Let $\overline{\mathcal K}_*$ be a compact subset of $\overline X(\frak d)$.
Let ${\mathcal K}_*$ be an $S^1$-invariant compact subset of $X(\frak d)$ such that
$\overline{\mathcal K}_* = {\mathcal K}_*/S^1$.
In Proposition \ref{purecover} we constructed a pure orbifold neighborhood $\overline U_*$ of ${\mathcal K}_*/S^1$.
\begin{lem}\label{S1tuki44}
There exists a pure orbifold neighborhood $U_*$ of ${\mathcal K}_*$ on which $S^1$ acts and
$U_*/S^1 = \overline U_*$.
\end{lem}
\begin{rem}
This lemma is somewhat loosely stated, since we did not define the notion of
$S^1$ action on pure orbifold neighborhood. The definition is: $U_{*}$ has a locally free effective smooth
$S^1$ action and all the structure maps commute with the $S^1$ action.
\end{rem}
\begin{proof}
We can prove this lemma by examining the proof of Proposition \ref{purecover}. Namely
$\overline U_*$ is obtained by gluing various Kuranishi charts and restricting it to suitable open subsets.
We take the inverse image of $U_{\frak p} \to \overline{U}_{\frak p}$ of those charts.
We can then glue and restrict them in the same way to obtain $U_*$.
We omit the detail.
\end{proof}
In Section \ref{sec:existenceofGCS} we then proceed to define a mixed orbifold neighborhood of $\overline X(D)$ for an ideal $D \subset \frak D$.
For an ideal $D \subset \frak D$ we put $D^{+1} = \{\frak d+1 \mid \frak d \in D\}$.
\begin{lem}\label{S1tuki442}
We assume that $\{\overline{\mathcal U}_{\frak d}\}$ together with other data provide the mixed orbifold neighborhood of $\overline X(D)$
obtained in Proposition \ref{existmixed}.
\par
Then we can take an $S^1$ equivariant mixed orbifold neighborhood $\{{\mathcal U}_{\frak d+1}\}$ (plus other data) on
$X(D^{+1})$ such that $\overline{\mathcal U}_{\frak d} = {\mathcal U}_{\frak d+1}/S^1$.
\end{lem}
\begin{proof}
This is proved again by examining the proof of Proposition \ref{existmixed}
and checking that the gluing process there can be lifted.
This is actually fairly obvious.
\end{proof}
We note that the chart of good coordinate system on $\overline X$ constructed in Part \ref{Part2}
is $\overline{\mathcal U}_{\frak d}$ and other data of good coordinate system
is obtained by the structure maps etc. of mixed orbifold neighborhood.
Therefore $ {\mathcal U}_{\frak d+1}$ becomes the required $S^1$ equivariant
good coordinate system of $X$.
The proof of Proposition \ref{S1goodcoordinate} is complete.
\end{proof}
\begin{lem}\label{lem214}
If the dimension of $X/S^1$ in the sense of Kuranishi structure is $-1$ then there exists an
$S^1$ equivariant multisection on the good coordinate system $\mathcal U_{\frak p}$
of $X$ whose zero set is empty.
\end{lem}
\begin{proof}
It suffices to define an appropriate notion of pull back of the multisection of
$\overline{\mathcal U}_{\frak p}$ to ones of $\mathcal U_{\frak p}$.
This is routine.
\end{proof}

\section{Floer's equation and its moduli space}
\label{Floersequation}

In this section we concern with the moduli space of solutions of
Floer's perturbed Cauchy-Riemann equation.
Such a moduli space appears in the proof of Arnold's conjecture of various kinds.
In the next section, we prove existence of $S^1$ equivariant Kuranishi structure
of such a moduli space in the case when our Morse function is time independent.

Let $H : X \times S^1 \to \R$ be a smooth function on a symplectic manifold $X$.
We put $H_t(x) = H(x,t)$ where $t \in S^1$ and $x \in X$.
The function $H_t$ generates the Hamiltonian vector field $\frak X_{H_t}$ by
$$
i_{\frak X_{H_t}}\omega = d{H_t}.
$$
We denote it by $\frak P(H)$ the set of the all 1-periodic orbits of the time dependent vector field $\frak X_{H_t}$.
We put
$$
\tilde{\frak P}(H) = \{(\gamma,w) \mid \gamma \in {\frak P}(H),\,\, u : D^2 \to X,\,\, u(e^{2\pi it}) = \gamma(t)\}/\sim,
$$
where $(\gamma,w) \sim (\gamma',w')$ if and only if $\gamma = \gamma'$ and
$$
\omega([w] - [w']) = 0,\quad
c_1([w] - [w']) = 0.
$$
Here $\omega$ is the symplectic form and $c_1$ is the first Chern class of $X$.
\begin{assump}\label{nondeg1}
All the  1-periodic orbits of the time dependent vector field $\frak X_{H_t}$ are non-degenerate.
\end{assump}
Following \cite{Flo89I}, we consider the maps $h : \R \times S^1 \to X$ that satisfy
\begin{equation}\label{Fleq}
\frac{\partial h}{\partial \tau} +  J \left( \frac{\partial h}{\partial t} - \frak X_{H_t} \right) = 0.
\end{equation}
Here $\tau$ and $t$ are the coordinates of $\R$ and $S^1 = \R/\Z$, respectively.
For $\tilde{\gamma}^{\pm} = ({\gamma}^{\pm},w^{\pm}) \in \tilde{\frak P}(H)$
we consider the boundary condition
\begin{equation}\label{bdlatinf}
\lim_{\tau \to \pm \infty}h(\tau,t) = \gamma^{\pm}(t).
\end{equation}
The following result due to Floer \cite{Flo89I} is by now well established.
\begin{prop}
We assume Assumption \ref{nondeg1}. Then for any solution $h$ of (\ref{Fleq})
with
$$
\int_{\R \times S^1}
\left\Vert\frac{\partial h}{\partial \tau}\right\Vert^2 d\tau dt < \infty
$$
there exists ${\gamma}^{\pm} \in {\frak P}(H)$ such that (\ref{bdlatinf})
is satisfied.
\end{prop}
Let $\tilde{\gamma}^{\pm} = ({\gamma}^{\pm},w^{\pm}) \in \tilde{\frak P}(H)$.
\begin{defn}\label{tildeMreg}
We denote by $\widetilde{\mathcal M}^{\text{\rm reg}}(X,H;\tilde{\gamma}^{-},\tilde{\gamma}^{+})$
the set of all maps $h : \R \times S^1 \to X$ that satisfy (\ref{Fleq}), (\ref{bdlatinf}) and
$$
w^- \# h \sim w^+.
$$
Here $\#$ is an obvious concatenation.
\par
The translation along $\tau \in \R$ defines an $\R$ action on $\widetilde{\mathcal M}^{\text{\rm reg}}(X,H;\tilde{\gamma}^{-},\tilde{\gamma}^{+})$.
This $\R$ action is free unless $\tilde{\gamma}^{-} = \tilde{\gamma}^{+}$.
We denote by ${\mathcal M}^{\text{\rm reg}}(X,H;\tilde{\gamma}^{-},\tilde{\gamma}^{+})$
the quotient space of this action.
\end{defn}
\begin{thm}{\rm (\cite[Theorem 19.14]{FOn})}\label{connectingcompactka}
We assume Assumption \ref{nondeg1}.
\begin{enumerate}
\item
The space
${\mathcal M}^{\text{\rm reg}}(X,H;\tilde{\gamma}^{-},\tilde{\gamma}^{+})$ has a compactification
${\mathcal M}(X,H;\tilde{\gamma}^{-},\tilde{\gamma}^{+})$.
\item
The compact space ${\mathcal M}(X,H;\tilde{\gamma}^{-},\tilde{\gamma}^{+})$ has an oriented Kuranishi structure with
corners.
\item
The codimension $k$ corner of ${\mathcal M}(X,H;\tilde{\gamma}^{-},\tilde{\gamma}^{+})$ is
identified with the union of
$$
\prod_{i=0}^{k} {\mathcal M}(X,H;\tilde{\gamma}_{i},\tilde{\gamma}_{i+1})
$$
over the $k+1$-tuples $(\tilde{\gamma}_{0},\dots,\tilde{\gamma}_{k+1})$
such that $\tilde{\gamma}_{0} = \tilde{\gamma}^{-}$, $\tilde{\gamma}_{k+1}=\tilde{\gamma}^{+}$
and $\tilde{\gamma}_{i} \in \tilde{\frak P}(H)$.
\end{enumerate}
\end{thm}
The proof is in Section \ref{calcu}.
\par
The main purpose of this section is to explain the proof of the next result.
We consider the case when $H$ is time independent.
In this case, Assumption \ref{nondeg1} implies that $H : X \to \R$ is a Morse function.
We also assume the  following:
\begin{assump}\label{nondeg2}
\begin{enumerate}
\item
The gradient vector field of $H$ satisfies the Morse-Smale condition.
\item
Any 1-periodic orbit of $\frak X_H$ is a constant loop.
(Namely it corresponds to a critical point of $H$.)
\end{enumerate}
\end{assump}
Condition (1) is satisfied for generic $H$. We can replace $H$ by $\epsilon H$
for small $\epsilon$ so that (2) is also satisfied.
\par
By assumption, elements of ${\frak P}(H)$ are constant loops. We write $\frak x \in X$ to denote its element.
We put
$$
\Pi = \frac{\text{\rm Im}~ (\pi_2(X) \to H_2(X;\Z))}{\text{Ker}(c_1) \cap \text{Ker}(\omega) \cap \text{\rm Im}~(\pi_2(X) \to H_2(X;\Z))}.
$$
Here we regard $c_1 : H_2(X;\Z) \to \Z$, $\omega : H_2(X;\Z) \to \R$.
An element of $\tilde{\frak P}(H)$ is regarded as a pair $(\frak z,\alpha)$, where $\frak z$ is a critical point of $H$ and
$\alpha\in \Pi$.
\par
We put
$$
{\mathcal M}^{\text{\rm reg}}(X,H;\frak z^{-},\frak z^{+};\alpha) = {\mathcal M}^{\text{\rm reg}}(X,H;(\frak z^-,\alpha^{-}),(\frak z^+,\alpha^{-}+\alpha)).
$$
It is easy to see that the right hand is independent of $\alpha^{-} \in \Pi$.
\par
Let ${\mathcal M}(X,H;\frak z^-,\frak z^+;\alpha)$ be its compactification as in Theorem \ref{connectingcompactka}.
\par
Let ${\mathcal M}^{\text{\rm reg}}(X,H;\frak z^{-},\frak z^{+};\alpha)^{S^1}$ be the fixed point set of the $S^1$ action
obtained by $t_0h(\tau,t) = h(\tau,t+t_0)$.
It is easy to see that this set is empty unless $\alpha = 0$ and in the case $\alpha = 0$ the fixed point set
${\mathcal M}^{\text{\rm reg}}(X,H;\frak z^{-},\frak z^{+};0)^{S^1}$ can be identified with
the set of gradient lines of $H$ joining $\frak z^{-}$ to $\frak z^{+}$.
This identification can be extended to their compactifications.
\begin{assump}\label{nondeg3}
\begin{enumerate}
\item
${\mathcal M}(X,H;\frak z^{-},\frak z^{+};0)^{S^1}$ is an
open subset of ${\mathcal M}(X,H;\frak z^{-},\frak z^{+};0)$.
Namely any solution of (\ref{Fleq}) which is sufficiently close to
an $S^1$ equivariant solution is $S^1$ equivariant.
\item
The moduli space ${\mathcal M}(X,H;\frak z^{-},\frak z^{+};0)$ is  Fredholm regular
at each point of ${\mathcal M}(X,H;\frak z^{-},\frak z^{+};0)^{S^1}$.

\end{enumerate}
\end{assump}
\begin{lem}\label{identifiedwithgrad}
Assumption \ref{nondeg3} is satisfied if we replace $H$ by $\epsilon H$ for a sufficiently small $\epsilon$.
\end{lem}
\begin{proof}
(2) is proved in  \cite[page 1038]{FOn}.
More precisely, it is proved there that for sufficiently small $\epsilon$ the following holds.
Let $\ell$ be a gradient line joining $\frak z^{-}$ to  $\frak z^{+}$ and let $h_{\ell}$ be the
corresponding element of ${\mathcal M}(X,\epsilon H;\frak z^{-},\frak z^{+};0)^{S^1}$.
We consider the deformation complexes of the gradient line equation at $\ell$ and
of the equation (\ref{Fleq}) at $h_{\ell}$.
The kernel and the cokernel of the former are contained in the
kernel and the cokernel of the later, respectively.
It is proved in \cite[page 1038]{FOn} that they  actually coincide each other if $\epsilon$ is
sufficiently small.
\par
Since $H$ satisfies the Morse-Smale condition, the element $\ell$ is Fredholm regular in the
moduli space of gradient lines. Therefore by the above mentioned result
the moduli space ${\mathcal M}^{\text{\rm reg}}(X,\epsilon H;\frak z^{-},\frak z^{+};0)$
is Fredholm regular at
$h_{\ell}$.
This implies (2).
(1) is a consequence of the same result and the implicit function theorem.
(We note that we can prove the same result at the point
${\mathcal M}(X,\epsilon H;\frak z^{-},\frak z^{+};\alpha)^{S^1}
\setminus {\mathcal M}^{\text{\rm reg}}(X,\epsilon H;\frak z^{-},\frak z^{+};\alpha)^{S^1}$
in the same way.)
\end{proof}
Thus, replacing $H$ by $\epsilon H$ if necessary, we may assume that $H$ satisfies Assumption \ref{nondeg3}.
We put
\begin{equation}\label{M0000}
{\mathcal M}_0(X,H;{\frak x}^{-},\tilde{\frak x}^{+};0) =
{\mathcal M}(X,H;{\frak x}^{-},\tilde{\frak x}^{+};0)
\setminus
{\mathcal M}(X,H;{\frak x}^{-},\tilde{\frak x}^{+};0)^{S^1}.
\end{equation}
Lemma \ref{identifiedwithgrad} implies that
${\mathcal M}_0(X,H;{\frak x}^{-},\tilde{\frak x}^{+};0)$
is open and closed in
${\mathcal M}(X,H;{\frak x}^{-},\tilde{\frak x}^{+};0)$.

\begin{thm}{\rm (\cite[page 1036]{FOn})}\label{existKuran}
If we assume Assumptions \ref{nondeg1}, \ref{nondeg2}  and \ref{nondeg3},
then the following holds.
\begin{enumerate}
\item
In case $\alpha \ne 0$ the Kuranishi structure on
${\mathcal M}(X,H;{\frak x}^{-},\tilde{\frak x}^{+};\alpha)$
can be taken to be $S^1$ equivariant.
\item
In case $\alpha = 0$ the same conclusion holds for  ${\mathcal M}_0(X,H;{\frak x}^{-},\tilde{\frak x}^{+};0)$.
\end{enumerate}
\end{thm}

\section{$S^1$ equivariant Kuranishi structure for the Floer homology
of time independent Hamiltonian}
\label{construS1equiv}
In this section we prove Theorem \ref{existKuran} in detail.
We begin with describing the compactification ${\mathcal M}(X,H;{\frak z}^{-},{\frak z}^{+};\alpha)$
of the moduli space ${\mathcal M}^{\rm reg}(X,H;\frak z^-,\frak z^+;\alpha)$.
Here we include the case $\alpha = 0$ and the $S^1$ fixed point since it will appear in the fiber product factor of the compactification.
\par
We consider $(\Sigma,z_-,z_+)$, a genus zero semistable curve with two marked points.
\begin{defn}\label{defn41}
Let $\Sigma_0$ be the union of the irreducible components of $\Sigma$ such that
\begin{enumerate}
\item $z_-,z_+ \in \Sigma_0$.
\item $\Sigma_0$  is connected.
\item $\Sigma_0$ is smallest among those satisfying (1),(2) above.
\end{enumerate}
We call $\Sigma_0$ the {\it mainstream} of $(\Sigma,z_-,z_+)$, or simply, of $\Sigma$.
An irreducible component of $\Sigma$ that is not contained in $\Sigma_0$ is called a {\it bubble component}.
\par
Let $\Sigma_a \subset \Sigma$ be an irreducible component of the mainstream.
If $z_- \notin \Sigma_a$ then there exists a unique singular point $z_{a,-}$ of $\Sigma$ contained in $\Sigma_a$
such that
\begin{enumerate}
\item
$z_-$ and $\Sigma_a \setminus\{z_{a,-}\}$ belong to the different connected components of
$\Sigma \setminus\{z_{a,-}\}$.
\item
$z_+$ and $\Sigma_a \setminus\{z_{a,-}\}$ belong to the same connected components of
$\Sigma \setminus\{z_{a,-}\}$.
\end{enumerate}
In case $z_- \in \Sigma_a$ we set $z_- = z_{a,-}$.
\par
We define $z_{a,+}$ in the same way.
\par
A {\it parametrization of the mainstream} of $(\Sigma,z_-,z_+)$ is
$\varphi = \{\varphi_a\}$, where $\varphi_a : \R \times S^1 \to \Sigma_a$
for each  irreducible component $\Sigma_a$ of the mainstream such that:
\begin{enumerate}
\item
$\varphi_a$ is a biholomorphic map $\varphi_a : \R \times S^1 \cong \Sigma_a \setminus \{z_{a,-},z_{a,+}\}$.
\item
$\lim_{\tau \to \pm \infty}\varphi_a(\tau,t) = z_{a,\pm}$.
\end{enumerate}
\end{defn}
\begin{defn}\label{defn210}
We denote by $\widehat{\mathcal M}(X,H;\frak z^-,\frak z^+,\alpha)$
the set of triples $((\Sigma,z_-,z_+),u,\varphi)$ satisfying
the following conditions:
\begin{enumerate}
\item
$(\Sigma,z_-,z_+)$ is a genus zero semistable curve with two marked points.
\item
$\varphi$ is a parametrization of the mainstream.
\item
$u : \Sigma \to X$ a continuous map from $\Sigma$ to $X$.
\item
If $\Sigma_a$ is an irreducible component of the mainstream and
$\varphi_a : \R \times S^1 \to \Sigma_a$ is as above then
the composition $h_a = u \circ \varphi_a$ satisfies the equation (\ref{Fleq}).
\item
If $\Sigma_a$  is a bubble component then $u$ is pseudo-holomorphic on it.
\item
$u(z_-) = \frak z^-$, $u(z_+) = \frak z^+$.
\item
$[u_*[\Sigma]] = \alpha$. Here $\alpha \in \Pi$.
\end{enumerate}
\end{defn}
\begin{defn}\label{3equivrel}
On the set $\widehat{\mathcal M}(X,H;\frak z^-,\frak z^+,\alpha)$  we define three
equivalence relations $\sim_1$, $\sim_2$, $\sim_3$ as follows.
\par
$((\Sigma,z_-,z_+),u,\varphi) \sim_1 ((\Sigma',z'_-,z'_+),u',\varphi')$ if and only if there exists
a biholomorphic map $v : \Sigma \to \Sigma'$ with the following properties:
\begin{enumerate}
\item
$u' = u\circ v$.
\item
$v(z_{-}) = z'_{-}$ and $v(z_{+}) = z'_{+}$.
In particular $v$ sends the mainstream of $\Sigma$ to the mainstream of $\Sigma'$.
\item
If $\Sigma_a$ is an irreducible component of the mainstream of $\Sigma$
and $v(\Sigma_a) = \Sigma'_a$, then we have
\begin{equation}\label{preserveparame}
v \circ \varphi_a = \varphi'_a.
\end{equation}
\end{enumerate}
\par
The equivalence relation $\sim_2$ is defined  replacing (\ref{preserveparame}) by
existence of $\tau_a$ such that
\begin{equation}\label{preserveparame2}
(v \circ \varphi_a)(\tau,t) = \varphi'_a(\tau+\tau_a,t).
\end{equation}
\par
The equivalence relation $\sim_3$ is defined by requiring only (1), (2) above. (Namely by removing
condition (3).)
\end{defn}
\begin{rem}
After taking the $\sim_3$ equivalence class, the data $\varphi$ does not remain.
Namely  $((\Sigma,z_-,z_+),u,\varphi) \sim_3 ((\Sigma,z_-,z_+),u,\varphi')$ for any $\varphi, \varphi'$.
\end{rem}
\begin{defn}\label{defn45}
We put
$$\aligned
\widetilde{\mathcal M}(X,H;\frak z^-,\frak z^+,\alpha) &=
\widehat{\mathcal M}(X,H;\frak z^-,\frak z^+,\alpha)/\sim_1, \\
{\mathcal M}(X,H;\frak z^-,\frak z^+,\alpha) &=
\widehat{\mathcal M}(X,H;\frak z^-,\frak z^+,\alpha)/\sim_2, \\
\overline{\overline{\mathcal M}}(X,H;\frak z^-,\frak z^+,\alpha) &=
\widehat{\mathcal M}(X,H;\frak z^-,\frak z^+,\alpha)/\sim_3.
\endaligned$$
We use $\frak p$ etc. to denote an element of ${\mathcal M}(X,H;\frak z^-,\frak z^+,\alpha)$ and
denote by $[\frak p]$ its equivalence class in
$\overline{\overline{\mathcal M}}(X,H;\frak z^-,\frak z^+,\alpha)$.
\par
Let $((\Sigma,z_-,z_+),u,\varphi)$ be an element of $\widehat{\mathcal M}(X,H;\frak z^-,\frak z^+,\alpha)$.
Suppose the mainstream of $\Sigma$ has $k$ irreducible components.
We add the bubble tree to the irreducible component of the mainstream where
it is rooted. We thus have obtained a decomposition
\begin{equation}\label{stedec}
\Sigma = \sum_{i=1}^k \Sigma_{i}.
\end{equation}
Here $z_- \in \Sigma_1$, $z_+ \in \Sigma_k$ and $\# (\Sigma_i \cap \Sigma_{i+1}) = 1$.
We call each summand of (\ref{stedec}) a {\it mainstream component}.
We put $z_{i+1} = \Sigma_i \cap \Sigma_{i+1}$ and call it $(i+1)$-th {\it transit point}.
We put $\frak z_i = u(z_i)$ and call it a {\it transit image}.
\par
Let $\frak p = ((\Sigma,z_-,z_+),u,\varphi) \in \widehat{\mathcal M}(X,H;\frak z^-,\frak z^+,\alpha)$
and let $\Sigma_i$ be one of its mainstream component.
We restrict $u$, $\varphi$ to $\Sigma_i$ and obtain $\frak p_i$.
We say that $\Sigma_i$ is a {\it gradient line component}
if $\frak p_i$ is a fixed point of the $S^1$ action.
\par
In the definition of  $\widehat{\mathcal M}(X,H;\frak z^-,\frak z^+,\alpha)$
we forget the map $u$ but remember only the homology class of $u\vert_{\Sigma_v}$
of each irreducible component and the images $u(z_i)$ of the transit points (that are critical
points of $H$). We then obtain a decorated moduli space of domain curves denoted by
$\widehat{\mathcal M}(\frak z^-,\frak z^+,\alpha)$.
We define equivalence relations $\sim_j$ ($j=1,2,3$) on it in the same way and
obtain
$\widetilde{\mathcal M}(\frak z^-,\frak z^+,\alpha)$,
${\mathcal M}(\frak z^-,\frak z^+,\alpha)$, and
$\overline{\overline{\mathcal M}}(\frak z^-,\frak z^+,\alpha)$.
For each element $\frak p$ of $\widehat{\mathcal M}(X,H;\frak z^-,\frak z^+,\alpha)$ etc.,
we denote by $\frak x_{\frak p}$ the element of $\widehat{\mathcal M}(\frak z^-,\frak z^+,\alpha)$ etc.
obtained by forgetting $u$ as above.
\par
For each $\frak p \in \widehat{\mathcal M}(X,H;\frak z^-,\frak z^+,\alpha)$ etc. or
$\frak x \in \widehat{\mathcal M}(\frak z^-,\frak z^+,\alpha)$ etc.,
we define a graph $\mathcal G_{\frak p}$ or $\mathcal G_{\frak x}$ in the same way
as in Section \ref{Graph}.
We include the data of the homology class of each component
and the images of the transit points in  $\mathcal G_{\frak p}$ (resp. $\mathcal G_{\frak x}$).
We  call $\mathcal G_{\frak p}$ (resp. $\mathcal G_{\frak x}$) the combinatorial type of $\frak p$ (resp. $\frak x$).
We denote by $\widehat{\mathcal M}(X,H;\frak z^-,\frak z^+,\alpha;\mathcal G)$ etc.
or $\widehat{\mathcal M}(\frak z^-,\frak z^+,\alpha;\mathcal G)$ etc. the
subset of the objects with combinatorial type $\mathcal G$.
\end{defn}
We consider the subset $\widehat{\mathcal M}^{\text{\rm reg}}(X,H;\frak z^-,\frak z^+,\alpha)$
of $\widehat{\mathcal M}(X,H;\frak z^-,\frak z^+,\alpha)$ consisting of all the elements
$((\Sigma,z_-,z_+),u,\varphi)$ such that $\Sigma = S^2$.
Let
$\widetilde{\mathcal M}^{\text{\rm reg}}(X,H;\frak z^-,\frak z^+,\alpha)$,
${\mathcal M}^{\text{\rm reg}}(X,H;\frak z^-,\frak z^+,\alpha)$,
$\overline{\overline{\mathcal M}}^{\text{\rm reg}}(X,H;\frak z^-,\frak z^+,\alpha)$
be the  $\sim_1$, $\sim_2$, $\sim_3$ equivalence classes of
$\widehat{\mathcal M}^{\text{\rm reg}}(X,H;\frak z^-,\frak z^+,\alpha)$, respectively.
\par
It is easy to see that $\widetilde{\mathcal M}^{\text{\rm reg}}(X,H;\frak z^-,\frak z^+,\alpha)$,
${\mathcal M}^{\text{\rm reg}}(X,H;\frak z^-,\frak z^+,\alpha)$
coincide with the ones in Definition \ref{tildeMreg}.
In particular
\begin{equation}\label{Rquotient}
{\mathcal M}^{\text{\rm reg}}(X,H;\frak z^-,\frak z^+,\alpha)
\cong
\widetilde{\mathcal M}^{\text{\rm reg}}(X,H;\frak z^-,\frak z^+,\alpha)/\R.
\end{equation}
Here the $\R$ action is obtained by the translation along the $\R$ direction of the source
and is free.
\par
Moreover we have
\begin{equation}
\overline{\overline{\mathcal M}}^{\text{\rm reg}}(X,H;\frak z^-,\frak z^+,\alpha)
\cong
{\mathcal M}^{\text{\rm reg}}(X,H;\frak z^-,\frak z^+,\alpha)/S^1.
\end{equation}
Here the $S^1$ action is obtained by the $S^1$ action of the source.
However we note that the fiber of the canonical map
\begin{equation}\label{C*action}
{\mathcal M}(X,H;\frak z^-,\frak z^+,\alpha)
\to
\overline{\overline{\mathcal M}}(X,H;\frak z^-,\frak z^+,\alpha)
\end{equation}
between the compactified moduli spaces may be bigger than $S^1$. In fact, if $(\Sigma,z_-,z_+)$ has $k$
mainstream components,
the fiber of $[(\Sigma,z_-,z_+),u,\varphi]$ in
${\mathcal M}(X,H;\frak z^-,\frak z^+,\alpha)$ is $(S^1)^k$ for the generic points.
\par
On the other hand there exists an $S^1$ action on
${\mathcal M}(X,H;\frak z^-,\frak z^+,\alpha)$
obtained by
$$
t_0\cdot [(\Sigma,z_-,z_+),u,\varphi]
= [(\Sigma,z_-,z_+),u,t_0\cdot\varphi]
$$
where
$$
(t_0\cdot \varphi)_a(\tau,t) = \varphi_a(\tau,t+t_0).
$$
\begin{defn}
$$
\overline{\mathcal M}(X,H;\frak z^-,\frak z^+,\alpha) =
{\mathcal M}(X,H;\frak z^-,\frak z^+,\alpha)/S^1.
$$
\end{defn}
The map (\ref{C*action}) factors through $\overline{\mathcal M}(X,H;\frak z^-,\frak z^+,\alpha)$.
\par
We can prove that $\overline{\mathcal M}(X,H;\frak z^-,\frak z^+,\alpha)$ is compact
in the same way as \cite[Theorem 11.1]{FOn}.
\par
In place of taking the quotient by the $\R$ action in (\ref{Rquotient}) we can
require the following balancing condition. (In other words we can
take a global section of this $\R$ action.)
\begin{defn}\label{balancedcond}
Let $((\Sigma,z_-,z_+),u,\varphi) \in \widetilde{\mathcal M}(X,H;\frak z^-,\frak z^+,\alpha)$.
Suppose that it has only one mainstream component.
We define a function $\mathcal A : \R \setminus \text{a finite set} \to \R$ as follows.
\par
Let $\tau_0 \in \R$. We assume $\varphi(\{\tau_0\} \times S^1)$ does not contain a root of the
bubble tree. (This is the way how we remove a finite set from  the domain of $\mathcal A$.)
Let $\Sigma_{{\rm v}_i}$, $i=1,\dots,m$ be the set of the irreducible components that is in a
bubble tree rooted on $\R_{\le \tau_0}\times S^1$. We define
\begin{equation}
\mathcal A(\tau_0)
=
\sum_{i=1}^m \int_{\Sigma_{{\rm v}_i}} u^*\omega
+ \int_{\tau=-\infty}^{\tau_0}\int_{t\in S^1} (u\circ\varphi)^*\omega
+ \int_{t\in S^1} H(u(\varphi(\tau_0,t))) dt.
\end{equation}
This is a nondecreasing function and satisifies
$$
\lim_{\tau\to-\infty} \mathcal A(\tau) = H(\frak z_-),
\qquad
\lim_{\tau\to+\infty} \mathcal A(\tau) = H(\frak z_+) + \alpha \cap \omega.
$$
We say $\varphi$ satisfies the {\it balancing condition} if
\begin{equation}
\lim_{\tau < 0
\atop \tau\to 0} \mathcal A(\tau)
\le
\frac{1}{2} \left(
H(\frak z_-) + H(\frak z_+) + \alpha \cap \omega
\right)
\le
\lim_{\tau > 0
\atop \tau\to 0} \mathcal A(\tau).
\end{equation}
In case $((\Sigma,z_-,z_+),u,\varphi) \in \widetilde{\mathcal M}(X,H;\frak z^-,\frak z^+,\alpha)$
we require the balancing condition in mainstream-component-wise.
\end{defn}
\begin{rem}\label{rem48}
We remark that there exists unique $\varphi$ satisfying the balancing condition in each of the $\R$ orbits,
except the following case: $u$ is constant on the (unique) irreducible component
in the mainstream. (In this case there must occur a nontrivial bubble.)
The uniqueness  breaks down in the case when there exists $\tau_0$ such that
$$
\sum_{i=1}^m \int_{\Sigma_{{\rm v}_i}} u^*\omega = \frac{1}{2} \omega \cap \alpha
$$
in addition.
(Here $\{{\rm v}_i\mid i=1,\dots,m\}$ is the bubbles associated to $\tau_0$ as in Definition \ref{balancedcond}.)
In such a case we replace $\mathcal A$ by the following regularized version
$$
\mathcal A'(\tau_0) = \frac{2}{\sqrt{\pi}}\int_{\R} e^{-(\tau-\tau_0)^2} \mathcal A(\tau) d\tau.
$$
The first derivative of $\mathcal A'$ is always strictly positive.
So there exists a unique $\varphi$ satisfying the modified
balanced condition (using $\mathcal A'$) in each $\R$ orbit.
\par
Note however the balancing condition will be mainly used later to define
canonical marked point. In the case there is a sphere bubble we will not
take a canonical marked point. So this remark is only for consistency of the terminology.
\end{rem}
We have thus defined a compactification of ${\mathcal M}^{\text{reg}}(X,H;\frak z^-,\frak z^+,\alpha)$.
We will construct a Kuranishi structure with corner on it.
The construction is mostly the same as the proof of
Theorem \ref{existsKura}, which  is a detailed version of the proof of \cite[Theorem 7.10]{FOn}.
Here we explain this proof in more detail than \cite{FOn}.
\par
We first remark that $\overline{\overline{\mathcal M}}(X,H;\frak z^-,\frak z^+,\alpha)$
does {\it not} have a Kuranishi structure in general. (Even in the case $\alpha \ne 0$.)
This is because there is an element in this moduli space whose
isotropy group is of positive dimension.
Namely if $\Sigma_i$ is a gradient line component, then
the biholomorphic $S^1$ action on the component $\Sigma_i$ is in the group of
automorphisms of this element of $\overline{\overline{\mathcal M}}(X,H;\frak z^-,\frak z^+,\alpha)$.
So a neighborhood of this element may not be a manifold
with corner.
\par
On the other hand, the $S^1$ action on ${\mathcal M}(X,H;\frak z^-,\frak z^+,\alpha)$
is always locally free (namely its isotropy group is a finite group) in case $\alpha \ne 0$.
In the case of $\alpha =0$ the $S^1$ action on ${\mathcal M}_0(X,H;\frak z^-,\frak z^+,0)$
is locally free.
\par
To define an obstruction bundle on the compactification we need to take an obstruction bundle data
in the same way as Definition \ref{obbundeldata}.
To keep consistency with the fiber product description of the boundary or corner
of  $\overline{\mathcal M}(X,H;\frak z^-,\frak z^+,\alpha)$
we will define it in a way invariant not only under the $S^1$ action but also
under the $(S^1)^k$ action on the part where there are $k$ irreducible components
in the mainstream.
\par
We also need to consider the case of gradient line component at the same time.
Note we assumed that the map $u$ is an immersion at the additional marked points
in Definition \ref{symstabili} (3)
(the definition of symmetric stabilization).
In case of gradient line component, there is no such point.
However since we assumed that the gradient flow of our Hamiltonian $H$ is Morse-Smale
and satisfying Assumption \ref{nondeg3} (2),
we actually do {\it not} need to perturb the equation on such a component.
So our obstruction bundle there is, by definition, a trivial bundle.
(And we do not need to stabilize such a component to define an  obstruction bundle.)
\par
Taking this into account we define an obstruction bundle data in our situation as the
following Definition \ref{obbundeldata}.
\par
\begin{defn}\label{symstab}
A {\it symmetric stabilization} of $((\Sigma,z_-,z_+),u,\varphi)$ is $\vec w$ such that
$\vec w \cap \Sigma_i$ is a symmetric stabilization
of $\Sigma_i$ in the sense of Definition \ref{symstabili}
if $\Sigma_i$
is not a gradient line component,
and $\vec w \cap \Sigma_i = \emptyset$ if $\Sigma_i$ is a gradient line component.
\end{defn}
\begin{defn}\label{canmarked}
Let $\frak p = ((\Sigma,z_-,z_+),u,\varphi)$
be as above. We assume $\Sigma_i$ is a gradient line component.
Note, then $\ell(\tau) = u(\varphi_i(\tau,t))$ is  a gradient line joining transit images $\frak z_i$ and $\frak z_{i+1}$.
There exists a unique $\tau_0$ such that
$$
H(\ell(\tau_0)) = \frac{1}{2}\left(
H(\frak z_i) + H(\frak z_{i+1})
\right).
$$
We put $w_i = \varphi_i(\tau_0,0)$, which we call the {\it canonical marked point} of the
gradient line component.
\end{defn}
\begin{rem}
We note that the pair $(\frak p,w_i)$ where $w_i$ is the canonical marked point
depends only on the $\sim_3$ equivalence class of $\frak p$ in the following sense.
Suppose $\frak p \sim_3 \frak p'$. We define $w'_i$ in the $i$-th mainstream component
of $\Sigma_{\frak p'}$ as above. Then there  exists an isomorphism $v : \Sigma_{\frak p} \to \Sigma_{\frak p'}$ satisfying Definition \ref{3equivrel}  (1), (2)
and such that $v(w_i) =w'_i$. This is because $u$ is $S^1$ equivariant on this irreducible component.
\par
On the other hand, if the homology class $u_*([\Sigma_i])$ is nonzero,
the pair $(\frak p,w_i)$
is not $\sim_3$ equivalent to $(\frak p,w'_i)$ in the above sense, where $w'_i = \varphi_i(\tau_0,t_0)$.
\end{rem}
We note that, if the $i$-th mainstream component $\Sigma_i$ consists of a gradient line and sphere bubbles,
then we put $\vec w$ only on the part of the sphere bubble.
Note the irreducible component that is the intersection of this mainstream component and the mainstream
is (source) stable since the root of the bubble is the third marked point.
\begin{defn}\label{obbundeldata1}
An {\it obstruction bundle data $\frak E_{\frak p}$ centered at}
$$
[\frak p] = [(\Sigma,z_-,z_+),u,\varphi] \in \overline{\overline{\mathcal M}}(X,H;\frak z^-,\frak z^+,\alpha)
$$
is the data satisfying the conditions described below.
We put $\frak x = (\Sigma,z_-,z_+)$. Let $\frak x_i$ be the $i$-th mainstream component. (It has two marked points.)
We put $\alpha_i = u_*[\Sigma_i] \in \Pi$.
\begin{enumerate}
\item
A symmetric stabilization $\vec w$ of $\frak p$.
We put $\vec w^{(i)} = \vec w \cap \frak x_i$.
\item
The same as Definition \ref{obbundeldata} (2).
\item
A universal family with coordinate at infinity of  $\frak x_{\frak  p} \cup \vec w \cup \vec w^{\rm can}$. Here we put
the canonical marked point (Definition \ref{canmarked}) for each gradient line component
and denote them by
$\vec w^{\rm can}$.
We require some additional condition (Condition \ref{condcoord} below) for the coordinate at infinity.
\item
The same as Definition \ref{obbundeldata} (4). Namely compact subsets $K^{\rm obst}_{\rm v}$ of $\Sigma_{\rm v}$.
(The support of the obstruction bundle.) We do not put $K^{\rm obst}_{\rm v}$  on the gradient line
components. In case the $i$-th mainstream component $\Sigma_i$ consists of
a gradient line and sphere bubbles,
then we put $K^{\rm obst}_{\rm v}$ only on the bubbles.
\item
The same as Definition \ref{obbundeldata} (5). Namely, finite dimensional complex linear subspaces $E_{\frak p,{\rm v}}(\frak y,u)$.
We do not put them on the gradient line components.
In case $i$-th mainstream component $\Sigma_i$ consists of the gradient line and sphere bubbles,
then we put them only on the bubbles.
\item
The same as Definition \ref{obbundeldata} (6) except the
differential operator there
\begin{equation}\label{ducomponents}
\aligned
\overline D_{u} \overline\partial : &L^2_{m+1,\delta}((\Sigma_{\frak y_{\rm v}},\partial \Sigma_{\frak y_{\rm v}});
u^*TX, u^*TL) \\
&\to
L^2_{m,\delta}(\Sigma_{\frak y_{\rm v}}; u^*TX \otimes \Lambda^{0,1})/E_{\frak p,{\rm v}}(\frak y,u)
\endaligned
\end{equation}
is replaced by the linearization of the equation (\ref{Fleq})
\item
The same as Definition \ref{obbundeldata} (7).
\item
We take a codimension 2 submanifold $\mathcal D_j$ for each of $w_j \in \vec w$ in the same way as Definition \ref{obbundeldata} (8).
We note that we do {\it not} take such submanifolds for the canonical marked points $\in \vec w^{\rm can}$.
(In fact since $u$ is not an immersion at the canonical marked points we can not choose such submanifolds.)
\end{enumerate}
\par
We require that the data $K^{\rm obst}_{\rm v}$, $E_{\frak p,{\rm v}}(\frak y,u)$ depend only on
the mainstream component
$\frak p_i = [(\Sigma_i,z_{i-1},z_i),u,\varphi]$
(where $z_i$ is the $i$-th transit point) that contains the $\rm v$-th irreducible component.
We call this condition {\it mainstream-component-wise}.
\end{defn}
The additional condition we assume in Item (3) above is as follows.
\begin{conds}\label{condcoord}
\begin{enumerate}
\item
Let $z_{i+1}$ be the $(i+1)$-th transit point, which is contained in $\Sigma_i$ and $\Sigma_{i+1}$.
Then the coordinate at infinity near $z_{i+1}$ coincides with the parametrization $\varphi_i$ or $\varphi_{i+1}$
up to the $\R \times S^1$ action. Namely it is  $(\tau,t) \mapsto \varphi_{i}(\tau+\tau_0,t+t_0)$
(resp. $\varphi_{i+1}(\tau+\tau_0,t+t_0))$ for some $\tau_0$ and $t_0$.
\item
Let $z$ be a singular point that is not a transit point and $\Sigma_{\rm v}$ an irreducible
component containing $z$.
Since $\Sigma_{\rm v}$ is a sphere there exists a biholomorphic map
$$
\phi : \Sigma_{\rm v} \cong \C \cup \{\infty\}
$$
such that $\phi(z) = 0$.
\par
Then the coordinate at infinity around $z$ is given as
$$
(\tau,t) \mapsto \phi^{-1}(e^{\pm 2\pi(\tau+\sqrt{-1}t)}),
$$
for some choice of $\phi$. Here $\pm$ depends on  the orientation of the edge corresponding to
$z$.
\end{enumerate}
\end{conds}
Note that the choice of coordinate at infinity satisfying the above condition is not unique.
\begin{rem}
In (2) above we make full use of the fact that our curve is of genus $0$. The construction developed in Part \ref{secsimple} and Part \ref{generalcase}
is designed so that it works in the case of arbitrary genus without change.
So we did not put this condition in Part \ref{secsimple} and Part \ref{generalcase}.
Condition \ref{condcoord} (2)  will be used to simplify the discussion on how to handle
the Hamiltonian perturbation in the gluing analysis. See Lemma \ref{gluehamexpondec}.
\end{rem}
\par
We can prove existence of an obstruction bundle data in the same way as Lemma \ref{existobbundledata}.
For example, we can choose the marked points $\vec w$ as follows:
We note that the restriction of $u$ to the irreducible component $\Sigma_{\rm v}$ is not homologous to zero
except the following two cases. So we can
find a point of $\Sigma_{\rm v}$ at which $u$ is an immersion and take it as an additional marked point.
\begin{enumerate}
\item
$\Sigma_{\rm v}$ is in the mainstream and is not a root of the sphere bubble.
\item
$\Sigma_{\rm v}$ is in the mainstream and is  a root of the sphere bubble.
\end{enumerate}
In Case (1), we take only the canonical marked point on this component.
In Case (2), this irreducible component is stable. So we do not take
additional marked points on this component.
Thus we can define $\vec w$.
\par\smallskip
We take and fix an obstruction bundle data for each element of $\overline{\overline{\mathcal M}}(X,H;\frak z^-,\frak z^+,\alpha)$.
\par
We defined the moduli spaces $\widehat{\mathcal M}(\frak z^-,\frak z^+,\alpha)$,
${\mathcal M}(\frak z^-,\frak z^+,\alpha)$,
and  $\overline{\overline{\mathcal M}}(\frak z^-,\frak z^+,\alpha)$ in Definition \ref{defn45}.
We add $\ell$ additional marked points on it
and
denote the moduli space of such objects as
$\widehat{\mathcal M}_{\ell}(\frak z^-,\frak z^+,\alpha)$,
${\mathcal M}_{\ell}(\frak z^-,\frak z^+,\alpha)$,
and  $\overline{\overline{\mathcal M}}_{\ell}(\frak z^-,\frak z^+,\alpha)$.
We denote by
$\widehat{\mathcal M}_{\ell}(\frak z^-,\frak z^+,\alpha;\mathcal G)$,
${\mathcal M}_{\ell}(\frak z^-,\frak z^+,\alpha;\mathcal G)$,
and  $\overline{\overline{\mathcal M}}_{\ell}(\frak z^-,\frak z^+,\alpha;\mathcal G)$
its subset so that its combinatorial type is $\mathcal G$.
(We include the datum on how
the additional marked points $\vec w$ in
$\mathcal G$ are distributed over the irreducible
components.)
\par
Let $((\Sigma,z_-,z_+),u,\varphi) \cup \vec w\cup \vec w^{\rm can} = \frak p \cup \vec w\cup \vec w^{\rm can}$ be as in Definition \ref{symstab} with decomposition (\ref{stedec}).
We put $\ell = \#(\vec w\cup \vec w^{\rm can})$ and
denote by  $\overline{\overline{\frak V}}(\frak p \cup \vec w\cup  \vec w^{\rm can})$  a neighborhood of
$\frak x_{\frak p} \cup \vec w\cup  \vec w^{\rm can}$ in $\overline{\overline{\mathcal M}}_{\ell}(\frak z^-,\frak z^+,\alpha;\mathcal G_{\frak p \cup \vec w\cup  \vec w^{\rm can}})$.
\par
In the same way as Definition \ref{def214} we define a map
\begin{equation}\label{doubleoverlinePhi}
\overline{\overline{\Phi}} : \overline{\overline{\frak V}}(\frak p \cup \vec w\cup  \vec w^{\rm can})
\times ((\vec T,\infty] \times S^1) \to \overline{\overline{\mathcal M}}_{\ell}(\frak z^-,\frak z^+,\alpha),
\end{equation}
that is an isomorphism onto an open neighborhood of $[\frak x_{\frak p} \cup \vec w\cup  \vec w^{\rm can}]$.
Here the notation $((\vec T,\infty] \times S^1)$ is similar to
Definition \ref{def29}.
\par
We next add the parametrization $\varphi$ of the mainstream to the map (\ref{doubleoverlinePhi})
and define its ${{\mathcal M}}_{\ell}(\frak z^-,\frak z^+,\alpha)$-version below.
\par
Now we define a manifold with corner $\widetilde D(k;\vec T_0)$ as follows.
We put
\begin{equation}\label{216}
\widetilde {\overset{\circ}D}(k;T_0)\\
=
\{(T_1,\dots,T_k) \in \R^{k} \mid T_{i+1} - T_i \ge T_{0,i}\}.
\end{equation}
We (partially) compactify $\widetilde {\overset{\circ}D}(k;\vec T_0)$ to $\widetilde {D}(k;\vec T_0)$ by admitting $T_{i+1} - T_i = \infty$
as follows.
We put $s'_i = 1/(T_{i+1}-T_i)$ then
$T_1$ and $s'_1,\dots,s'_{k-1}$ define another parameters.
So (\ref{216}) is identified with $\R \times \prod_{i=1}^k(0,1/T_{0,i}]$. We (partially) compactify it to $\R \times \prod_{i=1}^k[0,1/T_{0,i}]$.
By taking the quotients of
$\widetilde {\overset{\circ}D}(k;T_0)$ and
$\widetilde {D}(k;\vec T_0)$ by
the $\R$ action $T(T_1,\dots,T_k) = (T_1+T,\dots,T_k+T)$,
we obtain ${\overset{\circ}D}(k;\vec T_0)$ and $D(k;\vec T_0)$
respectively.
\par
Let ${{\frak V}}(\frak p \cup \vec w\cup  \vec w^{\rm can}) \subset {\mathcal M}_{\ell}(\frak z^-,\frak z^+,\alpha;\mathcal G
_{\frak p \cup \vec w\cup  \vec w^{\rm can}})$ be the inverse image of
$\overline{\overline{\frak V}}(\frak p \cup \vec w\cup  \vec w^{\rm can})$
under the projection
\begin{equation}\label{projectiontovarvar}
{\mathcal M}_{\ell}(\frak z^-,\frak z^+,\alpha;\mathcal G
_{\frak p \cup \vec w\cup  \vec w^{\rm can}})
\to
\overline{\overline{\mathcal M}}_{\ell}(\frak z^-,\frak z^+,\alpha).
\end{equation}
\begin{rem}
Note for an element $((\Sigma',z'_-,z'_+),\varphi') \cup \vec w'$, the marked points $w'_i$ that
correspond to the canonical marked points $\in \vec w^{\rm can}$ may not be canonical.
(Namely it may not be of the form $\varphi'(\tau_0,0)$ where $\tau_0$ is as in Definition \ref{canmarked}.)
\end{rem}
The space ${\mathcal M}_{\ell}(\frak z^-,\frak z^+,\alpha;\mathcal G
_{\frak p \cup \vec w\cup \vec w^{\rm can}})$ carries an $(S^1)^k$ action given by
$$
(t_1,\dots,t_k)(((\Sigma,z_-,z_+),u,\varphi) \cup \vec w\cup \vec w^{\rm can})
=
((\Sigma,z_-,z_+),u,\varphi') \cup \vec w\cup \vec w^{\rm can}
$$
where $\varphi = (\varphi_i)_{i=1}^k$
and  $\varphi' = (\varphi'_i)_{i=1}^k$ such that
$$
\varphi'_i(\tau,t) = \varphi_i(\tau,t+t_i).
$$
This action is locally free and the map (\ref{projectiontovarvar}) can be identified with the canonical
projection:
$$
{\mathcal M}_{\ell}(\frak z^-,\frak z^+,\alpha;\mathcal G
_{\frak p \cup \vec w\cup  \vec w^{\rm can}})
\to
{\mathcal M}_{\ell}(\frak z^-,\frak z^+,\alpha;\mathcal G
_{\frak p \cup \vec w\cup  \vec w^{\rm can}})/(S^1)^k.
$$
\par
It follows from this fact that
${{\frak V}}(\frak p \cup \vec w\cup  \vec w^{\rm can})$ is an open neighborhood of
the inverse image of $[\frak p]$ in ${\mathcal M}_{\ell}(\frak z^-,\frak z^+,\alpha;\mathcal G
_{\frak p \cup \vec w\cup  \vec w^{\rm can}})$.
\par
We now define:
\begin{equation}\label{418}
\overline{\Phi} : {{\frak V}}(\frak p \cup \vec w\cup  \vec w^{\rm can}) \times {D}(k;\vec T_0)
\times (\prod_{j=1}^m(T_{0,j},\infty] \times S^1)/\sim) \to {\mathcal M}_{\ell}(\frak z^-,\frak z^+,\alpha)
\end{equation}
that is a homeomorphism onto an open neighborhood of the  inverse image of
$[\frak p]$ in ${\mathcal M}_{\ell}(\frak z^-,\frak z^+,\alpha)$.
Here $\sim$ is as in Remark \ref{rem:161}.
\par
Let $\frak y \in {{\frak V}}(\frak p \cup \vec w\cup  \vec w^{\rm can})$ and $(\vec T,\vec \theta) \in \widetilde D(k;\vec T_0)
\times  (\prod_{j=1}^m(T_{0,j},\infty] \times S^1)/\sim)$. Note ${\frak V}(\frak p \cup \vec w \cup  \vec w^{\rm can})$ is the quotient space
of the $\R$ action. We represent the quotient by the slice
obtained by requiring the balancing condition
  (Definition \ref{balancedcond}).
\par
The element $\frak y$
comes with the coordinate around the $m$ singular points of $\Sigma$ that are not transit points.
(This is a part of the stabilization data centered at $\frak p$.)
We use the parameters in $(\prod_{j=1}^m(T_{0,j},\infty] \times S^1)/\sim$ to resolve those singular points.
\par
The rest of the parameter
$\vec T' = (T_1,\dots,T_k) \in D(k_1,k_2;\vec T_0)$ is used to resolve the
transit points as follows.
We consider the case this parameter $\vec T' $ is in
$\overset{\circ}D(k;\vec T_0)$.
Let us consider
$$
[-5T_i,5T_i]_i \times S^1_i
$$
and regard it as a subset of the domain of $\varphi_i : \R \times S^1 \to \Sigma_i$.
We define
$$
\varphi_0 : \bigcup_i ([-5T_i,5T_i]_i \times S^1_i)  \to \R \times S^1
$$
as follows. If $(\tau,t) \in [-5T_i,5T_i]_i \times S^1_i$ then
$$
\varphi_0(\tau,t) = (\tau+10 T_i,t).
$$
We use $\varphi_0\circ\varphi_i^{-1}$ to identify (a part of) $\Sigma_i$ with a subset of
$\R \times S^1$. We then use this identification to move bubble components (glued) and
marked points. So together with $\R \times S^1$ it gives a
marked Riemann surface. We thus obtain ${\frak Y}$.
The image of $\varphi_0$ is in the {\it core} of ${\frak Y}$ and the complement of
the core in the mainstream is the neck region.
\par
By taking the quotient with respect to the $\R$ action, we obtain:
$$
{\frak Y} = \overline{\Phi}(\frak y,\vec T,\vec \theta) \in {\mathcal M}_{\ell}(\frak z^-,\frak z^+,\alpha).
$$
We have thus defined
(\ref{418}).
\par\smallskip
We next consider $((\Sigma',z'_-,z'_+),u',\varphi')\cup \vec w'$ and define the $\epsilon$-close-ness of it from
$[\frak p \cup \vec w\cup  \vec w^{\rm can}] = [((\Sigma,z_-,z_+),u,\varphi) \cup \vec w\cup  \vec w^{\rm can}]$.
Here $((\Sigma',z'_-,z'_+),u',\varphi')$ is assumed to satisfy Definition \ref{defn210} (1)(2)(3)(6)(7)
and $\vec w'$ is the set of $\ell$ additional marked points.
We decompose
\begin{equation}\label{stedec2}
\Sigma' = \sum_{j=1}^{k'} \Sigma'_{j}
\end{equation}
into the mainstream components. Let $z'_j$ be the $j$-th transit point and
$\alpha'_j = u'_*([\Sigma'_{j}])$.
We assume there exists a map $i: \{1,\dots,k'\} \to \{1,\dots,k\}$
such that
\begin{enumerate}
\item[(a)]
$u(z_{i(j)}) = u'(z'_j)$.
\item[(b)]
$\sum_{i=i(j)}^{i(j+1)-1}\alpha_i = \alpha'_j$.
\end{enumerate}
Here $z_i$ is the $i$-th transit point of $\frak p$ and
$\alpha_i = u_*([\Sigma_{i}])$.
\par
\begin{defn}\label{epsiclose}
We say $((\Sigma',z'_-,z'_+),u',\varphi') \cup \vec w'$ is {\it $\epsilon$-close} to $[\frak p \cup \vec w\cup  \vec w^{\rm can}]$ if the
following holds.
\begin{enumerate}
\item
\begin{equation}\label{Phioverline2}
((\Sigma',z'_-,z'_+),u') \cup \vec w' = \overline{\overline{\Phi}}(\overline{\frak y},\vec T,\vec \theta)
\end{equation}
where $\overline{\frak y} \in \overline{\overline{\frak V}}(\frak p \cup \vec w\cup  \vec w^{\rm can})$.
And Definition \ref{epsiloncloseto} (1) holds.
\item
Definition \ref{epsiloncloseto} (2) holds.
\item
Definition \ref{epsiloncloseto} (3) holds.
\item
Definition \ref{epsiloncloseto} (4) holds. (Namely each of the component of $\vec T$ is $> 1/\epsilon$.)
\item
If $\Sigma_i$, $i=i(j),\dots,i(j+1)-1$ are all gradient line components,
then $\Sigma'_{j}$ is also a gradient line component
that is $\epsilon$-close to the union of the gradient lines $u\vert_{\Sigma_{i}}$, $i=i(j),\dots,i(j+1)-1$.
We also require that
$\vec w' \cap \Sigma'_{j}$ consists of $i(j+1) - i(j)$ points $z'_{i(j)+1},\dots,z'_{i(j+1)}$ such that
$$
\left\vert H(z_i) - \frac{1}{2}\left(
H(\frak z_i) + H(\frak z_{i+1})
\right)\right\vert < \epsilon
$$
for $i(j) \le i \le i(j+1)-1$.
\end{enumerate}
\end{defn}
If $((\Sigma',z'_-,z'_+),u',\varphi') \cup \vec w'$ is $\epsilon$-close to $[\frak p \cup \vec w\cup  \vec w^{\rm can}]$
and $\epsilon$ is sufficiently small,
we can define an obstruction bundle for $((\Sigma',z'_-,z'_+),u',\varphi')$ in a way similar to Definition \ref{Emovevvv} as follows.
\begin{defn}\label{obstructionbundle}
We consider the decomposition (\ref{stedec2}) of $\Sigma'$ into the
mainstream components and define $i(j)$ as in (a),(b) there.
We will define an obstruction bundle supported on each of $\Sigma'_j$.
\par
If $\Sigma'_j$ is a gradient line component, we set an obstruction bundle to be trivial on $\Sigma'_j$.
\par
Suppose $\Sigma'_j$ is not a gradient line component.
We remove all the  marked points $\vec w' \cap \Sigma'_j$ that correspond to
$\vec w^{\rm can}$.  We denote by $\vec w'_{0,j} \subseteq \vec w' \cap \Sigma'_j$  the
remaining marked points on $\Sigma'_j$.
It is easy to see that
$(\Sigma'_j;z'_j,z'_{j+1})\cup \vec w'_{0,j}$ is stable.
Let $\vec w_{0,j}\subseteq \vec w$ be the set of the marked points on $\Sigma$ corresponding to the marked points $\vec w'_{0,j}$.
\par
In the union $\bigcup_{i=i(j)}^{i(j+1) -1}\Sigma_i$, we shrink each of the
gradient line components $\Sigma_i$ to a point.
Let $\Sigma_{0,j}$ be the resulting semi-stable curve. Then $(\Sigma_{0,j};z_{i(j)},z_{i(j+1)}) \cup  \vec w_{0,j}$
is stable. It has a coordinate at infinity induced by one given in Definition \ref{obbundeldata} (3).
We remark that the union of the supports of the obstruction bundles in $\bigcup_{i=i(j)}^{i(j+1) -1}\Sigma_i$
may be regarded as subsets of $\Sigma_{0,j}$.
\par
We observe that $(\Sigma'_j;z'_j,z'_{j+1})\cup \vec w'_{0,j}$ is obtained from  $(\Sigma_{0,j};z_{i(j)},z_{i(j+1)}) \cup  \vec w_{0,j}$
by resolving singular points. Therefore using
the above mentioned coordinate at infinity we have a diffeomorphism from the
supports of the obstruction bundles in $\bigcup_{i=i(j)}^{i(j+1) -1}\Sigma_i$ onto open subsets of
$\Sigma'_j$, together with parallel transport to define an obstruction bundle
on $\Sigma'_j$, in the same way as Definition \ref{Emovevvv}.
\end{defn}
\begin{rem}\label{remakcannorole}
We remark that by construction the obstruction bundle that we defined on $((\Sigma',z'_-,z'_+),u',\varphi') \cup \vec w'$
is independent of the marked points $\in  \vec w'$ corresponding to the
canonical marked points $\vec w^{\rm can}$. This is important to see that our construction
is $S^1$ equivariant.
\par
In other words, the canonical marked points $\vec w^{\rm can}$ and the corresponding marked points
in $\vec w'$ do {\it not} play an important role in our construction.
We introduce them so that our terminology is as close to the one in Parts \ref{secsimple}--\ref{generalcase} as possible.
\end{rem}

In the same way as Corollary \ref{vindependenEcor}, we can show this obstruction bundle is independent of the equivalence relation
$\sim_i$ ($i=1,2,3$) that is defined in the same way as Definition \ref{3equivrel}.
We can define Fredholm regularity
and evaluation-map-transversality of such obstruction bundle in the same way as Definition \ref{fredreg},
Definition \ref{def:1720}, respectively.
Then an obvious analogue of Proposition
\ref{linearMV} is proved in the same way.
\begin{defn}\label{constrainttt}
Suppose $((\Sigma',z'_-,z'_+),u',\varphi') \cup \vec w'$ is as in Definition \ref{epsiclose}.
We say that it satisfies the {\it transversal constraint} if the following holds.
\begin{enumerate}
\item
Let $w'_i$ be one of the elements of $\vec w'$. If the corresponding
$w_i \in  \vec w\cup  \vec w^{\rm can}$ is contained in $\vec w$ then
$u'(w'_i)$ in contained in the
codimension 2 submanifold $\mathcal D_i$ that are given as a part of Definition \ref{obbundeldata1} (8).
\item
Suppose $\Sigma'_j$ is a gradient line component and  $w'_{i(j)+1},\dots,w'_{i(j+1)} = \vec w' \cap \Sigma'_{j}$. Then
$$
H(u(w'_i)) = \frac{1}{2}\left(
H(\frak z_i) + H(\frak z_{i+1})
\right)
$$
for $i(j)+1 \le i \le i(j+1)$.
\item
Let $w'_1,\dots,w'_n$ be the points in $\vec w'$ corresponding to $\vec w^{\rm can}$.
We may assume that they are all contained in the mainstream (by taking $\epsilon$ small.)
Then we require that the $S^1$ coordinate thereof are all $[0] \in \R/\Z = S^1$.
\end{enumerate}
\end{defn}
Then for each $\frak p \in \overline{\overline{\mathcal M}}(X,H;\frak z^-,\frak z^+,\alpha)$
we fix an obstruction bundle data $\frak C_{\frak p}$ centered at $[\frak p]$.
In particular we have $\vec w_{\frak p}$.
We choose $\epsilon_{\frak p}$ so that the conclusion of
Lemma \ref{transpermutelem} holds.
\par
For each $[\frak p] \in \overline{\overline{\mathcal M}}(X,H;\frak z^-,\frak z^+,\alpha)$
we denote by $\frak W_{\frak p}$ the subset of $\overline{\overline{\mathcal M}}(X,H;\frak z^-,\frak z^+,\alpha)$
consisting of $[\frak p']$ satisfying the following conditions.
There exists $\vec w'$ such that $\frak p' \cup \vec w'$ is $\epsilon_{\frak p}$-close to
$[\frak p \cup \vec w_{\frak p}\cup  \vec w^{\rm can}]$ and
$\vec w'$ satisfies the  transversal constraint.
\par
We then find a finite set $\frak C = \{\frak p_c\}$ such that
$$
\bigcup_{c} {\rm Int}\,\frak W_{\frak p_c} = \overline{\overline{\mathcal M}}(X,H;\frak z^-,\frak z^+,\alpha),
$$
where $\frak W_{\frak p_c}$ consists of elements $\frak p$ with the following property:
there exists $\vec w_{\frak p}$ such that $\frak p \cup \vec w_{\frak p}$ is $\epsilon_{\frak p_c}$-close to $\frak p_c
\cup \vec w_{\frak p_c}\cup \vec w^{\rm can}$ and $\vec w_{\frak p}$ satisfies the transversal constraint.
\par
We can construct such $\frak C = \frak C(\alpha)$
mainstream-component-wise in the following sense.
We decompose $\frak p = \cup \frak p_i$ into its mainstream components.
Then $\frak p \in \frak C(\alpha)$ if and only if $\frak p_i \in \frak C(\alpha_i)$\footnote{For the component $\alpha_i = 0$
we take a sufficiently dense subset of the moduli space of gradient lines and use it.} for each $i$.
As long as we consider a finite number of
$\alpha$'s, we can construct such $\frak C(\alpha)$ inductively over the energy of
$\alpha$.
\par
For each $[\frak p] \in \overline{\overline{\mathcal M}}(X,H;\frak z^-,\frak z^+,\alpha)$
we define the notion of stabilization data in the same way as
Definition \ref{stabdata}.
(In other words, we require Definition \ref{obbundeldata} (1)(2)(3)(8).)
\par
We put
\begin{equation}
\frak C_{[\frak p]} = \{ c\in \frak C(\alpha) \mid [\frak p] \in {\rm Int}\,\frak W_{\frak p_c}\}.
\end{equation}
We note that if $c \in \frak C_{[\frak p]}$ there exists $\vec w^{\frak p}_c$ such that
$\frak p\cup \vec w^{\frak p}_c$ is $\epsilon_c$-close to $\frak p_c \cup \vec w_c \cup \vec w^{\rm can}$
and $\vec w^{\frak p}_c$ satisfies the transversal constraint.
We take such $\vec w^{\frak p}_c$ and fix it.
\par
Let
$$
[\frak p] = [(\Sigma_{\frak p},z_{\frak p,-},z_{\frak p, +}),u_{\frak p},\varphi_{\frak p})] \in
\overline{\overline{\mathcal M}}(X,H;\frak z^-,\frak z^+,\alpha).
$$
\begin{defn}\label{defthickenmoduli}
We define a {\it thickened moduli space}
\begin{equation}\label{thicenmoduli}
{\mathcal M}_{(\ell_{\frak p},(\ell_c))}(X,H;\frak z^-,\frak z^+,\alpha;\frak p;\frak A;\frak B)_{\epsilon_0,\vec T_0}
\end{equation}
to be the set of $\sim_2$ equivalence classes of
$((\frak Y,u',\varphi'),\vec w'_{\frak p},(\vec w'_c))$ with the following properties.
(Here $\frak A \subset \frak B \subset \frak C_{[\frak p]}$.)
\begin{enumerate}
\item
$(\frak Y,u',\varphi') \cup \vec w'_{\frak p}$ is $\epsilon_0$-close to
$\frak p\cup \vec w_{\frak p} \cup \vec w^{\rm can}$.
Here $\ell_{\frak p} = \#\vec w'_{\frak p}$.
\item
$(\frak Y,u',\varphi') \cup \vec w'_{c}$ is $\epsilon_0$-close to
$\frak p\cup \vec w_c^{\frak p}$.
Here $\ell_{c} = \#\vec w'_{c}$.
\item
On the bubble we have
$$
\overline{\partial} u' \equiv  0 \mod \mathcal E_{\frak B}.
$$
Here $\mathcal E_{\frak B}$ is the obstruction bundle defined by  Definition \ref{obstructionbundle}
in the same way as Definition \ref{def:187}.
\item
On the $i$-th irreducible component of the mainstream we consider
$h'_i = u'\circ \varphi'_i$. Then it satisfies
$$
\frac{\partial h'_i}{\partial \tau}
+ J\left(
\frac{\partial h'_i}{\partial t} - \frak X_{H_t}
\right)
\equiv  0 \mod \mathcal E_{\frak B}.
$$
Here $\mathcal E_{\frak B}$ is as in (3).
\end{enumerate}
\end{defn}
The next lemma says that
${\mathcal M}_{(\ell_{\frak p},(\ell_c))}(X,H;\frak z^-,\frak z^+,\alpha;\frak p;\frak A;\frak B)_{\epsilon_0,\vec T_0}
$
carries the following $S^1$ action.
\begin{lem}\label{S1invarianceMM}
Suppose that $((\frak Y,u',\varphi'),\vec w'_{\frak p},(\vec w'_c))$
is an element of the moduli space
${\mathcal M}_{(\ell_{\frak p},(\ell_c))}(X,H;\frak z^-,\frak z^+,\alpha;\frak p;\frak A;\frak B)_{\epsilon_0,\vec T_0}
$
and $t_0 \in S^1$. Then
$((\frak Y,u',t_0\cdot\varphi'),t_0\vec w'_{\frak p},(t_0\vec w'_c))$
is an element of
${\mathcal M}_{(\ell_{\frak p},(\ell_c))}(X,H;\frak z^-,\frak z^+,\alpha;\frak p;\frak A;\frak B)_{\epsilon_0,\vec T_0}$.
Here
$
(t_0\cdot \varphi')(\tau,t) = \varphi'(\tau,t+t_0)
$
and $t_0\vec w'_c$ is defined as follows. $(t_0\vec w'_c)_i = (\vec w'_c)_i$ if it corresponds to a marked point in
$\vec w_{\frak p_c}$. If $(\vec w'_c)_i$ corresponds to a canonical marked point then
$(t_0\vec w'_c)_i = ((t_0\cdot \varphi') \circ(\varphi)^{-1})((\vec w'_c)_i)$.
\end{lem}
\begin{proof}
The only part of our construction which potentially breaks
the $S^1$ symmetry
is Definition \ref{constrainttt} (3). However as we remarked in Remark \ref{remakcannorole},
the marked points that correspond to the canonical marked points
do not affect the obstruction bundle. Therefore the $S^1$ symmetry is not broken.
\end{proof}
\begin{defn}
We denote by
$V_{(\ell_{\frak p},(\ell_c))}(X,H;\frak z^-,\frak z^+,\alpha;\frak p;\frak A;\frak B)_{\epsilon_0}$
the subset of ${\mathcal M}_{(\ell_{\frak p},(\ell_c))}(X,H;\frak z^-,\frak z^+,\alpha;\frak p;\frak A;\frak B)_{\epsilon_0,\vec T_0}
$ consisting of the elements with the same combinatorial type as $\frak p$.
(Compare \eqref{2193}.)
\end{defn}
In the same way as Lemma \ref{lem:191},
$V_{(\ell_{\frak p},(\ell_c))}(X,H;\frak z^-,\frak z^+,\alpha;\frak p;\frak A;\frak B)_{\epsilon_0}$
is a smooth manifold.
In the same way as Lemma \ref{S1invarianceMM} we can show that our space
$$V_{(\ell_{\frak p},(\ell_c))}(X,H;\frak z^-,\frak z^+,\alpha;\frak p;\frak A;\frak B)_{\epsilon_0}$$
has an $(S^1)^k$ action. Here $k$ is the number of mainstream components of $\frak p$.
Let $m$ be the number of singular points of $\Sigma_{\frak p}$ which are not transit points.
\par
Now we have the following analogue of
Theorem \ref{gluethm3}.
\begin{prop}\label{gluingpropham}
There exists a map
$$
\aligned
\text{\rm Glue} :
&V_{(\ell_{\frak p},(\ell_c))}(X,H;\frak z^-,\frak z^+,\alpha;\frak p;\frak A;\frak B)_{\epsilon_0}
\times \left(\prod_{i=1}^m (T_{0,i},\infty] \times S^1)\right)/\sim \times D(k;\vec T'_0)\\
&\to {\mathcal M}_{(\ell_{\frak p},(\ell_c))}(X,H;\frak z^-,\frak z^+,\alpha;\frak p;\frak A;\frak B)_{\epsilon_2}.
\endaligned
$$
Its image contains
${\mathcal M}_{(\ell_{\frak p},(\ell_c))}(X,H;\frak z^-,\frak z^+,\alpha;\frak p;\frak A;\frak B)_{\epsilon_3}$
for sufficiently small $\epsilon_3$.
\par
An estimate similar to Theorem \ref{exdecayT33} also holds.
\end{prop}
(The notation ${\mathcal M}_{(\ell_{\frak p},(\ell_c))}(X,H;\frak z^-,\frak z^+,\alpha;\frak p;\frak A;\frak B)_{\epsilon}$
is similar to one used in Theorem \ref{gluethm3}.)
\begin{proof}
The proof is mostly the same as the proofs of Theorems \ref{gluethm3}, \ref{exdecayT33}.
The only new point we need to discuss is the following.
\par
Our equation is the pseudo-holomorphic curve equation ((3) in Definition \ref{defthickenmoduli})
on the bubble but involves Hamiltonian vector field
((4) in Definition \ref{defthickenmoduli}) on the mainstream.
When we resolve the singular points, the bubble becomes the mainstream.
So we need to estimate the contribution of the (pull back by appropriate diffeomorphisms of the)
Hamiltonian vector field on such a
part. We need to do so by using the coordinate that is similar to those we used in
the proofs of Theorems \ref{gluethm3}, \ref{exdecayT33}.
We will use Lemma \ref{gluehamexpondec} below for this purpose.
\par
We put
$$
\Sigma_{\frak p} = \bigcup_{i=1}^k\Sigma_{i} \cup \bigcup_{\rm v} \Sigma_{\rm v}
$$
where $\Sigma_i$ are in the mainstream and $\Sigma_{\rm v}$ are the bubbles.
We note that each $\Sigma_{\rm v}$ is a two sphere $S^2$.
Let $z$ be a singular point contained in $\Sigma_{\rm v}$. According to Condition \ref{condcoord},
we have a disk $D_{z,\rm v} \subset \Sigma_{\rm v}$ centered at $z$ on which a
coordinate at infinity is defined. In case $z \in \Sigma_0$, we also have
a disk $D_{z,i} \subset \Sigma_{i}$ on which
a coordinate at infinity is defined.
\par
We also have
\begin{equation}\label{necktransti}
\varphi_i (((-\infty,5T_i] \cup [5T_{i+1},\infty))\times S^1) \subset \Sigma_i.
\end{equation}
The union of the images of (\ref{necktransti}) and the disks $D_{z,\rm v}$ and $D_{z,i}$
are the neck regions. Its complement is called the core. We write
$K_{\rm v}$ or $K_i$ the part of the core in the component $\Sigma_{\rm v}$,
$\Sigma_{i}$, respectively.
\par
Using the coordinate at infinity, we have an embedding
\begin{equation}\label{embeddingii}
i_{\rm v} : K_{\rm v} \to \Sigma', \qquad
i_{i} : K_{i} \to \Sigma'.
\end{equation}
Lemma \ref{gluehamexpondec}  provides an estimate of the maps (\ref{embeddingii}).
We put the metric on $K_{\rm v}$ regarding them as subsets of the sphere.
We put the metric on $K_i$ by regarding them as subsets of $\R \times S^1$.
(They are compact. So actually the choice of the metric does not matter.)
\par
We fix $\vec T$ (the lengths of the neck region when we glue and obtain $\Sigma'$.)
For each $\rm v$ we define $T_{\rm v}$ as follows.
Take a shortest path joining our irreducible component $\Sigma_{\rm v}$ to the
mainstream. Let $z_1,\dots,z_r$ be the singular points contained in this path.
Let $10T_{z_i}$ be the length of the neck region corresponding to the singular point $z_i$.
We put $T_{\rm v} = \sum T_{z_i}$.
\par
We observe that $i_{\rm v}(K_{\rm v})$ is in the mainstream if and only if $T_{\rm v}$ is finite.
\begin{lem}\label{gluehamexpondec}
There exist $C_{\ell}, c_{\ell} > 0$ with the following properties.
\begin{enumerate}
\item
Suppose that $i_{\rm v}(K_{\rm v})$ is in the mainstream.
We then regard
$$
i_{\rm v} : K_{\rm v} \to \R \times S^1.
$$
The $C^{\ell}$ norm of this map is smaller than $C_{\ell}e^{-c_{\ell}T_{\rm v}}$.
In particular the diameter of its image is smaller than $C_{1}e^{-c_{1}T_{\rm v}}$.
\item
We remark that the image of $i_i$ is in certain mainstream component.
So we may regard
$$
i_i : K_{i} \to \R \times S^1.
$$
Note $K_i \subset \R \times S^1$. Then $i_i$ extends to a biholomorphic map of the form
$
(\tau,t) \mapsto (\tau+\tau_0,t+t_0).
$
\end{enumerate}
\end{lem}
\begin{proof}
For simplicity of the notation we consider the case
$$
\Sigma_{\frak p} = (\R \times S^1) \cup S^2 = \Sigma_1 \cup \Sigma_{\rm v}.
$$
Let $T$ be the length of the neck region of $\Sigma'$. Then we have a canonical isomorphism
\begin{equation}\label{sigmaprimeidentified}
\Sigma'
= ((\R \times S^1) \setminus D^2_0) \cup ([-5T,5T]\times S^1) \cup (\C \setminus D^2).
\end{equation}
Here $D^2_0$ is an image of a small disk $\subset \R \times \R$ by the projection
$\R \times \R \to \R \times S^1$ and the disk $D^2$ in the third term is the
disk of radius 1 centered at origin.
The identification (\ref{sigmaprimeidentified}) is a consequence of Condition \ref{condcoord}.
\par
We have a biholomorphic map
$$
I : \Sigma'  \to \R \times S^1
$$
that preserves the coordinates at their two ends.
\par
We can take $I$ as follows.
Let $I_0 : D^2 \to D^2_0$ be the isomorphism that lifts to a homothetic embedding
$D^2 \to \R^2$.
\begin{enumerate}
\item
On $(\R \times S^1) \setminus D^2_0$, the map $I$ is the identity map.
\item
If $z \in \C \setminus D^2$ then
$$
I(z) = I_0(e^{-20\pi T}/z).
$$
\item
If $z = (\tau,t) \in [-5T,5T]\times S^1$ then
$$
I(z) = I_0(e^{-2\pi((\tau+5T)+\sqrt{-1}t)}).
$$
\end{enumerate}
Lemma \ref{gluehamexpondec} is immediate from this description.
\par
The general case can be proved by iterating a similar process.
\end{proof}
\begin{rem}\label{neckremark}
On the part of the neck region $[-4T,4T]\times S^1$ where we perform the gluing construction,
we can prove a similar estimate as Lemma \ref{gluehamexpondec} (1).
\end{rem}
Now we go back to the proof of Proposition \ref{gluingpropham}.
Lemma \ref{gluehamexpondec} (2) implies that on the mainstream the equation
Definition \ref{defthickenmoduli} (4) is preserved by gluing.
Lemma \ref{gluehamexpondec} (1) and Remark \ref{neckremark} imply that
the effect of the Hamiltonian vector field is small in the exponential order
on the other part. Therefore the presence of the Hamiltonian term
does not affect the proof of Theorems \ref{gluethm3}, \ref{exdecayT33} and
we can prove Proposition \ref{gluingpropham} in the same way as
in Section \ref{glueing}.
\end{proof}
We are now in the position to complete the proof of Theorem \ref{existKuran}.
We have defined the thickened moduli space and described it by the gluing map.
The rest of the construction of the Kuranishi structure is
mostly the same as the construction in Sections \ref{cutting}--
\ref{kstructure}.
We mention two points below. Except them there are nothing to modify.
\par\smallskip
\noindent(1)
We consider the process to put (transversal) constraint and forget the marked points.
This was done in Section \ref{cutting} to cut down the thickened moduli space
to an orbifold of correct dimension, which will become our Kuranishi neighborhood.
Here we use the constraint defined in Definition \ref{constrainttt}.
In the case when the marked point $w'_i$ corresponds to one of $\vec w_{\frak p}$ or $\vec w_{\frak p_c}$
that is not a canonical marked point, this process is exactly the same as in Section \ref{cutting}.
\par
In the case of the marked point $w'_i$ corresponds to a canonical marked point of $\frak p$ or $\frak p_c$ we use
Definition \ref{constrainttt} (2),(3). Note that these conditions determine the position of $w'_i$
on $\Sigma'$ uniquely.
On the hand, as we remarked in Remark \ref{remakcannorole}, the marked point $w'_i$
does not affect the obstruction bundle and hence the equations defining
our thickened moduli space. So the discussion of the process to put constraint and
forget such  a marked point is rather trivial.
\par\smallskip
\noindent(2)
Our thickened moduli space has an $S^1$ action. The gluing map we constructed
in Proposition \ref{gluingpropham} is obviously  $S^1$ equivariant.
(The obstruction bundle is invariant under the $S^1$ action as we remarked before.)
Therefore all the construction of the Kuranishi structure is
done in an $S^1$ equivariant way.
Note that we define the $S^1$ action on the thickened moduli space.
The smoothness of this action is fairly obvious.
\par
We note that the group $\Gamma_{\frak p}$ for $\frak p$ (Definition \ref{chertS1act}) in our case
of $\frak p = ((\Sigma,z_-,z_+),u,\varphi)$ consists of maps $v : \Sigma \to \Sigma$
that satisfies Definition \ref{3equivrel} (1)(2) for $((\Sigma',z'_-,z'_+),u,\varphi')
=  ((\Sigma,z_-,z_+),u,\varphi)$ and
$$
v \circ \varphi_a(\tau,t) = \varphi_a(\tau,t+t_0)
$$
for some $t_0 \in S^1$. The groups $\Gamma_{\frak p}$ and $S^1$ generate the group
$G_{\frak p}$.
\par\medskip
The proof of Theorem  \ref{existKuran} is now complete.
\qed
\begin{exm}\label{zeifertgeo}
Let $h : \R \times S^1 \to X$ be a solution of the equation (\ref{Fleq})
(without bubble). We assume that $h$ is injective.
We put
$$
h_k(\tau,t) = h(k\tau,kt) :  \R \times S^1 \to X.
$$
We define $\varphi : \R \times S^1 \to S^2 \setminus \{z_-,z_+\}$ by
$\varphi(\tau,t) = e^{2\pi(\tau+\sqrt{-1}t)}$.
We put $u_k = h_k \circ \varphi^{-1}$.
Then $((S^2,z_-,z_+),u_k,\varphi)$  is an element of ${\mathcal M}(X,H;\frak z^-,\frak z^+,\alpha)$.
Let $t_0\varphi(\tau,t) = \varphi(\tau,t+t_0)$.
It is easy to see that
$$
((S^2,z_-,z_+),u_k,\varphi) \sim_2 ((S^2,z_-,z_+),u_k,t_0\varphi)
$$
if and only if $t_0 = m/k$, $m\in \Z$.
\par
We can take the Kuranishi neighborhood of $\frak p = ((S^2,z_-,z_+),u_k,\varphi) $
of the form $V = S^1 \times V'$, on which the generator of the group $\Gamma_{\frak p} = \Z_k$ acts by
$(t,v) \mapsto (t+1/k,\psi(v))$ where $\psi : V' \to V'$ is not an identity map.
The $S^1$ action is by rotation of the first factor.
Thus the quotient $V/\Gamma_{\frak p}$ is a manifold and
$V/(\Gamma_{\frak p} \times S^1)$ is an orbifold.
See Example \ref{zeifert}.
\end{exm}
\begin{rem}
In this section we studied the moduli space
${{\mathcal M}}(X,H;\frak z^-,\frak z^+,\alpha)$ in the case when $H$ is a time independent
Morse function.
In an alternative approach (such as those \cite[Section 26]{fooospectr}),
we studied the case $H\equiv 0$ using Bott-Morse gluing.
Actually the discussion corresponding to this section is easier in the case of $H\equiv 0$.
In fact in case of $H \equiv 0$, we do not need to study the moduli space
of its gradient lines. (Some argument was necessary to discuss the moduli space of
gradient lines since its element has $S^1$ as an isotropy group. The main part of this
sedition is devoted to this point.)
\par
In \cite{FOn} we used the case $H$ is a time independent
Morse function rather than studying the case $H\equiv 0$ by Bott-Morse theory.
The reason is  that the chain level argument that we need to use for the case $H\equiv 0$ was
not written in detail at the time when  \cite{FOn} was written in 1996.
Now full detail of the chain level argument was written in \cite{fooo:book1}.
So  at the stage of 2012 (16 years after \cite{FOn} was written)
using the case $H\equiv 0$ to calculate
Floer homology of periodic Hamiltonian system is somewhat simpler to write up in detail rather than
using the case when $H$ is a time independent Morse function.
In this section however we focused on the case of time independent Morse function
and written up as much detail as possible
to  convince the readers
\end{rem}

\section{Calculation of Floer homology
of periodic Hamiltonian system}
\label{calcu}

We first prove
Theorem \ref{connectingcompactka}.
The proof is similar to the proof of Theorem \ref{existKuran}.
We indicate below the points where proofs are different.
\par
Let $(\Sigma,z_-,z_+)$ be as in Definition \ref{defn41}.
We define the notion of mainstream, mainstream component,
and transit point in the same way as in Section \ref{construS1equiv}.
\par
Let  $\tilde{\gamma}^{\pm} = ({\gamma}^{\pm},w^{\pm}) \in \tilde{\frak P}(H)$.
(Here $H$ is a time {\it dependent} periodic Hamiltonian.)
Let $((\Sigma,z_-,z_+),u,\varphi)$ be as in Definition \ref{defn210}
such that it satisfies
(1)(2)(4)(5) of  Definition \ref{defn210} and the
following three conditions:
\begin{enumerate}
\item[(3)$'$]
$u : \Sigma \setminus \{\text{transit points}\} \to X$ is a continuous map.
\item[(6)$'$]
There exist $\tilde\gamma_i = (\gamma_i,w_i) \in \tilde{\frak P}(H)$
for $i=1,\dots,k+1$
with $\tilde\gamma_1 = \tilde\gamma^-$,
$\tilde\gamma_{k+1} = \tilde\gamma^+$ such that
$$
\lim_{\tau\to -\infty} u(\varphi_i(\tau,t)) = \gamma_i(t),
\qquad
\lim_{\tau\to +\infty} u(\varphi_i(\tau,t)) = \gamma_{i+1}(t).
$$
Here $\varphi_i : \R \times S^1 \to \Sigma_i$ is the $i$-th component of
$\varphi$.
\item[(7)$'$]
$w_i \# u(\Sigma_i) \sim w_{i+1}$.
\end{enumerate}
We denote by $\widehat{\mathcal M}(X,H;\tilde\gamma^-,\tilde\gamma^+)$
the set of such $((\Sigma,z_-,z_+),u,\varphi)$.
We define equivalence relations $\sim_1$ and $\sim_2$ on it in the same way
as Definition \ref{3equivrel}.
(We do not use $\sim_3$ here.)
We then put
$$\aligned
\widetilde{\mathcal M}(X,H;\tilde\gamma^-,\tilde\gamma^+) &=
\widehat{\mathcal M}(X,H;\tilde\gamma^-,\tilde\gamma^+)/\sim_1, \\
{\mathcal M}(X,H;\tilde\gamma^-,\tilde\gamma^+) &=
\widehat{\mathcal M}(X,H;\tilde\gamma^-,\tilde\gamma^+)/\sim_2.
\endaligned$$
\par\smallskip
We next define a balancing condition.
Let
$((\Sigma,z_-,z_+),u,\varphi) \in \widehat{\mathcal M}(X,H;\tilde\gamma^-,\tilde\gamma^+)$
be an element with one mainstream component.
We define a function $\mathcal A : \R \setminus \text{a finite set} \to \R$ as follows.
\par
Let $\tau_0 \in \R$. We assume $\varphi(\{\tau_0\} \times S^1)$ does not contain a root of the
bubble tree. (This is the way how we remove a finite set from  the domain of $\mathcal A$.)
Let $\Sigma_{{\rm v}_i}$, $i=1,\dots,m$ be the set of the irreducible components that is in a
bubble tree rooted on $\R_{\le \tau_0}\times S^1$. We define
\begin{equation}
\aligned
\mathcal A(\tau_0)
=
\sum_{i=1}^m \int_{\Sigma_{{\rm v}_i}} u^*\omega
&+ \int_{\tau=-\infty}^{\tau_0}\int_{t\in S^1} (u\circ\varphi)^*\omega\\
&+ \int_{t\in S^1} H(t,u(\varphi(\tau_0,t))) dt
+ \int_{D^2} (w^-)^*\omega.
\endaligned
\end{equation}
Note that the action functional $\mathcal A_H$ is defined by
$$
\mathcal A_H(\tilde\gamma) =
\int_{t\in S^1} H(t,u(\gamma(t))) dt
+ \int_{D^2} (w)^*\omega.
$$
The function $\mathcal A(\tau_0) $ is nondecreasing  and satisifies
$$
\lim_{\tau\to-\infty} \mathcal A(\tau) = \mathcal A_H(\tilde\gamma^-),
\qquad
\lim_{\tau\to+\infty} \mathcal A(\tau) = \mathcal A_H(\tilde\gamma^+).
$$
We say $\varphi$ satisfies the {\it balancing condition} if
\begin{equation}
\lim_{\tau < 0
\atop \tau\to 0} \mathcal A(\tau)
\le
\frac{1}{2} \left(
\mathcal A_H(\tilde\gamma^-) + \mathcal A_H(\tilde\gamma^+)
\right)
\le
\lim_{\tau > 0
\atop \tau\to 0} \mathcal A(\tau).
\end{equation}
In case of general $((\Sigma,z_-,z_+),u,\varphi) \in \widehat{\mathcal M}(X,H;\tilde\gamma^-,\tilde\gamma^+)$
we consider the balancing condition in mainstream-component-wise.
\par
In the case there is a mainstream component $\Sigma_i$ such that
$\partial (u\circ \varphi_i)/\partial \tau = 0$ we can apply the method of Remark \ref{rem48}.
\par\smallskip
We next define the notion of canonical marked point.
Let
$\frak p = ((\Sigma,z_-,z_+),u,\varphi) \in {\mathcal M}(X,H;\tilde\gamma^-,\tilde\gamma^+)$.
Let $\Sigma_i$ be its mainstream component. We assume that there is no
sphere bubble rooted on it.
We are given a biholomorphic map
$\varphi_i : \R\times S^1 \to \Sigma_i \setminus \{z_i,z_{i+1}\}$
where $z_i,z_{i+1}$ are transit points on $\Sigma_i$.
We require $\varphi_i$ to satisfy the balancing condition.
Now we define the canonical marked point $w^{\rm can}_i$ on $\Sigma_i$ by
$
w_i = \varphi_i(0,0).
$
Let $\vec w^{\rm can}$ be the totality of all the canonical marked points on $\Sigma$.
\par
A symmetric stabilization of $\frak p = ((\Sigma,z_-,z_+),u,\varphi) \in {\mathcal M}(X,H;\tilde\gamma^-,\tilde\gamma^+)$
is $\vec w$ such that $\vec w \cap \Sigma_0 = \emptyset$ where $\Sigma_0$ is the
mainstream of $\Sigma$ and
$\vec w \cup\vec w^{\rm can}$ is the symmetric stabilization of $(\Sigma,z_-,z_+)$.
\begin{defn}\label{obbundeldata2}
An {\it obstruction bundle data $\frak E_{\frak p}$ centered at}
$$
\frak p = ((\Sigma,z_-,z_+),u,\varphi) \in {\mathcal M}(X,H;\tilde\gamma^-,\tilde\gamma^+)
$$
is the data satisfying the conditions described below.
We put $\frak x = (\Sigma,z_-,z_+)$. Let $\frak x_i$ be the $i$-th mainstream component. (It has two marked points.)
\begin{enumerate}
\item
A symmetric stabilization $\vec w $ of $\frak p$.
We put $\vec w^{(i)} = \vec w \cap \frak x_i$.
\item
The same as Definition \ref{obbundeldata} (2).
\item
A universal family with coordinate at infinity of $\frak x_{\frak  p} \cup \vec w \cup \vec w^{\rm can}$.
Here $\vec w^{\rm can}$ is the
canonical marked points as above.
We require Condition \ref{condcoord} for the coordinate at infinity.
\item
The same as Definition \ref{obbundeldata} (4). Namely compact subsets $K^{\rm obst}_{\rm v}$ of $\Sigma_{\rm v}$.
(We put $K^{\rm obst}_{\rm v}$ also on the mainstream.)
\item
The same as Definition \ref{obbundeldata} (5). Namely, finite dimensional complex linear subspaces $E_{\frak p,{\rm v}}(\frak y,u)$.
(We put $E_{\frak p,{\rm v}}(\frak y,u)$ also on the mainstream.)
\item
The same as Definition \ref{obbundeldata} (6) except the
differential operator there
\begin{equation}\label{ducomponents22}
\aligned
\overline D_{u} \overline\partial : &L^2_{m+1,\delta}((\Sigma_{\frak y_{\rm v}},\partial \Sigma_{\frak y_{\rm v}});
u^*TX, u^*TL) \\
&\to
L^2_{m,\delta}(\Sigma_{\frak y_{\rm v}}; u^*TX \otimes \Lambda^{0,1})/E_{\frak p,{\rm v}}(\frak y,u)
\endaligned
\end{equation}
is replaced by the linearization of the equation (\ref{Fleq})
\item
The same as Definition \ref{obbundeldata} (7).
\item
We take a codimension 2 submanifold $\mathcal D_j$ for each of $w_j \in \vec w$ in the same way as Definition \ref{obbundeldata} (8).
\end{enumerate}
\par
We require that the data $K^{\rm obst}_{\rm v}$, $E_{\frak p,{\rm v}}(\frak y,u)$ depend only on
$\frak p_i = [(\Sigma_i,z_{i-1},z_i),u,\varphi]$
(where $z_i$ is the $i$-th transit point) that contains the
$\rm v$-th irreducible component.
We call this condition {\it mainstream-component-wise}.
\end{defn}
We define ${\mathcal M}_{\ell}(\tilde\gamma^-,\tilde\gamma^+)$
in the same way as
${\mathcal M}_{\ell}(\frak z^-,\frak z^+,\alpha)$.
(Namely its element is $((\Sigma_i,z_{i},z_{i+1}),\varphi)$ together with $\ell$ additional marked points on
$\Sigma$,
elements $\tilde\gamma_i \in \tilde{\frak P}(H)$ assigned to each of the transit point, and
homology classes of each of the bubbles.)
We denote by ${\mathcal M}_{\ell}(\tilde\gamma^-,\tilde\gamma^+;\mathcal G)$
its subset consisting of elements with given combinatorial type $\mathcal G$.
\par
Let
$
\frak p = ((\Sigma,z_-,z_+),u,\varphi) \in {\mathcal M}(X,H;\tilde\gamma^-,\tilde\gamma^+)
$
and $\vec w \cup \vec w^{\rm can}$ be its symmetric stabilization.
We denote by $\frak V(\frak p \cup \vec w \cup \vec w^{\rm can})$
a neighborhood of $\frak p \cup \vec w \cup \vec w^{\rm can}$
in ${\mathcal M}_{\ell}(\tilde\gamma^-,\tilde\gamma^+;\mathcal G_{\frak p \cup \vec w \cup \vec w^{\rm can}})$.
\par
We can define
\begin{equation}\label{4182}
\overline{\Phi} : {{\frak V}}(\frak p \cup \vec w\cup  \vec w^{\rm can}) \times D(k;\vec T_0)
\times (\prod_{j=1}^m(T_{0,j},\infty] \times S^1)/\sim) \to {\mathcal M}_{\ell}(\tilde\gamma^-,\tilde\gamma^+)
\end{equation}
in the same way as (\ref{418}).
\par
We say $((\Sigma',z'_-,z'_+),u',\varphi') \cup \vec w'$ is
{\it $\epsilon$-close} to
$\frak p \cup \vec w\cup  \vec w^{\rm can}$, if
(1)(2)(3)(4) of Definition \ref{epsiclose} hold and if
\begin{enumerate}
\item[(5)$'$]
Let $w'_j = \varphi'_j(\tau_j,t_j) \in \vec w'$ be the marked point corresponding to the
canonical marked point $\in \vec w^{\rm can} \cap \Sigma_i$.
Then
$$
\left\vert
\mathcal A(\tau_j)
-
\frac{1}{2} \left(
\mathcal A_H(\tilde\gamma_i) + \mathcal A_H(\tilde\gamma_{i+1})
\right)
\right\vert < \epsilon
$$
and
$
\vert t_j - 0\vert < \epsilon.
$
\end{enumerate}
\par
If $((\Sigma',z'_-,z'_+),u',\varphi') \cup \vec w'$ is $\epsilon$-close to
$\frak p \cup \vec w\cup  \vec w^{\rm can}$
and the obstruction bundle data at $\frak p$ is given, then
they induce an obstruction bundle at $((\Sigma',z'_-,z'_+),u',\varphi') \cup \vec w'$
in the same way as Definition \ref{obstructionbundle}.
\begin{defn}\label{constrainttt2}
We say that $((\Sigma',z'_-,z'_+),u',\varphi') \cup \vec w'$ satisfies the {\it transversal constraint} if the following holds.
\begin{enumerate}
\item
The same as Definition \ref{constrainttt} (1).
\item
Let $w'_j = \varphi'_j(\tau_j,t_j) \in \vec w'$ be the marked point corresponding to the
canonical marked point $\in \vec w^{\rm can} \cap \Sigma_i$.
Then
$$
\mathcal A(\tau_j)
=
\frac{1}{2} \left(
\mathcal A_H(\tilde\gamma_i) + \mathcal A_H(\tilde\gamma_{i+1})
\right).
$$
\item
In the situation of (2) we have
$
t_j = [0].
$
\end{enumerate}
\end{defn}
For each $\frak p \in {\mathcal M}(X,H;\tilde\gamma^-,\tilde\gamma^+)$
we take a stabilization data $\frak C_{\frak p}$ centered at $\frak p$.
We also take $\epsilon_{\frak p}$ so that if
$((\Sigma',z'_-,z'_+),u',\varphi') \cup \vec w'$ is $\epsilon_{\frak p}$-close to
$\frak p \cup \vec w\cup  \vec w^{\rm can}$ then
Fredholm regularity (see Definition \ref{fredreg}) and evaluation map transversality (see Definition \ref{def:1720})
hold for
$((\Sigma',z'_-,z'_+),u',\varphi') \cup \vec w'$.
\par
Let $\frak W_{\frak p}$ be the set of elements $\frak p'$ with the following property:
there exists $\vec w'$ such that $\frak p' \cup \vec w'$ is $\epsilon_{\frak p}$-close to $\frak p \cup \vec w^{\frak p} \cup \vec w^{\rm can}$ and
$\vec w'$ satisfies the transversal constraint in the sense of Definition \ref{constrainttt2}.
\par
We use it to find a finite set $\frak C = \{\frak p_c\}$ such that
$$
\bigcup_{c} {\rm Int}\,\frak W_{\frak p_c} = {\mathcal M}(X,H;\tilde\gamma^-,\tilde\gamma^+).
$$
We then define
$$
\frak C_{\frak p} = \{ c \in \frak C \mid \frak p \in \frak W_{\frak p_c}\}.
$$
\begin{defn}\label{defthickenmoduli2}
We define a {\it thickened moduli space}
\begin{equation}\label{thicenmoduli22}
{\mathcal M}_{(\ell_{\frak p},(\ell_c))}(X,H;\tilde\gamma^-,\tilde\gamma^+;\frak p;\frak A;\frak B)_{\epsilon_0,\vec T_0}
\end{equation}
to be the set of $\sim_2$ equivalence classes of
$((\frak Y,u',\varphi'),\vec w'_{\frak p},(\vec w'_c))$ with the following properties.
(Here $\frak A \subset \frak B \subset \frak C_{\frak p}$.)
\begin{enumerate}
\item
$(\frak Y,u',\varphi') \cup \vec w'_{\frak p}$ is $\epsilon_0$-close to
$\frak p\cup \vec w_{\frak p} \cup \vec w^{\rm can}$.
Here $\ell_{\frak p} = \#\vec w'_{\frak p}$.
\item
$(\frak Y,u',\varphi') \cup \vec w'_{c}$ is $\epsilon_0$-close to
$\frak p\cup \vec w_c^{\frak p}$.
Here $\ell_{c} = \#\vec w'_{c}$.
\item
On the bubble we have
$$
\overline{\partial} u' \equiv  0 \mod \mathcal E_{\frak B}.
$$
Here $\mathcal E_{\frak B}$ is the obstruction bundle defined by  Definition \ref{obstructionbundle}
in the same way as Definition \ref{def:187}.
\item
On the $i$-th irreducible component of the mainstream we consider
$h'_i = u'\circ \varphi'_i$. Then it satisfies
$$
\frac{\partial h'_i}{\partial \tau}
+ J\left(
\frac{\partial h'_i}{\partial t} - \frak X_{H_t}
\right)
\equiv  0 \mod \mathcal E_{\frak B}.
$$
Here $\mathcal E_{\frak B}$ is as in (3).
\end{enumerate}
\end{defn}
\begin{defn}
We denote by
$V_{(\ell_{\frak p},(\ell_c))}(X,H;\tilde\gamma^-,\tilde\gamma^+;\frak p;\frak A;\frak B)_{\epsilon_0}$
the subset of the thickened moduli space ${\mathcal M}_{(\ell_{\frak p},(\ell_c))}(X,H;\tilde\gamma^-,\tilde\gamma^+;\frak p;\frak A;\frak B)_{\epsilon_0,\vec T_0}
$ consisting of elements with the same combinatorial type as $\frak p, \vec w, \vec w_c$.
(Compare \eqref{2193}.)
\end{defn}
The Fredholm regularity and evaluation map transversality imply that the space
$V_{(\ell_{\frak p},(\ell_c))}(X,H;\tilde\gamma^-,\tilde\gamma^+;\frak p;\frak A;\frak B)_{\epsilon_0}$
is a smooth manifold.
\begin{prop}\label{gluingpropham2}
There exists a map
$$
\aligned
\text{\rm Glue} :
&V_{(\ell_{\frak p},(\ell_c))}(X,H;\tilde\gamma^-,\tilde\gamma^+;\frak p;\frak A;\frak B)_{\epsilon_0}
\times \left(\prod_{i=1}^m (T_{0,i},\infty] \times S^1)\right)/\sim \times D(k;\vec T'_0)\\
&\to {\mathcal M}_{(\ell_{\frak p},(\ell_c))}(X,H;\tilde\gamma^-,\tilde\gamma^+;\frak p;\frak A;\frak B)_{\epsilon_2}.
\endaligned
$$
Its image contains
${\mathcal M}_{(\ell_{\frak p},(\ell_c))}(X,H;\tilde\gamma^-,\tilde\gamma^+;\frak A;\frak B)_{\epsilon_3}$
for sufficiently small $\epsilon_3$.
\par
An estimate similar to Theorem \ref{exdecayT33} also holds.
\end{prop}
The proof is the same as the proof of Proposition \ref{gluingpropham}.
Using Proposition \ref{gluingpropham2} we can prove Theorem
\ref{connectingcompactka}
in the same way as the last step of the proof of Theorem \ref{existKuran}.
\qed
\begin{rem}
In case $H$ in Theorem \ref{connectingcompactka} happens to be time independent,
the Kuranishi structure obtained by Theorem \ref{connectingcompactka} is different
from the one obtained by Theorem \ref{existKuran}.
In fact, during the proof of Theorem \ref{existKuran}
we chose a (sufficiently dense finite) subset of
$\overline{\overline{\mathcal M}}(X,H;\frak z^-,\frak z^+,\alpha)$
to define the obstruction bundle.
During the proof of Theorem \ref{connectingcompactka}
we chose a (sufficiently dense finite) subset of
${\mathcal M}(X,H;\tilde\gamma^-,\tilde\gamma^+)$ for the same purpose.
The elements of the moduli space $\overline{\overline{\mathcal M}}(X,H;\frak z^-,\frak z^+,\alpha)$
are $\sim_3$ equivalence classes
and the elements of the moduli space ${\mathcal M}(X,H;\tilde\gamma^-,\tilde\gamma^+)$
are $\sim_2$ equivalence classes.
\end{rem}
Using Theorem \ref{connectingcompactka}  we can define Floer homology
$HF(X,H)$ for time dependent 1-periodic Hamiltonian $H$ satisfying
Assumption \ref{nondeg1}. This construction (going back to Floer \cite{Flo88IV,Flo89I}, see also \cite{HoSa95,Ono95})  is well-established.
We sketch the construction here for completeness.
We use the universal Novikov ring
$$
\Lambda_0 = \left\{\left. \sum_{i=1}^{\infty} a_iT^{\lambda_i}  ~\right\vert~ a_i \in \Q, \,\, \lambda_i \in \R, \,\, \lim_{i\to \infty} \lambda_i = + \infty
\right\}.
$$
Let
$CF(X,H)$ be the free $\Lambda_0 $ module whose basis is identified with the set $\frak P(H)$.
\par
We take $E > 0$.
By Theorem \ref{connectingcompactka}, we obtained a system of Kuranishi structures on ${\mathcal M}(X,H;\tilde\gamma^-,\tilde\gamma^+)$
for each pair $\tilde\gamma^-,\tilde\gamma^+$ with $\mathcal A_H(\tilde\gamma^+) - \mathcal A_H(\tilde\gamma^-) < E$.
We take a system of multisections $\frak s$ on them that are compatible with the description of its boundary and corner
as in Theorem \ref{connectingcompactka} (3). Here we use the fact that
the obstruction bundle is defined mainstream-component-wise.
Note our Kuranishi structure is oriented.
We define
\begin{equation}
\partial^E [\gamma^-]
=
\sum_{\tilde\gamma^+, \mu(\tilde\gamma^+)-\mu(\tilde\gamma^- )= 1
\atop \mathcal A_H(\tilde\gamma^+) - \mathcal A_H(\tilde\gamma^-) < E}
 \#{\mathcal M}(X,H;\tilde\gamma^-,\tilde\gamma^+)^{\frak s}
 T^{\mathcal A_H(\tilde\gamma^+) - \mathcal A_H(\tilde\gamma^-)}[\gamma^+].
\end{equation}
Here we take a lift $\tilde\gamma^-$ of $\gamma^-$ to define the right hand side.
However we can show that the right hand side is independent of the choice of the lift.
The number $\mu(\tilde\gamma^+)$ is the Maslov index.
We have
$$
\dim {\mathcal M}(X,H;\tilde\gamma^-,\tilde\gamma^+)
= \mu(\tilde\gamma^+)-\mu(\tilde\gamma^-) -  1.
$$
Using the moduli space ${\mathcal M}(X,H;\tilde\gamma^-,\tilde\gamma^+)$
with
$\dim {\mathcal M}(X,H;\tilde\gamma^-,\tilde\gamma^+)
= 2$, we can prove
\begin{equation}
\partial^E \circ \partial^E \equiv 0  \mod T^E \Lambda_0
\end{equation}
in a well-established way.
Thus we define
$$
HF(X;H;\Lambda_0/T^E\Lambda_0) \cong H(CF(X,H) \otimes_{\Lambda_0}{\Lambda_0/T^E\Lambda_0},\partial^E).
$$
We can prove that $HF(X;H;\Lambda_0/T^E\Lambda_0)$ is independent of the choice of
Kuranishi structure and its multisection.
(See the proof of Theorem \ref{theorem10} later.)
We thus define
\begin{defn}
$$
HF(X;H;\Lambda_0) = \lim_{\longleftarrow} HF(X;H;\Lambda_0/T^E\Lambda_0).
$$
We also define
$$
HF(X;H) = HF(X;H;\Lambda_0) \otimes_{\Lambda_0}\Lambda
$$
where $\Lambda$ is the field of fractions of $\Lambda_0$.
\end{defn}
In fact, using the next lemma we can find a boundary operator $\partial$ on the full module
$CF(X,H)$ so that its homology is $HF(X;H;\Lambda_0)$.
\begin{lem}\label{homotopyextension}
Let $C$ be a  finitely generated free $\Lambda_0$ module and $E < E'$.
Suppose we are given
$\partial_E : C \otimes_{\Lambda_0} \Lambda_0/T^E\Lambda_0 \to C \otimes_{\Lambda_0} \Lambda_0/T^E\Lambda_0$,
$\partial_{E'} : C \otimes_{\Lambda_0} \Lambda_0/T^{E'}\Lambda_0 \to C \otimes_{\Lambda_0} \Lambda_0/T^{E'}\Lambda_0$
with $\partial_E \circ \partial_E = 0$, $\partial_{E'} \circ \partial_{E'} = 0$.
Moreover we assume $(C \otimes_{\Lambda_0} \Lambda_0/T^E\Lambda_0, \partial_{E'} \mod T^E\Lambda_0)$
is chain homotopy equivalent to $(C \otimes_{\Lambda_0} \Lambda_0/T^E\Lambda_0, \partial_{E})$.
Then we can lift $\partial_E$ to $\partial'_{E'} : C \otimes_{\Lambda_0} \Lambda_0/T^{E'}\Lambda_0 \to C \otimes_{\Lambda_0} \Lambda_0/T^{E'}\Lambda_0$
such that
$(C \otimes_{\Lambda_0} \Lambda_0/T^{E'}\Lambda_0, \partial_{E'})$
is chain homotopy equivalent to $(C \otimes_{\Lambda_0} \Lambda_0/T^{E'}\Lambda_0, \partial'_{E'})$.
\end{lem}
We omit the proof.
\begin{rem}
The method for taking projective limit $E \to \infty$ that we explained above is
a baby version of the one employed in \cite[Section 7]{fooo:book1}. (In \cite{fooo:book1}  the filtered $A_{\infty}$
structure is defined by using a similar method.)
In \cite{Ono95}, a slightly different way to go to projective limit was taken.
\par
For the main application, that is, to estimate the order of $\frak P(H)$ by Betti number, we actually do  not
need to go to the projective limit. See Remark \ref{limitremark}.
\end{rem}
Now we use the $S^1$ equivariant Kuranishi structure in Theorem \ref{existKuran}
to prove the next theorem.
\begin{thm}\label{theorem10}
For any time dependent 1-periodic Hamiltonian $H$ on a
compact manifold $X$ satisfying
Assumption \ref{nondeg1}, we have
$$
HF(X,H) \cong H(X;\Lambda)
$$
where the right hand side is the singular homology group with $\Lambda$
coefficients.
\end{thm}
\begin{proof}
Let $H'$ be a time {\it independent} Hamiltonian satisfying
Assumptions \ref{nondeg2}, \ref{nondeg3}.
We regard $H'$ as a Morse function and let ${\rm Crit}(H')$ be the set of
the critical points of $H'$.
We denote by $CF(X,H';\Lambda_0)$ the free $\Lambda_0$ module
with basis identified with ${\rm Crit}(H')$.
Let $\mu : {\rm Crit}(H') \to \Z$ be the Morse index.
For $\frak x^+,\frak x^- \in {\rm Crit}(H')$ with
$\mu(\frak x^+) - \mu(\frak x^-) = 1$ we define
\begin{equation}
\langle \partial \frak x^-,\frak x^+\rangle
= T^{H'(\frak x^+) - H'(\frak x^-)}\# \mathcal M(X,H',\frak x^-,\frak x^+;0),
\end{equation}
where $\# \mathcal M(X,H',\frak x^-,\frak x^+;0)$ is the number counted with orientation.
(Here $0$ denotes the equivalence class of zero in $\Pi$.)
By Assumptions \ref{nondeg2} this moduli space is smooth.)
It induces $\partial : CF(X,H';\Lambda_0) \to CF(X,H';\Lambda_0)$.
It is by now well established that $\partial\circ \partial = 0$.
We put
$$
HF(X,H';\Lambda_0) = \frac{\text{\rm Ker} \partial}{\text{\rm Im} \partial}.
$$
It is also standard by now that
$$
HF(X,H') = HF(X,H';\Lambda_0) \otimes_{\Lambda_0} \Lambda
$$
is isomorphic to the singular homology $H(X;\Lambda)$ of $\Lambda$ coefficients.
\par
We will construct a chain map from $CF(X,H';\Lambda_0)$ to
$CF(X,H;\Lambda_0)$.
Let $\mathcal H : \R \times S^1 \times X \to \R$ be a smooth map
such that
\begin{equation}
\mathcal H(\tau,t,x)
=
\begin{cases}
H'(x)    &\text{if $\tau \le -1$},
\\
H(t,x)    &\text{if $\tau \ge -1$}.
\end{cases}
\end{equation}
For a map $h : \R \times S^1 \to X$ we consider the equation
\begin{equation}\label{homomoeq}
\frac{\partial h}{\partial \tau}
+
J\left(
\frac{\partial h}{\partial t} - \frak X_{\mathcal H_{\tau,t}}
\right)
 = 0,
\end{equation}
where $\mathcal H_{\tau,t}(x) = \mathcal H(\tau,t,x)$.
\par
Given $\frak x \in \text{\rm Crit}(H')$ and $\tilde\gamma = (\gamma,w)
\in \tilde{\frak P}(H)$ we consider the set of the maps $h$ satisfying
(\ref{homomoeq}) together with the following boundary conditions.
\begin{enumerate}
\item
$
\lim_{\tau \to -\infty} h(\tau,t) = \frak x.
$
\item
$
\lim_{\tau \to +\infty} h(\tau,t) = \gamma(t).
$
\item
$[h] \sim [w]$. Here $h$ is regarded as a map from $D^2$
by identifying $\{-\infty\} \cup ((-\infty,+\infty] \times S^1)$ with $D^2$.
\end{enumerate}
We denote the totality of such $h$ by
${\mathcal M}^{\rm reg}(X,\mathcal H;\frak x,\tilde\gamma)$.
\begin{thm}\label{theorem51}
There exists a compactification ${\mathcal M}(X,\mathcal H;\frak x,\tilde\gamma)$
of ${\mathcal M}^{\rm reg}(X,\mathcal H;\frak x,\tilde\gamma)$,
which is Hausdorff.
\par
For each $E>0$ there exists a system of oriented Kuranishi structures  with corners
on ${\mathcal M}(X,\mathcal H;\frak x,\tilde\gamma)$
for $\mathcal A_H(\tilde \gamma) \le E$.
Its boundary is identified with the union of the following spaces, together with its
Kuranishi structures.
\begin{enumerate}
\item
$$
{\mathcal M}(X,H';\frak x,\frak x',\alpha) \times
{\mathcal M}(X,\mathcal H;\frak x',\tilde\gamma-\alpha)
$$
where $\frak x' \in {\rm Crit}(H')$, $\alpha \in \Pi$
and $\tilde\gamma- \alpha = (\gamma,w-\alpha)$.
\item
$$
{\mathcal M}(X,\mathcal H;\frak x,\tilde\gamma')
\times
{\mathcal M}(X,H,\tilde\gamma',\tilde\gamma)
$$
where $\tilde\gamma' \in \tilde{\frak P}(H)$.
\end{enumerate}
\end{thm}
\begin{proof}
The proof of Theorem \ref{theorem51} is
mostly the same as the proof of Theorems \ref{connectingcompactka}  and \ref{existKuran}.
So we mainly discuss the point where there is a difference in the proof.
\par
Let $(\Sigma,z_-,z_+)$ be a genus zero semi-stable curve with two marked points.
We fix one of the irreducible components in the mainstream and denote it by $\Sigma^0_0$.
We decompose
\begin{equation}
\Sigma = \Sigma^- \cup \Sigma^0 \cup \Sigma^+
\end{equation}
as follows. $\Sigma^0$ is the mainstream component containing $\Sigma^0_0$.
$\Sigma^-$ (resp. $\Sigma^+$) is the connected component of
$\Sigma \setminus \Sigma^0$ containing $z_-$ (resp. $z_+$).
We remark $\Sigma^0$ and/or $\Sigma^-$ may be empty.
\par
We consider $((\Sigma,z_-,z_+), \Sigma^0,u,\varphi)$
such that Definition \ref{defn210}
(1)(2)(4)(5) are satisfied. We assume moreover the following conditions.
\begin{enumerate}
\item[(3.1)$'$]
$\varphi\vert_{\Sigma^+} : \Sigma^+ \to X$ is continuous.
\item[(3.2)$'$]
We put either $\{ z_{0,-} \} = \Sigma^0 \cap \Sigma^-$
or $z_{0,-} = z_-$. (Here we put $z_{0,-} = z_-$ if $z_- \in \Sigma^0$.)
Then
$\varphi\vert_{\Sigma^0 \setminus \{z_{0,-}\}} : \Sigma^0 \setminus \{z_{0,-}\}\to X$ is continuous.
\item[(3.3)$'$]
The same condition as the condition (3)$'$ stated at the beginning of Section \ref{calcu} holds for $\Sigma^+$.
\item[(4.1)$'$]
$h\circ \varphi_i$ satisfies the equation (\ref{Fleq}) for $H'$ if $\Sigma_i \subset \Sigma^-$.
\item[(4.2)$'$]
$h\circ \varphi_0$ satisfies the equation (\ref{homomoeq}). Here $\varphi_0$ is a part
of $\varphi$ and is a parametrization of $\Sigma_0^0$.
\item[(4.3)$'$]
$h\circ \varphi_i$ satisfies the equation (\ref{Fleq}) for $H$ if $\Sigma_i \subset \Sigma^+$.
\item[(6.1)$'$]
Definition \ref{defn210} (6) is satisfied on $\Sigma^-$.
\item[(6.2)$'$]
Let $k-1$ be the number of the mainstream components in $\Sigma^+$.
Then there exist $\tilde\gamma_j  = (\gamma_j,w_j)
\in \tilde{\frak P}(H)$ for $j=1,\dots,k-1$. We have
$$
\lim_{\tau\to \infty} u(\varphi_0(\tau,t)) = \gamma_1(t).
$$
Here $\varphi_0$ is as in (4.2)$'$.
\item[(6.3)$'$]
Let $\Sigma^+_j$ ($j=1,\dots,k-1$) be the irreducible components in the mainstream
of $\Sigma^+$. Let $\varphi^+_j$ is a part of $\varphi$ and is a parametrization of
$\Sigma^+_j$. Then we have
$$
\lim_{\tau\to -\infty} u(\varphi^+_j(\tau,t)) = \gamma_j(t)
\qquad
\lim_{\tau\to \infty} u(\varphi^+_j(\tau,t)) = \gamma_{j+1}(t).
$$
Here $\gamma_{k} = \gamma$.
\item[(7)$'$]
$u_*([\Sigma]) = [w]$.
\end{enumerate}
We denote by $\widehat{\mathcal M}(X,\mathcal H;\frak x,\tilde\gamma)$ the set of all such
$((\Sigma,z_-,z_+), \Sigma^0,u,\varphi)$.
\par
We define three equivalence relations $\sim_1$,  $\sim_2$,  $\sim_3$
on  $\widehat{\mathcal M}(X,\mathcal H;\frak x,\tilde\gamma)$ as follows.
\par
The definition of $\sim_1$ is the same as Definition \ref{3equivrel}.
\par
We apply $\sim_2$ of Definition \ref{3equivrel} on $\Sigma^+$ and $\Sigma^-$
and $\sim_1$ of Definition \ref{3equivrel} on $\Sigma^0$.
This is $\sim_2$ here.
\par
We apply $\sim_2$ of Definition \ref{3equivrel} on $\Sigma^+$,
$\sim_3$ of Definition \ref{3equivrel} on  $\Sigma^-$
and $\sim_1$ of Definition \ref{3equivrel} on $\Sigma^0$.
This is $\sim_3$ here.
We  put
$$\aligned
\widetilde{\mathcal M}(X,\mathcal H;\frak x,\tilde\gamma) &=
\widehat{\mathcal M}(X,\mathcal H;\frak x,\tilde\gamma)/\sim_1, \\
{\mathcal M}(X,\mathcal H;\frak x,\tilde\gamma) &=
\widehat{\mathcal M}(X,\mathcal H;\frak x,\tilde\gamma)/\sim_2, \\
\overline{\overline{\mathcal M}}(X,\mathcal H;\frak x,\tilde\gamma) &=
\widehat{\mathcal M}(X,\mathcal H;\frak x,\tilde\gamma)/\sim_3.
\endaligned$$
\par
We define the notion of balancing condition on the mainstream component
in $\Sigma^-$ and $\Sigma^+$ in the same way as before.
(We do not define such a notion for $\Sigma^0_0$ since the equation
 (\ref{homomoeq}) is not invariant under the $\R$ action.)
\par
We next define the notion of canonical marked points.
For the mainstream components in $\Sigma^+$ or in $\Sigma^-$, the
definition is the same as before.
If $\Sigma^0$ contains a sphere bubble we do not define canonical marked points on it.
Otherwise the canonical marked point on this mainstream component is
$\varphi_0(0,0)$.
\par
Using this notion of canonical marked points we can define the notion of
obstruction bundle data for $[\frak p]\in \overline{\overline{\mathcal M}}(X,\mathcal H;\frak x,\tilde\gamma)$ in the same way as before.
(We put an obstruction bundle also on $\Sigma^0_0$.)
We take and fix an obstruction bundle data  for each of $[\frak p]\in \overline{\overline{\mathcal M}}(X,\mathcal H;\frak x,\tilde\gamma)$.
\par
We can use it to define a map similar to $\overline{\Phi}$ and $\overline{\overline{\Phi}}$
in the same way.
\par
We then use them to define the notion of $\epsilon$-close-ness in the same way.
\par
We next define transversal constraint.
Let $((\Sigma',z'_-,z'_+), \Sigma^{\prime 0},u',\varphi') \cup \vec w'$
is $\epsilon$ close to $((\Sigma,z_-,z_+), \Sigma^0,u,\varphi) \cup \vec w \cup \vec w^{\rm can}$.
We consider $w'_j \in \vec w'$.
If $w'_j$ is either in $\Sigma'_+$ or in $\Sigma'_-$ then the
definition of transversal constraint is the same as Definition \ref{constrainttt} or Definition \ref{constrainttt2}, respectively.
\par
Suppose  $w'_j \in \Sigma'_0$. If $w'_j$ corresponds to a marked point in $\vec w$ then
the definition of transversal constraint is the same as Definition \ref{constrainttt}.
We consider the case where $w'_j$ corresponds to a canonical marked point $w_i$.
There are three cases.
\begin{enumerate}
\item
$w_i \in \Sigma_0$. In this case the transversal constraint is $w'_j = \varphi_0(0,0)$.
\item
$w_i \in \Sigma_-$. Let $\Sigma_{-,i}$  be the mainstream component containing
$w_i$. (It is irreducible since $w_i$ is a canonical marked point.)
The transversal constraint first requires $w'_j = \varphi'_0(\tau_0,0)$ with $\tau_0 \le -1$.
Moreover it requires
$$
\aligned
&\int_{\Sigma_-} (u')^*\omega
+ \int_{\tau=-\infty}^{\tau_0} (u')^*\omega
+ H'(u'(\varphi_0(\tau_0,t)))
\\
&=
\frac{1}{2} \left(
{H'}(u(z_i)) + {H'}(u(z_{i+1}))
\right)
+
\int_{\Sigma_{\tau\le \tau(w_i)}} u^*\omega.
\endaligned
$$
Here $z_i$ and $z_{i+1}$ are transit points contained in $\Sigma_i$  and
$\Sigma_{\tau\le \tau_0}$ is defined as follows.
Let $w_i = \varphi_i(\tau_i,0)$.
We consider $\Sigma \setminus \{ \varphi_i(\tau_i,t) \mid t \in S^1\}$.
Then $\Sigma_{\tau\le \tau_0}$ is the connected component of it
containing $z_-$.
\item
$w_i \in \Sigma_+$. Let $\Sigma_{+,i}$  be the mainstream component containing
$w_i$. (It is irreducible since $w_i$ is a canonical marked point.)
The transversal constraint first requires $w'_j = \varphi'_0(\tau_0,0)$ with $\tau_0 \ge +1$.
Moreover it requires
$$
\aligned
&\int_{\Sigma_-} (u')^*\omega
+ \int_{\tau=-\infty}^{\tau_0} (u')^*\omega
+ \int_{t\in S^1}H(t,u'(\varphi_0(\tau_0,t)))dt
\\
&=
\frac{1}{2} \left(
\mathcal A_{H}(\tilde{\gamma_i}) + \mathcal A_{H}(\tilde{\gamma}_{i+1})
\right).
\endaligned
$$
Here the restriction of $u$ to $\Sigma_{+,i}$ gives an element of
$ {\mathcal M}(X,H;\tilde\gamma_i,\tilde\gamma_{i+1})$.
\end{enumerate}
For a point  $[\frak p]\in \overline{\overline{\mathcal M}}(X,\mathcal H;\frak x,\tilde\gamma)$
we define the notion of stabilization data in the same way as before.
\par
Now using this notion of transversal constraint and $\epsilon$-close-ness,
we define $\frak C_{[\frak p]} \subset \overline{\overline{\mathcal M}}(X,\mathcal H;\frak x,\tilde\gamma)$
for each $[\frak p]\in \overline{\overline{\mathcal M}}(X,\mathcal H;\frak x,\tilde\gamma)$.
We then define a finite set $\frak C$ such that
$$
\bigcup_{[\frak p_c] \in \frak C}\frak C_{[\frak p_c]} = \overline{\overline{\mathcal M}}(X,\mathcal H;\frak x,\tilde\gamma).
$$
We may assume that this choice is mainstream-component-wise in the same sense as before.
\par
We use the choice of $\frak C$ together with obstruction bundle data
we can define an obstruction bundle in the same way as before.
We use it  to define a thickened moduli space.
The rest of the proof is the same as the proof of Theorems \ref{connectingcompactka}  and \ref{existKuran}.
\end{proof}

\begin{lem}\label{512}
There exists a constant $\frak E^-(\mathcal H)$ depending only on $\mathcal H$ with the following properties.
If ${{\mathcal M}}(X,\mathcal H;\frak x,\tilde\gamma)$ is nonempty then
\begin{equation}
\mathcal A_{H}(\tilde\gamma) \ge H'(\frak x) -\frak E^-(\mathcal H).
\end{equation}
\end{lem}
This lemma is classical. See  \cite[(2.14)]{Ono95}.
(It was written also in \cite[Lemma 9]{Floerhofer}.)
\begin{rem}
The optimal estimate for $\frak E^-(\mathcal H)$ is $\frak E^-(\mathcal H) = \mathcal E^-(H-H')$, where
$$
\mathcal E^+(F) = \int_0^1 \max_{p\in X} F(t,p)dt,
\quad
\mathcal E^-(F) = -\int_0^1 \inf_{p\in X} F(t,p)dt.
$$
See \cite[Proposition 3.2]{Oh}.
(See also \cite[Proposition 2.1]{usher:depth}.)
A Lagrangian version of a similar optimal estimate
was obtained by  \cite{Chek96}.
\par
We do not need this optimal estimate for the purpose of this note, but it
becomes important to study spectral invariant.
\end{rem}
We now define
$$
\Phi_E : CF(X,H';\Lambda_0) \to  CF(X,H;\Lambda_0)
$$
as follows. We consider $\frak x$ and $\tilde{\gamma}$ so that:
\begin{enumerate}
\item[(a)] The (virtual) dimension of
${{\mathcal M}}(X,\mathcal H;\frak x,\tilde\gamma)$ is $0$.
\item[(b)]
$
\mathcal A_{H}(\tilde\gamma) \le H'(\frak x) + E.
$
\end{enumerate}
We then put
$$
\langle\Phi_E(\frak x),\tilde\gamma\rangle = T^{\mathcal A_{H}(\tilde\gamma) - H'(\frak x)
+ \frak E^-(\mathcal H)}\#{{\mathcal M}}(X,\mathcal H;\frak x,\tilde\gamma)^{\frak s}.
$$
Here we take and fix a system of multisections $\frak s$
of the moduli space
${{\mathcal M}}(X,\mathcal H;\frak x,\tilde\gamma)$ that is transversal to zero,
compatible with the description of the boundary, and
satisifies the inequality (b) above. We use it to define the right hand side.
\par
We note that the exponent of $T$ in the right hand side is nonnegative because of
Lemma \ref{512}.
\begin{lem}\label{513}
$$
\partial_E \circ \Phi_E - \Phi_E \circ \partial \equiv 0 \mod T^E\Lambda_0.
$$
\end{lem}
\begin{proof}
We use the case of moduli space ${{\mathcal M}}(X,\mathcal H;\frak x,\tilde\gamma)$
satisfying (a) above and has virtual dimension $1$.
It boundary is described as in Theorem \ref{theorem51}.
The case (2) there, counted with sign, gives $\partial_E \circ \Phi_E$.
\par
We consider the case (1) there.
We need to consider the case when virtual dimension of
${\mathcal M}(X,H';\frak x,\frak x',\alpha)$ is zero.
Using $S^1$ equivariance of our Kuranishi structure and multisection,
and Lemma \ref{lem214}, we find that such
${\mathcal M}(X,H';\frak x,\frak x',\alpha)$ is an empty set after perturbation,
unless $\alpha = 0$. In the case $\alpha =0$
we can prove,  in the same way as above, that
$$
{\mathcal M}(X,H';\frak x,\frak x',0) =
{\mathcal M}(X,H';\frak x,\frak x',0)^{S^1}
$$
after perturbation.
Therefore case (1) gives $\Phi_E \circ \partial$.
The proof of Lemma \ref{513} is complete.
\end{proof}
We have thus defined a chain map
\begin{equation}
\Phi_E : CF(X,H';\Lambda_0/T^E\Lambda_0) \to  CF(X,H;\Lambda_0/T^E\Lambda_0).
\end{equation}
We next put
$
\mathcal H'(\tau,t,x) = \mathcal H(-\tau,t,x).
$
We use it in the same way to define the moduli space
${{\mathcal M}}(X,\mathcal H';\tilde\gamma,\frak x)$.
This moduli space has an oriented Kuranishi structure with corners and its
boundary is described in a similar way as Theorem \ref{theorem51}.
(The proof of this fact is the same as the proof of Theorem \ref{theorem51}.)
If it is nonempty then we have
\begin{equation}
H'(\frak x) \ge  \mathcal A_{H}(\tilde\gamma) - \frak E^+(\mathcal H').
\end{equation}
\begin{rem}
Here $\frak E^+(\mathcal H')$ is certain constant depending only on $\mathcal H'$.
The optimal value is $\mathcal E^+(H-H')$.
\end{rem}
We put
$$
\langle\Psi_E(\tilde\gamma),\frak x\rangle = T^{H'(\frak x) -  \mathcal A_{H}(\tilde\gamma)
+ \frak E^+(\mathcal H')}\#{{\mathcal M}}(X,\mathcal H;\frak x,\tilde\gamma)^{\frak s}
$$
and obtain a chain map
\begin{equation}
\Psi_E : CF(X,H;\Lambda_0/T^E\Lambda_0) \to  CF(X,H';\Lambda_0/T^E\Lambda_0) .
\end{equation}
\begin{lem}\label{CC;-}
We may choose
$\frak E^+(\mathcal H')$, $\frak E^-(\mathcal H)$ and  such that
the following holds for $\frak E = \frak E^+(\mathcal H')+\frak E^-(\mathcal H)$.
\begin{enumerate}
\item
$\Psi_E\circ \Phi_E$ is chain homotopic to $[\frak x] \mapsto T^{\frak E}[\frak x]$.
\item
$\Phi_E\circ \Psi_E$ is chain homotopic to $[\tilde{\gamma}] \mapsto T^{\frak E}[\tilde{\gamma}]$.
\end{enumerate}
\end{lem}
\begin{proof}
For $S>0$ we define $\rho_S : \R \to [0,1]$ such that
$$
\rho_S(\tau) =
\begin{cases}
1 &\text{if $\vert\tau\vert < S-1$} \\
0 &\text{if $\vert\tau\vert \ge S$},
\end{cases}
$$
and put
$$
\mathcal H^S(t,x) = \mathcal H(\rho_S(\tau,x)).
$$
For $\frak x_{\pm} \in \text{\rm Crit}(H')$ and $\alpha \in \Pi$, we use
the perturbation of Cauchy-Riemann equation by the Hamilton vector field of
$H_{S,\tau,t}$ to obtain a moduli space ${{\mathcal M}}(X,\mathcal H^S;\frak x_-,\frak x_+;\alpha)$
in the same way. Its union for $S \in [0,S_0]$ also has a Kuranishi structure
whose boundary is given as in Theorem \ref{theorem51} (1), (2)
and ${{\mathcal M}}(X,\mathcal H^S;\frak x_-,\frak x_+;\alpha)$ with $S=0,S_0$.
\par
We consider the case $S=0$. In this case the equation for ${{\mathcal M}}(X,\mathcal H^S;\frak x_-,\frak x_+;\alpha)$
is $S^1$ equivariant. Therefore it has an $S^1$ equivariant Kuranishi structure
that is free for $\alpha \ne 0$.
For $\alpha = 0$ we obtain an $S^1$ equivariant Kuranishi structure
on ${{\mathcal M}}_0(X,\mathcal H^{S=0};\frak x_-,\frak x_+;0)$.
Therefore, by counting the moduli spaces of virtual dimension $0$ we have identity. (It becomes $[\frak x] \mapsto T^{\frak E}[\frak x]$
because of the choice of the exponent in the definition.)
\par
The case $S = S_0$ with $S_0$ huge gives the composition $\Psi_E\circ \Phi_E$.
\par
(1) now follows from a cobordism argument.
\par
The proof of (2) is similar.
\end{proof}
Using \cite[Proposition 6.3.14]{fooo:book1} we have
$$
HF(X,H';\Lambda_0) \cong (\Lambda_0)^{\oplus b'} \oplus \bigoplus_{i=1}^{m'} \frac{\Lambda_0}{T^{\lambda'_i}\Lambda_0}
$$
here $\lambda'_i$, $i=1,\dots,m'$ are positive real numbers.
It implies
$$
HF(X,H';\Lambda) \cong (\Lambda)^{\oplus b'}.
$$
We remark $H(X,H';\Lambda) \cong H(X;\Lambda)$ where the right hand side is the singular homology.
(Note $H'$ is a time independent Hamiltonian that is a Morse function on $X$.)
\par
Similarly we have
$$
HF(X,H;\Lambda_0) \cong (\Lambda_0)^{\oplus b} \oplus \bigoplus_{i=1}^{m} \frac{\Lambda_0}{T^{\lambda_i}\Lambda_0}
$$
and
$$
HF(X,H;\Lambda) \cong (\Lambda)^{\oplus b}.
$$
We take $E$ sufficiently larger than $\frak E$ and $\lambda_i$, $\lambda'_i$. Then
we can use Lemma \ref{CC;-} to show $b = b'$.
(See  \cite[Subsection 6.2]{bidisk} for more detailed proof of more precise results in a related situation.)
The proof of Theorem \ref{theorem51} is now complete.
\end{proof}
\begin{rem}\label{limitremark}
\begin{enumerate}
\item
The argument of the last part of the proof of Theorem \ref{theorem51} shows that
to prove the inequality (the homology version of Arnold's conjecture)
\begin{equation}\label{arnold}
\#\frak P(H) \ge \sum \text{\rm rank} H_k(X;Q)
\end{equation}
for periodic Hamiltonian system with non-degenerate closed orbit,
we do not need to use projective limit.
We can use $HF(X,H;\Lambda_0/T^E\Lambda_0)$ for sufficiently large
but fixed $E$. (Such $E$ depends on  $H$ and $H'$.)
\item
During the above proof of the isomorphism $H(X,H;\Lambda) \cong H(X;\Lambda)$ we did not construct
an isomorphism
among them but showed only the coincidence of their ranks. Actually we can construct the following
diagram
\begin{equation}\label{Diag}
\CD
CF(X,H';\Lambda_0/T^{E'}\Lambda_0)
@>{\Phi_{E'}}>>CF(X,H;\Lambda_0/T^{E'}\Lambda_0)
\\
 @V{}VV  @V VV \\
CF(X,H';\Lambda_0/T^{E}\Lambda_0)
@>{\Phi_{E}}>>CF(X,H;\Lambda_0/T^{E}\Lambda_0)\endCD
\end{equation}
for $E<E'$. Here the vertical arrows are composition of reduction modulo $T^E$ and chain homotopy equivalence.
(Note this map is not reduction modulo $T^{E}$. In fact the two chain complexes that are
the target and the domain of the vertical arrows, are constructed by using different Kuranishi structures
and multisections.)
We can prove that Diagram \ref{Diag} commutes up to chain homotopy.
Then using a lemma similar to Lemma \ref{homotopyextension} we can extend $\Phi_E$ to a
chain map $CF(X,H';\Lambda_0) \to CF(X,H;\Lambda_0)$.
We then can prove that it is a chain homotopy equivalence.
This argument is a baby version of one developed in \cite[Section 7.2]{fooo:book1}.
\end{enumerate}
\end{rem}
\begin{rem}\label{lagrangeproof}
As we mentioned at the beginning of Part \ref{S1equivariant}, there is an alternative (third) proof of (\ref{arnold})
which does {\it not} use $S^1$ equivariant Kuranishi structure, and which works
for an arbitrary compact symplectic manifold $X$.
Let $H$ be a time dependent Hamiltonian whose 1 periodic orbits are all non-degenerate.
Let $\varphi : X \to X$ be the
symplectic diffeomorphism that is time one map of the time dependent Hamiltonian vector field associated to $H$.
We consider a symplectic manifold $(X\times X,\omega \oplus -\omega)$.
The graph
$$
L(\varphi) = \{(x,\varphi(x)) \mid x \in X\}
$$
is a Lagrangian submanifold in $X\times X$.
Since the inclusion induces an injective homomorphism $H(L(\varphi)) \cong H(X) \to H(X\times X)$ in the homology
groups, \cite[Theorem C, Theorem 3.8.41]{fooo:book1} implies that the Lagrangian Floer cohomology between $L(\varphi)$ with itself
is defined, (after an appropriate bulk deformation).
Again since $L(\varphi) \to X\times X$ induces an injective homomorphism in the homology,
the spectral sequence in \cite[Theorem D (D.3)]{fooo:book1} degenerates at $E_2$ stage
and implies
$$
HF((L(\varphi),b,\frak b),(L(\varphi),b,\frak b);\Lambda) \cong H(L(\varphi);\Lambda) = H(X;\Lambda).
$$
(Here $b$ is an appropriate bounding cochain and $\frak b$ is an appropriate bulk class.)
Since $L(\varphi)$ is Hamiltonian isotopic to the diagonal $X$, \cite[Theorem 4.1.5]{fooo:book1}
implies
$$
HF((L(\varphi),b,\frak b),(L(\varphi),b,\frak b);\Lambda) \cong HF((L(\varphi),b,\frak b),(X,b',\frak b);\Lambda).
$$
Note $L(\varphi) \cap X \cong \frak P(H)$ and the intersection is transversal. (This is a consequence of
nondegeneracy of the periodic orbits.) Therefore
the rank of the Floer cohomology $HF((L(\varphi),b,\frak b),(X,b',\frak b);\Lambda)$ is not greater than
the order of $\frak P(H)$. The formula (\ref{arnold})  follows.
\par
Note in the above proof we use injectivity of $H(L(\varphi)) \to H(X\times X)$ to show that
$L(\varphi)$ (which is Hamiltonian isotopic to the diagonal) is unobstructed (after bulk deformation).
Alternatively we can use the involution $X \times X \to X\times X$, $(x,y) \mapsto (y,x)$ to
prove unobstructedness of the diagonal. (See \cite{fooo:inv}.)
\end{rem}
\par\newpage
\part{An origin of this article and some related matters}
\label{origin}

\section{Background of this article}

This article is more than 200 pages long and claims
{\it no} new results. Since this is unusual
for the article of this length posted in the e-print server,
we explain why such an article is written in this part.
\par
This article concerns the
foundation of the virtual fundamental chain/cycle technique, especially
of the one represented by the Kuranishi structure introduced by the first and
the fourth named authors in \cite{FuOn99I, FOn} and slightly rectified by the present authors
in \cite[Appendix A.1]{fooo:book1}.
There are various articles in the electronic network (including the e-print server) that
express negative opinions on the soundness of the
foundation of the virtual fundamental chain/cycle technique, not just
complaining about lack of details of the proofs.
We also have information that several people have
mentioned negative opinion on the soundness of the
foundation of the virtual fundamental chain/cycle technique
in various occasions such as lectures or talks in
conferences or workshops.
\par
However the authors of the present article
were not directly asked or pointed out explicit gaps or
objections in our writing before such negative opinion was expressed in
the articles for public display. In addition there have been presented
some lectures or talks in which such negative opinions
on virtual fundamental chain/cycle technique with the authors
of the present article not being present.
\par
We point out that the foundation of the virtual fundamental chain/cycle technique
was established in published and refereed journal papers
(\cite{FOn,LiTi98,LiuTi98,Rua99,Si}) by various authors. Most of such articles were
written in the year 1996 that is 16 years ago.
Various articles (\cite{ChiMo,CMRS,operad,fooo:book1,HWZ,joyce,joyce2,Liu,LiuTi98,LuT,McDuff07,Sie96}) on the foundation
of the virtual fundamental chain/cycle
have also been written and/or published.
\par
During these 16 years, virtual fundamental chain/cycle technique
have been applied for numerous purposes of different nature and by
many authors   successfully. 
Moreover in the course of generating these applications, all the
technical incompleteness or inconsistency in its foundation
have been exposed and corrected through the subsequent underpinning and
systematic usages of virtual fundamental chain/cycle techniques
in many different ways.
\par
Nevertheless the fact that various negative opinions
on the soundness of the foundation of the virtual fundamental chain/cycle technique
spread out in public has caused significant reservation
for various mathematicians, especially younger generation of mathematicians,
to use this technique for their purposes.
We are afraid that this will cause certain delay of the
development of the relevant mathematics as a whole.
\par
On the other hand, when these negative opinion on the foundation of the virtual fundamental
chain/cycle technique were expressed, most of them were not written in the way that
explicitly specify and/or pin down the point of their concern.
This fact and the fact that we were not directly asked questions or objections,
\footnote{
When the book \cite{fooo:book1}
was published, we corrected or supplied more detail
on Kuranishi structure about all the points we were directly told
or we found ourselves.
After \cite{fooo:book1} was published, we have not directly heard of
problems or objections about our writings until March 2012.
There are a few exceptions:
(1) A few people asked us to write more
about analytic issue (more than
we wrote in  \cite[Section A1.4]{fooo:book1}).
We know indirectly that
such demands exist among other people too.
(2) The first named author had some e-mail discussions with
D. Joyce.
(3) D. McDuff also asked some questions to us.}
made it very difficult for us to try to eliminate
such obstructions to the applications of virtual fundamental chain/cycle technique.
\par
Finally, there was set up the occasion of a group of invited mathematicians
discussing the soundness of the foundation of the virtual fundamental chain/cycle
technique in the google group named `Kuranishi' (its administrator is H. Hofer) through
which questions asked directly to the present authors.
Especially, K. Wehrheim sent us a list of questions on March 13, 2012
which she thought were problematic in the existing literature.
We took this opportunity and did our best to change the
above mentioned situation so that other mathematicians can also
freely use virtual fundamental chain/cycle technique.
\par
For this purpose, we temporarily halted most of our other on-going joint projects
and have been concentrated in preparing answers to Wehrheim's questions
in the google group `Kuranishi' in great detail as complete as possible and also in
replying to all the questions asked further there.
The pdf files \cite{Fu1}, \cite{FOn2}, \cite{fooo:ans3}, \cite{fooo:ans34} \cite{fooo:ans5}
were uploaded to the google group `Kuranishi' for this purpose during this discussion.
They form the main parts of this article (after minor modification
\footnote{There are certain corrections especially to \cite{FOn2},
which becomes Part \ref{Part2} of this article. We mentioned them
as Remarks \ref{rem520}, \ref{thankmac}.}).
\par
After we had uploaded those files that we thought answered  all the questions
raised by Wehrheim, an article \cite{MW1} of McDuff and Wehrheim was posted in the arXiv.
The article contains some objections or negative comments about the soundness of the
foundation of the virtual fundamental chain/cycle technique, especially of the
Kuranishi structure laid out in \cite{FOn}, \cite{fooo:book1}. Unfortunately,
the article does not even mention the presence of
our replies \cite{Fu1}, \cite{FOn2}, \cite{fooo:ans3}, \cite{fooo:ans34} \cite{fooo:ans5}
to those objections and criticisms.
As far as we are concerned, those objections had been already replied and confuted
in our files \cite{Fu1}, \cite{FOn2}, \cite{fooo:ans3}, \cite{fooo:ans34} \cite{fooo:ans5}.
We will explain in Section \ref{HowMWiswrong} more explicitly on these points.

\section{Our summary of the discussion at google group
`Kuranishi'}

In this section we present a summary of our discussion in the
google group
`Kuranishi' {\it from our point of view}.
There are several discussions at the beginning of this google group,
but we skip them since we were not directly involved in that
discussion.
The discussion we were directly involved in are
related to the questions by K. Wehrheim.
On March 13 2012, she sent a series of questions about
Kuranishi structures.
The first of them is:
\par\medskip

{\bf 1.)} Please clarify, with all details, the definition of a Kuranishi structure. And could you confirm that a special case of your work proves the following?
\begin{enumerate}
\item
The Gromov-Witten moduli space $\mathcal M_1(J,A)$ of $J$-holomorphic curves of genus 0, fixed homology class $A$, with 1 marked point has a Kuranishi structure.
\item
For any compact space $X$ with Kuranishi structure and continuous map $f:X\to M$ to a manifold $M$ (which suitably extends to the Kuranishi structure), there is a well defined $f_*[X]^{\rm vir}\in H_*(M)$.
\end{enumerate}
\par\medskip
We replied in \cite{Fu1} as follows:
\par\medskip
{\bf Question 1}
(1) Yes.
(2) We need to assume $f$ to be strongly continuous and Kuranishi structure has
tangent bundle\footnote{We need to take the version of \cite{fooo:book1} not
of \cite{FOn} for the
definition of the existence of tangent bundle.} and orientation.
\par\medskip
This question is rather formal so there is nothing more to explain here.
The next two questions are:
\par\medskip
{\bf 2.)}  The following seeks to clarify certain parts in the definition of Kuranishi structures and the construction of a cycle.
\begin{enumerate}
\item
What is the precise definition of a germ of coordinate change?
\item
What is the precise compatibility condition for this germ with respect to different choices of representatives of the germs of Kuranishi charts?
\item
What is the precise meaning of the cocycle condition?
\item
What is the precise definition of a good coordinate system?
\item
How is it constructed from a given Kuranishi structure?
\item
Why does this construction satisfy the cocycle condition?
\end{enumerate}
\par
{\bf 3.)}  Let $X$ be a compact space with Kuranishi structure and good
coordinate system. Suppose that in each chart the isotropy group $\Gamma_p=\{{\rm id}\}$
is trivial and $s^\nu_p:U_p\to E_p$ is a transverse section. What further conditions
on the $s^\nu_p$ do you need (and how to you achieve them) in order to ensure that the
perturbed zero set $X^\nu=\cup_p (s^\nu_p)^{-1}(0) / \sim$ carries a global triangulation, in particular
\begin{enumerate}
\item $X^\nu$ is compact,
\item $X^\nu$ is Hausdorff,
\item $X^\nu$ is closed, i.e.\ if $X^\nu=\bigcup_n \Delta_n$ is a triangulation
then $\sum_n f(\partial \Delta_n) = \emptyset$.
\end{enumerate}

Some part of Question 2 concerns the notion of
germ of Kuranishi neighborhood.
We explain this issue in Subsection \ref{gernkuranishi}.
As we explain there we do not use  the notion of
germ of Kuranishi neighborhood
in \cite{fooo:book1}.
Definition \ref{Definition A1.5} in this article
is the same as \cite[Definition A1.5]{fooo:book1}.
So at the time of this question asked, the point related to the notion of `germ of Kuranishi coordinate'
had been already corrected in \cite{fooo:book1}.
\par
The most important statement in our answer
\cite{Fu1} to Question 3 is:
\par\medskip
No further condition on $s^{\nu}_p$ is necessary if it is close enough  to original Kuranishi map
and if we shrink $U_p$'s during the construction.
\par\medskip
In \cite{FOn2}, we explained construction of the good coordinate system and
virtual fundamental chains based on the definition of \cite{fooo:book1} and
confirmed the above statement in detail.
It provides an answer to
Questions 2) and 3).
(The outline of that construction is
given in \cite{Fu1}.)
A few typos were pointed out and were immediately corrected.
Then the discussion at the google group on
this point stopped for a while.
\footnote{We remark that Hausdorffness was a point mentioned
by several people in the talks by McDuff and Wehrheim
at Institute for advanced study in March and April 2012.}
\par
In the mean time, we continued posting other parts \cite{fooo:ans3}, \cite{fooo:ans34},
\cite{fooo:ans5} of our answer.
After we had finished posting the last answer file \cite{fooo:ans5},
which we thought had answered all the questions raised by K. Wehrheim,
a version of the manuscript by McDuff and Wehrheim was posted to the google group `Kuranishi',
which soon appeared in the e-print arXiv as the article
\cite{MW1}.
Then discussions on Questions 2), 3) were restarted in the google group `Kuranishi'.
McDuff asked several questions and then pointed out an error in \cite{FOn2}.
We acknowledged it was an error and at the same time posted its correction.
This correction is used in Section \ref{sec:existenceofGCS} of this article.
We acknowledge McDuff for pointing out this error in Remark \ref{thankmac}.
In the presence of Example 6.1.14 given in \cite{MW1}, we
wrote the proof of Hausdorffness and compactness, which appeared
in Sections 2 and 4 of \cite{FOn2}, more carefully.  It then became Sections \ref{defgoodcoordsec},
\ref{sec:existenceofGCS} and \ref{gentoplem} of this article.
Thus Part \ref{Part2} of present article answers
Questions 2 and 3 in great detail.
Namely
the answer to
Question 2 (1) is
Definition \ref{Definition A1.3},
the answer to
Question 2 (2) is
explained in Subsection \ref{gernkuranishi},
the answer to Question 2 (3) is
Definition \ref{Definition A1.5},
the answer to Question 2 (4) is
Definition \ref{goodcoordinatesystem},
the answer to Question 2 (5) (6)
are given in Section \ref{sec:existenceofGCS}.
The answer to Question 3 (1) is
Lemma \ref{cuttedmodulilem},
the answer to Question 3 (2) is
Corollary \ref{corollary521},
and the answer to Question 3 (3)
is Lemma \ref{cycleproperties}.

As we mentioned in the introduction, according to our opinion
the points appearing in Questions 2 and 3 are of technical nature.
\par\medskip
The next question was:
\par\medskip
{\bf 4.)} For the Gromov-Witten moduli space $\mathcal M_1(J,A)$ of $J$-holomorphic
curves of genus 0 with 1 marked point, suppose that $A\in H_2(M)$ is primitive so
that $\mathcal M_1(J,A)$ contains no nodal or multiply covered curves.
\begin{enumerate}
\item
Given two Kuranishi charts $(U_p, E_p, \Gamma_p=\{{\rm id}\}, \ldots)$ and
$(U_q, E_q, \Gamma_q=\{{\rm id}\},\ldots)$ with overlap at $[r]\in\mathcal M_1(J,A)$,
how exactly is a sum chart $(U_r,E_r,\ldots)$ with $E_r \simeq E_p\times E_q$ constructed?
\item
How are the embeddings $U_p \supset U_{pr} \hookrightarrow U_r$ and $U_q \supset U_{qr}
\hookrightarrow U_r$ constructed?
\item
How is the cocycle condition proven for triples of such embeddings?
\end{enumerate}
This question addresses only a very special case in the standard practice of the common researches in the
field.\footnote{In fact
we do not need to use virtual fundamental cycle in the case
appearing in this question. Since the pseudo-holomorphic curve
appearing in such a moduli space is automatically somewhere injective,
the transversality can be achieved by taking generic $J$.}
It appears to us that the only point to mention is about how we
associate the obstruction space $E(u')$ to an unknown map $u'$ (equipped with
various other data).
Namely the Kuranishi neighborhood is the set of the solutions
of the equation
$
\overline\partial u' \equiv 0 \mod E(u').
$
\par
The answer which was in the first named author's google post \cite{Fu1}
was as follows. (We slightly change the notations below so that it is
consistent with the argument in Part \ref{generalcase} of this article.)
\par\medskip
We cover the given moduli space by a finite number of
sufficiently small closed sets $\frak M_{\frak p_c}$ of the mapping space
each of which is `centered at' a point $\frak p_c \in \frak M_{\frak p_c}$ that is represented
by a stable map $((\Sigma_c,\vec z_c),u_c)$.
($\Sigma_c$ is a (bordered) Riemann surface and $\vec z_c = (z_{c,1},
\dots,z_{c,m})$ are marked points. (Interior or boundary marked points.)
$u_c : \Sigma_c \to X$ is a pseudo-holomorphic map.)
We fix a subspace $E_c$ of $\Gamma(\Sigma_c;u_c^*TX\otimes \Lambda^{0,1})$
as in (12.7) in page  979 of \cite{FOn}.
For $\frak p = ((\Sigma_{\frak p},\vec z_{\frak p}),u_{\frak p})$
we collect $E_c$ for all $c$ with ${\frak p} \in \frak M_{\frak p_c}$
and the sum of them is $E_{\frak p}$.
The Kuranishi neighborhood of ${\frak p}$ is the set of solutions of

\begin{equation}\label{definingeq}
\overline{\partial} u'  \equiv 0 \mod E_{\frak p}.
\end{equation}

We will discuss the way how we identify $E_c$ to a
subset of  $\Gamma(\Sigma';(u')^*TX\otimes \Lambda^{0,1})$
in case $((\Sigma',\vec z'),u')$ is close to $((\Sigma_{\frak p},\vec z_{\frak p}),u_{\frak p})$.

When we fix $E_c$ we also fix
finitely many additional marked points $\vec w_{c}
= (w_{cj})$ where $w_{cj} \in \Sigma_c$, $j=1,\dots,k_c$
at the same time
and take transversals $\mathcal D_{cj}$ as in  \cite[Appendix]{FOn}.
We take sufficiently many marked points so that after adding those marked points
$(\Sigma_c,\vec z_c\cup \vec w_{c})$ becomes stable.

We consider $((\Sigma',\vec z'),u')$. For each $c$ we
add marked points $\vec w'_c = (z'_{cj})$, $w'_{cj}
\in \Sigma$, $j=1,\dots,k_i$ to $(\Sigma',\vec z')$
so that
$w'_{cj}$ is on   $\mathcal D_{cj}$ .
We add them to obtain $(\Sigma',\vec z'\cup \vec w'_c)$ that
becomes stable, for each $c$.
We require that it is close to $(\Sigma_c,\vec z_c \cup \vec w_{c})$
in Deligne-Mumford moduli space (or its bordered version).
Then we obtain a diffeomorphism
(outside the neck region) from $\Sigma_c$ to  $\Sigma'$
which sends $\vec z_c\cup \vec w_{c}$ to $\vec z' \cup \vec w'_c$,
preserving the enumeration. (See \cite[Lemma 13.18 and Appendix]{FOn}.)
Using this diffeomorphism and
(complex linear part of) the parallel transport on $X$
(the symplectic manifold) with respect to the Levi-Civita connection
along the minimal geodesic joining $u'(w'_{cj})$ with $u_c(w_{cj})$
(where $w_{cj} \in \Sigma_c$ is identified with $w'_{cj}
\in \Sigma'$ by the above mentioned diffeomorphism),
we send $E_c$ to a set of the sections of $(u')^*TX\otimes \Lambda^{0,1}$ on $\Sigma'$.
We do it for each of $c$.
(In other words the stabilization we use {\it depends} on $c$.)
Thus each of $E_c$ is identified with a subspace of
$\Gamma(\Sigma';(u')^*TX\otimes \Lambda^{0,1})$.
We take its sum and that is $E_{\frak p}$ at $(\Sigma',\vec z',(\vec w'_c))$.
\par
We have thus made sense out of (\ref{definingeq}).
\par
An important point here is that the subspace $E_c
\subset \Gamma(\Sigma';(u')^*TM\otimes \Lambda^{0,1})$ at $(\Sigma',\vec z',(\vec w'_c))$
is {\it independent} of ${\frak p}$ as far as $(\Sigma',\vec z')$ is close to ${\frak p}$.
The data we use for stabilization is chosen on ${\frak p}_c$
(not on ${\frak p}$) once and for all.
This is essential for the cocycle condition to hold\footnote{Since Equation  (\ref{definingeq})
makes sense in the way independent of $\frak p$ it seems possible to
simply take the union of its solution space to obtain
some Hausdorff metrizable space. That can play a role as the metric space
that contains all the Kuranishi neighborhoods. However
we insist that we should {\it not} build the general theory of Kuranishi structure
under the assumption of the existence of such an ambient space, since it will spoil the
flexibility of the definition of the general story we have.}.

(2),(3).
Once (1) is understood the coordinate change $\underline\phi_{{\frak q}{\frak p}}$ is just a map which send an en element
$((\Sigma',\vec z'),u')$ to the same element.
So the cocycle condition is fairly obvious. { [End of an answer in \cite{Fu1}]}
\par\medskip
We were then asked further details to provide and responded to their requests
by the posts \cite{fooo:ans3} and \cite{fooo:ans34}.
Those answers discuss much more general case than the case asked
in Question 4. We provided this general answer
because there was also the demand for providing more details of the gluing analysis.
\par
The case asked in Question 4 is of course contained
in \cite{fooo:ans3} and \cite{fooo:ans34} as a special case.
More specifically, our answers to the questions in Question 4 are given as follows:
\par
\begin{enumerate}
\item \cite[Prosition 2.125]{fooo:ans34} = Proposition \ref{chartprop}.
(It is based on \cite[Definition 2.60]{fooo:ans34} + \cite[Definition 2.63]{fooo:ans34}
+ \cite[Definition 2.119]{fooo:ans34}
=  Definitions \ref{defEc}, \ref{defthickened}, \ref{defVVVVV}.)

\item The proof is completed in \cite[Subsection 2.10 (2.385)]{fooo:ans34}
= Section \ref{kstructure} (\ref{eq23777}),
based on technical lemma \cite[Lemma 2.163]{fooo:ans34} =Lemma \ref{lem2143}.
The main part of the construction is
\cite[Propositions 2.131 and 2.152]{fooo:ans34} = Propositions \ref{prop2117} and \ref{prop21333}.

\item  The proof is completed in \cite[Subsection 2.10 (page 119)]{fooo:ans34}
= Section \ref{kstructure}.
The main part of the proof is
\cite[Lemma 2.145 and Proposition 2.158]{fooo:ans34} = Lemma \ref{lem125} and Proposition \ref{compaticoochamain}.
\end{enumerate}

There was the following additional question asked by Wehrheim on March 23.
\par\medskip
Q4: ``Assuming that I understand your construction of a single Kuranishi
chart in [FOn],\footnote{[FOn]
is \cite{FOn} in the reference of this article.}
I would like to see all the details of constructing
this simplest sum chart (which in [FOn] are rather scattered, as you
said). In particular, I would like to see the very explicit Fredholm
setup for the set of solutions of (0.5) ... i.e. what is the Banach
manifold? What is the Fredholm operator? Why is it smooth / transverse?
Or ... if you work with gluing maps and infinitesimal local slices ...
maybe you could give the construction in a series of lemmas without
proofs. Then I can ask about specific proofs.
In fact, in that case I could ask directly to get a proof of
injectivity / surjectivity / continuity of inverse  for the map $\psi$
from the "complement of tangent space to automorphism group in the
domain of gluing map" to the moduli space. I couldn't see much of an
explanation in [FOn] in this case without Deligne-Mumford parameters.''

\par\medskip
Our reply to it (which was sent at the same time as \cite{fooo:ans34})
was as follows:
\par\medskip
The Banach manifold we use to define smooth (or $C^m$) structure
is one of $C^m$ maps from the compact subset of Riemann surface to
$X$. (See the proof of \cite[Lemma 2.133]{fooo:ans34} = Lemma \ref{2120lem}
and \cite[Section 3.2]{fooo:ans34}
= Section \ref{toCinfty}.)
\footnote{More precisely we take a product with an
appropriate Deligne-Mumford moduli space (that includes the gluing parameter).}
We may use also the space of $L^2_m$ maps.
Because of elliptic regularity there is no difference which one we use.

The Fredholm operator appearing in the gluing construction is \cite[(2.247)]{fooo:ans34}
= (\ref{lineeqstep0vvv}).
This is not exactly the linearization operator of the equation.
We still have extra coordinate while solving equation.
Correct moduli space is obtained by cutting the solution space
down by putting some constraints.

The injectivity / surjectivity / continuity of inverse of the map $\psi$
is proved in \cite[Proposition 2.102]{fooo:ans34} = Proposition \ref{charthomeo}.
(More precisely, injectivity in \cite[Lemma 2.106]{fooo:ans34}=Lemma  \ref{injectivitypp},
surjectivity in  \cite[Lemma  2.105]{fooo:ans34}
= Lemma \ref{setisopen} and continuity in  \cite[Lemma 2.108]{fooo:ans34} = Lemma \ref{injectivitypp3}
are proved respectively.)
Its proof uses injectivity / surjectivity of the thickened moduli space
that is a part of \cite[Theorem 2.70]{fooo:ans34} = Theorem \ref{gluethm3}.
(The proof of that injectivity / surjectivity
is given in \cite[Subsection 1.5]{fooo:ans34} = Section \ref{surjinj}.)

Automorphism group of the domain is killed by adding
marked points to the source.
\par\medskip
After the discussion at the google group resumed,
we were asked to explain the following two points by
Wehrheim and McDuff (There were a few other points. But the two points below were
the major points of their concern of that time, we think).
When we set up the equation  (\ref{definingeq}), we use the additional
marked points $w'_{cj}$ on the source $\Sigma'$ of the map $u'$ there.
Therefore the thickened moduli space also involves the
added marked points and has dimension bigger than the given virtual
dimension by $2 \times \#\{w'_{cj}\}$. To kill off this extra dimension, we put the constraint
\begin{equation}\label{transversalconstr}
u'(w'_{cj}) \in \mathcal D_{cj}
\end{equation}
on the choice of $w'_{cj}$.
\par\medskip
\begin{enumerate}
\item[(A)]  The question was
whether this constraint equation is transversal or not.
\end{enumerate}
\par\medskip
We replied that the transversality of this equation had been
proved in the course of the proof of Lemma \ref{transstratasmf} etc.
\par
The right hand side $E_{\frak p}$ of the equation  (\ref{definingeq}) actually depends
on $u'$. Let us write it as $E_{\frak p}(u')$.
\begin{enumerate}
\item[(B)]
The question was
whether this forms the smooth vector bundle when we vary the associated
parameter space of $u'$'s and etc. (in an appropriate function space).
\end{enumerate}
\par\medskip
Our answer (Aug. 2) was as follows:
\par\medskip

This point is remarked in \cite[footnote 17 p. 77]{fooo:ans34} =
footnote 23 in Section \ref{glueing}
and also in \cite[Remark 1.39]{fooo:ans3} = Remark \ref{remark127}.
Also during the proof of gluing analysis,
we need to take (second) derivative of the projection to the
obstruction bundle. The projection to the
obstruction bundle is calculated in Formula \cite[(1.39)]{fooo:ans3} =(\ref{form118}).
We take its derivative in \cite[(1.40)]{fooo:ans3} = (\ref{DEidef}).
Existence of this (and higher) derivative is obvious
from the explicit formula \cite[(1.39)]{fooo:ans3} =(\ref{form118}).

In the proof of \cite[Lemma 1.22]{fooo:ans3} = Lemma \ref{mainestimatestep13}, which
is the key estimate for the Newton's method to work in gluing analysis,
the estimate of the second derivative of the projection
to the obstruction bundle becomes necessary.
More explicitly, it appears in \cite[(1.58)]{fooo:ans3} = (\ref{155ff}).
\par
However since we were specifically asked to provide more details we
replied to it. We reproduce our reply
in Subsection \ref{smoothness}.
\par\medskip

The fifth question of Wehrheim was:
\par\medskip
{\bf 5.)} How is equality of Floer and Morse differential for the Arnold conjecture proven?
\begin{enumerate}
\item
Is there an abstract construction along the following lines: Given a compact topological
space $X$ with continuous, proper, free $S^1$-action, and a Kuranishi structure for $X/S^1$
of virtual dimension $-1$, there is a Kuranishi structure for $X$ with $[X]^{\rm vir}=0$.
\item
How would such an abstract construction proceed?
\item
Let $X$ be a space of Hamiltonian Floer trajectories between critical points of index difference $1$,
in which breaking occurs (due to lack of transversality).
How is a Kuranishi structure for $X/S^1$ constructed?
\item
If the Floer differential is constructed by these means, why is it chain homotopy equivalent to
the Floer differential for a non-autonomous Hamiltonian?
\end{enumerate}

The reply in \cite{Fu1} was as follows:
\par\medskip

(1)(2) I do not think it is possible in completely abstract setting.
At least I do not know how to do it.
In a geometric setting such as the one appearing in page 1036 \cite{FOn},
Kuranishi structure is obtained by specifying the choice of
the obstruction space $E_p$ for each $p$.
We can take $E_p$ in an $S^1$ equivariant way so
the Kuranishi structure on the quotient $X/S^1$ is obtained.
And it is a quotient of one on $X$.
$S^1$ equivariant multisection can be constructed in an abstract setting
so if the quotient has virtual dimension $-1$ the zero set is empty.

(3) We can take a direct sum of the obstruction bundles, the support of which is
disjoint from the points where two trajectories are glued.
In the situation of (1) the obstruction bundle is $S^1\times S^1$ equivariant.
The symmetry is compatible with the diagonal $S^1$ action nearby.

(4) It is  \cite[Theorem 20.5]{FOn}.
\par\medskip
We then were requested to explain more detail.
We  sent \cite{fooo:ans5} to the google group,
which contains the required detail. It becomes Part \ref{S1equivariant} of this article.
\par
There were no further questions or objections to the points concerning
Question 5)
so far in the google group.
\par\medskip
At the time of writing this article (Sep. 10) all the
questions or objections asked in the google group `Kuranishi'
were answered, supplemented or confuted by us.

\section{Explanation of various specific points in the theory}

In this section we mention some of the points which we were asked
on the foundation of the virtual fundamental chain/cycle technique
and explain why they do not affect the rigor
of the virtual fundamental chain/cycle technique.
We have already discussed most of them in the main body of this article.
The discussion of this section mentions
other method to the problem also and
sometimes discuss some special case to clarify the idea.

\subsection{A note on the germ of the Kuranishi structure}
\label{gernkuranishi}

In \cite{FOn} the notion of `germ of Kuranishi neighborhood' is
defined as follows.
Let $(V_p, E_p, \Gamma_p, \psi_p, s_p)$,
$(V'_p, E'_p, \Gamma'_p, \psi'_p, s'_p)$
 be a Kuranishi neighborhoods centered at
$p \in X$ be as in Definition \ref{Definition A1.1}.
We say that they are equivalent
if there is a third Kuranishi neighborhood
$(V''_p, E''_p, \Gamma''_p, \psi''_p, s''_p)$
together with:
\begin{enumerate}
\item
Isomorphisms
$$
\Gamma_p \cong \Gamma''_p \cong \Gamma'_p.
$$
\item
Equivariant open embeddings
$$
V''_p \to V_p, \qquad
V''_p \to V'_p.
$$
\item
Fiberwise isomorphism bundle maps
$$
V''_p \times E''_p \to V_p \times E_p,
\quad
V''_p \times E''_p \to V'_p \times E'_p.
$$
which cover the open embeddings in (2)
and are equivariant.
\item
They are compatible with $\psi_p$, $s_p$
in an obvious sense.
\end{enumerate}
The equivariant class is called a germ of Kuranishi neighborhood.
It was used for a similar germ version of coordinate change.
Together with compatibility condition which is a `germ version' of
Definition \ref{Definition A1.5} (2), the definition of Kuranishi structure was given in \cite{FOn}.
\par
The main trouble of this definition is as follows.\footnote{This point
was mentioned by the first named author on 19th March 2012 in his
post to the google group Kuranishi.}
Using the open embedding of (2) above
we obtain a diffeomorphism between neighborhoods of $o_p$
in $V_p$ and in $V'_p$.
However such a diffeomorphism is {\it not} unique.
The germ of such diffeomorphisms at $o_p$ is not unique either.
(The compatibility with $\psi_p$ implies that its
restriction to the zero set of the Kuranishi map
$s_p$ is unique.)
\par
We then consider the coordinate change.
Suppose we have a coordinate change $(\hat\phi_{pq},\phi_{pq},h_{pq})$
from a Kuranishi neighborhood
$(V_q, E_q, \Gamma_q, \psi_q, s_q)$
centered at $q$ to  a Kuranishi neighborhood
$(V_p, E_p, \Gamma_p, \psi_p, s_p)$
centered at $p$.
(Here $q \in \text{\rm Im}\psi_p$.)
\par
We next take another pair of Kuranishi neighborhoods
$(V'_q, E'_q, \Gamma'_q, \psi'_q, s'_q)$
and
$(V'_p, E'_p, \Gamma'_p, \psi'_p, s'_p)$
of $q$ and $p$ respectively.
\par
We assume that
$(V'_q, E'_q, \Gamma'_q, \psi'_q, s'_q)$
and
$(V'_p, E'_p, \Gamma'_p, \psi'_p, s'_p)$
are equivalent to
$(V_q, E_q, \Gamma_q, \psi_q, s_q)$
and
$(V_p, E_p, \Gamma_p, \psi_p, s_p)$,
respectively.
\par
Then if we choose the diffeomorphisms
such as item (2) above between
\begin{enumerate}
\item[(a)]
A neighborhood of $o_p$ in $V_p$ and
a neighborhood of $o_p$ in $V'_p$,
\item[(b)]
A neighborhood of $o_q$ in $V_q$ and
a neighborhood of $o_q$ in $V'_q$,
\end{enumerate}
we can use them together with
$(\hat\phi_{pq},\phi_{pq},h_{pq})$ to find a
coordinate change $(\hat\phi'_{pq},\phi'_{pq},h'_{pq})$ from
$(V'_q, E'_q, \Gamma'_q, \psi'_q, s'_q)$
to
$(V'_p, E'_p, \Gamma'_p, \psi'_p, s'_p)$.
\par
The problem lies in the fact that
this induced coordinate change $(\hat\phi'_{pq},\phi'_{pq},h'_{pq})$
{\it does} depend on the choice of
the diffeomorphisms (a)(b) above.
It does depend on it even if we
consider a small neighborhood of $o_q$.
\par
As a consequence, it is hard to state the
compatibility condition between two different coordinate changes
in the way independent of the choice of the representative of the
germs of the Kuranishi neighborhood.
\par
{So the meaning of the statement that coordinate changes are compatible
which uses the language of germs is ambiguous from the writing in \cite{FOn}.
This is indeed a mathematical error of \cite{FOn}.}
{\begin{rem}
At the time of our writing of \cite[Appendix]{fooo:book1} however, we became aware of the
danger of the notion of germs of Kuranishi neighborhoods and etc.
This was the reason why we rewrote the definition of
the Kuranishi structure therein so that it does not involve the notion of germs.
However our understanding of the above mentioned problem at that time was not as
complete as now. That is the reason why we did not mention it in \cite{fooo:book1}.
\end{rem}}
\par
A correct way of clarifying this point, which was adopted in \cite{fooo:book1} and in this article, is
as follows. We take and fix representatives of Kuranishi neighborhoods
for each $p$. We also fix a choice of coordinate changes between
the Kuranishi neighborhoods, which are representatives of coordinate changes
between the Kuranishi neighborhoods of each $p$ and $q \in \psi_p(s_p^{-1}(0))$.
Here we fix a coordinate change {\it not} its germ.
Then the compatibility between coordinate changes
has definite meaning without ambiguity.
\par
Note that
the representatives of Kuranishi neighborhoods and coordinate changes
were taken and fixed in the proof of the existence of the
good coordinate system in \cite{FOn},
as mentioned explicitly in \cite[(6.14.2), (6.19.2), (6.19.4)]{FOn}.
Also the Kuranishi neighborhood and coordinate changes between them
appearing in the good coordinate system are the representatives
but not the germs. (This point is emphasized in  \cite[Remark 6.2]{FOn}.)
By this reason the proof of existence of the
good coordinate system in \cite[Lemma 6.3]{FOn} is correct notwithstanding the error mentioned above,
if we take the definition of Kuranishi structure from \cite{fooo:book1}.
The proof of the existence of good coordinate system
which is Theorem \ref{goodcoordinateexists} is given in Section \ref{sec:existenceofGCS}
uses the same idea as the proof of \cite[Lemma 6.3]{FOn}, but contains more detail.
\par\medskip
Another issue was mentioned by
D. Joyce on March 19 in the google group `Kuranishi' and
also during the E-mail discussion with the first named author.
\par
Let us consider a germ of coordinate change from a Kuranishi neighborhood
of $q$ to one of $p$.
 The issue is that we can't regard $p$ and $q$ as fixed: since we can
always make $V_p$ smaller in the germ, for fixed $p$ and $q$ you can always
make $V_p$, $V_q$ smaller in their germs so that they do not overlap, and
there is nothing to transform.\footnote{This sentence is copied
from D. Joyce's post on March 19 to the google group Kuranishi.}
\par
It seems that we can go around this problem by carefully
choosing the way of defining the notion of the germ of
coordinate change etc.
However since we have already modified the definition and eliminate the
notion of germ from our definition because of the first point,
we do not discuss this second issue any more.
\par\medskip
We remark that in his theory \cite{joyce2},
Joyce pushed the sheaf-theoretic point of views to its limit
and, we believe, his approach gives more thorough and systematic
answers to this problem. On the other hand, our purpose here
is to provide the shortest rigorous way of constructing
virtual fundamental cycles/chains in the situation where we apply the moduli space
of pseudo-holomorphic curve etc. to symplectic geometry, mirror symmetry and etc..
The way taken in \cite{fooo:book1} and in this article
lies rather at the opposite end of \cite{joyce2} in this respect.
We restrict our attention to effective orbifolds and then restrict
maps between them to embeddings.
In this way, we can study orbifolds and maps between them
in a way as close as possible to the way we study manifolds and smooth maps between them.
In \cite{joyce2}  the most general case of
orbifolds and morphisms between them are included.
To handle such general case, the language of 2-category is systematically used
in \cite{joyce2}.
\footnote{As we mentioned before this seems to be the correct way
to study Kuranishi space, {\it as its own}, which is {\it not} our purpose here.}

Part \ref{Part2} of this article
contains the construction of virtual fundamental chain/cycle and
the proof of its basic property. The proof we provide
is, in our opinion, strictly more detailed than is required
in the standard research articles.
Nevertheless it is much shorter than
other articles discussing similar points using different approach.

\subsection{Construction of the Kuranishi structure on the
moduli space of pseudo-holomorphic curve}
\label{subsec342}

One of the objections to the virtual fundamental chain/cycle technique, which we (mainly indirectly) heard of
before 2009 and sometimes after that,
is about the construction (or existence) of the (smooth) Kuranishi
structure on the moduli space of pseudo-holomorphic curves in various situations.
\par
Most of this problem is of analytic nature.
The main point we heard of is about the gluing (or stretching the neck) construction
and the smoothness of the resulting coordinate changes between the Kuranishi neighborhoods obtained
by the gluing (or stretching the neck) construction.
\par
We mentioned this point in the introduction and have provided one analytical way of constructing
\emph{smooth} Kuranishi structure with corners in detail
in Parts \ref{secsimple} and \ref{generalcase} of the present article.
\par
Here we discuss this point again at the same time mentioning another more
geometric way of constructing such structure.
\par
There is a well-established technique of extracting the moduli space as a manifold with boundary
in certain circumstances. It was used by Donaldson in gauge theory
(in his first paper \cite{Don83} to show that 1 instanton moduli of ASD connections on 4 manifold $M$
with $b_2^+ = 0$ has $M$ as a boundary.) According to this method one takes certain parameter
(that is, the degree of concentration of the curvature for the
case of ASD equation and the parameter $T$ in the situation of Part \ref{secsimple}).
We consider the submanifold where that parameter $T$ is sufficiently large, say $T_0$.
We throw away everything of the part where $T > T_0$.
Then the slice $T=T_0$ becomes the boundary of the `moduli space' we obtain.
It was more detailed
in a book by Freed and Uhlenbeck \cite{freedUhlen} in the gauge theory case.
Abouzaid used this technique in his paper \cite{Abexotic} about exotic spheres in $T^*S^n$,
including the case of corners.
At least as far as the results in \cite{FOn} are concerned we can
use this technique and prove all the results in \cite{FOn} since we have only to study the moduli
space of virtual dimension 0 and 1 only.
In other words we can use something like Theorem \ref{gluethm1}
for a large fixed $T$, but does not need to estimate the $T$ derivative or
study the behavior of the moduli space at $T=\infty$.
\par
{The reason is as follows. Suppose that the moduli space in consideration
carries the Kuranishi structure of virtual dimension $1$ or $0$. Then when
we consider the corners of codimension $2$ or higher
the restriction of the moduli space to that corner has negative virtual dimension.
So after a generic multivalued perturbation the zero set on the corner becomes
empty. So all we need is to extend multivalued perturbation to its neighborhood.
We remark that $C^0$ extension is enough for this purpose.}
\par
For the case of moduli space of virtual dimension 1, after generic perturbation
we have isolated zeros of the perturbed Kuranishi map on the boundary.
So, for large $T_0$, Theorem \ref{gluethm1} or its analogue
implies that set of zeros on the `boundary $T=T_0$' has one-one
correspondence with that of the actual boundary ($T=\infty$).
Because of this,  we do not need to carefully examine what happens
in a neighborhood of the set $T=\infty$.
All we need is to extend this given perturbation at $T=T_0$ to the interior.
\par
We also remark that it is unnecessary to show differentiability
of the coordinate change at this boundary.
This is because the zero set at the boundary is isolated.
\par
This argument is good enough to establish all the
results in \cite{FOn}.
\par
As we mentioned explicitly in \cite[page 978 line 13]{FOn} our argument
there, in analytic points, is basically the same as in \cite{McSa94}.
(Let us remark however the proof of `surjectivity' that is written in
\cite[Section 14]{FOn} is slightly different from one in \cite{McSa94}.)
So the novelty of \cite{FOn} does {\it not} lie in the analytic point but
in its general strategy, that is
\begin{enumerate}
\item Define some general notion of `spaces' that contain
various moduli spaces of pseudo-holomorphic curves as examples and work out transversality
issue in that abstract setting.
\item  Use multivalued abstract perturbation, which we call
multisection.
\end{enumerate}
\par\medskip
When we go beyond that and prove results such as those
we had proved in \cite{fooo:book1},
we need to study the moduli spaces of higher virtual dimension
and study chain level intersection theory.
In that case we are not sure whether the above mentioned technique is enough.
(It may work. But we did not think
enough about it.)
This is not the way we had taken in \cite{fooo:book1}.
\par
Our method in \cite{fooo:book1} was to use exponential decay
estimate (\cite[Lemma A1.59]{fooo:book1}) and to use $s = 1/T$ as the
coordinate on the normal direction to the stratum to define
smooth coordinate of the Kuranishi structure.  {Here $T$ is the gluing
parameter arising from the given analytic coordinate $z$ centered at the
puncture associated to the marked point.  More specifically, we have
$T = -\ln |z|$ (and hence $s = e^{-1/|z|}$.)}.

We refer \cite[Subsection A1.4]{fooo:book1}  and
\cite[Subection 7.1.2]{fooo:book1} where this construction is written.
\par
In Parts \ref{secsimple} and  \ref{generalcase}, we provide more details of the way how to use
alternating method to construct smooth chart at infinity
following the argument in \cite[Subsection A1.4]{fooo:book1}.

\subsection{Comparison with the
method of \cite{FOn} Sections 12-15 and of \cite[Appendix]{FOn}}
\label{comparizon}
\par
During the construction of the Kuranishi neighborhood of each
point in the moduli space of pseudo-holomorphic curve,
we need to kill the automorphism of the source curve
and fix the parametrization of it.
\footnote{This subsection is mostly the
copy of our post to the google group Kuranishi
on August 3.}
Let us call this process normalization in this subsection.
\par
In \cite{FOn} we provided two different normalizations.
One is written in Sections 12-15 and the other is
in Appendix.
Here is some explanation of the difference
between two approaches. There are two steps for this normalization process.
\smallskip
\par
{\bf Step (1) } To fix a coordinate or parametrization of the source.

{\bf Step (2) } To kill the ambiguity of the parametrization in Step (1).
\smallskip
\par
The difference between two techniques mainly lies
in Step (2).

The techniques of the appendix is certainly {\it not}
our invention.
It was used by many people as we quoted
in  Remark \ref{rem161}.
We believe it was already a standard method
when we wrote it in 1996.

In both techniques, we need to fix a parametrization
of the source curve. (Step (1).)
In case we construct a Kuranishi chart locally,
this is not a serious matter because we fix a parametrization
of the source at the center of the chart and use it to fix a parametrization nearby.
Transferring obstruction bundles centered at some points to
another point (on which we want to define an obstruction bundle)
is more serious issue. This point is discussed in \cite[Section 15]{FOn}.

In both techniques we used
the method to stabilize the
curve by {\it adding  marked points}.
In \cite[Sections 12-15]{FOn}
it is written in page 989.
There we add marked points
to reduce the isotropy group to a finite group.
(It is the isomorphism in (13.19).)
In the appendix it is written in page 1047.
\par
The difference of the two techniques lies in the way how we kill the ambiguity.
(Step (2) above.)
Namely, there are more parameters than the expected dimension of
Kuranishi neighborhood.
\par
In \cite[Sections 12-15]{FOn} it is done
by {\it requiring the local minimality
of the function} `meandist' (15.5).
Namely we require the map $u'$ to be as close to $u$
by the given parametrization.
(Here $u$ is the map part of the object that is the center of the chart and
$u'$ is the map part of the object that is a general element of the Kuranishi
neighborhood).
This is a version of the technique called the {\it center of mass technique}
which was discovered by K. Grove and H. Karcher in Riemannian
geometry \cite{GrKa}.
Probably using
this technique in this situation was new and
was not so standard.
\par
On the other hand, in \cite[Appendix]{FOn}, the ambiguity is
killed by {\it putting codimension $2$ submanifolds
and requiring that the marked points to land on those submanifolds}.
(As we mentioned above this was already a standard technique at the stage of 1996).
\par
Both techniques work.
But later we mainly use the technique given in the appendix.
For example, in \cite[p424]{fooo:book1} we wrote as follows.
\par\bigskip
For the case where $(\Sigma, \vec{z})$ is unstable, we use Theorem 7.1.44 and proceed
as follows.  (The argument below is a copy of that in
Appendix [FuOn99II]\footnote{The reference [FuOn99II] in \cite{fooo:book1} is just \cite{FOn} in the current article.}.
There
is an alternative argument which will be similar to Section 15 [FuOn99II].)  We add some
interior marked points $\vec{z}^+$ so that $(\Sigma, \vec{z}, \vec{z}^+)$ becomes stable.
We may also assume the following
(7.1.48.1) Any point $z_i^+$ lie on a component where $w$ is non-trivial.
(7.1.48.2)  $w$ is an immersion at each point of $\vec{z}^+$.

By the same reason as in the appendix [FuOn99II], we can make this assumption without
loss of generality.  We choose $Q_i^{2n-2} \subset M$ (a submanifold of codimension 2) for
each $z_i^+ \in \vec{a}^+$ such that $Q_i^{2n-2}$ intersect with $w(\Sigma)$ transversally at
$w(z_i^+)$.  .....
\par\medskip
There is a similar sentence in page 566 of the year-2006 preprint version of our book.
\par
Let us add a few words about the use of `slice' in \cite[Section 12]{FOn}.
There we first consider the space of maps with the parameter of the
source fixed.
Then we obtain a family of solutions (of nonlinear Cauchy-Riemann equation modulo
obstruction bundle) parametrized by
a finite dimensional space. (that is $V^+_{\sigma}$ which appears in the three lines above from \cite[Theorem 12.9]{FOn}.)
This part is \cite[Proposition 12.23]{FOn}.

This space is not the correct Kuranishi neighborhood since it has
too many parameters. (The extra parameters correspond to the automorphism
group of the source that
is {\it finite dimensional}, though.)
We take a slice at this stage.
Namely we use $V'_{\sigma}$ in place of $V^+_{\sigma}$.
(This is  \cite[Lemma 12.24]{FOn} that comes after \cite[Proposition 12.23]{FOn}.)
So when we take a slice, the space in question is {\it already} a solution space
of the elliptic PDE and so consists of smooth maps.
The reparametrization etc. is obviously smooth.

From this point of view, the situation is essentially different
from gauge theory.
When we consider the set of solutions of ASD equation (without gauge fixing condition),
we will get some {\it infinite dimensional space} and
its element {\it may not be smooth} because of the lack of ellipticity.
Dividing it by the infinite dimensional gauge transformation group is indeed a
nontrivial analytic problem.
So usually people study the process of gauge fixing and
solving ASD equation at the same time.
Then analysis is certainly an issue at that stage.
\par\medskip
Note that the action of the group of diffeomorphisms is mentioned at the beginning of \cite[Section 15]{FOn}.
If one reads this part conscientiously, one finds that it occurs during the
heuristic explanation why some approach (especially the same approach as gauge theory) has a
trouble and is {\it not} taken in \cite[Section 15]{FOn}.
The infinite dimensional group of diffeomorphisms appears only here
in \cite{FOn} and so it does not appear
in the part where the actual proof is performed. It is written in \cite{FOn}, line 17- 13
from the bottom on page 999 :
\par\bigskip
\dots
one may probably be able to prove a slice theorem \dots.
However, because of
the trouble we mentioned above, we do not use this infinite dimensional space and work more directly without using infinite dimensional manifold.
\par\bigskip
Here it is written clearly that we do not use a slice theorem in
\cite[Sections 12-15]{FOn}.
\footnote{We are informed that several people discuss a problem of the foundation
of virtual fundamental chain/cycle technique based on the fact that
the action of the group of diffeomorphisms ${\rm Diff}(\Sigma)$ on $L^p_1(\Sigma,X)$ is not differentiable.
Sometimes it is said they can not take slice to this group action because of the
lack of differentiability. As we mentioned here we {\it never} take slice of this group action
in our approach via the Kuranishi structure.}

Since in this article we provided the detail using the technique
killing the ambiguity by transversals
and not by the center of mass, and to write both techniques in such a detail ($>$100 pages) is too much,
we do not think it is necessary to discuss this comparison much longer,
for our purpose. Our purpose is to clarify the soundness of the virtual fundamental cycle or chain technique
based on Kuranishi structure and multisection.

\subsection{Smoothness of family of obstruction spaces $E(u')$}
\label{smoothness}

In the google group Kuranishi, we were asked on the smoothness of
the obstruction bundle $E_{\frak p}(u')$ which appears in the right hand side of
(\ref{definingeq}) for example, by McDuff. (Here $u'$ runs in an appropriate $L^2_{m}$ space of maps.)
The discussion in the google group ends up with the agreement that this
statement (that is $u' \mapsto E_{\frak p}(u')$ is a smooth family of vector subspaces) is correct.

Here we reproduce the pdf files that we uploaded to the google group Kuranishi
on August 10 and 12. They explain the proof in detail
in the particular case of the moduli space  $\mathcal M_1(A)$ of holomorphic maps
$u : S^2 \to X$ and of homology class $A$ that is primitive.
(Namely there is no nonconstant holomorphic maps $u_1, u_2 : S^2 \to X$
such that $u_{1*}[S^2] + u_{2*}[S^2] = A$.)

Here $1$ in the suffix means that we consider one marked point (= $0$).
In other words, we identify two maps $u,u'$ if there exists a biholomorphic map
$v : S^2 \to S^2$ such that $u\circ v = u'$ and $v(0) =0$.
In this case there is neither a bubble nor a nontrivial automorphism.
(\emph{This is the case to which \cite{MW1} restrict themselves}.)

We first explain construction of the Kuranishi chart in
this special case. {It is our opinion that these two posts essentially
take care of all the cases in application the content of \cite{MW1} can handle.}
\par\medskip
{\bf [The post on August 10]}
\par
Suppose one has
$u_{c(1)} : S^2 \to X$ and
$u_{c(2)} : S^2 \to X$,
for which we take obstruction bundles $E_1$ and $E_2$.
\footnote{$E_1$ is a finite dimensional vector space of smooth sections
of $u_{c(1)}^*TX \otimes \Lambda^{01}$. We assume the support of its element is
away from singular or marked points.}
Also we take $D_{i,1}$, $D_{i,-1}$,  which are codimension 2 submanifolds of $X$
that are transversal to $u_{c(i)}$ at $1 \in S^2$ and $-1 \in S^2$, respectively.
(We assume $u_{c(i)}$ is an immersion at $\pm 1$.)

Let $u : S^2 \to X$
be a third map to which we want to transfer $E_1$ and $E_2$.
(In other words $u$ is a pseudo-holomorphic map which will be
the center of the Kuranishi chart we are constructing.)

As an assumption
$u g_1$ is $C^0$ close to $u_{c(1)}$
and
$u g_2$ is $C^0$ close to $u_{c(2)}$.

Here $g_1,g_2 \in Aut(S^2,0)$.
(We are studying $\mathcal M_1(A)$ and $0 \in S^1$ is the marked point.
So $g_i(0) = 0$.)

But $g_1$ and $g_2$ can be very far a way from each other.

We require $u(g_i(1)) \in D_{i,1}$, $u(g_i(-1)) \in D_{i,-1}$.
This condition (together with $C^0$ closed-ness) determine $g_1$, $g_2$ uniquely.
(Actually $u g_i$  is $C^1$ close to $u_{c(i)}$ since they are pseudo-holomorphic.
We use this $C^1$ close-ness to show the uniqueness.
This point was discussed before already.)

The map $u$ is $C^0$ close to both  $u_{c(1)} g_1^{-1}$  and $u_{c(2)}
g_2^{-1}$.

The obstruction spaces $E_i$ are vector spaces of smooth sections of $u_{c(i)}^* TX \otimes
\Lambda^{01}$.
$g_i^{*}$ transforms them to a vector space of sections of
$(u_{c(i)} g_i^{-1})^* TX \otimes \Lambda^{01}$,
which we write $g_i^{*}E_i$.

We use parallel transport to send $g_i^{*}E_i$ to
a vector subspace of sections of
$u^*TX \otimes \Lambda^{01}$.

In case we want to construct a Kuranishi chart centered at $u$
we proceed as follows.
We choose $D_{2}, D_{-2}$ codimension 2 submanifolds of $X$ which
intersect transversally with $u$ at $u(2)$ and $u(-2)$ respectively.

We consider $u'$ which is a smooth map $S^2 \to X$ that is $C^{10}$ close to $u$
(Namely $\vert u-u'\vert_{C^{10}} < \epsilon_0$.). We do not assume $u'$ is
pseudo-holomorphic
and will transfer the obstruction bundles to $(u')^*TX \otimes \Lambda^{01}$.

We use four more marked points $w_{c(i),1}$, $w_{c(i),-1}$.
We assume $w_{c(i),1} \in S^2$ is $\epsilon_{0}$-close to $g_i(1)$ and
$w_{c(i),-1} \in S^2$ is $\epsilon_{0}$-close to $g_i(-1)$.
\footnote{ $w_{c(1),1} = w_{c(2),1}$ may occur. But it does not
cause problem. See the paragraph starting at the end of
page 55. (The paragraph start with `A technical point to take care of is ...'.}
($\epsilon_0$ is a small number depending on $u$.)

There exists $g_i' \in Aut(S^2,0)$  such that
$g'_i(1) = w_{c(i),1}$, $g'_i(-1) = w_{c(i),-1}$.
Such $g_i'$ is unique and $g'_i - g_i$ is small.
(More precisely it is estimated by $o(\epsilon_0\vert \epsilon_{c(i)})$.
Here $o(\epsilon_0\vert \epsilon_{c(i)})$ depends on $\epsilon_0, \epsilon_{c(i)}$
and goes to zero as $\epsilon_0$ goes to zero for each fixed $\epsilon_{c(i)}$.)

We use $g'_0$, $g'_1$ (which is determined by $u', w_{c(i),\pm 1}$ in the above way)
to obtain $E_i(u';\vec w)$.
(Here $\vec w'$ is 4 points $(w_{c(i),\pm 1})$.)
Namely we use parallel transport from $(g'_i)^{*}E_i$ to  $(u')^*TX \times \Lambda^{01}$, that is possible by
$C^0$ close-ness.
\footnote{$u$ is $\epsilon_{c_(i)}$ close to $u_{c_i}g_{c_i}$. $u'$ is $\epsilon_0$
close to $u$, and $g'_i$ is $o(\epsilon_0\vert \epsilon_{c(i)})$.
All are $C^0$ sense. We first choose $\epsilon_{c(i)}$ small
and then choose $\epsilon_0$ so that the sum
$\epsilon_{c(i)} + \epsilon_0 + o(\epsilon_0\vert \epsilon_{c(i)})$
is small enough.}

We may assume $E_1(u';\vec w) \cap E_2(u';\vec w) = \{0\}$, since we may
assume
$$
E_1(u;(g_i(\pm1))) \cap E_2(u;(g_i(\pm1))) = \{0\}
$$
by
Lemma \ref{transbetweenEs}
and $E_1(u';\vec w) \cap E_2(u';\vec w) = \{0\}$ is an open condition
(with respect to $C^0$ topology.)
(Note in case $u' = u$ and $\vec w = (g_i(\pm1))$ we have $g'_i = g_i$.)

Now we apply the implicit function theorem
to the equation

\begin{equation}\label{mainaquation34}
\overline{\partial} u' \equiv 0 \mod E_1(u';\vec w) \oplus E_2(u';\vec w)
\end{equation}
\noindent
to obtain the thickened moduli space.
(It is the set of $(u',\vec w)$ satisfying this equation.)
It is a smooth manifold.
\footnote{The surjectivity of the linearized equation is OK,
since it is OK at $u$ and we may choose $\epsilon_0$ depending on $u$.}
\footnote{In case $[u] \in \mathcal M_1(A)$ is a stable map with
nontrivial automorphism then it is an orbifold.}
The dimension of the thickened moduli space is greater
than the correct dimension of Kuranishi neighborhood.
The difference is 12. ($4 = \dim Aut(S^2,0)$ and each $w_{i,\pm 1}$
gives 2 extra dimension.)\footnote{We correct misprint in our google post here.}

We cut down the dimension of the thickened moduli space by imposing the constraints

\begin{equation}\label{trans}
\aligned
&u'(2) \in D_2, \quad u'(-2) \in D_{-2},\\
&u'(w_{i,1}) \in D_{i,1}, \quad u'(w_{i,-1}) \in D_{i,-1}.
\endaligned
\end{equation}

The condition in the first line kills the $ Aut(S^2,0)$ ambiguity of $u'$
and the condition in the second line
kills the ambiguity of the choice of $\vec w$.\footnote{We correct misprint in our google post here.}

These are 12 independent equations.
($12 = 6 \times \text{codim} D$). After cutting down the dimension of the
thickened moduli space by (\ref{trans}) we obtain the required Kuranishi neighborhood.
(Transversality of the equation (\ref{trans}) is OK
by choosing $\epsilon_0$ small, since it is OK at $u$.)
\par\medskip
This PDF file was an answer to a question of K.Wehrheim.
Then McDuff asked a question about the smoothness of the
right hand side of (\ref{mainaquation34}).

The next part is a reproduction of our answer to it
\par\medskip
{\bf [The post on August 12]}
\par
The equation
$$
\overline{\partial} u' \equiv 0  \mod E_1(u';\vec w) \oplus E_2(u';\vec w)
$$
is an elliptic PDE.
More precisely it is a family of elliptic PDE with parameter $\vec w$
and extra  $m = \dim E_1 + \dim E_2$
parameters $a_i$ we explain below. We rewrite the equation to
$$
\overline{\partial} u' + \sum a_i e_i(u',\vec w) = 0
$$
where $e_i(u',\vec w)$ is a basis of $E_1(u';\vec w) \oplus E_2(u';\vec w)$.
(This is an equation for $u'$. Its parameters  are
$\vec w$ and $a_i$.)

The coefficient of this elliptic PDE depends smoothly
on the parameters $\vec w$, $a_i$.
(See the argument below which shows that
$e_i(u',\vec w)$ is smooth both in $u'$ and $\vec w$.)

So the solution space with $\vec w$ (and $a_i$) moving consists of a
smooth manifold (if the surjectivity of its linearized equation is OK).
Moreover the projection to the parameter space (especially $\vec w$)
and all the evaluation maps are smooth on the solution space.
This is a standard fact in the theory of elliptic PDE.
\par\medskip
Finally let us explain the smoothness of $e_i(u',\vec w)$
with respect to $u'$ and $\vec w$.
Note $e_i(u',\vec w)$ is a member of a basis of
$E_j(u;\vec w)$ for $j=1,2$.
The section $e_i(u',\vec w)$ is defined from a basis $e_i$ of $E_j$
in a way explained before.
We explain it again below in a slightly different way so that
the smoothness of $e_i(u',\vec w)$ becomes obvious.
Hereafter we consider the case $j=1$.

We put $v = u_{c(1)}g_1^{-1} :  S^2 \to X$.
(Note $u'$ is close to $v$ is $C^0$ sense.)
Let us take an open set $\mathcal V$ of $S^2 \times X$
that is
$$
\mathcal V = \{ (z,x) \mid d(v(z),x) < 2\epsilon_0 + 2\epsilon_{c(1)}\}.
$$
We define a vector bundle on $\mathcal V$ by
$$
\frak E = \pi_2^*TX \otimes \pi_1^* \Lambda^{01}.
$$
($\pi_1$,$\pi_2$ are projections from $S^2 \times X$
to the  first and second factors.)
Let $G$ be a small neighborhood of $g_1$ in $Aut(S^2,0)$.
We pull back the bundle $\frak E$ to $\mathcal V \times G$
by the projection $\mathcal V \times G \to \mathcal V$.
We denote it by $\frak E \times G$.
Let us define a smooth section $\hat e_i$ of it as follows.
Let $g'_1 \in G$. We consider the composition
$u_{c(1)} \circ (g'_1)^{-1}$ and define
$$
 (u_{c(1)} \circ (g'_1)^{-1})^+: S^2 \to \mathcal V
$$
by
$$
z \mapsto (z,(u_{c(1)} \circ (g'_1)^{-1})(z))
$$
We identify $S^2$ with the image $(u_{c(1)} \circ (g'_1)^{-1})^+(S^2)$.
The restriction of $\frak E \times G$ to this $S^2$ is
$(u_{c(1)} \circ (g'_1)^{-1})^*TX \times \Lambda^{01}$.
To this bundle we transform the section $e_i$
of $u_{c(1)}^*TX \times \Lambda^{01}$ by $(g'_1)^{*}$.
(This is what we did before.)
We then extend it to a section on $\mathcal V \times \{g'_1\}$ by the
parallel transportation in the $X$ direction along a geodesic
with respect to an appropriate connection of $X$.
\par
Thus for each $g'_1$ we have a section of $\frak E\times \{g'_1\}$ on
$\mathcal V\times \{g'_1\}$.
Moving $g'_1$, we have a section of the bundle $\frak E \times G$
on $\mathcal V \times G$.
We denote it by $\hat e_i$.
It is obvious that $\hat e_i$ is smooth.

Let $W$ be the parameter space of $\vec w$
that is a finite dimensional manifold.
Let $g'_1(\vec w)$ be the biholomorphic map
sending $\pm 1$ to $w_{1,\pm 1}$.
It depends smoothly on $\vec w$.
We put
$$
\hat u' : S^2 \times W \to \mathcal V \times G
$$
by
$$
(z,\vec w) \mapsto (z,u'(z),g'_1(\vec w))
$$
\par
$e_i(u',w)$ coincides with the composition  $\hat e_i \circ \hat u'$.
So this is smooth in $u'$ and $\vec w$.

\subsection{A note about \cite[Section 12]{FOn}}
\label{sec12FO}
\par
This is a note related to
\cite[Remark 4.1.3]{MW1}.
This remark we think is related to
\cite[Lemma 12.24 and Proposition 12.25]{FOn}.
\par\medskip
Following \cite{MW1} we restrict our explanation
to the case of $\mathcal M_{1}(A)$,
the moduli space of pseudo-holomorphic sphere with one marked point and
of homology class $A$ such that there are no nonconstant
pseudo-holomorphic spheres
$u_1,u_2$ with $u_{1*}([S^2])   + u_{2*}([S^2]) = A$.
(Therefore all the elements of $\mathcal M_{1}(A)$ are somewhere injective.)
\par
Let us describe the point which we think is the concern of  \cite[Remark 4.1.3]{MW1}.
We put $G = \text{Aut}(S^2,0)$ the group of biholomorphic maps
$v : S^2 \to S^2$ such that $v(0) = 0$.
Let
$u : S^2 \to X$ be a pseudo-holomorphic map of homology class $A$.
We are going to construct a Kuranishi neighborhood centered at $[u]
\in \mathcal M_{1}(A)$.
\par
We take $E$ that is a finite dimensional space of
smooth sections of $u^*TX \otimes \Lambda^{01}$.
We assume that
the image of
$$
D_u\overline\partial : \Gamma(S^2;u^*TX) \to \Gamma(S^2;u^*TX\otimes \Lambda^{01})
$$
together with $E$ is $\Gamma(S^2;u^*TX\otimes \Lambda^{01})$.
\par
We consider the operator
\begin{equation}\label{01}
\overline {D_u}\overline\partial : \Gamma(S^2;u^*TX) \to \Gamma(S^2;u^*TX\otimes \Lambda^{01})/E.
\end{equation}
\par
Let $u' : S^2 \to X$ be a map which is $C^{10}$ close to $u$.
We define $E(u') \in \Gamma(S^2;(u')^*TX\otimes \Lambda^{01})$ by parallel transport from $u$.
(More precisely we use parallel transport of the tangent bundle $TX$ along the geodesic
joining $u(z)$ and $u'(z)$ for each $z \in S^2$.)
\par
This map $u' \mapsto E(u')$ is {\it not} invariant of $G$ action.
\par
Let $V$ be the set of the solution of
\begin{equation}\label{Enonequiveqqq}
\overline\partial u'  \equiv 0 \mod E(u'),
\end{equation}
such that $u'$ is $\epsilon$-close to $u$ in $C^{10}$ norm.
Implicit function theorem and assumption implies that
$V$ is a smooth manifold.
\par
On $V$ we have a vector bundle whose fiber at $u'$ is $E(u')$.
This is a smooth vector bundle.
We have a section $s$ of it such that $s(u') = \overline\partial u'$.
\par
$s^{-1}(0)$ maps to $\mathcal M_{1}(A)$.
However this is not injective since $G$ action is not killed.
\par
Let us identify $T_uV$ with the kernel of (\ref{01}).
\par
We consider the Lie algebra $T_eG$.
Since $u\circ g \in V$ for all $g$ we can embed $T_eG \subset T_uV$.
Let $V'$ be a submanifold of $V$ such that $u\in V'$ and
\begin{equation}\label{3502}
T_eV' \oplus T_eG = T_eV.
\end{equation}
This $V'$ is the slice appearing in \cite[Section 12]{FOn}.
(Note this slice is taken after we obtained a finite dimensional
space.)
We prove the following:
\begin{prop}\label{prof11aaa}
$u' \mapsto [u']$ induces a homeomorphism between
$
(V' \cap B_{\epsilon}V) \cap s^{-1}(0)
$
and a neighborhood of $[u]$ in $\mathcal M_{1}(A)$
for sufficiently small $\epsilon$.
\end{prop}
This proposition is a special case of \cite[Lemma 12.24 and Proposition 12.25]{FOn}.
\begin{proof}
Let $G_0$ be a small neighborhood of identity in $G$.
(We will describe how small it is later.)
Let $L^2_m(S^2,X)$ be the space of $L^2_m$ maps
to $X$ from $S^2$. We take $m$ huge.
$L^2_m(S^2,X)$ is a Hilbert manifold and
$V$ is its smooth submanifold of finite dimension.
Let $N_VL^2_m(S^2,X)$ be a tubular neighborhood of $V$ in
$L^2_m(S^2,X)$ and $\Pi : N_V(L^2_m(S^2,X)) \to V$  the projection.
Let $V_0$ be a relatively compact neighborhood of $u$ in $V$.
We take $G_0$ small such that if
$u'\in V_0$ and $g \in G_0$ then $u'\circ g \in N_VL^2_m(S^2,X)$.
We define
$$
F : V_0 \times G_0 \to V
$$
by
$$
F(u',g) = \Pi(u'\circ g).
$$
We remark that if $s(u') = 0$ then $F(u',g) = u'\circ g$.
$F$ is a smooth map since $V_0$ consists of smooth maps.
In fact  the map
$$
V_0 \times G_0 \to L^2_m(S^2,X)
$$
defined by $(u,g)\mapsto u\circ g$ is smooth.

\begin{lem}\label{lemma2adddd}
There exists a neighborhood $V_2$ of $u$ in $V$ with the following property.
If $u' \in V_2$ there exists $g \in G_0$ such that
$F(u',g) \in V'$.
\end{lem}
\begin{lem}\label{lemma3adddd}
There exist $\epsilon$ and $G_0$ such that the following holds.
If $u' \in V'  \cap B_{\epsilon}V$ and
$F(u',g) \in V' $, then $g = 1$.
\end{lem}
\begin{proof}[Proof of lemmas \ref{lemma2adddd},\ref{lemma3adddd}]
For each $u' \in V_0$ we put
$$
G_0u' = \{ F(u',g)  \mid g \in G_0\}.
$$
(Note  $u,g \mapsto F(u,g)$ is not a group action. So $G_0u'$ is not an orbit.)
\par
We may replace $G_0$ by a smaller neighborhood of identity and take a small neighborhood $V_{00}$ of $u$
so that $g \mapsto F(u',g)$ is a smooth embedding of $G_0$ if
$u' \in V_{00}$.
By assumption (\ref{3502}) the submanifold $G_0u$ intersects transversally with $V'$ at $u$ in $V$. So we may replace $G_0$
by a smaller neighborhood of identity again such that
$G_0u \cap V' = \{u\}$.
\par
Now since $u' \mapsto G_0u'$ is a smooth family of smooth submanifolds,
$G_0u'$ intersects transversally to $V'$ at one point if $u'$ is
sufficiently close to $u$. This implies Lemmas \ref{lemma2adddd},\ref{lemma3adddd}
\end{proof}

If $[u'']$ is in a small neighborhood of $[u]$ in
$\mathcal M_1(A)$ then by definition we may replace $u''$ and assume
$u'' \in V_2$. Then  there exists $g\in G_0$ such that $u' = F(u'',g)$ is in $V'$
by Lemma \ref{lemma2adddd}.
Since $u''$ is pseudo-holomorphic $u' = u''\circ g$.
Namely $[u'] = [u'']$.
Thus  $
(V' \cap B_{\epsilon}V) \cap s^{-1}(0)
$ goes to a neighborhood of $[u]$.
\par
The injectivity of this map is immediate from Lemma \ref{lemma3adddd}.
\par
It is easy to see that this map from
$(V' \cap B_{\epsilon}V) \cap s^{-1}(0)$ to an open set of $\mathcal M_1(A)$
is continuous.
The proof that it is an open mapping is the same as the proof of the fact that $
(V' \cap B_{\epsilon}V) \cap s^{-1}(0)
$ is a neighborhood of $[u]$.
\end{proof}
Note the reason why Proposition \ref{prof11aaa} is not completely
trivial lies on the fact that equation (\ref{Enonequiveqqq}) is not
$G$ invariant.
In \cite[Section 15]{FOn}, we made a difference choice of $E(u')$ using
center of mass technique so that $u\mapsto E(u)$ is
$G$ equivariant. Then Proposition \ref{prof11aaa} is trivial to prove for that
choice of $E(u')$.

\section{Confutations against the criticisms made in \cite{MW1} on the
foundation of the virtual fundamental chain/cycle technique}
\label{HowMWiswrong}

In an article \cite{MW1} there are several criticisms
on the earlier references on the foundation of
virtual fundamental chain/cycle technique.
This section provides our confutations against those criticisms.
\par
We think such confutations are necessary by the
following reason. There are various researches in progress
based on virtual fundamental chain/cycle technique by various people.
The authors of \cite{MW1} mention a plan to write a replacement
of a part of the results in the existing literature.

However, \cite{MW1} provides only very beginning of their plan
and \cite{MW1} concerns only the case where the pseudo-holomorphic curves discussed are
automatically somewhere injective, which is not applicable to any of the on-going
researches we mentioned above.

Based on their writing, it appears that the authors of \cite{MW1} do not plan to
study the chain level argument based on Kuranishi structure.
See \cite[Page 5 Line 12-14 from bottom]{MW1}.
Since the chain level argument is essential in many of the on-going researches,
we need something more than their planned `replacement'.

Therefore leaving the criticisms of \cite{MW1}
un-refuted would cause serious confusion among the researchers in the relevant field.
Hence we have decided to provide our confutations against the criticisms
displayed in \cite{MW1} as public as possible to the degree of the article \cite{MW1}.
\par
It is our understanding that an important and basic agreement of
the mathematical research is that researchers are free to use
the results of the published research papers
(with appropriate citations) unless explicit and specified gap or
problem was pointed out to the author but
the author failed to provide a reasonable answer or correction
on that particular point. This agreement is a part of the foundation of the refereeing system
of the mathematical publications on which the whole mathematical community much depends.
\par
By this reason we explain in detail where the misunderstanding behind the criticisms
of \cite{MW1} lies in, and then make our confutations against them word by word.
We do so only to the criticisms of \cite{MW1}
directed against the papers written by the present authors.
In \cite{MW1} there are criticisms to
other versions of virtual fundamental cycle/chain technique.
We found several problems there also.
However we restrict our confutations only to the criticisms
directed to the papers of the present authors.
This is because the present authors do not
have thorough knowledge of the other versions of virtual fundamental
cycle/chain technique, compared to that of their own.

The page numbers etc. below are those of the version of \cite{MW1}
appeared as arXiv:1280.1340v1 on Aug. 7 2012.
We note that the date Aug. 7 was after
we had posted all our detailed answers
\cite{Fu1}, \cite{FOn2}, \cite{fooo:ans3} \cite{fooo:ans34} and \cite{fooo:ans5} to
K. Wehrheim's questions.
The small letters are used for the quote from \cite{MW1}.
\par
Note similar criticisms appear repeatedly in \cite{MW1}. In such a case we repeat the
same answer.
Although some of the quotations below may not be direct criticisms, we supplement them in order to clarify
our mathematical points.

\begin{enumerate}
\item\label{MW1} Page 2 Line 13-14
\par
{\small
while some topological and analytic
issues were not resolved}
\par
We will clarify below that all such issues have been resolved.
\item\label{MW2}
Page 2 Line 11-14 from the bottom and
\par
{\small
The main analytic issue in each regularization approach is in the construction of
transition maps for a given moduli space, where one has to deal with the lack of
differentiability of the reparametrization action on infinite dimensional function spaces
discussed in Section 3.}
\par
Such issue does not cause any problem in our proof as we will explain
below. See items (\ref{MW3}),(\ref{MW13}),(\ref{MW18}),(\ref{MW21})-(\ref{MW28}).
\item\label{MW3}
Page 4 Line 16-20 from the bottom
\par
{\small
The issue here is the lack of differentiability of the reparametrization action,
which enters if constructions are transferred between different infinite dimensional local
slices of the action, or if a differentiable Banach manifold structure on a quotient
space of maps by reparametrization is assumed.
}
\par
Such an issue is irrelevant to our approach and is not present in our approach, since
infinite dimensional slice was never used.
We explained this point in the second half of Subsection \ref{comparizon}.
\item\label{MW4}
Page 4 line 10 from below:
\par
{\small However, in making these constructions explicit,
we needed to resolve ambiguities in the definition
of Kuranishi structure, concerning the precise meaning of germ of coordinate
changes and the cocycle condition discussed in Section 2.5.}
\par
This point had been already corrected in \cite{fooo:book1}.
We explained this point in Subsection \ref{gernkuranishi}.
\item\label{MW5}
Page 4 last 4 lines:
\par
{\small One issue that we will touch on only briefly in Section 2.2
is the lack of smoothness of the standard gluing constructions,
which affects the smoothness of the Kuranishi charts near the nodal
or broken curves.
}
\par
This points had already been discussed in \cite{fooo:book1}.
More detail had been given in \cite{fooo:ans3}  and \cite{fooo:ans34}.
They are basically the same as Parts \ref{secsimple} and \ref{generalcase}
of this article.
\item\label{MW6}
Page 4 last line - Page 5 4th line:
\par
{\small
A more fundamental topological issue is the necessity to ensure that
the zero set of a transversal perturbation is not only smooth, but also
remains compact as well as Hausdorff, and does not acquire boundary
in the regularization.
These properties nowhere addressed in the literature as far as we are
aware, are crucial for obtaining a global triangulation and thus a fundamental
homological class.
}
\par
These points had been addressed in \cite{Fu1,FOn2}.
They had been sent to the authors of \cite{MW1} as a reply to the question
raised by the very person who wrote :
`These properties nowhere addressed in the literature as far as we are
aware'.
\par
According to the opinion of present authors, this point is not a fundamental issue but
only a technical point. We leave each reader to see our proof
given in Part \ref{Part2} and form his/her own opinion about it.
Anyway correctness of the construction of virtual fundamental chain/cycle
is not affected at all
whether this point (which had already been resolved) is fundamental or not.
\item\label{MW7}
Page 5 4 - 6 line:
\par
{\small
Another topological issue is the
necessity of refining the cover by Kuranishi charts to a `good cover'
in which the unidirectional transition maps allow for an iterative construction of
perturbation.}
\par
This was proved in \cite[Lemma 6.3]{FuOn99I}.
Responding to the request of an author of \cite{MW1} more detail
of its proof had been provided in \cite{Fu1,FOn2}.
\item\label{MW78} Page 5 Line 20-30
\par
{\small
The case of moduli spaces with boundary given as the fiber product of other moduli
spaces, as required for the construction of $A_{\infty}$-structures, is beyond the scope of our
project.
It has to solve the additional task of constructing regularizations that respect
the fiber product structure on the boundary. This issue, also known as constructing coherent
perturbations, has to be addressed separately in each special geometric setting,
and requires a hierarchy of moduli spaces which permits one to construct the regularizations
iteratively. In the construction of the Floer differential on a finitely generated
complex, such an iteration can be performed using an energy filtration thanks to the
algebraically simple gluing operation. In more `nonlinear' algebraically settings, such
as $A_{\infty}$ structures, one needs to artificially deform the gluing operation, e.g. when dealing
with homotopies of data [Se].}
\par
This paragraph might intend to put some negative view
on the existing literature which constructed $A_{\infty}$ structure,
especially on \cite{fooo:book1}.
However since no mathematical problem or difficulty
in the existing literature is
mentioned
it is impossible for us to do anything other than ignoring
this paragraph.
\item\label{MW8}
Page 5 the last paragraph - Page 6 the first paragraph:
\par
{\small
Another fundamental issue surfaced when we tried to understand how Floer's proof
of the Arnold conjecture is extended to general symplectic manifolds using abstract
regularization techniques. In the language of Kuranishi structures, it argues that a
Kuranishi space $X$ of virtual dimension $0$, on which $S^1$ acts such that the fixed points
$F\subset X$ are isolated solutions, allows for a perturbation whose zero set is given by the
fixed points. At that abstract level, the argument in both [FO] and [LiuT]\footnote{
They are \cite{FOn} and \cite{LiuTi98} in the reference of this article} is that
a Kuranishi structure on $(X\setminus F)/S^1$ (which has virtual dimension 􀀀$-1$) can be pulled
back to a Kuranishi structure on $X\setminus F$ with trivial fundamental cycle. However, they
give no details of this construction.}
\par
Such detail, \cite{fooo:ans5} = Part \ref{S1equivariant} of this article,
had been given to the authors of \cite{MW1}.  \cite{fooo:ans5} is a reply  to the question
raised by the very person who wrote :
`However, they give no details of this construction'.

\item\label{MW9}
Page 9 line 8-9:
\par
{\small
The gluing analysis is a highly nontrivial Newton iteration scheme and should have
an abstract framework that does not seem to be available at present.
}
\par
We do not understand why {\it abstract frame work} {\bf should} be used.
Our gluing analysis is based on the study of the concrete geometric
situation where we perform our gluing construction.
To find some abstract formulation of gluing construction
can be an interesting research. However
it is {\it not required} in order to confirm the correctness of the gluing constructions in
various particular geometric cases.
Such gluing constructions had been used successfully in
gauge theory and in pseudo-holomorphic curve
by many people in the last 30 years. We wonder
whether the authors of \cite{MW1} question the soundness of all
those well established results or not.
If not the authors of \cite{MW1} should explain the reason why
abstract frame work should be used in this particular case and not in the other cases.
\item\label{MW10}
Page 9 line 9-10:
\par
{\small
In particular,
it requires surjective linearized operators, and so only applies after perturbation.
}
\par
In the construction of Kuranishi neighborhood we
modify the (nonlinear Cauchy-Riemann) equation
$
\overline\partial u' = 0
$
slightly and use
$
\overline\partial u' \equiv 0 \mod E(u').
$
In other words the surjective linearized operators
are obtained by introducing obstruction space $E(u')$ and
{\it not} by perturbation.
In other words, surjectivity of linearized operators
is used {\it before} perturbation.
\par
The perturbation via a generic choice of multisections
starts {\it after} finite dimensional reduction.
\item
Page 9\label{item12} line 11-12:
\par
{\small
Moreover, gluing of noncompact spaces requires uniform quadratic estimates, which
do not hold in general.
}
\par
The proof of uniform quadratic estimates is certainly necessary in all the situations
where gluing (stretching the neck) argument are used.
(Both in gauge theory and the study of pseudo-holomorphic curve for example.)
The way to handle it had been well established  more than 20 years ago.
\footnote{For example in \cite[Blowing up the metric; page 121-127]{freedUhlen} an estimate of
$L^4$ norm in terms of $L^2_1$ norm is discussed, in the case of one forms on noncompact $4$ manifolds
with cylindrical end.}
It is proved in our situation (in the same way as many of the other situations) as follows.
The domain curve $\Sigma_T$
(we use the notation of Part \ref{secsimple}) is a union of core and neck region.
The core consists a compact family of compact spaces and is independent of the
gluing parameter $T$.
Therefore `uniform quadratic estimates' is obvious there.
On the other hand, the length of the neck region is unbounded,
which is the noncompactness of \cite{MW1}'s concern we suppose.
However, on the neck region we have an exponential decay estimate
(See \cite[Lemma 11.2]{FOn} = Proposition \ref{neckaprioridecay} of this article.) and the
neck region is of cylindrical type. Therefore even the length of the neck region is
unbounded we have a {\it uniform}  quadratic estimates.
We also remark  that the weighted Sobolev norm
we introduced in (\ref{normformjula5}) and \cite[Subsection 7.1.2]{fooo:book1}  is designed so that
the norm of the right inverse of the linearized operator
becomes uniformly bounded. (Namely it is bounded by a constant
independent of $T$.)
\item\label{MW11}
Page 9 line 12-13:
\par
{\small
Finally, injectivity of the gluing map does not follow from
the Newton iteration and needs to be checked in each geometric setting.
}
\par
The classical proof of injectivity had been reviewed \cite{fooo:ans3} =
Part \ref{secsimple}. See Section \ref{surjinj}.
Maybe this proof is more popular in gauge theory community than in
symplectic geometry community.
(See \cite{Don86I} for example.)
\item\label{MW12}
Page 9 line 3-7:
\par
{\small
Each setting requires a different, precise definition of a Banach space of perturbations.
Note in particular that spaces of maps with compact support in a given
open set are not compact. The proof of transversality of the universal section is very
sensitive to the specific geometric setting, and in the case of varying $J$ requires each
holomorphic map to have suitable injectivity properties.}
\par
It seems that `a Banach space of perturbations' that they allude to here is the obstruction space
$E_{\frak p}(u')$. We always choose it so that it becomes a {\it finite} dimensional space of smooth sections.
(See for example \cite[(12.7.4)]{FOn} and Definition \ref{obbundeldata} (5) of present article.)
It is a finite dimensional space and so is  a `Banach space'.
(We do not think it is a good idea to introduce infinite dimensional space here.
It is because such infinite dimensional space is harder to control.)
It is very hard to understand the rest of this writing at least for us.
The word `universal section' is not defined, for example.
\item\label{MW13}
Page 11 lines 2-5 from the bottom.
\par
{\small
This differentiability issue was not apparent in [FO, LiT]\footnote{[FO] is \cite{FOn}
and [LiT] is \cite{LiTi98} in the reference of this article.}, but we encounter the same obstacle in
the construction of sum charts; see Section 4.2. In this setting, it can be overcome by working with
special obstruction bundles, as we outline in Section 4.3.}
\par
This `obstacle' seems to be related to the following point:
To associate an obstruction space $E_{\frak p}(u')$
for $u' : \Sigma' \to X$
we use the vector space
$E_c$ for various $c$, that is a space of sections $u_c^*TX \otimes \Lambda^{01}$
where $u_c : \Sigma_c \to X$. So
we need to use a diffeomorphism between a support $K_c$ of the elements of
$E_c$ and a subset of $\Sigma'$.
Namely we transform $E_c$ by using this diffeomorphism and a parallel transport on $X$.
\par
Note in \cite{FOn} as well as in all our articles,
the subspace $E_c$ is a finite dimensional space of
smooth sections. So clearly no issue appears by transforming them
by diffeomorphism.
It seems that `working with special obstruction bundles'
is nothing but this choice of \cite[(12.7.4)]{FOn} etc., that is, `a finite dimensional space of
smooth sections'.

\item\label{MW14}
Page 12 lines 16-20.
\par
{\small
In principle, the construction of a continuous gluing map should always be possible along the lines of \cite{McSa94},
though establishing the quadratic estimates is nontrivial in each setting. However, additional arguments
specific to each setting are needed to prove surjectivity, injectivity, and openness of the gluing map.}
\par
The method of \cite{McSa94} was used in \cite{FOn}.\footnote{It is quoted in \cite[page 984]{FOn} as follows:
\par
The proof is again a copy of McDuff-Salamon's in [47] with some minor modifications to handle the existence of the obstruction and moduli parameter.
\par [47] is \cite{McSa94} in the reference of
this article.}
The quadratic decay estimate was classical as we already mentioned in
 (\ref{item12}).
Because of our choice of `special obstruction bundles'  (\ref{MW13})),
there appears no new point arising from the addition of obstruction bundles.\footnote{We choose our obstruction bundle so that the
supports of its elements are away from nodes. So the obstruction bundle
affects our equation only at the core where noncompactness of the
source does not appear. At the neck region the equation is genuine
pseudo-holomorphic curve equation and is the same as one studied in
 \cite{McSa94}.}
The proof of surjectivity is in \cite[ Section 14]{FOn}.
More detail is in Section \ref{surjinj} of this article = \cite[Section 1.5]{fooo:ans3}, together
with the proof of injectivity (that is similar to one of surjectivity)
and openness of the gluing map (that follows from surjectivity
proven in  Section \ref{surjinj} as shown in Section \ref{cutting} (Lemma \ref{setisopen}
= \cite[Lemma 2.105]{fooo:ans34})).
Those proofs were classical already at the stage of 1996.

\item\label{MW15}  Page 12 lines 20-23.
\par
{\small Moreover, while homeomorphisms to their image suffice for the geometric regularization approach,
the virtual regularization approaches all require stronger differentiability of the gluing map; e.g. smoothness in
[FO,FOOO,J]\footnote{[FO,FOOO,J] are
\cite{FOn}, \cite{fooo:book1}  \cite{joyce} in the reference of this article.
}
}
\par
For the purpose of \cite{FOn} differentiability of the gluing map is required only in
its version with $T$ (the gluing parameter) fixed, as we explained in
Subsection \ref{subsec342}.
\par
The `stronger differentiability of the gluing map' had been
proved in \cite[Proposition A1.56]{fooo:book1}. More detail  of this
proof \cite{fooo:ans3} \cite{fooo:ans34}
(= Part \ref{secsimple} and \ref{generalcase} of this article) had been written
and sent to the members of the google group
`Kuranishi', which include the authors of \cite{MW1}.
\item\label{MW16}  Page 12  Remark 2.2.3
\par
{\small
None of [LiT, LiuT, FO, FOOO]
\footnote{They are \cite{LiTi98,LiuTi98,FOn,fooo:book1}
in the reference of this article.} give all details for the construction of a gluing map. In particular, [FO, FOOO]\footnote{They are \cite{FOn,fooo:book1}
in the reference of this article.} construct gluing maps with image in a space of maps,
but give few details on the induced map to the quotient space, see Remark 4.1.3. For closed nodal curves,
[McS, Chapter 10] constructs continuous gluing maps in full detail, but (at least in the first edition)
does not claim that the glued curves depend differentiably on the gluing parameter $a\in \C$ as $a \to  0$.
By rescaling $\vert a \vert$, it is possible to establish more differentiability for $a \to  0$.
}
\par
The analytic detail of the gluing map had been given in
\cite[Sections 7.1 and Sections A1.4]{fooo:book1}.
Even more detail of the same argument
(\cite{fooo:ans3} and \cite{fooo:ans34} that are Parts \ref{secsimple} and \ref{generalcase} of
this article) was sent
to the members of the google `Kuranishi' that include the authors of \cite{MW1}.
\par
In \cite[Remark 4.1.3]{MW1}, the authors of \cite{MW1} explained why {\bf their} approach
fails in the setting of \cite{FOn} Section 12-15.
Therefore it is irrelevant to {\bf our} approach.
See Items (\ref{MW29})-(\ref{MW293}) and Subsection \ref{sec12FO}.
\par
The book \cite{McSa94} which we quote in \cite{FOn} indeed does not claim
the differentiability on the gluing parameter.
As we explained in Subsection \ref{subsec342},
the differentiability of the gluing map with gluing parameter fixed,
is enough to establish all the results of \cite{FOn}.
This is because we only need to study the moduli space of virtual
dimension 1 and 0 for the purpose of \cite{FOn}.
We remark that this had been mentioned already in
\cite[page 782 line 6-8 from below]{fooo:book1} as follows.
\par\medskip
We remark that smoothness of coordinate change was not used in [FuOn99II]
\footnote{This is \cite{FOn} in the reference of this article.}, since only 0 and 1 dimensional
moduli spaces was used there. In other words, in Situation 7.2.2 mentioned in \S
7.2, we do not need it.
\par\medskip

\item\label{MW17} Page 16 - 19, Beginning of Subsection 2.5
\par
There are discussions about germ of Kuranishi neighborhood and its coordinate
change there. As we explained in Subsection  \ref{gernkuranishi}
this point had already been corrected in \cite{fooo:book1}.

\item\label{MW18} Page 19 lines 13-18.
\par
{\small
The first nontrivial step is to make sure that these representatives were chosen sufficiently small
for coordinate changes between them to exist in the given germs of coordinate changes. The second
crucial step is to make specific choices of representatives of the coordinate changes such that the
cocycle condition is satisfied. However, [FO, (6.19.4)] does not address the need to choose specific,
rather than just sufficiently small, representatives.
}
\par
This point is related to the notion of `germs of Kuranishi structure'.
So it had been already corrected in \cite{fooo:book1}.
The definition of Kuranishi structure in  \cite{fooo:book1} includes
the choice of representatives of the coordinate changes such that the cocycle condition is satisfied.

\item\label{MW19} Page 20 lines 18-21
\par
{\small
The basic issues in any regularization are that we need to make sense of the equivalence relation
and ensure that the zero set of a transverse perturbation is not just locally smooth (and hence can
be triangulated locally), but also that the transition data glues these local charts to a compact
Hausdorff space without boundary. }
\par
See item (6).

\item\label{MW20} Page 25 last two lines - Page 26 first line
\par
{\small
It has been the common understanding that by stabilizing the domain or working in finite dimensional
reductions one can overcome this differentiability failure in more
general situations. }
\par
This common understanding is absolutely correct.
(See Items (\ref{MW21}) - (\ref{MW26}).)
Indeed we had done so.

\item\label{MW21} Page 27 last paragraph - Page 28 second line.
\par
{\small
It has been the common understanding that virtual regularization techniques deal with the differentiability failure of the reparametrization action by working in finite dimensional reductions, in which the action is smooth. We will explain below for the global obstruction bundle approach, and in Section 4.2 for the Kuranishi structure approach, that the action on infinite dimensional spaces nevertheless needs to be dealt with in establishing compatibility of the local finite dimensional reductions. In fact, as we show in Section 4, the existence of a consistent set of such finite dimensional reductions with finite isotropy groups for a Fredholm section that is equivariant under a nondifferentiable group action is highly nontrivial. For most holomorphic curve moduli spaces, even the existence of not necessarily compatible reductions relies heavily on the fact that, despite the differentiability failure, the action of the reparametrization groups generally do have local slices. However, these do not result from a general slice
construction for Lie group actions on a Banach manifold, but from an explicit geometric construction using transverse slicing conditions.}
\par
The `highly nontrivial problem' mentioned in
`that is equivariant under a nondifferentiable group action is highly nontrivial'
in the above quote seems to be related to the issue written in Item (\ref{MW13}).
It had been resolved in the way explained there.
\par
The discussion of \cite[Section 4.2]{MW1} seems to be very much similar
to a special case of our discussion in Part \ref{generalcase} of this
article (= \cite{fooo:ans34}).
(See Item (\ref{MW27}).)
So it seems to us that \cite[Section 4.2]{MW1} also {\it supports} the
common understanding `stabilizing the domain or working in finite dimensional reductions one can overcome this differentiability failure in more
general situations' as Part \ref{generalcase} of this
article does.
\par
As we explained in Subsection \ref{comparizon}
we never used slice theorem for such
`general slice
construction for Lie group actions on a Banach manifold'
and had been used `explicit geometric construction using transverse slicing conditions'.
\par
{So the description of this part of \cite{MW1} is based on their misunderstanding of the Kuranishi
structure approach and presumption arising from their experience with other on-going projects}.
We emphasize that the action on infinite dimensional spaces
{\it never} enter in our approach.
\item\label{MW22}
Page 32 16 and 17 lines from bottom.
\par
{\small
Using such slices, the differentiability issue of reparametrizations still appears in many guises:}
\par
Let us explain below (Items (\ref{MW23}) - (\ref{MW26})) why {\it none} of those guises appear in our construction.
\item\label{MW23}
{\small
(i) The transition maps between different local slices - arising from different choices of fixed
marked points or auxiliary hypersurfaces are reparametrizations by biholomorphisms that vary with
the marked points or the maps. The same holds for local slices arising from different reference surfaces,
unless the two families of diffeomorphisms to the reference surface are related by a fixed diffeomorphism, and thus fit into a single slice.}
\par
The family of diffeomorphisms appearing here
is applied (as reparametrization) to the set of solutions of elliptic PDE. (Namely
after solving equations, $\overline{\partial}u' \equiv 0 \mod E_{\frak p}(u')$.)
By elliptic regularity they are smooth families of smooth maps.
So reparametrization does not cause any problem.
\item\label{MW24}
{
\small
(ii) A local chart for $\frak R$ near a nodal domain is constructed by gluing the components of the nodal
domain to obtain regular domains. Transferring maps from the nodal domain to the nearby regular domains
involves reparametrizations of the maps that vary with the gluing parameters.}
\par
Near a nodal domain the smoothness of coordinate change is
more nontrivial than the case (i) above.
\cite[Lemma A1.59]{fooo:book1} (which is generalized to Propositions
 \ref{changeinfcoorprop}  and \ref{reparaexpest} of this article
= \cite[Propositions 2.19,2.23]{fooo:ans34})
had been prepared for this purpose and had been used to resolve this point.
See the proof of Lemma \ref{2120lem} for example.
\item\label{MW25}
{\small
(iii) The transition map between a local chart near a nodal domain and a local slice of regular domains is
given by varying reparametrizations. This happens because the local chart produces a family of Riemann surfaces
that varies with gluing parameters, whereas the local slice has a fixed reference surface.
}
\par
The same answer as item (\ref{MW24}) above.
\item\label{MW26}
{
\small
(iv) Infinite automorphism groups act on unstable components of nodal domains.
}
\par
This is the reason why we need to add marked points to such components.
\item\label{MW27} Page 33 Lines 6-8
\par
{\small
We show in Remark 3.1.5 and Section 4.2 that these issues are highly nontrivial to deal with in abstract regularization approaches.}
\par
It is not so clear for us what `abstract regularization approaches' means.
It might be related to taking slice of infinite dimensional group.
We never do it as explained in  Subsection \ref{comparizon}.
\par
What is written in Section 4.2 \cite{MW1} is similar
to (a special case of) what we had written in \cite[Appendix] {FOn}, \cite[page 8-9, the answer to
Question 4]{Fu1}, \cite{fooo:ans34}
( = Part \ref{generalcase} of this article), as we
show by examples below.
Therefore as far as the correctness of the mathematical statements appearing here concerns,
the opinion of the authors of \cite{MW1} are likely to
coincide with ours.
\begin{enumerate}
\item
The conditions given in page 39 Lines 5 -12 from the bottom.
\par
These are similar to the special case of Definition \ref{obbundeldata} of obstruction bundle data.
\item
The discussion of the last part of page 42 of \cite{MW1}  (Item 1 there).
\par
The discussion there seems to be related to the
smoothness of $u' \mapsto E(u')$ that was explained in the posts on  Aug. 12, which is
reproduced in Section \ref{smoothness}.
\item
The Item 2 in page 43.
\par
The first part of this discussion is related to Lemma \ref{transbetweenEs}.
The rest seems to be similar to the construction of the obstruction
bundle we described in `the post of Aug. 10', which we reproduced in Section \ref{smoothness}.
\footnote{
We remark those posts on Aug. 10 and Aug. 12 are extracts (or adaptations to
a special case) of our earlier posts  \cite{fooo:ans34}.
We sent them to the google group `Kuranishi'
according to the request of the authors of \cite{MW1}.}
\item Page 44 Line 6-8.
\par
{\small
This construction is so canonical that coordinate
changes between different sum charts exist essentially automatically, and satisfy the
weak cocycle condition.}
\par
This sentence is very much similar to the following
which appears as a part of \cite{Fu1}.
So there seems to be an agreement concerning this point.
\par\medskip
Once (1) is understood the coordinate change $\phi_{{\frak q}{\frak p}}$ is just a map which send an en element
$((\Sigma,\vec z),u)$ to the same element.
So the cocycle condition is fairly obvious.
\item Page 46 \cite[(4.3.3)]{MW1}.
\par
The choice of the obstruction space in \cite[(4.3.3)]{MW1} is the same as
one which appeared
in \cite[(12.7.4)]{FOn}.
\end{enumerate}
\item\label{MW29}
Page 38-39,
Remark 4.1.3
\par
In \cite[Remark 4.1.3]{MW1}, the authors of \cite{MW1} explained why {\it their} approach
fails in the setting of \cite{FOn} Section 12-15.
Therefore it is irrelevant to {\it our} approach.
We show it by several examples below.

\item\label{MW291}
Page 38 Line 16-17,
\par
{\small
The above proof translates the construction of basic Kuranishi charts
in [FO]\footnote{This is \cite{FOn} in the reference of this article.} in the absence of nodes and Deligne-Mumford parameters into a formal setup.
}
\par
The authors of \cite{FOn} do not agree that this is a translation of the construction of \cite{FOn}.

\item\label{MW292}
Page 38 last 4 lines and Page 39 first line.
\par
{\small
[FO] construct the maps $\hat s$ and $\hat{\psi}$  on a {\it thickened}
Kuranishi domain" analogous to $\hat W_{f}$ and thus need to make the same restriction to
an infinitesimal local slice" as in Lemma 4.1.2. Again, the argument for injectivity
of $f$ given in Lemma 4.1.2 does not apply due to the differentiability failure of the
reparametrization action of $G = G_{\infty}$ discussed in Section 3.1.}
\par
The authors of \cite{MW1} claim here that why the proof of \cite[Lemma 4.1.2]{MW1} {\it they gave} fails in the
situation of \cite[Section 12]{FOn}.
This claim has no relation to our proof.
See Subsection \ref{sec12FO} for the correct proof of \cite[Lemma 12.24]{FOn}
in the situation of \cite{MW1}.
\item\label{MW292}
Page 39 Line 6-8
\par
{\small
The claim that $\psi_f$ has open image in $\sigma^{-1}(0)/G$ is analogous to [FO, 12.25], which
seems to assume that $\hat U_f$ is invariant under $G_{\infty}$ to assert
``$\hat\psi(s􀀀^{-1}(0) \exp_f(W_f)) = \hat{\psi}(\hat s^{-1}(0))$''. }
\par
When first and 4th named authors wrote \cite{FOn}, they of course were aware of the fact that
the choice of the obstruction bundle $E$ in \cite[Section 12]{FOn} is not invariant of the action of
automorphism group. This is {\it the} reason why \cite[Section 15 and appendix]{FOn} was written.
\footnote{
We remark that the automorphism group (that is written $G$ in the above quote from \cite{MW1})
acts on the zero set of Kuranishi map (that is written as $\hat s^{-1}(0)$ in the above quote.).
Therefore the set $V'_{\sigma} \cap s_{\sigma}^{-1}(0)/\text{\rm Aut} (\sigma)$ (that is
$\hat\psi(s^{-1}(0))$ in the above quote, which appears in \cite[Lemma 12.24 and Proposition 12.25]{FOn} is well defined. (The group $\text{\rm Aut}(\sigma)$ is the
finite group that is a automorphisms of stable {\it map}.)
}
More explicitly it is written in \cite[Page 1001 Line 18-19 from the bottom]{FOn} that
\par\medskip
The trouble
here is that $E_{\tau_i}$
is {\it not} invariant by the ``action'' of $Lie(Aut(\Sigma_{\tau_i}))_0$.
\par\medskip
(Note `not' was italic in \cite{FOn}.)
The proof of \cite[Proposition 12.25]{FOn} in the situation of \cite{MW1} without using $G$ invariance of obstruction
bundle
is in Subsection \ref{sec12FO}.
\item\label{MW293} Page 39 Line 8-10.
\par
{\small
However,
$G_{\infty}$-invariance of $\hat U_f$ requires $G_{\infty}$-equivariance of $\hat E$, i.e. an equivariant
extension of $E_f$ to the infinite dimensional domain $\hat{\mathcal V}$.
A general construction
of such extensions does not exist due to the differentiability failure of the $G_{\infty}$-action.
}
\par
It is not clear for us what `general construction' means. It might mean
`a construction in some abstract setting without using the
properties of explicit geometric setting'.
We never tried to find such a `general construction'.
Two explicit geometric constructions of such extensions are given in \cite{FOn}.
One in Section 15 the other in Appendix.

\item\label{MW28} Page 39 Lines 14-17 from the bottom.
\par
{\small
The differentiability issues in the above abstract construction of Kuranishi charts
can be resolved, by using a geometrically explicit local slice $\mathcal B_f \subset
\widehat{\mathcal B}^{k,p}$ as in (3.1.3).
(This is mentioned in various places throughout the literature, e.g. [FO, Appendix], but
we could not find the analytic details discussed here.)
}
\par
Such detail had been written in
\cite{fooo:ans34} and posted to the google group `Kuranishi'
of which the authors of \cite{MW1} are members. \cite{fooo:ans34} is a reply to a question raised by the very person
who wrote `we could not find the analytic details'.

\item\label{MW30} Page 53 Line 14-17 from the bottom.
\par
{\small
Note that we crucially use the triviality of the isotropy groups, in particular
in the proof of the cocycle condition. Nontrivial isotropy groups cause additional indeterminacy,
which has to be dealt with in the abstract notion of Kuranishi structures.}
\par
The detail of the existence of Kuranishi structure
of the moduli space {\it without using triviality of the isotropy groups} is
explained in detail in Part \ref{generalcase}  of this article.
\item\label{MW31}
Page 54
\par
{\small
(iii) In view of Sum Condition II′ and the previous remark, one cannot expect any two given basic Kuranishi
charts to have summable obstruction bundles and hence be compatible. This requires a perturbation of the
basic Kuranishi charts, which is possible only when dealing with a compactified solution space,
since each perturbation may shrink the image of a chart.
}
\par
This might be related to Lemma \ref{transbetweenEs}. (The proof of this lemma is easy.)
\item\label{MW32} Page 54
\par
{\small
(iv) This discussion also shows that even a simple moduli space such as $\mathcal M_1(A,J)$
does not have a canonical Kuranishi structure. Hence the construction of invariants from this
space also involves constructing a Kuranishi structure on the product cobordism
$\mathcal M_1(A,J) \times [0,1]$ intertwining any two Kuranishi structures for $\mathcal M_1(A,J)$
arising from different choices of basic charts and transition data.}
\par
The cobordism of Kuranishi structure and its application
to the well-defined-ness of virtual fundamental class
had been discussed in \cite[Lemmas 17.8,17.9]{FOn2} etc.
\item\label{MW33}
Page 75 13 -1 5
\par
{\small
However, in order to obtain a VMC from a Kuranishi structure, we either need to require the strong
cocycle condition, or make an additional subtle shrinking construction as in (i) that crucially uses the additivity condition.
}
\par
A detailed explanation of the construction of VMC from a Kuranishi structure {\it without using additivity condition}
is given in Part \ref{Part2} of this article.

\item\label{MW34}
Page 75 Line 20-22 from the bottom.
\par
{\small
The proof of existence in [FO, Lemma 6.3] is still based on notions of germs and addresses neither
the relation to overlaps nor the cocycle condition. }
\par
It is explained in Subsection \ref{gernkuranishi} why
the proof of [FO, Lemma 6.3] is {\it not} based on the notion of germs.
The detail of this proof is given in Section \ref{sec:existenceofGCS}.

\end{enumerate}

\end{document}